\documentclass [11pt]{amsart}
\usepackage {amsmath, amssymb, amscd, mathrsfs, mathtools, graphicx, stmaryrd, enumerate, color, url,epsf}
\usepackage[all,knot,arc]{xy}
\usepackage{hyperref}
\usepackage[pdftex,lmargin=1in,rmargin=1in,tmargin=1in,bmargin=.95in]{geometry}
\SetSymbolFont{stmry}{bold}{U}{stmry}{m}{n}

\newtheorem {theorem}{Theorem}[section]
\newtheorem {lemma}[theorem]{Lemma}
\newtheorem {proposition}[theorem]{Proposition}
\newtheorem {corollary}[theorem]{Corollary}
\newtheorem {conjecture}[theorem]{Conjecture}
\newtheorem {definition}[theorem]{Definition}

\theoremstyle{remark}
\newtheorem {remark}[theorem]{Remark}
\newtheorem {example}[theorem]{Example}

\newcommand\eps{\varepsilon}

\renewcommand\th{^{\text{th}}}

\newcommand\Z{\mathbb{Z}}
\newcommand\R{\mathbb{R}}
\newcommand\Q{\mathbb{Q}}
\def\E{\mathbb{E}}

\def\H{\mathbb{H}}
\def\FF{\mathfrak{F}}
\newcommand \bH {\overline{\H}}
\newcommand\T{\mathbb{T}}
\newcommand\Ta{\mathbb{T}_\alpha}
\newcommand\Tb{\mathbb{T}_\beta}
\newcommand\Tg{\mathbb{T}_\gamma}
\newcommand\Td{\mathbb{T}_\delta}
\newcommand\Tbp{\mathbb{T}_{\beta'}}
\newcommand\Tap{\mathbb{T}_{\alpha'}}

\newcommand\Xs{\mathbb{X}}
\newcommand\Os{\mathbb{O}}

\newcommand\Rect{\mathrm{Rect}}

\def\Sym{\mathrm{Sym}}

\def\Os{\mathbb O}

\def\zz {{\mathbb{Z}}}
\def\rr {{\mathbb{R}}}
\def\cc {{\mathbb{C}}}
\def\qq {{\mathbb{Q}}}

\def\TT {{\mathcal{T}}}

\def\mo {{\mathcal{A}}}
\def\kk {{\mathbf{k}}}
\def\del {{\partial}}
\def\delt{{\mathfrak{d}}}
\def\Sphere {{\mathcal{S}}}
\def\BHyper {{\overline{\Hyper}}}
\def\bRing{\bar{\mathcal{R}}}
\def\bRingbig{\bar{\mathfrak{R}}}
\def\bHyper{\bar{\Hyper}}
\def\tt {{\mathfrak{t}}}
\def\hey {\! \sslash \!}
\def\ss {{\mathfrak{s}}}
\def\nx {{\mathfrak{x}}}
\def\ux {{\mathfrak{u}}}
\def\A {{\mathscr{A}}}
\def\V {{\mathcal{V}}}

\def\C {{\mathcal{C}}}
\def\bC {\bar C}
\def\tC{\tilde{C}}
\def\tD{\tilde D}
\def\zed {{\mathcal{Z}}}
\def\lk {{\operatorname{lk}}}
\def\spc {{\operatorname{Spin^c}}}
\newcommand\SpinC\spc

\def\fin\qedhere
\def\pr {{\text{pr}}}
\DeclareMathOperator{\id}{id}
\DeclareMathOperator{\End}{End}

\def\from {{\leftarrow}}

\def\M {\mathcal {M}}
\def\Nod {\mathcal {N}}

\def\Zs {\mathbb{Z}}
\def\Ws {\mathbb{W}}
\def\s{\mathbf s}
\def\t{\mathbf t}
\def\a {\mathbf a}
\def\x{\mathbf x}
\def\c {\mathbf c}

\def\bu{\bar{u}}
\def\y{\mathbf y}
\def\z{\mathbf z}
\def\tx{\tilde{\x}}
 
\def\vs {\mathbf v}

\def\Chain {\mathfrak A}
\def\Am {\mathbf A^{\! \!-}}
\def\Cc {\mathfrak{C}}
\def\Cint{C_{\operatorname{int}}}
\def\Ccint {\Cc_{\operatorname{int}}}
\def\Bain {\mathfrak B}
\def\Ring {\mathcal R}
\def\K {\mathcal K}
\def\orK {\vec K}
\def\half {\frac{1}{2}}
\newcommand\orL{\vec{L}}
\newcommand\orN{\vec{N}}
\newcommand\orM{\vec{M}}

\def\cs {\mathbf c}
\def\S{\mathbf S}
\newcommand\alphas{\mbox{\boldmath$\alpha$}}
\newcommand\betas{\mbox{\boldmath$\beta$}}
\newcommand\etas{\mbox{\boldmath$\eta$}}
\newcommand\deltas{\mbox{\boldmath$\delta$}}
\newcommand\gammas{\mbox{\boldmath$\gamma$}}
\newcommand\thetas{\mbox{\boldmath$\theta$}}
\newcommand\lambdas{\mbox{\boldmath$\lambda$}}

\def\interior{\mathrm {int}}
\def\ws{\mathbf w}
\def\zs{\mathbf z}

\def\b{\mathbf{b}}
\def\gr{\mathrm{gr}}
\def\Field{\mathbb F}
\newcommand\Filt{A}

\newcommand\EmptyRect{\Rect^\circ}

\def\play {\operatorname{pl}}

\def\Per {\mathcal{P}}

\def\zero{\mathbf{0}}
\def\dd {{\mathbf {d}}}
\def\N {(\zz_{\geq 0})}
\def\U {{\mathbf{U}}}
\def\uu {{\mathcal{U}}}
\def\ff {{\mathbb{F}}}
\def\F {{\mathcal{F}}}
\def\G {{G}}

\def\ns {\mathbf{n}}
\def\he {{\mathcal{H}}}
\def\Hyper {{\mathcal{H}}} 
\def\S {{\mathbf{S}}}
\def\TR {\mathscr{T}}
\def\SR {\mathscr{S}}

\def\bd {{\langle \langle\tilde \Lambda \rangle \rangle}}
\def\twist {{\{ \tilde \Lambda \}}}
\def\dtwist {{\{ \{ \tilde \Lambda \} \}}}
\def\p {{\mathbf{p}}}
\def\Pr {{\mathcal{I}}}
\def\I {\Pr}
\def\De {{D}}

\def\D {{\mathcal{D}}}
\def\CC {{\mathscr{C}}}
\def\DD {{\mathscr{D}}}
\def\Disk {{\mathbb{D}}}

\def\hyp {{\mathscr{H}}}

\def\ori {{\mathfrak{o}}}
\def\ttt {{\mathbf{t}}}
\def\tTheta{\tilde{\Theta}}
\def\can {{\operatorname{can}}}
\def\mix {{\operatorname{mix}}}
\def\tot {{\operatorname{tot}}}

\def\Re {{\operatorname{Re }}}
\def\dim {{\operatorname{dim}}}
\def\mod {\operatorname{mod}}
\def\ker {{\operatorname{Ker}}}
\def\im {{\operatorname{Im}}}
\def\coker{{\operatorname{Coker}}}
\def\Id  {{\operatorname{Id}}}

\def\lra {{\longrightarrow}}
\def\fin\qedhere
\def\from {{\leftarrow}}
\def\In {\operatorname{In}}
\def\Out {\operatorname{Out}}
\def\Destab{\operatorname{Destab}}

\newcommand{\leftexp}[2]{{\vphantom{#2}}^{#1}{#2}}
\newcommand\bHF{\mathbf{HF}}
\newcommand\HFm{\mathbf{HF}^-}
\newcommand\CFm{\mathbf{CF}^-}
\newcommand\HFinf{\mathbf{HF}^\infty}
\newcommand\CFinf{\mathbf{CF}^\infty}

\newcommand\CFLm{\mathbf{CFL}^-}
\newcommand\HFmd{\mathbf{HF}^{-, \delta}}
\newcommand\CFmd{\mathbf{CF}^{-, \delta}}
\newcommand\CFmdd{\mathbf{CF}^{-, \delta \from \delta'}}
\newcommand{\iHF}{\mathit{HF}}
\newcommand{\iCF}{\mathit{CF}}
\newcommand{\HF}{\mathit{HF}}
\newcommand{\CF}{\mathit{CF}}
\newcommand\Cdd{\C^{-, \delta \from \delta'}}

\def\CFLI{\mathbf{CFL}^{\infty}}

\def\Ringbig{\mathfrak{R}}
\def\Algbig{\mathfrak{R}^Y}
\def\tU{\widetilde{U}}
\def\ms{\mathbf{m}}
\def\Trans{\Xi}

\def\as{\mathbf{a}}
\def\Ybig{\mathcal{Y}}

\begin {document}

\title{Heegaard Floer homology and integer surgeries on links}

\author[Ciprian Manolescu]{Ciprian Manolescu}
\address {Department of Mathematics, Stanford University\\ 
Stanford, CA 94305}
\email {cm5@stanford.edu}

\author [Peter Ozsv\'ath]{Peter Ozsv\'ath}
\address {Department of Mathematics, Princeton University\\ 
Princeton, NJ 08544}
\email {petero@math.princeton.edu}

\begin {abstract}
Let $L$ be a link in an integral homology three-sphere. We give a description of the Heegaard Floer homology of integral surgeries on $L$ in terms of some data associated to $L$, which we call a complete system of hyperboxes for $L$. Roughly, a complete systems of hyperboxes consists of chain complexes for (some versions of) the link Floer homology of $L$ and all its sublinks, together with several chain maps between these complexes. Further, we introduce a way of presenting closed four-manifolds with $b_2^+ \geq 2$ by four-colored framed links in the three-sphere. Given a link presentation of this kind for a four-manifold $X$, we then describe the Ozsv\'ath-Szab\'o mixed invariants of $X$ in terms of a complete system of hyperboxes for the link. Finally, we explain how a grid diagram  produces a particular complete system of hyperboxes for the corresponding link.
\end {abstract}

\maketitle

\newpage
\tableofcontents
\section {Introduction}

Heegaard Floer homology is a tool for studying low-dimensional
manifolds, introduced by Zolt{\'a}n Szab{\'o} and the second author.
In the case of closed three-manifolds, a 
genus $g$ Heegaard diagram naturally endows the $g$-fold symmetric product of
the Heegaard surface with a pair of half-dimensional tori; and a
suitable adaptation of Lagrangian Floer homology in the symmetric
product (relative to the Heegaard tori) turns out to depend
only on the underlying three-manifold. This three-manifold invariant
is constructed in \cite{HolDisk,HolDiskTwo}; related invariants of four-dimensional cobordisms are  constructed in \cite{HolDiskFour, Zemke};
invariants for knots and links in three-manifolds are developed in 
\cite{Knots,RasmussenThesis,Links}.  Of particular interest to us here are the mixed invariants of closed
four-manifolds defined in \cite{HolDiskFour}: they can detect exotic
smooth structures and, in fact, are conjecturally identical to the
Seiberg-Witten invariants \cite{Witten}.

The knot Floer homology groups from \cite{Knots,RasmussenThesis} are
closely related to the (closed) Heegaard Floer homology groups of
three-manifolds obtained as surgeries on the knot. Indeed,
both~\cite{Knots,RasmussenThesis} showed that the filtered knot Floer
complex contains enough information to recover the Heegaard Floer
homologies of all sufficiently large surgeries on the respective
knot. This is used as a stepping-stone to reconstruct the Heegaard
Floer homology of arbitrary surgeries on a knot in~\cite{IntSurg} and
\cite{RatSurg}.

Since every closed three-manifold can be obtained by surgery on a link
in the three-sphere, a natural question is whether the results from
\cite{IntSurg} admit a generalization for links. The goal of this
paper is to present such a generalization.

Let $K \subset Y$ be a knot in an integral homology
three-sphere. Recall that the knot Floer homology of $K$ is
constructed starting from a Heegaard diagram for $Y$ that has two
basepoints $w$ and $z$, which specify the knot. One can build Floer
homology groups by counting pseudo-holomorphic curves in the symmetric
product of the Heegaard surface in various ways. For example, one can
require the support of the curves to avoid $z$, and at the same time
keep track of the intersections with $w$ by powers of a $U$ variable:
this gives rise to knot Floer homology. Alternatively, one can define
complexes $A_s^+$ in which one keeps track the intersection number of the
curve with both $w$ and $z$, in a
way that depends on the value of an auxiliary parameter $s \in \zz$. 
When $s \gg 0$, the complex
$A_s^+$ corresponds to ignoring $z$ completely, and keeping track of
$w$ via a $U$ variable. When $s \ll 0$, we have the reverse: $A_s^+$
corresponds to ignoring $w$ completely, and keeping track of $z$. Note
that whenever $|s|\gg 0$, the homology $H_*(A_s^+)$ is the Heegaard
Floer homology of $Y$, regardless of the sign of $s$ (and in particular,
they are independent of the knot $K\subset Y$). 
The intermediate complexes $A_s^+$ for $s \in
\zz$, however, contain nontrivial information about the knot.
Indeed, according to~\cite{Knots,RasmussenThesis}, the complexes
$A^+_s$ (as $s$ varies) capture the Floer homology groups of the three-manifolds
obtained by sufficiently large surgeries on $K$. Moreover, these complexes are
basic building blocks of the constructions from \cite{IntSurg,RatSurg}.
 
Consider now an oriented link $\orL$ in an integral homology
three-sphere $Y$. The analogue of knot Floer homology was defined in
\cite{Links} and is called {\em link Floer homology}. The construction starts
with a Heegaard diagram $\Hyper^L$ for $Y$ that has several $w$ and
$z$ basepoints, specifying the link. Let $L_1, \dots, L_\ell$ be the
components of $L$. Following \cite{Links}, we consider the affine
lattice $\H(L)$ over $\zz^\ell$ defined by
$$ \H(L) = \bigtimes_{i=1}^\ell \H(L)_i, \ \ \H(L)_i = \frac{\lk(L_i, L - L_i)}{2} + \Z ,$$
where $\lk$ denotes linking number. By keeping track of the basepoints
in various ways, we can define generalized Floer chain complexes 
$\Am(\Hyper^L, \s), \s \in \H(L)$, which are the analogues of the
groups $A_s^+$ for knots with two basepoints. In the case where the diagram $\Hyper^L$ is link-minimal (that is, it has only one $w$ and one $z$ basepoint on each link component), the complexes $\Am(\Hyper^L, \s)$ are made of free modules; in general, this is not the case, and we will have to construct some resolutions $\Chain^-(\Hyper^L, \s)$ of $\Am(\Hyper^L, \s)$. We will then use the groups $\Chain^-(\Hyper^L,
\s)$ to reconstruct the Heegaard Floer homology of integer surgeries
on $L$. Note that, to keep in line with the conventions in
\cite{Links}, we will phrase our construction in terms of the $\iHF^-$
rather than the $\iHF^+$ version of Heegaard Floer homology.  We also ignore sign issues and work over the field $\ff = \zz/2\zz$.  

Further, for technical reasons, we find it useful to use a slightly
different variant of $\iHF^-$ than the one defined in \cite{HolDisk}: we
complete the groups with respect to the $U$ variables, so that they
become modules over the power series ring $\ff[[U]]$, compare
also~\cite{KMBook}. We denote this completed version by $\HFm$; it has
the following technical advantage over $\iHF^-$. The uncompleted version
$\iHF^-$ is functorial under cobordisms equipped with $\spc$ structures,
but it is not functorial under cobordisms per se, whereas the
completed version $\HFm$ is. In particular, $\HFm$ satisfies surgery exact
triangles (analogous to those in~\cite{HolDiskTwo}) just like $\iHF^+$.

Fix a framing $\Lambda$ for the link $\orL$. For a component $L_i$ of
$L$, we let $\Lambda_i$ be its induced framing, thought of as an
element in $H_1(Y - L)$. The latter group can be identified with
$\zz^\ell$ via the basis of oriented meridians for $\orL$. Given a
sublink $M \subseteq L$, we let $\Omega(M)$ be the set of all possible
orientations on $M$. For $\orM \in \Omega(M)$, we let $I_-(\orL,
\orM)$ denots the set of indices $i$ such that the component $L_i$ is
in $M$ and its orientation induced from $\orM$ is opposite to the one
induced from $\orL$. Set
$$ \Lambda_{\orL, \orM} = \sum_{i \in I_-(\orL, \orM)} \Lambda_i \in H_1(Y - L) \cong \zz^\ell.$$

Let $Y_\Lambda(L)$ be the three-manifold obtained from $Y$ by surgery
on the framed link $(L, \Lambda)$. The input that we use to
reconstruct $\HFm(Y_\Lambda(L))$ is called a {\em complete system of hyperboxes}
for the link $\orL$. The precise definition is given in
Section~\ref{sec:complete}.  Roughly, a complete system $\Hyper$
consists of Heegaard diagrams $\Hyper^{L'}$ representing all possible
sublinks $L' \subseteq L$, together with some additional data that
produces maps
$$\Phi^{\orM}_\s: \Chain^-(\Hyper^{L'}, \s) \to \Chain^-(\Hyper^{L'-M}, \psi^{\orM}(\s)),$$
for any $M \subseteq L' \subseteq L, \s \in \H(L')$, and $\orM\in
\Omega(M)$. Here, $\psi^{\orM}: \H(L') \to \H(L' - M)$ are natural
reduction maps. The orientation $\orM$ comes into play as follows:
starting with $\Hyper^{L'}$, we delete the $w$ basepoints
corresponding to components $L_i \subseteq M$ with $i \in I_-(\orL,
\orM)$, and delete the $z$ basepoints corresponding to the remaining
components $L_i \subseteq M$. The resulting diagram represents the
link $L - M$, and the complete system gives us a sequence of steps
that relate it to $\Hyper^{L-M}$, a diagram that also represents
$L-M$. The map $\Phi^{\orM}_\s$ is constructed by following that
sequence of steps. (In the case where the diagram is not link-minimal, the definition of $\Phi^{\orM}_\s$ also involves certain transition maps, which will be described in Section~\ref{sec:general}.)

Define
\begin {equation}
\label {eq:chl_zero}
 \C^-(\Hyper, \Lambda) = \bigoplus_{M \subseteq L} \prod_{\s \in \H(L)}  \Chain^-(\Hyper^{L - M}, \psi^{M}(\s) ),
\end {equation} 
where $\psi^{M}$ simply means $\psi^{\orM}$ with $\orM$ being the orientation induced from the one on $\orL$.
Equip $\C^-(\Hyper, \Lambda) $ with a boundary operator as follows. For $\s \in \H(L) $ and $\x \in  \Chain^-(\Hyper^{L - M}, \psi^{M}(\s))$, set
 \begin {eqnarray}
 \label{eq:DHL} 
\D^-(\s, \x) &=& \sum_{N \subseteq L - M} \sum_{\orN \in \Omega(N)} (\s + \Lambda_{\orL, \orN}, \Phi^{\orN}_{\psi^{M}(\s)}(\x)) \\
&\in&  \bigoplus_{N \subseteq L - M} \bigoplus_{\orN \in \Omega(N)}  \Chain^-(\Hyper^{L-M-N}, \psi^{M \cup \orN} (\s))   \subseteq \C^-(\Hyper, \Lambda). \notag
\end {eqnarray}

Our main result is: 

\begin {theorem}
\label {thm:FirstSurgery}
Fix a complete system of hyperboxes $\Hyper$ for an oriented,
$\ell$-component link $\orL$ in an integral homology three-sphere $Y$,
and fix a framing $\Lambda$ of $L$. There is an isomorphism of
homology groups:
\begin {equation}
\label {eq:surgery0}
 H_*(\C^-(\Hyper, \Lambda)) \cong  \HFm_*(Y_\Lambda(L)),
 \end {equation}
 where $\HFm$ is the
 completed version of Heegaard Floer homology over the power series
 ring $\ff[[U]]$.
\end {theorem}

Similar results hold for other variants of Heegaard Floer homology: $\widehat{\iHF}, \iHF^+$, and $\HFinf$ (the last being a completed version of $\iHF^{\infty}$). 

We refer to the isomorphism \eqref{eq:surgery0} as the {\em link surgery formula}. A more detailed description of this formula (but still without many technicalities) is given in Section~\ref{sec:ex}. We focus there on knots and two-component links, and go over the example of surgeries on the Hopf link. The complete statement of the link surgery formula, in the particular case of complete systems that are link-minimal (that is, the Heegaard diagrams in the system are link-minimal) is found in Section~\ref{sec:statement}. The statement in full generality, for arbitrary complete systems, is given in Section~\ref{sec:general}.

The proof of Theorem~\ref{thm:FirstSurgery} also gives a way of
describing the maps from $\HFm(Y)$ to $\HFm(Y_\Lambda(L))$ induced by
the surgery cobordism. More generally, let $W$ be a cobordism between
two connected three-manifolds $Y_1$ and $Y_2$, such that $W$ consists
of two-handle additions only. We can then find a link $\orL \subset
S^3$ with a framing $\Lambda$ and a sublink $L' \subseteq L$, such
that surgery on $L'$ (with framing specified in $\Lambda$) 
produces $Y_1$, and the framed link specified by $L - L'$, thought
of as a subset of $Y_1$, exactly corresponds to the cobordism $W$ going to $Y_2$. Using
this set-up, we can describe the map on Heegaard Floer homology
induced by the cobordism $W$ in terms of a complete system of
hyperboxes for $L$.

Refining this idea, we arrive at a similar description of a
non-trivial smooth, closed four-manifold invariant. Let $X$ be a
closed, oriented four-manifold with $b_2^+(X) \geq 2$. The
constructions from \cite{HolDiskFour} associated to every $\spc$
structure $\ss$ on $X$ an invariant called the {\em mixed invariant}
$\Phi_{X, \ss} \in \ff$. (The original definition was over $\zz$ and
involved a homology action, but we ignore this extra structure in
the present paper.) The manifold $X$ can be presented in terms of a
link as follows. Delete two four-balls from $X$ to obtain a cobordism
$W$ from $S^3$ to $S^3$. Then, split this cobordism into four-parts
$$ W= W_1 \cup_{Y_1} W_2 \cup _{Y_2} W_3 \cup_{Y_3} W_4,$$
such that $W_1$ consists of one-handles only, $W_2$ and $W_3$ of two-handles, and $W_4$ of three-handles; further, we arrange so that $Y_2$ is an admissible cut in the sense of \cite[Definition 8.3]{HolDiskFour}. Next, find a framed link $(\orL \subset S^3, \Lambda)$ that splits as a disjoint union $L_1 \cup L_2 \cup L_3$, such that surgery on $L_1$ produces $Y_1$, and surgery on $L_2$ and $L_3$ corresponds to the cobordisms $W_2$ and $W_3$, respectively. We refer to the data $(\orL = L_1 \cup L_2 \cup L_3, \Lambda)$ as a {\em cut link presentation} for $X$.

\begin {theorem}
\label {thm:FirstMixed}
Let $X$ be a closed four-manifold with $b_2^+(X) \geq 2$, with a cut link
presentation $(\orL = L_1 \cup L_2 \cup L_3, \Lambda)$. One
can describe the mixed invariants $\Phi_{X, \ss}, \ss \in \spc(X)$ with
coefficients in $\ff=\zz/2\zz$ in terms of the framing $\Lambda$ and a complete system of hyperboxes for
$\orL$.
\end {theorem}

The advantage of presenting the three- and four-manifold invariants in
terms of link Floer complexes is that the latter are
understood better. Indeed, there exist several combinatorial
descriptions of knot Floer homology, see
\cite{MOS,SarkarWang,CubeResolutions}. We focus on the description in
\cite{MOS}, in terms of grid diagrams, which has the advantage that it
extends to all versions of link Floer homology for links in $S^3$. It
turns out that a grid diagram for a link gives rise to a complete
system of hypercubes for that link  (provided that the grid has at least one free marking---see Section~\ref{sec:grid} for the exact condition). Hence, we can apply
Theorems~\ref{thm:FirstSurgery} and \ref{thm:FirstMixed} to obtain a
description of the Heegaard Floer invariants for three- and
four-manifolds in terms of counts of holomorphic curves on symmetric
products of grid diagrams. In \cite{MOT}, we use this result to
describe the Heegaard Floer invariants in purely combinatorial terms.

It would be interesting to see whether the procedure of constructing three- and four-manifold invariants from data associated to links, as presented in this paper, can be applied to other settings. Indeed, there are several homological invariants for links in $S^3$, that have much in common with Floer homology, but for which it is unknown whether they admit extensions to three- and four-manifolds. We are referring in particular to the link homologies constructed by Khovanov \cite{Khovanov}, and Khovanov and Rozansky \cite{KR1,KR2}.

In a different direction, the maps induced by surgery cobordisms can be used to construct another link invariant, a {\em link surgery spectral sequence} (see~\cite[Theorem 4.1]{BrDCov} and \cite[Theorem 5.2]{BaldwinSpectral}). The methods giving Theorem~\ref{thm:FirstSurgery} can also be used to give a description
of this spectral sequence in terms of complete systems of hyperboxes; see Theorem~\ref{thm:SpectralSequence} below for the precise statement.

The organization of this paper is as follows. In Section~\ref{sec:conventions} we explain the difference between the completed theories $\HFm$ and $\HFinf$ used in this paper, and the versions $\iHF^-, \iHF^\infty$ originally defined in \cite{HolDisk}. In Section~\ref{sec:cxes} we define the generalized Floer chain complexes $\Chain^-(\Hyper^L, \s)$, and explain their invariance properties.  We also define maps between generalized Floer complexes by counting $J$-holomorphic polygons.  Section~\ref{sec:ex} contains an overview of the link surgery formula, including a comparison with the knot surgery formula from \cite{IntSurg}, and an explicit computation for the case of the Hopf link. Sections~\ref{sec:conventions} through \ref{sec:ex} are sufficient for the reader who wants to reach a working understanding of the main result, without going into the proof or even into the full details of the statement.

In Section~\ref{sec:hyperco} we discuss some homological algebra that is needed throughout the rest of the paper: we introduce the notion of a hyperbox of chain complexes, and describe several operations on hyperboxes. Section~\ref{sec:triangles} contains some analytical results about the behavior of holomorphic polygon maps under a move called quasi-stabilization. Section~\ref{sec:03} contains similar analytical results for a move called index zero/three link stabilization. Section~\ref{sec:hyperHeegaard} is devoted to building up the definition of a complete system of hyperboxes for a link. In Section~\ref{sec:statement} we then give a more precise statement of Theorem~\ref{thm:FirstSurgery}, in the case of link-minimal complete systems. The proof of  Theorem~\ref{thm:FirstSurgery} for link-minimal complete systems occupies Section~\ref{sec:proof}, and is based on a truncation procedure explained in Section~\ref{sec:t} together with a surgery long exact sequence discussed in Section~\ref{sec:exact}. We then turn to general complete systems: In Section~\ref{sec:general} we give both the statement and the proof of Theorem~\ref{thm:FirstSurgery} for arbitrary complete systems.

In Section~\ref{sec:beyond} we present the extensions of Theorem~\ref{thm:FirstSurgery} to the other versions of Heegaard Floer homology, and to the invariants associated to cobordisms; we also prove Theorem~\ref{thm:FirstMixed} and discuss the link surgeries spectral sequence. Finally, in Section~\ref{sec:grid}, we explain how certain grid diagrams for links in $S^3$ give rise to complete systems of hyperboxes. Many of the diagrams in the resulting complete systems are quasi-stabilized, and therefore the results from Section~\ref{sec:triangles} can be used to simplify the description of the surgery complex in terms of grids.

\medskip \textbf {Acknowledgements.} We are thankful to Dylan Thurston for useful suggestions regarding the material in Section~\ref{sec:hyperco}, to Andr{\'a}s Juh{\'a}sz, Matthew Hedden, Sucharit Sarkar, and Zolt{\'a}n Szab{\'o} for helpful conversations, and to John Baldwin, Damek Davis, Haofei Fan, Jonathan Hales, Tye Lidman, Yajing Liu, Yi Ni, Faramarz Vafaee, Ian Zemke for comments on previous versions of the paper. We are particularly grateful to the anonymous referee for many useful suggestions, and to Ian Zemke for pointing out the need for trajectory $\alpha$-slides in Section~\ref{sec:moves}. The first author would also like to thank the mathematics department at the University of Cambridge for its hospitality during his stay there in the spring of 2009.

The first author was supported by NSF grants DMS-0852439, DMS-1104406, DMS-1402914, and a Royal Society University Research Fellowship. The second author was supported by NSF grants DMS-0804121, DMS-1258274 and a Guggenheim Fellowship.

\section {Conventions}
\label {sec:conventions}
Throughout this paper we work with Floer homology groups
with base field $\ff = \zz/2\zz$. 

Let $Y$ be a closed, connected, oriented $3$-manifold. We consider
the Heegaard Floer homology groups $\widehat{\iHF}, \iHF^+, \iHF^-, \iHF^\infty$
defined in~\cite{HolDisk}, \cite{HolDiskTwo}. These are modules over
the polynomial ring $\ff[U]$. (In the case of $\widehat{\iHF}$, the
action of $U$ is trivial.)

Let $\HFm$ and $\HFinf$ denote the completions of $\iHF^-, \iHF^\infty$
with respect to the maximal ideal $(U)$ in the ring $\ff[U]$. Since
completion is an exact functor, we can alternatively think of $\HFm$
as the homology of the complex $\CFm$ with the same generators as
$\iCF^-$, but whose coefficient ring is the formal power series ring
$\ff[[U]] $ rather than $\ff[U]$. Similarly, $\HFinf$ is the homology
of this same complex, whose base ring is now the field ring of
semi-infinite Laurent polynomials $\ff[[U, U^{-1}]$ (rather than
  $\ff[U, U^{-1}]$, as in the construction of the usual
  $\iHF^\infty$). 
  
  When $\ss$ is a torsion $\spc$ structure, recall that $\iHF^-(Y, \ss)$ is equipped with an
  absolute $\qq$-grading, such that $U$ drops grading by $2$; cf.~\cite{AbsGraded}. We define the $i^{\text{th}}$ graded piece of $\HFm(Y,\ss)$ to be the same as that of $\iHF^-(Y, \ss)$. This does not quite induce a grading on $\HFm(Y,\ss)$ in the usual sense, because $\HFm(Y,\ss)$ is not the direct sum of its graded parts, but rather the completion of the direct sum. Still, by a slight abuse of terminology, we will refer to this structure as an absolute $\qq$-grading on $\HFm(Y,\ss)$. The same goes for $\HFinf(Y, \ss)$.
  
  When $\ss$ is non-torsion, $\iHF^-(Y,
  \ss)$ and $\iHF^\infty(Y, \ss)$ admit relative $\zz/2k\zz$-gradings, for
  suitable $k$ (depending on $\ss$). In this case, they induce true relative $\zz/2k\zz$-gradings on      
  $\HFm(Y,\ss)$ and $\HFinf(Y, \ss)$, characterized by the fact that
  each generator (i.e. intersection point between
  totally real tori) has the same grading as it does when it is
  thought of as a generator of $\iCF^-(Y,\ss)$; and further, that multiplication
  by $U$ drops the grading by $2$.
  
One can define cobordism maps and mixed invariants as in
\cite{HolDiskFour}, using $\HFm$ and $\HFinf$ rather than $\iHF^-$ and
$\iHF^{\infty}$. This new setting is parallel to the one developed by
Kronheimer and Mrowka in \cite{KMBook} in the context of gauge
theory. It has the advantage that $\HFm$ and $\HFinf$ now admit exact
triangles analogous to those for $\widehat{\iHF}, \iHF^+$ from
\cite[Section 9]{HolDiskTwo}. Further, whereas for the definition of
$\iHF^-$ and $\iHF^\infty$ one needs to use strongly admissible Heegaard
diagrams as in \cite[Definition 4.10]{HolDisk}, in order to define
$\HFm$ and $\HFinf$ it suffices to consider weakly admissible
diagrams. Indeed, Lemma 4.13 in \cite{HolDisk} shows that the
differentials of $\CFm$ and $\CFinf$ are finite whenever the
respective Heegaard diagrams are weakly admissible.

More generally, whenever we discuss versions of Heegaard and link
Floer homology that were defined originally over polynomial rings
$\ff[U_1, \dots, U_\ell]$, in this paper we use their completions,
which are modules over formal power series rings $\ff[[U_1, \dots,
U_\ell]]$.

One could also define versions of ${\widehat{\iHF}}$ and $\iHF^+$ using the
completed ring; but since those are generated by complexes on which the
action of multiplication by $U$ is nilpotent on each generator, the resulting 
invariants
coincide with the versions defined over $\Field[U]$. In particular,
in the completed context, we have an exact sequence for 
any closed, oriented three-manifold $Y$
\begin{equation}
\label{eq:HFminfp}
\begin{CD}
...@>>>\HFm(Y)@>>>\HFinf(Y)@>>>\iHF^+(Y)@>>>... \ \ ,
\end{CD}
\end{equation}
where $\iHF^+$ is the Heegaard Floer homology group
from~\cite{HolDisk}.

Note that, when $Y$ is a three-manifold and $\ss$ is a torsion $\spc$
structure on $Y$, the groups $\iHF^-(Y, \ss)$ and $\iHF^\infty(Y, \ss)$
are determined by $\HFm$ and $\HFinf$, respectively. Indeed, we have
$\HFm_{\geq i}(Y, \ss) \cong (\iHF^-)_{\geq i}(Y, \ss)$ for any given
degree $i \in \qq$, and since the groups (and their module structure)
are determined by their truncations, the claim follows.  Similar
remarks apply to $\HFinf$ and $\iHF^\infty$.

For non-torsion $\spc$ structures $\ss$, there is some loss of
information when passing from $\iHF^-(Y,\ss)$ and
$\iHF^{\infty}(Y,\ss)$ to their completed 
analogues. For example, when $Y=S^1 \times S^2$, let $h$ be a
generator of $H^2(Y; \ss) \cong \Z$, and $\ss_k$ the $\spc$ structure
with $c_1(\ss) = 2kh$. Then
$$ \iHF^-(S^1 \times S^2, \ss_k) \cong  \ff[U]/(U^k - 1), \ \  \iHF^\infty(S^1 \times S^2, \ss_k) \cong  \ff[U, U^{-1}]/(U^k - 1),$$
so
$$ \HFm(S^1 \times S^2, \ss_k) \cong \HFinf(S^1 \times S^2, \ss_k) =0.$$

In general, for any non-torsion $\SpinC$ structure $\ss$ we have that
\begin {equation}
\label {eq:infty0}
\HFinf(Y, \ss) = 0.
\end {equation}
Indeed, Lemma 2.3 in \cite{HolDiskSymp} says that
$(1-U^N)\iHF^{\infty}(Y, \ux) = 0$ for some $N \geq 1$. Since $1-U^N$ is
invertible as a power series, after taking the completion we get
\eqref{eq:infty0}.  Consequently, Equation~\eqref{eq:infty0} combined with exactness in the sequence~\eqref{eq:HFminfp} gives a (grading-preserving) isomorphism
\begin {equation}
\label {eq:infty}
\HFm(Y, \ss) \cong \iHF^+(Y, \ss),
\end {equation}
for any non-torsion $\spc$ structure $\ss$.

\newpage

\section{Generalized Heegaard Floer complexes for links}
\label {sec:cxes}

We define here some complexes associated to a Heegaard diagram for a
link. As we shall see in Section~\ref{sec:large}, these are the complexes which govern
large surgeries on links.

\subsection {Heegaard diagrams}
\label {sec:hed}
Combining the constructions of \cite{Links} and \cite{MOS},
define a {\em multi-pointed Heegaard diagram} to be data of
the form $\he=(\Sigma, \alphas , \betas, 
\ws, \zs)$, where:
\begin {itemize} 
\item $\Sigma$ is a closed, oriented surface of genus $g$; 
\item $\alphas = \{\alpha_1, \dots,
\alpha_{g+k-1}\}$ is a collection of disjoint, simple closed curves on
  $\Sigma$ which span a $g$-dimensional lattice of $H_1(\Sigma; \zz)$,
  hence specify a handlebody $U_{\alpha}$; the same goes for
  $\betas= \{ \beta_1, \dots, \beta_{g+k-1} \}$, which specify a handlebody $U_{\beta}$;
\item $\ws= \{ w_1, \dots, w_k\} $ and $ \zs = \{z_1, \dots, z_m\}$ (with $k \geq m$) are collections of points on $\Sigma$ with the
  following property. Let $\{A_i\}_{i=1}^k$ be the connected
  components of $\Sigma - \alpha_1 - \dots - \alpha_{g+k-1}$ and
  $\{B_i\}_{i=1}^k$ be the connected components of $\Sigma - \beta_1 -
  \dots - \beta_{g+k-1}$. Then there is a permutation $\sigma$ of
  $\{1, \dots, m\}$ such that $w_i \in A_i \cap B_i$ for $i=1, \dots,
  k$, and $z_i \in A_i \cap B_{\sigma(i)}$ for $i=1, \dots, m$.
\end {itemize}

We do not take the orderings of the curves and basepoints to be part of the data of the Heegaard diagram; rather, we just ask for such orderings to exist, so that the conditions above are satisfied.

A Heegaard diagram $\he$ describes a closed, connected, oriented
$3$-manifold $Y = U_{\alpha} \cup_{\Sigma} U_{\beta}$, and an oriented
link $\orL \subset Y$ (with $\ell \leq m$ components), obtained as
follows. For $i=1, \dots, m$, we join $w_i$ to $z_i$ inside $A_i$ by
an arc which we then push by an isotopy into the handlebody
$U_{\alpha}$; then we join $z_i$ to $w_{\sigma(i)}$ inside $B_i$ by an
arc which we then push into $U_{\beta}$. The union of these arcs (with
the induced orientation) is the link $\orL$. We then say that $\he$ is
a multi-pointed Heegaard diagram representing $\orL \subset Y$. Note
that the definition we work with here is more general than the notion
of a multi-pointed Heegaard diagram from~\cite{Links}, as we allow
here for more than two basepoints per link component; moreover, we are
allowing for extra basepoints of type $w$ which are not thought of as
belonging to a link component. We refer to $w_{m+1}, \dots, w_k$ as {\em free basepoints}. The set of free basepoints is denoted $\ws^{\operatorname{free}}$.

In order to define the chain complexes associated to a Heegaard
diagram $\he$ (as below, Section~\ref{sec:chains}), we need to
require that it is {\em generic}, i.e. the alpha and beta curves
intersect each other transversely. Further, we should require that 
it is {\em admissible} in
the sense of \cite[Definition 2.2]{MOS}. More precisely: 

\begin {definition}
Let $\he =(\Sigma, \alphas, \betas, \ws, \zs)$ be a multi-pointed Heegaard diagram.

$(a)$ A {\em region} in $\he$ is the closure of a connected component of $\Sigma -  (\alpha_1 \cup \dots \cup \alpha_{g+k-1} \cup \beta_1 \cup  \dots \cup \beta_{g+k-1})$; 

$(b)$ A {\em periodic domain} in $\he$ is a two-chain $\phi$ on $\Sigma$ obtained as a linear combination of regions (with integer coefficients), such that the boundary of $\phi$ is a linear combination of $\alpha$ and $\beta$ curves, and the local multiplicity of $\phi$ at every $w_i \in \ws$ is zero. 

$(c)$ The diagram $\he$ is called {\em admissible} if every non-trivial periodic domain has some positive local multiplicities and some negative local multiplicities.
\end {definition}

Admissibility, as defined in part (c), is equivalent to the requirement that the
underlying diagram $(\Sigma, \alphas, \betas, \ws)$ representing $Y$
is weakly admissible; see \cite[Definition 4.10]{HolDisk}. 
As mentioned in Section~\ref{sec:conventions}, since in this paper we will use
coefficients in power series rings, there is no need to impose the
strong admissibility condition from \cite[Definition 4.10]{HolDisk}.

It is helpful to introduce some terminology for more restrictive classes of Heegaard diagrams:

\begin{definition}
\label{def:linkminimal}
 A Heegaard diagram $(\Sigma,
\alphas, \betas, \ws, \zs)$ is called {\em link-minimal} if $m=\ell$; that is, each link
  component has only two basepoints. 
\end{definition}
 
 \begin{definition}
 \label {def:minimal}
 A Heegaard diagram $(\Sigma,
\alphas, \betas, \ws, \zs)$ for a nonempty link is called {\em minimally-pointed} if $k=m=\ell$; that is, each link  component has only two basepoints, and there are no free basepoints. If the diagram represents an empty link, we call it {\em minimally-pointed} if $k=1$ and $m=0$.
 \end {definition}
 
 Clearly, a minimally-pointed diagram is link-minimal. Minimally-pointed diagrams were the ones originally studied in \cite{HolDisk} and \cite{Links}.
 
 We can impose an even more restrictive condition:
 
  \begin{definition}
 \label {def:hbasic}
 A Heegaard diagram $(\Sigma,
\alphas, \betas, \ws, \zs)$ is called {\em basic} if it is minimally-pointed and, further, for each $i=1, \dots, \ell$, the basepoints $w_i$ and $z_i$ (which determine one of the link components) lie on each side of a beta curve $\beta_i$, and are not separated by any alpha curves.
 \end {definition}

Basic Heegaard diagrams for knots were constructed in \cite[Section 2.2]{Knots}, as the doubly-pointed diagrams associated to marked Heegaard diagrams. To see that every link admits a basic Heegaard diagram, we can start with a minimally-pointed diagram; then, for each $i=1, \dots, \ell$, we add a $1$-handle with one foot near $w_i$ and one near $z_i$, and choose a path between these feet that does not intersect any alpha curves. Further, we add alpha curves along the chosen paths and the cores of the handles, and beta curves along the co-cores of the handles. See Figures~\ref{fig:hopf1} and  \ref{fig:hopf2} for an example of how this is done in the case of a diagram for the Hopf link.
 
 \medskip
 
More generally, we will also need to consider {\em multi-pointed Heegaard multi-diagrams}: data of the form $(\Sigma, \{\etas^i\}_{i=0}^l, \ws, \zs)$, such that each $(\Sigma, \etas^i, \etas^j, \ws, \zs)$ (for $0 \leq i \neq j \leq l$) is a multi-pointed Heegaard diagram. A multi-diagram is called {\em generic} if all the curves in it intersect transversely, and no three curves intersect at one point. Finally, we have the following admissibility condition; compare \cite[Definition 8.8 and Section 8.4.2]{HolDisk}:
 
 \begin {definition}
 \label{def:mpd}
A {\em multi-periodic domain} $\phi$ on a multi-pointed Heegaard multi-diagram $(\Sigma,  \{\etas^i\}_{i=0}^l, $ $\ws, \zs)$ is a two-chain on $\Sigma$ that is a linear combination of the connected components of $\Sigma \setminus (\cup_i \etas^i)$, with integer coefficients, such that $\del \phi$ is a linear combination of $\eta$ curves, and the local multiplicity of $\phi$ at every $w_i \in \ws$ is zero. The Heegaard multi-diagram is called {\em admissible} if every multi-periodic domain has some positive local multiplicities and some negative local multiplicities.
\end {definition}

\subsection {Generalized Floer complexes}
\label {sec:chains}

Let $\orL$ be an oriented $\ell$-component link in an integral
homology sphere $Y$. We denote its components by $L_1, \dots,
L_{\ell}$. Let $\he = (\Sigma, \alphas, \betas, \ws, \zs)$ be an
admissible, generic, multi-pointed Heegaard diagram for $\orL$.  As in the
previous subsection, $k$ denotes the number of basepoints in $\ws$,
and $m$ denotes the number of basepoints in $\zs$, so that $\ell\leq m
\leq k$.
 
In the Introduction we defined sets
$$ \H(L)_i = \frac{\lk(L_i, L - L_i)}{2} + \Z \subset \Q, \ \ \H(L) = \bigtimes_{i=1}^\ell \H(L)_i,$$  
where $\lk$ denotes linking number. Let us also set
$$ \bH(L)_i = \H(L)_i \cup \{-\infty, + \infty\}, \ \ \bH(L) = \bigtimes_{i=1}^{\ell} \bH(L)_i.$$

\begin {remark}
\label {rem:h1}
More invariantly, we could think of $\H(L)$ as an affine lattice over $H_1(Y - L; \Z)$, see \cite{Links}. The latter group is identified with $\Z^\ell$ using the oriented meridians of each component. Furthermore, $\H(L)$ itself is canonically identified with the space of $\spc$ structures on $Y$ relative to the link $\orL$, see \cite[Section 8.1]{Links}.
\end {remark}

The Heegaard diagram determines tori
$$\Ta = \alpha_1 \times \dots \times \alpha_{g+k-1}, \ \Tb=\beta_1 \times \dots \times \beta_{g+k-1} \subset \Sym^{g+k-1}(\Sigma).$$ 

For $\x, \y \in \Ta \cap \Tb$, we let $\pi_2(\x, \y)$ be the set of homology classes of Whitney disks from $\x$ to $\y$ relative to $\Ta$ and $\Tb$, as in \cite{HolDisk}. For each homology class of disks $\phi \in \pi_2(\x, \y)$, we denote by $n_{w_j}(\phi)$ and $n_{z_j} (\phi) \in \Z$ the multiplicity of $w_j$ (resp. $z_j$) in the domain of $\phi$. Further, we let $\mu(\phi)$ be the Maslov index of $\phi$.

Any intersection point $\x \in \Ta \cap \Tb$ has a Maslov grading $M(\x) \in \Z$ and an Alexander multi-grading given by 
$$A_i(\x) \in \H(L)_i, \ i \in \{1, \dots, \ell \}.$$ 

Let $\Ws_i$ and $\Zs_i$ be the set of indices for the $w$'s (resp. $z$'s) belonging to the $i\th$ component of the link. We then have
$$ A_i(\x) - A_i(\y) = \sum_{j \in \Zs_i} n_{z_j}(\phi) - \sum_{j \in \Ws_i} n_{w_j}(\phi),$$
where $\phi$ is any class in $\pi_2(\x, \y)$. Note that this formula determines the Alexander grading up to an overall integer shift. The absolute Alexander grading can be pinned down by either relating intersection points in $\Ta \cap \Tb$ to relative $\spc$ structures on $(Y, \orL)$ via vector fields, or (more simply) by using the symmetry of link Floer homology. See \cite[Sections 4 and 8]{Links} and \cite[Section 2.1]{MOS} for details.

Following \cite{Links} and \cite{MOS}, we define the {\em link Floer complex} $\CFLm(\he)$ as follows. We let $\CFLm(\he)$ be the free module over $\Ring(\he) = \Field[[U_1, \dots, U_k]]$ generated by $\Ta \cap \Tb$, and equipped with the differential: 
\begin{equation}
  \label{def:DefHFLm}
  \partial \x = \sum_{\y \in \Ta \cap \Tb} \sum_{\substack{\phi \in \pi_2(\x, \y)\\ \mu(\phi)=1 }} \# (\M(\phi)/\R) \cdot U_{1}^{n_{w_1}(\phi)} \cdots   U_{k}^{n_{w_k}(\phi)}\y. 
\end{equation}
Here, $\M(\phi)$ is the moduli space of pseudo-holomorphic curves (solutions to Floer's equation) in the class $\phi$, and $\R$ acts on $\M(\phi)$ by translations. Note that $\M(\phi)$ depends on the choice of a suitable path of almost complex structures on the symmetric product. We suppress the almost complex structures from notation for simplicity. 

The Maslov grading $M$ produces the homological grading on $\CFLm(\he)$, with each $U_i$ decreasing $M$ by two. Furthermore, each Alexander grading $A_i$ defines a filtration on $\CFLm(\he)$, with $U_i$ decreasing the filtration level $A_i$ by one, and leaving $A_j$ constant for $j \neq i$.

Of course, if we ignore the Alexander filtrations, then $\CFLm(\he)$ is simply a Heegaard Floer complex representing the underlying three-manifold $Y$, and we could write it as $\CFm(\Sigma, \alphas, \betas, \ws)$ or  $\CFm(\Ta, \Tb, \ws).$ However, we use the notation $\CFLm(\he)$ when we want to view it as an $\H(L)$-filtered complex. 

Given $\s = (s_1, \dots, s_\ell) \in \H(L)$, we define the {\em  generalized Heegaard Floer complex}
$$\Am(\he, \s) = \Am (\he, s_1, \dots, s_\ell)=\Am(\Ta,\Tb,\s)$$ 
to be the subcomplex of $\CFLm(\he)$ generated by elements $\x \in \Ta \cap \Tb$ with $A_i(\x) \leq s_i$ for all $i=1, \dots, \ell$.

\begin {remark}
\label {remark:knot1}
When $L=K$ is a knot and the diagram $\he$ has only two basepoints
(one $w$ and one $z$), the complex $\CFLm(\he)$ is the completion of the subcomplex $C\{i\leq 0\}$ of the knot Floer complex $CFK^{\infty}(Y, K)$, in the notation of \cite{Knots}. Further, the generalized Floer complex $\Am(\he, s)$ is the completion of the subcomplex $A_s^- = C\{\max(i, j-s) \leq 0\}$. Note
that in \cite{IntSurg}, the formula for integral surgeries along a
knot was stated in terms of $\iHF^+$. It involved a similar complex
$A_s^+ = C\{\max(i, j-s) \geq 0\}$.
\end {remark}

\begin {remark}
\label {rem:gen}
The restriction to the underlying three-manifold $Y$ being an integral
homology sphere is not essential. Indeed, with minor modifications,
one can define generalized Heegaard Floer complexes for arbitrary
null-homologous links in a three-manifold. Let us explain the only
other version of the construction relevant to this paper, namely when
$\he$ represents an unlink $\orL$ (of $m$ components) inside a ball in
$Y = \#^n(S^1 \times S^2)$. The generators $\x \in \Ta \cap \Tb$ are
then grouped in equivalence classes according to the $\spc$ structures
on $Y$. We focus attention on only those $\x$ that lie in the torsion
$\spc$ structure $\ss_0$. They have an Alexander grading $A(\x) \in
\zz^m$, and generate a link Floer complex $\CFLm(\he; \ss_0)$. For $\s \in \zz^m$, we then define $\Am(\he, \s)$ by the same procedure as above, but using only generators in
$\ss_0$. 
\end {remark}

So far we have defined $\Am(\Hyper, \s)$ for $\s = (s_1, \dots, s_\ell) \in \H(L)$. We can extend this definition to  all $\s \in \bH(L)$; that is, we can let some of the $s_i$ be $\pm \infty$. Indeed, if we allow some $s_i$ to be $+\infty$, we can use the same definition as before, and the condition $A_i(\x) \leq s_i$ is vacuous. In particular, when all $s_i$ are $+ \infty$, the complex $\Am(\Hyper, \s)$ is just the completed Heegaard Floer complex
$$ \CFm(\Sigma, \alphas, \betas, \ws) = \CFm(\Ta, \Tb, \ws)=\CFm(Y)$$ 
 representing the underlying three-manifold $Y$.

If we allow some of the $s_i$ to be $-\infty$, we cannot use the same definition, since $A_i(\x) \leq -\infty$ is impossible. Instead, we let $\Am(\he, \s)$ be generated by those $\x$ with $A_i(\x) \leq s_i$ for the values of $i$ with $s_i \neq -\infty$, whereas for the values of $i$ that are $-\infty$, we replace the coefficient $U_i^{n_{w_i}(\phi)}$ in \eqref{def:DefHFLm} with $U_i^{n_{z_i}(\phi)}$. In particular, when all $s_i$ are $- \infty$, the complex $\Am(\Hyper, \s)$ is just the completed Heegaard Floer complex
$$ \CFm(\Sigma, \alphas, \betas, \zs \cup \ws^{\operatorname{free}}) = \CFm(\Ta, \Tb, \zs \cup \ws^{\operatorname{free}}),$$ 
again representing the three-manifold $Y$.

\subsection{Heegaard moves}
\label{sec:hmoves}
 
Let $\he = (\Sigma, \alphas, \betas, \ws, \zs)$ and $\he' =
(\Sigma', \alphas', \betas', \ws', \zs')$ be two 
multi-pointed Heegaard diagrams representing the same link $\orL
\subset Y$.
If this is the case, we say that $\he$ and $\he'$ are {\em
equivalent}. It is important to note that we are fixing the 3-manifold
$Y$, not just its diffeomorphism class, and we are also fixing the
link $\orL$, not just its isotopy class. Similarly, the Heegaard surfaces
$\Sigma$ and $\Sigma'$ are thought of as embedded in $Y$.

\begin{figure}
\begin{center}
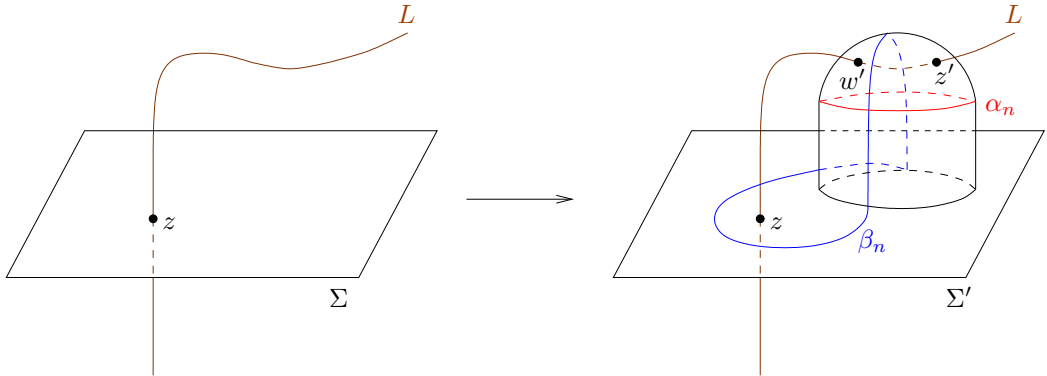
\end{center}
\caption {{\bf An index zero/three link stabilization.} The surface $\Sigma'$ is obtained from $\Sigma$ by deleting a disk and adding a cap, i.e., taking the connected sum with a sphere inside $Y$. This is the same picture as \cite[Figure 5]{MOS}, but we emphasize the fact that the construction is done inside the fixed 3-manifold $Y$, with a fixed link $L$.}
\label{fig:Stab03}
\end{figure}

We say that $\he$ and $\he'$ are related by:

\begin{enumerate}
\item a {\em 3-manifold isotopy} if there is a self-diffeomorphism $\phi: Y \to Y$ isotopic to the identity such that $\phi(L) = L$, $\phi(\Sigma) = \Sigma'$, and $\phi$ takes all the curves and basepoints on $\Sigma$ into the corresponding ones on $\Sigma'$; 
\item an {\em $\alpha$-curve isotopy}  if $\Sigma=\Sigma', \ws=\ws', \zs = \zs', \betas=\betas'$, and the curves in $\alphas'$ are obtained from those in $\alphas$ by an ambient isotopy of $\Sigma$ supported away from all the basepoints; a {\em $\beta$-curve isotopy} is similar, with the roles of $\alphas$ and $\alphas'$ switched with those of $\betas$ and $\betas'$;
\item an {\em $\alpha$-handleslide} if  $\Sigma=\Sigma', \ws=\ws', \zs = \zs', \betas=\betas'$, a curve $\alpha'_i \in \alphas'$ is obtained by handlesiding the respective curve  $\alpha_i \in \alphas$ over another curve in $\alphas$ (with the handleslide being done away from the basepoints), and the other curves in $\alphas'$ are isotopic (by an isotopy away from the basepoints) to the respective curves in $\alphas$; a {\em $\beta$-handleslide} is similar, with the roles of $\alphas$ and $\alphas'$ switched with those of $\betas$ and $\betas'$;
\item an {\em index one/two stabilization} if $\Sigma'$ is obtained by taking the connected sum (inside $Y$) of $\Sigma$ and a genus one surface with one alpha and beta curve intersecting transversely at a point; the connected sum is done away from the link $L$ and the basepoints, so that $\ws = \ws', \zs =\zs' $;
\item a {\em free index zero/three stabilization} if $\Sigma' = \Sigma, \alphas' = \alpha \cup \{\alpha_{n}\}, \betas' = \betas \cup \{\beta'\}, \ws' = \ws \cup \{w'\}, \zs' = \zs$, where the new curves $\alpha'$ and $\beta'$ intersect each other in two points, do not intersect any of the other curves, and both bound disks containing the new free basepoint $w'$. This kind of stabilization was called simple in \cite[Section 6.1]{Links}. See \cite[Figure 3]{Links} for a picture; 
\item an {\em  index zero/three link stabilization} if $\Sigma'$ is obtained from $\Sigma$ by taking the connected sum (inside $Y$) with a sphere, in the neighborhood of a $z$ basepoint, as shown in Figure~\ref{fig:Stab03}. This introduces the new curves $\alpha'$ and $\beta'$, as well as an  additional $(w, z)$ pair, denoted $(w', z')$ in the picture. 
\end{enumerate}

The moves above and their inverses are called {\em Heegaard moves}. (The inverses of the stabilization moves are called {\em destabilizations}.) An $\alpha$-curve isotopy or an $\alpha$-handleslide (resp. a $\beta$-curve isotopy or $\beta$-handleslide) is called {\em admissible} if the corresponding Heegaard triple diagram $(\Sigma, \alphas, \betas, \alphas', \ws, \zs)$ (resp. $(\Sigma, \alphas, \betas, \betas', \ws, \zs)$) is admissible. Admissible curve isotopies, admissible handleslides, together with all other moves (i) and (iv)-(vi) and their inverses are called {\em admissible Heegaard moves}.

Combining the arguments in \cite{HolDisk}, \cite[Proposition 3.9]{Links} and \cite[Lemma 2.4]{MOS}, we obtain the following:

\begin {proposition}
\label {prop:equiv}
Any two equivalent, generic, admissible, multi-pointed Heegaard diagrams for $\orL$ can be related by a sequence of admissible Heegaard moves. 
\end {proposition}

\begin{remark}
In fact, one can check that 3-manifold isotopies are redundant in the list (i)-(vi) above, that is, they can be obtained as combinations of the other moves. For example, a 3-manifold isotopy away from all the curves and the basepoints can be viewed as an index one/two stabilization, followed by an index one/two destabilization. Nevertheless, we keep (i) in the list because we will  refer to this move in some of our future discussions, and it makes the exposition clearer. 
\end{remark}

We have the following refinement of Proposition~\ref{prop:equiv}:

\begin {proposition}
\label {prop:RefinedMoves} 
Let $\he$ and $\he'$ be two equivalent, generic, admissible, multi-pointed Heegaard diagrams representing the same link $\orL$ with $\ell$ components. Suppose the diagram $\he$ is link-minimal, and that $\he'$ has $m$ basepoints of type $z$. For $\he'$, let us choose one basepoint of type $w$ on each component of $L$. The remaining $m-\ell$ basepoints of type $w$ on $\he'$ that are not free are called {\em subsidiary}.

Then, the diagram $\he'$ can be obtained from $\he$ by a sequence of admissible Heegaard moves that includes exactly $m-\ell$ index zero/three link stabilizations, with each of these stabilizations introducing a subsidiary basepoint. \end {proposition}

\begin{proof}
This follows from the arguments in \cite[proof of Lemma 2.4]{MOS}.
\end{proof}
 
For future reference, let us introduce some more terminology. 
If two Heegaard diagrams $\he$ and $\he'$ have the same underlying
Heegaard surface
$\Sigma$, and their 
collections of curves $\betas$ and $\betas'$  are
related by isotopies and handleslides only (supported away from the
basepoints), we say that $\betas$ and $\betas'$ are {\em strongly
equivalent}.

\begin {definition}
\label {def:ab}
Consider two multi-pointed Heegaard diagrams $\he =
(\Sigma, \alphas, \betas, \ws, \zs)$ and $\he' = (\Sigma', \alphas',
\betas', \ws', \zs')$, representing the same link $\orL \subset Y$.

$(a)$ The diagrams $\he$ and $\he'$ are called {\em strongly equivalent} if $\Sigma
= \Sigma', \ws=\ws', \zs = \zs'$, the curve collections $\alphas$ and
$\alphas'$ are strongly equivalent, and $\betas$ and $\betas'$ are
strongly equivalent as well. In other words, $\he$ and $\he'$ should differ by the Heegaard moves $(ii)$ and $(iii)$ only.

$(b)$ We say that $\he$ and $\he'$ differ by a {\em surface isotopy} if $\Sigma = \Sigma'$, and there is a self-diffeomorphism $\phi: \Sigma \to \Sigma$ isotopic to the identity and supported away from the link $\orL$, and such that $\phi$ takes the curves and basepoints in $\he$ to the corresponding curves and basepoints in $\he'$. If $\he$ and $\he'$ are surface isotopic, we write $\he \cong \he'$.
\end {definition}

Note that any surface isotopy can be embedded into an ambient 3-manifold isotopy. For this reason, we did not include it in the list of Heegaard moves (i)-(vi). However, surface isotopies will play an essential role in the definition of complete systems of hyperboxes in Section~\ref{sec:complete}. Observe also that, by definition, surface isotopies keep the basepoints on $\orL$ fixed. If a surface isotopy kept all the basepoints fixed, it could simply be viewed as a composition of curve isotopies. Therefore, the main role of surface isotopies is to allow a way of moving the free basepoints on $\Sigma$, while keeping $\Sigma$ fixed. 

\subsection{Invariance}
\label{sec:Inv}
Just as in the usual case of Heegaard Floer homology (or link Floer homology), there is an invariance statement for generalized Floer complexes for links. Before stating it, we need some preliminaries.

We start with a well-known algebraic lemma. We will use this lemma and variants of it throughout the paper, so we include its proof for completeness.
\begin{lemma}
\label{lem:ap1}
Let $R$ be any $\ff$-algebra, and let $C$ be a free complex over $R[[U_1]]$. Let $C_-$ and $C_+$ be copies of $C$. Consider the mapping cone complex
$$ C' = \mathit{Cone}\bigl(C_-[[U_2]] \xrightarrow{U_1-U_2} C_+[[U_2]] \bigr).$$
 Then, the complexes $C$ and $C'$ are chain homotopy equivalent over $R[[U_1]]$.
\end{lemma}

\begin{proof}
Let $\S$ be a set of free generators for $C$, over $R[[U_1]]$. For each $\x \in \S$, we denote by $\x_-$ and $\x_+$ the corresponding generators of $C_-$ resp. $C_+$. Then, $C'$ is free over $R[[U_1, U_2]]$, with generators $\x_-$ and $\x_+$ for each $\x \in \S$.  

Consider the morphisms over $R[[U_1]]$ given by
$$ \iota: C \to C', \  \ \iota(\x) = \x_+$$
and
$$ \rho: C' \to C, \ \ \rho(U_2^n \x_-)=0, \ \rho(U_2^n \x_+) = U_1^n \x.$$
We claim that $\iota$ and $\rho$ are homotopy inverses. Indeed, $\rho \circ \iota$ is the identity, whereas $\iota \circ \rho$ is chain homotopic to the identity via the homotopy
$$ H: C' \to C', \ \ H(U_2^n \x_-)=0, \ H(U_2^n \x_+) = \frac{U_1^n-U_2^n}{U_1-U_2} \x_-.$$
Here, by the fraction $ \frac{U_1^n-U_2^n}{U_1-U_2} $ we meant the polynomial expression
$$ U_1^{n-1} + U_1^{n-2}U_2 + \dots + U_1U_2^{n-2} + U_2^{n-1}.$$
This completes the proof.\end{proof}

\begin{remark}
\label{rem:ap1}
The map $\iota$ is a morphism over $R[[U_1, U_2]]$, which induces an isomorphism on homology. Thus, if we let $U_2$ act on $C$ just as $U_1$ does, we find that $C$ and $C'$ are quasi-isomorphic over $R[[U_1, U_2]]$. However, they are not normally chain homotopy equivalent over this bigger ring.
\end{remark}

Next, following \cite[Section 2]{Links}, if $C_*$ and $C'_*$ are complexes with a filtration by a $\zz^\ell$-affine space such as $\he(L)$, we say that $C_*$ and $C'_*$ are filtered chain homotopy equivalent if there are filtered chain maps $f: C_* \to C'_*, \ g: C'_* \to C_*$ and filtered chain homotopies $H_1: C_* \to C'_{*+1}, H_2:  C'_* \to C_{*+1}$ such that $g \circ f-\id = \del \circ H_1 + H_1 \circ \del$ and $f \circ g-\id = \del \circ H_2 + H_2 \circ \del$. 

Recall that we defined $\CFLm(\he)$ and $\Am(\he, \s)$ as complexes over the ring $\Ring(\he)=\ff[[U_1, \dots, U_k]]$, with one variable for each $w$ basepoint. Assuming $L$ is non-empty, let us also consider the ring
$$ \Ring(L)=\ff[[\tU_1, \dots, \tU_\ell]],$$
where each variable $\tU_i$ is associated to a link component $L_i \subseteq L$. We define a {\em natural inclusion} of $\Ring(L)$ into $\Ring(\he)$ to be a homomorphism (of unitary rings) that takes each $\tU_i$ to one of the $U_j$ variables such that the basepoint $w_j$ is on the link component $L_i$. Note that when the diagram $\he$ is link-minimal, then there is a unique natural inclusion $\Ring(L) \hookrightarrow \Ring(\he)$, and that inclusion is in fact an isomorphism.

If $L$ is empty, we let $\Ring(L)=\ff[[\tU]]$, and a natural inclusion is a homomorphism that takes $\tU$ to any of the $U_i$ variables.

Given a natural inclusion $\Ring(L) \hookrightarrow \Ring(\he)$, we can view $\CFLm(\he)$ and $\Am(\he, \s)$ as complexes over $\Ring(L)$, by restriction of scalars.

We are now ready to state the invariance theorem.

\begin{theorem}
  \label{thm:LinkInvariance}
  Let $\he$ and $\he'$ be two generic, admissible, multi-pointed
  Heegaard diagram for the same oriented link $\orL \subset Y$. Pick $\s \in \bH(L)$. Pick also natural inclusions $\Ring(L) \hookrightarrow \Ring(\he)$ and $\Ring(L) \hookrightarrow \Ring(\he')$, and view $\CFLm(\he)$, $\CFLm(\he')$, $\Am(\he, \s)$ and $\Am(\he', \s)$ as complexes over $\Ring(L)$, according to those inclusions. Then: 
  \begin{enumerate}[(a)]
\item The complexes $\CFLm(\he)$ and $\CFLm(\he')$ are filtered chain homotopy equivalent  over $\Ring(L)$.
\item The complexes $\Am(\he, \s)$ and $\Am(\he', \s)$ are chain homotopy equivalent over $\Ring(L)$.
\end{enumerate}
\end{theorem}

Note that part (b) is a simple algebraic consequence of part (a). Furthermore, part (a) follows by combining the arguments from~\cite[Theorem 4.4]{Links} and \cite[Proposition 2.5]{MOS}. However, some aspects of this proof will be of particular importance for the present paper. We give below an outline of the proof with special emphasis on those aspects.

\begin{proof}[Beginning of the proof of Theorem~\ref{thm:LinkInvariance}]
We need to verify invariance under each of the Heegaard moves (i)-(vi). Let us discuss these moves briefly in turn.

The proof of invariance under curve isotopies and index one/two stabilizations is completely similar to the one in \cite{HolDisk}. 

To check invariance under a 3-manifold isotopy, one needs to pull back the almost complex structure on the symmetric product under the isotopy, and then use the contractibility of the space of compatible almost complex structures. (This is the same argument as the one used in \cite[Theorem 6.1]{HolDisk}.)

Invariance under handleslides can also be proved along the lines in \cite{HolDisk}, but one has to define suitable (filtered) polygon maps between the link Floer complexes: the details of their construction are explained in Section~\ref{sec:polygon}. 

Let us now discuss the remaining two Heegaard moves (v)-(vi), which change the number of basepoints. 

Suppose that we change a  diagram $\he$ (with $k$ basepoints of type $w$) by a free index zero/three  stabilization introducing an additional 
free basepoint $w'$. Let $\he'$ be the new diagram.  
By the same argument as in \cite[proof of Proposition 6.5]{Links}, the new link Floer complex $\CFLm(\he')$ is the mapping cone
\begin{equation}
\label{eq:cone01}
 \CFLm(\he)[[U_{k+1}]] \xrightarrow{U_{j} - U_{k+1}} \CFLm(\he)[[U_{k+1}]], 
 \end{equation}
 for some $j \in \{1, \dots, k\}$. The Alexander filtrations on both the domain and the target of \eqref{eq:cone1} are induced from those on $ \CFLm(\he)$, with $U_{k+1}$ keeping the filtration level fixed (since it corresponds to a free basepoint). By Lemma~\ref{lem:ap1}, we have that $\CFLm(\he')$ is chain homotopy equivalent to $\CFLm(\he)$ over the previous ring $\Ring(\he)$ and hence over $\Ring(L)$, assuming that the natural inclusion $\Ring(L) \hookrightarrow \Ring(\he')$ factors through the natural inclusion into $\Ring(\he')$. Moreover, the maps constructed in the proof of Lemma~\ref{lem:ap1} are filtered, so $\CFLm(\he)$ and $\CFLm(\he')$ are filtered chain homotopy equivalent.

Given this, we get that $\Am(\he', \s)$ and $\Am(\he, \s)$ are chain homotopy equivalent. In fact, $\Am(\he', \s)$ is  the mapping cone
 \begin{equation}
 \label{eq:cone02}
  \Am(\he, \s)[[U_{k+1}]] \xrightarrow{U_{j} - U_{k+1}} \Am(\he, \s)[[U_{k+1}]]. 
   \end{equation}

Finally, suppose that we change a  diagram $\he$ by an index zero/three link stabilization introducing an additional $(w, z)$ pair on a link component $L_i \subseteq L$. By the argument in \cite[Proposition 2.3]{MOS}, the new link Floer complex $\CFLm(\he')$ is also of the form \eqref{eq:cone01}, albeit in this case the $A_i$ filtration on the domain of \eqref{eq:cone01} is shifted by $-1$ from the original one on $ \CFLm(\he)$. Nevertheless, since $U_{k+1}$ also decreases the $A_i$ filtration level by one, we get that $\CFLm(\he)$ is still filtered chain homotopy equivalent to $\CFLm(\he)$. This implies that $\Am(\he', \s)$ and $\Am(\he, \s)$ are chain homotopy equivalent, although we should point out that in this case $\Am(\he', \s)$ is not of the form \eqref{eq:cone02}. Once again, the chain homotopy equivalence is over $\Ring(L)$, and we need to assume that there is a factorization of natural inclusions 
$\Ring(L) \hookrightarrow \Ring(\he) \hookrightarrow \Ring(\he')$.

This completes the proof of invariance over $\Ring(L)$, with respect to {\em some} choices of natural inclusions $\Ring(L) \hookrightarrow \Ring(\he)$ and $\Ring(L) \hookrightarrow \Ring(\he')$. To see that invariance holds with respect to any such choices, we can make use of Proposition~\ref{prop:RefinedMoves}; this allows us to choose the basepoints that are introduced when we do the index zero/three link stabilizations. 
 \end {proof}
 
 We are grateful to Andr{\'a}s Juh{\'a}sz and Sucharit Sarkar for conversations about the following point:
 
 \begin {remark}
 \label {rem:canonical}
 Part (b) of Theorem~\ref{thm:LinkInvariance} implies that the homology $H_*(\Am(\he, \s))$ (for basic diagrams $\he$) is an invariant of $\orL \subset Y$ and $\s$, up to isomorphism. However, this is not so up to canonical isomorphism. The problem is that when relating a Heegaard diagram to itself via equivalence moves that include 3-manifold isotopies, the basepoints may trace homotopically non-trivial loops during those isotopies. Nevertheless, the generalized Floer homology $H_*(\Am(\he, \s))$ (for basic diagrams) is an invariant of $\orL \subset Y$, the level $\s$ and the basepoints $\ws, \zs \subset L$, up to canonical isomorphism. This can be proved along the lines of \cite[Section 6]{OzsvathStipsicz}.  
 \end {remark}

\subsection {Polygon maps}
\label {sec:polygon}

In this section we define certain maps between generalized Floer complexes that count holomorphic polygons. For the sake of clarity, we start with the case of triangles, which (along with quadrilaterals) are  used in the part of the proof of Theorem~\ref{thm:LinkInvariance} that deals with invariance under handleslides.  

Let $\he=(\Sigma, \alphas, \betas, \ws, \zs)$ be a generic, admissible,  
multi-pointed Heegaard diagram representing a link $\orL \subset Y$,
as in Section~\ref{sec:hed}. Let $\gammas$ be a new set of attaching curves for $\Sigma$, which is strongly equivalent to $\betas$. We assume that the Heegaard multi-diagram $(\Sigma, \alphas, \betas, \gammas, \ws, \zs)$ is generic and admissible. Note that the diagram $\he''  = (\Sigma, \alphas, \gammas, \ws, \zs)$ represents $\orL \subset Y$, whereas $\he'= (\Sigma, \betas, \gammas, \ws, \zs)$ represents an unlink in the connected sum of several $S^1 \times S^2$. More precisely, with $g, m, k$ as in Section~\ref{sec:hed}, we have that $\he'$ represents the unlink of $m$ components inside $\#^{g}(S^1 \times S^2)$.

We have link Floer complexes $\CFLm(\he)$ and $\CFLm(\he'')$ filtered by $\H(L)$, and (according to Remark~\ref{rem:gen}) a complex $\CFLm(\he'; \ss_0)$ filtered by $\zz^m$. We also have generalized Floer complexes $\Am(\he, \s)$ and $\Am(\he'', \s)$ for $\s \in \bH(L)$, and $\Am(\he', \s')$ for $\s' \in \zz^m$.

Note that there is a projection $\zz^m \longrightarrow \zz^\ell$, which takes $\s' = (s'_1, \dots, s'_m)$ to $\bar \s' = (\bar s'_1, \dots, \bar s'_\ell)$, where
$$ \bar s'_i = \sum_{j \in \Ws_i}  s'_j.$$
Since $\bH(L)$ is an affine space over $\zz^\ell$, for any $\s \in \bH(L)$ and $\s' \in \zz^m$  we can make sense of the expression $\s + \bar \s' \in \bH(L)$.

For $\x \in \Ta \cap \Tb, \x' \in \Tb \cap \Tg, \y \in \Tg \cap \Ta$, we let $\pi_2(\x, \x', \y)$ be the space of homotopy classes of  Whitney triangles connecting $\x, \x', \y$, as in \cite[Section 8.1.2]{HolDisk}. We define a triangle map
\begin{equation}
\label{eq:Triangle}
 f_{\alpha \beta \gamma}: \CFLm(\he) \otimes \CFLm(\he'; \ss_0) \to \CFLm(\he'')
 \end{equation}
by the formula 
$$ f_{\alpha \beta \gamma}(\x \otimes \x') = \sum_{\y \in \Ta \cap \Tg} \sum_{\{\phi \in \pi_2(\x, \x', \y)|\mu(\phi)=0 \}} \#( \M(\phi)) \cdot U_{1}^{n_{w_1}(\phi)} \cdots   U_{k}^{n_{w_k}(\phi)}  \y. $$

Here, $\M(\phi)$ is the moduli space of pseudo-holomorphic triangles in the class $\phi$, which is required to have Maslov index $\mu(\phi) =0$. It is straightforward to check that $f_{\alpha \beta \gamma}$ is a filtered chain map, where the filtration on the domain is by $\H(L)$, given by adding up the filtrations on $\CFLm(\he)$ and $\CFLm(\he'; \ss_0)$ via $(\s, \s') \to \s + \bar \s'.$ 

Hence, for  $\s=(s_1, \dots, s_\ell) \in \bH(L)$ and $\s'=(s'_1, \dots, s'_\ell) \in \zz^\ell$, we obtain triangle maps 
$$\Am(\he, \s) \otimes \Am(\he', \s') \to \Am(\he'', \s + \bar \s'),$$
which we still denote by $ f_{\alpha \beta \gamma}$.

This construction of triangle maps can be extended to more general polygon maps; compare \cite[Section 4.2]{BrDCov}. Let $(\Sigma,  \{\etas^j\}_{j=0}^r, \ws, \zs)$ be an admissible, generic, multi-pointed  Heegaard multi-diagram, such that each $\etas^j$ is a $(g+k-1)$-tuple of attaching circles. For simplicity, we assume that the curve collections $\etas^j$ come in two equivalence classes, such that the collections in the same class are strongly equivalent, and the diagrams formed by curve collections in different classes represent a link $\orL$ in an integral homology sphere $Y$.  We can then define Floer groups $\Am(\T_{\eta^i}, \T_{\eta^j}, \s)$, where each $\s$ is either in $\zz^m$ or in $\bH(L)$. In the former case we have a well-defined projection $\bar \s \in \zz^\ell$; in the latter, we use the notation $\bar \s$ to simply denote $\s$.

Suppose $r \geq 1$, and that $\etas^0$ and $\etas^r$ are in different equivalence classes, so that $(\Sigma, \etas^0, \etas^r, \ws, \zs)$ represents the link in the integral homology sphere. Then, we can define linear maps
\begin {equation}
\label {eq:feta}
 f_{\eta^0, \dots, \eta^r} : \bigotimes_{j=1}^r \Am(\T_{\eta^{j-1}},\T_{\eta^j}, \s_j) \to \Am(\T_{\eta^0}, \T_{\eta^r}, \bar \s_1 + \dots+\bar \s_r)
 \end {equation}
$$ f_{\eta^0, \dots, \eta^r}(\x_1 \otimes \dots \otimes \x_r)= \sum_{\y \in \T_{\eta^{0}} \cap \T_{\eta^r}}\sum_{\{\phi \in \pi_2(\x_1, \dots, \x_r, \y)| \mu(\phi)=2-r\} } \# (\M(\phi)) \cdot  U_{1}^{n_{w_1}(\phi)} \cdots   U_{k}^{n_{w_k}(\phi)} \y,$$ 
given by counting isolated pseudo-holomorphic $(r+1)$-gons in $\Sym^{g+k-1}(\Sigma)$, with edges on $\T_{\eta^0}, \dots$, $\T_{\eta^r}$, and with specified vertices. Here, the Maslov index $\mu(\phi)$ denotes the expected dimension of the space of pseudo-holomorphic polygons in the class $\phi$, where the domain (a disk with $r+1$ marked points) has a fixed conformal structure. Since the moduli space of conformal structures on the domain has dimension $(r-2)$, the Maslov index equals the expected dimension of $\M(\phi)$ minus $(r-2)$, where $\M(\phi)$ is the space of all pseudo-holomorphic $(r+1)$-gons in the class $\phi$. We warn the reader that this definition of $\mu(\phi)$ is different from the one in \cite[Section 4.2]{BrDCov}, where $\mu(\phi)$ was simply the expected dimension of $\M(\phi)$. Our definition of $\mu(\phi)$ coincides with that used by Sarkar in \cite[Section 4]{SarkarIndex}. It has the advantage that it makes the Maslov index additive under the natural juxtaposition maps.

When $r=1$, by $f_{\eta^0, \eta^1}$ we simply mean the differential $\del$ for a generalized Floer complex $\Am(\T_{\eta^0}, \T_{\eta^1}, \s)$. 

The maps $f_{\eta^0, \dots, \eta^r}$ can also be defined when all the curve collections $\etas^0, \dots, \etas^r$ are strongly equivalent. The definition is completely analogous, except there is no need for the bars on the values $\s$; the image should be in $\Am(\T_{\eta^0}, \T_{\eta^r}, \s_1 + \dots+ \s_r)$, where $\s_1 + \dots + \s_r \in \zz^m$. 

For simplicity,  we will ignore the subscripts on the maps $f_{\eta^0, \dots, \eta^r}$, and denote them all by $f$. The maps $f$ satisfy a generalized associativity property, which can be written as
\begin {equation}
\label {eq:f}
 \sum_{0 \leq i < j \leq r} f(\x_1, \dots, \x_i, f(\x_{i+1}, \dots, \x_j), \x_{j+1}, \dots, \x_r ) = 0,
 \end {equation}
for any $\x_i \in \Am(\T_{\eta^{j-1}}, \T_{\eta^j}, \s_j), j=1, \dots, r$. Compare Equation (9) in \cite{BrDCov}.

\begin{proof}[Completion of the proof of Theorem~\ref{thm:LinkInvariance}]
As mentioned in Section~\ref{sec:hmoves}, the triangle and quadrilateral maps defined here are used to prove the part of Theorem~\ref{thm:LinkInvariance} that deals with invariance under handleslides. Indeed, one can follow the arguments of \cite[Section 9]{HolDisk} almost verbatim. The only difference
is that we use link Floer complexes, and we have to make sure that our maps are filtered. In particular, if $\betas$ and $\gammas$ are curve collections that differ from each other by an elementary handleslide as in \cite[Figure 9]{HolDisk}, instead of $\iCF^-(\betas, \gammas, \ss_0)$ we use $\Am(\Tb, \Tg, \zero)$, where $\zero=(0, \dots, 0)$ denotes the zero vector. Since all the $2^{g+k-1}$ generators are in Alexander grading zero, it is easy to see that $\Am(\Tb, \Tg, \zero)$ has trivial differential. Therefore, its homology is  that of a $(g+k-1)$-dimensional torus, over the base ring $\Ring(\he) = \Field[[U_1, \dots, U_k]]$. When constructing the triangle maps corresponding to this handleslide, we use the top degree generator $\Theta \in \Am(\Tb, \Tg, \zero)$. Observe that, when $\s' =\zero$, the map $f_{\alpha \beta \gamma}$ from \eqref{eq:Triangle} produces a filtered map 
$$\CFLm(\he)  \to \CFLm(\he''), \ \  \x \mapsto f_{\alpha \beta \gamma}(\x \otimes \Theta)$$
and hence chain maps from $\Am(\he, \s)$ to $\Am(\he'', \s)$. These can be shown to be chain homotopy equivalences by using (filtered) quadrilateral maps, along the lines of \cite[Section 9]{HolDisk}.   
It should be noted that in \cite[Section 9]{HolDisk} it was only proved that the maps induced by handleslides are quasi-isomorphisms. However, the stronger statement that they are chain homotopy equivalences follows from the existence of a chain homotopy between the triangle map given by a small isotopy, and a continuation map; see \cite[proof of Proposition
11.4]{LipshitzCyl},  \cite[proof of Theorem
6.6]{OzsvathStipsicz}, and the proof of Lemma~\ref{lemma:kan} below.
\end{proof}

\subsection{An alternative description} 
\label{sec:alternative}
In general, the generalized Floer complex $\Am(\he, \s)$ is not free as an $\Ring(\he)$-module. In Section~\ref{sec:general} we will construct a resolution of $\Am(\he, \s)$, denoted $\Chain^-(\he, \s)$. For now, we will focus on the special case when the Heegaard diagram $\he$ is link-minimal (cf. Definition~\ref{def:linkminimal}). In that case, $\Am(\he, \s)$ admits an alternative description as a free module $\Chain^-(\he, \s)$, as follows.

Since we assume that $\he$ is link-minimal, we have $\Ws_i = \{w_i\}$ and $\Zs_i =\{z_i\}$ for each $i=1, \dots, \ell$. We may also have some free basepoints $w_{\ell+1}, \dots, w_k$. For $i \in \{1, \dots, \ell \}$, $s \in  \bH(L)_i$, and $\phi \in \pi_2(\x, \y)$, we set
$$ E^i_s(\phi) = \begin {cases} 
n_{w_i}(\phi) & \text{ if } \Filt_i(\x) \leq s, \Filt_i(\y) \leq s \\
s-\Filt_i(\x) + n_{z_i}(\phi) & \text{ if } \Filt_i(\x) \leq s, \Filt_i(\y) \geq s \\
n_{z_i}(\phi) & \text{ if } \Filt_i(\x) \geq s, \Filt_i(\y) \geq s \\
\Filt_i(\x) - s + n_{w_i}(\phi) & \text{ if } \Filt_i(\x) \geq s, \Filt_i(\y) \leq s.
\end {cases} $$
 
Equivalently, we can write
\begin {eqnarray}
\label {eiz}
E^i_s(\phi) &=&  \max(s-A_i(\x), 0) - \max(s-A_i(\y), 0)+ n_{z_i}(\phi) \\
\label {eiw}
&=& \max(A_i(\x) -s, 0) - \max(A_i(\y) -s,0) +  n_{w_i}(\phi).
\end {eqnarray}

In particular, $E^i_{-\infty}(\phi) = n_{z_i}(\phi)$ and $E^i_{+\infty}(\phi) =n_{w_i}(\phi)$.

Given $\s = (s_1, \dots, s_\ell) \in \bH(L)$, we define the complex
$\Chain^-(\he, \s)$ as the free module over $\Ring(\he) = \Field[[U_1, \dots, U_k]]$ generated by $\Ta \cap \Tb$, and equipped with the differential: 
\begin{equation}
  \label{def:DefD}
  \partial \x = \sum_{\y \in \Ta \cap \Tb} \sum_{\substack{\phi \in \pi_2(\x, \y)\\ \mu(\phi)=1 }} \# (\M(\phi)/\R) \cdot U_{1}^{E^1_{s_1}(\phi)} \cdots  U_{\ell}^{E^\ell_{s_\ell}(\phi)}\cdot  U_{\ell+1}^{n_{w_{\ell+1}}(\phi)}\cdots U_{k}^{n_{w_k}(\phi)}\y. 
\end{equation}
We also put a $\zz$-grading $\mu_\s$ on $\Chain^- (\he, \s)$ by setting
\begin {equation}
\label {eq:ms}
 \mu_\s(\x) = M(\x) - 2\sum_{i=1}^\ell \max(A_i(\x)-s_i, 0),
  \end{equation}
at least when none of the values $s_i$ is $-\infty$. When some of the values $s_i$ are $-\infty$, we replace the corresponding expressions $ \max(A_i(\x)-s_i, 0)$ by $A_i(\x)$ in the formula \eqref{eq:ms}.

\begin{lemma}
\label{lem:AmC}
The complexes $\Am(\he, \s)$ and $\Chain^-(\he, \s)$ are isomorphic over $\Ring(\he)$.
\end{lemma}

\begin{proof}
When none of the values $s_i$ is $-\infty$, the isomorphism $\Am(\he, \s) \xrightarrow{\cong} \Chain^-(\he, \s)$ is given by
\begin {equation}
\label {eq:mapsto}
 \x \mapsto  U_1^{\max(A_1(\x)-s_1, 0)}\cdots U_\ell^{\max(A_\ell(\x)-s_\ell, 0)}  \x.
 \end{equation}
When some of the values $s_i$ are $-\infty$, we replace $ \max(A_i(\x)-s_i, 0)$ by $A_i(\x)$ in the above formula.
\end{proof}

\begin{remark}
Although we could identify the complexes $\Am(\he, \s)$ and $\Chain^-(\he, \s)$ using the isomorphism from Lemma~\ref{lem:AmC}, in practice we will use the notation $\Am(\he, \s)$ when we think of the subcomplex of $\CFLm(\he)$, and $\Chain^-(\he, \s)$ when we think of the alternative description above, as a free complex.   
\end{remark}

\begin{remark}
When the Heegaard diagram $\he$ is not link-minimal, we have $m$ variables $U_i$, but only $\ell$ Alexander gradings. If we set the variables $U_i$ corresponding to the same link component to be equal to each other, then one can define a complex by a formula similar to \eqref{def:DefD}. However, this complex would not be well-behaved with Heegaard move (vi) from Section~\ref{sec:hmoves}, the index zero/three link stabilization. Indeed, under this move, the generalized Floer complex $\Chain^-_*(\he, \s)$ would turn into the direct sum $\Chain^-_{*+1}(\he, \s') \oplus \Chain^-_*(\he, \s)$, where $\s'$ differs from $\s$ by $-1$ in the spot corresponding to the link component that is being stabilized.
\end{remark}

Recall from Remark~\ref{rem:gen} that the complexes $\Am(\he, \s)$ can also be defined when $\he$ represents an unlink in $\#^n(S^1 \times S^2)$. When $\he$ is also link-minimal, we can construct isomorphic free complexes $\Chain^-(\he, \s)$ in that context as well. Furthermore, the polygon maps from \eqref{eq:feta} can be re-written in terms of the free complexes as
\begin {equation}
\label {eq:fetaFree}
 f_{\eta^0, \dots, \eta^r} : \bigotimes_{j=1}^r \Chain^-(\T_{\eta^{j-1}},\T_{\eta^j}, \s_j) \to \Chain^-(\T_{\eta^0}, \T_{\eta^r}, \s_1 + \dots+ \s_r)
 \end {equation}
$$ f_{\eta^0, \dots, \eta^r}(\x_1 \otimes \dots \otimes \x_r)= \sum_{\y \in \T_{\eta^{0}} \cap \T_{\eta^r}}\sum_{\{\phi \in \pi_2(\x_1, \dots, \x_r, \y)| \mu(\phi)=2-r\} } \# (\M(\phi)) \cdot  \U^{E_{\s_1, \dots, \s_r}(\phi)} \y,$$ 
where we use the notation
$$ \U^{E_{\s_1, \dots, \s_r}}(\phi)=U_1^{E^1_{s_{1,1}, \dots, s_{r,1}}(\phi)} \dots U_{\ell}^{E^{\ell}_{s_{1,\ell}, \dots, s_{r,\ell}}(\phi)} \cdot  U_{\ell+1}^{n_{w_{\ell+1}}(\phi)}\cdots U_{k}^{n_{w_k}(\phi)},$$
with
$$ \s_j=(s_{j,1}, \dots, s_{j, \ell}), \ j=1, \dots, r$$
and
\begin{align*} E^i_{s_{1,i}, \dots, s_{r,i}}(\phi) &= \sum_{j=1}^r \max(s_{j,i}-A_i(\x_j), 0) - \max(\sum_{j=1}^r s_{j,i} - A_i(\y), 0) + n_{z_i}(\phi) \\
&= 
\sum_{j=1}^r \max(A_i(\x_j)-s_{j,i}, 0) - \max(A_i(\y)-\sum_{j=1}^r s_{j,i}, 0) + n_{w_i}(\phi).
\end{align*}

\subsection {Reduction}
\label {sec:reduction}
Suppose that $M$ is a sublink of $L = L_1 \amalg \dots \amalg L_\ell$. We choose an orientation on $M$ (possibly different from the one induced from $\orL$), and denote the corresponding oriented link by $\orM$. We let $I_+(\orL, \orM)$ (resp. $I_-(\orL, \orM)$) be the set of indices $i$ such that the component $L_i$ is in $M$ and its orientation induced from $\orL$ is the same as (resp. opposite to) the one induced from $\orM$. 

For $i \in \{1, \dots, \ell\}$, we define a projection map $p_i^{\orM} : \bH(L)_i \to \bH(L)_i $ by
$$ p_i^{\orM}(s) = 
\begin{cases}
+\infty & \text{ if } i \in I_+(\orL, \orM), \\
-\infty & \text{ if } i \in I_-(\orL, \orM), \\
s & \text{ otherwise.}
\end {cases}$$

Then, for $\s = (s_1, \dots, s_\ell) \in \bH(L)$, we set
 $$p^{\orM} (\s)= \bigl(p_1^{\orM}(s_1), \dots, p_\ell^{\orM}(s_\ell)\bigr).$$

 Set $N = L-M$.  We define a map 
\begin{equation}
\label{eq:psiorm}
 \psi^{\orM} : \bH(L) \longrightarrow \bH(N)
 \end{equation}
as follows. The map $\psi^{\orM}$ depends only on the summands $\bH(L)_i$ of $\bH(L)$ corresponding to $L_i \subseteq N$. Each of these $L_i$'s appears in $N$ with a (possibly different) index $j_i$, so there is a corresponding summand $\bH(N)_{j_i}$ of $\bH(N)$.  
We then set
\begin {equation}
\label {eq:psy}
 \psi^{\orM}_i : \bH(L)_i \to \bH(N)_{j_i}, \ \ s_i \to s_i - \frac{\lk(L_i, \orM)}{2},
 \end {equation}
where $L_i$ is considered with the orientation induced from $L$, while $\orM$ is with its own orientation. We then define $\psi^{\orM}$ to be the direct sum of the maps $\psi_i^{\orM}$, pre-composed with projection  to the relevant factors.

\begin {remark}
If we view $\bH(L)$ as a lattice over $H_1(S^3 -L)$, see Remark~\ref{rem:h1}, we can describe the map $\psi^{\orM}$ as
$$ \psi^{\orM}(\s) = \s - \frac{[\orM]}{2},$$
where we denote elements in $H_1(S^3 -L)$ the same as their inclusions into $H_1(S^3 - N)$.   
\end {remark}

\begin {definition}
\label {def:reduce}
Let $\he$ be a multi-pointed  Heegaard diagram representing a link $\orL \subset Y$. Let $M \subseteq L$ be a sublink, with an orientation $\orM$ (not necessarily the one induced from $\orL$). The {\em reduction} of $\he$ at $\orM$, denoted $r_{\orM} \he$, is the  Heegaard diagram for $\orL - M$ obtained from $\he$ as follows:  first, we delete the basepoints $z$ from all components of $M \subseteq L$ oriented the same way in $\orL$ as in $\orM$; second, we delete the basepoints $w$, and relabel the basepoints $z$ as $w$,  from all components of $M$ oriented the opposite way in $\orL$ as in $\orM$. 
\end {definition}

\begin{remark}
The reduction of a link-minimal Heegaard diagram is link-minimal. By contrast, the reduction of a minimally pointed diagram is usually not minimally pointed (because we introduce new free basepoints).
\end{remark}

Using the interpretation of $\H(L)$ as a space of relative $\spc$ structures (cf. Remark~\ref{rem:h1}) it follows from \cite[Section 3.7]{Links} that there is an identification:
\begin {equation}
\label {eq:red_lmn0}
\begin {CD}
\Am(\he, p^{\orM}(\s)) @>{\cong}>> \Am(r_{\orM}(\he), \psi^{\orM}(\s)).
\end {CD}
\end {equation}

When the diagram $\he$ is link-minimal, we can write this identification in terms of the free complexes from Section~\ref{sec:alternative}:
 \begin {equation}
\label {eq:red_lmn}
\begin {CD}
\Chain^-(\he, p^{\orM}(\s)) @>{\cong}>> \Chain^-(r_{\orM}(\he), \psi^{\orM}(\s)).
\end {CD}
\end {equation}

\subsection{Inclusion maps}
\label{sec:inclusions}

Let $\orL \subset Y$ and $M \subseteq L$ be as in the previous subsection. Let also $\he$ be a generic, admissible, multi-pointed Heegaard diagram for $\orL$. If $M$ is equipped with  the orientation induced from $\orL$, then for $\s \in \bH(L)$ we have a natural inclusion map between subcomplexes of $\CFLm(\he)$:
$$ \Pr^{M}_\s : \Am(\he, \s) \to \Am (\he, p^{M}(\s)).$$

If $\he$ is link-minimal, we can write $\Pr^{M}_\s$ in terms of the alternative description of generalized Floer complexes from Section~\ref{sec:alternative}. In that setting, we have the formula
$$ \Pr^{M}_\s : \Chain^-(\he, \s) \to \Chain^- (\he, p^{M}(\s)),$$
\begin{equation}
\label{eq:promo}
 \Pr^{M}_\s \x =  \prod_{\{i \mid L_i \subseteq  L-M\} } U_{i}^{\max(A_i(\x) - s_i, 0)}  \cdot \x.
 \end{equation}

Now suppose that $M$ is equipped with an arbitrary orientation $\orM$, not necessarily the one from $\orL$. We use the notation $I_{\pm}(\orL, \orM)$ from the previous subsection. If $\he$ is link-minimal then, by analogy with \eqref{eq:promo}, we define an inclusion map
$$ \Pr^{\orM}_\s : \Chain^-(\he, \s) \to \Chain^- (\he, p^{\orM}(\s))$$
by
\begin {equation}
\label {eq:proj}
\Pr^{\orM}_\s \x =  \prod_{i \in I_+(\orL, \orM)} U_{i}^{\max(A_i(\x) - s_i, 0)} \cdot \prod_{i \in I_-(\orL, \orM)} U_{i}^{\max(s_i - A_i(\x), 0)} \cdot \x. 
\end {equation}

Alternatively, in terms of the original description, we have
\begin {equation}
\label {eq:proj3}
\Pr^{\orM}_\s : \Am(\he, \s) \to \Am(\he, p^{\orM}(\s)), \ \ \Pr^{\orM}_\s \x =  \prod_{i \in I_-(\orL, \orM)} U_{i}^{s_i - A_i(\x)} \cdot \x. 
\end {equation}

Note that the expressions \eqref{eq:promo}, \eqref{eq:proj} and \eqref{eq:proj3} are well-defined only when the exponents of the $U_{i}$ variables involved there are finite. That is, we need to require that $s_i \neq -\infty$ for all $i \in I_+(\orL, \orM)$, and $s_i \neq +\infty$ for all $i \in I_-(\orL, \orM)$. These conditions will always be satisfied when we consider inclusion maps in this paper.

Observe that $\Pr^{\orM}_\s$ is a chain map, and it shifts the grading \eqref{eq:ms} by a definite amount:
\begin {equation}
\label {eq:mumu}
 \mu_{p^{\orM}(\s)}(\Pr^{\orM}_\s(\x)) =  \mu_\s(\x) - 2\sum_{i \in I_-(\orL, \orM)} s_i.
\end {equation}

\begin {remark}
\label {remark:knot2}
When $L=K$ is a knot, the maps $\Pr^{\vec K}_s, \Pr^{-\vec K}_s$ correspond to the inclusions $v_{s}^-$ and $h_{s}^-$ of $A_s^-$ into the subcomplexes $C(\{i \leq 0\})$ and $C(\{j \leq s \})$, respectively; compare Remark~\ref{remark:knot1}. In \cite{IntSurg}, there are analogous maps $v^+_s, h^+_s : A_s^+ \to B^+ \cong \iCF^+(Y)$.
\end {remark}

\begin{remark}
\label{rem:notLM}
When $\he$ is not link-minimal, there does not seem to be a natural way of defining maps $ \Pr^{\orM}_\s : \Am(\he, \s) \to \Am(\he, p^{\orM}(\s))$ for arbitrary orientations $\orM$. This is one of the reasons why, in Section~\ref{sec:general}, we will construct resolutions $\Chain^-(\he, \s)$ of $\Am(\he, \s)$; the resolutions will come with maps $\Pr^{\orM}_\s : \Chain^-(\he, \s) \to \Chain^- (\he, p^{\orM}(\s))$ that generalize \eqref{eq:proj}.
\end{remark}

\section {An overview of the link surgery formula}
\label {sec:ex}

According to the link surgery formula (Theorem~\ref{thm:FirstSurgery}), the Heegaard Floer homology groups $\HFm_*$ of a link surgery $Y_\Lambda(L)$ can be computed from the surgery complex $\C^-(\Hyper, \Lambda)$. Constructing the surgery complex is rather technical, and the full details will be postponed until Section~\ref{sec:statement}. Until then, the purpose of the present section is to give the reader a basic understanding of the link surgery formula, and to show how it can be applied to compute some examples.  

\subsection{Preliminaries}
\label{sec:prelim}
As mentioned in the Introduction, the input data for the surgery complex is called a complete system of hyperboxes $\Hyper$ for the link $\orL$. In this section, we will focus on the simplest type of complete system, called {\em basic}. Complete systems of hyperboxes will be defined in Section~\ref{sec:hyperHeegaard}, and basic ones in Subsection~\ref{sec:basic}. For now, it suffices to note that a basic system exists for every link, and is determined by the choice of a single Heegaard diagram $\Hyper^L$ for $L$. The diagram $\Hyper^L$ has to be basic in the sense of Definition~\ref{def:hbasic}; that is, it should have no free basepoints, and every link component should intersect the Heegaard surface in exactly two basepoints (one $w_i$ and one $z_i$), situated on each side of a beta curve. 

Apart from $\Hyper^L$, the other diagrams appearing in a basic complete system $\Hyper$ are the reductions
$$ \Hyper^{L-M} := r_{M} (\Hyper^L),$$
obtained from $\Hyper^L$ by deleting the $z$ basepoints on the sublink $M$; cf. Section~\ref{sec:reduction}. (In our case, $M$ is given the orientation induced from $\orL$.) Note that all the diagrams $\Hyper^{L-M}$ are link-minimal.

To the diagrams $ \Hyper^{L-M}$ we associate generalized Floer complexes $\Chain^-(\Hyper^{L - M}, \s)$, as in Section~\ref{sec:alternative}. These are modules over the ground ring
$$\Ring := \Ring(\Hyper^{L}) = \Field[[U_1, \dots, U_{\ell}]].$$

With this in mind, recall from \eqref{eq:chl_zero} that, as a $\Ring$-module, the surgery complex is the infinite direct product
\begin{equation}
\label{eq:CHL}
 \C^-(\Hyper, \Lambda) = \bigoplus_{M \subseteq L} \prod_{\s \in \H(L)}  \Chain^-(\Hyper^{L - M}, \psi^{M}(\s) ).
 \end{equation}
Here, $\psi^{M}$ is the map from \eqref{eq:psiorm}, with $M$ again having the orientation induced from $\orL$. 

To simplify the notation somewhat, we denote a typical term in the chain complex \eqref{eq:CHL} by
\begin {equation}
\label {eq:ces}
 \C^\eps_\s = \Chain^-(\Hyper^{L-M}, \psi^M(\s)), 
 \end {equation}
where $\eps=\eps(M)= (\eps_1, \dots, \eps_{\ell}) \in \{0,1\}^\ell$ is such that $L_i \subseteq M$ if and only if $\eps_i = 1$. 

Furthermore, as noted in \eqref{eq:DHL}, the differential on the complex $ \C^-(\Hyper, \Lambda)$ is given by
$$\D^-(\s, \x) = \sum_{N \subseteq L - M} \sum_{\orN \in \Omega(N)} (\s + \Lambda_{\orL, \orN}, \Phi^{\orN}_{\psi^{M}(\s)}(\x)),$$
for $\s \in \H(L)$ and $\x \in \Chain^-(\Hyper^{L-M}, \psi^M(\s))$. In this formula, the maps 
$$\Phi^{\orN}_{\psi^{M}(\s)}: \Chain^-(\Hyper^{L-M}, \psi^M(\s)) \to \Chain^-(\Hyper^{L-M-N}, \psi^{M\cup N}(\s))$$
are constructed from polygon maps of the type considered in Section~\ref{sec:polygon}; cf. also Equation~\eqref{eq:fetaFree}. We shall not give their exact definition here, but we note that when $N = \emptyset$, the map $\Phi^{\orN}_{\psi^{M}(\s)}$ is just the usual differential on $\Chain^-(\Hyper^{L-M}, \psi^M(\s))$, counting holomorphic disks.

Let us denote a typical summand in the differential $\D^-$ by
\begin {equation}
\label {eq:dees}
\D^{\eps, \eps'}_{\eps^0, \s} = \Phi^{\orN}_{\psi^M(\s)} : \C^{\eps^0}_\s \to \C^{\eps^0 + \eps}_{\s + \eps' \cdot \Lambda},
 \end {equation}
where $\eps^0=\eps(M), \eps=\eps(N)$, and $\eps' \in \{0,1\}^{\ell}$  is such that $i \in I_-(\orL, \orN)$ if and only if $\eps'_i = 1$. The dot product $\eps' \cdot \Lambda$ denotes the vector $\sum \eps'_i\Lambda_i$. Furthermore, as in Section~\ref{sec:reduction}, $I_-(\orL, \orN)$ is the set of indices $i$ such that the orientation of $L_i \subset M$ is different in $\orN$ from the one in $\orL$.

Whenever we drop a subscript or superscript from our notation, we will mean the direct product (or sum, as the case may be) over all possible values of that subscript or superscript. For example, $\C^\eps = \prod_\s \C^\eps_\s$, and $\C = \oplus_\eps \C^\eps = \C^-(\Hyper, \Lambda)$. 

In practice, if one wants to compute the Heegaard Floer homology groups $\HFm(Y_\Lambda(L), \ux)$ using Theorem~\ref{thm:FirstSurgery}, it is helpful to replace the infinite direct product from \eqref{eq:CHL} with a finite one. We refer to this procedure as {\em horizontal truncation}. It can be done along the same lines as the corresponding argument in the case of surgery on knots, from  \cite[Section 4.1]{IntSurg}. There are many ways of doing horizontal truncation. A general procedure will be given in Section~\ref{sec:truncate}, but until then, we will explain a few ad hoc variants in the examples below.

The organization of the rest of Section~\ref{sec:ex} is as follows. In Section~\ref{sec:knots} we specialize the link surgery formula to the case of knots, and explain the relation to the already-existing knot surgery formula from \cite{IntSurg}. In Section~\ref{sec:unknot} we explain how horizontal truncation works in the simple case of $+1$ surgery on the unknot, and we also justify our use of direct products rather than direct sums. In Section~\ref{sec:twoC}, we turn our attention to two-component links, and describe the surgery complex in more detail in that case. Finally, in Section~\ref{sec:hopf}, we show how the truncation technique can be used to arrive at an explicit computation in the case of surgeries on the Hopf link.

\subsection{Knots}
\label{sec:knots}

Let $\orL = \orK$ be an oriented knot in a homology three-sphere $Y$. In this case, only two Heegaard diagrams appear in the surgery complex: a basic diagram $\Hyper^{K}$ for $K$ (with one $w$ basepoint and one $z$ basepoint), and its reduction $\Hyper^{\emptyset}$, a diagram for $Y$ that is obtained from $\Hyper^{K}$ by deleting $z$. The lattice $\H(L)$ is just $\zz$, and $\H(\emptyset)$ consists of the single element $0$. Further, let $U=U_1$ be the unique variable in the ground ring $\Ring=\Field[[U]]$, and let $m \in \zz$ be the surgery coefficient that specifies the framing $\Lambda$. As a $\Ring$-module, the surgery complex $\C =  \C^-(\Hyper, m)$  takes the form
$$ \C= \prod_{s \in \zz} \C^0_\s \oplus \prod_{s \in \zz} \C^1_s.$$
Here, $\C^0_s$ is the generalized Floer complex $\Chain^-(\Hyper^K, s)$, and each $\C^1_s$ is a copy of the complex
$$\Chain^-(\Hyper^{\emptyset}, 0) = \CFm(\Hyper^{\emptyset}),$$
whose homology is $\HFm(Y)$.

With regard to the differential, this is composed of the following terms: $\D^{0,0}_{0,s}$ (which is the usual differential on $\C^0_s$), $\D^{0,0}_{1,s}$ (which is the usual differential on $\C^1_s$), as well as
$$ \D^{1,0}_{0,s} = \Phi^{K}_s : \C^0_s \to \C^1_s$$
and
$$ \D^{1,1}_{0,s} = \Phi^{-K}_s : \C^0_s \to \C^1_{s+m}.$$

Thus, the complex $\C$ can be viewed as the mapping cone of the map
\begin{equation}
\label{eq:MCone}
  \prod_{s \in \zz} \C^0_\s \ \to \ \prod_{s \in \zz} \C^1_s, \ \ \ (s, \x) \mapsto (s, \Phi^{K}_s) + (s+m, \Phi^{-K}_s).
  \end{equation}

Let us be more specific about what the maps $ \Phi^{K}_s$ and $\Phi^{-K}_s$ are in this case. The map $ \Phi^{K}_s$ is just the inclusion $\I^{K}_s$ defined in Section~\ref{sec:inclusions}, whose codomain $\Chain^-(\Hyper^K, +\infty)$ is identified with $\Chain^-(\Hyper^{\emptyset}, 0)$ by \eqref{eq:red_lmn}. 

The map $\Phi^{-K}_s$ is also constructed from the corresponding inclusion $\I^{-K}_s$, but note that the codomain of $\I^{-K}_s$ is the complex $\Chain^-(\Hyper^K, -\infty) \cong \Chain^-(r_{-K}(\Hyper^K), 0)$. The diagram $r_{-K}(\Hyper^K)$ is obtained from $\Hyper^K$ by deleting $w$ and relabeling $z$ as the new $w$. This diagram represents the underlying three-manifold $Y$, just like $\Hyper^{\emptyset}$, so it can be related to the latter by a sequence of Heegaard moves. In fact, in a basic system, we use a specific set of Heegaard moves, which will be described in Section~\ref{sec:basic}. The Heegaard moves induce a chain homotopy equivalence from $\Chain^-(r_{-K}(\Hyper^K))$ to $\Chain^-(\Hyper^{\emptyset}, 0)$. We define $\Phi^{-K}_s$ to be the composition of $\I^{-K}_s$ with this homotopy equivalence.

The mapping cone \eqref{eq:MCone} is strongly reminiscent of the mapping cone 
$$ \bigoplus_{s\in\zz} A_s^+ \ \to \ \bigoplus_{s\in \zz} B_s^+, \ \ (s, \x) \mapsto (s, v_s^+(\x)) + (s+m, h_s^+(\x))$$
that appears in the knot surgery formula from \cite{IntSurg}. The main difference is that our complex $\C$ computes the completed, minus flavor of Heegaard Floer homology, $\HFm(Y_m(K))$, whereas the one in \cite{IntSurg} computed the plus flavor $\HF^+(Y_m(K))$. As noted in Remark~\ref{remark:knot1}, the complexes $\C^0_s$ are the completions of $A_s^- = C\{\max(i, j-s) \leq 0\}$, in the notation of \cite{IntSurg}. Similarly, the complexes $\C^1_s = \CFm(Y)$ are the analogues of $B_s^+=\CF^+(Y)$. Moreover, as explained in Remark~\ref{remark:knot2}, the maps $\Phi^{K}_s$ and $\Phi^{-K}_s$ are the analogues of $v_s^+$ and $h_s^+$ from \cite{IntSurg}.

\subsection {Remarks on direct products and the unknot}
\label {sec:unknot}

The reader who is familiar with the knot surgery formula from \cite{IntSurg} may wonder why our definition of the complex $\C$ in \eqref{eq:CHL} involved a direct product rather than a direct sum. The results in \cite{IntSurg} were phrased in terms of direct sums, but they only applied to $\iHF^+$ and $\widehat{\iHF}$. In the case of $\HFm$ direct sums do not give the right answer. This can be seen even in the simple case of $+1$ surgery on the unknot $\vec U$ in $S^3$. 

Specifically, let us consider a genus one Heegaard diagram for the unknot, with one alpha curve and one beta curve intersecting transversely at a single point, and two basepoints $w$ and $z$. This is a basic Heegaard diagram, to which we can associate a basic complete system of hyperboxes. The complex $\C$ for $+1$ surgery is
\begin{equation}
\label{eq:Unknot}
 \C = \bigoplus_{\eps \in \{0,1\}} \prod_{s \in \zz} \C^{\eps}_s.
 \end{equation}
Here, if we use the alternative description of the complexes $\C^0_s = \Chain^-(\Hyper^K, s)$ from Section~\ref{sec:alternative}, then each such complex is freely generated by the single intersection point, so it is a copy of $\ff[[U]]$. The complexes $\C^0_s$ are also copies of $\ff[[U]]$. Let $a_s \in \C^0_s$ and $b_s \in \C^1_s$ be the generators of each piece. 

The differential $\D$ is obtained as the sum of maps
$$ \Phi^{\vec U}_s: \C^{0}_s \longrightarrow \C^{1}_s \ \ \text{and} \ \ \Phi^{-\vec U}_s: \C^{0}_s \longrightarrow \C^{1}_{s+1}.$$

The maps $\Phi^{\vec U}_s$ are inclusion maps $\Pr^{\vec U}_s$, which are given by multiplication by a certain power of $U$; cf. Equation~\eqref{eq:proj}. The maps $\Phi^{-\vec U}_s$ are inclusion maps $\Pr^{-\vec U}_s$ post-composed with a homotopy equivalence, which in this case is the identity.  Thus, using \eqref{eq:proj}, we obtain:
\begin {equation}
\label {eq:phiu}
 \Phi^{\vec U}_s = \begin{cases}
1 & \text{if } s \geq 0 \\
U^{-s} & \text{if } s \leq 0,
\end {cases} \ \ \  \ \  \ \  \ \  \ \ 
\Phi^{-\vec U}_s = \begin{cases}
U^s & \text{if } s \geq 0 \\
1 & \text{if } s \leq 0.
\end {cases}
\end{equation}
(Compare Section 2.6 in \cite{IntSurg}.) The homology of the complex $\C$ is then isomorphic to $\ff[[U]]$, being freely generated by the element  
$$ \sum_{s \in \zz}  U^{|s|(|s|-1)/2} a_s .$$

The Heegaard Floer homology of $+1$ surgery on the unknot is $\HFm(S^3) \cong \ff[[U]]$, so Theorem~\ref{thm:FirstSurgery} gives the right answer. However, if instead of the direct product in \eqref{eq:Unknot} we had used a direct sum, the homology of the resulting complex would have been a more complicated $\ff[[U]]$-module.  Indeed, the map $\D$ would then have nontrivial cokernel in the new $\C^1=\bigoplus \C^1_s$. The cokernel would be generated as a $\ff[[U]]$-module by classes $[b_i], i \in \zz$, subject to the relations:
$$[b_0]= [b_1] = U[b_{-1}] = U[b_{2}] = U^3[b_{-2}] = U^3[b_3]=\dots$$ 
The basic reason why direct sums are not suitable for the Surgery Theorem for $\HFm$ is that they do not behave as well with respect to filtrations as direct products do. As we shall see soon, filtrations play an important role in the truncation procedure, which in turn is essential for the proof of Theorem~\ref{thm:FirstSurgery}.

Let us fix some terminology. For us, a {\em filtration} $\F$ on a $\Ring$-module $\mo$ is a collection of $\Ring$-submodules $\{\F^i(\mo) \mid i \in \zz \}$ of $\mo$ such that $\F^i(\mo) \subseteq \F^j(\mo)$ for all $i \leq j$.  A filtration is called {\em bounded above} if  $\F^i(\mo) = \mo$ for $i \gg 0$, and {\em bounded below} if $\F^i(\mo) = 0$ for $i \ll 0$. A filtration is called {\em bounded} if it is both bounded above and bounded below.

If $\mo$ is equipped with a differential $\del$ that turns it into  a chain complex, we say that the chain complex $(\mo, \del)$ is filtered by $\F$ if $\del$ preserves each submodule $\F^i(\mo)$. The {\em associated graded} complex $\gr_{\F} \mo$ is defined as
\begin {equation}
\label {eq:assocgraded}
 \gr_{\F} (\mo) = \bigoplus_{i \in \zz} \bigl ( \F^i(\mo)/\F^{i-1}(\mo) \bigr ),
 \end {equation}
equipped with the differential induced from $\F$.

If $\F$ is a bounded filtration on a chain complex $(\mo, \del)$, a standard result from homological algebra says that if $\gr_{\F}(\mo)$ is acyclic, then $\mo$ is acyclic as well. (Note that this can fail for filtrations that are not bounded.)

A standard way to construct bounded filtrations is as follows. If $\mo$ is freely generated over $\Ring$ by a collection of generators $G$, a bounded map $ \F: G \to \zz$
defines a bounded filtration on $\mo$ by letting $\F^i(\mo)$ be the submodule generated by the elements $g\in G$ with $\F(g) \leq i$. 

Suppose now that we have a direct product of $\Ring$-modules  
$$\mo = \prod_{s \in S} \mo_s,$$ 
indexed over a countable set $S$. Suppose further that each $\mo_s$ is
a free module over $\Ring$ with a set of
generators $G_s$. Assume that $\mo$ is equipped with a differential
$\del$.
(Typically, in our examples, each $G_s$ is finite, and 
the differential $\del$ is 
defined on each term so that it is locally finite, as in
Definition~\ref{def:locfin}.)

In this situation, an assignment
$$\F : \bigcup_{s \in S} G_s \to \zz$$
specifies bounded filtrations on each $\mo_s$. Together these produce a filtration on $\mo$ given by
$$\F^i(\mo) = \prod_{s \in S} \F^i(\mo_s).$$
This filtration on $\mo$ is generally neither bounded above nor bounded below. It is bounded above provided  that there exists $i \gg 0$ such that $\F^i(\mo_s) = \mo_s$ for all $s$; that is, if $\F(g) \leq i$ for all $g \in G_s, s \in S$. We have the following:

\begin {lemma}
\label {lem:acyclicAG}
Consider a module $\mo = \prod_{s \in S} \mo_s$, where each $\mo_s$ is a freely generated over $\Ring$ by a  set of generators $G_s$. Suppose $\F: \bigcup_{s \in S} G_s \to \zz$ defines a filtration on $\mo$ that is bounded above. Further, suppose $\mo$ is equipped with a differential $\del$, and that the associated graded complex $\gr_{\F}(\mo)$ is acyclic. Then $\mo$ itself is acyclic.
\end {lemma}

\begin {proof}
Consider the spectral sequence associated to $\F$, whose $E^1$ term is $H_*(\gr_{\F}(\mo)) = 0$. According to \cite[Theorem 5.5.10]{Weibel}, the spectral sequence converges to $H_*(\mo)$ if the the filtration is complete, exhaustive, and the spectral sequence is regular. Completeness of the filtration means that $\mo= \varprojlim \mo/{\F}^i\mo$, which is true because $\mo$ is constructed as a direct product. (Note that for a complete filtration we have $\cap_i \F^i\mo=0$.) Exhaustiveness is a weaker condition than being bounded above. Regularity (as defined in \cite[Definition 5.2.10]{Weibel}) is automatic when the $E^1$ term is zero, because the higher differentials have to be zero as well.  
\end {proof}

One can apply Lemma~\ref{lem:acyclicAG} to calculate $H_*(\C^-(\vec U, 1))$, where $\C = \C^-(\vec U, 1)$ is the surgery complex for $+1$ surgery on the unknot, considered above. (Of course, this calculation can also be done directly.) As a module, the complex $\C$ is the direct product of all $\C^\eps_s$, over $\eps \in \{0,1\}$ and $s \in \zz$. Each term in the direct product is freely generated by an element $a_s$ or $b_s$. We have two filtrations $\F_0$ and $\F_1$ on $\C$ defined by:
$$ \F_0(a_s) = \F_0(b_s) = -s, $$
$$ \F_1(a_s) = \F_1(b_{s+1}) = s.$$

The map $\Phi_s^{\vec U}$ preserves the $\F_0$-grading and decreases $\F_1$ by one, whereas $\Phi_s^{-\vec U}$ decreases $\F_0$ by one and preserves $\F_1$. Neither $\F_0$ nor $\F_1$ are bounded above. However, we can consider the following subcomplexes of $\C$:
$$ \C_{>0} = \prod_{s>0} (\C^0_s \oplus \C^1_s) \ \ \ \text{and} \ \ \ \C_{<0} =  \prod_{s<0} (\C^0_s \oplus \C^1_{s+1}).$$

Then the restriction of $\F_0$ to $\C_{>0}$ and the restriction of $\F_1$ to $\C_{<0}$ are both bounded above. Further, in the associated graded of these subcomplexes the differential cancels out all the terms in pairs. Applying Lemma~\ref{lem:acyclicAG} we get that
$$ H_*(\C_{>0}) \cong H_*(\C_{<0}) \cong 0.$$

The quotient complex $\C/(\C_{>0} \oplus \C_{<0})$ is simply $\C^0_0 \cong \ff[[U]]$ with a trivial differential. From the corresponding long exact sequence we get that
$$ H_*(\C) \cong H_*(\C^0_0) \cong \ff[[U]],$$
as expected. 

Although Lemma~\ref{lem:acyclicAG} can be successfully used to compute the homology of the complex associated to $+1$ surgery on the unknot, if we try to do $-1$ surgery instead, we will need to work with filtrations that are bounded below, rather than above. As the following example shows, the exact analogue of Lemma~\ref{lem:acyclicAG} fails for bounded below filtrations. 

\begin {example}
Let $\mo = \prod_{s\geq 0} \mo_s$, where each $\mo_s$ is a free module over $\ff[[U]]$ of rank two, with generators $a_s$ and $b_s$. We equip $\mo$ with the locally finite differential $\del$ given by $\del a_0=b_0$,  $\del a_s = b_s + b_{s-1}$ for $s > 0$, and $\del b_s = 0$, for all $s$. Let $\F$ be the filtration on $\mo$ defined by $\F(a_s)=\F(b_s)=s$. Then $\F$ is bounded below but not bounded above, and $\gr_\F \mo$ is acyclic. However, $\mo$ itself is not acyclic, because $\sum_s a_s$ is in the kernel of its differential.
\end {example}

Nevertheless, we can obtain the desired result by imposing an additional condition on our filtration:
 
\begin{definition}
\label{def:utame}
Suppose we have a module $\mo = \prod_{s \in S} \mo_s$, where each $\mo_s$ is a freely generated over $\Ring=\ff[[U_1, \dots, U_\ell]]$ by a finite set of generators $G_s$. Suppose $\mo$ is equipped with a locally finite differential $\del$, and with a filtration induced by an assignment $\F: \bigcup_{s \in S} G_s \to \zz$. For $i \in \zz$, let $\mo^{[i]} \cong \F^i(\mo)/\F^{i-1}(\mo)$ be the free $\Ring$-module generated by those $g \in \cup_s G_s$ with $\F(g) =i$. Then, we say that the filtration $\F$ is {\em $U$-tame} if there exists $k \in \zz$ such that in the restriction of $\del$ to $\mo/\F^k \mo$, all  terms that drop the filtration level contain at least a nontrivial power of $U_j$ for some $j$; that is, for all $i > k$ we have
$$ \del(\mo^{[i]}) \subset \mo^{[i]} \oplus \bigl( (U_1, \dots, U_\ell) \cdot \F^{i-1}(\mo) \bigr).$$
\end {definition}

\begin{lemma}
\label{lem:acyclicAG2}
Let $\mo$ be a chain complex as in Definition~\ref{def:utame}, with a $U$-tame filtration $\F$. Suppose that the associated graded complex $\gr_{\F}(\mo)$ is acyclic. Then $\mo$ itself is acyclic. 
\end {lemma}

\begin {proof}
Let $k \in \zz$ be as in Definition~\ref{def:utame}. The complex $\F^k(\mo)$ (with the restriction of the filtration $\F$) satisfies the hypotheses of Lemma~\ref{lem:acyclicAG}, so it is acyclic.  Thus, it suffices to prove that the quotient complex $\mo'=\mo/\F^{k}(\mo)$ is acyclic.

If we view $\mo'$ as a complex over $\ff$ (rather than over $\Ring$), we can equip it with a second filtration $\uu$, given by the total degree of the $U_j$ variables. Precisely, if $g \in \cup_s G_s$ is a generator with $\F(g_s) > k$, we set
$$\uu(U_1^{m_1} \dots U_{\ell}^{m_\ell}\cdot g) = -(m_1+\dots+m_\ell).$$

Because of the original $\Ring$-structure on $\mo'$, we see that $\uu$ is indeed a filtration on $\mo'$ (viewed as an $\ff$-complex). The $U$-tameness assumption implies that the differential on the associated graded $\gr_\uu \mo'$ preserves each $\F$-filtration level $\mo^{[i]}$, $i > k$. Thus, $\gr_\uu \mo'$ splits as a direct sum of complexes of the form $\gr_\uu \mo^{[i]}$ for $i > k$.

We claim that each $\gr_\uu \mo^{[i]}$ is acyclic for $i > k$. We know that $\mo^{[i]}$ itself (being a summand of $\gr_{\F}(\mo)$) is acyclic. The short exact sequence 
$$ 0 \longrightarrow \mo^{[i]} \xrightarrow{U_1} \mo^{[i]} \longrightarrow \mo^{[i]}/U_1 \mo^{[i]} \longrightarrow 0$$
induces a long exact sequence in homology, showing that $ \mo^{[i]}/U_1 \mo^{[i]} $ is acyclic. By iterating this argument, we get that $ \mo^{[i]}/(U_1, \dots, U_\ell) \mo^{[i]}$ is acylic. Since  $\gr_\uu \mo^{[i]}$ is just a direct sum of complexes isomorphic to $ \mo^{[i]}/(U_1, \dots, U_\ell) \mo^{[i]}$, it must be acyclic as well.
We deduce from this that $\gr_\uu \mo'$ is acyclic. 

It now follows that $\mo'$ is acyclic, by Lemma~\ref{lem:acyclicAG}
applied to the filtration $\uu$. (Note that $\uu$ is bounded above by
$0$.)  Consequently  $\mo$ is also acyclic.
\end {proof}

As an application of Lemma~\ref{lem:acyclicAG2}, let us compute the homology of the complex $\C=\C^-(\vec U, -1)$, associated to $-1$ surgery on the unknot. Just as in the case of $+1$ surgery, our new complex $\C$ is the direct product (as a module) of all $\C^\eps_s$, over $\eps \in \{0,1\}$ and $s \in \zz$. Each term in the direct product is freely generated by one element, denoted $a_s$ (for $\C^0_s$) and $b_s$ (for $\C^1_s$). The maps $\Phi_s^{\vec U}$ and $\Phi_s^{-\vec U}$ are still given by the formulas \eqref{eq:phiu}, but now $\Phi_s^{-\vec U}$ takes $\C^0_s$ to $\C^1_{s-1}$. Explicitly, we have
$$\del a_s = \begin{cases}
b_s + U^s b_{s-1} & \text{if } s \geq 0 \\
U^{-s}b_s + b_{s-1} & \text{if } s \leq 0.
\end {cases}, \ \ \ \ \ \ \ \ \ \ \del b_s = 0.$$

Consider the following quotient complexes of $\C$:
$$ \C_{>0} = \prod_{s>0} (\C^0_s \oplus \C^1_s) \ \ \ \text{and} \ \ \ \C_{<0} =  \prod_{s<0} (\C^0_{s+1}  \oplus \C^1_{s}).$$

On $\C_{>0}$, we define a filtration $\F_0$ by
$$ \F_0(a_s) = \F_0(b_s) = s.$$
This filtration is not bounded above, but it is $U$-tame. Since its associated graded is acyclic, we get that $\C_{>0}$ is acyclic by Lemma~\ref{lem:acyclicAG2}.

Similarly, on $\C_{<0}$ we have a $U$-tame filtration $\F_1$ given by the assignment
$$ \F_1(a_s) = \F_1(b_{s-1}) = -s.$$
Its associated graded is acyclic, so $\C_{<0}$ is acyclic as well. We deduce that $\C$ is quasi-isomorphic to its subcomplex $\C^1_0$, which has homology $\ff[[U]]$.

\medskip
The calculations above (for $+1$ and $-1$ surgery on the unknot) are the simplest examples of horizontal truncation, and serve as models for the general case. In Section~\ref{sec:hopf} we will explain horizontal truncation in a more complicated example, that of the Hopf link.

\subsection{Two-component links}
\label{sec:twoC}
We now consider the case of a two-component link $\orL = L_1 \cup L_2 \subset Y$. We are interested in surgery on $L$ with surgery coefficients $p_1$ on $L_1$ and $p_2$ on $L_2$. The framing $\Lambda$ is given by the matrix
$$ \Lambda = \begin{pmatrix}
p_1 & c\\
c & p_2
\end{pmatrix},$$
where $c$ is the linking number between $L_1$ and $L_2$. We denote by $\Lambda_1$ and $\Lambda_2$ the two rows of $\Lambda$.

The affine lattice $\H(L)$ associated to $L$ is $\zz \times \zz$ if $c$ is even, and $(\tfrac{1}{2} + \zz) \times (\tfrac{1}{2} + \zz)$ if $c$ is odd. The lattices $\H(L_1)$ and $\H(L_2)$ are copies of $\zz$, and $\H(\emptyset) = 0$.

A basic complete system $\H$ for $\orL$ involves a basic diagram $\Hyper^L$ for $L$ (with four basepoints: $w_1$ and $z_1$ on $L_1$, and $w_2$ and $z_2$ on $L_2$), as well as its reductions 
$ \Hyper^{L_1}$ (obtained from $\H^L$ by deleting $z_2$), $\Hyper^{L_1}$ (obtained from $\Hyper^L$ by deleting $z_1$), and $\Hyper^{\emptyset}$ (obtained by deleting both $z_1$ and $z_2$). 

The surgery complex $\C$ is a ``double mapping cone'' of the form
\begin{equation}
\label{eq:Co2}
\xymatrix{
\displaystyle \prod_{\s \in \H(L)} \C^{00}_{\s}\ar[d] \ar[r] \ar[dr] & \displaystyle \prod_{\s \in \H(L)} \C^{10}_{\s} \ar[d] \\
\displaystyle \prod_{\s \in \H(L)} \C^{01}_{\s}\ar[r]  & \displaystyle \prod_{\s \in \H(L)} \C^{00}_{\s}
 }
\end{equation}
Here, the horizontal arrows consist of maps of the form $\Phi^{\pm L_1}_{\s}$, the vertical arrows of maps $\Phi^{\pm L_2}_{\s}$, and the diagonal of maps $\Phi^{\pm L_1 \cup \pm L_2}_\s$ and $\Phi^{\pm L_1 \cup \mp L_2}_\s$. Here, $+L_i$ (resp. $-L_i$) means $L_i$ with the orientation induced from (resp. opposite to) $\orL$. There are also maps $\Phi^{\emptyset}_\s$ (not shown in the diagram above), which represent the differentials on each complex $\C^{\eps}_{\s}$. 

The values of $\s$ in the targets of each map are shifted by an amount depending on the type of map and the framing $\Lambda$. Precisely, whenever we have a negatively oriented component $-\orL_i$ in the superscript of a map $\Phi_\s$, we add the vector $\Lambda_i \in \zz \times \zz$ to $\s$. Note that all the maps preserve the equivalence class of $\s$ in the quotient 
$$ \H(L) / \langle \Lambda_1, \Lambda_2 \rangle.$$
This quotient can be naturally identified with the space of $\spc$ structures on the  surgered manifold $Y_\Lambda(L)$; cf. Remark~\ref{rem:h1} and \cite[Section 3.7]{Links}. Therefore, the complex $\C$ splits as a direct sum, according to these $\spc$ structures. In fact, there is a more refined version of Theorem~\ref{thm:FirstSurgery}, which establishes an isomorphism between the homology of each of these summands and the corresponding $\HFm(Y_{\Lambda}(L), [\s])$; see Theorem~\ref{thm:surgery} in Section~\ref{sec:statement}. 

Let us now say a few words about the maps in \eqref{eq:Co2}. Those of the form $\Phi^{L_i}_{\s}$
are given by inclusions $\I^{L_i}_\s$, as in Section~\ref{sec:inclusions}. Those of the form $\Phi^{-L_i}_{\s}$ are the compositions of inclusions $\I^{-L_i}_\s$ with certain ``descent maps'' $D^{L_i}_{p^{L_i}(\s)}$. The general definition of descent maps will be given in Equation~\eqref{eq:des}, but in the case at hand they are chain homotopy equivalences induced by Heegaard moves between the respective Heegaard diagrams (just as in the case of knots, discussed in Section~\ref{sec:knots}). 

Finally, the diagonal maps in \eqref{eq:Co2} are chain homotopies between the corresponding compositions of vertical and horizontal maps. For example, $\Phi^{L_1 \cup - L_2}_\s$ is a chain homotopy between $\Phi^{L_1}_{\psi^{-L_2}(\s)} \circ \Phi^{- L_2}_\s$ and $ \Phi^{- L_2}_{\psi^{L_1(\s)}} \circ \Phi^{L_1}_{\s}$. These chain homotopies are constructed as compositions of inclusion and descent maps. In the case of a basic system for two component links, such as the one considered here, it turns out that most of the descent maps, and hence the corresponding diagonal maps in \eqref{eq:Co2}, are zero.\footnote{The reason why those descent maps vanish is because, in the terminology of Section~\ref{sec:hyperco} below, they are obtained by compressing hyperboxes of size $\dd =(d_1,d_2)$, where at least one of the $d_i$ is zero. See the formulas \eqref{eq:deo1o2} and \eqref{eq:Phi}.} Specifically, the only diagonal map that may be non-trivial is $\Phi^{- L_1 \cup - L_2}_\s$. This is the composition of the inclusion $\I^{- L_1 \cup - L_2}_\s$ with a descent map $D^{- L_1 \cup - L_2}_\s$.

The descent map  $D^{- L_1 \cup - L_2}_\s$ is a chain homotopy between two maps induced by Heegaard moves between the diagrams $r_{-L}(\Hyper^L)$ and $\Hyper^{\emptyset}$ (that is, between the diagrams obtained from $\Hyper^L$ by deleting the $w$ basepoints, resp. deleting the $z$ basepoints). While the maps induced by Heegaard moves usually involve counts of holomorphic triangles, the chain homotopy $D^{- L_1 \cup - L_2}_\s$ involves also counts of holomorphic quadrilaterals. In practice, these counts are harder to get a handle on, but in some cases (such as the example discussed in the next section) we can simplify the complex $\C$ so that we do not have to compute $D^{- L_1 \cup - L_2}_\s$ at all.

\subsection {Surgeries on the Hopf link}
\label {sec:hopf}

\begin{figure}
\begin{center}
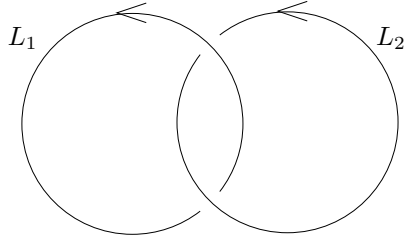
\end{center}
\caption {{\bf The positive Hopf link.}
}
\label{fig:hopf}
\end{figure}

Let $\orL = L_1 \cup L_2 \subset S^3$ be the positive Hopf link shown in Figure~\ref{fig:hopf}. Its link Floer homology was computed in \cite[Section 12]{Links}. Let $p_1, p_2$ be two integers, $p_1 p_2 \neq \pm 1$. Then $(p_1,p_2)$-surgery on $\orL$ produces the lens space $L(p_1p_2-1, p_1)$, which is a rational homology sphere admitting $p_1p_2-1$ different $\spc$ structures. For any such $\spc$ structure $\ss$, the corresponding Heegaard Floer homology group is $\HFm(S^3_{p_1,p_2}(\orL), \ss) \cong \ff[[U]]$ as a relatively $\Z$-graded module, see \cite[Proposition 3.1]{HolDiskTwo}. In this section we show how Theorem~\ref{thm:FirstSurgery} can be used to recover this calculation. We assume for simplicity that $p_1, p_2 \geq 2$.

\begin{figure}
\begin{center}
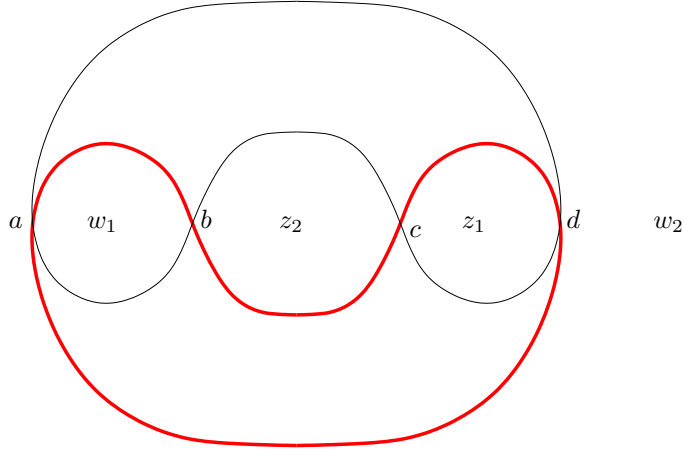
\end{center}
\caption {{\bf A Heegaard diagram for the Hopf link.}
The thicker (red) curve is $\alpha$, while the thinner (black) curve is $\beta$.}
\label{fig:hopf1}
\end{figure}

\begin{figure}
\begin{center}
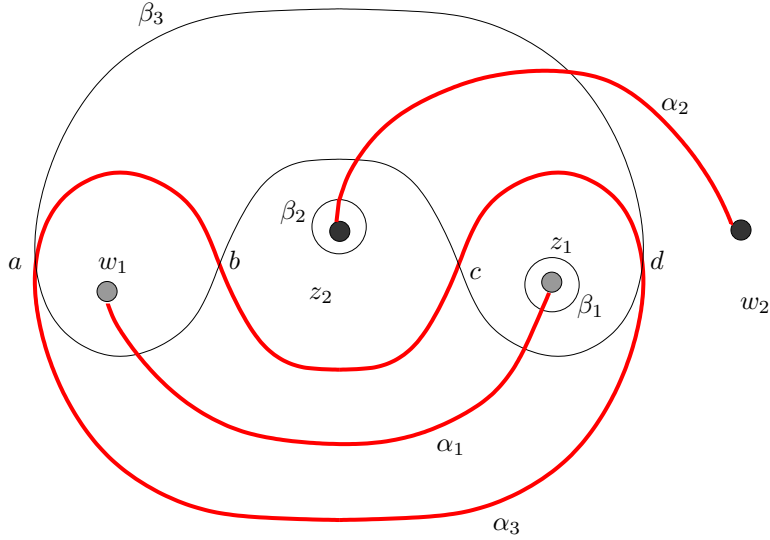
\end{center}
\caption {{\bf Another Heegaard diagram for the Hopf link.}
This picture is obtained from Figure~\ref{fig:hopf1} by stabilizing twice and doing some handleslides. It has the advantage that the basepoints come in pairs $(w_1, z_1)$ and $(w_2, z_2)$, with $w_i$ and $z_i$ on each side of the curve $\beta_i$.
}
\label{fig:hopf2}
\end{figure}

We use a basic system of hyperboxes for the Hopf link $\orL$. To construct it, we start with the genus zero Heegaard diagram for $\orL$ from Figure~\ref{fig:hopf1}. We then stabilize it (and do some handleslides) to end up with the basic Heegaard diagram  pictured in Figure~\ref{fig:hopf2}, in which $w_i$ and $z_i \ (i=1,2)$ are separated by the new beta curves $\beta_i$. The old $\alpha$ and $\beta$ curves are denoted $\alpha_3$ and $\beta_3$, respectively.

In the genus zero Heegaard diagram in Figure~\ref{fig:hopf1}, the intersection $\alpha \cap \beta$ consists of four points, denoted $a, b, c$ and $d$. The index one holomorphic disks correspond to bigons. There are twelve such bigons: one from $b$ to $a$ containing $w_1$, one from $b$ to $c$ containing $z_2$, one from $d$ to $c$ containing $z_1$, on from $d$ to $a$ containing $w_2$; two other bigons go from $a$ to $b$ and contain $z_1$, but they cancel each other in the Floer complex, so for all our purposes they can be ignored. Similarly, there are two bigons from $c$ to $d$ containing $w_1$, two bigons from $a$ to $d$ containing $z_2$, and two bigons from $c$ to $b$ containing $w_2$; all of these can be ignored.

In the genus two $\alpha$-$\beta$ Heegaard diagram from Figure~\ref{fig:hopf2} (which is the diagram denoted $\Hyper^L$ as part of the basic system), the tori $\Ta = \alpha_1 \times \alpha_2 \times \alpha_3$ and $\Tb = \beta_1 \times \beta_2 \times \beta_3$ again intersect  each other in four points. Indeed, since $\beta_1$ intersects a single $\alpha$ curve, namely $\alpha_1$, and that intersection consists of a single point, that point must contribute to any intersection in $\Ta \cap \Tb$. Similar remarks apply to $\alpha_1 \cap \beta_2$. Therefore, the intersection points in $\Ta \cap \Tb$ are determined by their component from $\alpha_3 \cap \beta_3$. We denote them still by $a, b,c$ and $d$, using the obvious correspondence with the intersections in Figure~\ref{fig:hopf1}.

The index one holomorphic disks between $\Ta$ and $\Tb$ in Figure~\ref{fig:hopf1} are also in one-to-one correspondence with those in Figure~\ref{fig:hopf2}. Indeed, each bigon from Figure~\ref{fig:hopf1} corresponds to an annular domain in Figure~\ref{fig:hopf2}. These annular domains are of the same kind as those considered in \cite[proof of Lemma 3.4]{HolDiskTwo}, where it is proved that they support exactly one holomorphic representative (modulo $2$). Therefore, when building $\alpha$-$\beta$ Floer chain complexes from Figure~\ref{fig:hopf2}, we may just as well look at the simpler Figure~\ref{fig:hopf1} and count the corresponding bigons.

The chain complexes we build from $\Hyper^L$ are $\Chain^-(\Hyper^L, \s)$, for $\s \in \H(L)$. We use the alternate definition for these complexes given in Section~\ref{sec:alternative}. Note that:
 
 $$\H(L) = \bigl(\frac{1}{2} + \Z\bigr ) \times \bigl (\frac{1}{2} + \Z \bigr).$$

The Alexander gradings $(A_1, A_2)$ of $a, b, c, d$ are $(\half, \half), (-\half, \half), (-\half, -\half)$ and $(\half, - \half)$, respectively. Therefore, the formulas for the exponents $E^i_{s_i}(\phi)$ that appear in the definitions of $\Chain^-(\Hyper^L, \s)$ depend only on the signs of $s_i$. 

More precisely, for each $\s = (s_1, s_2) \in \H(L)$, the complex $\Chain^-(\Hyper^L, \s)$ is freely generated over $\ff[[U_1, U_2]]$ by $a, b, c$ and $d$. When $s_1, s_2 > 0$, the differential $\del$ on the respective complex counts powers of $U$ according to the multiplicities of $w_1, w_2$ and ignores $z_1, z_2$. We get the following complex, denoted $\Chain^{++}:$
$$ \Chain^{++} :  \  \  \del a = \del c = 0, \ \del b = U_1a + c, \  \del d = U_2 a + c. $$

When $s_1 > 0, s_2 < 0$, we use $w_1$ and $z_2$, and ignore $z_1, w_2$. We get the complex
$$   \Chain^{+-} :   \  \ \del a = \del c = 0, \ \del b = U_1a + U_2c, \  \del d = a + c. $$

When $s_1 < 0, s_2 > 0$, we use $z_1$ and $w_2$ and obtain
 $$  \Chain^{-+} : \  \ \del a = \del c = 0, \ \del b = a + c, \  \del d = U_2 a + U_1c. $$

Finally, when $s_1 < 0, s_2 < 0$, we use $z_1$ and $z_2$ and obtain
$$   \Chain^{--} :  \  \ \del a = \del c = 0, \ \del b =a + U_2c, \  \del d = a + U_1c. $$

We now turn our attention to the other three Heegaard diagrams in the basic system $\Hyper$, namely $\Hyper^{L_1}, \Hyper^{L_2}$ and $\Hyper^{\emptyset}$. Note that
$$ \H(L_1) = \H(L_2) = \Z, \ \ \ \H(\emptyset) = 0$$
and we have
$$ \psi^{\pm L_1} : \H(L) \to \H(L_2), \ \ \psi^{\pm L_1}(s_1, s_2) = s_2 \mp \half,$$
$$ \psi^{\pm L_2} : \H(L) \to \H(L_1), \ \ \psi^{\pm L_1}(s_1, s_2) = s_1 \mp  \half.$$

The diagram $\Hyper^{L_1}$ is obtained from $\Hyper^{L}$ by deleting $z_2$. Let us study $\Chain^-(\Hyper^{L_1}, s)$ for $s \in \Z$. The four generators $a,b,c,d$ have Alexander gradings $0,-1,-1,0$, respectively. Thus, the complex $\Chain^-(\Hyper^{L_1}, s)$ is isomorphic to $\Chain^{++}$ for $s \geq 0$ and to $\Chain^{-+}$ for $s < 0$.

Similarly, $\Hyper^{L_2}$ is obtained from $\Hyper^{L}$ by deleting $z_1$. The complex  $\Chain^-(\Hyper^{L_2}, s)$ is isomorphic to $\Chain^{++}$ for $s \geq 0$ and to $\Chain^{+-}$ for $s < 0$.

Lastly, $\Hyper^{\emptyset}$ is obtained from $\Hyper^{L}$ by deleting both $z_1$ and $z_2$, and the corresponding complex $\Chain^-(\Hyper^{\emptyset}, 0)$ is a copy of $\Chain^{++}$.

The surgery coefficients $p_1, p_2$ on the two components of $L$ describe a framing $\Lambda$ of $L$. Let $(1,0), (0,1)$ be the generators of $H_1(S^3 - L) \cong \Z^2$ corresponding to the meridians of $L_1$ and $L_2$, respectively. Since the linking number between $L_1$ and $L_2$ is $1$, the framings of the components are $\Lambda_1 = (p_1, 1)$ and $\Lambda_2 = (1,p_2 )$.

Let us now describe the full complex $\C^-(\Hyper, \Lambda)$, whose homology is presumed to produce 
$$\HFm(S^3_{\Lambda}(\orL)) \cong \ff[[U]]^{\oplus (p_1p_2-1)}.$$ 

As an $\Ring$-module, $\C = \C^-(\Hyper, \Lambda)$ is the direct product of complexes $\C_\s$ over $\s=(s_1, s_2) \in \H(L) = (\Z+ \half)^2$, where 
$$ \C_\s = \Chain^-(\Hyper^{L}, \s)  \oplus \Chain^-(\Hyper^{L_1}, s_1 - 1/2) \oplus \Chain^-(\Hyper^{L_2}, s_2 - 1/2) \oplus \Chain^-(\Hyper^{\emptyset}, 0).$$

Following \eqref{eq:ces}, we denote the four terms in the direct sum above by $\C^{00}_\s, \C^{01}_\s, \C^{10}_\s, \C^{11}_\s$, in this order; for simplicity, we write the superscript as $00$ rather than $(0,0)$, etc. Recall that each of the terms is freely generated over $\ff[[U_1, U_2]]$ by four generators $a, b, c, d$. In order to be able to tell these generators apart, we denote them by $a^{\eps_1\eps_2}_\s, b^{\eps_1\eps_2}_\s$, etc. when they live in $\C^{\eps_1 \eps_2}_\s$.

\begin{figure}
\begin{center}
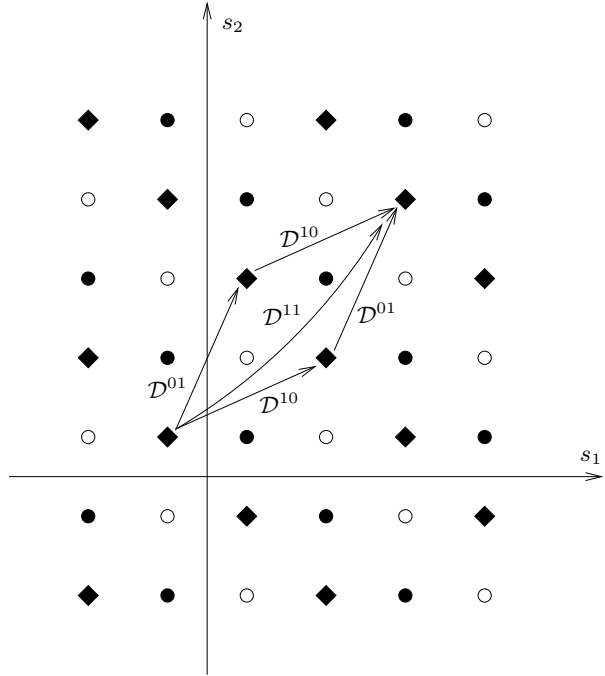
\end{center}
\caption {{\bf The complex $\C^-(\Hyper, \Lambda)$ for $p_1 = p_2 = 2$.}
We show here the lattice $\H(L)$, as a union of various icons in the plane: black dots, white dots, and black diamonds. Each type of icon corresponds to a particular $\spc$ structure on the surgered manifold.
We also show how various parts of the differential $\D^-$ act on the lattice. Not shown is $\D^{00}$, which simply preserves each icon.  Note that all parts of $\D^-$ preserve the type of the icon.
}
\label{fig:D-}
\end{figure}

The differential $\D^-: \C^-(\Hyper, \Lambda) \to \C^-(\Hyper, \Lambda)$ also splits as a sum of four terms
$$ \D^- = \D^{00} + \D^{01} + \D^{10} + \D^{11},$$
where $\D^{\eps'_1\eps'_2}$ maps $\C_\s$ to $\C_{\s + \eps'_1 \Lambda_1 + \eps'_2 \Lambda_2}$. We have chosen here to drop $\eps^0, \eps$ and $\s$ from \eqref{eq:dees}. A graphical representation of the differentials $\D^{\eps'_1\eps'_2}$ is given in Figure~\ref{fig:D-}. Note that the equivalence relation on $\H(L)$ generated by 
$$ (s_1, s_2) \sim (s_1 + p_1, s_2 + 1), \ \ (s_1, s_2) \sim (s_1 + 1, s_2 + p_2)$$ 
breaks $\H(L)$ into $p_1p_2 - 1$ equivalence classes, corresponding to the $\spc$ structures on the  surgered manifold $S^3_{\Lambda}(\orL)$.

More precisely, $\D^{00}$ consists of the differentials $\del = \Phi^{\emptyset}_\s$ on the chain complexes $\C_\s^{\eps_1 \eps_2}$ themselves, plus the cross-terms
$$
\Phi^{L_1}_\s : \C_\s^{00} \to \C_\s^{10} , \ \ \ \ \ \  
 \Phi^{L_2}_\s : \C_\s^{00} \to \C_\s^{01}, $$
  
 $$\Phi^{L_1}_{\psi^{L_2}(\s)} : \C_\s^{01} \to \C_\s^{11}, \ \ \ \ \ \  \Phi^{L_2}_{\psi^{L_1}(\s)} : \C_\s^{10} \to \C_\s^{11},
$$
 and
$$ \Phi^L_\s : \C_\s^{00} \to \C_\s^{11}. $$
Note that $\Phi^L_\s$ is a chain homotopy between $\Phi^{L_1}_{\psi^{L_2}(\s)} \circ \Phi^{L_2}_\s $ and $\Phi^{L_2}_{\psi^{L_1}(\s)} \circ \Phi^{L_1}_\s$. In fact, as noted in Section~\ref{sec:twoC}, the map $\Phi^L_\s$ vanishes in the case of a basic system.

The term $\D^{10}$ is simpler. It consists of the maps
$$ \Phi^{-L_1}_{(s_1, s_2)} : \C^{00}_{(s_1, s_2)} \to \C^{10}_{(s_1 + p_1, s_2+1)},$$
$$ \Phi^{-L_1}_{\psi^{L_2}(s_1, s_2)} : \C^{01}_{(s_1, s_2)} \to \C^{11}_{(s_1 + p_1, s_2+1)},$$
and
$$\Phi^{(-L_1)\cup L_2}_{(s_1, s_2)} : \C^{00}_{(s_1, s_2)} \to \C^{11}_{(s_1 + p_1, s_2+1)}.$$
The chain homotopy $\Phi^{(-L_1)\cup L_2}_{(s_1, s_2)}$ is again zero.

A similar description applies to $\D^{01}$. Finally, the term $\D^{11}$ is the simplest of all, consisting only of maps of the type
$$ \Phi^{-L}_{(s_1, s_2)} : \C^{00}_{(s_1, s_2)} \to \C^{11}_{(s_1 + p_1 + 1, s_2 + p_2 +1)}. $$

\begin{figure}
\begin{center}
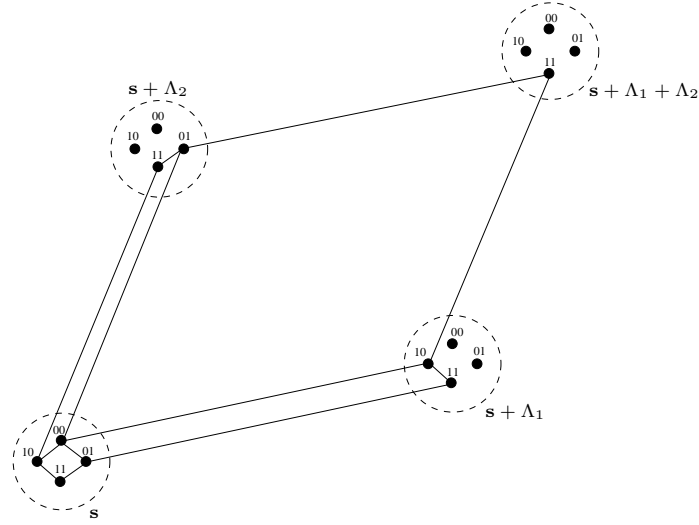
\end{center}
\caption {{\bf Four locally-defined filtrations.} This is a more in-depth look at Figure~\ref{fig:D-}. Each dashed circle corresponds to one of the icons from Figure~\ref{fig:D-}, and the four bullets inside a dashed circle are the four summands $\C^{\eps}_\s$. The values of $\eps$ are noted near each bullet, and the values of $\s$ near each dashed circle. Each of the four parallelograms represents a summand (namely, the one containing $\C^{00}_\s$) in the associated graded complex of one of the four filtrations $\F_{00}, \F_{10}, \F_{01}$ and $\F_{11}$.
}
\label{fig:filtered}
\end{figure}

Let us introduce four locally-defined filtrations 
$$ \F_{\omega_1 \omega_2}, \ \omega_1, \omega_2 \in \{0,1\}$$
on the complex $\C^-(\Hyper, \Lambda)$. These filtrations will play an important role in calculating the homology. 
The differential on the associated graded of $ \F_{\omega_1 \omega_2}$ will be denoted $\tilde \D^{\omega_1 \omega_2}$. These associated graded complexes are shown schematically in Figure~\ref{fig:filtered}.

Recall that the generators of $\C^-(\Hyper, \Lambda)$ (as a direct product) are of the form
$$ g^{\eps_1, \eps_2}_{(s_1, s_2)}, \ \ \ g \in \{a, b, c, d\},  \  \eps_1, \eps_2 \in \{0, 1\}, \ s_1, s_2 \in \Z + \half. $$ 
 
The first filtration $\F_{00}$ is defined on generators by
$$\F_{00}( g^{\eps_1, \eps_2}_{(s_1, s_2)}) = \min\{-s_1 , -s_2\}.$$

The differential $\D^-$ either preserves or decreases the filtration level. In the associated graded, the only visible part of $\D^-$ is $\tilde \D^{00}=\D^{00}$. Consequently, the associated graded splits as a direct product of terms of the form $(\C_\s, \D^{00})$. 

The next filtration $\F_{10}$ is defined on generators by
$$ \F_{10}( g^{\eps_1, \eps_2}_{(s_1, s_2)}) = s_1 - (p_1-1)s_2 - \eps_1.$$

Again, the differential $\D^-$ either preserves or decreases the filtration level. In the associated graded, the only visible part of $\D^-$ is $\tilde \D^{10}$, which is the sum of $\D^{10}$ and the parts of $\D^{00}$ that preserve $\eps_1$. Consequently, the associated graded splits as a direct product of terms of the form:
\begin {equation}
 \label {eq:dtilde10}
\xymatrix{
 \C^{00}_{(s_1, s_2)}\ar[d] \ar[r] \ar[dr] & \C^{10}_{(s_1 + p_1, s_2+1)} \ar[d] \\
 \C^{01}_{(s_1, s_2)}\ar[r]  & \C^{11}_{(s_1 + p_1, s_2 + 1)}
 }
\end {equation}

We similarly have a filtration $\F_{01}$ given by
$$  \F_{01}( g^{\eps_1, \eps_2}_{(s_1, s_2)}) = s_2 - (p_2-1)s_1 - \eps_2.$$

The differential $\tilde \D^{01}$ on its associated graded consists of $\D^{01}$ and the parts of $\D^{00}$ that preserve $\eps_2$. There is a direct product splitting of the associated graded analogous to \eqref{eq:dtilde10}:
\begin {equation}
 \label {eq:dtilde01}
\xymatrix{
 \C^{00}_{(s_1, s_2)}\ar[d] \ar[r] \ar[dr] & \C^{10}_{(s_1 , s_2)} \ar[d] \\
 \C^{01}_{(s_1+1, s_2+p_2)}\ar[r]  & \C^{11}_{(s_1 + 1, s_2 + p_2)}
 }
\end {equation}

The last filtration $\F_{11}$ is defined on generators by
$$\F_{11}( g^{\eps_1, \eps_2}_{(s_1, s_2)}) = \min\{ s_1 - p_1 \eps_1 - \eps_2, s_2 - \eps_1 - p_2 \eps_2\}.$$

The corresponding differential $\tilde \D^{11}$ on the associated graded is the sum of the following terms: the differentials $\del = \Phi^{\emptyset}_\s$ on the chain complexes $\C_\s^{\eps_1 \eps_2}$ themselves, plus the cross-terms
 $$ \Phi^{-L_1}_{(s_1, s_2)} : \C^{00}_{(s_1, s_2)} \to \C^{10}_{(s_1 + p_1, s_2+1)}, \ \ \
  \Phi^{-L_1}_{\psi^{L_2}(s_1, s_2)} : \C^{01}_{(s_1, s_2)} \to \C^{11}_{(s_1 + p_1, s_2+1)},$$
 $$ \Phi^{-L_2}_{(s_1, s_2)} : \C^{00}_{(s_1, s_2)} \to \C^{01}_{(s_1 + 1, s_2+p_2)}, \ \ \
 \Phi^{-L_2}_{\psi^{L_1}(s_1, s_2)} : \C^{10}_{(s_1, s_2)} \to \C^{11}_{(s_1 + 1, s_2+p_2)},$$
and
$$ \D^{11} = \Phi^{-L}_{(s_1, s_2)} : \C^{00}_{(s_1, s_2)} \to \C^{11}_{(s_1 + p_1 + 1, s_2 + p_2 +1)}.$$

Consequently, the associated graded of $\F_{11}$ splits as a direct product of terms of the form: 
\begin {equation}
 \label {eq:dtilde11}
\xymatrix{
 \C^{00}_{(s_1, s_2)}\ar[d] \ar[r] \ar[dr] & \C^{10}_{(s_1 + p_1, s_2+1)} \ar[d] \\
 \C^{01}_{(s_1 + 1, s_2+p_2)}\ar[r]  & \C^{11}_{(s_1 + p_1 + 1, s_2 + p_2 + 1)}
 }
\end {equation}

Let us now turn to the computation of the homology of  $\C^-(\Hyper, \Lambda)$. This complex has a subcomplex 
$$ \C_{\geq 0} = \prod_{\{(s_1, s_2)| \max(s_1, s_2) > 0\}} \C_{(s_1, s_2)}.$$ 

\begin {lemma}
\label {lemma:>0}
$H_*(\C_{\geq 0}) = 0$. 
\end {lemma}

\begin {proof}
We use the restriction of the filtration $\F_{00}$ to $\C_{\geq 0}$. Note that this restriction is bounded above.  In light of Lemma~\ref{lem:acyclicAG}, it suffices to show that the homology of the associated graded groups $H_*(\C_{(s_1, s_2)}, \D^{00})$ vanishes, whenever $s_1 > 0$ or $s_2 > 0$. 

Let us first consider the case $s_1, s_2 > 0$. Then $\C_{(s_1, s_2)}$ consists of four copies of the same complex, $\Chain^{++}$, related by the five cross-terms in the description of $\D^{00}$ above. Among these, $ \Phi^{L_1}_\s,  \Phi^{L_2}_\s , \Phi^{L_1}_{\psi^{L_2}(\s)}$ and $\Phi^{L_2}_{\psi^{L_1}(\s)}$ are all isomorphisms, corresponding to the identity on the complex $\Chain^{++}$, while the chain homotopy $\Phi^L_\s$ is trivial.  Thus, the complex $\C_\s$ can be described as
$$\begin {CD}
\Chain^{++} @>{\cong}>> \Chain^{++} \\
@V{\cong}VV @VV{\cong}V \\
\Chain^{++} @>{\cong}>> \Chain^{++}
\end {CD} $$
This is clearly acyclic.

Next, let us consider the case $s_1 > 0, s_2 < 0$. Then $\C_{(s_1, s_2)}$ consists of two copies of $\Chain^{+-}$, namely $\C^{00}_{(s_1, s_2)}$ and $\C^{10}_{(s_1, s_2)}$, and two copies of $\Chain^{++}$, namely $\C^{01}_{(s_1, s_2)}$ and $\C^{11}_{(s_1, s_2)}$. Thus $\C_\s$ can be described as
$$\begin {CD}
\Chain^{+-} @>{\cong}>> \Chain^{+-} \\
@VVV @VVV \\
\Chain^{++} @>{\cong}>> \Chain^{++}
\end {CD} $$

Even though the vertical maps are not isomorphisms, the horizontal ones are identities. This suffices to show that $\C_\s$ is acyclic. Indeed, the filtration $\F_{01}$ restricts to a filtration on $\C_\s$, whose associated graded differential consists of the two horizontal maps above. Therefore, this associated graded is acyclic, and so is $\C_\s$.

Similar remarks apply to the case $s_1 < 0, s_2 > 0$. This shows that $\C_{\geq 0}$ is acyclic.
\end {proof}

Lemma~\ref{lemma:>0} implies that the homology of $\C^-(H, \Lambda)$ is the same as that of its quotient complex
$$ \C_{< 0} = \prod_{s_1, s_2 < 0} \C_{(s_1, s_2)}.$$ 

Next, we show that a large part of the complex $\C_{< 0}$ is also acyclic. Pick small $\zeta_1, \zeta_2 > 0$ linearly independent over $\qq$. Consider the parallelogram $P_\rr$ in the plane with vertices 
$$ (-\zeta_1, -\zeta_2), \ (-\zeta_1 - p_1, -\zeta_2 - 1), \ (-\zeta_1 - 1, -\zeta_2 - p_2), \ (-\zeta_1 - p_1 -1, -\zeta_2 - p_2 - 1). $$
 
If $\zeta_1, \zeta_2$ are sufficiently small, the parallelogram $P_\rr$ contains a unique representative from each equivalence class in $\H(L)$, i.e. from each $\spc$ structure on the surgered manifold. Set $P = P_\rr \cap \H(L)$.
 
For $\omega_1, \omega_2 \in \{0,1\}$, let $Q_{\omega_1 \omega_2}$ be the quadrant in the plane given by
$$  Q_{\omega_1 \omega_2} = \{(s_1, s_2) \in \H(L) | s_1 < -\zeta_1 +(\omega_1-1) p_1 + (\omega_2-1), \ s_2 < -\zeta_2 +(\omega_1 -1) +(\omega_2-1)p_2\}.$$ 

\begin{figure}
\begin{center}
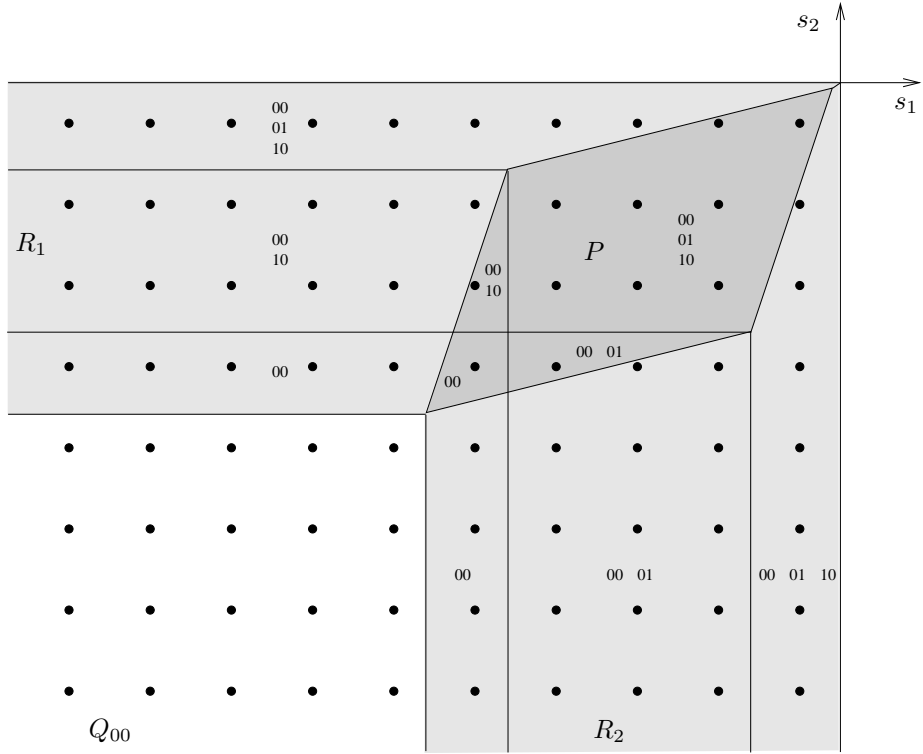
\end{center}
\caption {{\bf The complex $\C_{<0}$ for $p_1= 4, p_2 =3$.}
The parallelogram $P$ is darkly shaded, and the two regions $R_1$ and $R_2$ are more lightly shaded. Each dot represents an element of $Q_{11} \subset \H(L)$. The boundaries of $Q_{00}, Q_{01}, Q_{10}, Q_{11}$ and $P$ split the lower left quadrant in eleven regions. In each of these regions we mark the values $\eps_1\eps_2$ for which the respective groups $\C^{\eps_1 \eps_2}_\s$ are part of the complex $\C_R=\C_{<0}/\C_{<R}$.
}
\label{fig:parallelo}
\end{figure}

The complement $R=Q_{11} \setminus Q_{00}$ consists of the union of $P$ and two other regions: one to the left of $P$, which we denote by $R_1$, and one below $P$, which we denote by $R_2$. These regions and the parallelogram $P$ are shown in Figure~\ref{fig:parallelo}. 

Consider the following submodule of $\C_{<0}:$
$$ \C_{<R} = \prod_{\s \in Q_{00}} \bigl(\C_\s^{00}  \oplus \C_{\s+\Lambda_1}^{10} \oplus \C_{\s+\Lambda_2}^{01} \oplus \C_{\s+\Lambda_1 + \Lambda_2}^{11}    \bigr) .$$

It is straightforward to check that $\C_{<R}$ is a subcomplex of $\C_{<0}$. We denote the corresponding quotient complex by $\C_R$.

\begin {lemma}
\label {lemma:r}
$H_*(\C_{<R}) = 0$.
\end {lemma}

\begin {proof}
We use the restriction of the filtration $\F_{11}$ to the subcomplex $\C_{<R}$. Since this restriction is bounded above, by Lemma~\ref{lem:acyclicAG} it suffices to prove that the associated graded groups are acyclic. 

Recall that the differential of the associated graded of $\F_{11}$ is denoted $\tilde \D^{11}$; see the description of the terms of $\tilde \D^{11}$ before the diagram ~\eqref{eq:dtilde11}. As discussed in Section~\ref{sec:twoC}, the maps $\Phi$ appearing in $\tilde \D^{11}$ are compositions of descent maps $\De$ with inclusion maps $\Pr$. Since $s_1, s_2 \leq -1/2$ and the generators all have Alexander gradings at least $-1/2$, from the description 
\eqref{eq:proj} of the inclusion maps we see that in our case the inclusion maps are the identity. 

Turning our attention to the descent maps, one could in principle compute them explicitly by counting holomorphic polygons. However, it is not necessary to do so. All we need to know is that the descent maps at one component only (be it $-L_1$ or $-L_2$), given by counting holomorphic triangles, induce isomorphisms on homology. This is true because they correspond to the natural triangle maps between strongly equivalent Heegaard diagrams.

The associated graded of $\C_{<R}$ with respect to $\F_{11}$ breaks down into a direct product of factors of the form ~\eqref{eq:dtilde11}, one for each $\s \in Q_{00}$. The fact that the vertical and horizontal maps in \eqref{eq:dtilde11} are quasi-isomorphisms implies that the respective factors are acyclic. The claim follows.
\end {proof}

We denote the quotient complex of $\C_{<R} \subset \C_{<0}$ by $C_R$. Lemma~\ref{lemma:r} implies that
$$ H_*(\C_{< 0}) = H_*(\C_R).$$

Let $\C_P \subset \C_R$ be the submodule
$$ \C_P = \prod_{\s \in P} \C_\s^{00}.$$

This is a quotient complex of $\C_R$. The respective subcomplex $\C_{R_1 \cup R_2}$ splits as a direct sum of two complexes
$$ \C_{R_1} =\Bigl( \prod_{\s \in R_1}  \C_\s \oplus \prod_{\s \in P} \C_\s^{10} \Bigr) \cap \C_R, $$
$$  \C_{R_2} = \Bigl( \prod_{\s \in R_2} \C_\s  \oplus \prod_{\s \in P} \C_\s^{01} \Bigr) \cap \C_R.$$

\begin {lemma}
\label {lemma:cr12}
$ H_*(\C_{R_1} ) = H_*(\C_{R_2} ) = 0$.
\end {lemma}

\begin {proof}
Consider the restriction of the filtration $\F_{10}$ to the complex $\C_{R_1}$. This restriction is bounded above. The respective associated graded splits as a direct product of mapping cone complexes of two possible kinds:
$$ \Chain^{--} \cong \C_\s^{00}  \xrightarrow{\Phi_\s^{-L_1}} \C_{\s+\Lambda_1}^{10} \cong \Chain^{+-},$$
for $\s \in Q_{01} \cap R_1$, and
$$ \Chain^{--} \cong \C_\s^{00}  \xrightarrow{\Phi_\s^{L_2}} \C_\s^{01} \cong \Chain^{-+},$$
for $\s = (s_1, -1/2) \in R_1$.  

Mapping cone complexes of the first kind are acyclic by the same reasoning as in Lemma~\ref{lemma:r}, because $\Phi_\s^{-L_1}$ is a quasi-isomorphism.

Let us study a mapping complex of the second kind. To compute $\Phi_\s^{L_2}$, note that it is the composition of a descent map $\De$ and an inclusion map $\Pr$; the former is the identity and the latter is multiplication by suitable powers of $U_2$.  More precisely, in terms of the generators $a, b, c,d $ of $\Chain^{--}, \Chain^{-+}$, the map $\Phi_\s^{L_2}$ is given by
 $$ a \to  U_2a, \ \ b \to U_2b, \ \ c \to c, \ \ d \to d.$$
This induces an isomorphism on homology, which implies that the respective mapping cone complex is acyclic. Hence $H_*(\C_{R_1}) = 0$.

The proof that $\C_{R_2}$ is acyclic is similar, but uses the filtration $\F_{01}$ instead of $\F_{10}$.
\end {proof}

Putting together Lemmas~\ref{lemma:>0}, \ref{lemma:r} and ~\ref{lemma:cr12} we obtain that the homology of the full complex $\C^-(\Hyper, \Lambda)$  is the same as the homology of the complex
$$ \C_P = \prod_{\s \in P} \C_\s^{00}.$$

For $s_1, s_2 < 0$ we have $\C_\s^{00} \cong \Chain^{--}$, whose homology is easily seen to be isomorphic to $\ff[[U_1, U_2]]/(U_1 - U_2) \cong \ff[[U]]$. Since there are $p_1p_2 - 1$ lattice points from $\H(L)$ inside the parallelogram $P$, we obtain
$$ H_*(\C^-(\Hyper, \Lambda)) \cong \ff[[U]]^{\oplus (p_1p_2 - 1)},$$
as expected.

\subsection{Remarks about the general case}
For an $\ell$-component link $\orL \subset Y$, the surgery complex \eqref{eq:CHL} takes the form of an $\ell$-dimensional hypercube, similar to the double complex \eqref{eq:Co2}. At the vertices of the hypercube we have complexes $\Chain^-(\Hyper^{L-M}, \psi^M(\s))$, which are resolutions of the generalized Floer complexes $\Am(\Hyper^{L-M}, \psi^M(\s))$. When the diagrams $\Hyper^{L-M}$ are link-minimal  (that is, have only two basepoints on each link component), the complexes $\Chain^-(\Hyper^{L-M}, \psi^M(\s))$ were constructed in Section~\ref{sec:alternative}, and are isomorphic to $\Am(\Hyper^{L-M}, \psi^M(\s))$. In the general case, the resolutions $\Chain^-(\Hyper^{L-M}, \psi^M(\s))$ will be defined in Section~\ref{sec:general}.

Describing the maps $\Phi^{\orN}_{\psi^M(\s)}$ that are the building blocks of the differential requires additional work, even in the link-minimal case (such as for a basic system). In the link-minimal case, the maps $\Phi^{\orN}_{\psi^M(\s)}$ are the composition of inclusion maps $\I$ (of the type defined in Section~\ref{sec:inclusions}), and certain descent maps $D$. The simplest descent maps (those that appear along the edges of the hypercube) are equivalences induced by Heegaard moves. The descent maps that are needed for the higher-dimensional faces of the hypercube are chain homotopies between the edge maps, and higher homotopies between these. They will be defined by counting holomorphic polygons. Note that for the equivalences induced by Heegaard moves, we need to compose the maps corresponding to a single move. Similarly, for the higher descent maps, we need to compose higher homotopies from hypercubes. The homological algebra needed to describe this kind of compositions will be developed in Section~\ref{sec:hyperco}.

The general link surgery formula is stated not just for a basic system, but for a general complete system of hyperboxes. (We need this more general formulation in the context of grid diagrams, for example; see Section~\ref{sec:grid}.) Complete systems will be defined in Section~\ref{sec:hyperHeegaard}. In a complete system, the diagram $\Hyper^M$ associated to $M \subset L$ is not necessarily the reduction of $\Hyper^L$ at $M$. Rather, for any orientation $\orM$ of $M$, as part of the complete system we specify a sequence of Heegaard moves that relate $r_{\orM}(\Hyper^L)$ to $\Hyper^M$. These Heegaard moves are then used to construct the descent maps. Furthermore, in the case where the diagrams are not link-minimal, to get the descent maps we will compose the maps induced by Heegaard moves with certain other transition maps; we refer to Section~\ref{sec:general} for more details.

\section {Hyperboxes of chain complexes and compression}

\label {sec:hyperco}

In this section we develop some homological algebra that is essential for the statement of the surgery theorem. All the vector spaces we consider are over $\ff = \zz/2\zz$. 

When $f$ is a function, we denote its $n\th$ iterate by $ f^{\circ n}$, i.e. $f^{\circ 0} = id, \ f^{\circ 1} = f, \ f^{\circ (n+1)} = f^{\circ n} \circ f$.  

\subsection {Hyperboxes of chain complexes}
\label {sec:hyperv}
An $n$-dimensional {\em hyperbox} is a subset of $\rr^n$ of the form $[0, d_1] \times \dots \times [0, d_n]$, where $d_i \geq 0, i=1, \dots, n$. We will assume that  $\dd = (d_1, \dots, d_n) \in \N^n$ is a collection of nonnegative integers. We then let $\E(\dd)$ be the set of points in the corresponding hyperbox with integer coordinates, i.e.
$$  \E(\dd) =  \{ \eps = (\eps_1, \dots, \eps_n) \ | \  \eps_i \in \{0,1, \dots, d_i\}, \ i=1, \dots, n \}. $$

In particular, $ \E_n =  \E(1, \dots, 1) = \{0,1\}^n$ is the set of vertices of the $n$-dimensional unit hypercube. 

For $\eps = (\eps_1, \dots, \eps_n)\in \E(\dd)$, we set
$$ \| \eps \| = \eps_1 + \dots + \eps_n.$$

We can view the elements of $\E(\dd)$ as vectors in $\R^n$. There is a partial ordering on $\E(\dd)$, given by $\eps' \leq \eps \iff \forall i, \ \eps'_i \leq \eps_i$. We write $\eps' < \eps$ if $\eps' \leq \eps$ and $\eps' \neq \eps$.

For $i=1, \dots, n$, let $\tau_i \in \E_n$ be the $n$-tuple formed of $n-1$ zeros and a single one, where the one is in position $i$. Then, for any $\eps \in \E(\dd)$ we have
$$ \eps = \eps_1 \tau_1 + \dots + \eps_n \tau_n.$$

\begin {definition}
\label {def:hyperbox}
An {\em $n$-dimensional hyperbox of chain complexes} of size $\dd \in \N^n$ consists of a collection of $\Z$-graded vector spaces 
$$ (C^{\eps})_{\eps \in \E(\dd)}, \ \  C^{\eps} = \bigoplus_{* \in \Z} C^{\eps}_*,$$ 
together with a collection of linear maps
$$ \De^{\eps}_{\eps^0}  : C_*^{\eps^0} \to C_{*-1+\| \eps \| }^{\eps^0 + \eps},$$
one map for each $\eps^0 \in \E(\dd)$ and  $\eps \in \E_n $ such that $\eps^0 + \eps \in \E(\dd)$. The maps are required to satisfy the relations
\begin {equation}
\label {eq:d2}
 \sum_{\eps' \leq \eps}  \De^{\eps - \eps'}_{\eps^0 + \eps'} \circ \De^{\eps'}_{\eps^0}   = 0,
 \end {equation}
for all $\eps^0 \in \E(\dd), \eps \in \E_n$ such that $\eps^0 + \eps \in \E(\dd)$. 
\end {definition}

Given a hyperbox of chain complexes as above, we denote 
$$ C =  \bigoplus_{\eps \in \E(\dd)} C^{\eps} $$
and define linear maps $ \De^{\eps} : C \to C$, by setting them on generators to 
$$\De^{\eps}(x) =  \begin {cases} \De^{\eps}_{\eps^0} (x) & \text{ for } x \in C^{\eps^0} \text { with } \eps^0 + \eps \in \E(\dd), \\
0 & \text{ for } x \in C^{\eps^0} \text { with } \eps^0 + \eps \not \in \E(\dd).
\end {cases}$$

We denote a typical hyperbox of chain complexes by $H = \bigl( (C^{\eps})_{\eps \in \E(\dd)}, (\De^{\eps})_{\eps \in \E_n} \bigr); $\footnote{Observe that the same notation $\eps$ is used for both elements of $\E(\dd)$ and elements of $\E_n$. Of course, $\E_n$ can be viewed a subset of $\E(\dd)$.} the maps $\De^{\eps}_{\eps^0}$ are implicitly taken into account in the direct sums $\De^{\eps}$.  Sometimes, by abuse of notation, we let $\De^{\eps}$ stand for any of its terms   $\De^{\eps}_{\eps^0}$. If $\dd = (1, \dots, 1)$, we say that $H $ is a {\em hypercube of chain complexes}.

Observe that a $0$-dimensional hyperbox of chain complexes is simply a chain complex, while a $1$-dimensional hyperbox with $\dd = (d)$ consists of chain complexes $C^{(i)}, i=0, \dots, d$, together with a string of chain maps
\begin {equation}
\label {eq:string}
 C^{(0)} \xrightarrow{\De^{(1)}} C^{(1)}  \xrightarrow{\De^{(1)}}  \dots  \xrightarrow{\De^{(1)}} C^{(d)}. 
 \end {equation}

To give another example, a $2$-dimensional hypercube is a diagram of complexes and chain maps
$$\begin {CD}
\label {eq:square}
C^{(0,0)} @>{\De^{(1,0)}}>> C^{(1,0)} \\
@V{\De^{(0,1)}}VV @VV{\De^{(0,1)}}V \\
 C^{(0,1)} @>{\De^{(1,0)}}>> C^{(1,1)}
\end {CD}$$
together with a chain homotopy 
\begin {equation}
\label {eq:chh}
\De^{(1,1)} : C^{(0,0)} \to C^{(1,1)}
\end {equation}
between $\De^{(1,0)} \circ \De^{(0,1)}$ and $\De^{(0,1)} \circ \De^{(1,0)}$.

In general, if $ \bigl( (C^{\eps})_{\eps \in \E(\dd)}, (\De^{\eps})_{\eps \in \E_n} \bigr)$ is an $n$-dimensional hyperbox, then $(C^{\eps}, D^{(0,\dots, 0)})$ are chain complexes. Along the edges of the hyperbox we see strings of chain maps  $\De^{\tau_i}, \ i=1, \dots, n$. In fact, let us imagine the hyperbox $[0, d_1] \times \dots \times [0, d_n]$ to be split into $d_1d_2 \dots d_n$ unit hypercubes. Then along each edge of one of these hypercubes we see a chain map. Along the two-dimensional faces of the unit hypercubes we have chain homotopies, and along higher-dimensional faces we have higher homotopies.

Observe that when $H = (C^{\eps}, \De^{\eps})_{\eps \in \E_n}$ is a hypercube, we can also form a total complex $H_\tot = (C_*, \De)$, where the grading on $C$ is given by
$$ C_* =  \bigoplus_{\eps \in \E(\dd)} C^{\eps}_{*+\|\eps\|} $$
and the boundary map is the sum $\De = \sum \De^\eps$. 

\subsection {Compression} 
\label {sec:compression1}

Let $H =  \bigl( (C^{\eps})_{\eps \in \E(\dd)}, (\De^{\eps})_{\eps \in \E_n} \bigr)$ be an $n$-dimensional hyperbox of chain complexes. We will explain how to construct from $H$ an $n$-dimensional hypercube $\hat H = (\hat C^{\eps}, \hat \De^{\eps})_{\eps \in \E_n}$. The process of turning $H$ into $\hat H$ will be called {\em compression}.

The simplest example of compression is when $n=1$, and $H$ is a string of chain complexes and chain maps as in \eqref{eq:string}. Then compression is composing the maps. Precisely, the compressed hypercube $\hat H$ consists of the complexes $\hat C^{(0)} = C^{(0)}$ and $\hat C^{(1)} = C^{(d)}$, linked by the chain map 
$$\bigl( \De^{(1)}\bigr)^{\circ d} = \De^{(1)} \circ \dots \circ \De^{(1)} \ : C^{(0)} \longrightarrow C^{(d)}.$$

For general $n$ and $\dd = (d_1, \dots, d_n)$, the compressed hypercube $\hat H$ has at its vertices the same complexes as those at the vertices of the original hyperbox $H$:
$$ \hat C^{(\eps_1, \dots, \eps_n)} = C^{(\eps_1d_1, \dots, \eps_n d_n)}, \ \eps = (\eps_1, \dots, \eps_n) \in \E_n.$$
    
Further, along each edge of $\hat H$ we should see the composition of the respective edge maps in $H, $ i.e. 
$$ \hat \De^{\tau_i} = \bigl(\De^{\tau_i}\bigr)^{\circ d_i}.$$

The construction of the maps corresponding to the higher-dimensional faces of $\hat H$ is rather involved, and will occupy Sections~\ref{sec:songs} -\ref{sec:compression2}. For now, to give a flavor of the respective formulae, let us explain the simplest nontrivial case, namely $n=2$.

When $n=2$ and $\dd = (d_1, d_2)$, the hyperbox $H$ is a rectangle split into $d_1d_2$ unit squares. Along the horizontal edges we have chain maps denoted $f_1=\De^{(1,0)}$ and along the vertical edges we have chain maps denoted $f_2=\De^{(0,1)}$. Further, each unit square carries a chain homotopy $f_{\{1,2\}}=\De^{(1,1)}$ between $f_1 \circ f_2$ and $f_2 \circ f_1$. Then, on the edges of the compressed hypercube $\hat H$ we have maps
$$ \hat \De^{(1,0)} =  f_1^{\circ d_1}, \ \  \hat \De^{(0,1)} =  f_2^{\circ d_2}.$$
For the diagonal map $\hat  \De^{(1,1)}$ we choose
\begin {equation}
\label {eq:deo1o2}
\hat \De^{(1,1)} = \sum_{j_1=1}^{d_1} \sum_{j_2=1}^{d_2} f_1^{\circ(j_1-1)} \circ f_2^{\circ (j_2-1)} \circ f_{\{1,2\}} \circ f_2^{\circ (d_2 - j_2)} \circ f_1^{\circ (d_1 - j_1)}.
\end {equation}

\begin{figure}
\begin{center}
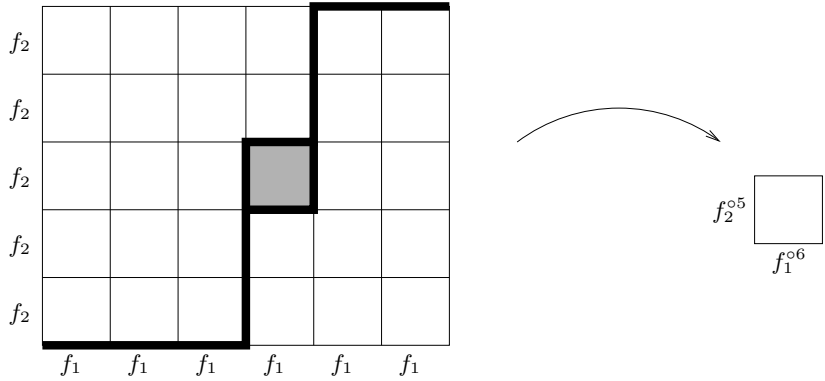
\end{center}
\caption {{\bf Compression of a rectangle into a square.}
This is the graphical representation of the term of the form $f_1^{\circ 2} \circ f_2^{\circ 2} \circ f_{\{1,2\}}  \circ f_2^{\circ 2} \circ f_1^{\circ 3}$ appearing in the sum~\eqref{eq:deo1o2}, with $d_1 = 6, d_2 = 5, j_1 = j_2 =3$. Each unit segment which is part of the thick line corresponds to a map $f_1$ or $f_2$, while the shaded square is the chain homotopy $f_{\{1,2\}}$. Taking the sum of all the terms in \eqref{eq:deo1o2} corresponds to filling up the whole rectangle with $30$ unit squares, and represents a chain homotopy between $f_1^{\circ 6} \circ f_2^{\circ 5}$ and $f_2^{\circ 5} \circ f_1^{\circ 6}$. 
}
\label{fig:compression}
\end{figure}

It is easy to check that $\hat \De^{(1,1)}$ is a chain homotopy between $ \hat \De^{(1,0)} \circ \hat \De^{(0,1)} =f_1^{\circ d_1} \circ f_2^{\circ d_2}$ and $\hat \De^{(0,1)} \circ \hat \De^{(1,0)} =f_2^{\circ d_2} \circ f_1^{\circ d_1}$. See Figure~\ref{fig:compression} for a pictorial interpretation.

\subsection {The algebra of songs}
\label {sec:songs}

Let $X$ be a finite set. (Typically, $X$ will be a subset of the set of nonnegative integers.)

\begin {definition}
An {\em $X$-valued song} is a finite, ordered list of {\em items}, where each item can be either a {\em note}, i.e. an element of $X$, or a {\em harmony}, i.e. a subset of $X$. 
\end {definition}

For example, if $X=\{1,2,3\}$, a typical $X$-valued song is written as
$$ s = (213\{2,3\}2\{\}12\{3\}).$$
The song $s$ has nine items, six being notes and three being
harmonies. Note that we allow the empty harmony $\{ \}$, and that we
distinguish between the note $3$ and the one-element harmony $\{3\}$.
Also, our convention is to write songs between parantheses. For
example, if $A \subseteq X$ is a harmony, by $(A)$ we mean the
one-item song made of that harmony.

Let $\tilde S(X)$ be the $\ff$-vector space freely generated by all
$X$-valued songs. Since songs form a monoid under concatenation (with
the empty song as the identity element), this induces the structure of
a non-commutative algebra on $\tilde S(X)$.

\begin {definition}
\label {def:songs}
The {\em algebra of $X$-valued songs}, denoted $S (X)$, is the quotient of $\tilde S(X)$ by the ideal $I(X)$ generated by the following relations:
\begin {itemize}
\item For any note $x \in X$,
\begin {equation}
\label {eq:unu}
 (x \{ \}) = (\{ \} x)
 \end {equation}
 and
 \begin {equation}
\label {eq:unu'}
  (x\{x\}) = (\{ x \}x).
 \end {equation}
 \item
For any note $x \in X$ and song $s \in \tilde S(X)$,
\begin {equation}
\label {eq:doi}
(x\{x\}sx) + (xs\{x\}x) = (xs) + (sx).
\end {equation}
\item For any harmony $A \subseteq X$,  
\begin {equation}
\label {eq:trei}
 \sum_{B \subseteq A} (B)(A \setminus B) = 0.
 \end {equation}
\end {itemize}
\end {definition}

\begin {remark}
The reader may wonder why we defined the algebra $S(X)$ this way. The motivation behind our choice of the relations ~\eqref{eq:unu}-\eqref{eq:trei}  is that they do not affect the playing of songs, as defined in Section~\ref{sec:plays} below. See Lemma~\ref{lemma:plays} for the relevant result.
\end {remark}

Let $Y = X \cup \{y\}$ be the set obtained from $X$ by adding a new note $y$. We define an operation on songs:
$$ \psi_y : \tilde S(X) \to \tilde S(Y),$$
as follows. For a note $x \in X$, we set
$$\psi_y(x) = (xy\{x,y\}yx).$$
For a harmony $A \subseteq X$, let $\Pi(A)$ be the set of all ordered decompositions $(A_1, \dots, A_k)$ of $A$ into a disjoint union
$$A = A_1 \amalg A_2 \amalg \dots \amalg A_k,$$
where $k \geq 0$ and all $A_i$'s are nonempty. For $A \neq \{ \}$, we set
$$ \psi_y(A) = \sum_{(A_1,\dots, A_k) \in \Pi(A)} (y(A_1\cup \{y\})y(A_2 \cup \{y\})y \dots y(A_k \cup \{y\})y),$$
while for $A = \{ \}$ we set $\psi_y(A) = (y)$.

So far we defined $\psi_y$ only on one-item songs, consisting of either one note or one harmony. We extend it to arbitrary songs by requiring it to act as a derivation:
$$ \psi_y(s_1s_2) = \psi_y(s_1)s_2 + s_1\psi_y(s_2).$$
For example,
$$ \psi_3(2\{1,2\}) = (23\{2,3\}32\{1,2\}) + (23\{1,2,3\}3) + (23\{1,3\}3\{2,3\}3) + (23\{2,3\}3\{1,3\}3).$$

Finally, we extend $\psi_y$ to all of $\tilde S(X)$ by requiring it to be linear.

\begin {lemma}
\label {lemma:sxy}
The operation $\psi_y$ descends to a linear map between $S(X)$ and $S(Y)$.
\end {lemma}

\begin {proof}
We need to check that when we apply $\psi_y$ to the relations \eqref{eq:unu}-\eqref{eq:trei} from Definition~\ref{def:songs}, we obtain relations that hold true in $S(Y)$, i.e., lie in the ideal $I(Y)$.

Let us first look at the relation~\eqref{eq:unu}. The claim is that $\psi_y (x \{ \}) = \psi_y(\{ \} x)$. Indeed, we have
\begin {align} 
\label {eq:newunu}
\psi_y (x \{ \}) + \psi_y(\{ \} x) &= (xy\{x, y\}yx\{\}) + (xy) +  (\{\} xy\{x, y\}yx) + (yx) \\ 
&= (xy(\{x, y\}\{\} + \{ \} \{x, y\}) yx) + (xy) + (yx) \notag \\
&= (xy(\{x\}\{y\} + \{y\} \{x\}) yx) + (xy) + (yx) \notag \\ 
&= (x(y\{x\}\{y\}y + y\{y\} \{x\} y)x) + (xy) + (yx) \notag \\ 
&= (x(y\{x\} + \{x\}y)x) +( xy)+(yx) \notag \\ 
&= (xy\{x\}x) + (x\{x\}yx )+ (xy)+(yx) \notag \\ 
&= 0. \notag
\end {align} 
To get the second equality in~\eqref{eq:newunu} we used \eqref{eq:unu}, namely the fact that $\{\}$ commutes with $x$ and $y$. To get the third equality we applied \eqref{eq:trei} for $A = \{x, y\}$, while to get the fifth and seventh equalities we applied \eqref{eq:doi}.

The similar result for Equation~\eqref{eq:unu'} is simpler. It suffices to apply \eqref{eq:unu'} twice and\eqref{eq:doi} once:
\begin {align}
\psi_y(x\{x\}) + \psi_y (\{ x \}x) &= (xy\{x, y\}yx\{x\} ) + (\{x\}xy\{x, y\}yx) + (xy\{x, y\}y) + (y\{x, y\}yx) 
\\
&=
 (xy\{x, y\}y\{x\}x) + (x\{x\}y\{x, y\}yx) + (xy\{x, y\}y) + (y\{x, y\}yx) \notag \\
&= 0. \notag
\end {align}

Here is the analogous result for Equation~\eqref{eq:doi}:
\begin {align} 
\label {eq:newdoi}
 & \phantom{==} \psi_y (x\{x\}sx) + \psi_y(xs\{x\}x) +\psi_y (xs) + \psi_y(sx) \\ 
&=  (xy\{x,y\}yx\{x\}sx) + (xy\{x,y\}ysx) + (x\{x\}\psi_y(s)x) + (x\{x\}sxy\{x,y\}yx) +\notag \\
& \phantom{==} (xy\{x,y\}yxs\{x\}x) + (xsy\{x,y\}yx) + (x\psi_y(s)\{x\}x) + (xs\{x\}xy\{x,y\}yx) + \notag \\
& \phantom{==} (xy\{x,y\}yxs) + (x\psi_y(s)) + (\psi_y(s)x) + (sxy\{x,y\}yx) \notag \\
&= \bigl( (xy\{x,y\}yx\{x\}sx)+  (xy\{x,y\}yxs\{x\}x) + (xy\{x,y\}ysx)+  (xy\{x,y\}yxs) \bigr)  +\notag \\
& \phantom{==}\bigl((x\{x\}sxy\{x,y\}yx)+  (xs\{x\}xy\{x,y\}yx) +(xsy\{x,y\}yx)+ (sxy\{x,y\}yx)\bigr)  + \notag \\
& \phantom{==} \bigl(  (x\{x\}\psi_y(s)x)+ (x\psi_y(s)\{x\}x)+(x\psi_y(s)) + (\psi_y(s)x) \bigr)=0. \notag
\end {align}

In the last step, the four terms in each of the large parantheses cancel each other out by applying Equation \eqref{eq:doi}.

Lastly, we prove that Equation~\eqref{eq:trei} holds true after applying $\psi_y$.
Let us introduce the following notational shortcut: if $A$ is a subset of $X$, we denote by $\tilde A = A \cup \{y\} \subseteq Y$.

For any $A \subseteq X$, we have
\begin {align}
\label {eq:bab}
\psi_y \Bigl (\sum_{B \subseteq A} (B)(A \setminus B) \Bigr) &=  \sum_{B \subseteq A} \psi_y(B) \cdot (A \setminus B) + (A \setminus B) \cdot \psi_y(B) \\
&= \sum_{B \subseteq A} \sum_{(B_1, \dots, B_k) \in \Pi(B)} \Bigl( (y\tilde B_1y\dots y\tilde B_k y (A \setminus B)) + ((A \setminus B)  y\tilde B_1y\dots y\tilde B_k y) \Bigr) \notag \\
&= \sum_{(A_1, \dots, A_k) \in \Pi(A)} \Bigl( 
(y\tilde A_1 y \dots y\tilde A_{k-1}y A_k) + 
(A_1 y\tilde A_2 y \dots y\tilde A_k y)\Bigr) + \notag \\
& \phantom{==} \sum_{(A_1, \dots, A_k) \in \Pi(A)} \Bigl( (y \tilde A_1 y \dots y\tilde A_k y \{ \}) +
(\{ \} y \tilde A_1 y \dots y\tilde A_k y)
 \Bigr). \notag
\end {align}

The last equality was obtained by splitting the summation on the second line into terms with $B\neq A$ and $B=A$.

The final expression in \eqref{eq:bab} is a sum of two terms, where each term is a summation over the elements of $\Pi(A)$. Using the fact that $\{ \}$ and $y$ commute, the second summation (the one appearing on the very last line of \eqref{eq:bab}) is seen to equal:
\begin {equation}
\label {eq:secondsum}
 \sum_{(A_1, \dots, A_k) \in \Pi(A)} \sum_{i=1}^k \Bigl( (y\tilde A_1y \dots y \tilde A_{i-1}y) \cdot 
\bigl( (\{ \}\tilde A_i) + (\tilde A_i \{ \}) \bigr) \cdot(y\tilde A_{i+1}y  \dots y\tilde A_k y) \Bigr).
\end {equation}

Applying \eqref{eq:trei} to $\tilde A_i$ we see that the paranthesis in the middle of the summation term in ~\eqref{eq:secondsum} equals
$$  (\{ \}\tilde A_i) + (\tilde A_i \{ \}) = \sum_{\substack{B \subseteq \tilde A_i  \\ B \neq \emptyset, \tilde A_i}} (B)(\tilde A_i \setminus B) = 
(A_i\{y\}) + (\{y\}A_i) + \sum_{\substack{B \subseteq A_i  \\ B \neq \emptyset, A_i}}
\Bigl( (\tilde B)(A_i \setminus B)+ (B)(\tilde A_i \setminus B)  \Bigr).
$$  

Plugging this back into \eqref{eq:secondsum}, we obtain that  \eqref{eq:secondsum} equals
\begin {equation}
\label {eq:sumaa}
  \sum_{(A_1, \dots, A_k) \in \Pi(A)} \sum_{i=1}^k \Bigl( (y\tilde A_1y \dots y \tilde A_{i-1}) \cdot 
\bigl( (y\{y \} A_i y) + (y A_i \{ y\}y) \bigr) \cdot(\tilde A_{i+1}y  \dots y\tilde A_k y) \Bigr) + 
\end {equation}
$$  \sum_{(A_1, \dots, A_k) \in \Pi(A)} \sum_{i=1}^{k-1} (y\tilde A_1y \dots y \tilde A_{i-1} A_i y  \dots y\tilde A_k y)  + \sum_{(A_1, \dots, A_k) \in \Pi(A)} \sum_{i=1}^{k-1}  (y\tilde A_1y \dots y A_{i-1} \tilde A_i y  \dots y\tilde A_k y) .$$

Here, in the last two summations, we changed notation so that $(A_1, \dots, A_{i-1}, B, A_i \setminus B, A_{i+1}, \dots, A_k)$ is renamed $(A_1, \dots, A_k)$.

Now, applying \eqref{eq:doi} to the middle paranthesis in the first summation in \eqref{eq:sumaa}, we can replace $(y\{y \} A_i y) + (y A_i \{ y\}y)$ with $(yA_i) + (A_i y)$. Consequently, most of the terms in that first summation cancel out with terms in the second and third summations in \eqref{eq:sumaa}. The only remaining terms are some corresponding to $i=1$ and $i=k$. More precisely, we get that \eqref{eq:sumaa} equals
$$\sum_{(A_1, \dots, A_k) \in \Pi(A)} \Bigl( 
(y\tilde A_1 y \dots y\tilde A_{k-1}y A_k) + 
(A_1 y\tilde A_2 y \dots y\tilde A_k y)\Bigr). $$ 

This exactly corresponds to the first summation in the final expression in \eqref{eq:bab}. Hence, we obtain
$$\psi_y \Bigl (\sum_{B \subseteq A} (B)(A \setminus B) \Bigr) = 0,$$
as desired.
\end {proof}

\subsection {Symphonies}

\begin {definition}
\label {def:symph}
Let $X$ be a finite, totally ordered set, and $m(X)$ the maximal element in $X$. The {\em symphony $\alpha(X) \in \tilde S(X)$ on the set $X$} is defined, recursively, by 
$$ \alpha(\emptyset) = (\{\}), \ \ \alpha(X) = \psi_{m(X)}\bigl(\alpha(X \setminus \{m(X)\}) \bigr).$$

We call $\alpha_n = \alpha(\{1, 2, \dots, n\})$ the {\em $n\th$ standard symphony}.
\end {definition}

For example, 
\begin {eqnarray*}
\alpha_1 &=& \psi_1(\{ \}) = (1),  \\
\alpha_2 &=& \psi_2(1) = (12\{1,2\}21), \\
\alpha_3 &=& \psi_3(12\{1,2\}21) = (123\{1,2,3\}321) + (123\{2,3\}32\{1,2\}21) + (12\{1,2\}23\{2,3\}321) + \\ & \; & (12\{1,2\}213\{1,3\}31)+  (13\{1,3\}312\{1,2\}21) + (123\{1,3\}3\{2,3\}321) + (123\{2,3\}3\{1,3\}321).
\end {eqnarray*}

Computer experimentation shows that $\alpha_4$ is a linear combination of $97$ different songs, and $\alpha_5$ a linear combination of $2051$ different songs. In general, a song $s$ that appears with nonzero multiplicity in $\alpha_n$ is easily seen to satisfy the following two conditions. Let $s$ consist of $k$ notes and $l$ harmonies, and let $h_1, \dots, h_l$ be the cardinalities of each harmony. Then we have:
\begin {equation}
\label {eq:h1}
k=2n+l-1, \ \ \ \sum_{i=1}^l h_i = n+l-1.
\end {equation}

Of course, not every song that satisfies \eqref{eq:h1} appears in the formula for $\alpha_n$.

\begin {lemma}
\label {lemma:symph}
For any finite, totally ordered set $X$, we have the following relation in $S(X)$:
\begin {equation}
\label {eq:symph}
\sum_{Y \subseteq X} \alpha(Y)\alpha(X \setminus Y) = 0.
\end {equation}
\end {lemma}

\begin {proof}
Induction on the cardinality $n$ of $X$. For $n=0$ the corresponding relation $\{\}\{\} = 0$ is Equation \eqref{eq:trei} for $A = \{ \}$, while for $n=1$ the corresponding relation $(x\{\}) = (\{\}x)$ is Equation~\eqref{eq:unu}. 

Let $X$ be a set of cardinality $n \geq 2$, and denote $m = m(X)$. Suppose that \eqref{eq:symph} is true for all sets of cardinality $< n$, and, in particular, for $X' = X \setminus \{m\}$. Then:
\begin {eqnarray*}
\sum_{Y \subseteq X} \alpha(Y)\alpha(X \setminus Y) &=& \sum_{Y \subseteq X'} \bigl( \alpha(Y)\alpha(X' \cup \{m\}  \setminus Y) +  \alpha(Y \cup \{m \}) \alpha(X' \setminus Y) \bigr) \\
&=& \sum_{Y \subseteq X'} \Bigl( \alpha(Y)\psi_m \bigl(\alpha(X'   \setminus Y)\bigr) +  \bigl( \psi_m(\alpha(Y ) \bigr) \alpha(X' \setminus Y) \Bigr) \\
&=& \psi_m \bigl(\sum_{Y \subseteq X'} \alpha(Y)\alpha(X' \setminus Y) \bigr) \\
&=& \psi_m(0) = 0.
\end {eqnarray*} \end {proof}

\subsection {Hypercubical collections}
\label {sec:plays}

Let $X$ be a finite set and $(A, +, *)$ be a (possibly non-commutative) algebra over $\ff$. Given an element $\A \in A$, we denote by $\A^{*j} = \A * \dots * \A$ its $j\th$ power. In particular, $\A^{*0}=1$ is the unit.

\begin {definition}
An $n$-dimensional {\em hypercubical collection} in the algebra $A$, modeled on $X$,  is a collection $\A$ composed of elements $\A_Z  \in A$, one for each $Z \subseteq X$, which are required to satisfy the relations
\begin {equation}
\label {eq:wzw}
 \sum_{Z' \subseteq Z} \A_{Z'}  * \A_{ Z \setminus Z' } = 0,
 \end {equation}
for any $Z \subseteq X$.
\end {definition}  
  
\begin {example}
\label {ex:cde}
Let $H=(C^\eps, \De^\eps)$ be a hyperbox of chain complexes as in Section~\ref{sec:hyperv}. Choose $X = \{1, 2, \dots, n\}$. For any $Z \subseteq X$, we can define an element $\zeta(Z) = (\zeta(Z)_1, \dots, \zeta(Z)_n) \in \E_n = \{0,1\}^n$, by
\begin {equation}
\label {eq:zeta}
\zeta(Z)_i = \begin {cases} 
1 & \text{if } i \in Z, \\
0 & \text{otherwise.}
\end {cases}
\end {equation}
Then  
$$  \A_Z = \De^{\zeta(Z)}$$
form a hypercubical collection in the algebra $\End (C)$, modeled on $X$.
\end {example}  
 
 \begin {definition}
 \label {def:playing}
Let $X$ be a finite set and let $\dd = (d_x)_{x \in X}$ be a collection of positive integers indexed by $X$.
Let $\A = \{A_ Z \}_{Z \subseteq X}$ be a hypercubical collection in an algebra $A$. Let also $s$ be an $X$-valued song, spelled out as a sequence of items
$$ s = (x^1_1\dots x^1_{r_1} \{y^1_1, \dots, y^1_{t_1} \}  x^2_1\dots x^2_{r_2} \{y^2_1, \dots, y^2_{t_2} \}  \dots \dots x^l_1\dots x^l_{r_l} \{y^l_1, \dots, y^l_{t_l} \}  x^{l+1}_1\dots x^{l+1}_{r_{l+1}} ), $$
where some of the $r_i$'s can be zero.

The result of {\em playing the song $s$ to the hypercubical collection $\A$, in the register $\dd$,} is the algebra element $\play_\A^{\dd}(s) \in A$ defined by: 
\begin {equation}
\label {eq:big}
 \play_\A^{\dd}(s)= \sum  \A_{\{x^1_1\}}^{* j^1_1} * \dots * \A_{ \{x^1_{r_1}\}}^{*j^1_{r_1}} *  \A_{\{y^1_1, \dots, y^1_{t_1}\} }  *   \A_{\{x^2_1\}}^{* j^2_1} * \dots * \A_{ \{x^2_{r_2}\}}^{*j^2_{r_2}} *  \A_{\{y^2_1, \dots, y^2_{t_2}\} }  * \cdots
 \end {equation}
$$  \cdots *  \A_{\{x^l_1\}}^{* j^l_1} * \dots * \A_{ \{x^l_{r_l}\}}^{*j^l_{r_l}} *  \A_{\{y^l_1, \dots, y^l_{t_l}\} }  *   \A_{\{x^{l+1}_1\}}^{* j^{l+1}_1} * \dots * \A_{ \{x^{l+1}_{r_{l+1}}\}}^{*j^{l+1}_{r_{l+1}}},
$$
where the sum is over all possible powers  $j^\sigma_u \geq 0$ satisfying, for each $x \in X$,
$$ \sum_{\sigma=1}^{l+1} \sum_{u=1}^{r_\sigma} \delta_{x, x^\sigma_u} j^\sigma_u + \sum_{\sigma=1}^l
 \sum_{u=1}^{t_\sigma} \delta_{x, y^\sigma_u} = d_x.$$
Here $\delta_{i, j}$ is the Kronecker delta symbol.
\end {definition} 

It is worth giving some examples of playing. First, note that, unless every $x \in X$ appears at least once in the song $s$ (either as a note or as part of a harmony), we have $ \play_\A^{\dd}(s)=0$. 

If $X = \{1\}$ and $\dd=(d)$, then $\A$ consists of two algebra elements $\A_{\emptyset}$ and $\A_{\{1\}}$. Playing the song $(1)$ to $\A$ yields the $d\th$ power $\A_{\{1\}}^{*d}$. On the other hand, playing the song $(\{1\})$ to $\A$ yields $\A_{\{1\}}$ when $d=1$ and $0$ otherwise. 

If $X = \{1,2\}$ and $\dd=(d_1, d_2)$, then $\A$ consists of four algebra elements $\A_{\emptyset}, \A_{\{1\}}, 
\A_{\{2\}}$ and $\A_{\{1,2\}}$. Playing the second standard symphony $\alpha_2 = (12\{1, 2\}21)$ to $\A$ yields
$$ \play_{\A}^{(d_1, d_2)}(\alpha_2) =\sum_{\substack{j_1^1 + j^2_2 = d_1 -1 \\ j_2^1 + j_1^2 = d_2 -1}} \A_{\{1\}}^{*j_1^1} * \A_{\{2\}}^{* j_2^1} * \A_{\{1,2\}} * \A_{\{2\}}^{* j^2_1} * \A_{\{1\}}^{* j_2^2}, $$
or, equivalently,
$$ \play_{\A}^{(d_1, d_2)}(\alpha_2) =\sum_{j_1=1}^{d_1} \sum_{j_2=1}^{d_2} \A_{\{1\}}^{*(j_1-1)} * \A_{\{2\}}^{* (j_2-1)} * \A_{\{1,2\}} * \A_{\{2\}}^{* (d_2 - j_2)} * \A_{\{1\}}^{* (d_1 - j_1)}.  $$

(Compare Equation~\eqref{eq:deo1o2} and Figure~\ref{fig:compression}.) 

\medskip

\begin {lemma}
\label {lemma:plays}
Let $\A$ be a hypercubical collection in an algebra $A$, modeled on a set $X$, and pick a series of nonnegative integers $\dd = (d_x)_{x \in X}$. Then, the operation of playing songs to $\A$ in the register $\dd$ descends to a linear map $\play^\dd_\A : S(X) \to A$.
\end {lemma}

\begin {proof}
We extend the playing of songs linearly to a map $\tilde S(X) \to A$. In order to show that it descends to $ S(X)$, we need to check that the relations ~\eqref{eq:unu}-\eqref{eq:trei} hold true after playing them. In fact, since playing is not multiplicative, one needs to check that these relations, when multiplied on the left and right with arbitrary songs, still hold true after playing.

For example, let us look at the relation ~\eqref{eq:unu}, namely $(x\{\}) = (\{\}x)$, for any $x \in X$. The claim is that the maps $\play_{(s_1x\{\}s_2)}$ and $ \play_{(s_1\{\}xs_2)}$ are equal, for any songs $s_1, s_2$. This is true, because $x \{ \}$ corresponds to taking a power of $\A_{\{x\}}$ in the big summation in \eqref{eq:big}, followed by the factor $\A_\emptyset$. Equation \eqref{eq:wzw} for $Z = \{x\}$ imples that $\A_\emptyset$ commutes with $\A_{\{x\}}$, so it also commutes with its power.

The relation \eqref{eq:unu'} holds true after playing because both $(x\{x\})$ and $(\{x \}x)$ correspond to taking an arbitrary, but nontrivial, power of $\A_{\{x\}} $ in the big summation in \eqref{eq:big}.

The relation~\eqref{eq:doi} holds true after playing because the left hand side is a sum of two terms, both roughly of the form $xsx$, except that in the first we impose the condition that the exponent of $\A_{\{x\}}$ terms is nonzero on the left of $s$, and in the second that it is nonzero on the right of $s$. Therefore, most of the terms obtained after playing cancel in pairs. The only remaining ones are those obtained by playing either $sx$ or $xs$ (with no $\A_{\{x\}}$ powers on the left and right, respectively).

Finally, the fact that \eqref{eq:trei} holds true after playing boils down to Equation~\eqref{eq:wzw} in the definition of a hypercubical collection.
\end {proof}

Let $\A$ be a hypercubical collection modeled on a set $X$. For any $X' \subseteq X$,  the subcollection  
composed of $\A_Z$ for $Z \subseteq X'$ is a hypercubical collection modeled on $X'$, which we denote by $\A|X'$. Further, if $\dd = (d_x)_{x \in X}$ is a series of nonnegative integers, by picking only the terms $d_x$ for $x \in X'$ we obtain a new series, denoted $\dd|X'$.

The following is a straightforward consequence of the definition of playing, taking into account Lemma~\ref{lemma:plays}:
\begin {corollary}
\label {cor:plays}
Let $\A$ be a hypercubical collection in an algebra $A$, modeled on a set $X$.  Let $X = X' \amalg X''$ be a decomposition of $X$ as a disjoint union.  Then, we have
$$ \play_{\A|X'}^{\dd|X'} (s') *  \play_{\A|X''}^{\dd|X"} (s'')  =  \play_{\A}^{\dd} (s's''),$$
for any $s' \in S(X'), s'' \in S(X'')$.  
\end {corollary}

Suppose now that $X$ is a finite, totally ordered set. We equip all $Z \subseteq X$ with the induced total ordering.

Let $\A$ be a hypercubical collection in an algebra $A$, modeled on $X$. Pick $\dd =(d_x)_{x \in X}, \ d_x > 0$ and, for any $Z \subseteq X$, define
$$\A^\dd_Z  = \play_{\A|Z}^{\dd|Z} (\alpha(Z)),$$
where $\alpha(Z)$ is the symphony on $Z$ from Definition~\ref{def:symph}.

\begin {lemma}
\label {lemma:ntom}
The elements $ \{ \A^\dd_Z \}_{Z \subseteq X}$ form a new hypercubical collection $\A^\dd$ in $A$.
\end {lemma}

\begin {proof}
We need to check that, for any $Z \subseteq X$,
\begin {equation}
\label {eq:ww}
\sum_{Z' \amalg Z'' = Z}  \A^\dd_{Z'} * \A^\dd_{Z''}  = 0.
\end {equation}

Indeed, the left hand side in \eqref{eq:ww} equals
\begin {equation}
\label {eq:www}
\sum_{Z' \amalg Z'' = Z}  \play_{\A|Z'}^{\dd |Z'} (\alpha(Z')) *  \play_{\A|Z''}^{\dd |Z''} (\alpha(Z'')) = \sum_{Z' \amalg Z'' = Z}  \play_{\A|Z}^{\dd |Z} (\alpha(Z') \alpha(Z'')) = 
\end {equation}
\begin {flushleft}
$ \hskip2.5cm =  \play_{\A|Z}^{\dd |Z} \bigl(\sum_{Z' \amalg Z'' = Z}  \alpha(Z') \alpha(Z'')\bigr)  = \play^{\dd |Z}_{\A|Z} (0)  = 0$.
\end {flushleft}
The first equality in \eqref{eq:www} is a consequence of Corollary~\ref{cor:plays}, the second of linearity (Lemma~\ref{lemma:plays}), and the third of Lemma~\ref{lemma:symph}.
 \end {proof}

\subsection {Back to compression}
\label {sec:compression2}

Let  $H =  \bigl( (C^{\eps})_{\eps \in \E(\dd)}, (\De^{\eps})_{\eps \in \E_n} \bigr)$ be an $n$-dimensional  hyperbox of chain complexes as in Section~\ref{sec:hyperv}. In Section~\ref{sec:compression1} we advertised the construction of a compressed hypercube $\hat H = (\hat C^{\eps}, \hat \De^{\eps})_{\eps \in \E_n}$, with
$$ \hat C^{(\eps_1, \dots, \eps_n)} = C^{(\eps_1d_1, \dots, \eps_n d_n)}.$$

We are now ready to explain the exact construction of the maps  $\hat \De^{\eps}$. Let $A$ be the algebra $\End(C)$ under composition, where $C = \oplus_{\eps \in \E(\dd)} C^{\eps}$. As mentioned in Example~\ref{ex:cde}, the maps $\De^{\eps}= \De^{\zeta(Z)}= \A_Z$ form a hypercubical collection $\A$ in $A$, modeled on $X = \{1, \dots, n\}$. (Note that every $\eps \in \E_n$ can be written as $\zeta(Z)$, for a unique $Z \subseteq X$.) 

For $Z \subseteq \{1, \dots, n\}$, set
\begin {equation}
\label {eq:tildeDe}
\hat \De^{\zeta(Z)} = \A^{\dd}_Z = \play_{\A|Z}^{\dd|Z} (\alpha(Z)). 
\end {equation} 
 
For example, when $Z = \{i\}, \eps = \tau_i$ for some $i\in \{1, \dots, n\}$, the map along the corresponding edge of the hypercube is $\hat \De^{\tau_i} = \play_{\A|\{i\}}^{(d_i)} (i) = (\De^{\tau_i})^{\circ d_i}$, as noted in Section~\ref{sec:compression1}.  For $n=2$, by playing the symphony on a set of two elements, we recover formula \eqref{eq:deo1o2} for $\hat \De^{(1,1)}$. 

\begin {proposition}
\label {prop:compressed}
$\hat H = (\hat C^{\eps}, \hat \De^{\eps})_{\eps \in \E_n}$ is a hypercube of chain complexes.
\end {proposition}

\begin {proof}
The relations \eqref{eq:d2} are a direct consequence of Lemma~\ref{lemma:ntom}. The fact that $\hat \De^{\eps}$ changes grading by $\| \eps \|-1$ (as required in the definition of a hypercube) follows from the similar property for the maps $\De^{\eps}$, together with the second relation in ~\eqref{eq:h1}, which is satisfied by all the terms appearing in a symphony on a set of size $n$.
\end {proof}

\subsection {Chain maps and homotopies}
\label {sec:chmaps}
Let 
$$\leftexp{0}{H} =  \bigl( (\leftexp{0}{C}^{\eps})_{\eps \in \E(\dd)}, (\leftexp{0}{\De}^{\eps})_{\eps \in \E_n} \bigr),  \ \  \leftexp{1}{H} =  \bigl( (\leftexp{1}{C}^{\eps})_{\eps \in \E(\dd)}, (\leftexp{1}{\De}^{\eps})_{\eps \in \E_n} \bigr)$$  
be two hyperboxes of chain complexes, having the same size $\dd \in \N^n$. Let $(\dd, 1) \in \N^{n+1}$ be the sequence obtained from $\dd$ by adding $1$ at the end.

\begin {definition}
\label {def:chmap}
A {\em chain map} $F : \leftexp{0}{H} \to \leftexp{1}{H}$  is a collection of linear maps 
$$F_{\eps^0}^\eps: \leftexp{0}{C}^{\eps^0}_* \to \leftexp{1}{C}^{\eps^0 + \eps} _{*+ \| \eps \|},$$
satisfying 
\begin{equation}
\label{eq:DF}
\sum_{\eps' \leq \eps} \bigl( \De^{\eps - \eps'}_{\eps^0 + \eps'} \circ F^{\eps'}_{\eps^0}  +   F^{\eps - \eps'}_{\eps^0 + \eps'} \circ \De^{\eps'}_{\eps^0} \bigr)  = 0,
\end{equation}
for all $\eps^0 \in \E(\dd), \eps \in \E_n$ such that $\eps^0 + \eps \in \E(\dd)$. 
\end {definition}

In other words, a chain map between the hyperboxes $\leftexp{0}{H} $ and $\leftexp{1}{H}$ is an $(n+1)$-dimensional hyperbox of chain complexes, of size $(\dd, 1)$, such that the sub-hyperbox corresponding to $\eps_{n+1}=0$ is $\leftexp{0}{H}$ and the one corresponding to $\eps_{n+1}=1$ is $\leftexp{1}{H}$. The maps $F$ are those maps $\De$ in the new hyperbox that increase $\eps_{n+1}$ by $1$.

Note that, in the particular case when $\dd = (0, \dots, 0)$, so that $\leftexp{0}{H} $ and $\leftexp{1}{H}$ are ordinary chain complexes, the notion of chain map coincides with the usual one. Also note that, when $\dd = (1, \dots, 1)$, so that $\leftexp{0}{H} $ and $\leftexp{1}{H}$ are hypercubes, a chain map $F$ induces an ordinary chain map $F_\tot$ between the corresponding total complexes $\leftexp{0}{C}_\tot$ and  $\leftexp{1}{C}_\tot$.

The identity chain map $\Id: H \to H$ is defined to consist of the identity maps $\Id_{\eps^0}^\eps$ when $\eps = (0, \dots, 0)$, and zero for other $\eps$.

\begin {lemma}
\label {lemma:chmap}
A chain map $F$ between hyperboxes $\leftexp{0}{H} $ and $\leftexp{1}{H}$ induces a natural chain map $\hat F$ between the compressed hypercubes  $\leftexp{0}{\hat H} $ and $\leftexp{1}{\hat H}$.
\end {lemma}

\begin {proof}
As mentioned above, the information in $F$ can be used to build a new hyperbox of size $(\dd, 1)$  composed of
$\leftexp{0}{H} $ and $\leftexp{1}{H}$. Compressing this bigger hyperbox gives the required map $\hat F$.
\end {proof} 

\begin {definition}
Let  $F : \leftexp{0}{H} \to \leftexp{1}{H}$ and  $\G : \leftexp{1}{H} \to \leftexp{2}{H}$ be chain maps between hyperboxes of the same size. Their composite $\G \circ F$ is defined to consist of the maps
$$ (\G \circ F)_{\eps^0}^{\eps} = \sum_{\{\eps'| \eps' \leq \eps \}}   \G_{\eps^0+\eps'}^{\eps-\eps'} \circ F_{\eps^0}^{\eps'}.$$
\end {definition}

\begin {definition}
Let $F, G$ be two chain maps between hyperboxes $\leftexp{0}{H} $ and $\leftexp{1}{H}$. A {\em chain homotopy} between $F$ and $\G$  is a collection $\Psi$ of linear maps 
$$ \Psi_{\eps^0}^\eps: \leftexp{0}{C}^{\eps^0}_* \to \leftexp{1}{C}^{\eps^0 + \eps} _{*+1+  \| \eps \|},$$
for all satisfying 
$$F_{\eps^0}^\eps - \G_{\eps^0}^\eps = \sum_{\eps' \leq \eps} \Bigl(  \De^{\eps - \eps'}_{\eps^0 + \eps'} \circ \Psi^{\eps'}_{\eps^0}  +   \Psi^{\eps - \eps'}_{\eps^0 + \eps'} \circ \De^{\eps'}_{\eps^0}   \Bigr) ,$$
for $\eps^0 \in \E(\dd), \eps \in \E_n$ such that $\eps^0 + \eps \in \E(\dd)$. 
\end {definition}

Note that we can also interpret a chain homotopy as a bigger hyperbox, namely an $(n+2)$-dimensional one of size $(\dd, 1, 1)$, where $\eps_{n+1} = 0$ and $\eps_{n+1} =1$ are the $(n+1)$-dimensional hyperboxes corresponding to $F$ and $\G$, respectively, and the maps in the new direction (from $\eps_{n+2}=0$ to $\eps_{n+2} =1$) are the identities (preserving $\eps_{n+1}$) and $\Psi$ (increasing $\eps_{n+1}$ by one).

We define a {\em chain homotopy equivalence} between hyperboxes as a chain map that has an inverse up to chain homotopy. The following follows from the same kind of argument as Lemma~\ref{lemma:chmap}:

\begin {lemma}
\label {lemma:chainhe}
If $F :  \leftexp{0}{H} \to \leftexp{1}{H}$ is a chain homotopy equivalence, then the compressed map $\hat F$ is also a chain homotopy equivalence.
\end {lemma}

Observe that if $F :  \leftexp{0}{H} \to \leftexp{1}{H}$ is a chain homotopy equivalence between hypercubes, the map $F_{\tot}$ between the respective total complexes is an ordinary chain homotopy equivalence.

\begin {definition}
\label {def:quasiiso}
A chain map $F :  \leftexp{0}{H} \to \leftexp{1}{H}$ is called a {\em quasi-isomorphism of hyperboxes} if its components 
$$F^{(0, \dots, 0)}_{\eps^0} :   \leftexp{0}{C}^{\eps^0} \to \leftexp{1}{C}^{ \eps^0} $$
induce isomorphisms on homology, for all $\eps^0 \in \E(\dd)$. 
\end {definition}

We note that a chain homotopy equivalence of hyperboxes is a quasi-isomorphism. Further, if $F :  \leftexp{0}{H} \to \leftexp{1}{H}$ is a quasi-isomorphism between hypercubes, the total map $F_{\tot}$ is an ordinary quasi-isomorphism.

\subsection {Elementary enlargements}
\label {sec:ele}

For future reference, we introduce here a simple operation on hyperboxes, called {\em elementary enlargement}. 

Let $H =  \bigl( (C^{\eps})_{\eps \in \E(\dd)}, (\De^{\eps})_{\eps \in \E_n} \bigr)$ be a hyperbox of chain complexes, of size $\dd \in \N^n$. Pick $k \in \{1, \dots, n\}$ and $j \in \{0, 1, \dots, d_k\}$. Define $\dd^+ = ({d}^+_1, \dots, {d}^+_n )\in \N^n$ by $\dd^+ = \dd + \tau_k$, i.e.
$${d}^+_i = \begin {cases} d_i & \text{if  } i \neq k\\ d_i + 1& \text{if  } i=k.\end {cases}$$

We construct a new hyperbox
$$H^+ =  \bigl( ({C}^{+,\eps})_{\eps \in \E(\dd^+)}, ({\De}^{+,\eps})_{\eps \in \E_n} \bigr)$$
by replicating (i.e. introducing a new copy of) the complexes in positions with $\eps_k = j$. The new copy will be in position $\eps_k = j+1$, and everything with higher $\eps_k$ is shifted one step to the right. The two identical copies (which can be thought of as sub-hyperboxes) are linked by the identity chain map. Precisely, we set
$$ C^{+, \eps} = \begin {cases} C^\eps & \text{if  } \eps_k \leq j \\ C^{\eps - \tau_k} & \text{if  } \eps_k \geq j+1 \end {cases} $$
and
$$ \De^{+, \eps}_{\eps^0} =  
\begin {cases} 
\De^\eps_{\eps^0} & \text{if  }  \eps^0_k + \eps_k  \leq j \\
\Id & \text{if } \eps^0_k = j,  \eps=\tau_k \\
0 & \text{if } \eps^0_k = j,  \eps_k = 1, \|\eps\| \geq 2 \\
\De^{\eps}_{\eps^0- \tau_k} & \text{if  } \eps^0_k > j.
 \end {cases} $$

We then say that $H^+$ is obtained from $H$ by an elementary enlargement at position $(j, k)$. The following is easy to check from the definitions:

\begin {lemma}
\label {lemma:ci}
Let $H^+$ be an elementary enlargement of a hyperbox of chain complexes $H$. Then the compressed hypercubes $\hat H$ and $\hat H^+$ are identical.
\end {lemma}

\subsection {Canonical inclusions}
\label {sec:caninc}
We describe here yet another construction that will be useful to us later. We restrict to the case of hypercubes (since this is all we need), but everything can also be done in the more general context of hyperboxes.

\begin {definition}
\label {def:canhyper}
Let $(K_*, \del)$ be a chain complex. The $n$-dimensional {\em canonical hypercube} $H(K, n)$ associated to $(K_*, \del)$ is defined to consist of the vector spaces 
$$ C^\eps_* = K_*, \ \eps \in \E_n,$$
together with the maps
$$ D_{\eps^0}^\eps = \begin {cases}
\del &\text{if } \|\eps\| = 0,\\
\Id &\text{if } \|\eps \| =1,\\
0 & \text{if } \|\eps \| \geq 2.\end {cases}$$  
\end {definition}
 
Now suppose that $H = (C^\eps, D ^\eps)_{\eps \in \E_n}$ is an arbitrary $n$-dimensional hypercube of chain complexes. Our aim is to construct a chain map
$$ F^\can_H : H(C^{(0,\dots, 0)}, n) \to H,$$
which will be called the {\em canonical inclusion}.

When $n=0$ this is simply the identity map. When $n=1$, the hypercube $H$ consists of a single chain map $f = D^{(1)}: C^{(0)} \to C^{(1)}$ between two chain complexes. The canonical inclusion is then the square:
$$ \begin {CD}
C^{(0)} @>{\Id}>> C^{(0)} \\
@V{\Id}VV @VV{f}V \\
C^{(0)} @>{f}>> C^{(1)}
\end {CD}$$
with the diagonal map being trivial.

For general $n$, we construct the canonical inclusion as a composition of $n$ different chain maps as follows. For $i \in \{0, \dots, n\}$ and $\eps=(\eps_1, \dots, \eps_n) \in \E_n$ we let $\eps[\leq i] \in \E_n$ be the multi-index obtained from $\eps$ by changing all entries above $i$ to be zero. In other words,
$$ \eps[\leq i]_j = \begin{cases}
\eps_j & \text{if } j \leq i, \\
0 & \text{if } j > i.
\end {cases}
$$
Similarly, we let $\eps[>i]$ to be obtained from $\eps$ by setting all entries less than or equal to $i$ to be zero. 

We define a hypercube $H[i]$ to consist of the chain groups
$$ C[i]^\eps = C^{\eps[\leq i]}, \ \eps \in \E_n,$$
and the maps
$$ D[i]_{\eps}^{\eps'-\eps} = \begin {cases}
D_{\eps[\leq i]}^{(\eps' - \eps)[\leq i]} &\text{if } \eps[>i] = \eps'[>i],\\
\Id &\text{if } \eps[\leq i] = \eps'[\leq i] \text{ and } \|\eps'[>i] - \eps[>i] \| =1,\\
0 & \text{otherwise.}
\end {cases} $$

Note that $H[0]$ is the canonical hypercube $H(C^{(0,\dots, 0)}, n)$, while $H[n] = H$. 

For $i=1, \dots, n$, we define chain maps
$$ F[i]: H[i-1] \to H[i]$$
to consist of 
$$ F[i]_{\eps}^{\eps' - \eps} =  \begin {cases}
D_{\eps[\leq (i-1)]}^{\eps'[\leq i]-\eps[\leq (i-1)]} &\text{if } \eps_i =1, \ \eps[>i] = \eps'[>i],\\
\Id &\text{if } \eps = \eps' \text{ and } \eps_i =0, \\
0 & \text{otherwise.}
\end {cases} $$

The canonical inclusion is then:
$$ F^\can_H = F[n] \circ F[n-1] \circ \dots \circ F[1]. $$

\section {Quasi-stabilizations}
\label {sec:triangles}

In this section we introduce a new move that relates certain equivalent, multi-pointed Heegaard diagrams, called quasi-stabilization. Basically, a quasi-stabilization is the composition of a free index zero/three stabilization and some handleslides.

Our goal is to study how the polygon maps on Heegaard Floer complexes behave under this move. There are two motivations for this. First, the behavior of polygon maps under ordinary (free) index zero/three stabilizations (which can be viewed as particular examples of quasi-stabilizations) is one of the inputs in the construction of complete systems of hyperboxes in Section~\ref{sec:hyperHeegaard}, as well as in the proof of the Surgery Theorem~\ref{thm:surgery}. Second, the more general quasi-stabilizations are needed in Section~\ref{sec:grid}, where they appear in the context of grid diagrams.

For concreteness, we will first describe the results in the case of triangles (for general quasi-stabilizations). Then we will explain how similar arguments can be used to study higher polygon maps. At the end we will specialize to the case of ordinary index zero/three stabilizations.

\subsection {The set-up}
\label {sec:setup}

Let $\Hyper = (\Sigma, \alphas, \betas, \ws, \zs)$ be a multi-pointed Heegaard diagram, as in Section~\ref{sec:hed}. Suppose $\Hyper$ represents a link $\orL$ in an integral homology sphere $Y$.  Fix $\s \in \bH(L)$, so that we have a well-defined generalized Floer complex $\Am(\Hyper,\s)=\Am(\Ta, \Tb, \s)$.

Let $g$ be the genus of $\Sigma $ and $d$ the number of  alpha (or beta) curves. We assume that $\beta_1 \in \betas$ bounds a disk containing a free basepoint $w_1$, and that the only alpha curve intersecting $\beta_1$ is $\alpha_1$, which does so at two points $x$ and $x'$. On the other hand, $\alpha_1$ can intersect other beta curves. We also require that the homology class $[\alpha_1]$ lies in the span of the other $[\alpha_i]$ in $H_1(\Sigma)$.
 
Let $\bar \Hyper = (\Sigma,  \bar \alphas, \bar \betas, \bar \ws, \zs)$ be the diagram obtained from $\Hyper$ by deleting $\alpha_1, \beta_1, $ and $w_1$. We then say that $\bar \Hyper$ is obtained from $ \Hyper$ by {\em quasi-destabilization}. The reverse process is called {\em quasi-stabilization}.

\begin {remark}
By handlesliding $\alpha_1$ over other alpha curves, we can arrange so that it does not intersect any beta curve except $\beta_1$. The resulting diagram is then a usual free index zero/three stabilization of $\bar \Hyper$; see \cite{Links}. 
\end {remark}

\begin{figure}
\begin{center}
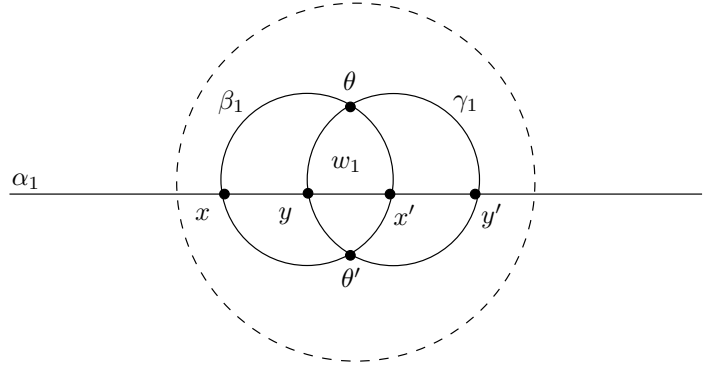
\end{center}
\caption {{\bf Quasi-stabilization.} We show here a part of the triple Heegaard diagram described in Section~\ref{sec:setup}.  In Section~\ref{sec:degeneration}, we will stretch the complex structure along the dashed curve.}
\label{fig:trio}
\end{figure}

Consider now an extra collection of $d$ attaching curves  $\gammas$ on $\Sigma$, such that $\gamma_1 \in \gammas$ 
has the same properties as $\beta_1$: it bounds a disk containing $w_1$, and the only alpha curve that it intersects is $\alpha_1$, with the respective intersection consisting of two points $y$ and $y'$. Furthermore, we assume that $\gamma_1$ is a small Hamiltonian translate of $\beta_1$, and intersects $\beta_1$ in two points $\theta$ and $\theta'$, as in Figure~\ref{fig:trio}. We assume that the relative positions of $x, x', y, y', \theta$ and $\theta'$ are exactly as in the figure.

Let $ \bar \gammas$ be the collection of curves obtained from $ \gammas$ by removing $ \gamma_1$. Then $(\Sigma, \alphas, \gammas, \ws, \zs)$ is a quasi-stabilization of $(\Sigma, \bar \alphas, \bar \gammas, \bar \ws, \zs)$. For any $\x \in \Ta \cap \Tb$, the intersection $\x \cap \alpha_1 $ is either $x$ or $x'$. We denote by $\bar \x \in \T_{\bar \alpha} \cap  \T_{\bar \beta}$ the generator obtained from $\x$ by deleting the point in $\x \cap \alpha_1$. Similarly, for $\y \in \Ta \cap \Tg$, there is a corresponding generator $\bar \y$ in $\T_{\bar \alpha}\cap  \T_{\bar \gamma}$, obtained by deleting $y$ or $y'$. 

Pick an intersection point $\thetas$ in $\Tb \cap \Tg$ such that $\theta \in \thetas$. We have a map
$$ F: \Am(\Ta, \Tb,\s) \to \Am(\Ta, \Tg, \s), \ \ \ F(\x) = f(\x \otimes \thetas),$$ 
which counts index zero pseudo-holomorphic triangles with one vertex at $\thetas$, as in Section~\ref{sec:polygon}.

Set $\bar \thetas = \thetas - \{\theta\} \in \T_{\bar \beta} \cap \T_{\bar \gamma}$. There is a corresponding triangle map in the quasi-destabilized diagram:
$$ \bar F: \Am(\T_{\bar \alpha} , \T_{\bar \beta}, \s) \to \Am(\T_{\bar \alpha},  \T_{\bar \gamma}, \s), \ \ \ \bar F(\bar \x) = f(\bar \x \otimes \bar \thetas).$$ 

Next, we define a map $G: \Am(\Ta, \Tb, \s) \to \Am(\Ta,
\Tg, \s)$ as follows. Writing the coefficients of $\bar F$ as 
$n_{\bar \x,\bar y}$, so that  for $\x \in \Ta \cap \Tb$, 
$$\bar F(\bar \x) = \sum_{\bar \y \in \T_{\bar \alpha}\cap  \T_{\bar \gamma}} n_{\bar \x, \bar \y} \bar \y,$$ we set
$$G( \x) = \sum_{\bar \y \in \T_{\bar \alpha}\cap  \T_{\bar \gamma}} n_{\bar \x, \bar \y} (\bar \y \cup (\x \cap \alpha_1)).$$ 
 
\begin {proposition}
\label {prop:degen}
For a quasi-stabilized triple Heegaard diagram $(\Sigma, \alphas, \betas, \gammas, \ws, \zs)$ as above, and suitable almost complex structures on the symmetric products, the maps $F$ and $G$ coincide.
\end {proposition}

The proof of Proposition~\ref{prop:degen} will occupy Sections~\ref{sec:lipshitz}--\ref{sec:degeneration}.

\subsection {Cylindrical formulations}
\label {sec:lipshitz}

We recall Lipshitz's cylindrical formulation of Heegaard 
Floer homology \cite{LipshitzCyl}, see also \cite[Section 5.2]{Links}. Instead of holomorphic strips $[0,1]\times \rr \to \Sym^d(\Sigma)$ with boundaries on $\Ta$ and $\Tb$, Lipshitz considers pseudo-holomorphic maps from a Riemann surface $S$ (with boundary) to the target
$$  W = \Sigma \times [0,1] \times \rr.$$

The four-manifold $W$ admits two natural projection maps
$$ \pi_\Sigma: W \longrightarrow \Sigma \ \text{ and } \ \ \pi_\Disk: W \longrightarrow [0,1] \times \rr.$$

The notation $\pi_\Disk$ refers to the unit disk $\Disk \subset \cc$, which can be viewed as the conformal compactification of $[0,1] \times \rr$, obtained by adding the points $\pm i$.  

We equip $W$ with an almost complex structure $J$ translation invariant in the $\rr$-factor, and such that $\pi_\Disk$ is a pseudo-holomorphic map. Further, we ask for $J$ to be tamed by a natural split symplectic form on $W$. Typically, we choose $J$ to be a small perturbation of a split complex structure $j_\Sigma \times j_\Disk$, where $j_\Sigma$ and $j_\Disk$ are complex structures on $\Sigma$ and $[0,1] \times \rr$, respectively. Sometimes (for example, to ensure positivity of intersections) it will be convenient to require $J$ to be {split on $U$}, that is, split on $U \times [0,1] \times \rr$, where $U \subset \Sigma$ is an open subset. 

\begin {definition}
\label {def:annoy}
An {\em annoying curve} is a pseudo-holomorphic curve in $W$ contained in a fiber of $\pi_\Disk$.
\end {definition}
 
To define the differential on the cylindrical Heegaard Floer complex, Lipshitz uses pseudo-holomorphic maps
$$ u: S \longrightarrow W = \Sigma \times [0,1] \times \rr$$
with the following properties:
\begin {itemize}
\item $S$ is a Riemann surface with boundary and $2d$ punctures on its boundary, of two types: $d$ ``positive'' punctures $\{p_1, \dots, p_d\}$ and $d$ ``negative'' punctures $\{q_1, \dots, q_d\}$;
\item $u$ is a smooth embedding;
\item $u(\del S) \subset (\alphas \times \{1\} \times \rr) \cup (\betas \times \{0\} \times \rr)$;
\item $u$ has finite energy;
\item For each $i$, $u^{-1}(\alpha_i \times \{1\} \times \rr)$ and $u^{-1}(\beta_i \times \{0\} \times \rr)$ consist of exactly one component of $\del S -\{p_1, \dots, p_d, q_1, \dots, q_d\}$;
\item No components of the image $u(S)$ are annoying curves;
\item Any sequence of points in $S$ converging to $q_i$ resp. $p_i$ is mapped under $\pi_\Disk$ to a sequence of points whose second coordinate converges to $-\infty$ resp. $+\infty$.
\end {itemize}

Curves of this kind are called {\em cylindrical flow lines}. Any cylindrical flow line $u$ can be extended to a map $\bar u$ from the closure of $S$ to the compactification $\Sigma \times \Disk$. The image of this extension contains the points $x \times \{-i\}$ and $y \times \{i\}$, for $x \in \x, y \in \y$, where $\x, \y \in \Ta \cap \Tb \subset \Sym^d(\Sigma)$. We then say that $u$ connects $\x$ to $\y$.

To every cylindrical flow line $u: S \to W$ one can associate a strip $\tilde u: [0,1] \times \rr \to \Sym^d(\Sigma)$ with boundaries on $\Ta $ and $\Tb$, by setting $\tilde u(z) = \pi_\Sigma((\pi_\Disk \circ u)^{-1}(z))$. Thus, cylindrical flow lines can be organized according to moduli spaces $\M(\phi)$, indexed by homology classes  $\phi \in \pi_2(\x, \y)$ for the corresponding Whitney disks. Moreover, in \cite[Appendix A]{LipshitzCyl}, Lipshitz identifies the moduli spaces of cylindrical flow lines in a class $\phi$ with the respective moduli spaces of pseudo-holomorphic strips (ordinary flow lines), for suitable almost complex structures, in the case when the Maslov index $\mu(\phi)$ is one. It follows that the Heegaard Floer complex can be defined just as well by counting cylindrical flow lines instead of pseudo-holomorphic strips.
 
When studying degenerations of cylindrical flow lines (for example, in the proof that $\del^2 = 0$ in the cylindrical setting), we also encounter maps of the following kind:
 
\begin {definition}
Consider a Riemann surface $S$ with boundary and $d$ punctures $\{p_1, \dots, p_d\}$ on its boundary. A (cylindrical) {\em boundary degeneration} is a pseudo-holomorphic map $ u:  S \to \Sigma \times (- \infty, 1] \times \rr$ which has finite energy, is a smooth embedding, sends $\del S$ into $\alphas \times \{1\} \times \rr$, contains no component in the fiber of the projection to $(-\infty, 1] \times \rr$, and has the property that each component of $u^{-1}(\alpha_i \times \{1\} \times \rr)$ consists of exactly one component of $\del S \setminus \{p_1, \dots, p_d\}$. A similar definition can be made with $\beta$ playing the role of $\alpha$, and using the interval $[0, \infty)$ instead of $(-\infty, 1]$.
\end {definition}

Note that for a boundary degeneration $u$, the points at infinity must be mapped
to a fixed $\x \in \Ta$. The boundary degenerations with endpoint $\x$ can be organized into moduli spaces $\Nod(\psi)$ according to homology classes $\psi \in \pi_2^\alpha(\x) \cong H_2(\Sigma, \alphas)$. 

\begin{figure}
\begin{center}
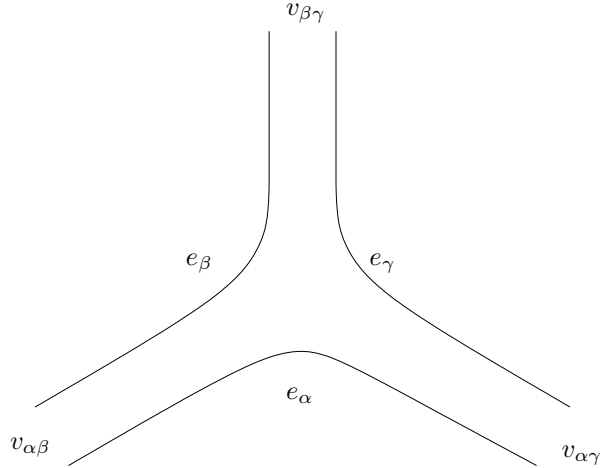
\end{center}
\caption {{\bf The triangular region $\Delta$.} This region is conformally equivalent to a triangle with punctures at the vertices.}
\label{fig:Delta}
\end{figure}

Next, let us recall from Section~\ref{sec:polygon} that when one has three collections of curves $\alphas, \betas, \gammas$ on fixed Heegaard
surface with marked basepoints $(\Sigma, \ws)$,  one can define a map
$$ f = f_{\alpha \beta \gamma}: \Am(\Ta, \Tb, \s) \otimes \Am(\Tb, \Tg, \zero) \longrightarrow \Am(\Ta, \Tg, \s)$$
by counting index zero pseudo-holomorphic triangles in $\Sym^d(\Sigma)$, with boundaries on $\Ta, \Tb$ and $\Tg$. These maps admit a cylindrical formulation, too, see \cite[Section 10]{LipshitzCyl}. Indeed, consider a contractible subset $\Delta \subset \cc$ as in Figure~\ref{fig:Delta}, with three boundary components $e_\alpha, e_\beta$ and $e_\gamma$, and three infinite ends $v_{\alpha\beta}, v_{\beta \gamma}, v_{\alpha \gamma}$, all diffeomorphic to $[0,1] \times (0, \infty)$. Setting
$$ W_\Delta = \Sigma \times \Delta,$$
note that there are natural projections $ \pi_\Sigma, \pi_\Delta$ to the two factors. We equip $W_\Delta$ with an almost complex structure having properties analogous to those of the almost complex structure on $W$. 

We then consider pseudo-holomorphic maps
$$ u: S \longrightarrow W_\Delta$$
with the following properties:
\begin {itemize}
\item $S$ is a Riemann surface with boundary and $3d$ punctures $p^{\alpha\beta}_i, p^{\beta \gamma}_i, p^{\alpha \gamma}_i,  \  i\in \{1, \dots, d\}$, on the boundary $\del S$;
\item $u$ is a smooth embedding;
\item $u(\del S) \subset (\alphas \times e_\alpha) \cup (\betas \times e_\beta) \cup (\gammas \times e_\gamma)$;
\item $u$ has finite energy;
\item For each $i=1, \dots, d$ and $\sigma \in \{\alpha, \beta, \gamma \}$, the preimage  $u^{-1}(\sigma_i \times e_\sigma)$ consists of exactly one component of the punctured boundary of $S$;
\item No components of the image $u(S)$ are annoying curves;
\item Any sequence of points in $S$ converging to $p^{\alpha\beta}_i$ (resp. $p^{\beta \gamma}_i, p^{\alpha \gamma}_i$) is mapped under $\pi_\Delta$ to a sequence of points converges towards infinity in the strip-like end $v_{\alpha \beta}$ (resp. $v_{\beta \gamma}, v_{\alpha \gamma}$).
\end {itemize}

Maps $u$ like this are called {\em cylindrical triangles}. They can be organized into moduli spaces $\M(\phi)$ according to homology classes $\phi \in \pi_2(\x,  \y, \zs)$, for $\x \in \Ta \cap \Tb, \y \in \Tb \cap \Tg, \zs \in \Ta \cap \Tg$. The moduli space of such maps in classes $\phi$ with $\mu(\phi) = 0$ can be identified with the moduli space of ordinary pseudo-holomorphic triangles in $\Sym^d(\Sigma)$, as  used in the definition of the map $f= f_{\alpha \beta \gamma}$. (This can be proved by the same methods as in \cite[Section 13]{LipshitzCyl}.) It follows that $f$ can be described in terms of counts of cylindrical triangles.

Similar descriptions can be given to the higher polygon maps from Section~\ref{sec:polygon}.

 \subsection {Domains}
 \label {sec:domains}
Let $\phi \in \pi_2(\x, \y)$ be a homology class of Whitney disks in a Heegaard diagram $(\Sigma, \alphas, \betas, \ws)$. The curves $\alphas$ and $\betas$ split the Heegaard surface into several connected components $R_1, \dots, R_r$, which we call {\em regions}.  A {\em domain} $\D$ on the Heegaard diagram is by definition a  linear combination of regions, with integer coefficients. The class $\phi$ has an associated {\em domain} $\D(\phi)$, see \cite[Section 3.5]{HolDisk}:
$$ \D(\phi) = \sum_{i=1}^r (\phi \cdot [\{z_i\} \times \Sym^{d-1}(\Sigma)])   R_i,$$
where $z_i$ is a point chosen in the interior of the region $R_i$, and $\cdot$ denotes intersection product.
  
Let $\D$ be a domain and $x \in \alpha_i \cap \beta_j$ an intersection point, for some $i, j$. A neighborhood of $x$ is split by $\alpha_i$ and $\beta_j$ into four quadrants. Two of the four quadrants have the property that as we move counterclockwise around $x$, we first see $\alpha_i$ on their boundary and then $\beta_j$; for the other two quadrants, we first see $\beta_j$ and then $\alpha_i$. Let $m^{\alpha \beta}(\D, x)$ the sum of the multiplicities of $\D$ in the two quadrants of the first type, and $m^{\beta \alpha}(\D, x)$ the sum of multiplicities in the other two quadrants.

Given a point $x  \in \alpha_i \cap \beta_j$ and a $d$-tuple $\x \in \Ta \cap \Tb$, we set 
$$\delta(\x, x) = \begin {cases} 1 & \text{if } x \in \x, \\
0 & \text{otherwise.}
\end {cases}
$$

\begin {definition}
\label {def:acc}
Let $\D$ be a domain on  $(\Sigma, \alphas, \betas)$ and $\x, \y \in \Ta \cap \Tb$. The domain $\D$ is said to be {\em acceptable for the pair $(\x, \y)$} if, for every $i, j =1, \dots, d$ and $x \in \alpha_i \cap \beta_j$, we have
\begin {equation}
\label {eq:corners}
m^{\alpha \beta}(\D, x) -  m^{\beta \alpha}(\D, x) = \delta(\x, x) - \delta(\y, x).
\end {equation}
\end {definition}

The proof of the following lemma is straightforward:
\begin {lemma}
\label {lemma:acc}
A domain $\D$ is acceptable for the pair $(\x, \y)$ if and only if it is of the form $\D(\phi)$ for some $\phi \in \pi_2(\x, \y)$.
\end {lemma}

We now turn to the Maslov index $\mu(\phi)$, which is
the expected dimension of the moduli space of pseudo-holomorphic representatives of $\phi \in \pi_2(\x, \y)$. The Maslov index can be calculated in terms of the domain $\D=\D(\phi)$ using the following formula due to Lipshitz \cite[Corollary 4.3]{LipshitzCyl}:

\begin {equation}
\label {eq:lipshitz1}
 \mu(\phi) = e(\D) + \sum_{x \in \x} n_x(\D) + \sum_{y \in \y} n_y(\D).
\end {equation} 

Here, $n_p(\D)$ denotes the average multiplicity of $\D$ in the four quadrants around a point $p$, while $e(\D)$ is the Euler measure of the domain, as defined in \cite{LipshitzCyl}.  
 
Now consider a boundary degeneration class $\psi \in \pi_2^\alpha(\x)$. Its domain is then an {\em $\alpha$-periodic domain} $\Per = \Per(\psi) \in H_2(\Sigma, \alphas)$, i.e. a linear combination of components of $\Sigma - \alphas$. In fact, there is a one-to-one correpsondence between periodic domains and classes in $ \pi_2^\alpha(\x)$. The respective Maslov index is given by
\begin {equation}
\label {eq:muper1}
 \mu(\psi) = e(\Per) + 2\sum_{x \in \x} n_x(\Per).
 \end {equation} 
 
We also have an alternate characterization, see \cite[Lemma 5.4]{Links}. Recall that we have a basepoint $w_i$ in each component of $\Sigma - \alphas$. Then:
\begin {equation}
\label {eq:muper2}
 \mu(\psi) = 2 \sum_{i=1}^{d-g+1} n_{w_i}(\Per).
\end {equation}  

Next, we turn to homology classes of triangles. Let $(\Sigma, \alphas, \betas, \gammas, \ws, \zs)$ be a triple Heegaard diagram, with each curve collection consisting of $d$ curves. By {\em regions} we now mean the connected components of $\Sigma \setminus (\alphas \cup \betas \cup \gammas)$. Given a homology class $\phi \in \pi_2(\x, \y, \z)$, for $\x \in \Ta \cap \Tb, \y \in \Tb \cap \Tg, \zs \in \Ta \cap \Tg$,
 its domain $\D = \D(\phi)$ is defined as before. We have analogues of Definition~\ref{def:acc} and Lemma~\ref{lemma:acc}:

\begin {lemma}
\label {lemma:acc3}
 The necessary  and sufficient conditions for a domain $\D$ to be of the form $\D(\phi)$ for some $\phi \in \pi_2(\x, \y, \z)$ is that $\D$ is {\em acceptable} for the triple $(\x, \y, \z)$, that is, it should satisfy:
\begin {eqnarray}
\label {eq:corners1}
m^{\alpha \beta}(\D, x) -  m^{\beta \alpha}(\D, x) &=& \delta(\x, x)   \text{ for } x \in \alpha_i \cap \beta_j, \\
\label {eq:corners2}
m^{\beta \gamma}(\D, y) -  m^{\gamma \beta}(\D, y) &=& \delta(\y, y)  \text{ for } y \in \beta_i \cap \gamma_j, \\
\label {eq:corners3}
m^{\gamma \alpha}(\D, z) -  m^{ \alpha \gamma}(\D, z) &=& \delta(\z, z) \text{ for } z \in \gamma_i \cap \alpha_j.
\end {eqnarray}
\end {lemma}

The vertex multiplicities of a domain are defined  as in the case of bigons. We can similarly define the Euler measure. Further, we let $a(\D)$ denote the intersection $\del \D \cap \alphas$, viewed as a $1$-chain on $\Sigma$, supported on $\alphas$. Similarly we define $b(\D) = \del \D \cap \betas$ and $c(\D) = \del \D \cap \gammas$. We let $a(\D) . c(\D)$ denote the average of the four algebraic intersection numbers between $a'(\D)$ and $c(\D)$, where $a'(\D)$ is a small translate of $a(\D)$ in any of the four ``diagonal'' directions off $\alphas$, such that no endpoint of $a(\D)$ lies on $\gammas$, and no endpoint of $c(\D)$ lies on $\alphas$. We could similarly define $b(\D) . a(\D)$ or $c(\D). b(\D)$. Sarkar \cite[Theorem 4.1]{SarkarIndex} proved the following formula for the index of holomorphic triangles:\footnote{Note that our conventions are different from those in \cite{SarkarIndex}. There, the sides of a holomorphic triangle were on $\alpha, \beta, \gamma$ in counterclockwise order, whereas we have them clockwise. Thus, Sarkar's expression $a(\D).c(\D)$ differs from ours by a negative sign. The formula \eqref{eq:sarkar} above is written with our conventions.}
\begin {equation}
\label {eq:sarkar}
\mu(\phi) =e(\D) +  \sum_{x \in \x} n_x(\D) + \sum_{y \in \y} n_y(\D) -  a(\D).c(\D)  - \frac{d}{2}.
\end {equation} 
 
Sarkar also gave a generalization of this formula to higher polygons. Suppose we have curve collections $\etas^i, i=0, \dots, l$ on a pointed surface $(\Sigma, \ws, \zs)$, such that each collection consists of $d$ curves. We then consider a homotopy class of $(l+1)$-gons $\phi \in \pi_2(\x^0, \dots, \x^l)$, where $\x^i \in \T_{\eta^i} \cap \T_{ \eta^{i+1}}$ for $i< l$ and $\x^l \in \T_{\eta^0} \cap \T_{\eta^l}$. We can define $\D(\phi)$ as before, and we have acceptability conditions similar to \eqref{eq:corners1}-\eqref{eq:corners3}. We let $a^i(\D) = \del \D \cap \etas^i$. With our conventions, Theorem 4.1 in \cite{SarkarIndex} says:
 \begin {equation}
\label {eq:sarkar2}
\mu(\phi) =e(\D) +  \sum_{x \in \x^0} n_x(\D) + \sum_{y \in \x^1} n_y(\D) -  a^0(\D) . \sum_{j=2}^l a^j(\D) - \sum_{j > k > 1} a^j(\D) . a^k(\D)  - \frac{d(l-1)}{2}.
\end {equation} 
 
Here, the Maslov index $\mu(\phi)$ is as defined in Section~\ref{sec:polygon}. 

\subsection {Convergence and gluing for the moduli spaces of triangles}
\label {sec:gluing}

Let $(\Sigma^1, \alphas^1, \betas^1, \gammas^1)$ and $(\Sigma^2, \alphas^2, \betas^2, \gammas^2)$ be two triple Heegaard diagrams. (For the purposes of this subsection, we can ignore the basepoints.) For $i=1,2$, we let $d_i$ be the number of curves in the collection $\alpha^i$ (or $\beta^i$, or $\gamma^i$), and $g_i \leq d_i$ the genus of $\Sigma^i$.

Consider an extra simple closed curve $\alpha^1_s$ on $\Sigma^1$ that is disjoint from the other curves in $\alphas^1$ and lies in their homological span. Set
$$ \alphas^{1+} = \alphas^1 \cup \{\alpha^1_s\}.$$
 
Pick also one of the curves in $\alphas^2, $ say $\alpha^2_1$, and call it $\alpha^2_s$. The subscript $s$ stands for ``special.'' 

Pick points $p_i \in \alpha^i_s, i=1,2$, that do not lie on any of the beta or gamma curves. We form the connected sum $\Sigma = \Sigma_1 \# \Sigma_2$ at $p_1$ and $p_2$, of genus $g = g_1 + g_2$. By joining each of the two ends of $\alpha^1_s$ with at $p_1$ with the respective end of $\alpha^2_s$ at $p_2$, we obtain a new curve $\alpha_s = \alpha_s^1 \# \alpha_s^2$ on $\Sigma$. We set
$$ \alphas = \alphas^1  \cup \bigl( \alphas^2 - \{\alpha^2_s\} \bigr) \cup \{\alpha_s\}.$$

This is a collection of $d=d_1 + d_2$ attaching curves on $\Sigma$. We can also form collections $\betas = \betas^1 \cup \betas^2$ and $\gammas = \gammas^1 \cup \gammas^3$. Together, they turn $\Sigma$ into a triple Heegaard diagram, which we call the {\em special connected sum} of $(\Sigma^1, \alphas^1, \betas^1, \gammas^1)$ and $(\Sigma^2, \alphas^2, \betas^2, \gammas^2)$.

\begin{figure}
\begin{center}
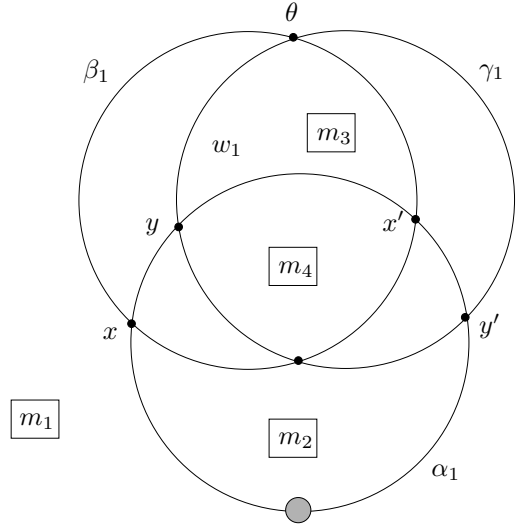
\end{center}
\caption {{\bf The sphere $\Sphere$.} This is the one-point compactification of the disk bounded by the dashed curve in Figure~\ref{fig:trio}. It could be viewed as a triple Heegaard diagram, except it is missing a basepoint. The compactification point $p^2$ is shown by a gray dot. The four boxes show the multiplicities of a triangular domain inside the corresponding regions.}
\label{fig:ess}
\end{figure}

\begin {example}
\label {ex:ess}
A triple Heegaard diagram $(\Sigma, \alphas, \betas, \gammas)$ as in Section~\ref{sec:setup} can be viewed as the special connected sum of the diagram $(\Sigma, \bar \alphas, \bar \betas, \bar \gammas)$ with the genus zero diagram $(\Sphere, \alpha_1, \beta_1, \gamma_1)$ shown in Figure~\ref{fig:ess}. The notation is as in Section~\ref{sec:setup}. 
\end {example}

Given a special connected sum of triple Heegaard diagrams, note that intersection points $\x^1 \in \T_{\alpha^1} \cap \T_{\beta^1}$ and $\x^2 \in \T_{\alpha^2} \cap \T_{\beta^2}$ give rise to an intersection point $\x^1 \times \x^2 \in  \T_{\alpha} \cap \T_{\beta}$. Conversely, any $\x \in  \T_{\alpha} \cap \T_{\beta}$ is of this form, because none of the points on $\alpha^1_s \cap \beta^1_i, i=1, \dots, d_1$, can be part of $\x$. Similar remarks apply to generators of the form $\y=\y_1 \times \y_2 \in \Tb \cap \Tg$ and $\z = \z_1 \times \z_2 \in \Tg \cap \Ta$.

Given a homology class $\phi \in \pi_2(\x, \y, \z)$ on the special connected sum, we denote by  $m_1 = m_1(\D)$ and $m_2= m_2(\D)$ are the multiplicities of $\D$ on each side of the curve $\alpha_s$, near the connected sum neck. (See Figure~\ref{fig:ess} for an example.)

\begin {lemma}
\label {lemma:split}
For a triple Heegaard diagram obtained as a special connected sum as above, pick $\x = \x^1 \times \x^2 \in  \T_{\alpha} \cap \T_{\beta}, \y=\y_1 \times \y_2 \in \Tb \cap \Tg$ and $\z = \z_1 \times \z_2 \in \Tg \cap \Ta$. Then, there is a natural surjective map:
\begin {equation}
p:  \pi_2(\x, \y, \z) \to \pi_2(\x^2, \y^2, \z^2).
\end {equation}

Furthermore, a choice of a domain $\Per \in H_2(\Sigma, \alphas^{1+})$ whose boundary contains $\alpha^1_s$ with multipicity one gives rise to a map
$$ \sigma: \pi_2(\x^2, \y^2, \z^2) \to \pi_2(\x, \y, \z)$$
such that $p \circ \sigma = \id$.
\end {lemma}

\begin {proof}
We use the identification between homology classes $\phi \in \pi_2(\x, \y, \z)$ and acceptable domains $\D = \D(\phi)$, see Lemma~\ref{lemma:acc3}. An acceptable domain $\D$ on the special connected sum gives rise to an acceptable domain $\D_2$ on $\Sigma_2$ with respect to $(\x^2, \y^2, \z^2)$, by restriction. This produces the map $p$.

Given $\Per$ is as in the statement of the lemma, the desired map $s$ takes an acceptable domain $\D_2$ on $\Sigma_2$ to the acceptable domain  $\D = \D_2 + (m_1(\D_2) - m_2(\D_2)) \Per$ on the special connected sum. \end {proof}

Let $\phi \in \pi_2(\x, \y, \z)$ have a domain $\D$. We define an equivalence relation on pairs $(\phi^1, \Per)$, where $\phi^1 \in \pi_2(\x^1, \y^1, \z^1)$ and $\Per \in H_2(\Sigma, \alphas^{1+})$ has $\alpha^1_s$ with multiplicity $m_1(\D) - m_2(\D)$ on its boundary. Two pairs $(\phi^1_1, \Per_1)$ and $(\phi^1_2, \Per_2)$ are set to be equivalent if $\phi^1_1 + \Per_1 = \phi^1_2 + \Per_2$, as two-chains on $\Sigma^1$. By identifying homology classes with acceptable domains (as in the proof of Lemma~\ref{lemma:split}), we see that every $\phi$ determines a unique such equivalence class $\phi^{1+}$. We set
$$ \mu(\phi^{1+}) = \mu(\phi^1) + \mu(\Per),$$
for any $(\phi^1, \Per) \in \phi^{1+}$.

\begin {lemma}
\label {lemma:mumu}
Let  $\phi \in \pi_2(\x, \y, \z)$ be a homology class of triangles in a triple Heegaard diagram obtained by special connected sum, as above. Let $\phi^2 \in \pi_2(\x^2, \y^2, \z^2)$ be its restriction to $\Sigma^2$, and $\phi^{1+}$ the corresponding equivalence class of pairs on $\Sigma^1$. Then:
$$ \mu(\phi) = \mu(\phi^{1+}) + \mu(\phi^2) - m_1(\D) - m_2(\D),$$
where $\D = \D(\phi)$ is the domain of $\phi$.
\end {lemma}

\begin {proof}
In Sarkar's formula \eqref{eq:sarkar}, all terms except $e(\D)$ are additive under the special connected sum. When adding up the Euler measures, we have to subtract $m_1(\D)+m_2(\D)$ because doing the special connected sum involves deleting two disks, each made of two bigons. Two of these four bigons have multiplicity $m_1(\D)$, the other two $m_2(\D)$, and the Euler measure of a bigon is $1/2$.  
\end {proof}

We now proceed to study holomorphic triangles on a special connected sum. We will use Lipshitz's cylindrical formulation from Section~\ref{sec:lipshitz}. 

Note that if a homology class $\phi$ (of cylindrical flow lines, boundary degenerations, triangles, etc.) admits pseudo-holomorphic representatives, the principle of positivity of intersections implies that the domain $\D(\phi)$ is a linear combination of regions with only nonnegative coefficients:
$$ \D(\phi) \geq 0,$$
see  \cite[Lemma 3.2]{HolDisk}.

In addition to the cylindrical flow lines, boundary degenerations and triangles from Section~\ref{sec:lipshitz}, when discussing special connected sums we will also need to study some new objects:

\begin {definition}
\label {def:annoy3}
Consider a Heegaard surface $\Sigma$ and a collection of attaching circles $\alphas$ on $\Sigma$. An {\em annoying $\alpha$-degeneration} is a holomorphic curve $u: S \to W_\Delta = \Sigma \times \Delta$ such that $S$ is a connected Riemann surface with boundary and punctures on the boundary, and there exists an unpunctured component $\del_0 S$ of the boundary $\del S$ satisfying $u(\del_0 S) \subset \alphas \times e_\alpha$. Here, $\Delta$ is as in Figure~\ref{fig:Delta}.
\end {definition}

\begin {lemma}
An annoying $\alpha$-degeneration $u: S \to W_\Delta$ is an annoying curve in the sense of Definition~\ref{def:annoy}, that is, all of $S$ is mapped to a fixed point $p \in e_\alpha$ under $\pi_\Delta \circ u$.
\end {lemma}

\begin {proof} Let $D(S)$ be the double of $S$ taken along the component $\del_0 S \subseteq \del S$, and $D(\Delta)$ the double of $\Delta$ along $e_\alpha$. We can extend $\pi_\Delta \circ u$ to a holomorphic map $f: D(S) \to D(\Delta)$ using Schwartz reflection. Since $\del_0S$ is compact and $e_\alpha$ is not, there exists some $z_0 \in \del_0 S$ with $f'(z_0) = 0$. If $f$ were not constant, it would have a branch point of order $k \geq 2$ at $z_0$. This contradicts the local model near $z_0$, which is that of a holomorphic function $f$ mapped to $\cc$ such that $\Re\  f(z) > 0$ for $\Re \ z > 0$.
\end {proof}

Annoying $\alpha$-degenerations can be organized according to their domains, which are relative homology classes $\Per \in H_2(\Sigma, \alphas)$. The domain $\Per$ must be nonnegative. Even though there are no transversality results for annoying curves, compare \cite[Section 3]{LipshitzCyl}, one can still define the Maslov index $\mu(\Per)$ according to the formula \eqref{eq:muper1}, by treating $\Per$ as in the case of usual boundary degenerations. 

We will mostly be interested in annoying $\alpha$-degenerations on the first surface $\Sigma^1$ that is part of the special connected sum. The degenerations will be taken with respect to the collection of curves $\alphas^{1+}$ that includes $\alpha^1_s$. For an annoying $\alpha$-degeneration of this form, with domain $\Per$, and any $\x \in \T_{\alpha^1}$, we have the formula \eqref{eq:muper1}. However, in  \eqref{eq:muper1} there are no vertex multiplicity contributions from the curve $\alpha^1_s$, so Equation~\eqref{eq:muper2} needs to be modified accordingly. Precisely, if we place basepoints $w_1, \dots, w_{d^1 - g^1}$ in all components of $\Sigma^1 - \alphas^{1+}$ except the two that have $\alpha^1_s$ on their boundary, we obtain
\begin {equation}
\label {eq:muannoy}
\mu(\Per) = m_1(\Per) + m_2(\Per) + 2 \sum_{i=1}^{d^1 - g^1} n_{w_i}(\Per).
\end {equation}

This has the following consequence:
\begin {lemma}
\label {lemma:peru}
Let $u$ be an annoying $\alpha$-degeneration in $(\Sigma^1, \alphas^{1+})$, with domain $\Per$. Then $\mu(\Per) \geq 0$, with equality if and only if $\Per = 0$.
\end {lemma}

We now turn to studying how cylindrical triangles in $W_{\Delta} = \Sigma \times \Delta$ relate to those in $W_\Delta^1 = \Sigma^1 \times \Delta$ and $W_\Delta^2 = \Sigma^2 \times \Delta$ when we do a special connected sum. Pick almost complex structures $J^1$ and $J^2$ on $W_\Delta^1$ and $W_\Delta^2$ and disk neighborhoods $D^1, D^2$ of $p^1$ in $\Sigma^1$, resp. $p^2$ in $\Sigma^2$. We assume that $J^1$ and $J^2$ are split near $D^1, D^2$. For $T > 0$, we form 
the connected sum
$$ \Sigma(T) = (\Sigma^1 - D^1) \# \bigl ( [-T-1, T+1] \times S^1 \bigr) \# (\Sigma^2 - D^2) $$
by inserting a long cylinder, using the identifications $\del D^1 \cong \{-T\} \times S^1$ and $\del D^2 = \{ T \} \times S^1$. We construct an almost complex structure $J(T)$ on 
$$ W_\Delta(T) = \Sigma(T) \times \Delta$$ 
by extending $J^1, J^2$ on the two sides, and using a split complex structure on the cylinder. The quantity $T$ is called the {\em neck-length}. 

By a {\em broken triangle} in a homology class $\phi$, we mean the juxtaposition of a cylindrical triangle with some cylindrical flow lines and ordinary boundary degenerations, such that the sum of all their homology classes (as a two-chain on the Heegaard surface) is $\phi$.

We then have the following convergence result:

\begin {proposition}
\label {prop:converges}
Suppose we have a special connected sum, with the notations above. Consider a homology class $\phi \in \pi_2(\x, \y, \z)$. Let $\phi^2 \in \pi_2(\x^2, \y^2, \z^2)$ be its restriction to $\Sigma^2$, and $\phi^{1+}$ the corresponding equivalence class of pairs on $\Sigma^1$. Suppose that the moduli space $\M(\phi)$ of cylindrical triangles is nonempty for a sequence of almost complex structures $J(T_i)$ with $T_i \to \infty$. Then, the moduli space of broken holomorphic triangles in the class $\phi^2$ is nonempty. Further, there exists a representative $(\phi^1, \Per)$ of the equivalence class $\phi^{1+}$ such that  the moduli space of broken holomorphic triangles in the class $ \phi^1$ is nonempty, and there exist some annoying $\alpha$-degenerations with domains that sum to $\Per$.
\end {proposition}

\begin {proof} The proof is similar to that of the second part of Theorem 5.1 in \cite{Links}, and is based on Gromov compactness, compare also \cite[Sections 7, 8, 10]{LipshitzCyl}. In the limit $T_i \to \infty$, the sequence of holomorphic triangles must have a subsequence converging to some holomorphic objects on $\Sigma^1$ and $\Sigma^2$. On $\Sigma^2$ the object is a broken triangle $u^2$. The only new twist is that when $\pi_\Sigma \circ u^2$ maps a point of the boundary of the domain to the connected sum point $p^2 \in \Sigma^2$, on the other side (i.e. on $\Sigma^1$) an annoying $\alpha$-degeneration must appear in the limit. In the end on $\Sigma^1$ we obtain a union of a broken triangle and 
some annoying $\alpha$-degenerations.
\end {proof}

There is also a gluing result:

\begin {proposition}
\label {prop:gluing}
Consider a homology class $\phi \in \pi_2(\x, \y, \z)$ on a  triple Heegard diagram obtained by special connected sum. Let $\phi^2 \in \pi_2(\x^2, \y^2, \z^2)$ be the restriction of $\phi$ to $\Sigma^2$, and $\phi^{1+}$ the  equivalence class of pairs obtained by restricting $\phi$ to $\Sigma^1$. Suppose that $\phi^{1+}$ contains a representative of the form $(\phi^1, 0)$, with $\phi^1 \in \pi_2(\x^1, \y^1, \z^1)$. Further, suppose that  
 $d_2 > g_2, \ \mu(\phi^1) = 0,\  \mu(\phi^2) = 2m$, and the domain $\D$ of $\phi$ has $m_1(\D) = m_2(\D) =m$, so that $\mu(\phi) =0$, see Lemma~\ref{lemma:mumu}. Consider the maps
 $$ \rho^1: \M(\phi^1) \longrightarrow \Sym^m(\Delta) \ \text{ and } \rho^2  : \M(\phi^2) \longrightarrow \Sym^m(\interior(\Delta) \cup e_\alpha) $$
 where
 $$ \rho^i(u) = \pi_\Delta ((\pi_\Sigma \circ u^i)^{-1}(\{ p^i \})).$$
 If the fibered product
 $$ \M(\phi^1) \times_{\Sym^m(\Delta)} \M(\phi^2) = \{ u^1 \times u^2 \in \M(\phi^1) \times \M(\phi^2) | \rho^1(u^1) = \rho^2(u^2) \}$$
 is a smooth manifold, then this fibered product can be identified with the moduli space $\M(\phi)$, for sufficiently large neck-length.
\end {proposition}

\begin {proof}
We claim that, for sufficiently large neck-length, if $u: S \to \Sigma$ is a holomorphic representative of $\phi$, no point of $\del S$ is mapped to $p^2$ under $\pi_\Sigma \circ u$. Indeed, if such points existed, in the limit when $T \to \infty$ we would get a broken triangle in a class $\psi^1$ and one (or more) annoying $\alpha$-degenerations on $\Sigma^1$, summing up to a class $\Per \in H_2(\Sigma^1, \alphas^{1+})$. We must have $\mu(\psi^1) \geq 0$ because of the existence of a holomorphic representative, $\mu(\Per) > 0$ by Lemma~\ref{lemma:peru}, and $\mu(\psi^1) + \mu(\Per) = 0$ because $(\psi^1, \Per) \sim (\phi^1, 0)$. This is a contradiction, so our claim was true. 

By the definition of the fibered product, if a holomorphic triangle $u^2 \in \M(\phi^2)$ is such that $\rho^2(u^2) \not\subset \Sym^m(\interior(\Delta))$, that triangle cannot appear in the fibered product. With these observations in mind, the rest of the proof is completely analogous to that of the third part of Theorem 5.1 in \cite{Links}. Basically, the index conditions forbid the presence of flow lines and boundary degenerations as part of broken triangles in 
the relevant moduli space. Further, the hypothesis $d_2 > g_2$ is used to exclude the presence of sphere bubbles on the $\Sigma^2$ side. One can then use the gluing arguments from \cite{LipshitzCyl}, applied to triangles.
\end {proof}

\subsection {A degeneration argument}
\label {sec:degeneration}

We now return to the setting of Section~\ref{sec:setup}. We view the triple Heegaard diagram $(\Sigma, \alphas, \betas, \gammas)$ as the special connected sum of the diagram $(\Sigma, \bar \alphas, \bar \betas, \bar \gammas)$ with the genus zero diagram $(\Sphere, \alpha_1, \beta_1, \gamma_1)$ from Figure~\ref{fig:ess}, see Example~\ref{ex:ess}.

\begin {lemma}
\label {lemma:tris}
Let $\psi \in \pi_2(a, \theta, b)$ be a homology class of triangles in $\Sphere$, with $a \in \{x, x'\}$ and $b\in \{y, y'\}$. Let $m_1, m_2, m_3, m_4$ be the local multiplicities  of the domain of $\psi$ in the regions marked as such in Figure~\ref{fig:ess}. Then:
$$ \mu(\psi) = m_1 + m_2 + m_3 + m_4. $$
\end {lemma}

\begin {proof}
The claimed equality is true when $\psi$ is the index zero triangle in $\pi_2(x, \theta, y)$. Any other class $\psi$ is related to this by the juxtaposition of a linear combination of embedded index one $\alpha$-$\beta$ or $\alpha$-$\gamma$ bigons (i.e. homology classes of flow lines) and index two disks (i.e. homology classes of boundary degenerations on $\alpha$, $\beta$ or $\gamma$). It is straightforward to check the equality for these disks and bigons.
\end {proof}

\begin {proof}[Proof of Proposition~\ref{prop:degen}] We seek to understand the moduli space of triangles $\M(\phi)$ for $\phi \in \pi_2(\x, \thetas, \y)$ with $\mu(\phi) =0$. 

Suppose $\M(\phi) \neq \emptyset$, for any sufficiently large neck-length $T$. Let $\psi=\phi^2  \in \pi_2(a, \theta, b)$ be the restriction of $\phi$ to $\Sphere$, where $a \in \{x, x'\}$ and $b\in \{y, y'\}$. Let also $\phi^{1+}$ be the equivalence class of pairs which is the restriction of $\phi$ to $(\Sigma, \bar \alphas, \bar \betas, \bar \gammas)$. Using Proposition~\ref{prop:converges}, there must be a pair $(\bar \phi, \Per) \in \phi^{1+}$ that admits holomorphic representatives. Hence, $\bar \phi$ and $\Per$ are nonnegative domains, and $\mu(\phi^{1+}) \geq 0$.

On the other hand, by Lemmas~\ref{lemma:mumu} and \ref{lemma:tris}, we have
$$ 0 = \mu(\phi) = \mu(\phi^{1+}) + \mu(\psi) - m_1 - m_2 = \mu(\phi^{1+}) + m_3 + m_4.$$

Since all the terms on the right hand side are nonnegative, we deduce that $\mu(\phi^{1+}) = m_3 = m_4=0$.
The fact that $\mu(\phi^{1+}) =0$ together with Lemma~\ref{lemma:peru} implies that there can be no annoying $\alpha$-degenerations: $\Per = 0, \ \mu( \bar \phi) = 0$ and $m_1 = m_2$. Denote by $m$ the common value $m_1 = m_2$. We are now able to apply Proposition~\ref{prop:gluing} to obtain an identification:
\begin {equation}
\label {eq:fibered}
\M(\phi) \cong \M(\bar \phi) \times_{\Sym^m(\Delta)} \M(\psi).
\end {equation}

The fact that $m_1 = m_2$ and $m_3 =m_4 = 0$ implies that $\psi$ must be a class in either $\pi_2(x, \theta, y)$ or $\pi_2(x', \theta, y')$. Without loss of generality, let us consider $\psi \in \pi_2(x, \theta, y)$. We have $\mu(\psi) = 2m$. From the proof of Proposition~\ref{prop:gluing} we know that for any $u \in \M(\psi)$, 
$$ \rho(u) = \pi_\Delta ((\pi_\Sigma \circ u)^{-1}(\{p^2\})) $$
lies in $\Sym^m(\interior(\Delta))$, that is, it does not contain any points on $e_\alpha$. Given $\p \in \Sym^m(\Delta)$, set
$$ \M(\psi, \p) = \{u \in \M(\psi) | \rho(u) = \p \}.$$

Define
$$ M(\p) = \sum_{\{\psi \in \pi_2(x, \theta, y) | m_1(\psi) = m_2(\psi) =m, m_3(\psi) = m_4(\psi) = 0 \}} \# \M(\psi, \p).$$

A Gromov compactness argument shows that $M(\p)$ is independent of $\p$, modulo $2$, compare \cite[Lemma 6.4]{Links}.  By taking the limit as $\p$ consists of $m$ distinct points, all approaching the edge $e_\beta$ of $\Delta$ with spacing at least $T $ between them, with $T \to \infty$, we obtain that the respective contributions $M(\p)$ are splicings of the index zero triangle class $\psi_0 \in \pi_2(x, \theta, y)$ and $m$ $\beta$-boundary degenerations of index two, compare \cite[Lemma 6.4]{Links}. There is a unique possible class  of $\beta$-boundary degenerations of index two with $m_3 = m_4 = 0$, namely the exterior of the curve $\beta_1$ from Figure~\ref{fig:ess}. For this class the count of pseudo-holomorphic representatives (modulo the two-dimensional group of automorphisms) is $1 (\mod 2)$, see \cite[Theorem 5.5]{Links}. Moreover, the index zero triangle class $\psi_0$ has a unique pseudo-holomorphic representative. It follows that 
$$M(\p) \equiv 1 \ (\mod 2).$$
 
A similar equality holds for the sum of contributions from classes $\psi \in \pi_2(x', \theta, y')$ with $ m_1(\psi) = m_2(\psi) =m, m_3(\psi) = m_4(\psi) = 0$. Combining these observations with \eqref{eq:fibered}, we deduce that
$$ \sum_{\phi \in \pi_2(\x, \thetas, \y)| \mu(\phi) = 0\}} \M(\phi) \equiv \sum_{\bar \phi \in \pi_2(\bar \x, \bar \thetas, \bar \y)| \mu(\bar \phi) = 0\}}  \M(\bar \phi) \ \ (\mod 2),$$
for sufficiently large neck-length. This implies that the triangle maps $F$ and $G$ are the same.
\end {proof}

\subsection {Higher polygons}
\label{sec:HigherP}

\begin{figure}
\begin{center}
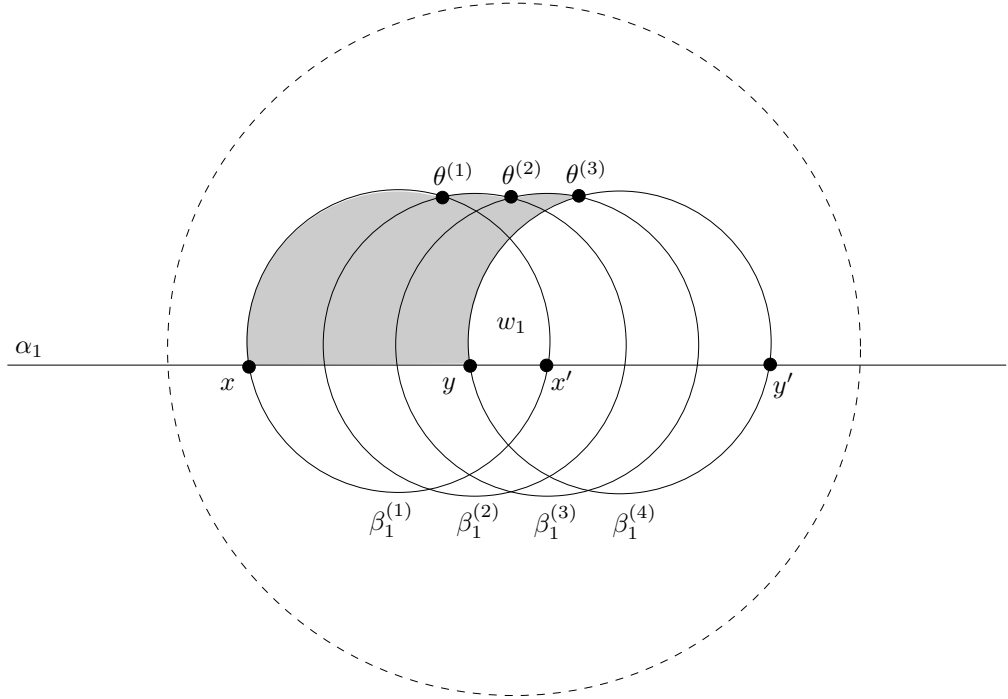
\end{center}
\caption {{\bf Several curve collections.} This is the analogue of Figure~\ref{fig:trio} for higher polygons. We show here the case $l=4$. The shaded domain is a pentagon of index zero.}
\label{fig:several}
\end{figure}

We now turn to a generalization of
Proposition~\ref{prop:degen}. Suppose we have $l \geq 2$ collections
of attaching curves $\alphas, \betas^{(1)}, \betas^{(2)}, \dots,
\betas^{(l)}$ on a multi-pointed Heegaard surface $(\Sigma,
\ws, \zs)$, such that each diagram $(\Sigma, \alphas,
\betas^{(i)}, \ws, \zs)$ is the quasi-stabilization of a diagram
$(\Sigma, \bar \alphas, \bar \betas^{(i)}, \bar \ws, \zs)$,
obtained by adding curves $\alpha_1, \beta^{(i)}_1$ and the basepoint
$w_1$. We also assume that, for every $i \neq j$, the curves
$\beta^{(i)}_1$ and $\betas^{(j)}_1$ differ by a small Hamiltonian
isotopy, and intersect each other in two points. See
Figure~\ref{fig:several}.

Let $\{x, x'\} = \alpha_1 \cap \beta_1^{(1)}$ and $\{y, y'\} = \alpha_1 \cap \beta^{(l)}_1$, with $x$ and $y$ to the left of $x'$ and $y'$, as in Figure~\ref{fig:several}. Let also $\theta^{(i)} \in \beta_1^{(i)} \cap \beta^{(i+1)}_1$ be the upper intersection point between the two curves. 

 There is a unique index zero $(l+1)$-gon class $\psi_0$ on the sphere $\Sphere$ with vertices at $x, \theta^{(1)}, \dots, \theta^{(l-1)}, y$, such that $\psi_0$ has a positive domain: see Figure~\ref{fig:several}. The moduli space of holomorphic representatives for $\psi_0$ is $(l-2)$-dimensional, corresponding to all possible lengths of the cuts at $\theta^{(1)}, \dots, \theta^{(l-2)}$.
  
\begin {lemma}
\label {lemma:conf}
 For a fixed, generic conformal structure on the domain (which is a disk with $l+1$ marked points on its boundary), the number of pseudo-holomorphic representatives of $\psi_0$ is one $(\text{mod } 2)$. 
\end {lemma}

\begin {proof}
We use induction on $l$. The case $l=3$ was treated in  \cite[proof of Theorem 4.7]{BrDCov}.
For $l \geq 4$, let $J_0$ be a generic conformal structure on the domain. Take a generic, smooth path of conformal structures $(J_t)_{t \in [0,1)}$, which starts at $J_0$ and limits (as $t \to 1$) to a degenerate conformal structure, corresponding to taking the length of one of the cuts starting at $\theta^{(l-2)}$ to infinity. Thus, in the limit $t \to 1$ the domain degenerates into the union of an $l$-sided polygon and a triangle. Let $\M_t= \M_t(\psi_0)$ be the moduli space of pseudo-holomorphic representatives of $\psi_0$ with the conformal structure $J_t$ on the domain. By Gromov compactness and generic transversality, the union 
$$ \M = \bigcup_{t \in [0,1]} \M_t$$
is a one-dimensional compact manifold with boundary $\M_0 \cup \M_1$. 
By the inductive hypothesis, the cardinality of $\M_1$ is odd; hence, the same must be true for $\M_0$.  
\end {proof}

For $i=1, \dots, l-1$, pick $\bar \thetas^{(i)} \in \T_{\bar \beta^{(i)}} \cap \T_{\bar \beta^{(i+1)}}$ and set 
$$ \thetas^{(i)}=
\bar \thetas^{(i)} \cup \{\theta^{(i)}\} \in  \T_{ \beta^{(i)}} \cap \T_{ \beta^{(i+1)}}.$$ 

Just as in Section~\ref{sec:triangles}, there is a map
$$ F: \Am(\Ta, \T_{ \beta^{(1)}}, \s) \to \Am(\Ta, \T_{ \beta^{(l)}}, \s), \ \ F(\x) = f(\x \otimes \thetas^{(1)} \otimes \dots \otimes \thetas^{(l-1)}),$$
this time given by counting pseudo-holomorphic $(l+1)$-gons of index $2-l$, see Section~\ref{sec:polygon}.

For  any $\x \in \Ta \cap \T_{ \beta^{(1)}}$ and $\y \in \Ta \cap  \T_{ \beta^{(l)}}$, we can eliminate their intersections with $\alpha_1$ to obtain generators $\bar \x   \in \T_{\bar \alpha} \cap \T_{\bar \beta^{(1)}}$ and $\y \in \T_{\bar \alpha} \cap  \T_{\bar \beta^{(l)}}$. We can define a map $\bar F$ by counting $(l+1)$-gons in the destabilized diagram with $l-1$ fixed vertices at $\bar \thetas^{(1)}, \dots, \bar \thetas^{(l-1)}$. If $\bar F(\bar \x) = \sum_{\bar \y \in \T_{\bar \alpha}\cap  \T_{\bar  \beta^{(l)}}} n_{\bar \x, \bar \y} \bar \y$, we set
$$G( \x) = \sum_{\bar \y \in \T_{\bar \alpha}\cap  \T_{\bar  \beta^{(l)}}} n_{\bar \x, \bar \y} (\bar \y \times (\x \cap \alpha_1)).$$ 
 
\begin {proposition}
\label {prop:polydegen}
For $\Sigma, \alphas,  \betas^{(1)}, \dots, \betas^{(l)}, \thetas^{(1)}, \dots, \thetas^{(l-1)}$ as above, and suitable almost complex structures on the symmetric products, the maps $F$ and $G$ coincide.
\end {proposition}

\begin {proof}
The arguments are completely analogous to the ones for $l=2$. We insert a long neck along the dashed curve from Figure~\ref{fig:several}. We thus view $\Sigma$ as the special connected sum of the sphere $\Sphere $ and the quasi-destabilized diagram, except now each has $l+1$ collections of attaching curves. Let $\phi \in \pi_2(\x,  \thetas^{(1)}, \dots, \thetas^{(l-1)}, \y)$ be homology class of $(l+1)$-gons, of index $2-l$, that admits pseudo-holomorphic representatives. In the limit when the neck-length $T \to \infty$, the class $\phi$ splits into homology classes of broken $(l+1)$-gons $\psi$ on the sphere $\Sphere$, $\bar \phi$ on the quasi-destabilized diagram, and a class $\Per$ of  annoying $\alpha$-degenerations. Lemma~\ref{lemma:mumu}, Equation~\eqref{eq:muannoy} and Lemma~\ref{lemma:tris} still hold true, and therefore we have
\begin {equation}
\label {eq:pit}
2-l = \mu(\phi) = \mu(\bar \phi) + \mu(\Per) + m_3 + m_4 \geq (2-l) + 0 + 0 + 0. 
\end {equation}

Hence, $\Per = 0$ (so there are no annoying $\alpha$-degenerations), $m_3 = m_4 = 0$, and we end up with a fibered product description of $\M(\phi)$ analogous to \eqref{eq:fibered}. On the $\Sphere$ side, we can use a limiting process to ensure that the holomorphic representatives of $\psi$ are  splicings of an index zero $(l+1)$-gon and several $\beta^{(1)}$-boundary degenerations. One possibility for the index zero $(l+1)$-gon is that it lies in the class  $\psi_0$ analyzed in Lemma~\ref{lemma:conf}. In the fibered product description the conformal structure of the domain of a pseudo-holomorphic $(l+1)$-gon in  $\psi_0$ is specified by the $(l+1)$-gon on the quasi-destabilized diagram. By Lemma~\ref{lemma:conf}, the number of the holomorphic representatives of $\psi_0$ with this constraint is $1$ (mod $2$). A similar discussion applies to the index zero $(l+1)$-gon class with vertices at 
$x', \theta^{(1)}, \dots, \theta^{(l-1)}, y'$. The fibered product description then implies the identification of the two maps $F$ and $G$.
\end {proof}

The reader may wonder what happens in the case $l=1$, that is, how are flow lines in a diagram $\bar \Hyper$ related to flow lines in its quasi-stabilization $\Hyper$. This question is more difficult, because if we try to degenerate along the special connected sum neck, we can no longer avoid the presence of annoying $\alpha$-degenerations. Indeed, in  \eqref{eq:pit} the inequality $\mu(\bar \phi) \geq 2-l = 1$ does not hold true, due to the existence of index zero flow lines on $\bar \Hyper$ (namely, trivial ones). This leaves open the possibility that $\mu(\Per) = 1$. 

Nevertheless, we make the following:
\begin {conjecture}
\label {conj:flows}
Let $\bar \Hyper= (\Sigma, \bar \alphas, \bar \betas,\bar \ws, \zs)$ be a Heegaard diagram, and $\Hyper = (\Sigma, \alphas, \betas, \ws,\zs)$ be its quasi-stabilization, as in Section~\ref{sec:setup}. Suppose $w_2$ is the second basepoint (apart from $w_1$) in the component of $\Sigma - \bar \alphas$ that contains the curve $\alpha_1$. Suppose the variables corresponding to $w_1$ and $w_2$ are $U_1$, resp. $U_2$. Then, for suitable almost complex structures, there is an identification between the Floer complex $\Am(\Hyper, \s)$ and the mapping cone complex 
 $$ \Am(\bar \Hyper, \s)[[U_1]] \xrightarrow{U_1 - U_2} \Am(\bar \Hyper, \s)[[U_1]].$$
 \end {conjecture} 
  
Note that this is a generalization of Proposition 6.5 in \cite{Links}, which dealt with ordinary index zero/three stabilizations. The difficulty in proving Conjecture~\ref{conj:flows} is the lack of available transversality and gluing results for annoying curves. In Section~\ref{sec:handleslid}, we will give a proof of the conjecture for  a particular class of Heegaard diagrams, using rather ad-hoc arguments.

\subsection{Stabilizations}
\label {sec:stabil}

We now specialize to the case of ordinary (free) index zero/three stabilizations. 

Let $(\Sigma,\alphas,\{\betas^{(i)}\}_{i=1}^{l},\ws,\zs)$ be as in
the previous section, and
$$ F: \Am(\Ta, \T_{ \beta^{(1)}}, \s) \to \Am(\Ta, \T_{ \beta^{(l)}}, \s), \ \ F(\x) = f(\x \otimes \thetas^{(1)} \otimes \dots \otimes \thetas^{(l-1)})$$
the corresponding polygon map.
Assume further that the curves $\beta^{(i)}_1$ is a small Hamiltonian translate of $\alpha_1$, intersecting it in two points. Thus,
$(\Sigma,\alphas,\betas^{(i)},\ws,\zs)$ is obtained from
$(\Sigma,{\bar\alphas},{\bar\betas}^{(i)},{\bar\ws},\zs)$
by an index zero/three stabilization, and
$$ {\bar F}: \Am(\T_{\bar\alpha}, \T_{ {\bar\beta}^{(1)}}, \s) \to \Am(\T_{\bar\alpha}, \T_{ {\bar\beta}^{(l)}}, \s), \ \ 
{\bar F}({\bar\x}) = f({\bar\x} \otimes {\bar\thetas}^{(1)} \otimes
\dots \otimes {\bar\thetas}^{(l-1)}).$$ 

We keep the notation from Conjecture~\ref{conj:flows}, with basepoints $w_1,w_2$ and variables $U_1, U_2$. Associate one $U$ variable to each $w$ basepoint, and let $\Ring$ be the power series ring generated over $\Field$ by all $U$ variables except $U_1$.

From the proofs of Theorem~\ref{thm:LinkInvariance} and Lemma~\ref{lem:ap1}, we know that there are destabilization maps
$$\Destab^{i} \colon \Am(\Sigma,{\alphas},{\betas}^{(i)},\s) \longrightarrow \Am(\Sigma,{\bar\alphas},{\bar\betas}^{(i)},\s),$$
defined by
$$ U_1^n({\bar \x} \times x') \mapsto 0, \ \ \
U_1^n({\bar \x} \times x)\mapsto U_2^n \bar \x,$$
and these maps are equivalences over the ring $\Ring$.

\begin{proposition}
\label {prop:StabPolygon}
The equivalences induced by destabilization commute with the polygon maps, that is,
 $$  \Destab^{l} \circ  F={\bar F}\circ \Destab^1.$$
\end{proposition}

\begin {proof} Use Proposition~\ref{prop:polydegen}.
\end {proof}

Versions of Proposition~\ref{prop:StabPolygon} are true for ordinary
Heegaard Floer complexes, as well. In that
case, we can consider an admissible multi-diagram
$(\Sigma,\alphas,\{\betas^{(i)}\}_{i=1}^{l},\ws)$, 
where again the curves $\beta^{(i)}_1$ are small isotopic translates of $\alpha_1$, so that we can form the index
zero/three destabilizations, and corresponding maps 
\begin{align*}
 { F}_o: \CFm(\T_{\alpha}, \T_{{\beta}^{(1)}},\ws) &\to \CFm(\Ta, \T_{ {\beta}^{(l)}},\ws), \\
 {\bar F}_o: \CFm(\T_{\bar \alpha}, \T_{{\bar \beta}^{(1)}}, \bar \ws) &\to \CFm(\Ta, \T_{ {\bar \beta}^{(l)}}, \bar \ws).
\end{align*}
(Note that the Heegaard diagrams appearing here can represent
arbitrary three-manifolds.)

We have analogous destabilization maps $\Destab_o$.

\begin{proposition}
\label {prop:StabPolygonH}
In the case of ordinary Heegaard Floer complexes,
the equivalences induced by destabilization commute with the polygon maps,
that is, $$  \Destab_o^{l} \circ { F}_o={\bar F}_o\circ \Destab^1_o.$$
\end{proposition}

\begin {proof} 
  This is analogous to the proof of  Proposition~\ref{prop:StabPolygon}.
\end {proof}

\section{Index zero/three link stabilizations}
\label{sec:03}

In this section we study one of the Heegaard moves that appeared in Section~\ref{sec:hmoves}, namely index zero/three link stabilization. The behavior of holomorphic disks under this move was discussed in \cite{MOS} and \cite{Zemke}. We will review those results, and then (in Section~\ref{sec:hpol}) we extend them to the case of more general holomorphic polygons. Furthermore, in Section~\ref{sec:final}, we will study what happens to holomorphic polygons under a type of strong equivalence that relates different kinds of stabilizations.

The analytical input for this section comes from the work of Zemke in \cite{Zemke}.

\subsection{Algebraic preliminaries}
In Section~\ref{sec:Inv} we established Lemma~\ref{lem:ap1}, which said that a free complex $C$ over $R[[U_1]]$ is chain homotopy equivalent to the mapping cone
\begin{equation}
\label{eq:conefree}
\mathit{Cone}\bigl(C_-[[U_2]] \xrightarrow{U_1-U_2} C_+[[U_2]] \bigr).
\end{equation}
That cone appeared naturally in the context of (both free and link) index zero/three stabilizations; see the proof of Theorem~\ref{thm:LinkInvariance}.

To better understand link index zero/three stabilizations, we will also need the following variant of Lemma~\ref{lem:ap1}.

\begin{lemma}
\label{lem:ap2}
Let $R$ be an $\ff$-algebra, and let $C$ be a free complex over $R[[U_1]]$. Let $C^{U_1 \to U_2}_{\pm}$ be two copies of the free complex over $R[[U_2]]$ obtained from $C$ by replacing the variable $U_1$ with $U_2$. Then, the complexes $C$ and
$$ C' = \mathit{Cone}\bigl(C^{U_1 \to U_2}_+[[U_1]] \xrightarrow{U_1-U_2} C^{U_1 \to U_2}_-[[U_1]] \bigr)$$
are chain homotopy equivalent over $R[[U_1]]$.
\end{lemma}

\begin{proof}
Let us fix a set of free generators $\S$ for $C$. We write the differential on $C$ as
$$\del \x = \sum_{\y \in \S} c(\x, \y) U_1^{n(\x, \y)} \y,$$
where $c(\x, \y) \in \{0,1\}$ and $n(\x, \y)$ are nonnegative integers. Further, given a generator $\x \in \S$, we denote by $\x_+$ and $\x_-$ the corresponding generators in $C_+^{U_1 \to U_2}$ and $C_-^{U_1 \to U_2}$.

We construct chain maps between $C$ and $C'$ in both directions, such that they are module maps over $R[[U_1]]$ and they are homotopy inverses. 

We let $\rho: C' \to C$ be given by
$$ \rho (U_2^m \x_-) = U_1^m \x, \ \ \rho(U_2^m \x_+) = 0.$$
In the opposite direction, we define $\iota: C \to C'$ by
$$ \iota(\x) = \x_- + \sum_{\y \in \S} c(\x, \y)  \frac{U_1^{n(\x, \y)} - U_2^{n(\x, \y)}}{U_1 - U_2} \y_+.$$

One can check that $\rho$ and $\iota$ commute with the differentials, and $\rho \circ \iota$ is the identity. Moreover, if we define 
$$ H: C' \to C', \ \ H(U_2^n \x_+)=0, \ H(U_2^n \x_-) = \frac{U_1^n-U_2^n}{U_1-U_2} \x_+,$$
then we find that $\iota \circ \rho$ is chain homotopic to the identity, via the homotopy $H$.
\end{proof}

\subsection{Holomorphic disks}
\label{sec:hd03}
\begin{figure}
\begin{center}
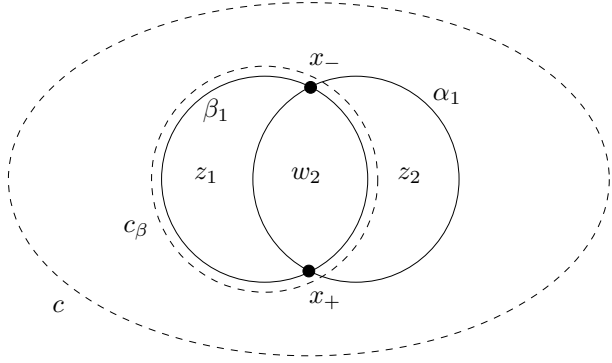
\end{center}
\caption {An index zero/three link stabilization.}
\label{fig:stabilize0}
\end{figure}

As defined in Section~\ref{sec:hmoves}, an index zero/three link stabilizations is a local move on a Heegaard diagram, in the neighborhood of an existing basepoint of type $z$. The move introduces a new alpha curve, a new beta curve, and a basepoint of each type; see Figure~\ref{fig:Stab03}. For convenience, we include Figure~\ref{fig:Stab03} again here, as Figure~\ref{fig:stabilize0},  with the curves $\alpha_n$ and $\beta_n$ relabeled as $\alpha_1$ and $\beta_1$, and the basepoints $z, w', z'$ relabeled as $z_1, w_2, z_2$. (Not shown is, for example, the basepoint $w_1$, which is separated from $z_1$ by a number of beta curves.) 

We will denote the initial diagram by $\bHyper=(\bar \Sigma, \bar \alphas, \bar \betas, \bar \ws, \bar \zs)$, and the stabilized diagram by $\Hyper=(\Sigma, \alphas, \betas, \ws, \zs)$. We let $L_1 \subseteq L$ be the link component on which the stabilization is performed.

We are interested in the relation between the holomorphic disks in $\Hyper$ to those in $\bHyper$. If we ignore the $z$ basepoints, our picture is just that of a free index zero/three stabilization, and an analysis of the holomorphic disks was done in \cite[Section 6]{Links}; cf. also \cite[Section 2]{MOS}. The analysis is based on viewing $\Hyper$ as the connected sum of $\bHyper$ (with $z_1$ deleted) and a sphere $\Sphere$ containing the new curves (i.e., containing what is shown in Figure~\ref{fig:stabilize0}), and then degenerating the Heegaard surface along the curve $c$ in Figure~\ref{fig:stabilize0}, by taking the length of the connected sum neck to infinity. Further, the connected sum point $p$ is then taken to be close to one of the circles $\alpha_1$ or $\beta_1$ in Figure~\ref{fig:Stab03}. 

The analysis in \cite{Links} and \cite{MOS} did not distinguish between the disks that cross $z_1$ and those that cross $z_2$. A more complete characterization of the holomorphic disks was given in  \cite[Lemmas 14.3 and 14.4]{Zemke}. This is based on a different degeneration. We still stretch the neck along the dashed curve $c$ in Figure~\ref{fig:stabilize0}, but then, instead of taking the connected sum point close to a curve, we stretch along the dashed curve $c_{\beta}$. Note that stretching along $c_{\beta}$ is very similar to the degeneration used for quasi-stabilizations in Section~\ref{sec:triangles}. 

We state the results from  \cite[Lemmas 14.3 and 14.4]{Zemke} here, with slightly different conventions: $z_1$ and $z_2$ playing the roles of $w'$ and $w$ in \cite[Lemma 14.3]{Zemke}, and $\alpha_1$ and $\beta_1$ are switched.

\begin{proposition}[Zemke \cite{Zemke}]
\label{prop:ZemkeDisks}
For generic almost complex structures associated to sufficiently large neck stretching along $c$ and $c_{\beta}$, we have that, with regard to the counts of holomorphic disks $( \mathit{modulo }\ 2)$ in $\Hyper$,
\begin{itemize}
\item The only disks of Maslov index one in a class $\phi \in \pi_2(\x \times x_-, \y \times x_-)$ from $\Hyper$ appear when $n_{z_1}(\phi)=n_{w_2}(\phi)=0$, in which case they are in one-to-one correspondence with the holomorphic disks of Maslov index one in the corresponding class $\bar \phi \in \pi_2(\x, \y)$ in $\bHyper$. Here, $\phi$ is obtained from $\bar \phi$ by taking connected sum (along $c$) with some copies of the complement of the disk bounded by $\beta_1$ in Figure~\ref{fig:stabilize0}, so that $n_{z_2}(\phi) = n_{z_1}(\bar \phi)$;
\item The same description applies to the holomorphic disks of Maslov index one in classes $\phi \in \pi_2(\x \times x_+, \y \times x_+)$;
\item The only holomorphic disks of Maslov index one in a class $\phi \in \pi_2(\x \times x_-, \y \times x_+)$ from $\Hyper$ appear when $\x = \y$ and the domain of $\phi$ is the bigon in Figure~\ref{fig:stabilize0} containing $w_2$; moreover, in that case there is only one such disk;
\item The only holomorphic disks of Maslov index one in a class $\phi \in \pi_2(\x \times x_+, \y \times x_-)$ from $\Hyper$ appear when $\x = \y$ and $n_{z_1}(\phi)=n_{w_2}(\phi)=n_{z_2}(\phi)=0$; moreover, for each $\x$, the total count of such disks (over all $\phi$) is $1$ $( \text{mod  }2)$, with the disk being in a class $\phi$ that has $n_{w_1}(\phi)=0$ and goes over no other basepoints.
\end{itemize}
\end{proposition}

Proposition~\ref{prop:ZemkeDisks} allows us to relate the various Heegaard Floer complexes associated to $\Hyper$ and $\bHyper$. For example, if we ignore the $z$ basepoints, we are in the setting of a free index zero/three stabilization. Consider the complex $C = \CFm( \T_{\bar \alpha}, \T_{\bar \beta},\bar \ws).$ Then, the corresponding complex $\CFm(\T_{\alpha}, \T_{\beta}, \ws)$ in the stabilized diagram is identified with
\begin{equation}
\label{eq:firstcone}
 C_-[[U_2]] \xrightarrow{U_2-U_1} C_+[[U_2]],
 \end{equation}
where we denoted by $C_{\pm}$ the copies of $C$ that are obtained by including the basepoint $x_{\pm}$ in $\Sphere$. This cone appeared in \eqref{eq:conefree}, and explains why Lemma~\ref{lem:ap1} is relevant to our situation. In particular, Lemma~\ref{lem:ap1} implies the following result (which we have already seen as part of the proof of Theorem~\ref{thm:LinkInvariance}).
 \begin{corollary}
 \label{cor:rhonormal}
In the setting of Figure~\ref{fig:stabilize0}, the destabilization map $$ \rho: \CFm(\T_{\alpha}, \T_{\beta}, \ws) \to \CFm(\T_{\bar \alpha}, \T_{\bar \beta}, \bar \ws)$$
given by
$$\rho( U_2^m (\x \times x_+)) = U_1^m \x, \ \ \rho(U_2^m (\x \times x_-))= 0$$
is a chain homotopy equivalence over $\bar \Ring$ (the ring with one $U$ variable for each $w$ basepoint in the diagram $\bar \Hyper$).

Furthermore, the equivalence $\rho$ is filtered with respect to the Alexander filtrations, and therefore gives rise to equivalences
$$ \rho: \Am(\Hyper, \s) \to \Am(\bar \Hyper, \s).$$
\end{corollary}

Note that $\Am(\Hyper, \s)$ does not admit a mapping cone description similar to \eqref{eq:firstcone}. However, such descriptions do apply to, for example, generalized link Floer complexes that involve only the $w$ basepoints on the component $L_1$, but may involve both $w$'s and $z$'s on the other components.
\medskip

\begin{figure}
\begin{center}
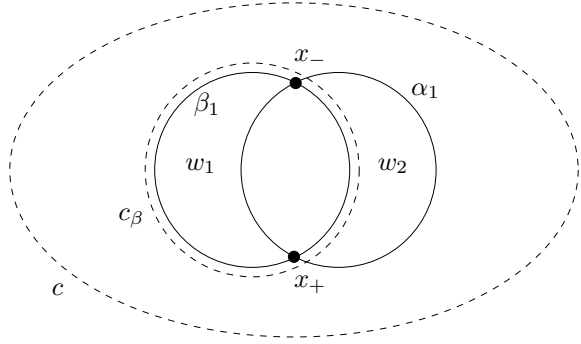
\end{center}
\caption {A variant of the free index zero/three stabilization.}
\label{fig:variantfree}
\end{figure}

Also relevant to link stabilizations is Lemma~\ref{lem:ap2}. Indeed, let us ignore the $w$ basepoints and look at the complexes constructed with $z$ basepoints. We re-label the $z$ basepoints on $L_1$ as $w$'s, so that our notation to be consistent with that from Section~\ref{sec:hed}. Then, Figure~\ref{fig:stabilize0} becomes Figure~\ref{fig:variantfree}. The move in Figure~\ref{fig:variantfree}  consists in introducing two new curves and a new basepoint $w_2$ in the neighborhood of $w_1$. This is a variant of a free index zero/three stabilization: it differs from that by a strong equivalence (cf. Section~\ref{sec:final} below), and in fact it could replace the ordinary free index zero/three stabilizations in the list of Heegaard moves in Section~\ref{sec:hmoves}.

If we now denote $C= \CFm(\T_{\bar \alpha}, \T_{\bar \beta}, \bar \ws),$ 
then Proposition~\ref{prop:ZemkeDisks} implies that the corresponding complex $\CFm(\T_{\alpha}, \T_{\beta}, \ws)$ in $\Hyper$ is
\begin{equation}
\label{eq:secondcone}
 C^{U_1 \to U_2}_+[[U_1]] \xrightarrow{U_1-U_2} C^{U_1 \to U_2}_-[[U_1]].
 \end{equation}

 Lemma~\ref{lem:ap2} has the following:
 \begin{corollary}
 \label{cor:rhovariant}
For the move shown in Figure~\ref{fig:variantfree}, the destabilization map $$ \rho: \CFm(\T_{\alpha}, \T_{\beta}, \ws) \to \CFm(\T_{\bar \alpha}, \T_{\bar \beta}, \bar \ws),$$
\begin{equation}
\label{eq:rhova}
\rho( U_2^m (\x \times x_-)) = U_1^m \x, \ \ \rho(U_2^m (\x \times x_+))= 0
\end{equation}
is a chain homotopy equivalence over $\bar \Ring$.
\end{corollary}

\begin{remark} The result of Corollary~\ref{cor:rhovariant} also holds if instead of $\CFm$ we work with generalized  link Floer complexes that involve only the $z$ basepoints on the component $L_1$ (now re-labeled as $w$'s), but may involve both $w$'s and $z$'s on the other components.
\end{remark}

\subsection{Changes in the almost complex structures}
\label{sec:changeJ}
In the previous subsection we worked with a generic family (call it $J_{\beta})$ of almost complex structures associated to neck stretching along the curves $c$ and $c_{\beta}$. We could just as well stretch along $c$ and a curve $c_{\alpha}$ around $\alpha_1$; let us denote the corresponding family of almost complex structures by $J_{\alpha}$. 

With $J_{\alpha}$, the analogue of Proposition~\ref{prop:ZemkeDisks} holds, with the difference that the disks in a class $\phi \in \pi_2(\x \times x_-, \y \times x_-)$ or $ \pi_2(\x \times x_+, \y \times x_+)$ have $n_{z_2}(\phi)=0$ instead of $n_{z_1}(\phi)=0$; and the classes $\phi$ are obtained from $\bar \phi$ by taking connected sum with some copies of the complement of the disk bounded by $\alpha_1$. 
Note that the description of $\CFm(\T_{\alpha}, \T_{\beta}, \ws)$ as the cone \eqref{eq:firstcone} is still true, and so is Corollary~\ref{cor:rhonormal}.

\begin{figure}
\begin{center}
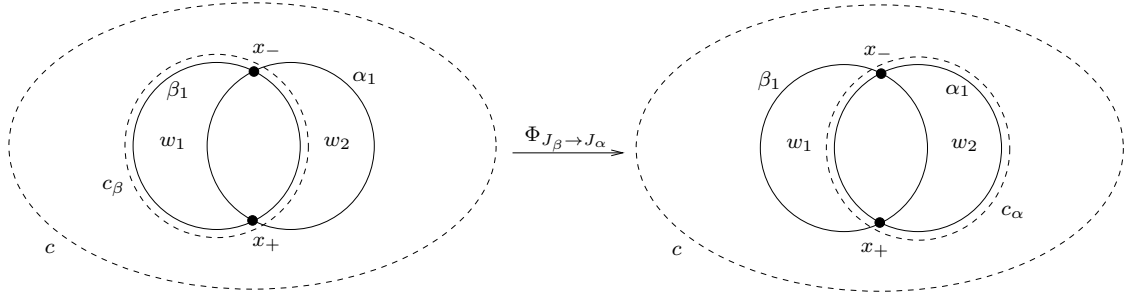
\end{center}
\caption {A change of almost complex structures from Figure~\ref{fig:variantfree}.}
\label{fig:changeJ}
\end{figure}

On the other hand, once we ignore the $w$ basepoints and relabel $z$'s as $w$'s, the description of $\CFm(\T_{\alpha}, \T_{\beta}, \ws)$ from \eqref{eq:secondcone} changes. In the situation depicted on the right hand side of Figure~\ref{fig:changeJ}, with $J_{\alpha}$, the complex $\CFm(\T_{\alpha}, \T_{\beta}, \ws)$ is simply the cone
\begin{equation}
\label{eq:thirdcone}
 C_+[[U_2]] \xrightarrow{U_2-U_1} C_-[[U_2]].
\end{equation}
Corollary~\ref{cor:rhovariant} still holds, with the $\rho$ map defined by the same formula \eqref{eq:rhova}. 

One can interpolate between $J_{\beta}$ and $J_{\alpha}$ by a family of almost complex structures, all with sufficiently large connected sum neck along $c$. This induces a map on Floer complexes
\begin{equation}
\label{eq:Jba}
 \Phi_{J_{\beta} \to J_{\alpha}} : \CFm(\T_{\alpha}, \T_{\beta}, \ws, J_{\beta}) \to \CFm(\T_{\alpha}, \T_{\beta}, \ws, J_{\alpha}).
\end{equation}
See Figure~\ref{fig:changeJ}. The change of almost complex structure map was studied in \cite[Lemma 14.5]{Zemke}. We state the result below with our conventions. 

Here, and later in the paper, we will write maps between mapping cones as a $2 \times 2$ matrices; each cone will be viewed (as a vector space) as the direct sum of its domain and target, in this order. 

\begin{proposition}[Zemke \cite{Zemke}]
\label{prop:changeJ}
For sufficiently stretched almost complex structures, the map \eqref{eq:Jba} can be written as
$$ \Phi_{J_{\beta} \to J_{\alpha}} = \begin{pmatrix}
\id & * \\
0 & \id
\end{pmatrix}$$
or, in more expanded form,
$$
\xymatrix{
C^{U_1 \to U_2}_+[[U_1]] \ar[r]^-{\id} \ar[d]_{U_1 - U_2} &  C_+[[U_2]] \ar[d]^{U_1 - U_2} \\
C^{U_1 \to U_2}_-[[U_1]]  \ar[r]^-{\id} \ar[ur]^{*} & C_-[[U_2]].
}
$$
\end{proposition}

\begin{remark}
The upper right entry $\ast$ in $ \Phi_{J_{\beta} \to J_{\alpha}} $ can be computed explicitly; see \cite[Remark 14.6]{Zemke}. However, we do not need to know it for our purposes.
\end{remark}

\begin{corollary}
\label{cor:Jba}
The map \eqref{eq:Jba} is related to the projections $\rho$ of the form \eqref{eq:rhova} as follows:
$$ \rho \circ  \Phi_{J_{\beta} \to J_{\alpha}} = \rho.$$
\end{corollary}

\begin{proof}
This follows immediately from Proposition~\ref{prop:changeJ}. Indeed, recall that the maps $\rho$ take the domains of the cones to zero, and act on the targets by taking both $U_1$ and $U_2$ to $U_1$.
\end{proof}

We also need to consider the change of almost complex structures map for a free index zero/three stabilization done in the neighborhood of a basepoint $w_1$; see Figure~\ref{fig:changeJfree}. This will be useful to us in Section~\ref{sec:final} below. To be in agreement with the notation there, we will denote the unstabilized diagram by $(\bar \Sigma, \bar \alphas, \bar \gammas, \bar \ws)$, and the stabilized one by $(\Sigma, \alphas, \gammas, \ws)$, with $\alphas = \bar \alphas \cup \{\alpha_1\}$, $\gammas =\bar \gammas \cup \{\gamma_1\}$, and $\ws=\bar \ws \cup \{w_2\}$. 
We consider almost complex structures $J_{\alpha}$ stretched along $c$ and $c_{\alpha}$, and $J_{\gamma}$ stretched along $c$ and $c_{\gamma}$. We then interpolate between them using a family of almost complex structures sufficiently stretched along $c$. This induces a map
\begin{equation}
\label{eq:Jag}
 \Phi_{J_{\alpha} \to J_{\gamma}}: \CFm(\T_{\alpha}, \T_{\gamma}, \ws, J_{\alpha}) \to \CFm(\T_{\alpha}, \T_{\gamma}, \ws, J_{\gamma}).
 \end{equation}

\begin{figure}
\begin{center}
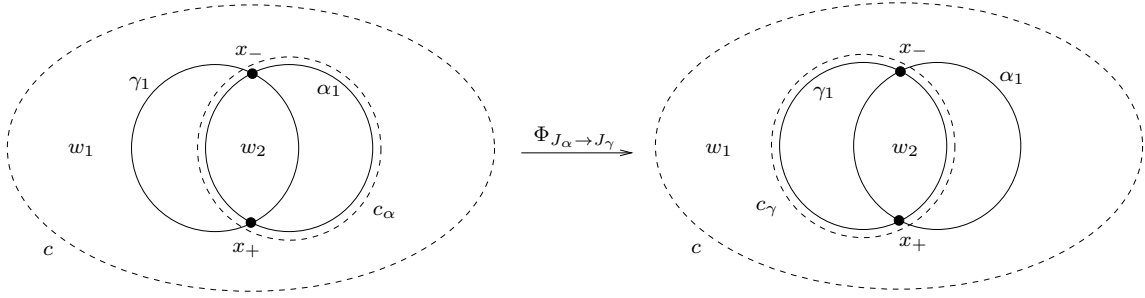
\end{center}
\caption {A change of almost complex structures for free index zero/three stabilizations.}
\label{fig:changeJfree}
\end{figure}

Note that, if we denote $C= \CFm(\T_{\bar \alpha}, \T_{\bar \gamma}, \bar \ws)$, then both the domain and the target of $ \Phi_{J_{\alpha} \to J_{\gamma}}$ can be identified with the cone
\begin{equation}
\label{eq:fourthcone}
C_-[[U_2]] \xrightarrow{U_2-U_1} C_+[[U_2]].
\end{equation}

The following result comes from Proposition 14.22 in \cite{Zemke}.

\begin{proposition}[Zemke \cite{Zemke}]
\label{prop:changeJfree}
For sufficiently stretched almost complex structures, the map \eqref{eq:Jag} can be written as
$$ \Phi_{J_{\alpha} \to J_{\gamma}} = \begin{pmatrix}
\id & 0 \\
* & \id
\end{pmatrix}$$
or, in more expanded form,
$$
\xymatrix{
C_-[[U_2]] \ar[r]^-{\id} \ar[d]_{U_1 - U_2} \ar[dr]^{*} &  C_-[[U_2]] \ar[d]^{U_1 - U_2} \\
C_+[[U_2]]  \ar[r]^-{\id}  & C_+[[U_2]].
}
$$
\end{proposition}

Once again, the lower right term $\ast$ is computed explicitly in \cite[Proposition 14.22]{Zemke}, but we do not need to know that formula here.

We are interested in the relation of the map \eqref{eq:Jag} with the projections $\rho$. 

\begin{corollary}
\label{cor:Jag}
There is a chain homotopy 
$$ \rho \circ  \Phi_{J_{\alpha} \to J_{\gamma}} \simeq \rho.$$
\end{corollary}

\begin{proof}
Note that, since now the lower right term is nonzero, we do not have $ \rho \circ  \Phi_{J_{\alpha} \to J_{\gamma}} = \rho$ on the nose, as in Corollary~\ref{cor:Jba}. On the other hand, we can look at the homotopy inverses to $\rho$. From the proof of Lemma~\ref{lem:ap1}, such homotopy inverses are given by the formula
$$ \iota: \CFm(\T_{\bar \alpha}, \T_{\bar \gamma}, \bar \ws) \to \CFm(\T_{\alpha}, \T_{\gamma}, \ws), \ \ \iota(\x)= \x \times x_+.$$

Proposition~\ref{prop:changeJfree} implies that $\iota = \Phi_{J_{\alpha} \to J_{\gamma}} \circ \iota$, on the nose. Since $\iota \circ \rho \simeq \id$ and $\rho \circ \iota = \id$, we have 
$$ \rho \circ  \Phi_{J_{\alpha} \to J_{\gamma}} \simeq  \rho \circ  \Phi_{J_{\alpha} \to J_{\gamma}} \circ \iota \circ \rho=  \rho \circ \iota \circ   \rho =  \rho,$$
as desired.
\end{proof}

\subsection{Holomorphic polygons}
\label{sec:hpol}
Our next goal is to study the behavior of holomorphic polygons under an index zero/three link stabilization. The discussion will be modelled on Section~\ref{sec:triangles}, where we did a similar study for quasi-stabilizations.

\begin{figure}
\begin{center}
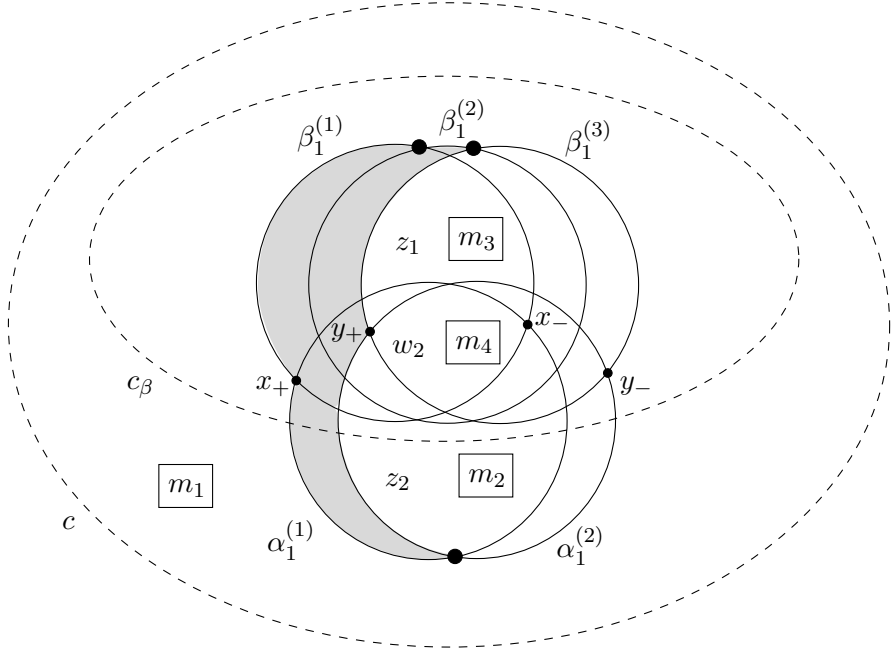
\end{center}
\caption {General polygons under an index zero/three stabilization. The black dots represent the $\theta$ intersection points. The shaded region is the domain of a pentagon class $\psi_0$.}
\label{fig:Polygon03}
\end{figure}

Consider a Heegaard multi-diagram
$$ (\bar\Sigma, \bar \alphas^{(1)}, \dots, \bar \alphas^{(k)}, \bar \betas^{(1)}, \dots, \bar \betas^{(l)}, \bar \ws, \bar \zs),$$
and let 
$$(\Sigma, \alphas^{(1)}, \dots, \alphas^{(k)}, \betas^{(1)}, \dots, \betas^{(l)}, \ws, \zs)$$ 
be its index zero/three link stabilization, in the sense that each pair $(\alphas^{(i)}, \betas^{(j)})$ is related to $(\bar \alphas^{(i)}, \bar \betas^{(j)})$ by an index zero/three link stabilization, introducing two new basepoints $w_2$ and $z_2$ near $z_1$, and two new curves $\alpha^{(i)}_1$, $\beta^{(j)}_1$. Furthermore, we want the new curves $\alpha^{(i)}_1$ to be Hamiltonian translates of one another, intersecting at two points each; and the same for the new curves $\beta^{(j)}_1$. See Figure~\ref{fig:Polygon03} for an example, with $k=2$ and $l=3$.

Mark the top degree intersection points $\theta_\alpha^{(i)} \in \alpha^{(i)}_1 \cap \alpha^{(i+1)}_1$, for $i=1, \dots, k-1$, and $\theta_\beta^{(j)} \in \beta^{(j)}_1 \cap \beta^{(j+1)}_1$, for $j=1, \dots, l-1$. Pick $\thetas_\alpha^{(i)} \in \T_{\alpha^{(i)}} \cap \T_{\alpha^{(i+1)}}$ containing $\theta_\alpha^{(i)}$, and pick $\thetas_\beta^{(j)} \in \T_{\beta^{(j)}} \cap \T_{\beta^{(j+1)}}$ containing $\theta_\beta^{(j)}$. Let $\bar \thetas_\alpha^{(i)}=\thetas_\alpha^{(i)} - \{ \theta_\alpha^{(i)} \}$ and $\bar \thetas_\beta^{(j)}= \thetas_\beta^{(j)} - \{ \theta_\beta^{(j)} \}$ be the corresponding intersection points in the destabilized diagram.

We want to study holomorphic $(k+l)$-gons with boundaries on $\alphas^{(k)}, \dots, \alphas^{(1)}, \betas^{(1)}, \dots, \betas^{(l)}$, in this clockwise order, and with vertices at $\thetas_\alpha^{(i)}, \thetas_\beta^{(j)}$ as well as arbitrary $\x \in  \T_{\alpha^{(1)}} \cap \T_{\beta^{(1)}}$ and $\y \in  \T_{\alpha^{(k)}} \cap \T_{\beta^{(l)}}$. A class $\phi$ of such holomorphic $(k+l)$-gons can be viewed as the connected sum of a class $\bar \phi$ on the destabilized diagram, and a class $\psi$ of polygons on the sphere $\Sphere$. Furthermore, given $\bar \phi$ and intersection points $r_1 \in \alpha^{(1)}_1 \cap \beta^{(1)}_1, r_2 \in \alpha^{(k)}_1 \cap \beta^{(l)}_1$, we define 
$\Phi(\bar \phi, r_1, r_2)$ as the set of classes $\phi=(\bar \phi, \psi)$ satisfying $n_{w_2}(\phi)=0$, and with $\psi$ going from $r_1$ to $r_2$. This set is nonempty only when $(r_1, r_2)=(x_+, y_+)$ or $(r_1, r_2)=(x_-, y_-)$, in which case it consists in splicing of a standard polygon $\psi_0$ on $\Sphere$ (an example is the shaded region in Figure~\ref{fig:Polygon03}) and some beta boundary degenerations.

\begin{proposition}
\label{prop:PolyDestab03}
Choose generic almost complex structures $J_{\beta}$ associated to sufficiently large connected sum neck along $c$, and also sufficiently stretched along a curve $c_{\beta}$ enclosing the $\beta^{(j)}_1$ curves, as in Figure~\ref{fig:Polygon03}. Then:
\begin{enumerate}[(a)]
\item
When the components of $(\x, \y)$ on $\Sphere$ are the pair $(x_+, y_-)$, the counts of holomorphic $(k+l)$-gons in classes 
$$\phi \in \pi_2(\thetas_\alpha^{(k-1)}, \dots, \thetas_\alpha^{(1)}, \x, \thetas_\beta^{(1)}, \dots, \thetas_\beta^{(l-1)}, \y)$$ 
are always zero;
\item Suppose $\x = \bar \x \times r_1$ and $\y = \bar \y \times r_2$, with $(r_1, r_2)=(x_+, y_+)$ or $(x_-, y_-)$. Then, the counts of holomorphic $(k+l)$-gons in classes 
$$\phi \in \pi_2(\thetas_\alpha^{(k-1)}, \dots, \thetas_\alpha^{(1)}, \x, \thetas_\beta^{(1)}, \dots, \thetas_\beta^{(l-1)}, \y)$$ 
can be nonzero only when $n_{w_2}(\phi)=0$. Further, for 
 $$\bar \phi \in \pi_2(\bar \thetas_\alpha^{(k-1)}, \dots,\bar \thetas_\alpha^{(1)}, \bar \x, \bar \thetas_\beta^{(1)}, \dots, \bar \thetas_\beta^{(l-1)}, \bar \y),$$  if we define $\Phi(\bar \phi, r_1, r_2)$ as the set of classes $\phi=(\bar \phi, \psi)$ satisfying $n_{w_2}(\phi)=0$, and with $\psi$ going from $r_1$ to $r_2$, then
\begin{equation}
\label{eq:emfi}
 \# \M(\bar \phi) = \sum_{\phi \in \Phi(\bar \phi, r_1, r_2)} \# \M(\phi) \ \ (\text{mod } 2).
 \end{equation}
\end{enumerate}
\end{proposition}

\begin{proof}
We use similar arguments to those in the proof of Proposition~\ref{prop:polydegen}.

In the limit as we stretch the connected sum neck along $c$, holomorphic polygons become splicings of (possibly broken) holomorphic polygons on the two sides. The homology classes of polygons split also, as $\phi = (\bar \phi, \psi)$, where $\bar \phi$ is on $\bar \Sigma$ and $\psi$ on $\Sphere$.

Unlike for quasi-stabilizations, now the curve $\alpha_1$ does not go through the point $p$, and hence cannot serve as boundary for the domains on $\bar \Sigma$. Hence, $\bar \phi$ is now an ordinary class of polygons (rather than a pair of a polygon and an $\alpha$ boundary degeneration), so the analysis becomes somewhat simpler. 

Consider a few domain multiplicities for $\phi$, shown in Figure~\ref{fig:Polygon03}:
$$m_1=n_p(\phi), \ m_2 = n_{z_2}(\phi), \ m_3=n_{z_1}(\phi), \ m_4=n_{w_2}(\phi).$$

The analogue of Lemma~\ref{lemma:mumu} gives the index of the class $\phi$:
$$ \mu(\phi) = \mu(\bar \phi) + \mu(\psi) - 2m_1.$$

There is also an analogue of Lemma~\ref{lemma:peru}, which gives the index of the component $\psi$:
$$\mu(\psi)= (3-k-l) + m_1+m_2+m_3+m_4.$$

Hence,
\begin{equation}
\label{eq:muofphi}
 \mu(\phi)= \mu(\bar \phi) + m_2+m_3+m_4 - m_1.
 \end{equation}

Further, by the relations for a domain to be acceptable (cf. Lemma~\ref{lemma:acc3}), we have
$$ m_3 - m_4 = m_1 - m_2 + \delta,$$
with 
$$\delta= \begin{cases}
1 &\text{if } (r_1, r_2) = (x_+, y_-), \\
0 & \text{if } (r_1,r_2) = (x_+, y_+) \text{ or } (x_-, y_-),\\
-1 &\text{if } (r_1, r_2) = (x_-, y_+).
\end{cases}$$

From here and \eqref{eq:muofphi} we get
\begin{equation}
\label{eq:muofphinew}
\mu(\phi) = \mu(\bar \phi) + 2m_4 + \delta.
\end{equation}

We are interested in classes $\phi$ with $\mu(\phi)=3-k-l$, so that we can count rigid holomorphic $(k+l)$-gons. If such a class $\phi$ contains holomorphic representatives, we must have $\mu(\bar \phi) \geq 3-k-l$. The relation~\eqref{eq:muofphinew} implies that $2m_4 + \delta \leq 0.$ Since $m_4 \geq 0$, this disallows the case $\delta =1$; that is, we cannot have $(r_1, r_2) = (x_+, y_-)$. Part (a) is proved.

Part (b) deals with the case $\delta =0$. Then, $\mu(\bar \phi) \geq 3-k-l = \mu(\phi)$ and \eqref{eq:muofphinew} imply that
$$\mu(\bar \phi) =0, \ \ m_4=0.$$ 

This shows that $\phi \in \Phi(\bar \phi, r_1, r_2)$. Observe that the classes in  $\Phi(\bar \phi, r_1, r_2)$ are splicings of a standard polygon class $\psi_0$ on $\Sphere$ (an example is the shaded region in Figure~\ref{fig:Polygon03}) and some disk classes (namely, exteriors of the curves $\alpha_1^{(i)}$ or $\beta_1^{(j)}$ in Figure~\ref{fig:Polygon03}).

We have convergence and gluing results entirely similar to those in the proof of Proposition~\ref{prop:polydegen}; compare Propositions~\ref{prop:converges} and \ref{prop:gluing}. From there we obtain the desired relation \eqref{eq:emfi} between the holomorphic polygons in the class $\bar \phi$ and those in the possible classes $\phi$.
\end{proof}

\begin{remark}
If we had a single set of alpha curves (and several betas), then the result of Proposition~\ref{prop:PolyDestab03} would follow more directly from the analysis in Section~\ref{sec:HigherP}; see also \cite[Lemma 14.25]{Zemke}. In that case, we would actually get a stronger conclusion; for example, that the holomorphic polygon counts are zero for $(r_1, r_2)=(x_-, y_+)$ as well.
\end{remark}

\begin{remark}
In the proof of Proposition~\ref{prop:PolyDestab03} we actually did not use that the almost complex structures are stretched along $c_{\beta}$; only along $c$. However, we need to stretch along $c_{\beta}$ in order to identify holomorphic disks, and thus express the domain and target of  polygon maps in terms of mapping cones---for Propositions~\ref{prop:PolyDestab03a} and ~\ref{prop:PolyDestab03b} below.
\end{remark}

We present below two consequences of Proposition~\ref{prop:PolyDestab03}. They are similar in spirit to Propositions~\ref{prop:StabPolygonH} and \ref{prop:StabPolygon}. 

First, we consider the polygon maps on ordinary Floer complexes $\CFm$, and study their behavior with respect to the variant of index zero/three stabilization shown in Figure~\ref{fig:variantfree}. The relevant picture is obtained from Figure~\ref{fig:Polygon03} by deleting $w_2$ and relabeling $z_1$ and $z_2$ as $w_1$ and $w_2$. In the stabilized diagram, we have a map
$$F: \CFm(\T_{\alpha^{(1)}}, \T_{ \beta^{(1)}}) \to \CFm(\T_{\alpha^{(k)}}, \T_{ \beta^{(l)}}), \ \ F(\x) = f(\thetas_\alpha^{(k-1)}\otimes \dots \otimes \thetas_\alpha^{(1)} \otimes \x \otimes \thetas_\beta^{(1)} \otimes \dots \otimes \thetas_\beta^{(l-1)}),$$
whereas in the original (destabilized) diagram, we have
$$ \bar F:  \CFm(\T_{\bar \alpha^{(1)}}, \T_{\bar \beta^{(1)}}) \to \CFm(\T_{\bar \alpha^{(k)}}, \T_{ \bar \beta^{(l)}}), \ \ F(\bar \x) = f(\bar \thetas_\alpha^{(k-1)} \otimes \dots \otimes \bar \thetas_\alpha^{(1)} \otimes \bar \x \otimes \bar \thetas_\beta^{(1)} \otimes \dots \otimes \bar \thetas_\beta^{(l-1)}).$$

\begin{proposition}
\label{prop:PolyDestab03b}
In the setting of Figure~\ref{fig:variantfree}, the polygon maps induced on $\CFm$ commute (on the nose) with the projections $\rho$ from Corollary~\ref{cor:rhovariant}; that is, we have
$$ \rho \circ F =\bar F \circ \rho.$$ 
\end{proposition}

\begin{proof}
The map $\rho$ takes generators containing $x_+$ or $y_+$ to zero, and otherwise sends both $U_1$ and $U_2$ to $U_1$. In view of this, the result follows readily from Proposition~\ref{prop:PolyDestab03}: When $\bar \y$ appears with some coefficient $U_1^n$ in the polygon map (on the destabilized diagram) applied to some $\bar \x$, then $\bar \y \times y_-$ appears with a coefficient $U_1^{n_1} U_2^{n_2}$ in the polygon map applied to $\bar \x \times x_-$, such that $n_1 + n_2=n$. This implies commutation with $\rho$.
\end{proof}

\begin{remark} Proposition~\ref{prop:PolyDestab03b} easily extends to maps between generalized link Floer complexes, where these complexes may use both types of basepoints on the link components that are not involved in the stabilization move.
\end{remark}

Secondly, we can look at an index zero/three stabilization as in Figure~\ref{fig:Polygon03}, and consider the polygon maps induced on generalized link Floer complexes:
\begin{equation}
\label{eq:FA}
F: \Am(\T_{\alpha^{(1)}}, \T_{ \beta^{(1)}}, \s) \to \Am(\T_{\alpha^{(k)}}, \T_{ \beta^{(l)}}, \s), \ \ F(\x) = f(\thetas_\alpha^{(k-1)}\otimes \dots \otimes \thetas_\alpha^{(1)} \otimes \x \otimes \thetas_\beta^{(1)} \otimes \dots \otimes \thetas_\beta^{(l-1)})
\end{equation}
and
\begin{equation}
\label{eq:bFA}
 \bar F:  \Am(\T_{\bar \alpha^{(1)}}, \T_{\bar \beta^{(1)}}, \s) \to \Am(\T_{\bar \alpha^{(k)}}, \T_{ \bar \beta^{(l)}}, \s), \ \ F(\bar \x) = f(\bar \thetas_\alpha^{(k-1)} \otimes \dots \otimes \bar \thetas_\alpha^{(1)} \otimes \bar \x \otimes \bar \thetas_\beta^{(1)} \otimes \dots \otimes \bar \thetas_\beta^{(l-1)}).
\end{equation}

We are interested in whether these maps commute with the projections $\rho$. In Corollary~\ref{cor:rhonormal}, the maps $\rho$ on $\Am(\Hyper, \s)$ are induced from those on $\CFm$, using the Alexander filtration. Further, homotopy inverses $\iota$ to the $\rho$ maps  on $\CFm$ can be found from the proof of Lemma~\ref{lem:ap1}; they are given by the formulas
$$ \iota(\x)= \x \times x_+, \ \ \iota(\y) = \y \times y_+.$$
The same formulas must give homotopy inverses for the filtered versions of $\rho$, on $\Am(\Hyper, \s)$. 

Proposition~\ref{prop:PolyDestab03} now shows that the polygon maps commute with $\iota$, on the nose:
\begin{equation}
\label{eq:iotaf}
F \circ \iota = \iota \circ \bar F.
\end{equation}
In the case when $F$ and $\bar F$ are triangle maps (i.e., $k=1, l=2$ or $k=2, l=1$), the same argument as in the proof of Corollary~\ref{cor:Jag} shows that these triangle maps commute with the homotopy inverses to $\iota$, which are the maps $\rho$:
\begin{equation}
\label{eq:rhoff}
 \rho \circ F \simeq \bar F \circ \rho.
 \end{equation}

For $k + l \geq 4$, however, the polygon maps $F$ and $\bar F$ are not chain maps, so a chain homotopy as in \eqref{eq:rhoff} would not make sense. Nevertheless, we can state a similar result by considering hypercubes of chain complexes, as in Section~\ref{sec:hyperco}. Suppose that 
$$C = (C^{\eps}, D^{\eps})_{\eps \in \E_n}$$ is a hypercube, where all the complexes $C^{\eps}$ are generalized link Floer complexes of the form $\Am(\T_{\alpha{(\eps)}}, \T_{\beta{(\eps)}}, \s)$, for some collections of curves $\alphas(\eps), \betas(\eps)$ depending on $\eps$, and the maps $D^{\eps}$ are polygon maps of the form \eqref{eq:FA} above, with $k+l = \|\eps\|+2.$ (Such hypercubes will appear throughout this paper; cf. Section~\ref{sec:hyperfloer} below.) Furthermore, suppose that all $\alphas(\eps)$ and $\betas(\eps)$ are obtained from collections $\bar \alphas(\eps)$ and $\bar \betas(\eps)$ by index zero/three stabilizations in the same spot. We have a hypercube 
$$\bar C= (\bar C^{\eps}, \bar D^{\eps})_{\eps \in \E_n}$$
in the destabilized diagram, where the maps $\bar D^{\eps}$ are of the form \eqref{eq:bFA}.

The analogue of \eqref{eq:rhoff} in this general situation is the following.

\begin{proposition}
\label{prop:PolyDestab03a}
Consider two hypercubes of chain complexes $\bar C$ and $C$ as above, related by an index zero/three link stabilization, with almost complex structures chosen as in Proposition~\ref{prop:PolyDestab03}. Then, there exists a chain map (as in Definition~\ref{def:chmap})
$$ \Psi: C \to \bar C$$
such that its components $\Psi^{\zero}_{\eps}: C^{\eps} \to \bar C^{\eps}$ are the projections $\rho$ from Corollary~\ref{cor:rhonormal}.
\end{proposition}

\begin{proof}
 For the components of $\Psi$ that increase the value of $\|\eps\|$ by one, we choose the chain homotopies in \eqref{eq:rhoff}; then, Equation~\eqref{eq:rhoff} says that the condition for $\Psi$ to be a chain map, Equation~\eqref{eq:DF}, is satisfied along two-dimensional faces. 

The higher components of $\Psi$ are constructed inductively in the dimension of the faces. This is done using Equation \eqref{eq:iotaf} and arguments similar to those in the proof of Corollary~\ref{cor:Jag}, but applied to chain maps between hypercubes. 
\end{proof}

\subsection{A strong equivalence}
\label{sec:final}

\begin{figure}
\begin{center}
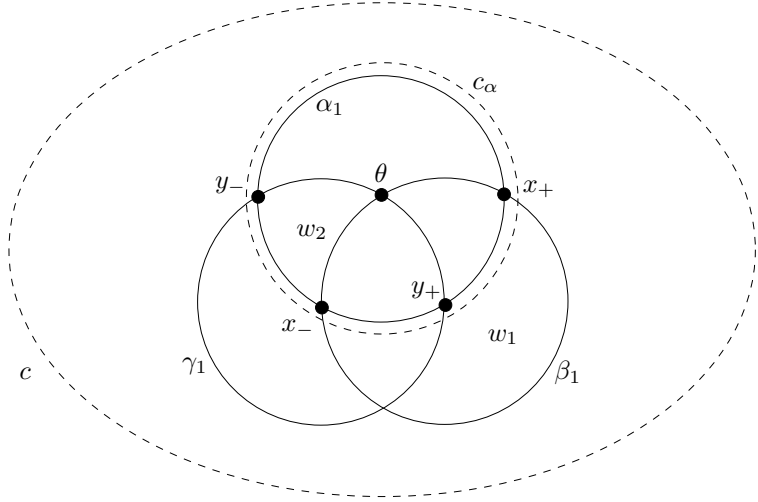
\end{center}
\caption {A strong equivalence between different kinds of stabilizations.}
\label{fig:StabTriple0}
\end{figure}

We now turn to discussing a move that will appear naturally in the context of index zero/three link stabilizations for complete systems of hyperboxes, in Section~\ref{sec:moves}. This move is pictured in Figure~\ref{fig:StabTriple0}, and consists in replacing the curve $\beta_1$ by $\gamma_1$ in that diagram, as well as replacing all the other beta curves (not shown) with small isotopic translates of themselves, intersecting them in two points. We denote by $\Hyper=(\Sigma, \alphas, \betas, \ws)$ the initial diagram, by $\Hyper'=(\Sigma, \alphas, \gammas, \ws)$ the final one, and by $\bar \Hyper=(\bar \Sigma, \bar\alphas, \bar\betas, \bar\ws)$ and  $\bar \Hyper'=(\bar \Sigma, \bar\alphas, \bar\gammas, \bar\ws)$ the destabilized diagrams, with a single $w$ in place of Figure~\ref{fig:StabTriple0}. Note that $\bar \Hyper$ and $\bar \Hyper'$ are strongly equivalent, and in fact surface isotopic, in the sense of Definition~\ref{def:ab}. As usual, we view $\Sigma$ as the connected sum of $\bar \Sigma$ and a sphere $\Sphere$ that contains Figure~\ref{fig:StabTriple0}.

Observe that $\Hyper$ is exactly the variant of free index zero/three stabilization pictured in Figure~\ref{fig:variantfree}, whereas $\Hyper'$  is an ordinary free index zero/three stabilization of $\bar \Hyper'$. Our move from $\Hyper$ to $\Hyper'$ is a strong equivalence, cf. Definition~\ref{def:ab} (a), and it can be viewed as handlesliding $\beta_1$ over several other beta curves---those on the boundary of the component of $\bar \Sigma \setminus \betas$ that contains $w_1$.

For each curve $\beta_i$ not shown in the figure, let $\gamma_i$ be its translate, and let $\theta_i \in \beta_i \cap \gamma_i$ be the intersection point that gives higher homological grading. We denote by $\thetas \in \T_{\beta} \cap \T_{\gamma}$ the generator formed by all $\theta_i$ together with the point $\theta$ in Figure~\ref{fig:StabTriple0}. Let $\bar \thetas  = \thetas -\{\theta\} \in \T_{\bar \beta} \cap \T_{\bar \gamma}$ be the corresponding generator in the destabilized diagram.

To compute Floer complexes, we will use almost complex structures for sufficiently large neck stretching along the curves $c$ and $c_{\alpha}$ in Figure~\ref{fig:StabTriple0}. With these almost complex structures, recall from Section~\ref{sec:changeJ} that the Heegaard Floer complexes $\CFm(\Hyper, J_{\alpha})$ and $\CFm(\Hyper', J_{\alpha})$ can be described as the mapping cones \eqref{eq:thirdcone} and \eqref{eq:fourthcone}, respectively. Thus, they are chain homotopy equivalent to $\CFm(\bar \Hyper)$ resp. $\CFm(\bar \Hyper')$, via the projections $\rho$.

At the level of Heegaard Floer complexes, the strong equivalence from $\bar \Hyper$ to $\bar \Hyper'$ induces a map
$$ \bar F: \CFm(\bar \Hyper) \to \CFm (\bar \Hyper')$$
given by counting holomorphic triangles with one vertex at $\bar \thetas$. 

Similarly, the strong equivalence from $\Hyper$ to $\Hyper'$ induces a map
$$ F_{\alpha} : \CFm(\Hyper, J_{\alpha}) \to \CFm(\Hyper', J_{\alpha})$$
given by counting holomorphic triangles with one vertex at $\thetas$. 

The map $F_{\alpha}$ is a version of the transition map computed in \cite[Proposition 14.8]{Zemke}.  Translated into our setting, his result reads:

\begin{proposition}[Zemke \cite{Zemke}]
\label{prop:ZemkeTransition}
For generic almost complex structures $J_{\alpha}$ associated to a sufficiently large necks along the curves $c$ and $c_{\alpha}$, the map $F_{\alpha}$ is given by the $2 \times 2$ matrix
\begin{equation}
\label{eq:fmatrix}
 F_{\alpha} = \begin{pmatrix} \bar F & 0 \\ 0 &  \bar F\end{pmatrix},
 \end{equation}
or, in a more expanded form,
$$
\xymatrix{
C_+[[U_2]] \ar[r]^-{\bar F} \ar[d]_{U_1 - U_2} &  C'_-[[U_2]] \ar[d]^{U_1 - U_2} \\
C_-[[U_2]]  \ar[r]^-{\bar F}  & C'_+[[U_2]],
}$$
where $C=\CFm(\bar \Hyper)$ and $C'=\CFm(\bar \Hyper')$.
\end{proposition}

\begin{remark}
Since the $\bar \betas$ and $\bar \gammas$ curves are obtained from each other by small isotopies, the triangle map $\bar F$ is chain homotopic to the nearest point map; see \cite[proof of Theorem 6.6]{OzsvathStipsicz}. On homology, we can think of the map induced by $\bar F$  as the identity.
\end{remark}

\begin{corollary}
\label{cor:PreRhoMove}
Under the hypotheses of Proposition~\ref{prop:ZemkeTransition}, we have a commutative diagram
$$\xymatrix{
\CFm(\Hyper, J_{\alpha}) \ar[r]^{F_{\alpha}}  \ar[d]^{\rho} & \CFm(\Hyper', J_{\alpha}) \ar[d]^{\rho} \\
\CFm(\bar \Hyper) \ar[r]^{\bar F} & \CFm(\bar \Hyper').
}$$
\end{corollary}

\begin{proof}
This is an immediate consequence of Proposition~\ref{prop:ZemkeTransition}.
\end{proof}

While in Proposition~\ref{prop:ZemkeTransition} we worked with almost complex structures $J_{\alpha}$ stretched along $c$ and $c_{\alpha}$, in this paper we will mostly work in the setting of Section~\ref{sec:hd03}; that is, with almost complex structures $J_{\beta}$, stretched along $c$ and $c_{\beta}$ (or $J_{\gamma}$, stretched along $c$ and $c_{\gamma})$. Recall from \eqref{eq:secondcone} that $\CFm(\Hyper, J_{\beta})$ can be identified with the cone
$$C^{U_1 \to U_2}_+[[U_1]] \xrightarrow{U_1-U_2} C^{U_1 \to U_2}_-[[U_1]].$$
Further, $\CFm(\bar \Hyper, J_{\gamma})$ can be identified with the cone from \eqref{eq:firstcone}, namely
$$ C_-[[U_2]] \xrightarrow{U_1-U_2} C_+[[U_2]].$$

We can go from $\CFm(\Hyper, J_{\beta})$ to $\CFm(\Hyper', J_{\gamma})$ by first changing the almost complex structures from $J_{\beta}$ to $J_{\alpha}$ as in Proposition~\ref{prop:changeJ}, then applying the triangle map $F_{\alpha}$ from Proposition~\ref{prop:ZemkeTransition}, and then changing the almost complex structure from $J_{\alpha}$ to $J_{\gamma}$ as in Proposition~\ref{prop:changeJfree}; that is, we define
$$ F: \CFm(\Hyper, J_{\beta}) \to \CFm(\Hyper', J_{\gamma}), $$
$$ F=  \Phi_{J_{\alpha} \to J_{\gamma}} \circ F_{\alpha} \circ \Phi_{J_{\beta} \to J_{\alpha}}.$$

Our main interest is the relation of $F$ with the destabilization maps $\rho$ from Corollaries~\ref{cor:rhonormal} and \ref{cor:rhovariant}. 

\begin{corollary}
\label{cor:RhoMove}
The following diagram commutes up to a chain homotopy (represented by the diagonal map): 
$$\xymatrix{
\CFm(\Hyper, J_{\alpha}) \ar[r]^F  \ar[d]^{\rho} \ar[dr] & \CFm(\Hyper', J_{\gamma}) \ar[d]^{\rho} \\
\CFm(\bar \Hyper) \ar[r]^{\bar F} & \CFm(\bar \Hyper').
}$$
\end{corollary}

\begin{proof}
Put together the results of Corollaries~\ref{cor:Jba}, \ref{cor:PreRhoMove} and \ref{cor:Jag}.
\end{proof}

We now turn to discussing how higher polygon maps interact with the move in Figure~\ref{fig:StabTriple0}. 

Specifically, in Figure~\ref{fig:StabTriple0}, we will introduce several isotopic translates of each curve. Suppose we have a Heegaard multi-diagram
$$ (\bar\Sigma, \bar \alphas^{(1)}, \dots, \bar \alphas^{(k)}, \bar \betas^{(1)}, \dots, \bar \betas^{(l)}, \bar \gammas^{(1)}, \dots, \bar \gammas^{(m)}, \bar \ws, \zs),$$
and let 
$$(\Sigma, \alphas^{(1)}, \dots, \alphas^{(k)}, \betas^{(1)}, \dots, \betas^{(l)}, \gammas^{(1)}, \dots, \gammas^{(m)}, \ws, \zs)$$ 
be obtained from it by adding a new basepoint $w_2$ and new curves 
$$\alpha^{(i)}_1, i=1, \dots, k; \ \ \beta^{(i)}_1, i=1, \dots, l; \ \  \gamma^{(i)}_1, i=1, \dots, m$$ that are translates of $\alpha_1$, $\beta_1$ and $\gamma_1$ from Figure~\ref{fig:StabTriple0}. Furthermore, we want the new curves $\alpha^{(i)}_1$ to intersect one another at two points each, and similarly for the new beta and gamma curves. 

We mark the top degree intersection points $\theta_\alpha^{(i)} \in \alpha^{(i)}_1 \cap \alpha^{(i+1)}_1$, for $i=1, \dots, k-1$, and similarly $\theta_\beta^{(i)}$ and $\theta_\gamma^{(i)}$. Pick $\thetas_\alpha^{(i)} \in \T_{\alpha^{(i)}} \cap \T_{\alpha^{(i+1)}}$ containing $\theta_\alpha^{(i)}$, and let $\bar \thetas_\alpha^{(i)}=\thetas_\alpha^{(i)} - \{ \theta_\alpha^{(i)} \}$. We pick  $\thetas_\beta^{(i)} $ containing $\theta_\beta^{(i)}$ and $\thetas_\gamma^{(i)}$ containing $\theta_\gamma^{(i)}$, and define $\bar \thetas_\beta^{(i)}$ and $\bar \thetas_\gamma^{(i)}$ similarly. 

Furthermore, we assume that the curves $\bar \betas^{(l)}$ and $\bar \gammas^{(1)}$ are small isotopic translates of each other, intersecting in two points. We let $\bar \thetas$ be the top intersection point in $\T_{\bar \beta^{(l)}} \cap \T_{\bar \gamma^{(1)}}$, and then set $\thetas = \bar \thetas \times \theta$, where $\theta \in \beta^{(l)}_1\cap \gamma^{(1)}_1$ is as in Figure~\ref{fig:StabTriple0}.

We now count holomorphic $(k+l+m)$-gons with boundaries on $\alphas^{(k)}, \dots, \alphas^{(1)}, \betas^{(1)}, \dots, \betas^{(l)}, \gammas^{(1)}$, $\dots, \gammas^{(m)}$, in this clockwise order, and with vertices at $\thetas_\alpha^{(i)}, \thetas_\beta^{(i)}, \thetas, \thetas_\gamma^{(i)}$ as well as arbitrary $\x \in  \T_{\alpha^{(1)}} \cap \T_{\beta^{(1)}}$ and $\y \in  \T_{\alpha^{(k)}} \cap \T_{\gamma^{(m)}}$. If we use almost complex structures $J_{\alpha}$ stretched along $c$ and $c_{\alpha}$ (where $c_{\alpha}$ encloses the curves $\alpha^{(i)}_1$), the holomorphic polygon  counts give rise to a map
$$ G_{\alpha}: \CFm(\T_{\alpha^{(1)}}, \T_{\beta^{(1)}}, J_{\alpha}) \to \CFm(\T_{\alpha^{(k)}}, \T_{\gamma^{(m)}}, J_{\alpha}).$$
There is a corresponding map in the destabilized diagram
$$\bar G: \CFm(\T_{\bar \alpha^{(1)}}, \T_{\bar \beta^{(1)}}) \to \CFm(\T_{\bar \alpha^{(k)}}, \T_{\bar \gamma^{(m)}}).$$

We can also consider almost complex structures $J_{\beta}$ stretched along $c$ and $c_{\beta}$ (where $c_{\beta}$ encloses $\beta^{(1)}_1$), and $J_{\gamma}$  stretched along $c$ and $c_{\gamma}$ (where $c_{\gamma}$ encloses $\gamma^{(m)}_1$). By pre- and post-composing with change of almost complex structure maps, we can define
$$ G : \CFm(\T_{\alpha^{(1)}}, \T_{\beta^{(1)}}, J_{\beta}) \to \CFm(\T_{\alpha^{(k)}}, \T_{\gamma^{(m)}}, J_{\gamma}),$$
$$G =  \Phi_{J_{\alpha} \to J_{\gamma}} \circ G_{\alpha} \circ \Phi_{J_{\beta} \to J_{\alpha}}.
$$

The analysis done for holomorphic triangles in \cite[Proposition 14.8]{Zemke} extends to the case of higher polygons (see also the proof of Proposition~\ref{prop:PolyDestab03}). We obtain the following result. 
\begin{proposition}
\label{prop:PolyDestabMove}
The higher polygon maps $G_{\alpha}$ commute with the projections $\rho$:
$$ \rho \circ G_{\alpha} = \bar G \circ \rho.$$
\end{proposition}

Observe that Proposition~\ref{prop:PolyDestabMove} implies an analogous result for the maps $G$ instead of $G_{\alpha}$, but with the commutation only up to homotopy, and assuming that the maps $G$ are chain maps:
$$\rho \circ G \simeq \bar G \circ \rho.$$
Of course, in general, the polygon maps $G$ are not chain maps, but rather fit into hypercubes of chain complexes similar to those in the setting of Proposition~\ref{prop:PolyDestab03a}.
The ``commutation up to homotopy'' can then be phrased in terms of the existence of a chain map between hypercubes, which has the maps $\rho$ along its edges.

\begin{remark}
The discussion in this section was for Heegaard Floer complexes, using only $w$ basepoints. However, if there are some $z$ basepoints (away from the sphere $\Sphere$, so that $w_1$ and $w_2$ are still free), the same results apply to generalized link Floer complexes.
\end{remark}

\section {Hyperboxes of Heegaard diagrams}
\label {sec:hyperHeegaard}

In this section we define the notion of a complete system of hyperboxes for a link. As advertised in the introduction, this notion is the basic input for the surgery theorem.

Let $\betas = (\beta_1, \dots, \beta_{g+k-1})$ and $\betas'= (\beta'_1, \dots, \beta'_{g+k-1})$ be two collections of curves on $(\Sigma, \ws, \zs)$. We will need the following terminology:

\begin {definition}
\label {def:approx}
Suppose that for any $i$, the curve $\beta'_i$ is obtained from $\beta_i$ by an isotopy, such that $\beta_i$ and $\beta'_i$ intersect each other transversely, in exactly two points, and do not intersect any of the other curves in the diagrams. If this is the case, we write $\betas \approx \betas'$, and we say that $\betas'$ {\em approximates} $\betas$. Replacing $\betas$ by $\betas'$ is called an {\em approximation}.
 \end {definition}

Fix a multi-pointed, generic, admissible Heegaard diagram
$(\Sigma,\alphas,\betas,\ws,\zs)$
for a link $L$, and $\s\in\bH(L)$.
If $\betas'$ approximates $\betas$ sufficiently closely, then there is
a nearest-point map from $\Am(\Ta,\Tb,\s)$ to
$\Am(\Ta,\Tbp,\s)$ taking each intersection point to the
corresponding nearest intersection point. 

\begin{lemma}
  \label{lem:NearestPoints}
  If $\betas'$ approximates $\betas$ sufficiently closely, then 
  the nearest-point map is an isomorphism of chain complexes.
\end{lemma}

\begin{proof}
  When the approximation is sufficiently small, the nearest point map,
  which is clearly an isomorphism of modules, coincides with a 
  continuation map obtained from varying the almost complex structure. The latter is always a chain map, see~\cite[Section 7.3]{HolDisk}.
\end{proof}

\subsection {$\beta$-Hyperboxes}
\label {sec:betah}
Let $\betas$ and $\betas'$ be collections of curves on $(\Sigma, \ws, \zs)$ that are strongly equivalent. Note that the Heegaard diagram $(\Sigma, \betas, \betas', \ws, \zs)$ is link-minimal, and represents an unlink inside the connected sum of several copies of $S^1 \times S^2$. If the pair $(\betas, \betas')$ is generic and admissible, we can define a generalized Floer chain complex $\Am(\Tb, \Tbp, \zero)$ as in Section~\ref{sec:chains}; see in particular Remark~\ref{rem:gen}. (Here, $\zero$ is the zero vector.) We will use the alternative description of that complex from Section~\ref{sec:alternative}, and thus denote it by $\Chain^-(\Tb, \Tbp, \zero)$. The resulting Floer homology $H_*(\Chain^-(\Tb, \Tbp, \zero))$ is the generalized Heegaard Floer homology of an unlink inside the connected sum of several copies of $S^1 \times S^2$, hence it equals the homology of a torus; see the completion of the proof of Theorem~\ref{thm:LinkInvariance} at the end of Section~\ref{sec:polygon}. As such, there is a well-defined maximal degree element $\theta_{\beta, \beta'} \in H_*(\Chain^-(\Tb, \Tbp, \zero))$. Observe that, if the surface $\Sigma$ has genus $g$, then this maximal degree is $\mu = g/2$. In the particular case when $\betas \approx \betas'$, there is also a canonical cycle (intersection point) $\Theta^\can_{\beta, \beta'}$ representing  $\theta_{\beta, \beta'}$. 

Recall the notation from Section~\ref{sec:hyperv}: in particular, pick $\dd =(d_1, \dots, d_n) \in \N^n, n \geq 0$, and consider the set of multi-indices $\E(\dd)$. We say that two multi-indices $\eps, \eps' \in \E(\dd)$ with $\eps \leq \eps'$ are {\em neighbors} if $\eps' - \eps \in \E_n$. Note that, in the definition of  a hyperbox of chain complexes (Definition~\ref{def:hyperbox}), we only have linear maps $\De^{\eps' - \eps} : C^{\eps} \to C^{\eps'}$ in the case when $\eps$ and $\eps'$ are neighbors.

The hyperbox $\E(\dd)$ can be viewed as a union of several unit hypercubes: if all $d_i$'s are nonzero, there are $d_1d_2 \dots d_n$ unit hypercubes of dimension $n$, whereas if some $d_i$'s are zero, we get a union of unit hypercubes of smaller dimension. More precisely, let $\dd^\circ = (d_1^\circ, \dots, d_n^\circ)$ consist of the values $d_i^\circ = \max(d_i - 1, 0)$, and let $n^\circ$ be the number of nonzero $d_i$'s. Then $\E(\dd)$ is the union of the unit hypercubes $\eps + \E_{n^\circ}$ for $\eps \in \dd^\circ$. 

 \begin {definition}
 \label {def:betae} 
An {\em empty $\beta$-hyperbox} of size $\dd \in \N^n$ on a fixed
surface with marked points $(\Sigma, \ws, \zs)$ consists of a
collection $\{\betas^\eps\}_{\eps\in\E(\dd)}$ of strongly equivalent
sets of 
attaching beta curves $\betas^\eps$, indexed by $\eps
\in \E(\dd)$. Further, we require that, for each unit
hypercube of the form $\eps + \E_{n^\circ}$ with $\eps \in \dd^\circ$, the
corresponding Heegaard multi-diagram 
$(\Sigma, \{\betas^{\eps'} \}_{\eps'\in(\eps+\E_{n^\circ})},\ws,\zs)$ 
is generic and admissible.
\end {definition}

\begin {definition}
\label {def:betah}
Let $\{\beta^\eps\}_{\eps \in \E(\dd)}$ be an empty $\beta$-hyperbox, consisting of diagrams on a surface of genus $g$. Set $\mu = g/2$. A {\em filling} $\Theta$ of  the hyperbox consists of chain elements $$\Theta_{\eps , \eps'} \in \Chain^-_{\mu + \|\eps' - \eps\| -1}(\T_{\beta^{\eps}}, \T_{\beta^{\eps'}}, \zero),$$ one for each pair  $(\eps, \eps')$ such that $\eps < \eps'$ and $\eps, \eps'$ are neighbors. The chains $\Theta_{\eps , \eps'}$ are required to satisfy the following conditions:
\begin {itemize}
\item When $\eps < \eps'$ and $\| \eps' - \eps \| = 1$ (i.e. $\eps, \eps'$ are the endpoints of an edge in the hyperbox), $\Theta_{\eps , \eps'}$ is a cycle representing the maximal degree element $\theta_{\beta^{\eps}, \beta^{\eps'}}$ in Floer homology;
\item
For any $\eps < \eps'$ such that $\eps, \eps'$ are neighbors, we have 
\begin {equation}
\label {eq:compatibility}
\sum_{l=1}^{\| \eps' - \eps \|} \sum_{\{\eps = \eps^0 < \dots < \eps^l = \eps'\}} f(\Theta_{\eps^0, \eps^1} \otimes \dots \otimes \Theta_{\eps^{l-1}, \eps^l}) = 0. 
\end {equation}
\end {itemize}

The data consisting of an empty $\beta$-hyperbox $\{\beta^\eps\}_{\eps \in \E(\dd)}$ and  a filling by $\{\Theta_{\eps, \eps'}\}$ is simply called a {\em $\beta$-hyperbox}.
\end {definition}

\begin {remark} \label {rem:beta}
The simplest kind of $\beta$-hyperbox is a $\beta$-hypercube,  i.e. one with $\dd = (1, \dots, 1)$, so that $\E(\dd) = \E_n$. Then the Heegaard multi-diagram consisting of all $2^n$ curve collections has to be admissible. Further, any two $\eps, \eps' \in \E_n$ with $\eps \leq \eps'$ are neighbors, so for any such $\eps, \eps'$ the pair $(\betas^\eps, \betas^{\eps'})$ needs to come equipped with a chain $\Theta_{\eps, \eps'}$, such that these chains satisfy \eqref{eq:compatibility}. 

This kind of $\beta$-hypercube is quite natural, and has already appeared implicitly in the Heegaard Floer literature, e.g. in \cite{BrDCov}. There are two reasons why we need the more general kind of hyperbox. The first is that 
we want to allow for some of the beta pairs in the hyperbox (those for which $\eps, \eps'$ are not neighbors) to form non-admissible diagrams. These non-admissible diagrams make an appearance in the construction of basic systems, see Section~\ref{sec:basic} below. The second reason is that general hyperboxes appear naturally in the context of grid diagrams, see Section~\ref{sec:completegrid} below.
\end {remark}

\begin {lemma}
\label {lemma:filling}
Any empty $\beta$-hyperbox admits a filling. Moreover, if a filling is partially defined on the $m$-skeleton of a $\beta$-hyperbox (i.e., the elements $\Theta_{\eps, \eps'}$ are defined only for neighbors $\eps, \eps'$ with $\|\eps' - \eps\| \leq m$, and satisfy the required conditions), it can be extended to a filling on the whole hyperbox.
\end {lemma}

\begin {proof}
We construct the chain elements $\Theta_{\eps, \eps'}$ inductively on $\|\eps' - \eps\|$. 

When $\|\eps' - \eps\| =1$, we choose arbitrary cycles representing the maximal degree elements in homology. 

When $\|\eps' - \eps\| >1$, suppose we want to define $\Theta_{\eps, \eps'}$ and we have defined elements $\Theta_{\gamma, \gamma'}$ whenever $\|\gamma' - \gamma\| < \|\eps' - \eps\|$, satisfying \eqref{eq:compatibility}. Set
$$ c= \sum_{l=2}^{\| \eps' - \eps \|} \sum_{\{\eps = \eps^0 < \dots < \eps^l = \eps'\}} f(\Theta_{\eps^0, \eps^1} \otimes \dots \otimes \Theta_{\eps^{l-1}, \eps^l}).$$ 

Using \eqref{eq:f}, we obtain:

\begin {eqnarray*}
\del c &=& \sum_{l=2}^{\| \eps' - \eps \|} \sum_{\{\eps = \eps^0 < \dots < \eps^l = \eps'\}} f(f(\Theta_{\eps^0, \eps^1} \otimes \dots \otimes \Theta_{\eps^{l-1}, \eps^l})) \\
&=& \sum_{l=2}^{\| \eps' - \eps \|} \sum_{\{\eps = \eps^0 < \dots < \eps^l = \eps'\}} \sum_{\substack{0 \leq i < j \leq l \\ (i, j) \neq (0, l)} } f(\Theta_{\eps^0, \eps^1} \otimes \dots f(\Theta_{\eps^{i}, \eps^{i+1}} \otimes \dots \otimes \Theta_{\eps^{j-1}, \eps^j} ) \dots \otimes \Theta_{\eps^{l-1}, \eps^l})) \\
&=& 0.
\end {eqnarray*}

Here, for the last equality we applied  \eqref{eq:compatibility} to the pair $(\eps^i, \eps^j)$. 

Thus, $c$ is a cycle in $\Chain^-_{\mu + \|\eps' - \eps\| -2}(\T_{\beta^{\eps}}, \T_{\beta^{\eps'}}, \zero)$. If $\|\eps' - \eps\| =2$, then $c$ is the sum of two terms of the form $f(\Theta_{\eps, \eps^1} \otimes \Theta_{\eps_1, \eps'})$. Both of these terms represent the nonzero homology class of maximal degree (compare \cite[proof of Lemma 4.5]{BrDCov}), and therefore their sum $c$ is a boundary. If  $\|\eps' - \eps\| > 2$, then the respective homology group $\Chain^-_{\mu + \|\eps' - \eps\| -2}(\T_{\beta^{\eps}}, \T_{\beta^{\eps'}}, \zero)$ is zero (being beyond the maximal degree $\mu$), so $c$ is again a boundary. Thus, in any case we can choose $\Theta_{\eps, \eps'}$ such that $\del \Theta_{\eps, \eps'}= c$. Then  \eqref{eq:compatibility} is satisfied for the pair $(\eps, \eps')$. \end {proof}

\subsection {Hyperboxes of strongly equivalent Heegaard diagrams}
\label{sec:hse}
We define an {\em $\alpha$-hyperbox} on $(\Sigma, \ws, \zs)$ to be the same as  a $\beta$-hyperbox, except we denote the collections of curves by $\alphas$'s and, for any neighbors $\eps < \eps'$, we are given elements $\Theta_{\eps', \eps} \in \Chain^-(\T_{\alpha^{\eps'}}, \T_{\alpha^\eps}, \zero)$ rather than in $ \Chain^-(\T_{\alpha^{\eps}}, \T_{\alpha^{\eps'}}, \zero)$. The compatibility relation \eqref{eq:compatibility} has to be modified accordingly:
\begin {equation}
\label {eq:compatibility2}
 \sum_{l=1}^{\| \eps' - \eps \|} \sum_{\{\eps' = \eps^0 > \dots > \eps^l = \eps\}} f(\Theta_{\eps^0, \eps^1} \otimes \dots \otimes \Theta_{\eps^{l-1}, \eps^{l}}) = 0. 
\end {equation}

Now suppose that we have $\dd=(d_1, \dots, d_n) \in \N^n$, for some $n \geq 0$, and we are given  maps
$$ r_i : \{1, \dots, d_i \} \to \{\alpha, \beta\}, \ \ i=1, \dots, n. $$ 

We can then assign to each edge $(\eps, \eps')$ of the hyperbox $\E(\dd)$ a symbol $r(\eps, \eps') \in \{\alpha, \beta\}$ as follows: if the edge is parallel to the $i\th$ axis and its projection to that axis is the segment $[j-1, j]$, we choose $r(\eps, \eps') = r_i(j)$. We refer to $ \{r_i\}_{i=1}^n$ as {\em bipartition maps}.

Given bipartition maps as above and $\eps \in \E(\dd)$, we set
$$ \eps_i^{\alpha} = \# \bigl (r_i^{-1}(\alpha) \cap \{1, \dots, \eps_i\} \bigr) ,$$
$$ \eps_i^{\beta} = \# \bigl( r_i^{-1}(\beta) \cap \{1, \dots, \eps_i\} \bigr),$$
for $i=1, \dots, n$.  These define $n$-tuples $\eps^{\alpha} \in \E(\dd^{\alpha}), \eps^{\beta} \in \E(\dd^{\beta})$, where $\dd^\alpha = (d_i^{\alpha})_{i=1}^n, \  \dd^\beta = (d_i^{\beta})_{i=1}^n$ 
(where, of course, $d_i^\alpha=\# r_i^{-1}(\alpha)$ and
$d_i^\beta=\# r_i^{-1}(\beta)$).
Note that $\eps_i^\alpha + \eps_i^\beta = \eps_i$ for any $i$.

\begin {definition}
\label {def:alphabeta}
Choose $\dd \in \N^n$ and bipartition maps $r_i, i=1, \dots, n$. A {\em hyperbox $\Hyper$ of strongly equivalent Heegaard diagrams} consists of an $\alpha$-hyperbox of size $\dd^\alpha \in \N^n$, and a $\beta$-hyperbox of size $\dd^\beta \in \N^n$, both on the same multi-pointed surface  $(\Sigma, \ws, \zs)$. These are required to satisfy the following conditions. For each multi-index $\eps \in \E(\dd)$, we can consider the Heegaard diagram $$\Hyper_\eps = (\Sigma,  \alphas^{\eps^{\alpha}}, \betas^{\eps^{\beta}}, \ws, \zs).$$
For each unit hypercube $(\eps + \E_{n^\circ}) \subseteq \E(\dd)$, the curve collections appearing in the diagrams $\Hyper_{\eps'}$, for all $\eps' \in  (\eps + \E_{n^\circ})$, are required to form a generic, admissible Heegaard multi-diagram.
\end {definition}

We can view some of the information in $\Hyper$ in the following way. The hyperbox $\Hyper$ has at each of its vertices a Heegaard diagram $\Hyper_\eps$, such that all of these are strongly equivalent, and each edge in the hyperbox corresponds to changing either the alpha or the beta curves. Further, we have elements $\Theta^{\alpha}_{\eps', \eps}= \Theta_{\eps', \eps}$ for $\eps < \eps'$ neighbors in $\E(\dd^\alpha)$, and 
$\Theta^{\beta}_{\eps, \eps'}= \Theta_{\eps, \eps'}$ for $\eps < \eps'$ neighbors in $\E(\dd^\beta)$. 
We usually refer to a hyperbox $\Hyper$ as going between the Heegaard diagrams $\Hyper_{(0,\dots, 0)}$ and $\Hyper_{(d_1, \dots, d_n)}$. 

\begin {remark}
\label {rem:switch} Suppose we have a hyperbox $\Hyper$ of strongly equivalent Heegaard diagrams representing a link $\orL \subset Y$. Then the hyperbox $\Hyper$ naturally gives rise to a {\em reduced hyperbox} $r_{\orM}(\Hyper)$, with the same size and bipartition maps, made of the reduced Heegaard diagrams $(r_{\orM}(\Hyper))_\eps = r_{\orM}(\Hyper_\eps)$, see Definition \ref{def:reduce}. Indeed, in $r_{\orM}(\Hyper)$ we can take the images of the $\Theta$-chain elements from $\Hyper$ under the corresponding inclusion maps defined in Section~\ref{sec:inclusions}. 
\end {remark}

\subsection {Hyperboxes of Floer complexes}
\label {sec:hyperfloer}
 Let $\Hyper$ be a hyperbox of strongly equivalent Heegaard diagrams, of size $\dd\in \N^n$ and with partition maps $r_i, i=1, \dots, n$. At each vertex we have an admissible Heegaard diagram $\Hyper_\eps$. Let us assume that $\Hyper_\eps$ represents a link $\orL$ in an integral homology sphere $Y$. Fix $\s \in \bH(L)$. We will construct an associated hyperbox of chain complexes as in Section~\ref{sec:hyperv}. To each vertex $\eps \in \E(\dd)$ we assign the generalized Floer chain complex
$$ C^\eps_* = \Am_*(\Hyper_\eps, \s) = \Am_*(\T_{\alpha^{\eps^\alpha}}, \T_{\beta^{\eps^\beta} }, \s),$$
compare Section~\ref{sec:chains}.

We define the linear maps $\De^{\eps'-\eps}_{\eps} : C^{\eps} \to C^{\eps'}$ by the formula
$$ \De^{\eps'-\eps}_{\eps}(\x) =\sum_{l, p} \sum_{\{\eps^{'\alpha} = \gamma^0 > \dots > \gamma^l = \eps^\alpha\}}  \sum_{\{\eps^\beta = \zeta^0 < \dots < \zeta^p = \eps^{'\beta}\}}  f(\Theta^\alpha_{\gamma^0, \gamma^1} \otimes \dots \otimes \Theta^\alpha_{\gamma^{l-1}, \gamma^l} \otimes \x \otimes \Theta^\beta_{\zeta^0, \zeta^1} \otimes \dots \otimes \Theta^\beta_{\zeta^{p-1}, \zeta^p} ).
$$

\begin {example}
When $\eps = \eps'$ the corresponding map $\De^0_\eps$ is just the differential $f = \del$ on the generalized Floer complex $\Am_*(\Hyper_\eps, \s)$. 
\end {example}

\begin {example}
\label {ex:edge}
When $\|\eps' - \eps\| = 1, $ i.e. $(\eps, \eps')$ is an edge in the hypercube then $\De^{\eps'-\eps}_{\eps} $ is given by counting holomorphic triangles with one specified vertex. Since that vertex is a theta cycle representing the maximal degree element in homology, the chain map  $\De^{\eps'-\eps}_{\eps} $ is a chain homotopy equivalence.
\end {example}

\begin {example}
When $\|\eps' - \eps\| = 2$, we distinguish two cases, according to whether the edges of the square between $\eps$ and $\eps'$ are marked with only one or both of the symbols $\alpha, \beta$. In the first case, assuming they are marked with $\beta$, the map $\De^{\eps'-\eps}_{\eps} $ is a sum of three different polygon maps: one counting triangles with one vertex at $\Theta^\beta_{\eps^\beta, \eps^{'\beta}}$ and two counting quadrilaterals with two vertices at  $\Theta^\beta_{\eps^\beta, \zeta}$ and  $\Theta^\beta_{\eps^\beta, \zeta}$, where $\zeta$ is an intermediate multi-index between $\eps^\beta$ and $\eps^{'\beta}$. (There are two possibilities for $\zeta$.) In the second case, when $\eps, \eps'$ are the vertices of a square with two edges marked with $\alpha$ and two with $\beta$, the map $\De^{\eps'-\eps}_{\eps} $ is given by counting quadrilaterals with two specified vertices at $\Theta^\beta_{\eps^\beta, \eps^{'\beta}}$ and $\Theta^\alpha_{\eps^{'\alpha}, \eps^{\alpha}}$. 
\end {example}

\begin {lemma}
\label {lemma:cde}
$\Am(\Hyper, \s) = \bigl( (C^{\eps})_{\eps \in \E(\dd)}, (\De^{\eps})_{\eps \in \E_n} \bigr)$ is a hyperbox of chain complexes.
\end {lemma}

\begin {proof}
We need to check Equation~\eqref{eq:d2}, i.e. that for any $\eps, \eps'' \in \E_n$ we have
\begin {equation}
\label {eq:d3}
 \sum_{\{\eps' | \eps < \eps' < \eps''\}}    \De^{\eps'' - \eps'}_{\eps'} \circ \De^{\eps' - \eps}_\eps = 0. 
 \end {equation}

Indeed, the summation in \eqref{eq:d3} equals
$$ \sum f(\Theta^\alpha_{\gamma^0, \gamma^1} \otimes \dots \otimes \Theta^\alpha_{\gamma^{i-1}, \gamma^i} \otimes f(\Theta^\alpha_{\gamma^{i}, \gamma^{i+1}} \otimes \dots \otimes  \Theta^\alpha_{\gamma^{l-1}, \gamma^l} \otimes \x \otimes \Theta^\beta_{\zeta^0, \zeta^1} \otimes \dots \otimes  \Theta^\beta_{\zeta^{j-1}, \zeta^j}) \otimes  \Theta^\beta_{\zeta^{j}, \zeta^{j+1}} \otimes \dots \otimes \Theta^\beta_{\zeta^{p-1}, \zeta^p} ),$$
where the sum is taken over all possible $l, p, i, j$ and multi-indices $\eps^{'\alpha} = \gamma^0 > \dots > \gamma^l = \eps^\alpha, \  \eps^\beta = \zeta^0 < \dots < \zeta^p = \eps^{'\beta}$. Applying Equation~\eqref{eq:f} we find that this sum further equals
$$ \sum f(\Theta^\alpha_{\gamma^0, \gamma^1} \otimes \dots \otimes f( \Theta^\alpha_{\gamma^{i-1}, \gamma^i} \otimes \dots \otimes \Theta^\alpha_{\gamma^{j-1}, \gamma^{j}}) \otimes \dots \otimes  \Theta^\alpha_{\gamma^{l-1}, \gamma^l} \otimes \x \otimes \Theta^\beta_{\zeta^0, \zeta^1} \otimes \dots \otimes  \Theta^\beta_{\zeta^{p-1}, \zeta^p} )+$$
$$ \sum f(\Theta^\alpha_{\gamma^0, \gamma^1} \otimes \dots \otimes  \Theta^\alpha_{\gamma^{l-1}, \gamma^l} \otimes \x \otimes \Theta^\beta_{\zeta^0, \zeta^1} \otimes \dots \otimes  f(\Theta^\beta_{\zeta^{i-1}, \zeta^i} \otimes  \dots \otimes \Theta^\beta_{\zeta^{j-1}, \zeta^{j}}) \otimes \dots \otimes \Theta^\beta_{\zeta^{p-1}, \zeta^p} ).$$

Both of these sums vanish. Indeed, let us fix $i, l-j, p$ and all $\gamma$'s and $\zeta$'s in the first sum, except for $\gamma^i, \dots, \gamma^{j-1}$. Then the corresponding sum of $f( \Theta^\alpha_{\gamma^{i-1}, \gamma^i} \otimes \dots \otimes \Theta^\alpha_{\gamma^{j-1}, \gamma^{j}})$ over $\gamma^i, \dots, \gamma^{j-1}$ is zero by Equation~\eqref{eq:compatibility2}. Similarly, by applying Equation~\eqref{eq:compatibility} we find that the second sum is also zero. 
\end {proof}

In the case where the hyperbox $\Hyper$ consists of link-minimal diagrams, we can use the alternative description of $\Am(\Hyper_\eps, \s)$ as free complexes $\Chain^-_*(\Hyper_\eps, \s)$, as in Section~\ref{sec:alternative}. We then denote the resulting hyperbox of chain complexes by 
$\Chain^-(\Hyper, \s)$.

\subsection {Moves on hyperboxes} \label {sec:movesh} In this section we describe a series of moves on hyperboxes of strongly equivalent Heegaard diagrams. In light of Lemma~\ref{lemma:filling}, we see that  hyperboxes are easy to construct, and thus are rather flexible objects.

Recall the list of Heegaard moves between multi-pointed Heegaard diagrams from Section~\ref{sec:hmoves}. Suppose now that $\Hyper$ is a hyperbox of strongly equivalent Heegaard diagrams. We have a similar list of {\em hyperbox Heegaard moves} on $\Hyper$:

\begin {enumerate}

\item A {\em 3-manifold isotopy} of $\Hyper$ consists of applying the same 3-manifold isotopy to all the Heegaard diagrams $\Hyper^\eps$ simultaneously, keeping the same partition maps and taking the $\Theta$-chain elements (which are linear combinations of collections of intersection points) to the corresponding linear combinations in the isotopic diagrams;

\item An {\em index one/two stabilization} of a hyperbox $\Hyper$ consists of a
simultaneous index one/two stabilization of all diagrams $\Hyper_\eps$, in the same
position. Note that if two collections of attaching curves are strongly equivalent, they remain so after the stabilization. With regard to the fillings, we need to pair each $\Theta$-chain element with the unique 
intersection point between the new $\alpha$ curve and new $\beta$ curve. The fact that the new $\Theta$-chain elements still satisfy the relations \eqref{eq:compatibility} follows from the argument of~\cite[Lemma~4.7]{HolDiskFour} generalized to polygons;

\item A {\em free index zero/three stabilization} of a hyperbox $\Hyper$ consists of a
simultaneous free index zero/three stabilization of all diagrams $\Hyper_\eps$, in the same position. We again note that if two collections of attaching curves are strongly equivalent, they remain so after the stabilization. With regard to the fillings, we pair each $\Theta$-chain element with the maximal degree intersection point between the new $\alpha$ and $\beta$ curves. The fact that the new $\Theta$-chain elements still satisfy the relations \eqref{eq:compatibility} follows from the first part of Proposition~\ref{prop:StabPolygon};

\item An {\em index zero/three link stabilization} of a hyperbox $\Hyper$ consists of a
simultaneous index zero/three link stabilization of all diagrams $\Hyper_\eps$, in the same position. The fillings are constructed as in (iv), by pairing with the maximal degree intersection point between the new curves. The relations \eqref{eq:compatibility} are then a consequence of Proposition~\ref{prop:PolyDestab03a};

\item Instead of curve isotopies and handleslides, we now have global shifts, defined as follows.  
Consider two hyperboxes $\Hyper$ and $\Hyper'$ having the same size $\dd \in \N^n$ and the same partition maps $r_i$. A {\em global shift} $S$ from $\Hyper$ to $\Hyper'$ is an $(n+1)$-dimensional hyperbox of size $(\dd, 1) \in \N^{n+1}$, such that its sub-hyperbox corresponding to $\eps_{n+1} = 0$ is $\Hyper$ and its sub-hyperbox corresponding to $\eps_{n+1} = 1$ is $\Hyper'$. Note that there are two kinds of global shifts, $\alpha$ and $\beta$, according to the value of the map $r_{n+1} : \{1\} \to \{\alpha, \beta\}$. Observe that, by definition, a global shift does not change the Heegaard surface with basepoints $(\Sigma, \ws, \zs)$; 

\item An {\em elementary enlargement}, to be defined later
(see Definition~\ref{def:Enlargement} below);

\item Inverses to the above. In particular, the inverse process to a stabilization is called {\em destabilization}, and the inverse of an elementary enlargement is called a {\em contraction}.
\end {enumerate}

Note that the moves (i)-(v) (and their inverses) preserve the size of the respective hyperbox. 

As a particular kind of 3-manifold isotopy, we define a {\em surface isotopy} of $\Hyper$ to consist in applying the same surface isotopy to all the Heegaard diagrams simultaneously, compare Definition~\ref{def:ab} (b). If $\Hyper$ and $\Hyper'$ are surface isotopic hyperboxes, we write $\Hyper \cong \Hyper'$.

Let us now focus on understanding global shifts further. Note that a global shift induces a chain map between the respective hyperboxes of generalized Floer chain complexes, compare Definition~\ref{def:chmap}. It turns out that this map is always a chain homotopy equivalence. Before proving this fact, we need some preliminaries.

\begin {definition}
Let $\Hyper$ and $\Hyper'$ be two hyperboxes of strongly equivalent Heegaard diagrams, such that either: 

(a) the corresponding alpha curves on $\Hyper$ and $\Hyper'$ coincide, whereas the corresponding beta curve collections approximate each other $(\beta^{\eps} \approx \beta^{'\eps}$) in the sense of Definition~\ref{def:approx}, or

(b) the corresponding beta curves coincide, and the corresponding alpha curve collections  approximate each other. 

Further, suppose that the approximations are sufficiently small (for Lemma~\ref{lem:NearestPoints} to hold) and suppose the $\Theta$-chains correspond to each other under the respective nearest point maps. We then say that $\Hyper'$ {\em approximates} $\Hyper$, and write $\Hyper \approx \Hyper'$. If we are in case (a), we call this an approximation of type $\beta$, and if we are in case (b), we call it of type $\alpha$.
\end {definition}

\begin {lemma}
\label {lemma:kan} 
Suppose $\Hyper \approx \Hyper'$, where the approximation is of type $\beta$. Then, there exists a canonical global shift between $\Hyper$ and $\Hyper'$, such that along the new edges we see the canonical elements (mentioned in Section~\ref{sec:betah})
$$\Theta^\can_\eps = \Theta^\can_{\betas^\eps, \betas^{'\eps}} \in \T_{\beta^\eps} \cap \T_{\beta^{'\eps}}.$$
\end {lemma}

\begin {proof}
For simplicity, we explain the construction in the case when the hyperboxes are one-dimensional of length one. Thus, $\Hyper$ consists of two curve collections $\betas^0$ and $\betas^1$ and a chain element $ \Theta$ relating them. Similarly, $\Hyper'$ consists of curve collections $\betas^{'0}$ and $\betas^{'1}$ and a chain element $ \Theta'$ that is the image of $\Theta$ under the nearest point map. Along the new edges we are required to place the two intersection points $ \Theta^\can_0$ and $\Theta^\can_1$.

We also have a cycle
$$ \Theta'' \in \Chain^-(\T_{\beta^0}, \T_{\beta^{1'}}, \zero)$$ 
that is the image of $\Theta$ under the nearest point map changing
$\betas^1$ into $\betas^{'1}$. (This is a cycle because of
Lemma~\ref{lem:NearestPoints}.)
Recall that there is a canonical homotopy between 
the continuation map from $\Chain^-(\T_{\beta^{0'}},\T_{\beta^{1'}}, \zero)\to
\Chain^-(\T_{\beta^{0}},\T_{\beta^{1'}}, \zero)$
and the triangle map $\x\mapsto f(\Theta_0^\can\otimes \x)$, see
\cite[proof of Proposition
11.4]{LipshitzCyl}, or \cite[proof of Theorem
6.6]{OzsvathStipsicz}. This chain homotopy is defined by counting holomorphic
bigons with a varying almost-complex structure and stretching of the neck.
Thus,
$$ f(\Theta^\can_0 \otimes \Theta') = \Theta'' + dH_0,$$
where $H_0$ is the image of $\Theta'$ under the canonical chain homotopy.
Similarly, we have
$$ f(\Theta \otimes \Theta^\can_1) = \Theta'' + dH_1,$$
where $H_1$ is the image of a map counting stretched bigons with one vertex at $\Theta$. We then place the chain element  $H_0 + H_1$ on the  two-dimensional face of our global shift. Relation ~\eqref{eq:compatibility} is satisfied.

If $\Hyper$ and $\Hyper'$ are one-dimensional (of arbitrary length), we place chain elements as above on all their two-dimensional faces. Further, the construction generalizes to hyperboxes of higher dimension, by considering continuation maps with several stretched necks.
\end {proof}

The global shift constructed in Lemma~\ref{lemma:kan} is called the {\em identity shift}. There exists a similar identity shift between hyperboxes that approximate each other through an approximation of type $\alpha$.
 
Now, suppose $\Hyper \approx \Hyper'$ and fix $\s \in \bH(L)$. Let $(C^\eps, D^\eps)$ and $(C'^\eps, D'^\eps)$ be the hyperboxes of generalized Floer complexes (for the value $\s$) associated to $\Hyper$ and $\Hyper'$, respectively. Note that there is a nearest point map from $C^\eps$ to $C'^{\eps}$, defined by taking each generator (intersection point) to its image under the approximation. This is a chain map provided that the almost complex structures are chosen in a compatible way. In fact, by adding zeros on the higher dimensional faces, we obtain a chain map between the respective hyperboxes of generalized Floer complexes.  We refer to this the nearest point map, too.

\begin {lemma}
\label {lemma:id}
An identity shift induces a map on the hyperboxes of Floer complexes which is  chain homotopic to the nearest point map. 
\end {lemma}

\begin {proof}
The triangle map induced by an approximation is chain homotopic to the nearest point map, see \cite[proof of Theorem 6.6]{OzsvathStipsicz}. The argument in \cite{OzsvathStipsicz} can also be applied to higher polygons, with the result that the chain homotopies between triangle and nearest point maps lift to chain homotopies between the respective hyperboxes.
\end {proof}

\begin {lemma}
\label {lemma:zzz}
The chain map $F(S)$ induced by any global shift $S$ on hyperboxes of Floer complexes is a chain homotopy equivalence.
\end {lemma}

\begin {proof}
Let $S$ be a global shift between two $n$-dimensional hyperboxes $\Hyper, \Hyper'$ of size $\dd$, so that $S$ has size $(\dd, 1)$. Let $\Hyper'', \Hyper'''$ be two hyperboxes that approximate $\Hyper$ and each other. 

We construct an $(n+2)$-dimensional hyperbox $\tilde S$ of size $(\dd, 1, 1)$ as follows. Its sub-hyperbox corresponding to $\eps_{n+2} = 0$ is $S$, so that $\eps_{n+1} = \eps_{n+2} = 0$ corresponds to $\Hyper$ and $\eps_{n+1}=1, \eps_{n+2} =0$ is $\Hyper'$. Its sub-hyperbox corresponding to $\eps_{n+1} = 1$ is any global shift between $\Hyper'$ and $\Hyper''$, for example the reverse $S^r$ of $S$ composed with the identity shift from $\Hyper$ to $\Hyper''$. Thus $\eps_{n+1} = \eps_{n+2} =1$ corresponds to $\Hyper''$. For $\eps_{n+1} =0, \eps_{n+2} =1$ we take the other hyperbox $\Hyper'''$ that approximates  $\Hyper$. For $\eps_{n+1}=0$ we choose the identity shift between $\Hyper$ and $\Hyper'''$, and for $\eps_{n+2}=1$ the identity shift between $\Hyper'''$ and $\Hyper''$. We then fill in the remaining faces of the hyperbox $\tilde S$ (those corresponding to increasing both $\eps_{n+1}$ and $\eps_{n+2}$ by $1$) with $\Theta$ chain elements  in an arbitrary way.

 On the level of Floer complexes, the hyperbox $\tilde S$ then produces a chain homotopy between $F(S^r) \circ F(S)$ and the identity.
\end {proof}

\begin{definition}
\label{def:Enlargement}
We define an {\em elementary enlargement} of a hyperbox $\Hyper$ of strongly equivalent Heegaard diagrams by analogy with the corresponding concept from Section~\ref{sec:ele}. Let $\dd \in \N^n$ be the size of $\Hyper$.  Pick $i_0 \in \{1, \dots, n\}$ and $j_0 \in \{0, 1, \dots, d_k\}$. Define $\dd^+ = \dd + \tau_k$ as in Section~\ref{sec:ele}. We also choose a symbol $\kappa = \alpha$ or $\beta$.

We construct a new hyperbox $\Hyper^+$ of size $\dd^+$, with bipartition maps $r_i^+ = r_i$ for $i\neq i_0$, and
$$ r_{i_0}^+( j) = \begin {cases} 
r_{i_0}(j) & \text{if } j \leq j_0 \\
\kappa & \text{if } j=j_0+1 \\
r_{i_0}(j-1) & \text{if } j \geq j_0+2.
\end {cases} $$

The hyperbox $\Hyper^+$ is obtained from $\Hyper$ by splitting the latter into two halves along the hyperplane $\eps_{i_0}=j_0$, adding a new sub-hyperbox at $\eps_{i_0}=j_0 +1$ that approximates the one at $\eps_{i_0}=j_0$, and inserting the identity global shift between them. We say that $H^+$ is obtained from $H$ by an elementary enlargement. The reverse process is called {\em elementary contraction}. 
\end {definition}

In Section~\ref{sec:ele} we defined a similar notion of elementary enlargement for hyperboxes of chain complexes. If $\Hyper^+$ is the elementary enlargement of a hyperbox of Heegaard diagrams $\Hyper$, let $\Chain^-(\Hyper, \s)$ and $\Chain^-(\Hyper^+, \s)$ be the respective hyperboxes of generalized Floer chain complexes, as in Section~\ref{sec:hyperfloer}. Let also $\Chain^-(\Hyper, \s)^+$ be the corresponding elementary enlargement of the hyperbox of chain complexes $\Chain^-(\Hyper, \s)$. By applying Lemma~\ref{lemma:id}, we see that $\Chain^-(\Hyper, \s)^+$ is chain homotopy equivalent to $\Chain^-(\Hyper^+, \s)$. 

\begin {lemma}
\label {lemma:moves}
Let $\Hyper$ and  $\Hyper'$ be two hyperboxes of strongly equivalent Heegaard diagrams representing the same oriented link $\orL \subset Y$, and having the same dimension. Then:

(a) The hyperboxes $\Hyper$ and $\Hyper'$ can be related by a sequence of hyperbox Heegaard moves. 

(b) Suppose further that the diagrams in $\Hyper$ are link-minimal, and the diagrams in $\Hyper'$ have $m$ basepoints of type $z$, and therefore $m-\ell$ subsidiary $w$ basepoints, as in Proposition~\ref{prop:RefinedMoves}. (We choose the same basepoints to be subsidiary in all the diagrams in $\Hyper'$.) Then, $\Hyper'$ can be obtained from $\Hyper$ using a sequence of hyperbox Heegaard moves that includes exactly $m-\ell$ index zero/three link stabilizations, with each of these stabilizations introducing a subsidiary basepoint.
\end {lemma}

\begin {proof}
For part (a), use 3-manifold isotopies, stabilizations, and destabilizations to transform the hyperboxes into two new ones that have the same underlying Heegaard surface $(\Sigma, \ws, \zs)$. After a few more index one/two stabilizations as in \cite[Lemma 2.4]{MOS}, we can also arrange so that the $\alpha$ curves in $\Hyper$ are strongly equivalent to the $\alpha$ curves in $\Hyper'$, and the $\beta$ curves in $\Hyper$ are strongly equivalent to the $\beta$ curves in $\Hyper'$. 
 We can then use elementary enlargements to arrange so that the hyperboxes $\Hyper$ and $\Hyper'$ have the same size and bipartition maps. Using the winding procedure from \cite[Section 5]{HolDisk}, we can find a sequence of empty hyperboxes (of the same size and with the same bipartition maps) that interpolate between the two: $\Hyper=\Hyper^0, \Hyper^1, \dots, \Hyper^p = \Hyper'$ such that each pair $(\Hyper^{j-1}, \Hyper^{j})$ forms an empty hyperbox of one dimension bigger (satisfying the required admissibility conditions). We then choose arbitrary fillings of the intermediate hyperboxes $\Hyper^j$, see Lemma~\ref{lemma:filling}, as well as arbitrary global shifts between $\Hyper^{j-1}$ and $\Hyper^j$, for $j=1, \dots, p$. The result is a sequence of global shifts relating $\Hyper$ and $\Hyper'$.

Part (b) follows by combining the strategy above with Proposition~\ref{prop:RefinedMoves}.
\end {proof}

Recall the invariance statement about generalized Floer complexes for links (Theorem~\ref{thm:LinkInvariance}). We have an analogous result for hyperboxes of strongly equivalent Heegaard diagrams representing a link:

\begin{proposition}
\label {prop:moves2}
Let  $\Hyper, \Hyper'$ be two hyperboxes of strongly equivalent Heegaard diagrams that represent the same oriented link $\orL \subset Y$, and have the same dimension. Let $\Am(\Hyper, \s)$ and $\Am(\Hyper', \s)$ be the hyperboxes of generalized Floer chain complexes associated to $\Hyper$ and $\Hyper'$, respectively, for some $\s \in \bH(L)$. Pick also natural inclusions $\Ring(L) \hookrightarrow \Ring(\he)$ and $\Ring(L) \hookrightarrow \Ring(\he')$, as in Theorem~\ref{thm:LinkInvariance}. Then, the hypercubes obtained from $\Am(\Hyper, \s)$ and $\Am(\Hyper', \s)$ by compression (cf. Section~\ref{sec:compression1}) are chain homotopy equivalent over $\Ring(L)$, in the sense of Section~\ref{sec:chmaps}.
\end {proposition}

\begin {proof}
By Lemma~\ref{lemma:moves}, it suffices to investigate the effect of hyperbox Heegaard moves on the corresponding (compressed) hyperboxes of generalized Floer complexes. Index one/two stabilizations and 3-manifold isotopies produce chain homotopy equivalences between the respective hyperboxes, see \cite{HolDisk} and \cite{HolDiskFour}.  The same goes for global shifts, according to  Lemma~\ref{lemma:zzz}. The chain homotopy equivalences then descend to the compressed hyperboxes, according to Lemma~\ref{lemma:chainhe}. Elementary enlargements also produce quasi-isomorphisms between the respective compressions, see Lemma~\ref{lemma:ci}. Free index zero/three stabilizations result in equivalences; cf. \cite[Proposition 6.5]{Links} and Proposition~\ref{prop:StabPolygon}. Index zero/three link stabilizations also result in equivalences, as a consequence of Corollary~\ref{cor:rhonormal} and Proposition~\ref{prop:PolyDestab03a}. 
\end {proof}

\subsection {Sublinks and hyperboxes}
\label {sec:sublinks}

Let $\orL \subset Y$ be an oriented link, and $M \subseteq L$ a sublink. We choose an orientation $\orM$ of $M$, not necessarily the one induced from $\orL$. We denote by $\orL - M$ the sublink $L-M$ with the orientation induced from $\orL$.

\begin {definition}
\label {def:pair}
A {\em hyperbox of Heegaard diagrams for the pair $(\orL, \orM)$} is an $m$-dimensional hyperbox $\Hyper$ of strongly equivalent Heegaard diagrams representing the link $\orL - M$, together with an ordering $M_1, \dots, M_m$ of the components of $M$. (Here $m$ is the number of such components.)
\end {definition}

\begin {remark} In particular, a hyperbox of Heegaard diagrams for a pair $(\orL, \emptyset)$ is simply a Heegaard diagram for $\orL$.
\end {remark}

Definition~\ref{def:pair} may appear mysterious at first. The intuition behind it is that a hyperbox $\Hyper$ for a pair $(\orL, \orM)$ corresponds to subtracting the link $\orM$ from $\orL$. Indeed, although this is not required by the definition, the initial diagram $\Hyper_{(0,\dots, 0)}$ of $\Hyper$  will always be obtained from a Heegaard diagram for $\orL$ by deleting some of its basepoints. The final diagram will represent $\orL - M$, and the ordering of the components is telling us the parts of the hyperbox associated to subtracting the respective components $M_i$ from $\orL$. Indeed, we can think of the ordering as a one-to-one correspondence between the coordinate axes of the hyperbox and the components of $M$, where the $i\th$ coordinate corresponds to $M_i$.

To make this more precise, let $\Hyper^{\orL, \orM}$ be a hyperbox for the pair $(\orL, \orM)$. Then, for each $M' \subseteq M$, we denote by $\Hyper^{\orL, \orM}(M')$ the Heegaard diagram $\Hyper^{\orL, \orM}_{\eps(M')}$, where $\eps(M')$ is the multi-index with components $\eps(M')_i, i=1, \dots, m$, given by
$$ \eps(M')_i = \begin {cases} d_i & \text{ if } M_i \subseteq M', \\
0 & \text{otherwise.} \end {cases}$$

Informally, we think of $\Hyper^{\orL, \orM}(M')$ as the intermediate step in the hyperbox obtained after subtracting $M'$ from $\orL$. In particular, the initial and final vertex of the hyperbox $\Hyper^{\orL, \orM}$ are $\Hyper^{\orL, \orM}(\emptyset)$ and $\Hyper^{\orL, \orM}(M)$, respectively. 

For every $M' \subseteq M$, there is a sub-hyperbox of $\Hyper^{\orL, \orM}$ going from $\Hyper^{\orL, \orM}(M')$ to $\Hyper^{\orL, \orM}(M)$. We denote it by $\Hyper^{\orL, \orM}(M', M)$. Note that $ \Hyper^{\orL, \orM}(M', M)$ is a hyperbox associated to the pair  $(\orL-M', \orM - M')$.  We denote also by $\Hyper^{\orL, \orM}(\emptyset, M')$ the sub-hyperbox of $\Hyper^{\orL, \orM}$ going from $\Hyper^{\orL, \orM}(\emptyset)$ to $\Hyper^{\orL, \orM}(M')$, which is complementary to $\Hyper^{\orL, \orM}(M', M)$.

We say that two hyperboxes $\Hyper^{\orL, \orM}$ and $\tilde \Hyper^{\orL, \orM}$ for the same pair $(\orL, \orM)$ are {\em isotopic} if they have the same same ordering of the components of $M$, and the underlying hyperboxes of strongly equivalent Heegaard diagrams are isotopic. If this is the case, we write $\Hyper^{\orL, \orM} \cong \tilde \Hyper^{\orL, \orM}$.

Up to now, the orientation for $M$ has not played any role in the definition of a hyperbox for a pair $(\orL, \orM)$. However, it played a role when we discussed reduction (see Definition ~\ref{def:reduce} and Remark~\ref{rem:switch}). Observe that if $\Hyper^{\orL, \orM}$ is a hyperbox for a pair $(\orL, \orM)$, its reduction at an oriented sublink $\orN \subseteq L- M$, denoted $r_{\orN}(\Hyper^{\orL, \orM})$, is a hyperbox for the pair $(\orL - N, \orM)$.

\begin{figure}
\begin{center}
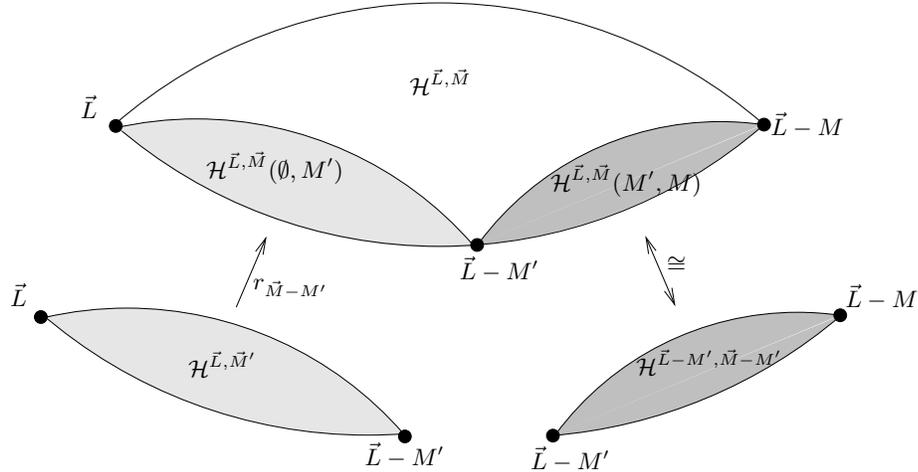
\end{center}
\caption {{\bf Compatibility between hyperboxes.} We represent here hyperboxes schematically by lenses. On top we have a hyperbox $
\Hyper^{\orL', \orM}$ and two sub-hyperboxes, shown as sub-lenses with different shadings. The compatibility condition (a) in Definition~\ref{def:comp1}  is the equality between the two lightly shaded lenses. Condition (b) is the equality between the two darkly shaded lenses.
}
\label{fig:hyper1}
\end{figure}

We are now ready to mention the following notions of compatibility for hyperboxes (see Figure~\ref{fig:hyper1}):

\begin {definition}
\label {def:comp1}
Let $\Hyper^{\orL, \orM}$ be a hyperbox for a pair $(\orL, \orM)$, and $\Hyper^{\orL, \orM'}$ be a hyperbox for the pair $(\orL, \orM')$, where $M'$ is a sublink of $M$ with the orientation $\orM'$ induced from $\orM$. Let also $\Hyper^{\orL-M', \orM-M'}$ be a hyperbox for the pair $(\orL-M', \orM-M')$.

$(a)$ We say that the hyperboxes $\Hyper^{\orL, \orM}$ and $\Hyper^{\orL, \orM'}$ are {\em compatible} if  
\begin {equation}
\label {eq:comp1}
r_{\orM-M'}(\Hyper^{\orL, \orM'}) = \Hyper^{\orL, \orM}(\emptyset, M').
\end {equation}

$(b)$ We say that the hyperboxes $\Hyper^{\orL, \orM}$ and $\Hyper^{\orL-M', \orM-M'}$ are {\em compatible} if there is a surface isotopy
\begin {equation}
\label {eq:comp2}
\Hyper^{\orL- M', \orM-M'} \cong \Hyper^{\orL,  \orM} (M', M). 
 \end {equation}
\end {definition}

\subsection {Complete systems of hyperboxes} 
\label {sec:complete}
In the following definition (and in the rest of this subsection), all sublinks of a link $L$ that are denoted $L'$ will come with the orientation $\orL'$ induced from $\orL$ and, similarly, all sublinks in a link $M$ that are denoted $M'$ will come with the orientation $\orM'$ induced from $\orM$.

\begin {definition}
\label {def:precomplete}
A {\em complete pre-system of hyperboxes} $\Hyper$ representing the
link $\orL$ consists of a collection of hyperboxes, subject
to certain compatibility conditions, as follows.
For each pair of subsets $M\subseteq L'\subseteq L$,
and each orientation $\orM\in\Omega(M)$, the complete pre-system
assigns a hyperbox $\Hyper^{\orL',\orM}$ for the pair 
$(\orL', \orM)$. 
Moreover, the hyperbox $ \Hyper^{\orL', \orM}$ is
required to be compatible with both $\Hyper^{\orL', \orM'}$ and $\Hyper^{\orL'- M',
\orM-M'}$.
\end {definition}

In particular, note that a complete pre-system contains hyperboxes of the form $\Hyper^{\orL', \emptyset}$, which are zero-dimensional; in other words, they consist of a single Heegaard diagram, which we denote $\he^{L'}$. The diagram $\he^{L'}$ represents the sublink $L' \subseteq Y$.

As previously mentioned in Section~\ref{sec:sublinks}, we think of the hyperbox $\Hyper^{\orL', \orM}$ as a way of de-stabilizing $\orL'$ at the components of $\orM$. Indeed, $\Hyper^{\orL', \orM}$ goes from $r_{\orM}(\he^{L'})$ (i.e. $\he^{L'}$ with half of the basepoints on $M$ deleted, according to the orientation $\orM$) to $\he^{L' - M}$.

\begin{figure}
\begin{center}
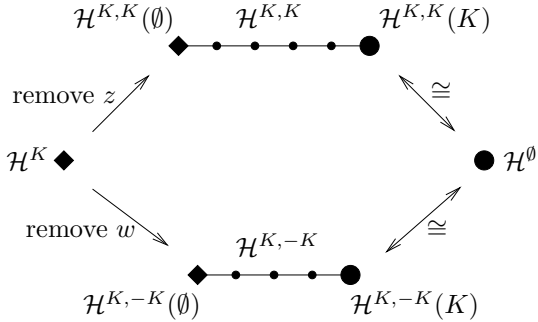
\end{center}
\caption {{\bf A complete pre-system of hyperboxes for a knot $\orK$.} For simplicity, we assume here that $\Hyper^K$ has no free basepoints, and exactly two linked basepoints ($w$ and $z$). The symbols $\cong$ indicate surface isotopies.
}
\label{fig:completeK}
\end{figure}

\begin {example}
\label{ex:completeknot}
Let $\orL = \orK \subset Y$ be a knot. Then, a complete pre-system of hyperboxes for $\orK$ consists of two multi-pointed Heegaard diagrams: $\he^{K}$ for $K$ and $\he^{\emptyset}$ for $Y$ itself, together with two one-dimensional hyperboxes of strongly equivalent Heegaard diagrams: one, $\he^{ K,  K}$, going from a diagram $\he^{ K,  K}(\emptyset)$ (which is just $\he^{K}$ with the $z$ points on $K$ removed) to some diagram $\he^{ K,  K}(K)$ (surface isotopic to $\he^{\emptyset}$); and another one, $\he^{ K, - K}$, going from the diagram $\he^{ K, - K}(\emptyset)$ (which is just $\he^{K}$ with the $w$ points removed), down to a diagram $\he^{ K, - K}(K)$ (surface isotopic to $\he^{\emptyset}$). This is illustrated in Figure~\ref{fig:completeK}.
\end {example}

In short, a complete pre-system of hyperboxes for a knot $K$ produces a multi-pointed diagram $\he^{K}$ for the knot, together with a way of relating the diagram $\he^K$ with the $z$ points on $K$ removed 
to the diagram $\he^K$ with the $w$ points on $K$ removed, via a certain sequence of strongly equivalent Heegaard diagrams and $\Theta$-cycles, plus some surface isotopies in the middle (moving the $w$ basepoints into the $z$ basepoints). Note that both of the diagrams that we relate in this fashion represent $Y$ itself. Observe also that the sequence of $\Theta$-cycles (and surface isotopies) induces a corresponding sequence of chain maps on Floer complexes as in Section~\ref{sec:hyperfloer}. These chain maps are all chain homotopy equivalences, see Example~\ref{ex:edge}. By composing them we obtain a chain homotopy equivalence between the initial and the final Floer complex. This is exactly the kind of structure that was used in \cite{IntSurg}, in the context of describing the Heegaard Floer homology of integer surgeries on knots. One should view complete pre-systems as a generalization of this structure to the case of links.

\begin{figure}
\begin{center}
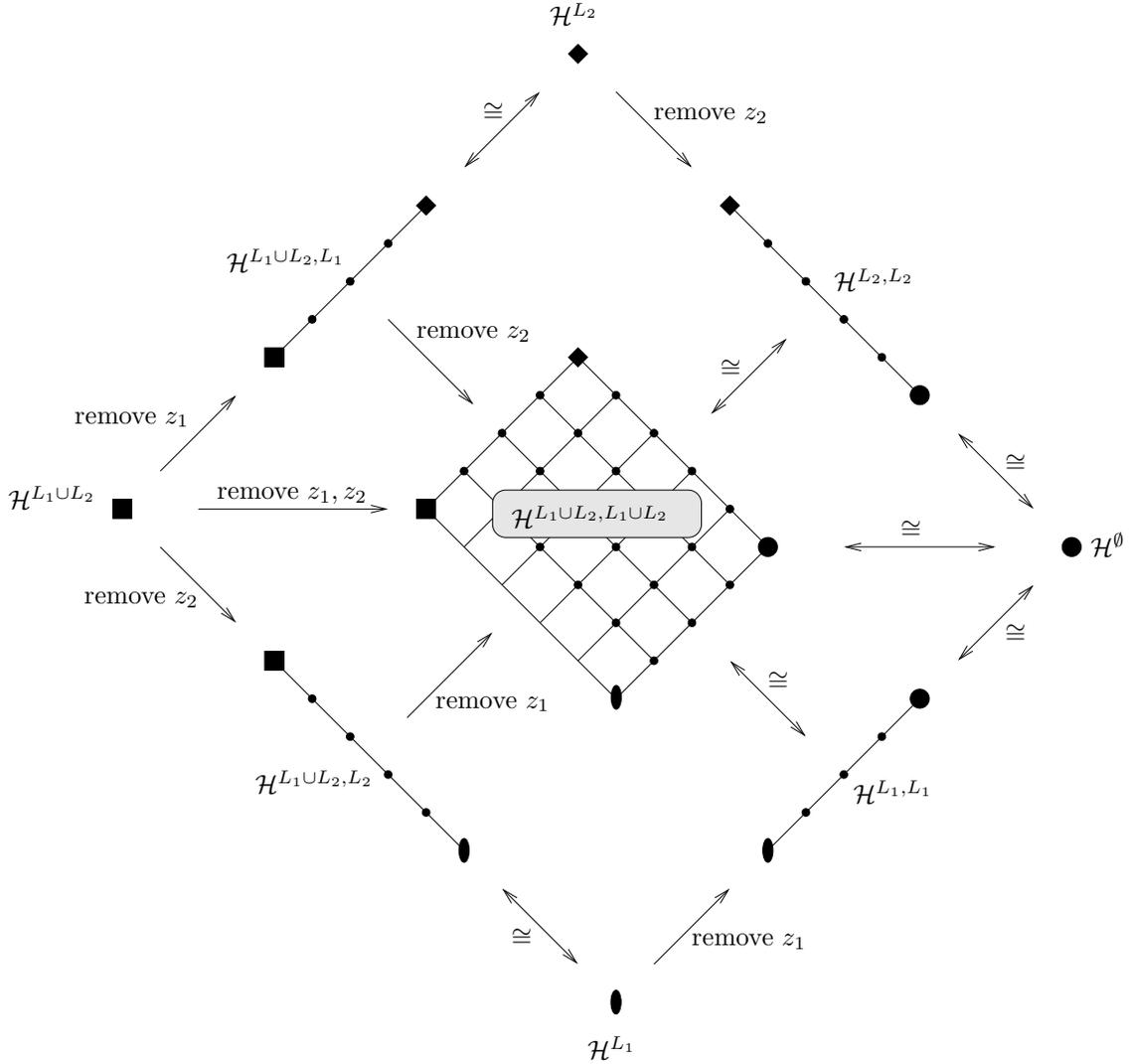
\end{center}
\caption {{\bf A part of a complete pre-system of hyperboxes for a link $\orL = L_1 \cup L_2$.} For simplicity, we assume that $\Hyper^L$ has exactly two basepoints $w_1$ and $z_1$ on $L_1$, and two basepoints $w_2$ and $z_2$ on $L_2$. We show here one quarter of the complete pre-system, consisting of the two-dimensional hyperbox $\Hyper^{L_1 \cup L_2,  L_1 \cup L_2}$, and all smaller hyperboxes related to it by compatibility conditions. 
}
\label{fig:completeL}
\end{figure}

\begin {example}
\label{ex:hyperboxes2}
When $\orL = L_1 \cup L_2$ is a link of two components, a complete pre-system of hyperboxes for $\orL$ consists of four zero-dimensional hyperboxes:
$$ \Hyper^{L_1 \cup L_2}, \Hyper^{L_1},  \Hyper^{L_2}, \Hyper^{\emptyset},$$
eight one-dimensional hyperboxes:
$$\Hyper^{L_1 \cup L_2,  L_1}, \Hyper^{L_1 \cup L_2,  -L_1}, \Hyper^{L_1 \cup L_2,  L_2}, \Hyper^{L_1 \cup L_2,  -L_2},$$
$$ \Hyper^{L_1, L_1}, \Hyper^{L_1, -L_1}, \Hyper^{L_2, L_2}, \Hyper^{L_2, -L_2},$$
and four two-dimensional hyperboxes:
$$\Hyper^{L_1 \cup L_2,  L_1 \cup L_2}, \Hyper^{L_1 \cup L_2,  -L_1 \cup L_2}, \Hyper^{L_1 \cup L_2,  L_1 \cup -L_2}, \Hyper^{L_1 \cup L_2,  -L_1 \cup -L_2}.$$

These hyperboxes are related by various compatibility conditions. Some of these conditions are illustrated in Figure~\ref{fig:completeL}.
\end {example}

In the case of a knot $K \subset S^3$, the complete pre-system only contains data for quasi-isomorphisms (in fact, chain homotopy equivalences) between various Heegaard Floer complexes of $Y=S^3$. At the level of homology, there is a unique $\ff[[U]]$-module isomorphism from $\HFm(Y) = \ff[[U]]$ to itself, namely the identity. However, in the case of links (or even for knots, but in more general homology spheres $Y$),  complete pre-systems give rise to quasi-isomorphisms between generalized Floer complexes of sublinks, whose homology can be complicated. Recall from Remark~\ref{rem:canonical} that, in general, the resulting isomorphisms on generalized Floer homology are not canonical. To ensure naturality, we need to control the paths traced by basepoints under the surface isotopies involved in the compatibility relations \eqref{eq:comp2}.

We start with the following:
\begin {definition}
\label {def:cs}
Let $\he= (\Sigma, \alphas, \betas, \ws, \zs)$ be a  multi-pointed Heegaard diagram for an oriented link $\orL \subset Y$, with a corresponding Heegaard splitting $Y = U_{\alpha} \cup_{\Sigma} U_{\beta}$. For each basepoint $w \in \ws \cap L$, we denote by $z(w) \in \zs$ the {\em successor of $w$ on $\orL$}, that is, the basepoint on the same component $\orL_j$ of $\orL$ as $w$, which appears just after $w$ as we go around $\orL_j$ according to its orientation. A {\em good set of trajectories}
$$ \cs = \{c_w | w \in \ws \cap L \}$$
for $\he$ consists of disjoint, smoothly embedded paths $$c_w: [0,1] \to \Sigma, \ \ \ c_w(0)=w, \ c_w(1) = z(w),$$ such that each path $c_w$ is homotopic (inside the handlebody $U_{\alpha}$) to the minimal oriented segment $l_w$ on $\orL$ going from $w$ to $z(w)$, by a homotopy whose interior avoids the link $L$. Further, we require that the image of each $c_w$ is disjoint from the free basepoints (i.e., the basepoints in $\ws \setminus L$).
\end {definition}

Observe that the alpha and the beta curves did not play an essential role in the definition above; we could have made the same definition by starting with a Heegaard splitting (and suitable basepoints) instead of a multi-pointed Heegaard diagram.

\begin{definition}
\label{def:surfacedata}
Let $\orL \subset Y$ be an oriented link. A {\em set of surface data} $(\Sigma, \ws, \zs, \cs)$ for $\orL$ consists of:
\begin{itemize}
\item an (oriented) embedded surface $\Sigma \subset Y$ that produces a Heegaard splitting $Y=U_{\alpha} \cup_{\Sigma} U_{\beta}$ with $\del U_{\alpha} = -\del U_{\beta} = \Sigma$; we assume $\Sigma$ to be transverse to $L$;
\item a decomposition of $L \cap \Sigma$ into alternating $w$ and $z$ basepoints; $\zs$ is the set of $z$ basepoints, and $\ws$ is the union of those basepoints and possibly some additional ones on $\Sigma \setminus L$; 
\item a good set of trajectories $\cs$ as in Definition~\ref{def:cs}.
\end{itemize}
\end{definition}

Consider now a complete pre-system of hyperboxes $\Hyper$ for $\orL \subset Y$. Note that all the hyperboxes $\Hyper^{\orL', \orM}$ in $\Hyper$ contain Heegaard diagrams with the same underlying Heegaard surface $\Sigma$. Let $\ws^{\orL', \orM}$ (resp. $\zs^{\orL', \orM}$) be the set of basepoints of type $w$ (resp. $z$) on the diagrams in the hyperbox $\Hyper^{\orL', \orM}$. In particular, we let $\ws^{L'} = \ws^{\orL', \emptyset}$ and $\zs^{L'} = \zs^{\orL', \emptyset}$. We also set $\ws = \ws^L$ and $\zs = \zs^L$.

Given a sublink $M \subseteq L$ with an orientation $\orM$, we let $M = M_+ \amalg M_-$, where $M_+$ (resp. $M_-$) consists of those components of $M$ that have the same (resp. opposite) orientation in $\orL$ as in $\orM$. Suppose $M$ is a sublink of some $L' \subseteq L$ (with the orientation $\orL'$ coming from $\orL$). Because of the compatibility condition \eqref{eq:comp1} (applied to $\orL'$ instead of $\orL$, and taking $M' = \emptyset$), we must have
$$ \ws^{\orL', \orM} =  (  \ws^{L'} \setminus M_-)   \cup ( \zs^{L'} \cap M_- )$$
and
$$ \zs^{\orL', \orM} = \zs^{L'} \cap (L' - M).$$

Thus, the basepoints on all hyperboxes are determined by those on the diagrams $\Hyper^{L'}$, for $L' \subseteq L$. To get further control on the basepoints (and the way they move under surface isotopies), 
consider a good set of trajectories $\cs = \{c_w\}$ for the initial diagram $\Hyper^L$ (which represents $\orL$). We say that the complete pre-system $\Hyper$ is {\em dependent} on the trajectory set $\cs$ if the following conditions are satisfied:
\begin {itemize}
\item For any $\orL' \subseteq \orL$, we have $\ws^{L'} = \ws$ and $\zs^{L'} = \zs \cap L'$, so that for any $M \subset L'$ with an orientation $\orM$,
$$ \ws^{\orL', \orM} =  (  \ws \setminus M_-)   \cup ( \zs \cap M_- ), \ \ \ \zs^{\orL', \orM} = \zs \cap (L' - M);$$
\item For any $M' \subseteq M \subseteq L' \subseteq L$, and any orientation $\orM$ of $M$, in the respective compatibility relation \eqref{eq:comp2}, which reads
\begin {equation}
\label {eq:comp22}
 \Hyper^{\orL'- M', \orM-M'} \cong \Hyper^{\orL',  \orM} (M', M)
 \end {equation}
we use a surface isotopy that moves each basepoint $w \in \ws \cap M'_-$ (appearing in the diagram on the left hand side of \eqref{eq:comp22}) to its successor $z(w) \in \zs \cap M'_-$ (appearing in the diagram on the right hand side of \eqref{eq:comp22}), exactly tracing the path $c_w$. Moreover, we require the surface isotopy to be supported in a small neighborhood of the path $c_w$, and in particular to fix all the basepoints in $$ \bigl(\ws \setminus M_-\bigr)  \cup \bigl(\zs \cap (M-M')_- \bigr) \cup \bigl(\zs \cap (L' - M) \bigr),$$ 
which appear on both sides of \eqref{eq:comp22}.
\end {itemize}

\begin {definition}
\label {def:complete}
A {\em complete system of hyperboxes} $(\Hyper, \cs)$ representing a link $\orL \subset Y$ consists of a complete pre-system of hyperboxes $\Hyper$ for $\orL$ together with a good set of trajectories $\cs$ for $\Hyper^L$ such that $\Hyper$ is dependent on $\cs$. (We usually drop $\cs$ from the notation, and refer to $\Hyper$ as a complete system of hyperboxes.)
\end {definition}

If $\Hyper$ is a complete system of hyperboxes for $\orL$ with a good set of trajectories $\cs$, and if $(\Sigma, \ws, \zs)$ is the multi-pointed Heegaard surface of the initial diagram $\Hyper^L$, we refer to $(\Sigma, \ws, \zs, \cs)$ as the {\em underlying surface data} for $\Hyper$. (Compare Definition~\ref{def:surfacedata}.)

\begin{definition}
\label{def:linkminimalsystem}
A complete system of hyperboxes $(\Hyper, \cs)$ for $\orL \subset Y$ is called {\em link-minimal} if $\Hyper^L$ (and hence all the Heegaard diagrams in $\Hyper$) are link-minimal, in the sense of Definition~\ref{def:linkminimal}.
\end{definition}

\begin{definition}
\label{def:mpsystem}
A complete system of hyperboxes $(\Hyper, \cs)$ for $\orL \subset Y$ is called {\em minimally-pointed} if the diagram $\Hyper^L$ is minimally pointed, in the sense of Definition~\ref{def:minimal}. (Note that this does not imply that the other diagrams in the system are minimally-pointed.)
\end{definition}

\begin{definition}
\label{def:purebsystem}
A complete system of hyperboxes $(\Hyper, \cs)$ for $\orL \subset Y$ is called {\em of pure $\beta$-type} if all its hyperboxes are pure $\beta$-hyperboxes, i.e. the corresponding bipartition maps $r_i$ take all the indices to $\beta$. (Hence, in such a system the $\alpha$ curves always stay fixed, and we have $\eps = \eps^\beta$ for any multi-index $\eps$.) 
\end{definition}

\subsection {Basic systems}
\label {sec:basic}

Let $\orL  \subset Y$ be a link. In this section we describe a special kind of minimally pointed complete system of pure $\beta$-type for $\orL$, which we call {\em basic}.

In a basic complete system, all hyperboxes of the form $\Hyper^{\orL', \orM}$ will be trivial, i.e. of size $(0, \dots, 0)$, when $\orM$ has the orientation induced from $\orL$ (or, equivalently, from $\orL'$). In particular, all diagrams $\he^{L'}=\he^{L', \emptyset}$ will simply be obtained from an initial diagram, $\he^{L}$, by deleting the $z$ basepoints on the components of $L - L'$. 

Let $\ell$ be the number of components of $L$. We choose the Heegaard diagram $\he^{L}$ to be basic (see Definition~\ref{def:hbasic}), that is, of genus $g$, with $g+ \ell - 1$ alpha curves and $g + \ell - 1$ beta curves, $\ell$ basepoints marked $w$, $\ell$ basepoints marked $z$, and such that the basepoints $w_i$ and $z_i, i =1 \dots, \ell$, lie on each side of a beta curve $\beta_i$, and are not separated by any alpha curves. Thus, the beta curves split the surface $\Sigma$ into $\ell$ components $\Sigma_1, \dots, \Sigma_{\ell}$, numbered such that $\Sigma_i$ contains the basepoints $w_i$ and $z_i$, and has both sides of $\beta_i$ as parts of its boundary.  We denote by $L_i$ the component of $L$ on which $w_i$ and $z_i$ lie. We construct a good set of trajectories $\cs =\{ c_{w_i} \}$ for $\he^L$ to consist of small paths $c_{w_i}$ from $w_i$ to $z_i$ that intersect $\beta_i$ once and do not intersect any of the other curves. We will take $(\Sigma, \{w_i\}, \{z_i\}, \cs)$ to be the underlying surface data for our complete system.

\begin{figure}
{
   \def\svgwidth{4in}
   %% Creator: Inkscape 1.0.2 (e86c8708, 2021-01-15), www.inkscape.org
%% Accompanies image file 'sigma1.pdf' (pdf, eps, ps)
\begingroup%
  \makeatletter%
  \providecommand\color[2][]{%
    \errmessage{(Inkscape) Color is used for the text in Inkscape, but the package 'color.sty' is not loaded}%
    \renewcommand\color[2][]{}%
  }%
  \providecommand\transparent[1]{%
    \errmessage{(Inkscape) Transparency is used (non-zero) for the text in Inkscape, but the package 'transparent.sty' is not loaded}%
    \renewcommand\transparent[1]{}%
  }%
  \providecommand\rotatebox[2]{#2}%
  \newcommand*\fsize{\dimexpr\f@size pt\relax}%
  \newcommand*\lineheight[1]{\fontsize{\fsize}{#1\fsize}\selectfont}%
  \ifx\svgwidth\undefined%
    \setlength{\unitlength}{484.87813058bp}%
    \ifx\svgscale\undefined%
      \relax%
    \else%
      \setlength{\unitlength}{\unitlength * \real{\svgscale}}%
    \fi%
  \else%
    \setlength{\unitlength}{\svgwidth}%
  \fi%
  \global\let\svgwidth\undefined%
  \global\let\svgscale\undefined%
  \makeatother%
  \begin{picture}(1,0.40867058)%
    \lineheight{1}%
    \setlength\tabcolsep{0pt}%
    \put(0,0){\includegraphics[width=\unitlength,page=1]{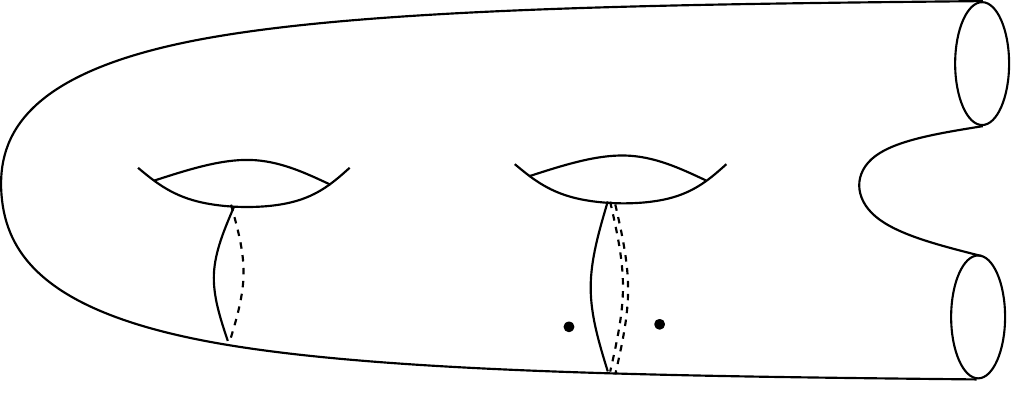}}%
    \put(0.5287993,0.07497971){\makebox(0,0)[lt]{\lineheight{1.25}\smash{\begin{tabular}[t]{l}$z_i$\end{tabular}}}}%
    \put(0.66303769,0.09915972){\makebox(0,0)[lt]{\lineheight{1.25}\smash{\begin{tabular}[t]{l}$w_i$\end{tabular}}}}%
    \put(0.58775721,0.00363533){\makebox(0,0)[lt]{\lineheight{1.25}\smash{\begin{tabular}[t]{l}$\beta_i$\end{tabular}}}}%
    \put(0.49136312,0.1112305){\makebox(0,0)[lt]{\lineheight{1.25}\smash{\begin{tabular}[t]{l}$\beta''_i$\end{tabular}}}}%
    \put(0,0){\includegraphics[width=\unitlength,page=2]{sigma1.pdf}}%
    \put(0.62499317,0.06581365){\makebox(0,0)[lt]{\lineheight{1.25}\smash{\begin{tabular}[t]{l}$c_{w_i}$\end{tabular}}}}%
  \end{picture}%
\endgroup%

}
\caption {{\bf The curves $\beta_i$ and $\beta_i''$.}
We show here a component $\Sigma_i$ of the complement of the beta curves in $\Sigma$.
}
\label{fig:Sigma_i}
\end{figure}

Let $\beta_i''$ be the curve obtained from $\beta_i$ by a small isotopy pushing $w_i$ into $z_i$, following the path $c_{w_i}$, as in Figure~\ref{fig:Sigma_i}. Note that, if we ignore $w_i$, then $\beta''_i$  can also be obtained from $\beta_i$ by handleslides supported in $\Sigma_i$ away from $z_i$, namely by handlesliding $\beta_i$ over all the other curves on the boundary of $\Sigma_i$; compare \cite[Proof of Proposition 7.1]{HolDisk}.

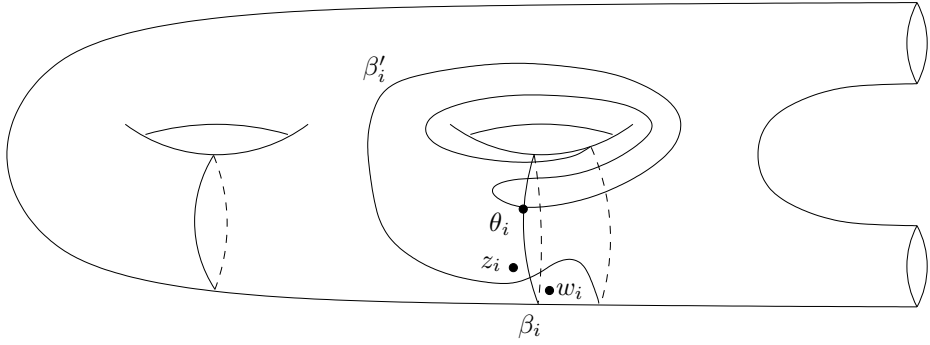
\begin{figure}
{
   \def\svgwidth{4in}
   %% Creator: Inkscape 1.0.2 (e86c8708, 2021-01-15), www.inkscape.org
%% Accompanies image file 'sigma2.pdf' (pdf, eps, ps)
\begingroup%
  \makeatletter%
  \providecommand\color[2][]{%
    \errmessage{(Inkscape) Color is used for the text in Inkscape, but the package 'color.sty' is not loaded}%
    \renewcommand\color[2][]{}%
  }%
  \providecommand\transparent[1]{%
    \errmessage{(Inkscape) Transparency is used (non-zero) for the text in Inkscape, but the package 'transparent.sty' is not loaded}%
    \renewcommand\transparent[1]{}%
  }%
  \providecommand\rotatebox[2]{#2}%
  \newcommand*\fsize{\dimexpr\f@size pt\relax}%
  \newcommand*\lineheight[1]{\fontsize{\fsize}{#1\fsize}\selectfont}%
  \ifx\svgwidth\undefined%
    \setlength{\unitlength}{484.87813058bp}%
    \ifx\svgscale\undefined%
      \relax%
    \else%
      \setlength{\unitlength}{\unitlength * \real{\svgscale}}%
    \fi%
  \else%
    \setlength{\unitlength}{\svgwidth}%
  \fi%
  \global\let\svgwidth\undefined%
  \global\let\svgscale\undefined%
  \makeatother%
  \begin{picture}(1,0.40867058)%
    \lineheight{1}%
    \setlength\tabcolsep{0pt}%
    \put(0,0){\includegraphics[width=\unitlength,page=1]{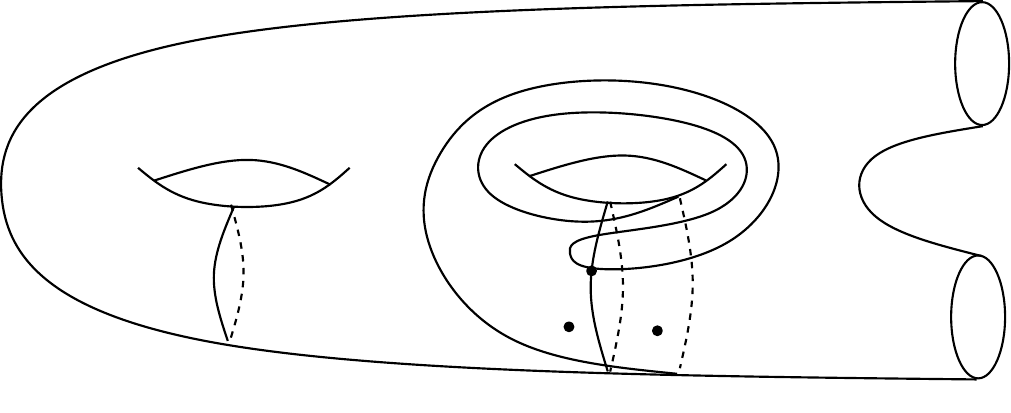}}%
    \put(0.5287993,0.07497971){\makebox(0,0)[lt]{\lineheight{1.25}\smash{\begin{tabular}[t]{l}$z_i$\end{tabular}}}}%
    \put(0.63003447,0.09469888){\makebox(0,0)[lt]{\lineheight{1.25}\smash{\begin{tabular}[t]{l}$w_i$\end{tabular}}}}%
    \put(0.5502981,0.11990965){\makebox(0,0)[lt]{\lineheight{1.25}\smash{\begin{tabular}[t]{l}$\theta_i$\end{tabular}}}}%
    \put(0.58775721,0.00363533){\makebox(0,0)[lt]{\lineheight{1.25}\smash{\begin{tabular}[t]{l}$\beta_i$\end{tabular}}}}%
    \put(0.40205953,0.2872198){\makebox(0,0)[lt]{\lineheight{1.25}\smash{\begin{tabular}[t]{l}$\beta'_i$\end{tabular}}}}%
  \end{picture}%
\endgroup%

}
\caption {{\bf The curves $\beta_i$ and $\beta_i'$.}
The diagram formed by them is admissible.
}
\label{fig:sigma}
\end{figure}

We would like to construct a complete system using the curves $\beta_i$ and $\beta_i''$. Of course, if the isotopy from $w_i$ to $z_i$ is supported just in a small neighborhood of $c_{w_i}$, then the curves $\beta_i$ and $\beta_i''$ are not transverse (they have a whole arc in common away from that neighborhood), so cannot be part of the same Heegaard diagram. This issue persists even if we slightly isotope $\beta_i''$ further so that its intersection with $\beta_i$ is transverse and consists of two points, as in Figure~\ref{fig:Sigma_i}. Indeed, let $B_i^1$ be the bigon between $\beta_i$ and $\beta_i''$ containing $z_i$ and $B_i^2$ the other bigon. Then $\Sigma_i - B_i^1 + B_i^2$ is a periodic domain in $\Sigma$ with only nonnegative multiplicities, and with zero multiplicity at the basepoint $z_i$. This is  a source of problems if one tries to construct a Heegaard diagram using the curves $\beta_i$ and $\beta_i''$ and the basepoint $z_i$, because such diagrams are not admissible. In order to fix this problem, we introduce an intermediate curve $\beta'_i$ as in Figure~\ref{fig:sigma}. Then there are no periodic domains as above between $\beta_i$ and $\beta_i'$, nor between $\beta_i'$ and $\beta_i''$.  

Let us now describe the hyperbox $\Hyper^{\orL, - \orL}$, which is the biggest hyperbox in our basic system. As mentioned above, the hyperbox $\Hyper^{\orL, -\orL}$ is a pure $\beta$-hyperbox. Its punctures are $z_i$ for $i=1, \dots, \ell$, but playing the role of $w$'s in the definition of a multi-pointed Heegaard diagram. The hyperbox has dimension $\ell$ and size $\dd = (2,2, \dots, 2)$.  For $\eps \in \E(\dd)=\{0,1,2\}^\ell$, the collection of curves $\betas^\eps$ is given by an aproximation:
$$ \beta^\eps_i \approx \begin {cases}
\beta_i & \text{ if }  \eps_i = 0 \\
\beta_i' & \text{ if } \eps_i = 1\\
\beta_i'' & \text{ if } \eps_i = 2,
\end {cases} $$ 
for $i \leq \ell$, and $ \beta^\eps_i \approx \beta_i$
for $i > \ell$. We also arrange so that $\beta_i^\eps \approx \beta_i^{\eps'}$ for any $\eps \neq \eps'$.

Note that with this choice of $\betas^\eps$, when two multi-indices $\eps$ and $\eps'$ are neighbors, we  never see a pair of curves that approximate $\beta_i$ resp. $\beta_i''$ (for the same $i$) in the Heegaard diagram $(\Sigma, \betas^{\eps}, \betas^{\eps'}, \ws, \zs)$. Hence, the admissibility hypothesis in the definition of the $\beta$-hyperbox is satisfied.  Further, there is a natural choice for the chain $\Theta_{\eps, \eps'}$ when $(\eps, \eps')$ is an edge in the hyperbox: namely, the respective intersection points of maximal degree. For example, if in the Heegaard diagram $(\Sigma, \betas^{\eps}, \betas^{\eps'}, \ws, \zs)$ we see the curves $\beta_i$ and $\beta_i'$, then the point marked $\theta_i$ in Figure~\ref{fig:sigma} would be part of the corresponding $\Theta$-chain. When $\eps < \eps'$ are neighbors but $\|\eps' - \eps\| \geq 2$, we set $\Theta_{\eps, \eps'} = 0$.
 
The fact that the $\Theta$-chains satisfy the compatibility relations \eqref{eq:compatibility} has a proof similar to those of \cite[Lemma 9.7]{HolDisk} and \cite[Lemma 4.3]{BrDCov}.

For future reference, we denote the $\beta$ curves in the $\Hyper^{\orL, - \orL}$ hyperbox by $$ \tilde \beta^\eps_i =\beta^\eps_i, \ \ \eps \in \{0,1,2\}^\ell,\  i \in \{1, \dots, g+\ell - 1\}.$$

Let us now describe an arbitrary hyperbox $\Hyper^{\orL', \orM}$ in the basic system. Let $M = M_+ \amalg M_-$ as in Section~\ref{sec:complete}, where $M_+$ (resp. $M_-$) consists of those components of $M$ that have the same (resp. opposite) orientation in $\orL$ as in $\orM$. Denote by $m$ (resp. $m_+, m_-$) the number of components in $M$ (resp. $M_+, M_-$). Order the components of $M$ according to their ordering as components of $L$:
$$ M = L_{i_1} \cup \dots \cup L_{i_{m}}, \ \ i_1 < \dots < i_{m}.$$

The hyperbox $\Hyper^{\orL', \orM}$ has size $\dd^{\orM} = (d^{\orM}_1, \dots, d^{\orM}_m)$, where 
$$d^{\orM}_j = \begin {cases} 0 & \text{if } L_{i_j} \subseteq M_+ \\  2 & \text {if } L_{i_j} \subseteq M_-. \end {cases}$$

Its diagrams all lie on the surface $\Sigma$ with punctures $\ws^{\orL', \orM}, \zs^{\orL', \orM}$ given by
$$ \ws^{\orL', \orM} = \{w_i | L_i \subseteq L - M_- \} \cup \{z_i | L_i \subseteq M_-\}$$
and
$$ \zs^{\orL', \orM} = \{z_i | L_i \subseteq L' - M \}.$$

Consider the injective map
$$ \lambda^{\orM}  : \E(\dd^{\orM}) \to \E(\dd^{-\orL}) = \{0,1,2\}^\ell,$$
given by
$$\Bigl( \lambda^{\orM}(\eps_1, \dots, \eps_{m})\Bigr)_i=
\begin {cases} 
\eps_j & \text{if } i=i_j \text{ for some } j \\   
0 & \text{otherwise.}
\end {cases}
$$

For the beta curves $\beta^\eps_i$ in the hyperbox $\Hyper^{\orL', \orM}$, we choose 
\begin{equation}
\label{eq:bt}
 \beta^\eps_i = \tilde \beta^{\lambda^{\orM}(\eps)}_i, \ \ \eps \in \E(\dd^{\orM}), \ i \in \{1, \dots, g+\ell - 1\}.
 \end{equation}

Finally, we let the $\Theta$-cycles on the edges be the respective intersection points of maximal degree, and we let the rest of the $\Theta$-chain elements to be zero. 

This completes the description of the hyperboxes $\Hyper^{\orL', \orM}$. (Note that when $M_- = \emptyset$, we indeed get a trivial hyperbox consisting of the initial Heegaard diagram $\Hyper^L$ with the points $z_i$ on $(L-L') \cup M$ removed.) The verification of the conditions in the definition of a complete system is an easy exercise.

\begin {definition}
\label {def:basic}
A complete system of hyperboxes $\Hyper$ representing the link $\orL$ is called {\em basic} if it is constructed as above. If this is the case, we refer to $\Hyper$ more simply as a {\em basic system} for the link $\orL$. 
\end {definition}

 To review, every basic system for $\orL$ is associated to a particular $2\ell$-pointed Heegaard diagram $(\Sigma, \alpha_1, \dots, \alpha_{g+\ell-1}, \beta_1 , \dots, \beta_{g+\ell-1}, w_1, \dots, w_\ell; z_1, \dots, z_\ell)$ with the property that, for each $i=1, \dots, \ell$, the base points $w_i$ and $z_i$ lie on each side of the curve $\beta_i$ and are not separated by any $\alpha$ curves. In addition, the construction of the basic system requires several choices (the curves $\beta'_i, \beta''_i$, and their approximations). For simplicity, we typically forget those choices and just say that the basic system is associated to a special Heegaard diagram as above.

\begin {remark}
Another complete system of hyperboxes for a link is described in Section~\ref{sec:completegrid}, using grid diagrams.
\end {remark}

\subsection{A more general construction}
While Section~\ref{sec:basic} was focused on basic systems, the same ideas can be used to construct complete systems from any Heegaard diagram and good set of trajectories.

\begin{proposition}
\label{prop:gencon}
Let $\Hyper^L=(\Sigma, \alphas, \betas, \ws, \zs)$ be a multi-pointed Heegaard diagram for a link $\orL \subset Y$, and let $\cs$ be a good set of trajectories for $\Hyper^L$. Then, there exists a complete system of hyperboxes $\Hyper$, which is of pure $\beta$-type, dependent on $\cs$, and has $\Hyper^L$ as its initial diagram.
\end{proposition}

\begin{proof}
We will construct $\Hyper$ so that its hyperboxes $\Hyper^{\orL', \orM}$ are of the same size as the ones defined for a basic system in Section~\ref{sec:basic}. Also, since $\Hyper$ is meant to be of pure $\beta$-type, the collections of alpha curves will be fixed to be the same collection $\alphas$ from $\Hyper^L$, in all the hyperboxes. 

We begin by describing the biggest hyperbox $\Hyper^{\orL, -\orL}$, which is of size $\dd = (2,2, \dots, 2)$.  Recall that each trajectory $c_w$ in $\cs$ gives a surface isotopy supported in a neighborhood of $c_w$, taking $w$ to its successor $z(w)$. For $\eps \in \E(\dd)=\{0,1,2\}^\ell$, let $\tau_{\eps}$ be the composition of all surface isotopies corresponding to $c_w$ for basepoints $w$ on link components $L_i$ with $\eps_i \neq 0$. 

If none of the values $\eps_i$ are $1$, we choose the collection of beta curves $\betas^{\eps}$ to be
$$ \betas^\eps := \tau_{\eps}(\betas).$$ 
Thus, when some value $\eps_i$ changes from $0$ to $2$, the beta curves change by an isotopy moving some basepoints, similarly to how $\beta_i''$ was obtained from $\beta_i$ in Section~\ref{sec:basic}.

We construct the remaining collections $\betas^\eps$ by induction on the number of values $1$ among the $\eps_i$. Let $\bar \eps \in \{0,2\}^\ell$ be obtained from $\eps$ by changing all values $1$ among $\eps_i$ to $2$. Then, we define $\betas^\eps$ by winding the curves $\betas^{\bar \eps} = \tau_{\eps}(\betas)$ sufficiently to achieve admissibility with the already defined collections of beta curves, in each hypercube. This is done as in \cite[Sections 5 and 8.4.2]{HolDisk}, and is similar to how we obtained the curve $\beta_i'$ from $\beta_i''$ in Section~\ref{sec:basic}.

After the beta curves are constructed, we fill in the hyperboxes in $\Hyper^{\orL, -\orL'}$ with $\Theta$-chain elements  using Lemma~\ref{lemma:filling}. We then define the other hyperboxes $\Hyper^{\orL', \orM}$ based on what we have in $\Hyper^{\orL, -\orL'}$, using the same formula \eqref{eq:bt} for the beta curves as in the basic case. We also let the $\Theta$-chain elements to be induced from those in $\Hyper^{\orL, -\orL'}$: the element for the pair $(\eps, \eps')$ in $\Hyper^{\orL', \orM}$ should be the same as the one for $(\lambda^{\orM}(\eps), \lambda^{\orM}(\eps'))$ in $\Hyper^{\orL, -\orL'}$.
\end{proof}

\subsection {Moves on complete systems of hyperboxes}
\label {sec:moves}

In Section~\ref{sec:movesh} we defined several moves on hyperboxes of strongly equivalent Heegaard diagrams. We can define an analogous list of {\em system Heegaard moves} between complete systems of hyperboxes.

First, note that {\em 3-manifold isotopies, index one/two stabilizations, free index zero/three stabilizations, global shifts} (and all their inverses) have straightforward extensions to complete systems of hyperboxes. A certain move of one of these types on a complete system consists of applying that type of move to all the hyperboxes in the system, in a way compatible with restrictions. The good set of trajectories should be taken into the corresponding one by the respective moves. The (index one/two, or free index zero/three) stabilizations are required to be done away from the trajectories $c_w$ in the good set, so that these trajectories are preserved. 

The three other system Heegaard moves in the list are elementary enlargements / contractions, index zero/three link stabilizations (and destabilizations), and trajectory $\alpha$-slides, which require more discussion.

An {\em elementary enlargement} of a complete system $\Hyper$ consists of picking a  component $M_0$ of $L$, with an orientation $\orM_0$, and doing compatible elementary enlargements of those hyperboxes $\Hyper^{\orL', \orM}$ with $\orM_0 \subseteq \orM$ (and the orientation on $\orM_0$ is the one induced from $\orM$); these elementary enlargements are all done along the coordinate axis corresponding to $M_0$. The good set of trajectories is unchanged. An {\em elementary contraction} is the inverse of an elementary enlargement.

An {\em index zero/three link stabilization} of a complete system $\Hyper$ is as follows. We do an index zero/three link stabilization of the initial diagram $\Hyper^L$ resulting in two additional basepoints $w'$ and $z'$ near an old basepoint $z$ on a component $L_i \subseteq L$, see Figure~\ref{fig:Stab03}. Let $\Sigma'$ be the new Heegaard surface, and $\ws', \zs'$ the new collections of basepoints. We let $c_{w'}$ be a short trajectory (inside the disk bounded by $\alpha'$) joining $w'$ to $z'$. Adding $c_{w'}$ to the good set of trajectories $\cs$ for $\Hyper$, we obtain a new trajectory set $\cs'$. We take $(\Sigma, \ws', \zs', \cs')$ to be the underlying surface data for the new stabilized complete system $\Hyper'$, which we now describe. Let $M \subseteq L' \subseteq L$ be sublinks, with $M = M_+ \amalg M_-$ having an orientation $\orM$ (that coincides with the orientation of $\orL$ exactly on $M_+$). To get from the hyperboxes in $\Hyper$ to the corresponding hyperboxes in $\Hyper'$, we do the following.

 If $L_i \subseteq L' - M_-$, we change the hyperbox $\Hyper^{\orL', \orM}$ by an index zero/three link stabilization performed in the same place as in $\Hyper^L$.
 
 If $L_i \subseteq L- L'$, note that the $z$ and $z'$ basepoints should disappear. Therefore, we change $\Hyper^{\orL', \orM}$ by a free index zero/three stabilization at the same location (followed by a 3-manifold isotopy to get from $\Sigma$ to $\Sigma'$).

 \begin{figure}
\begin{center}
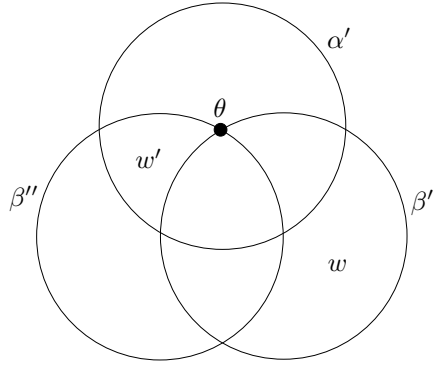
\end{center}
\caption {{\bf Index zero/three link stabilizations of complete systems.} Compared to the right hand side of Figure~\ref{fig:Stab03}, we deleted $w'$ and relabeled $z$ and $z'$ as $w$ and $w'$, respectively. The curve $\beta''$ can be obtained from $\beta'$ by handlesliding it over other beta curves.
}
\label{fig:StabTriple}
\end{figure}

 If $L_i \subseteq M_-$, note that the $w'$ basepoint disappears, and we relabel $z$ and $z'$ as $w$ and $w'$, respectively. We first construct a hyperbox $\check \Hyper^{\orL', \orM}$ from $\Hyper^{\orL', \orM}$ by taking (in each Heegaard diagram) the respective connected sum with a sphere, and adding the two new curves $\alpha'$ and $\beta'$ encircling $w$ and $w'$, as in Figure~\ref{fig:StabTriple}. The corresponding hyperbox $ (\Hyper')^{\orL', \orM}$ in the new system $\Hyper'$ is obtained from $\check \Hyper^{\orL', \orM}$ by increasing the length of its side in the $L_i$ direction by one, and adding on (at the end of each segment in the $L_i$ direction) new diagrams in which the curve $\beta'$ is replaced by $\beta''$ as in Figure~\ref{fig:StabTriple}. (Observe that replacing $\beta'$ with $\beta''$ is the strong equivalence discussed in Section~\ref{sec:final}, with different notation.) The $\Theta$-chain elements on the new edges in the  $L_i$ direction contain the point $\theta$ in the figure. After this modification, the sub-hyperbox $(\Hyper')^{\orL', \orM}(L_i, M)$ is  surface isotopic to the free index zero/three stabilization of $\Hyper^{\orL' - L_i, \orM - L_i}$, as it should be according to the compatibility relation \eqref{eq:comp2}. Note that the change from $\check \Hyper^{\orL', \orM}$ to $(\Hyper')^{\orL', \orM}$ can be realized as the composition of an elementary enlargement and a global shift (of hyperboxes). We do these moves in a compatible way on all  hyperboxes corresponding to pairs $(\orL', \orM)$ with $L_i \subseteq M_-$. This completes the description of the index zero/three link stabilization of $\Hyper$.
 
\begin{figure}
 \hskip1.5cm
{
   \def\svgwidth{5.7in}
   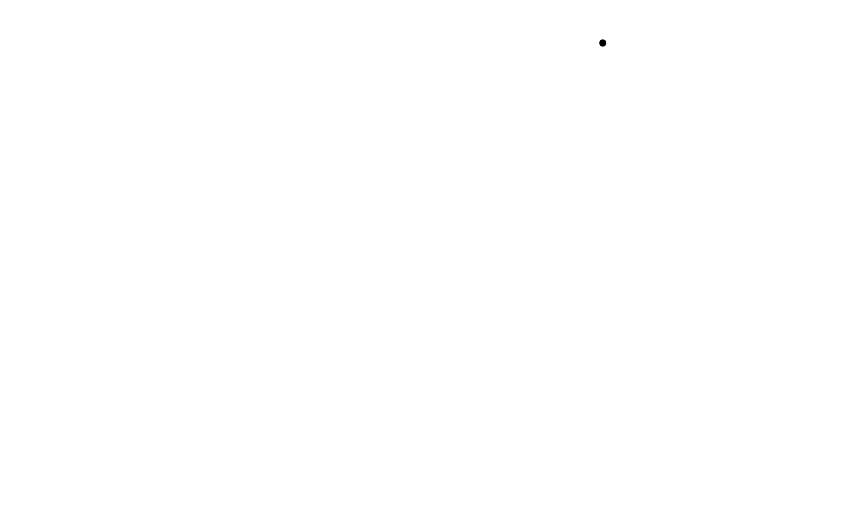
}
\caption{A trajectory $\alpha$-slide.}
\label{fig:aslide}
\end{figure}

 Finally, we have one more system Heegaard move, the {\em trajectory $\alpha$-slide}, which did not have an analogue in the list of moves on hyperboxes. This involves changing the good set of trajectories $\cs$. It is a move that we only perform on complete systems of pure $\beta$-type, i.e. those where the alpha curves are fixed;  see Definition~\ref{def:purebsystem}. Suppose $\Hyper$ is such a complete system, and let $\tau$ be an embedded arc connecting a basepoint $z_i$ to an alpha curve $\alpha_j$, and not intersecting any alpha curves in its interior. Observe that the trajectory $c_{w_i}$ does not intersect alpha curves either---otherwise these would be moved by the isotopy along $c_{w_i}$, but we assume them to be fixed. On the other hand, both $c_{w_i}$ and $\tau$ may intersect various beta curves, from different diagrams in the complete system. See Figure~\ref{fig:aslide} (a). 
 
The trajectory $\alpha$-slide is a move that transforms $\Hyper$ into a new complete system $\Hyper'$, where the trajectory $c_{w_i}$ is replaced by $c'_{w_i}$, the result of a handleslide of $c_{w_i}$ over $\alpha_j$ along the arc $\tau$, as in Figure~\ref{fig:aslide} (b). To motivate the description of $\Hyper'$, consider the relations 
\[
 \Hyper^{\orL'- M', \orM-M'} \cong \Hyper^{\orL',  \orM} (M', M)
 \]
from~\eqref{eq:comp2}, where the surface isotopy goes along $c_{w_i}$; this happens when $w_i \in M'_-$. Let us denote this surface isotopy by $\psi$. Then, if on the left hand side of the above equation we have a diagram with curves $\beta_k$, on the right hand side we will have the curves $\psi(\beta_k)$, as in Figure~\ref{fig:aslide} (c). Now, in the new complete system $\Hyper'$, we would like to have on the right hand side the curves $\psi'(\beta_k)$ instead, where $\psi'$ is the isotopy following $c'_{w_i}$, as in Figure~\ref{fig:aslide} (d). To achieve this, given any hyperbox $\Hyper^{\orL', \orM}$ with $w_i \in M_-$, we increase its side length corresponding to the component $L_i$ (the one containing $w_i$) by $2$, by letting $\psi(\beta_k)$ be followed by a new curve $\beta_k'$ and then by $\psi'(\beta_k)$; the intermediate curves $\beta_k'$ are chosen by twisting $\beta_k$ sufficiently to ensure admissibility. We then fill in the new parts of the hyperboxes using the procedure from Lemma~\ref{lemma:filling}. The result is a new complete system of pure $\beta$-type, which is our desired $\Hyper'$.
 
With the list of system Heegaard moves in place, we can now state the analogue of Lemma~\ref{lemma:moves} in this context:

\begin {proposition}
\label {prop:moves}
Let $\Hyper$ and $\Hyper'$ be two complete systems of hyperboxes representing the same oriented link $\orL \subset Y$.

$(a)$ The complete systems $\Hyper, \Hyper'$ can be related by a sequence of system Heegaard moves.

$(b)$ Suppose that $\Hyper$ is a link-minimal system, and that the diagram $(\Hyper')^L$ in $\Hyper'$ has $m$ points of type $z$, and hence $m-\ell$ subsidiary $w$ basepoints. Then, $\Hyper'$ can be obtained from $\Hyper$ by a sequence of system Heegaard moves that includes exactly $m-\ell$ index zero/three link stabilizations, with each of these stabilizations introducing a subsidiary basepoint.
\end {proposition}

Recall that each complete system $\Hyper$ has some underlying surface data $S=(\Sigma, \ws, \zs, \cs)$; compare Definition~\ref{def:surfacedata}. To prove Proposition~\ref{prop:moves}, we will relate the complete systems by first relating their underlying surface data. With the exceptions of global shifts and elementary enlargements / contractions (which do not change the surface data), each system move has a corresponding move at the level of the surface data itself. Thus, we distinguish the following types of {\em moves of surface data}:  3-manifold isotopies, index one/two stabilizations, index zero/three stabilizations (free or linked), and trajectory $\alpha$-slides. (In the case of trajectory $\alpha$-slides, to define these without reference to alpha curves, we simply require one of the trajectories in $\cs$ to slide over a simple closed curve that bounds a disk in the handlebody $U_{\alpha}$.)

\begin {lemma}
\label {lem:movesS}
Let $S$ and $S'$ be two sets of surface data for the same oriented link $\ell$-component $\orL \subset Y$. 
 
$(a)$ The sets $S$ and $S'$ can be related by a sequence of moves of surface data.

$(b)$ Suppose that $S$ is link-minimal (that is, it contains exactly $\ell$ basepoints of type $z$), and that $S'$  has $m$ points of type $z$, and hence $m-\ell$ subsidiary $w$ basepoints. Then, $S'$ can be obtained from $S$ by a sequence of moves of surface data that includes exactly $m-\ell$ index zero/three link stabilizations.
\end {lemma}

\begin{proof}
Let us first consider the case when both $S= (\Sigma, \ws, \zs, \cs)$ and $S'= (\Sigma', \ws', \zs', \cs')$ are minimally pointed; that is, they each have only $\ell$ basepoints of type $w$ and $\ell$ basepoints of type $z$.

In this case, after a 3-manifold isotopy, we can assume $S$ and $S'$ have the same basepoints $\ws=\ws'$ and $\zs=\zs'$, and in fact that the corresponding surfaces $\Sigma$ and $\Sigma'$ agree in a neighborhood of each basepoint (intersecting that neighborhood in the same small disk). Moreover, we can make the trajectories $c_w$ agree with $c_w'$ in small neighborhoods of the basepoints.

After this, to relate $S$ and $S'$, we use a Morse theoretic argument. We fix a Riemannian metric on $Y$. For each component $L_i$ of $L$ (with basepoints $w_i$ and $z_i$), we also fix a standard self-indexing Morse function $f_i$ in a small tubular neighborhood $\nu(L_i)$ of $L_i$. We ask for $f_i$ to have exactly two critical points, both on $L_i$, one of index $0$ and one of index $3$,  such that $L_i$ consists of two gradient flow trajectories between these two critical points, and  
$$f_i^{-1}(3/2)= \Sigma \cap \nu(L_i) = \Sigma' \cap \nu(L_i)$$
so that the alpha handlebody (from either $S$ or $S'$) intersects $\nu(L_i)$ in $f_i^{-1}[0, 3/2]$, and the beta handlebody intersects $\nu(L_i)$ in $f_i^{-1}[3/2, 3].$  

Furthermore, we can assume that $c_{w_i} \cap \nu(L_i) = c'_{w_i} \cap \nu(L_i)$ consists of two small arcs on $\Sigma$: an arc $a_i$ from $w_i$ to some point $w_i^*$ on $\del \nu(L_i)$, and another arc $b_i$ from some point $z_i^*$ on $\del \nu(L_i)$ to $z_i$. We let $c_i \subset c_{w_i}$ be the remaining arc from $w_i^*$ to $z_i^*$ along the trajectory $c_{w_i}$, after we remove $a_i$ and $b_i$. We define $c'_i \subset c'_{w_i}$ similarly.

We choose a pushoff $\gamma_i$ of the arc $L_i \cap f_i^{-1}[0, 3/2]$, so that $\gamma_i$ goes along the boundary $\del \nu(L_i)$, from $w_i^*$ to $z_i^*$. See Figure~\ref{fig:Morse}.
\begin{figure}
 \hskip1.5cm
{
   \def\svgwidth{5in}
   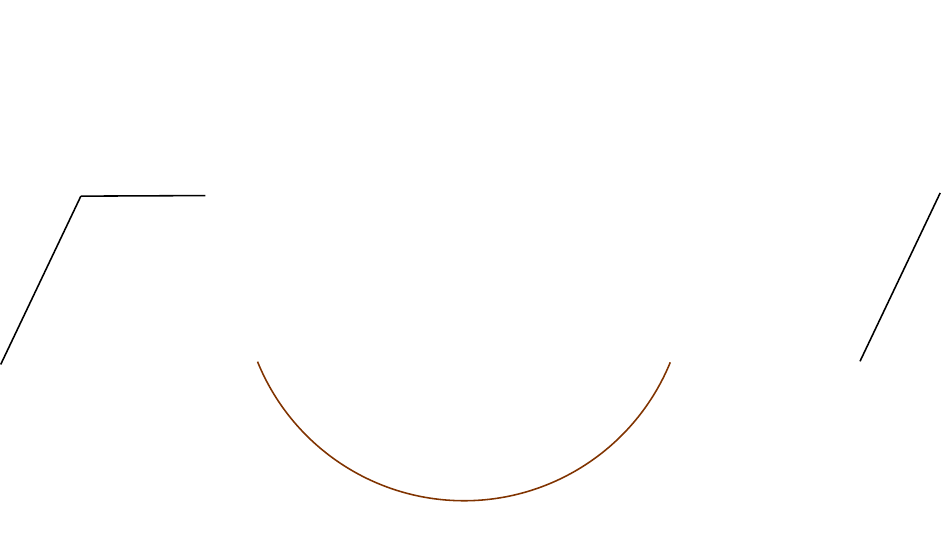
}
\caption{{\bf Constructing a Morse-Smale function.} We show here a neighborhood of a link component $L_i$, together with the path $c_i$. The arrows represent gradient flow lines taking $\gamma_i$ into $c_i$.}
\label{fig:Morse}
\end{figure}

We extend the functions $f_i$ to a single self-indexing Morse-Smale function $f$ on all of $Y$, by introducing new index $1$ and $2$ critical points (but no additional index $0$ or $3$ critical points) so that $\Sigma=f^{-1}(3/2)$. We ask that the arc $\gamma_i$ is taken by the forward gradient flow to the arc $c_i$ on $\Sigma$, from $w_i^*$ to $z_i^*$. To see that we can choose $f$ this way, note that, by assumption, $\gamma_i$ is homotopic to $c_i$ inside the alpha handlebody, and therefore, by Dehn's lemma, the union of $\gamma_i$ and $c_i$ bounds an embedded disk $D_i$ in that handlebody. We extend the functions $f_i$ from $\nu(L_i)$ to a neighborhood of that disk, so that the flow takes $\gamma_i$ to $c_i$, then extend them to the rest of the alpha handlebody (so that on $\Sigma$ it takes the value $3/2$), and finally to the beta handlebody. 

Next, we apply the same procedure to $\Sigma'$, obtaining a Morse-Smale function $g$ with $\Sigma'=g^{-1}(3/2)$, such that $c'_i$ is the image of $\gamma_i$ under the forward gradient flow.

We now relate the two functions $f$ and $g$ by a generic smooth homotopy $(F_t)_{0\leq t \leq 1}$,  with $F_0=f$, $F_1=g$, keeping $F_t$ fixed in the neighborhoods $\nu(L_i)$. Following the level set at $3/2$ along this homotopy, we can relate the surfaces $\Sigma$ and $\Sigma'$ by $3$-manifold isotopies and  index $1$/$2$ stabilizations and destabilizations; compare \cite[Proposition 2.2]{HolDisk} or \cite[Proposition 2.15]{Juhasz}. As for the trajectories $c_{w_i}$ and $c'_{w_i}$, they remain fixed near the basepoints throughout the homotopy, while the arc $c_i$ changes into $c'_i$ by following the images of $\gamma_i$ through the homotopy. In this generic one-parameter family of functions, at finitely many points it may happen that the forward flow takes $\gamma_i$ into an index $1$ critical point. Then, $c_i$ changes not by an isotopy on $\Sigma$, but rather by a trajectory $\alpha$-slide (over the intersection of $\Sigma$ with the ascending manifold of the index $1$ critical point). 

Thus, we have succeeded in relating any two minimally-pointed sets of surface data. If we have arbitrary surface data $S'$, we can find some minimally-pointed surface data $S$ and a sequence of 3-manifold isotopies and index zero/three (free and linked) stabilizations, that take $S$ into $S'$. Indeed, we can position these (de-)stabilizations around the trajectories $c_w$ and the free basepoints that we want to eliminate from $S'$, until we get to the minimally pointed setting. This completes the proof of part (a). For part (b), note that in the above process we used the minimal number of index zero/three link stabilizations.
\end{proof}

\begin {proof}[Proof of Proposition~\ref{prop:moves}]
First, let us consider the case when $\Hyper$ and $\Hyper'$ have the same underlying surface data $S=(\Sigma, \ws, \zs, \cs)$. If so, elementary enlargements and contractions can be used to make the sizes (and bipartition maps) of all corresponding hyperboxes to agree. Then, we use global shifts to relate the resulting complete systems, as we did in Lemma~\ref{lemma:moves} for hyperboxes. 

In view of the above observation, instead of relating two complete systems $\Hyper$ and $\Hyper'$ by moves, it suffices to consider their surface data $S$ and $S'$, use Lemma~\ref{lem:movesS} to relate these by a sequence of intermediate data  
\begin{equation}
\label{eq:fromStoS}
 S=S_0, \ S_1, \ S_2, \dots, S_n = S'
 \end{equation}
and, for each $k=0, \dots, n-1$, find \emph{some} complete system $\Hyper_k$ with data $S_k$, related by some system move to a complete system $\Hyper'_k$ with data $S_{k+1}$. We would then be able to relate $\Hyper$ to $\Hyper_0$ (because they have surface data $S_0$), each $\Hyper_k$ to $\Hyper'_k$ and to $\Hyper_{k+1}$ (as the latter two have the same surface data $S_{k+1}$), and finally $\Hyper'_{n-1}$ with $\Hyper'$. The combination of these system moves will enable us to get from $\Hyper$ to $\Hyper'$.

It remains to see that every move of surface data can be upgraded to a corresponding move of complete systems. We have already seen in Proposition~\ref{prop:gencon} that for any surface data, we can find a complete system of pure $\beta$-type having that data. If our move is a $3$-manifold isotopy, we just follow the same isotopy to change the complete system. If our move is a stabilization (of any kind), let $S$ be the data before the stabilization and $S'$ the one after. Up to isotopy, we can assume that the stabilization happens in a neighborhood of some point $p$, and therefore we can choose a complete system $\Hyper$ with underlying data $S$, such that we can apply the same stabilization to it (i.e., the curves involved in the diagrams in $\Hyper$ do not go through $p$); the result is a system $\Hyper'$ for $S'$. 

Finally, if our move is a trajectory $\alpha$-slide over some curve $\alpha_0$, we pick a complete system $\Hyper$ of pure $\beta$-type that contains $\alpha_0$ as one of its alpha curves. (This can be done because in the proof of Proposition~\ref{prop:gencon}, there was a lot of freedom in the choice of alpha curves: we just wanted them to represent the given handlebody, and we can choose such a set by starting with any curve bounding a disk in the handlebody.) Then, we let the new complete system $\Hyper'$ be obtained from $\Hyper$ by the given trajectory $\alpha$-slide.
 \end {proof}

\newpage \section {Statement of the surgery theorem for link-minimal complete systems}
\label {sec:statement}

We are now ready to state a version of Theorem~\ref{thm:FirstSurgery} for the case of link-minimal complete systems. We let $\orL = L_1 \amalg \dots \amalg L_\ell$ be an oriented link in an integer homology sphere $Y$. We suppose we are given a link-minimal complete system of hyperboxes $\Hyper$ for $\orL$, as in Definition~\ref{def:linkminimalsystem}. In particular, $\Hyper^L$ is a Heegaard diagram for $\orL$. We let $k$ be the number of $w$ basepoints in $\Hyper^L$ (or in any other diagram in $\Hyper$), so that $\Hyper^L$ has $k-\ell$ free basepoints.

Since all the diagrams we work with are link-minimal, we are free to use the description of generalized  Floer complexes $\Am(\cdot, \s)$ as free complexes $\Chain^-(\cdot, \s)$; cf. Section~\ref{sec:alternative}.

\subsection {Descent to a sublink}
\label {subsec:desublink}
Let $M \subseteq L$ be a sublink, with an orientation $\orM$, which can be different from the one coming from $\orL$, as in Section~\ref{sec:reduction}.

 Recall that in Section~\ref{sec:reduction} we defined projection maps
 $$ p^{\orM}:\bH(L) \to \bH(L),$$
as well as reduction maps
$$ \psi^{\orM}: \bH(L) \to \bH(L-M).$$
Fix $\s \in \bH(L)$. Since $\psi^{\orM}$ depends only on the components of $L$ that are not on $M$, we have 
$$\psi^{\orM}(p^{\orM}(\s)) = \psi^{\orM}(\s).$$

Consider the $m$-dimensional hyperbox $\Hyper^{\orL, \orM}$ from the complete system. This hyperbox consists of diagrams representing the link $L-M$. By the definition of the complete system, $\Hyper^{\orL, \orM}$ has to be compatible with $\Hyper^{\orL, \emptyset} = \Hyper^L$, i.e. the initial diagram in $\Hyper^{\orL, \orM}$ is obtained from $\Hyper^L$ by deleting the basepoints $z_i$ on components $L_i$ with $i \in I_+(\orL, \orM)$, deleting $w_i$ on the components $L_i$ with $i \in I_-(\orL, \orM)$, and also relabeling $z_i$ as $w_i$ for the components $L_i$ with $i \in I_-(\orL, \orM)$. In view of \eqref{eq:red_lmn}, we have an identification
\begin{equation}
\label{eq:idem}
 \Chain^-(\he^L, p^{\orM}(\s)) \cong \Chain^-(\he^{\orL, \orM}(\emptyset), \psi^{\orM}(\s)).
 \end{equation}

Let $\dd^{\orM} \in \N^m$ be the size of $\Hyper^{\orL, \orM}$. We have an associated hyperbox of generalized Floer complexes 
$$\Chain^-( \Hyper^{\orL, \orM}, \psi^{\orM}(\s)) = \bigl( (C^{\eps}_{\s})_{\eps \in \E(\dd^{\orM})}, (\De^{\eps}_\s)_{\eps \in \E_m} \bigr)$$ as defined in Section~\ref{sec:hyperfloer}. By compressing the hyperbox $(C^\eps_\s, \De^\eps_\s)$ we obtain a hypercube $(\hat C^\eps_\s, \hat \De^\eps_\s)_{\eps \in \E_m}$, see Sections \ref{sec:compression1}, ~\ref{sec:compression2}. We are only interested in the longest diagonal map in this hypercube, namely $\hat \De^{(1, \dots, 1)}_\s$, which we simply denote by $\hat \De^{\orM}_\s$. According to Equation~\eqref{eq:tildeDe}, we have
$$ 
\hat \De^{\orM}_\s = \play^{\dd^{\orM}}_{\{\De^\eps_\s\}} (\alpha_m),$$
where $\alpha_m$ is the $m\th$ standard symphony and $\play$ denotes the operation of playing songs, see Definitions \ref{def:symph}, \ref{def:playing}. 

Thus, for any $\s \in \bH(L)$, we have defined a {\em descent map}
$$  \hat\De^{\orM}_{\s} : \Chain^-(\he^L, p^{\orM}(\s)) \to \Chain^-(\he^{\orL, \orM}(M), \psi^{\orM}(\s)),$$
as a sum over compositions of polygon maps associated to various sub-hyperboxes of $\Hyper^{\orL, \orM}$. We used here the identification~\eqref{eq:idem}. Since $\hat\De^{\orM}_{\s}$ depends only on the projection $p^{\orM}(\s)$, we will write $\hat\De^{\orM}_{\s}$ as $\hat\De^{\orM}_{p^{\orM}(\s)}.$

\begin {example} 
\label {ex:edge2}
Suppose $M$ is a single link component $L_i$. Each edge of the one-dimensional hyperbox $\Hyper^{\orL, \orL_i}$ comes equipped with a corresponding theta chain element, and this gives a triangle map between the Floer homology groups associated to the initial and final Heegaard diagrams for that edge. The map $\hat \De^{\orL_i}_{p^{\orM}(\s)}$ is the composition of these triangle maps.  Note that all the triangle maps, and hence also $\hat \De^{\orL_i}_{p^{\orM}(\s)}$, are chain homotopy equivalences, see Example~\ref{ex:edge}. 
\end {example}

In a complete system of hyperboxes, we have an isotopy $\he^{\orL, \orM}(M) \cong \he^{L - M}$, which induces an identification between the respective Floer complexes. When we change the range of $\hat \De^{\orM}_{p^{\orM}(\s)}$ via this identification, we denote the resulting map by
\begin {equation}
\label {eq:des}
\De^{\orM}_{p^{\orM}(\s)}: \Chain^-(\he^L, {p^{\orM}(\s)}) \to \Chain^-(\he^{L-M}, \psi^{\orM}(\s)).
 \end {equation}

\begin{figure}
\begin{center}
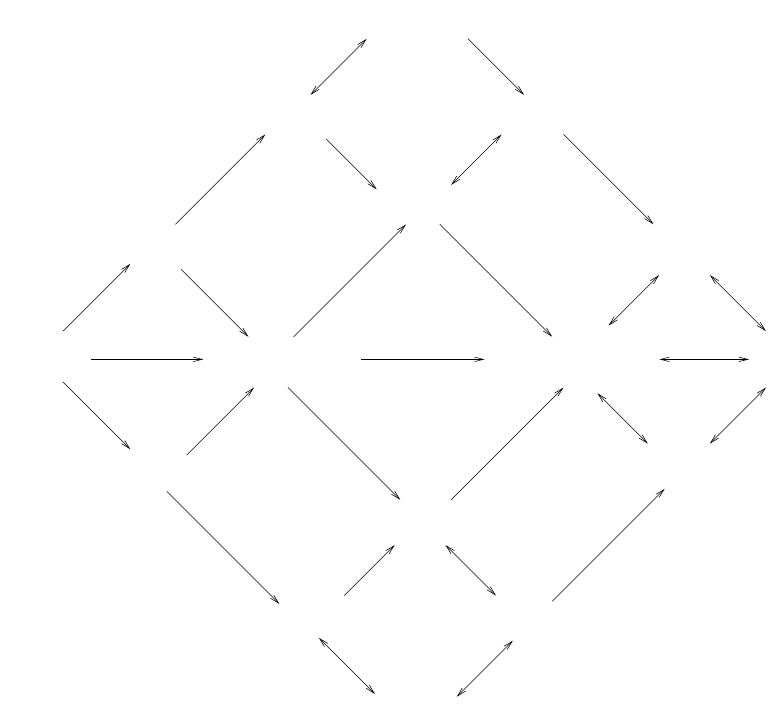
\end{center}
\caption {{\bf Chain maps coming from hyperboxes in a complete system.} We consider here a complete system $\Hyper$ for a link $\orL = L_1 \cup L_2$. Recall that one quarter of such a  complete system was illustrated in Figure~\ref{fig:completeL}. We show here the corresponding chain maps: to each reduction (deletion of basepoints) in the complete system we have an associated inclusion map $\Pr$, to each face of a hyperbox we have an associated descent map $\hat D$, and to each surface isotopy we have an associated isomorphism denoted by the symbol $\cong$. For simplicity, all the sublinks are taken here with their induced orientation from $\orL$, and we drop the arrows from the notation. There are similar maps corresponding to other hyperboxes in the complete system, in which some link components can appear with the opposite orientation. 
}
\label{fig:mapsL1}
\end{figure}

\begin{figure}
\begin{center}
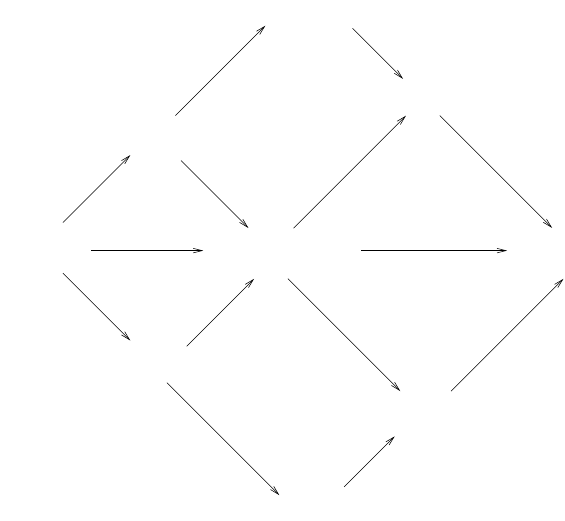
\end{center}
\caption {{\bf More chain maps coming from hyperboxes in a complete system.} This is the same as Figure~\ref{fig:mapsL1}, but we have composed the descent maps $\hat D$ with chain isomorphisms (coming from surface isotopies) to obtain the maps denoted $D$.
}
\label{fig:mapsL2}
\end{figure}

\begin {example}
In the case of a link with two components, some of the inclusion and descent maps are illustrated in Figures~\ref{fig:mapsL1} and \ref{fig:mapsL2}.
\end {example}

For any $\s \in \H(L)$, we now define a map
$$ \Phi^{\orM}_\s : \Chain^-(\he^L, \s) \to \Chain^-(\he^{L-M}, \psi^{\orM}(\s)),$$
\begin {equation}
 \label {eq:Phi}
 \Phi^{\orM}_\s = \De^{\orM}_{p^{\orM}(\s)} \circ \Pr^{\orM}_{\s}.
\end {equation}

Note that we can define similar maps if we replace $\orL$ by a sublink $\orL'$. By abuse of notation, we will always denote the maps corresponding to inclusion and descent at $\orM$ by $ \Pr^{\orM}_{\s}, \De^{\orM}_{p^{\orM}(\s)}, \Phi_\s^{\orM}$, even though their domains of definition may vary.

\begin {lemma}
\label {lemma:comm}
Let $M_1 , M_2 \subseteq L$ be two disjoint sublinks, with orientations $\orM_1$ and $\orM_2$. For any $\s \in \bH(L)$, we have: 
$$ \Pr^{\orM_2}_{\psi^{\orM_1}(\s)} \circ \De^{\orM_1}_{p^{\orM_1}(\s)} = \De^{\orM_1}_{p^{\orM_1 \cup \orM_2}(\s)} \circ \Pr^{\orM_2}_{p^{\orM_1}(\s)}.$$ 
\end {lemma}

\begin{proof}
For $\s = (s_1, \dots, s_\ell)$, let us denote by $s'_i$ the components of $\s' = \psi^{\orM_1}(\s)$, obtained from $s_i$ by subtracting half the linking number of $L_i$ and $M_1$, compare Equation~\eqref{eq:psy}.

Observe that both $ \De^{\orM_1}_{p^{\orM_1}(\s)}$ and $\De^{\orM_1}_{p^{\orM_1 \cup \orM_2}(\s)}$ are sums of compositions of polygon maps. The same polygons get counted in both maps, but with different powers of $U_i$. More precisely, suppose we have a chain of polygons relating intersection points $\x$ from $\he^L$ and $\y$ from $\he^{L-M_1-M_2}$. If $j$ is one of the colors, let $e_j$ resp. $e_j'$ be the exponent of $U_j$ in the coefficient of $\y$ in $\De^{\orM_1}_{p^{\orM_1}(\s)}(\x)$, resp. $\De^{\orM_1}_{p^{\orM_1 \cup \orM_2}(\s)}(\x)$. We have:
\begin {eqnarray*}
 e_j - e_j' &=& \sum_{\{i \in I_+(\orL, \orM_2)| \tau_i = j\}} \bigl( \max(A_i(\x) - s_i, 0) - \max(A_i(\y) - s'_i, 0) \bigr)+\\ & \; & \sum_{\{i \in I_-(\orL, \orM_2)| \tau_i = j\}} \bigl( \max(s_i - A_i(\x), 0) -  \max(s'_i - A_i(\y), 0) \bigr),
 \end {eqnarray*}
compare Section~\ref{sec:polygon}.

On the other hand, according to Equation~\eqref{eq:proj}, the map $\Pr^{\orM_2}_{p^{\orM_1}(\s)}$ contributes a power of $U_j$ with exponent
$$\sum_{\{i \in I_+(\orL, \orM_2)| \tau_i = j\}}  \max(A_i(\x) - s_i, 0) + \sum_{\{i \in I_-(\orL, \orM_2)| \tau_i = j\}} \max(s_i - A_i(\x), 0) $$
and the map $ \Pr^{\orM_2}_{\psi^{\orM_1}(\s)} $ a power of $U_j$ with exponent
$$ \sum_{\{i \in I_+(\orL, \orM_2)| \tau_i = j\}}  \max(A_i(\y) - s'_i, 0) + \sum_{\{i \in I_-(\orL, \orM_2)| \tau_i = j\}}   \max(s'_i - A_i(\y), 0) . $$
These contributions exactly cancel out the difference between $e_j'$ and $e_j$.
\end{proof}

\begin {proposition}
\label {prop:phiphi}
Choose a sublink $M \subseteq L$, and endow it with an orientation $\orM$. Then, for any $\s \in \bH(L)$, we have the following relation
\begin {equation}
\label {eq:summ}
 \sum_{\orM_1 \amalg \orM_2 = \orM} \Phi^{\orM_2}_{\psi^{\orM_1} (\s)} \circ \Phi^{\orM_1}_{\s} = 0,
 \end {equation}
where $\orM_1$ and $\orM_2$ are equipped with the orientations induced from $\orM$.
\end {proposition}

\begin {proof}
We have
\begin {eqnarray*}
 \sum_{\orM_1 \amalg \orM_2 = \orM} \Phi^{\orM_2}_{\psi^{\orM_1} (\s)} \circ \Phi^{\orM_1}_{\s} &=& \sum_{\orM_1 \amalg \orM_2 = \orM} \De^{\orM_2}_{p^{\orM_2}(\psi^{\orM_1}(\s))} \circ \Pr^{\orM_2}_{\psi^{\orM_1} (\s)}\circ \De^{\orM_1}_{p^{\orM_1}(\s)} \circ \Pr^{\orM_1}_{\s} \\
 & = & 
\sum_{\orM_1 \amalg \orM_2 = \orM} \De^{\orM_2}_{p^{\orM_2}(\psi^{\orM_1}(\s))} \circ  \De^{\orM_1}_{p^{\orM_2}(\psi^{\orM_1}(\s))} \circ \Pr^{\orM_2}_{p^{\orM_1} (\s)}\circ  \Pr^{\orM_1}_{\s} \\
&=&
\Bigl( \sum_{\orM_1 \amalg \orM_2 = \orM} \De^{\orM_2}_{p^{\orM_2}(\psi^{\orM_1}(\s))} \circ  \De^{\orM_1}_{p^{\orM_2}(\psi^{\orM_1}(\s))} \Bigr) \circ \Pr^{\orM_1 \cup \orM_2}_{\s} \\
&=& 0.
\end {eqnarray*}

Indeed, the second equality above follows from Lemma~\ref{lemma:comm}, and the last equality is a consequence of Proposition~\ref{prop:compressed}, together with the properties of a complete system of hyperboxes. 
\end {proof}

\subsection {The surgery theorem}
\label {subsec:surgery}

Let us fix a framing $\Lambda$ for the link $\orL$. For a component $L_i$ of $L$, we let $\Lambda_i$ be its induced framing, thought of as an element in $H_1(Y - L)$.

Given a sublink $N \subseteq L$, we  
let $\Omega(N)$ be the set of all possible orientations on $N$, as in the Introduction. For $\orN \in \Omega(N)$, we let 
$$ \Lambda_{\orL, \orN} = \sum_{i \in I_-(\orL, \orN)} \Lambda_i \in H_1(Y - L).$$

We view $\H(L) \subseteq H_1(Y - L; \Q)$ as an affine lattice over $H_1(Y-L)$ as in Remark~\ref{rem:h1}. Thus, if $\s \in \H(L)$, then $\s + \Lambda_{\orL, \orN} $ is also in $\H(L)$.

Now consider the $\Ring$-module
\begin {equation}
\label {eq:chl}
 \C^-(\Hyper, \Lambda) = \bigoplus_{M \subseteq L} \prod_{\s \in \H(L)}  \Chain^-(\Hyper^{L - M}, \psi^{M}(\s) ),
\end {equation} 
where $\psi^{M}$ simply means $\psi^{\orM}$ with $\orM$ being the orientation induced from the one on $\orL$.

Note that the definition of $\C^-(\Hyper, \Lambda)$ involves direct products; in fact, the direct sum in \eqref{eq:chl} can equally be thought of as a direct product, since it is finite. It is worth saying a few words about how one can define maps between direct products:

\begin {definition}
\label {def:locfin}
Let $S$ and $T$ be countable index sets, and $\mathcal{A} = \prod_{s \in S} \mathcal{A}_s, \mathcal{B} = \prod_{t \in T} \mathcal{B}_t$ direct products of modules over a commutative ring $R$. Suppose we are given module homomorphisms 
$$ F_{s,t} : \mathcal{A}_s \to \mathcal{B}_t,$$ 
for each $s \in S$ and $t \in T$. The collection of maps $\{F_{s, t} \}$ is called {\em locally finite} if for each $t \in T$, only finitely many $F_{s, t}$ are nonzero.  
\end {definition}

If $\{F_{s, t}\}$ is a locally finite collection of homomorphisms as in Definition~\ref{def:locfin}, we can  assemble them into a single homomorphism 
$$F: \mathcal{A} \to \mathcal{B}, \ \ F(\{a_s\}_{s \in S}) =  \bigl \{ \sum_{s \in S} F_{s,t}(a_s) \bigr \}_{t \in T}.$$

With this in mind, we equip the module $\C^-(\Hyper, \Lambda) $ with a boundary operator $\D^-$ as follows. For $\s \in \H(L) $ and $\x \in  \Chain^-(\Hyper^{L - M}, \psi^{M}(\s))$, we set
 \begin {eqnarray*} 
\D^-(\s, \x) &=& \sum_{N \subseteq L - M} \sum_{\orN \in \Omega(N)} (\s + \Lambda_{\orL, \orN}, \Phi^{\orN}_{\psi^{M}(\s)}(\x)) \\
&\in&  \bigoplus_{N \subseteq L - M} \bigoplus_{\orN \in \Omega(N)}  \Chain^-(\Hyper^{L-M-N}, \psi^{M \cup \orN} (\s))   \subseteq \C^-(\Hyper, \Lambda).
\end {eqnarray*}

This defines a locally finite collection of maps between the modules $\Chain^-(\Hyper^{L - M}, \psi^{M}(\s))$, producing a well-defined map 
$$ \D^- :  \C^-(\Hyper, \Lambda) \to \C^-(\Hyper, \Lambda).$$
According to Proposition~\ref{prop:phiphi}, $\C^-(\Hyper, \Lambda)$ is a chain complex. Note that  $\C^-(\Hyper, \Lambda)$ naturally breaks into a direct product of terms $\C^-(\Hyper, \Lambda, \ux)$, according to equivalence classes $\ux$ of the values $\s$. Here $\s_1$ and $\s_2$ are equivalent if they differ by an element in the (possibly degenerate) sublattice
$$ H(L, \Lambda) \subseteq H_1(Y -L),$$
generated by all possible $\Lambda_{\orL, \orN}$ or, equivalently, by the component framings $\Lambda_i \in H_1(Y - L)$. The space of equivalence classes is parametrized by the quotient
$$ \H(L) / H(L, \Lambda),$$
which can be naturally identified with the space of $\spc$ structures on the  surgered manifold $Y_\Lambda(L)$, see Remark~\ref{rem:h1} and \cite[Section 3.7]{Links}. 

Given a $\spc$ structure $\ux$ on $Y_\Lambda(L)$, we set 
\begin {equation}
\label {eq:deltux}
 \delt(\ux) = \gcd_{\xi \in H_2(Y_\Lambda(L); \zz)} \langle c_1(\ux), \xi\rangle,
 \end {equation}
where $c_1(\ux)$ is the first Chern class of the $\spc$ structure.  The Heegaard Floer homology $\HFm_*(Y_\Lambda(L), \ux)$ admits a relative $\zz/\delt(\ux)\zz$-grading, see \cite{HolDisk}. In Section~\ref{sec:gradings} we will construct a 
relative $\zz/\delt(\ux)\zz$-grading on the complex $\C^-(\Hyper, \Lambda, \ux)$ as well.

The Surgery Theorem (for link-minimal systems) then says:

\begin {theorem}
\label {thm:surgery}
Fix a link-minimal complete system of hyperboxes $\Hyper$ for an oriented, $\ell$-component link $\orL$ in an integral homology three-sphere $Y$, and fix a framing $\Lambda$ of $L$. Then, for any $\ux \in \spc(Y_\Lambda(L)) \cong  \H(L) / H(L, \Lambda)$, there is an isomorphism of relatively graded $\ff[[U]]$-modules
\begin {equation}
\label {eq:surgery1}
 H_*(\C^-(\Hyper, \Lambda, \ux), \D^-) \cong  \HFm_*(Y_\Lambda(L), \ux).
 \end {equation}
\end {theorem}

Note that the left hand side of~\eqref{eq:surgery1} is a priori an $\Ring$-module, where $\Ring = \ff[[U_1, \dots, U_k]]$. However, part of the claim of the theorem is that all $U_i$'s act the same, so we can think of it as an $\ff[[U]]$-module. 

The proof of Theorem~\ref{thm:surgery} will be given in Section \ref{sec:proof}.

\begin {remark}
In the case when the $\spc$ structure $\ux$ is torsion, one should be able to use the same techniques as in \cite{IntSurg} to obtain an isomorphism of absolutely graded groups, with a well-determined shift in grading between the two sides of \eqref{eq:surgery1}. However, we will not pursue this direction in the present paper. 
\end {remark}

\subsection {Gradings}
\label {sec:gradings}
As promised in the previous subsection, we proceed to construct a relative $\zz/\delt(\ux)\zz$-grading on the complex
$\C^-(\Hyper, \Lambda, \ux)$.  

Let us identify $H_1(Y-L)$ with $\zz^\ell$ as in Remark~\ref{rem:h1}. We view the framing $\Lambda$ as an $\ell$-by-$\ell$ symmetric matrix with columns $\Lambda_i$. The matrix element $c_{ij}$ in $\Lambda$ (for $i, j=1, \dots, \ell$) is the linking number between $L_i$ and $L_j$ when $i\neq j$, and the surgery coefficient $\lambda_i$ of $L_i$ when $i=j$.

Let $H(L, \Lambda)^\perp \subset \zz^\ell$ be the orthogonal complement to $H(L, \Lambda)$, that is,
$$  H(L, \Lambda)^\perp = \{\vs \in \zz^\ell |   \vs \cdot \Lambda_i = 0, \forall i \} = \bigl \{(v_1, \dots, v_\ell) \in \zz^\ell \big |   \sum_i v_i \Lambda_i = 0 \bigr\}.$$

There are natural identifications 
$$H^2(Y_{\Lambda}(L)) \cong H_1(Y_\Lambda(L)) \cong \zz^\ell/ H(L, \Lambda) ,$$ such that
$$ c_1([\s]) = [2\s],$$
for any $\s \in \spc(Y_{\Lambda}(L)) \cong \H(L)/H(L, \Lambda)$. This can be deduced from the formulas for the Chern class in \cite[Equation (24) and Lemma 3.13]{Links}, compare also Lemma~\ref{lem:spc4} below.

 Using Poincar\'e duality, we obtain a natural identification $H_2(Y_\Lambda(L)) \cong H(L, \Lambda)^\perp$. Hence, 
\begin {equation}
\label {eq:dgr}
\delt(\ux) = \gcd_{\vs \in H(L, \Lambda)^\perp} \sum_i 2s_iv_i,
\end {equation}
where we wrote $\vs = (v_1, \dots, v_\ell)$, and $\s=(s_1, \dots, s_\ell)$ is any element in the corresponding equivalence class $\ux \in  \H(L)/H(L, \Lambda)$.

\begin {remark}
\label {rem:even}
It is clear from \eqref{eq:deltux} that $\delt(\ux)$ is always even. One can also verify this using the description of $\delt(\ux)$ given in \eqref{eq:dgr}. Indeed, let $\lambdas = (\lambda_1, \dots, \lambda_\ell) = (c_{11}, \dots, c_{\ell \ell})$ be the diagonal vector of the framing matrix $\Lambda$. For $\s \in \H(L)$ and $\vs \in H(L, \Lambda)^\perp$, we have
$ 2\s \equiv \Lambda_1 + \dots + \Lambda_\ell - \lambdas \pmod 2$
and
$ (\Lambda_1 + \dots + \Lambda_\ell ) \cdot \vs = 0$,
so 
$$ (2\s) \cdot \vs  \equiv  \lambdas \cdot \vs \equiv  \sum_{i=1}^\ell c_{ii} v_i  \equiv \sum_{i=1}^\ell c_{ii} v_i^2  \equiv   \sum_{i=1}^\ell \sum_{j=1}^\ell c_{ij} v_i v_j  \equiv  \vs^T \Lambda \vs 
\equiv 0 \pmod 2.$$
\end {remark}

\begin {lemma}
\label {lemma:nu}
Fix an equivalence class $\ux \in  \H(L)/H(L, \Lambda)$. There exists a function $\nu: \ux \to \zz/\delt(\ux)\zz$ with the property that
\begin {equation}
\label {eq:nus}
 \nu(\s + \Lambda_i) \equiv \nu(\s) + 2s_i,
 \end {equation}
for any $i=1, \dots, \ell$ and $\s = (s_1, \dots, s_\ell) \in \ux$.
\end {lemma}

\begin {proof}
Pick some $\s^0=(s^0_1, \dots, s^0_\ell) \in \ux$. Any other $\s \in \ux$ is of the form $\s^0 + \sum a_i\Lambda_i$, for some $a_i \in \zz$. Set
$$ \nu(\s^0 + \sum_{i=1}^\ell a_i\Lambda_i) =  \sum_{i=1}^\ell 2a_i s^0_i + \sum_{i, j=1}^\ell a_ia_j c_{ij}. $$ 

There is an indeterminacy in expressing $\s$ as $\s^0 + \sum a_i\Lambda_i$, namely one can add an element in $H(L, \Lambda)^\perp$ to the vector $(a_1, \dots, a_\ell)$. It is easy to check that $\nu(\s)$ is independent (modulo $\delt(\ux)$) of how we express $\s$, and that \eqref{eq:nus} is satisfied.
\end {proof}

\begin {remark}
The function $\nu$ from Lemma~\ref{lemma:nu} is unique up to the addition of a constant. 
\end {remark}

Fix a function $\nu$ as in Lemma~\ref{lemma:nu}. Each factor $ \Chain^-(\Hyper^{L - M}, \psi^{M}(\s) )$ appearing in the complex $\C^-(\Hyper, \Lambda, \ux)$ admits a natural $\zz$-grading $\mu^M_\s = \mu_{\psi^M(\s)}$ as in \eqref{eq:ms}. We define a $\zz/\delt(\ux)\zz$-grading $\mu$ on $\C^-(\Hyper, \Lambda, \ux)$ as follows. For $\s \in \ux$ and $\x \in \Chain^-(\Hyper^{L - M}, \psi^{M}(\s) )$, we set
$$ \mu(\s, \x) = \mu^M_\s(\x) +\nu(\s) - |M|,$$
where $|M|$ denotes the number of components of $M$.

\begin {lemma}
The differential $\D^-$ on $\C^-(\Hyper, \Lambda, \ux)$ decreases $\mu$ by one modulo $\delt(\ux)$.
\end {lemma}

\begin {proof}
Use \eqref{eq:mumu}, \eqref{eq:Phi}, and the fact that $\De^{\orM}_{p^{\orM}(\s)}$ (being the longest map in an $|M|$-dimensional hypercube of chain complexes) changes grading by $|M| -1$.
\end {proof}

\section {Truncation}
\label {sec:t}

As noted near the end of Section~\ref{sec:prelim}, in order to calculate the Heegaard Floer homology groups $\HFm(Y_\Lambda(L), \ux)$ using Theorem~\ref{thm:surgery}, it is helpful to replace the infinite direct product from \eqref{eq:chl} with a finite one. This is called {\em horizontal truncation}, and we saw a few instances of it in Sections~\ref{sec:unknot} and \ref{sec:hopf}. In Section~\ref{sec:truncate} below, we will describe a general way of doing horizontal truncation, for surgery on an arbitrary link.

Even after horizontal truncation, the direct product in \eqref{eq:chl} is still an infinite-dimensional $\ff$-vector space, due to the fact that each term is a free module over a ring of power series. However, in Section~\ref{sec:algebra} we show that the power series ring can be replaced (essentially without any loss of information) by a finite-dimensional polynomial ring. This process is called {\em vertical truncation}, and is done by setting large powers of the $U$ variables to zero. By combining horizontal and vertical truncation, we can replace the right hand side of \eqref{eq:chl} with a finite-dimensional chain complex. In Section~\ref{sec:another}, we describe an alternate way of doing so, by applying a slightly different horizontal truncation to the vertically truncated complex. In Section~\ref{sec:twt}, we alter the combined truncation procedure from Section~\ref{sec:another} further, by adding certain ``crossover maps''. We call the result a {\em folded truncation} of the original complex. Folded truncation will play an important role in the proof of Theorem~\ref{thm:surgery} presented in Section~\ref{sec:proof}.

Before describing all these truncations, let us recall some notation from Section~\ref{sec:ex}. We denote a typical term in the surgery complex \eqref{eq:chl} by
%\label {eq:ces}
$$
 \C^\eps_\s = \Chain^-(\Hyper^{L-M}, \psi^M(\s)), 
$$
where $\eps=\eps(M) \in \E_\ell = \{0,1\}^\ell$ is such that $L_i \subseteq M$ if and only $\eps_i = 1$. A typical summand in the differential $\D^-$ is denoted
%\label {eq:dees}
$$
\D^{\eps, \eps'}_{\eps^0, \s} = \Phi^{\orN}_{\psi^M(\s)} : \C^{\eps^0}_\s \to \C^{\eps^0 + \eps}_{\s + \eps' \cdot \Lambda},
$$
where $\eps^0=\eps(M), \eps=\eps(N)$, and $\eps' \in \E_\ell$  is such that $i \in I_-(\orL, \orN)$ if and only if $\eps'_i = 1$. The dot product $\eps' \cdot \Lambda$ denotes the vector $\sum \eps'_i\Lambda_i$. Note that we always have $\eps' \leq \eps$.

Dropping a subscript or superscript from the notation means considering the direct product over all possible values of that subscript or superscript. For example, $\C^\eps = \prod_\s \C^\eps_\s$, and $\C = \oplus_\eps \C^\eps = \C(\Hyper, \Lambda)$. Observe that $(\C^\eps, \D^\eps)$ form a hypercube of chain complexes as defined in Section~\ref{sec:hyperv} (except it may have only a $\zz/2\zz$-grading, rather than a $\zz$-grading) and $\C$ is the total complex of this hypercube.

\subsection {Horizontal truncation}
\label {sec:truncate}
We now return to the setting of Section~\ref{sec:statement}, where $L$ is an arbitrary link, equipped with an orientation, a complete system of hyperboxes $\Hyper$, and a framing $\Lambda$. 

\begin {lemma}
\label {lemma:isos}
There exists a constant $b > 0$ such that, for any $i=1, \dots, \ell$, and for any sublink $M \subseteq L$ not containing the component $L_i$, the chain map
$$\Phi^{\orL_i}_{\psi^{\orM}(\s)}: \Chain^-(\Hyper^{L-M}, \psi^{\orM}(\s)) \to \Chain^-(\Hyper^{L-M-L_i}, \psi^{\orM \cup \orL_i}(\s))$$
induces an isomorphism on homology provided that either
\begin {itemize}
\item  $\s \in \H(L)$ is such that $s_i > b$, and $L_i$ is given the orientation  induced from $L$; or
\item $\s \in \H(L)$ is such that $s_i < -b$, and $L_i$ is given the orientation opposite the one induced from $L$. 
\end {itemize}
\end {lemma}

\begin {proof}
For $|s_i|$ sufficiently large, and $L_i$ oriented as in the lemma (according to the sign of $s_i$), the inclusion map $\Pr^{\orL_i}_{\psi^{\orM}(\s)}$ from \eqref{eq:proj} is the identity. Moreover, the descent map $ \De^{\orL_i}_{p^{\orL_i}(\psi^{\orM}(\s))}$ is a composition of maps along the edges of the corresponding hyperbox $\Hyper^{\orL-M, L_i}$, hence induces an isomorphism on homology, see Example~\ref{ex:edge2}. The conclusion then follows in light of Equation \eqref{eq:Phi}.
\end {proof}

Lemma~\ref{lemma:isos} is the key ingredient in truncating the complex $\C^-(\Hyper, \Lambda)$. Roughly, it allows the terms of this complex to cancel in pairs, whenever $\s \in \H(L)$ has at least one component $s_i$ with $|s_i| > b$. The result is that the homology of the complex $\C^-(\Hyper, \Lambda)$ can be computed by restricting to some terms corresponding to $\s$ lying in a compact subset  of $\H(L)$.

Of course, we need to be more explicit about how this is done. For simplicity, let us assume (for the moment) that the framing vectors $\Lambda_1, \dots, \Lambda_\ell$ are linearly independent in $\rr^\ell$. This is equivalent to asking for $Y_\Lambda(L)$ to be a rational homology sphere.

Let us first recall the case of surgery on knots, see \cite[Section 4.1]{IntSurg}. Then the framing coefficient $\Lambda$ is a nonzero integer.  Set
\begin {equation}
\label {eq:knottruncate}
 \C^-(\Hyper, \Lambda)\langle b \rangle = \bigoplus_{-b \leq s \leq b} \Chain^-(\Hyper^L, s) \oplus \bigoplus_{-b+\Lambda \leq s \leq b} \Chain^-(\Hyper^\emptyset, \psi^{L}(s)) \subset \C^-(\Hyper, \Lambda).
 \end {equation}

It is easy to check that $\C^-(\Hyper, \Lambda)\langle b \rangle $ is a subcomplex of $\C^-(\Hyper, \Lambda)$ for $\Lambda < 0$, and a quotient complex of $\C^-(\Hyper, \Lambda)$ when $\Lambda > 0$. In both cases, an application of  Lemma~\ref{lemma:isos} shows that $\C^-(\Hyper, \Lambda)$ and $ \C^-(\Hyper, \Lambda)\langle b \rangle$ are quasi-isomorphic.

Next, we turn to the case when the link $L$ has two components. We denote by $\lambda_1, \lambda_2$ the framings of $L_1, L_2$ (as compared to the framings coming from Seifert surfaces for those components), and by $c$ the linking number between $L_1$ and $L_2$. Thus, in terms of the standard basis of $H_1(Y-L) \cong \zz^2$, we have 
$$\Lambda_1 = (\lambda_1, c), \ \ \ \Lambda_2 = (c, \lambda_2).$$

Recall our assumption that $\Lambda_1 $ and $\Lambda_2$ are linearly independent, i.e. $\lambda_1\lambda_2 - c^2 \neq 0$.

As before, $\C^{\eps_1\eps_2}_\s$ will denote the term $ \Chain^-(\Hyper^{L - M}, \psi^{M}(\s) )$ appearing in \eqref{eq:chl}, where $\eps_i=1$ or $0$ depending on whether or not $L_i \subseteq M, i=1,2$. We say that $\C^{\eps_1\eps_2}_\s$ is supported at the point $\s = (s_1, s_2) \in \H(L) \subset \rr^2$.

\begin{figure}
\begin{center}
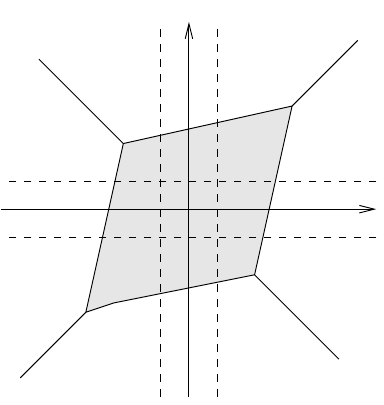
\end{center}
\caption {{\bf The parallelogram $Q$ and the regions $R_i$.} 
The parallelogram $Q$, which is (roughly) the support of the truncated complex  $\C^-(\Hyper, \Lambda)\langle b \rangle$, is shown shaded. 
}
\label{fig:quadrilateral}
\end{figure}

Let $b$ be the constant from Lemma~\ref{lemma:isos}. We seek to define a chain complex  $\C^-(\Hyper, \Lambda)\langle b \rangle $ quasi-isomorphic to  $\C^-(\Hyper, \Lambda)$, and composed of only finitely many of the terms  $\C^{\eps_1\eps_2}_\s$. For this purpose, we construct a parallelogram $Q$ in the plane, with edges parallel to the vectors $\Lambda_1, \Lambda_2$, and with vertices $P_1, P_2, P_3, P_4$ as in Figure~\ref{fig:quadrilateral}. We require the coordinates of $P_1$ to satisfy $x,y>b$, the coordinates of $P_2$ to satisfy $x<-b, y>b$, etc. Further, we choose half-lines $l_1, l_2, l_3, l_4$ with $l_i$ starting at $P_i$ and 
staying in the $i$th quadrant. For example, $l_1$ has to form an angle between $0$ and $\pi/2$ with the $x$ axis, etc. The half-lines $l_i$ split the complement $\rr^2 - Q$ into four regions, denoted $R_1, R_2, R_3, R_4$, see Figure~\ref{fig:quadrilateral}. We require that the edges of the parallelogram $Q$ and the half-lines $l_i$ miss the lattice $\H(L)$.

\begin{figure}
\begin{center}
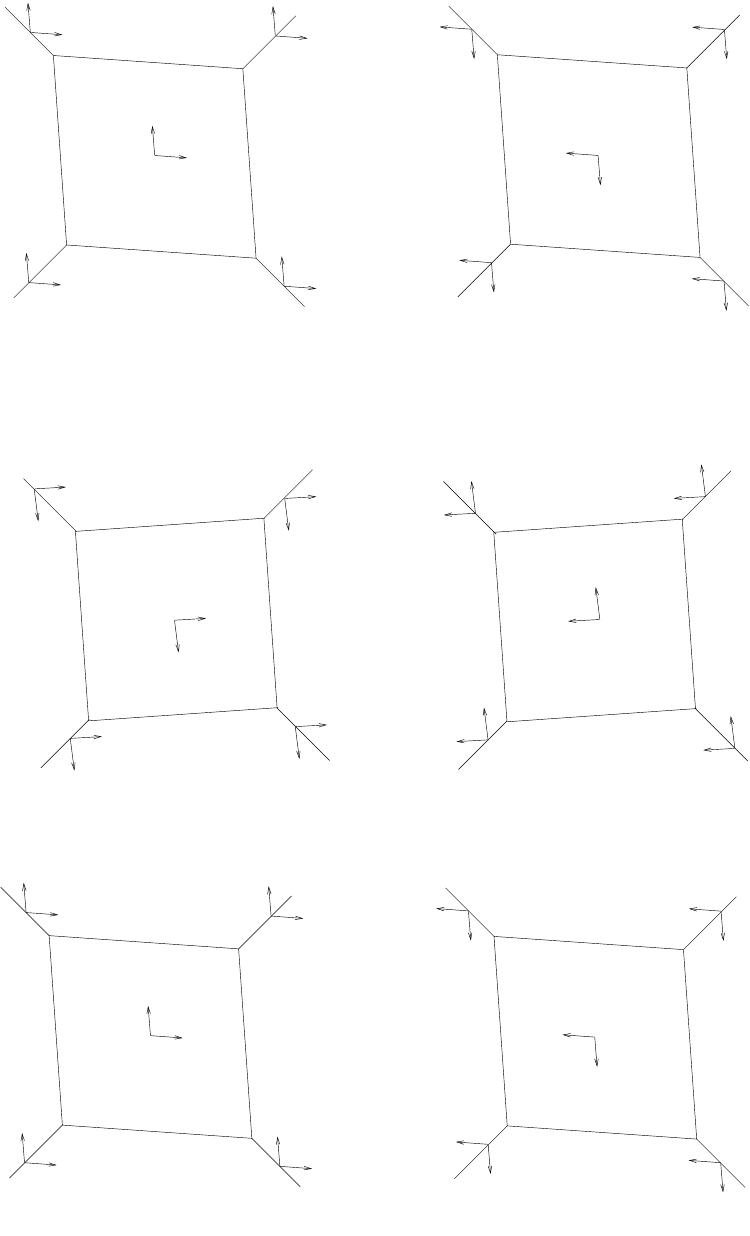
\end{center}
\caption {{\bf Horizontal truncation.} We show here the required slopes for the half-lines $l_i$, in relation to the vectors $\Lambda_1, \Lambda_2$. There are six cases.
}
\label{fig:six}
\end{figure}

We require that the slopes of the lines $l_1, l_2, l_3$ and $l_4$, compared to those of the vectors $\Lambda_1, \Lambda_2$,  are as indicated in Figure~\ref{fig:six}. We distinguish there six cases, according to the values $\lambda_1, \lambda_2$ and $c$. (Note that these cases cover all the possibilities for $\lambda_1 \lambda_2 - c^2 \neq 0$, with some overlap. When we are in an overlap situation, we are free to choose either setting.) For example, in Case I, we require both vectors $\Lambda_1$ and $\Lambda_2$ to point into $R_2$ when placed on the half-line $l_2$. 

By Lemma~\ref{lemma:isos}, the maps $\Phi^{L_1}$ appearing in the definition of the differential for $\C^-(\Hyper, \Lambda)$ induce isomorphisms on homology when restricted to terms supported in $R_1$. The same is true for the maps $\Phi^{L_2}$ supported in $R_2$, for the maps $\Phi^{-L_1}$ supported in $R_3$, and for the maps $\Phi^{-L_2}$ supported in $R_4$. (Here, as in Section~\ref{sec:hopf}, we let 
$L_1, L_2$ have the orientation induced from $L$, and we add a minus sign when we wish to indicate the opposite orientation.) These facts allow us to truncate the complex $\C^-(\Hyper, \Lambda)$ to a obtain a quasi-isomorphic one $\C_Q = \C^-(\Hyper, \Lambda)\langle b \rangle $, which is (roughly) supported in the parallelogram $Q$. This is obtained by taking successive subcomplexes and quotient complexes from  $\C^-(\Hyper, \Lambda)$, and cancelling out complexes supported in $R_1, R_2, R_3, R_4$. Some care has to be taken with what happens near the boundaries of these regions, so we proceed to do a case by case analysis. We present the first case in detail, and for the others we sketch the necessary modifications.

\medskip
\textbf{ Case I: } $\lambda_1, \lambda_2 > 0, \lambda_1 \lambda_2 - c^2 > 0$. Let $\C_{R_1 \cup R_2}$ be the subcomplex of $\C^-(\Hyper, \Lambda)$ consisting of those terms supported in $R_1 \cup R_2$. There is a filtration $\F_{00}$ on $\C_{R_1 \cup R_2}$ (analogous to the one used for the Hopf link in Section~\ref{sec:hopf}) such that $\F_{00}$ is bounded above, and in the associated graded we only see the differentials that preserve $\s$. The associated graded then splits into a direct sum $\C_{R_1} \oplus \C_{R_2}$, where $\C_{R_i}$ is the direct product of the terms supported in $R_i$. On $\C_{R_1}$ there is an additional filtration such that the differential of its associated graded consists only of maps of the form $\Phi^{L_1}$, which are isomorphisms on homology. Hence $H_*(\C_{R_1}) = 0$. Similarly $H_*(\C_{R_2}) = 0$, using an associated graded that leaves only the maps $\Phi^{L_2}$ in the differential. Putting these together, we deduce that $\C_{R_1 \cup R_2}$ is acyclic.

The quotient complex of $\C_{R_1 \cup R_2}$ is supported in $Q \cup R_3 \cup R_4$. Let us define a subcomplex of it, denoted $\C_{R_3 \cup R_4}$, to consist of those terms $\Chain^{\eps_1\eps_2}_{\s}$ with the property that $\s - \eps_1\Lambda_1 -\eps_2\Lambda_2 \in R_3 \cup R_4$. This is roughly supported in $R_3 \cup R_4$, although some terms spill over into $Q$. We define $\C_Q$ to be its quotient complex.

We claim that $\C_Q$ is quasi-isomorphic to the original complex $\C^-(\Hyper, \Lambda)$. For this, we need to show that $\C_{R_3 \cup R_4}$ is acyclic. In the region $R_3$, we would like to use the maps $\Phi^{-L_1}$ to cancel out terms in pairs. However, there exist a whole region of values $\s \in R_3$ such that $\s+\Lambda_2$ lands in $R_2$ rather than $R_3$. The direct product of $\C^{00}_\s$ and $\C^{01}_\s$ over the values $\s$ in that region forms a quotient complex of $\C_{R_3 \cup R_4}$; this quotient is acyclic, because the maps $\Phi^{L_2}$ (being close to the region $R_2$) make the terms $\C^{00}_\s$ and $\C^{01}_\s$  cancel out in pairs. The corresponding subcomplex $\C'_{R_3 \cup R_4}$ is quasi-isomorphic to $\C_{R_3 \cup R_4}$. Similarly, we can eliminate the terms $\C^{00}_\s$ and $\C^{10}_\s$ from $\C'_{R_3 \cup R_4}$ for those $\s \in R_4$ such that $\s + \Lambda_1 \in R_1$. The result is a quasi-isomorphic complex $\C''_{R_3 \cup R_4}$. Consider the associated graded of this complex with respect to a filtration $\F_{11}$ (analogous to the one used for the Hopf link in Section~\ref{sec:hopf}), such that the remaining differentials preserve $\s - \eps_1\Lambda_1 - \eps_2\Lambda_2$. The associated graded is acyclic, as it splits into a direct sum according to the regions $R_3$ and $R_4$, and those are acyclic by Lemma~\ref{lemma:isos}. We conclude that $H_*(\C_{R_3 \cup R_4}) = H_*(\C''_{R_3 \cup R_4}) = 0$. 

\medskip
\textbf{ Case II: } $\lambda_1, \lambda_2 < 0, \lambda_1 \lambda_2 - c^2 > 0$. This is similar to Case I, except $\C_{R_1 \cup R_2}$ and $\C_{R_3 \cup R_4}$ are quotient complexes, and $\C_Q$ is a subcomplex. Further, the filtrations we use are no longer bounded above; still, they are $U$-tame in the sense of Definition~\ref{def:utame}, and we can apply Lemma~\ref{lem:acyclicAG2}.

\medskip
\textbf{ Case III: } $\lambda_1 > 0, \lambda_2 < 0$. We define a subcomplex $\C_{R_1}$ composed of those terms $\C^{\eps_1\eps_2}_\s$ such that either $\s \in R_1$, or ($\s \in R_4, \eps_2=1$ and $\s - \Lambda_2 \in R_1$). We also define another subcomplex $\C_{R_3}$ composed of those terms $\C^{\eps_1\eps_2}_\s$ such that either $\s -\eps_1\Lambda_1 \in R_3$, or ($\s \in R_4, \eps_2=1$ and $\s - \eps_1 \Lambda_1 - \Lambda_2 \in R_3$). Both $\C_{R_1}$ and $\C_{R_3}$ are acyclic. The corresponding quotient complex admits two further acyclic quotient complexes: one, $\C_{R_2}$, consisting of $\C^{\eps_1\eps_2}_\s$ with $\s \in R_2$ and the other, $\C_{R_4}$, consisting of $\C^{\eps_1\eps_2}_\s$ such that $\s - \eps_2\Lambda_2 \in R_4$. We let $\C_Q$ be the resulting subcomplex. 

\medskip
\textbf{ Case IV: } $\lambda_1 < 0, \lambda_2 > 0$. Similar to Case III, but now $\C_{R_1}$ and $\C_{R_3}$ are quotient complexes, while $\C_{R_2}$ and $\C_{R_4}$ are subcomplexes.

\medskip
 \textbf{ Case V: } $c > 0, \lambda_1 \lambda_2 - c^2 < 0$. This is similar to Case I, except that in order to make the complex $\C_{R_3 \cup R_4}$ acyclic, we need to slightly modify its definition near the corners $P_2, P_4$ of the parallelogram $Q$. Specifically, when we define $\C_{R_3 \cup R_4}$, we still ask that $\s - \eps_1\Lambda_1 -\eps_2\Lambda_2 \in R_3 \cup R_4$, but from such pairs $(\s, \eps)$ we eliminate those with: (a) $\s \in Q, \s+\Lambda_1 \in R_2, \s - \Lambda_2 \in R_3$ and $\eps = (0,1)$; as well as those with: (b) $\s \in Q, \s + \Lambda_2 \in R_1, \s - \Lambda_1 \in R_4$ and $\eps=(1,0)$. 

\medskip
\textbf{ Case VI: } $c < 0, \lambda_1 \lambda_2 - c^2 < 0$. This is similar to Case II, except we need to make adjustments near the corners $P_2$ and $P_4$ as in Case V.

\medskip

A similar construction can be done for surgery on links with an arbitrary number of components $\ell$. Recall that we assume the framing vectors $\Lambda_1, \dots, \Lambda_\ell$ to be linearly independent in $\rr^\ell$ (i.e. $Y_\Lambda(L)$ is a rational homology sphere). Let $x_1, \dots, x_\ell$ be the coordinates in $\rr^\ell \cong H_1(Y - L; \rr)$. The coordinate hyperplanes $x_i = 0$ split $\rr^\ell$ into $2^\ell$ ``hyper-quadrants''. We construct a ``skewed hyperbox'' $Q$ (analogous to the parallelogram $Q$ in the case $\ell=2$), with one vertex in each hyper-quadrant, as follows. The vertices of $Q$ are $P_\eps, \eps \in \E_\ell = \{0,1\}^\ell$, with coordinates satisfying
$$ x_i(P_\eps) > b, \text{ if } \eps_i = 0,$$
$$ x_i(P_\eps) < -b, \text{ if } \eps_i = 1. $$

Thus, the skewed hyperbox $Q$ contains the hypercube $[-b, b]^\ell$, where $b$ is the constant in Lemma~\ref{lemma:isos}. Further, we require the edges of $Q$ to be parallel to the vectors $\Lambda_i$.

Let 
$$F_{i, \sigma}, \ i\in \{1, \dots, \ell\}, \sigma \in \{-1,1\} $$ 
be the hyperface of $Q$ that lies completely in the half-space given by $\sigma x_i > 0$. In other words, $F_{i, \sigma}$ has as vertices all  $P_\eps$ with $(-1)^{\eps_i} = \sigma$. 

 We can truncate the complex $\C^-(\Hyper, \Lambda)$ to obtain a quasi-isomorphic one $\C_Q$, roughly supported in the skewed hyperbox $Q$. This truncated complex $\C_Q=\C^-(\Hyper, \Lambda)\langle b \rangle$ is obtained from $\C^-(\Hyper, \Lambda)$ after canceling some acyclic subcomplexes and quotient complexes, one for each face $F_{i, \sigma}$. We use Lemma~\ref{lemma:isos} to show acyclicity, along the same lines as in the case $\ell = 2$. Note that, near the faces $F_{i, +1}$, the truncation is done exactly along the boundaries of $Q$, while near the faces $F_{i, -1}$, we allow some terms to spill in or out of $Q$, i.e. instead of requiring (locally) that $\s \in Q$, we have a requirement of the form $\s -  \eps_i \Lambda_i \in Q$. This allows for the cancellation of the terms outside $\C_Q$. We leave the verification of the details to the interested reader. 
  
Finally, let us turn to the case when the framing matrix $\Lambda$ is degenerate, so that $b_1(Y_\Lambda(L)) > 0$. Then one can still truncate each complex $\C^-(\Hyper, \Lambda, \ux)$, corresponding to a specific $\spc$ structure $\ux$ on $Y_\Lambda(L)$. The truncations is done in the same way as in the case $b_1(Y_{\Lambda}(L)) =0$, but we cut only in $\ell - b_1(Y_{\Lambda}(L))$ directions. 

 In principle, the complex $\C^-(\Hyper, \Lambda)$ is a direct product over $\C^-(\Hyper, \Lambda, \ux)$, and there are infinitely many $\spc$ structures $\ux$. Nevertheless, Theorem 2.3 in \cite{HolDiskTwo}, together with Equation~\ref{eq:infty}, implies that there are only finitely many $\spc$ structures $\ux$ for which $\HFm(Y_\Lambda(L), \ux) \neq 0$. In addition, one can find an a priori bound (in terms of a suitable Heegaard diagram) to determine the range of possible $\ux$ with nonzero Floer homology. By only taking those particular $\ux$ in the direct product $\prod \C^-(\Hyper, \Lambda, \ux)$ we obtain a chain complex quasi-isomorphic to the original $\C^-(\Hyper, \Lambda)$, which we can then truncate to arrive at a quasi-isomorphic finite direct product.
We denote the resulting chain complex with finite support by $\C^-(\Hyper, \Lambda)\langle b \rangle$, depending on $b \gg 0$.

\subsection {Vertical truncation} 
\label {sec:algebra}
This section is an analogue of Section 2.7 in \cite{IntSurg}, with $\HFm$ replacing $\iHF^+$, and with the use of possibly several $U$ variables, as well as non-torsion $\spc$ structures. 

Let $C$ be a chain complex over $\Ring = \ff[[U_1, \dots, U_\ell]]$, with a relative $\zz/2N\zz$-grading, where $N$ is a nonnegative integer and each $U_i$ has degree $-2$. Let $\delta$ be a positive integer. 
Let $\Ring^\delta$ be the quotient of $\Ring$ by the ideal generated by $U_i^{\delta}, i=1, \dots, \ell$. We then denote by $C^{\delta}$ be the complex $C \otimes_\Ring \Ring^\delta$. Further, if $F: C \to D$ is a map between chain complexes over $\Ring$, we denote by $F^\delta$ the corresponding map between $C^\delta$ and $D^\delta$. The procedure of replacing $C$ by $C^\delta$, or $F$ by $F^\delta$ is referred to as {\em vertical truncation}.

If $\delta' \geq \delta$, note that there is a natural projection $\Ring^{\delta'} \to \Ring^{\delta}$, which gives a map $C^{\delta'} \to C^{\delta}$. Set
\begin {equation}
\label {eq:dd}
H_*^{\delta \from \delta'}(C) = \im \bigl(H_*(C^{\delta'}) \longrightarrow H_*(C^{\delta}) \bigr) .
\end {equation}

\begin {definition}
A chain complex $C$ over $\Ring =\ff[[U_1, \dots, U_\ell]]$ is said to be of {\em torsion $\CFm$ type} if it admits an absolute $\qq$-grading, a relative $\zz$-grading, and it is quasi-isomorphic (over $\Ring$) to a finitely generated, free chain complex over $\Ring$.
\end {definition}

\begin {remark} 
The prototype of a chain complex of torsion $\CFm$ type is the Heegaard Floer complex $\CFm(Y, \ux)$, where $Y$ is a three-manifold and $\ux$ a torsion $\spc$ over $Y$; see \cite{HolDisk}, \cite{HolDiskTwo}, \cite{HolDiskFour}. Also, when $\ux$ is torsion, the complex $\C^-(\Hyper, \Lambda, \ux)$ constructed in Section~\ref{subsec:surgery} is of $\CFm$ type, because it is quasi-isomorphic to a (finitely generated, free) horizontally truncated complex as in Section~\ref{sec:truncate}.
\end {remark} 

Note that if $C$ is of torsion $\CFm$ type, then the homology $H_k(C)$ vanishes for $k \gg 0$. Further, as an $\Ring$-module, $H_*(C)$ admits a quotient module $H_{\geq k} = \oplus_{i \geq k} H_i(C),$ for any $i \in \zz$.

In the lemma below, the symbol $\cong$ denotes isomorphism of $\Ring$-modules.

\begin {lemma}
\label {lemma:tors}
(a) Let $C$ be a complex of torsion $\CFm$ type. Then, for any $k \in \qq$, there exists a constant $d$ such that for all integers $\delta \geq d$, we have $H_{\geq k}(C^\delta) \cong  H_{\geq k}(C)$.

(b) If $A, B$ are chain complexes of torsion $\CFm$ type satisfying $H_*(A^\delta) \cong H_*(B^\delta)$ for all $\delta \gg 0$, then $H_*(A) \cong H_*(B)$.

(c) More generally, suppose we are given integers $a, b \geq 0$. If $A, B$ are chain complexes of torsion $\CFm$ type satisfying 
$$H_*(A^\delta) \otimes \bigl( H_{*+2\delta -1}(S^{2\delta-1}) \bigr)^{\otimes a} \cong H_*(B^\delta)  \otimes \bigl( H_{*+2\delta -1}(S^{2\delta-1}) \bigr)^{\otimes b}$$ 
for all $\delta \gg 0$, then $H_*(A) \cong H_*(B)$. (The tensor products are taken over $\Field$.)
\end {lemma}

\begin {proof}
(a) This is similar to the proof of Lemma 2.7 in \cite{IntSurg}. We consider the short exact sequence
\begin {equation}
\label {eq:ses}
 0 \longrightarrow C \xrightarrow{U_1^{\delta}} C \longrightarrow C/U_1^{\delta}C \longrightarrow 0.
 \end {equation}

For a given $k$, we choose $d$ such that $H_i(C) = 0$ for $i \geq k+2d-1$. The induced long exact sequence then gives an isomorphism $H_{\geq k}(C) \cong H_{\geq k}(C/U_1^{\delta}C)$ for $\delta \geq d$. Iterate this argument, replacing $C$ with $C/U_1^{\delta}C$ and $U_1$ with $U_2$, then use $U_3$, etc.

Parts (b) and (c) follow from (a), combined with the fact that $H_{\geq k}(A) \cong H_{\geq k}(B)$ for all $k$ implies $H_*(A) \cong H_*(B)$. For (c), we also need to note that  
$ H_{*}(A^\delta)$ and $H_*(A^\delta) \otimes \bigl( H_{*+2\delta -1}(S^{2\delta-1}) \bigr)^{\otimes a}$ have the same quotient modules corresponding to degrees $\geq k$, for any $\delta$ sufficiently large compared to $k$; and that a similar statement applies to $B^\delta$. 
\end {proof}

\begin {remark}
\label {rem:tors}
If $C$ is finitely generated, we can estimate the  value of $d$ in Lemma~\ref{lemma:tors} (a) as follows: if $m$ is the maximal degree of the generators of $C_*$, we can choose $d > (m-k)/2$.
\end {remark}

\begin {definition}
A chain complex $C$ over $\Ring =\ff[[U_1, \dots, U_\ell]]$ is said to be of {\em non-torsion $\CFm$ type} if it admits a relative $\zz/2N\zz$-grading for some $N> 0$, and it is quasi-isomorphic (over $\Ring$) to a finitely generated, free chain complex over $\Ring$. 
\end {definition}

\begin {lemma}
\label {lemma:nt}
Let $Y$ be a three-manifold and $\ux$ a non-torsion $\spc$ structure over $Y$. Then the Heegaard Floer complex $\CFm(Y, \ux)$ is of non-torsion $\CFm$ type, and there exists $d \geq 0$ such that $U_i^d  \HFm(Y, \ux) = 0$ for all $i$. 
\end {lemma}

\begin {proof}
The fact that $\CFm(Y, \ux)$ is of non-torsion $\CFm$ type was established in \cite{HolDisk}.

For the second statement, note that all $U_i$'s act the same way on homology, see \cite{Links}. Let us denote their common action by $U$. We need to check that this action is nilpotent. Indeed, we have $\HFinf (Y, \ux) = 0$ by Equation~\eqref{eq:infty0}. Since $\HFinf$ is the ring of fractions of $\HFm$ with respect to $(U)$, the action of $U$ on $\HFm$ must be nilpotent.
\end {proof}

\begin {lemma}
\label {lemma:nontors}
(a) Let $C$ be a  complex of non-torsion $\CFm$ type. Suppose there exists $d \geq 0$ such that $U_i^d  H_*(C) = 0$ for all $i$. Then, for all integers $\delta \geq d$ and $\delta' \geq \delta + d$, we have an  isomorphism of relatively graded $\Ring$-modules:
\begin {equation}
\label {eq:s3}
H_*^{\delta \from \delta'}(C) \cong  H_*(C) .
\end {equation}

(b) Let $A, B$ be relatively $\zz/2N\zz$-graded chain complexes of non-torsion $\CFm$ type satisfying $H_*^{\delta \from \delta'}(A) \cong H_*^{\delta \from \delta'}(B)$ for all $\delta' \geq \delta \gg 0$. Suppose there exists $d \geq 0$ such that $U_i^d  H_*(A) = 0$ for all $i$. Then $H_*(A) \cong H_*(B)$, as relatively graded $\Ring$-modules.
\end {lemma}

\begin {proof}
(a) Note that the short exact sequences \eqref{eq:ses} for $\delta'$ and $\delta$ fit into a commutative diagram
$$
\begin {CD}
0 @>>> C @>{U_1^{\delta'}}>> C @>>> C/U_1^{\delta'}C @>>> 0 \\
@.  @V{U_1^{\delta'-\delta}}VV @| @VVV  @. \\
0 @>>> C @>{U_1^{\delta}}>> C @>>> C/U_1^{\delta}C @>>> 0
\end {CD}
$$

At the level of homology, this produces the commutative diagram
$$
\begin {CD}
0 @>>> H_*(C)  @>>> H_*(C/U_1^{\delta'}C) @>>> H_*(C) @>>> 0 \\
@.   @| @VVV  @VV{U_1^{\delta'-\delta}}V @. \\
0  @>>> H_*(C)  @>>> H_*(C/U_1^{\delta}C) @>>> H_*(C) @>>> 0
\end {CD}
$$

The third vertical arrow is zero by hypothesis. It follows that $\im \bigl (H_*(C/U_1^{\delta'}C) \longrightarrow  H_*(C/U_1^{\delta}C)\bigr) \cong H_*(C)$. By iterating this argument $\ell$ times, we obtain \eqref{eq:s3}.

(b) If $U_i^d$ annihilates $H_*(A)$, the long exact sequence in homology associated to \eqref{eq:ses} implies that $U_i^d$ also annihilates $H_*(A^\delta)$, for all $\delta \gg 0$. Since $H_*(A^\delta) \cong H_*(B^\delta)$ by hypothesis (where we chose $\delta' = \delta$), we get that $U_i^d H_*(B^\delta) = 0$ for $\delta \gg 0$. 

We claim that this implies  $U_i^d H_*(B) = 0$ for all $i$. Let us explain the argument in the case $\ell =1$. The long exact sequence in homology 
$$ \dots \to H_*(B) \xrightarrow{U_1^\delta} H_*(B) \longrightarrow H_*(B^\delta) \to \dots $$
implies that $M = U_1^d H_*(B)$ is in $U_1^\delta H_*(B) = U_1^{\delta - d} M$ for all large $\delta$. Since $M$ is a finitely generated module over the local Noetherian ring $\Ring = \ff[[U_1, \dots, U_\ell]]$, Krull's Theorem implies that $\cap_i  \ U_1^iM = 0$, see \cite[Corollary 10.19]{AtiyahM}. Thus $U_1^d H_*(B) = 0$. Iterating this argument $\ell$ times produces the same conclusion for an arbitrary number of $U_i$ variables.

The claim that $H_*(A) \cong H_*(B)$ now follows by applying (a) to both $A$ and $B$.
\end {proof}

\begin {lemma}
\label {lemma:ntors}
Let $C$ be a complex of non-torsion $\CFm$ type, such that all $U_i$ act the same on homology, and $U^d  H_*(C) = 0$ for some $d \geq 0$, where $U$ denotes the common $U_i$ action. Suppose we have $\dim_{\Field} (H_*(C^{\delta+1})) =  \dim_{\Field} (H_*(C^{\delta}))$ for some $\delta \geq 1$. Then $H_*^{\delta \from 2\delta }(C) \cong H_*(C)$.

\end {lemma}

\begin {proof}
Explicitly, the homology $H_*(C)$ is a direct sum $ \oplus_j (\ff[[U]]/U^{k_j})$, so that $\dim_{\Field} (H_*(C)) = \sum k_j$. The long exact sequence on homology associated to \eqref{eq:ses}, iterated $\ell$ times, implies that $\dim_{\Field}   (H_*(C^\delta)) = 2^\ell \sum \min (k_j, \delta)$. If $\dim_{\Field} (H_*(C^{\delta+1})) =  \dim_{\Field} (H_*(C^{\delta}))$, we must have $\delta \geq k_j$ for all $j$, which means that $U_i^\delta H_*(C) = 0$. The claim now follows from Lemma~\ref{lemma:nontors} (a).
\end {proof}

Let us apply this discussion to the complex $\C^-(\Hyper, \Lambda)$ from \eqref{eq:chl}. Suppose we understand the chain groups and differentials and we want to compute its homology. First, we decompose the complex into terms of the form $\C^-(\Hyper, \Lambda, \ux)$, according to $\spc$ structures $\ux$. Then, we apply the horizontal truncation from Section~\ref{sec:truncate} and get quasi-isomorphic complexes $\C^-(\Hyper, \Lambda, \ux)\langle b \rangle$. If $\ux$ is torsion, the complex $\C^-(\Hyper, \Lambda, \ux)\langle b \rangle$ is finite dimensional in each given degree, so we can compute its homology. Alternatively, we could replace it by a quasi-isomorphic, finite dimensional complex $\C^-(\Hyper, \Lambda, \ux)\langle b \rangle^\delta$, see Lemma~\ref{lemma:tors}, where $\delta$ can be estimated as in Remark~\ref{rem:tors}. If $\ux$ is non-torsion, the complex $\C^-(\Hyper, \Lambda, \ux)\langle b \rangle$ is of non-torsion $\CFm$ type, and its homology is annihilated by a power of the $U_i$'s, see Theorem~\ref{thm:surgery} and Lemma~\ref{lemma:nt}. We then start computing the homology of the complexes $\C^-(\Hyper, \Lambda, \ux)\langle b \rangle^\delta$, and let $\delta$ increase by one until we find that the complexes for $\delta$ and $\delta+1$ have the same total rank. By Lemma~\ref{lemma:ntors}, we have $H_*(\C^-(\Hyper, \Lambda, \ux)) \cong H_*^{\delta \from 2\delta}(\C^-(\Hyper, \Lambda, \ux))$. The latter homology group can be computed from finite dimensional complexes.

\subsection {A combined truncation} \label{sec:another} Pick $\delta > 0$ and consider the vertically truncated complex $\C^\delta = \C^-(\Hyper, \Lambda)^\delta$. We could apply the horizontal truncation procedure from Section~\ref{sec:truncate} to obtain a quasi-isomorphic, finite dimensional complex $\C^-(\Hyper, \Lambda)^\delta \langle b \rangle$. In this section we describe a different way of doing horizontal truncation. This new procedure cannot be applied directly to $\C^-(\Hyper, \Lambda)$; it is essential to do the vertical truncation by $\delta$ first. The new horizontal truncation is an intermediate step towards constructing the folded truncation in Section~\ref{sec:twt}.

We use the notation introduced at the beginning of Section~\ref{sec:t}, i.e. we denote by $\C_\s^{\eps, \delta}$ the factors of $\C^\delta$, and by $\D^{\eps, \eps', \delta}_{\eps^0, \s}$ the maps that form the differential. The property that distinguishes the vertically truncated complex $\C^\delta$ from $\C$ is the following:

\begin {lemma}
\label {lemma:bdelta}
Fix $\delta > 0$. Then, there is a constant $b^\delta > 0$ such that, for any $i=1, \dots, \ell$, the map $\D^{\eps, \eps', \delta}_{\eps^0, \s}$ is a quasi-isomorphism provided that either
\begin {itemize}
\item  $\eps = \tau_i$ (i.e. $\eps_i = 1$ and $\eps_j = 0$ for $j\neq i$), $\eps' = 0$, and $s_i > b^\delta$, or
\item $\eps = \eps' = \tau_i$ and $s_i < -b^\delta$; 
\end {itemize}
and, further, we have $\D^{\eps, \eps', \delta}_{\eps^0, \s}=0$ provided that either
\begin {itemize}
\item  $\eps_i= \eps'_i = 1$ and $s_i > b^\delta$, or
\item $\eps_i =1, \eps_i' =0$, and $s_i < -b^\delta$. 
\end {itemize}
\end {lemma}

\begin {proof}
For the first part of the statement, it suffices to make sure that $b^\delta \geq b$, where $b$ is the constant from Lemma~\ref{lemma:isos}. For the second statement (about the triviality of the respective maps), observe that, for example, $\eps_i = \eps'_i =1$ means that $i \in I_-(\orL, \orN)$, and we are asked to show that $\Phi^{\orN}_{\psi^{M}(\s)}$ is zero. This is true because by \eqref{eq:Phi}, one of the factors of $\Phi^{\orN}_{\psi^{M}(\s)}$ is the map $\Pr^{\orN}_{\psi^M(\s)}$. This ``inclusion'' is the zero map because it contains a large power of $U_i$ see \eqref{eq:proj}, and that power is set to zero in the vertical truncation. The case $\eps_i =1, \eps_i' =0$ is similar.
\end {proof}

Let us fix some $\zeta = (\zeta_1, \dots, \zeta_\ell) \in \rr^\ell$ such that the values $\zeta_i$ are very close to zero, and linearly independent over $\qq$. We let $P_\rr(\Lambda) \subset \rr^\ell$ be the hyper-parallelepiped with vertices
$$ \zeta + \frac{1}{2}(\pm \Lambda_1 \pm \Lambda_2 \pm \dots \pm \Lambda_\ell),$$
for all possible choices of signs. This is a fundamental domain for $\rr^\ell/H(L, \Lambda)$, where $H(L,\Lambda)$ is the lattice generated by the vectors $\Lambda_i$, as in Section~\ref{subsec:surgery}. Let $P(\Lambda)$ be the collection of points in the lattice $\H(L)$ that also lie in $P_\rr(\Lambda)$. Because of our choice of $\zeta$, there are no lattice points on the boundary of $P_\rr(\Lambda)$. Therefore,
\begin {equation}
\label {eq:plambda}
P(\Lambda) \cong \H(L)/H(L, \Lambda). 
\end {equation}

In terms of the standard basis of $H_1(Y -L)$, we write
$$ \Lambda_i = (c_{i1}, \dots, c_{i \ell}),$$
where $c_{ii} =\lambda_i$ is the framing coefficient on the component $L_i$, and $c_{ij}$ is the linking number between $\orL_i$ and $\orL_j$ for $i\neq j$.

Pick $m_i \gg 0$, for $i=1, \dots, \ell$, and let $\tilde \Lambda$ be obtained from $\Lambda$ by increasing the framing coefficients by $m_i$; that is,
$$\tilde \Lambda_i = (\tilde c_{i1}, \dots, \tilde c_{i \ell}),$$
with
$$ \tilde c_{ij} = \begin {cases}
c_{ii} + m_i & \text{if } i=j, \\
c_{ij} & \text{if } i \neq j. 
\end {cases} $$

For each $\eps \in \E_\ell$, consider the hyper-parallelepiped $P(\tilde \Lambda, \Lambda, \eps)_\rr \subset \rr^\ell$ with vertices
$$ \zeta + \frac{1}{2} \sum_{i=1}^\ell \bigl(\sigma_i \tilde \Lambda_i + (1-\sigma_i)\eps_i\Lambda_i \bigr),$$
over all possible choices of signs $\sigma_i \in \{\pm1\}$. Set 
$$P(\tilde \Lambda, \Lambda, \eps) = P(\tilde \Lambda, \Lambda, \eps)_\rr \cap \H(L).$$ For $\eps = (0, \dots, 0)$, we recover the old $P(\tilde \Lambda)$, while for $\eps = (1, \dots, 1)$, the hyper-parallelepiped $P(\tilde \Lambda, \Lambda, \eps)$ is a rectangular hyperbox of size $(m_1, \dots, m_\ell)$.

Set
\begin {equation}
\label {eq:cdbd}
 \C^\delta\bd = \bigoplus_{\eps \in \E_\ell} \bigoplus_{\s \in P(\tilde \Lambda, \Lambda, \eps)} \C^{\eps, \delta}_\s.
 \end {equation}

\begin{figure}
\begin{center}
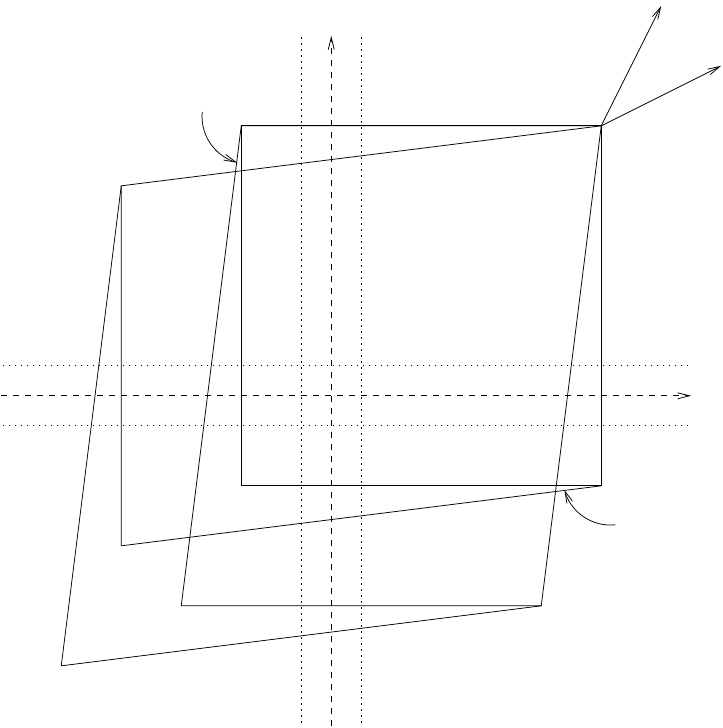
\end{center}
\caption {{\bf The complex $\C^\delta\bd$ for a two-component link.} We show here the combined truncation corresponding to surgery on a framed link with, say, $\Lambda_1 = (100, 200), \Lambda_2=(200, 100), m_1=m_2 = 600$ and $b^\delta = 50$. The four parallelograms in the picture are $P(\tilde \Lambda, \Lambda, \eps)$ for $\ell=2,\ \eps \in \{0,1\}^2$, and $\tilde \Lambda, \Lambda$ fixed. In each region we mark the values of $\eps$ such that the corresponding  $P(\tilde \Lambda, \Lambda, \eps)$ contains the region. This tells us which terms $\C^{\eps, \delta}_\s$ form the truncated complex $\C^\delta\bd$. The framing $\tilde \Lambda$ is sufficiently large compared to $b^\delta$, as explained in the proof of Proposition~\ref{prop:deltatruncate}.}
\label{fig:Ps}
\end{figure}

See Figure \ref{fig:Ps} for an example.

\begin {remark}
\label {rem:mmm}
The hyper-parallelepipeds $P(\tilde \Lambda, \Lambda, \eps)$ behave nicely with respect to the maps $\psi^M: \H(L) \to \H(L-M)$ as defined in Section~\ref{sec:reduction}. Indeed, consider the link 
$$M = \bigcup_{\{i|\eps_i = 1\}} L_i $$ with the orientation induced from $\orL$. Then $\psi^M$ takes $P(\tilde \Lambda, \Lambda, \eps)$ exactly to the hyper-parallelepiped $P(\tilde \Lambda|_{L-M})$, and it does so in an $m(M)$-to-one fashion, where
$$ m(M) = \prod_{\{i| L_i \subseteq M\}} m_i.$$
Therefore, we can re-write \eqref{eq:cdbd} as 
\begin {equation}
\label {eq:cdbd2}
 \C^\delta\bd = \bigoplus_{M \subseteq L} \bigoplus_{\s \in  P(\tilde \Lambda|_{L-M})} \oplus^{m(M)} \bigl( \Chain^{-, \delta}(\Hyper^{L-M}, \s) \bigr).
 \end {equation}
\end {remark}

\begin {proposition}
\label {prop:deltatruncate}
Fix $\delta > 0$. If we pick the values $m_i$ sufficiently large, the direct sum $\C^\delta\bd$, equipped with the restriction of the differential $\D^{-, \delta}$, forms a chain complex quasi-isomorphic to $(\C^\delta, \D^{-, \delta})$.
\end {proposition}

\begin {proof}
Let $P_\rr$ be a hyper-parallelepiped in $\rr^\ell$, with vertices $V^\sigma$, for $\sigma = (\sigma_1, \dots, \sigma_\ell) \in \{-1,1\}^\ell$. We assume that $\sigma_i V^\sigma_i > 0$ for all $\sigma$ and $i$; that is, each vertex lies in the hyper-quadrant in $\rr^\ell$ that corresponds to $\sigma$. For each
$$ \omega = (\omega_1, \dots, \omega_\ell) \in \{-1,0,1\}^\ell,$$
we define a subset $P_\rr[\omega] \subset \rr^\ell$ as follows. First, define a completion of $\omega$ to be a vector $\sigma \in \{-1, 1\}^\ell$ such that $\sigma_i = \omega_i$ whenever $|\omega_i | =1$. In other words, a completion of $\omega$ is a vector in which we replace the zero entries in $\sigma$ with $+1$ or $-1$. Let $P_\rr^\omega$ be the sub-parallelepiped of $P_\rr$ with vertices $V^\sigma$, where $\sigma$ runs over all possible completions of $\omega$. Further, given a vector with only nonnegative entries $\ttt= (t_1, \dots, t_\ell) \in [0, \infty)^\ell$, we define $\ttt * \omega$ to be the pointwise product
$$\ttt * \omega = (t_1\omega_1, \dots, t_\ell \omega_\ell) \in \rr^\ell. $$

Set:
$$ P_\rr[\omega] = \bigcup_{\ttt \in [0, \infty)^\ell} \bigl(P_\rr^\omega + \ttt * \omega\bigr).$$

Note that $P_\rr[(0, \dots, 0)] = P_\rr$. We have a decomposition:
$$ \rr^\ell = \bigcup_{\omega \in \{-1,0,1\}^\ell}  P_\rr[\omega].$$

See Figure~\ref{fig:nine} for the case $\ell = 2$. 

\begin{figure}
\begin{center}
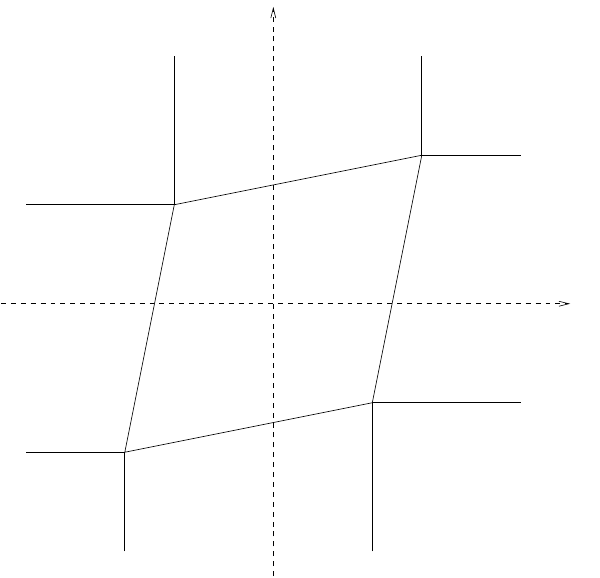
\end{center}
\caption {{\bf A decomposition of the plane into nine regions.}  
We show here an example of the subsets $ P_\rr[\omega] \subset \rr^2$, for $\omega \in \{-1,0,1\}^2$.
}
\label{fig:nine}
\end{figure}

Suppose there are no points in the lattice $\H(L)$ on the boundary of any $P_\rr[\omega]$. Then, letting $P[\omega] = P_\rr[\omega] \cap \H(L)$, we get a decomposition as a disjoint union
$$ \H(L) = \coprod_{\omega \in \{-1,0,1\}^\ell}  P[\omega].$$

We can apply this to any of the hyper-parallelepipeds $P(\tilde \Lambda, \Lambda, \eps)$ and obtain a decomposition of $\H(L)$ for each $\eps$. 

For $\omega \in \{-1,0,1\}^\ell$, consider the direct sum
$$\C^\delta\bd[\omega] = \bigoplus_{\eps \in \E_\ell} \prod_{\s \in P(\tilde \Lambda, \Lambda, \eps)[\omega]} \C^{\eps, \delta}_\s,$$
so that
\begin {equation}
\label {eq:cdsum}
\C^\delta = \bigoplus_{\omega \in \{-1,0,1\}^\ell}  \C^\delta\bd[\omega].
\end {equation}

We choose the values $m_i$ sufficiently large such that whenever $V$ is a vertex of some $P(\tilde \Lambda, \Lambda, \eps)$ with $\sigma_i V_i > 0$, we in fact have $ \sigma_i V_i > b^\delta$. (Here $b^\delta$ is the constant in Lemma~\ref{lemma:bdelta}.)

Starting from here, the idea is to use Lemma~\ref{lemma:bdelta} to show that the terms $\C^\delta\bd[\omega]$ in the direct sum \eqref{eq:cdsum} produce acyclic complexes for $\omega \neq (0, \dots 0)$. We eliminate them one by one, beginning with the ones with $|\omega| = \sum \omega_i = \ell$, then those with $|\omega|=\ell -1$, all the way to $|\omega|=1$, after which we are left only with $\C^\delta\bd[(0, \dots, 0)] = \C^\delta\bd$.

We start with $\omega$ such that  $|\omega| = \ell$, i.e. all values $\omega_i$ are $1$ or $-1$. Let $\varpi = (\varpi_1, \dots, \varpi_\ell)$ be the vector consisting of the values $\varpi_i = (1-\omega_i)/2$ for all $i$. Each $P(\tilde \Lambda, \Lambda, \eps)[\omega]$ is obtained from the hyper-quadrant determined by $\omega$ by translating by an amount that depends on $\eps$. More precisely, if $V$ is the unique vertex of  $P(\tilde \Lambda, \Lambda, 0)[\omega]$, then the vertex of $P(\tilde \Lambda, \Lambda, \eps)[\omega]$  is $V + (\varpi * \eps) \cdot \Lambda$.

The second part of Lemma~\ref{lemma:bdelta} implies that the maps $\De^{\eps, \eps', \delta}_{\eps^0, \s}$ are trivial unless $\eps' = \varpi * \eps$. It follows from here that $\C^\delta\bd[\omega]$ is a direct summand of $\C^\delta$ as a chain complex, i.e. it is preserved by the differential. Moreover, since by definition $\De^{\eps, \eps', \delta}_{\eps^0, \s}$ maps $\C^{\eps^0, \delta}_\s$ into $\C^{\eps^0 + \eps, \delta}_{\s+ \eps' \cdot \Lambda}$, we deduce that the differential on $\C^\delta\bd[\omega]$ preserves the quantity 
$$ \s - (\varpi * \eps) \cdot \Lambda.$$

Hence, the complex $\C^\delta\bd[\omega]$ splits into a direct sum of terms corresponding to values of $ \s - (\varpi * \eps) \cdot \Lambda$. Each such term is an $\ell$-dimensional hypercube that has on its edges maps of the form  $\De^{\eps, \eps', \delta}_{\eps^0, \s}$ with $\eps = \tau_i$. By the first part of Lemma~\ref{lemma:bdelta}, all these edge maps are quasi-isomorphisms. Therefore, the respective hypercube complexes are acyclic, and so is $\C^\delta\bd[\omega]$. 

We then proceed inductively on $\ell - |\omega|$. At each stage, we have a complex $\C^\delta\bd[\leq |\omega|]$ obtained from $\C^\delta$ by cancelling the terms with higher $|\omega|$. Each $\C^\delta\bd[\omega]$, for the given $|\omega|$, splits as a direct sum of subcomplexes and quotient complexes of $\C^\delta\bd[\leq |\omega|]$. They are acyclic by an application of Lemma~\ref{lemma:bdelta}. The claim follows.\end {proof}

\subsection {A folded truncation}
\label {sec:twt}
In this section we describe the variant of truncation which will be useful to us in Section~\ref{sec:proof}.
We keep the notation from the previous subsection, with $\C^\delta = \C^-(\Hyper, \Lambda)^\delta$ and $\tilde \Lambda$ as before. The folded truncation $\C^\delta\twist$ that we are about to describe will be isomorphic to $\C^\delta\bd$ as an $\Ring$-module, but its differential will contain  some additional maps.

For $\eps = \eps(M) \in \E_\ell$ and $\s \in \H(L)$, we introduce the notation $[\s]_{\eps}$ to denote the 
translate of $\s$ that lies in the hyper-parallelepiped $P(\tilde \Lambda, \Lambda,  \eps)$, where we allow translations by multiples of $\tilde \Lambda_i - \eps_i \Lambda_i$ for any $i=1, \dots, \ell$. In other words, we consider the tesselation of $\rr^\ell$ by translates of  $P(\tilde \Lambda, 
\Lambda, \eps)$, and ``fold'' all the hyper-parallelepipeds in this tesselation onto the fundamental domain $P(\tilde \Lambda, \Lambda, \eps)$ using translation. For simplicity, we write $\C^{\eps}_{[\s]}$ for the complex $\C^{\eps}_{[\s]_{\eps}}.$

Observe that if $\s \in \H(L)$ is sufficiently close to a face of $P(\tilde \Lambda, \Lambda,  \eps)$ (say, at a distance of no more than the sum of the absolute values of the entries of $\Lambda$) and the values $m_i$ were chosen sufficiently large, then $[\s]_{\eps}$ is also close to a (possibly different) face of $P(\tilde \Lambda, \Lambda,  \eps)$. Further, in this case in $\C^{\eps}_{\s} 
= \Chain^-(\Hyper^{L-M}, \psi^M(\s))$ all entries in $\psi^M(\s)$ are either sufficiently large or sufficiently small so that they can be replaced with $\pm \infty$, so $\C^{\eps}_{\s}$ is naturally isomorphic to $\Chain^-(\Hyper^\emptyset)$, using the data in the complete system $\Hyper$. The same is true of $\C^{\eps}_{[\s]}$, so we have an identification between $\C^{\eps}_{\s}$ and $\C^{\eps}_{[\s]}$.

With this in mind, we define
$$ \C\twist = \bigoplus_{\eps \in \E_\ell} \bigoplus_{\s \in P(\tilde \Lambda, \Lambda, \eps)} \C^{\eps }_\s,$$
 with the differential $\D^-\twist$ consisting of all maps
\begin {equation}
\label {eq:nouv}
 \D^{\eps, \eps'}_{\eps^0, \s}\twist = \Phi^{\orN}_{\psi^M(\s)} : \C^{\eps^0}_\s \to \C^{\eps^0 + \eps}_{[\s + \eps' \cdot \Lambda]}
 \end {equation}
for $\eps^0=\eps(M), \eps=\eps(N), \eps' \in \E_\ell$ such that $\eps^0 + \eps \in \E_\ell$ and $\eps' \leq \eps.$ 

Of course, it may happen that the value $\s$ is in the hyper-parallelepiped $P(\tilde \Lambda, \Lambda, \eps^0)$, but $\s + \eps' \cdot \Lambda$ is not in $P(\tilde \Lambda, \Lambda, \eps^0+ \eps)$. In such cases, note that $\s$ must be close enough to a face of $P(\tilde \Lambda, \Lambda, \eps^0+ \eps)$; hence, the target of the map  $ \D^{\eps, \eps'}_{\eps^0, \s}$, which is a priori $\C^{\eps^0 + \eps}_{\s + \eps' \cdot \Lambda}$, is identified with $C^{\eps^0 + \eps}_{[\s + \eps' \cdot \Lambda]}$, as explained above. The maps $ \D^{\eps, \eps'}_{\eps^0, \s} $ of this form are called {\em crossover maps}. 

We let $\C^{\delta}\twist$ be the vertical truncation (by $\delta$) of $\C\twist$. Then $\C^{\delta}\twist$ differs from the truncation  $\C^\delta\bd$ only in that its differential contains the crossover maps.

\begin{figure}
\begin{center}
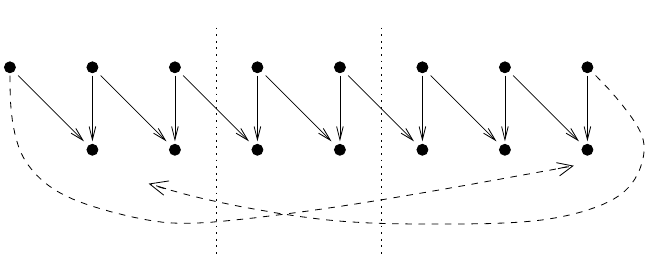
\end{center}
\caption {{\bf The folded truncation for surgery on a knot $K$.}  in this example $\Lambda =1$ and $b=8$. Each black dot represents a complex $\C^\eps_s$, each vertical arrow represents a map $\Phi^{+K}_s$, and each diagonal arrow represents a map $\Phi^{-K}_s$. The dashed arrows are crossover maps that appear in $\C \{2b\}$. They are set to zero in $\C^{\delta} \{2b \}$ because they originate outside the region between the two dotted lines---compare Lemma~\ref{lemma:bdelta}. }
\label{fig:crknot}
\end{figure}

\begin {example}
\label{ex:knott}
If $K$ is a knot with framing coefficient $\Lambda \in \zz$, recall from Equation~\eqref{eq:knottruncate} from Section~\ref{sec:truncate} that we have a horizontal truncation 
$$\C\langle b \rangle  =   \bigoplus_{-b \leq s \leq b} \C^0_s \oplus \bigoplus_{-b+\Lambda \leq s \leq b} \C^1_s  \subset \C.$$
If we also do vertical truncation by $\delta$, the resulting complex $\C^\delta \langle b \rangle$ is the same as the combined truncation $\C^\delta \bd$ from Section~\ref{sec:another}, with $\tilde \Lambda = 2b$. When $\Lambda = 1$, the truncation $\C^\delta \langle b \rangle = \C^\delta \langle \langle 2b \rangle \rangle$ is shown in Figure~\ref{fig:crknot}. In principle, in the folded truncation $\C^\delta \{ 2b \}$ we have to add the crossover maps $\Phi^{+K}_{-b}$ and $\Phi^{-K}_b$ shown by  dashed lines in Figure~\ref{fig:crknot}. (In our notation, for example, $\Phi^{+K}_{-b}$ is $\D^{\eps, \eps'}_{\eps^0, \s}\twist$ with $\eps^0 = 0, \eps=1, \eps'=0$ and $\s= -b$, so that $[\s+\eps'\cdot \Lambda]_{\eps^0+\eps} = [\s]_{\eps^0+\eps} = b$.)
However, observe that if $b$ is chosen sufficiently large compared to $b^\delta$, then the maps $\Phi^{+K}_{-b}$ and $\Phi^{-K}_b$ vanish, according to Lemma~\ref{lemma:bdelta}. Thus, in this case the folded truncation is the same as the combined truncation. (In fact, the same is true for any $\Lambda$, that is, for any integer surgery on a knot.)
\end {example}

\begin{figure}
\begin{center}
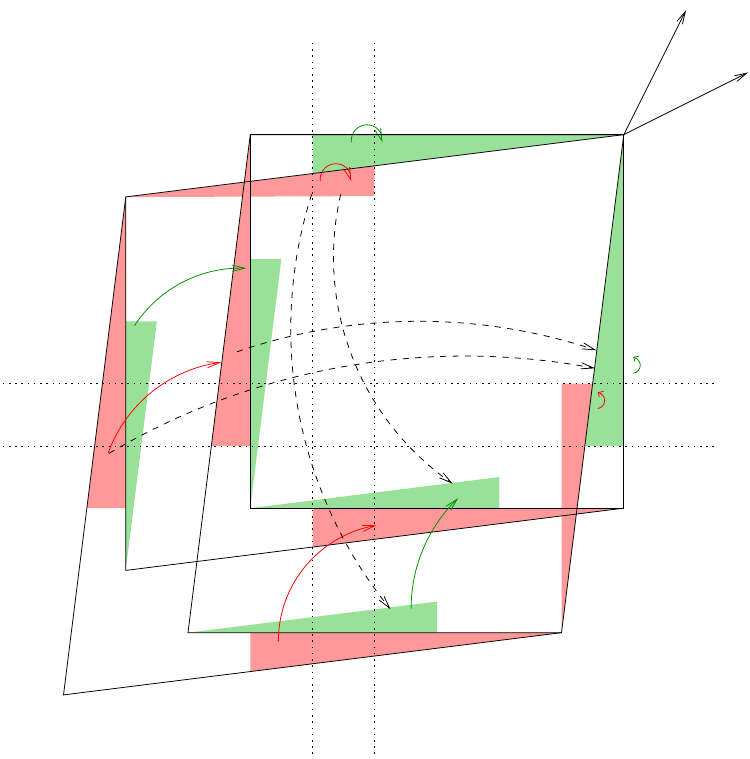
\end{center}
\caption {{\bf The folded truncation for surgery on a two-component link.}  We show here the same parallelograms as in Figure~\ref{fig:Ps}, corresponding to surgery on a link with $\Lambda_1=(200,100)$ and $\Lambda_2=(100, 200)$. Examples of crossover maps are given by dashed lines. They go from the quotient complex $ \C^\delta_{\Out}\twist$ (indicated by the red regions) to the subcomplex  $ \C^\delta_{\In}\twist$ (indicated by the green regions). In addition to the four crossover maps shown by dashed lines, there are four others: two of the form $\Phi^{+L_1}$ and two of the form $\Phi^{-L_2}$ (not shown here in order to avoid overcrowding the picture). In each green or red region we mark the values $\eps$ such that $\C^{\eps, \delta}_\s$ is part of the respective quotient complex or subcomplex. Further, near each green region we mark (in green square brackets) the edge of the parallelogram that blocks the  points in that region: $[1, +], [1, -], [2, +]$, or $[2, -]$---see the proof of  Proposition~\ref{prop:foldtruncate} for an explanation. The subcomplex $\C^\delta_{\In}\twist$ is acyclic because its differential contains quasi-isomorphisms (indicated by the green arrows) that cancel its terms in pairs. Similarly, $ \C^\delta_{\Out}\twist$ is acyclic because of the quasi-isomorphisms indicated by the red arrows.  
}
\label{fig:crossover}
\end{figure}

\begin {example}
When $L$ is a link of two or more components, there exist crossover maps that remain nontrivial even when vertically truncating by $\delta$. Examples are shown in Figure~\ref{fig:crossover}. For concreteness, consider the dashed arrow $\Phi^{+L_2}$ starting at the leftmost red triangle in the figure, and ending at the rightmost green triangle. In our notation, this is the crossover map $\D^{\eps, \eps'}_{\eps^0, \s}\twist$ for $\eps^0 = (0,0), \eps=(0,1), \eps'=(0,0)$, with $\s$ in the leftmost red triangle, and $[\s]_{\eps^0+\eps}$ in the rightmost green triangle. 
\end {example}

For future use, let us also write the definition of $\C^{\delta}\twist$ in a slightly different way. For $\eps=\eps(M) \in \E_\ell$, let $\TR^\eps$ be the (finite) polynomial ring over $\ff$ in variables $T_i$ for those $i$ such that $\eps_i = 1$ (that is, $L_i \subseteq M$), with relations $T_i^{m_i} = 1$. If $M = L_{i_1} \cup \dots \cup L_{i_p}$ and $\kk = (k_1, \dots, k_p) \in \zz^p$, we formally write
$$  T^\kk = \prod_{j=1}^p T_{i_j}^{k_j}.$$

Fix a reference value $\s^0=(s^0_1, \dots, s^0_\ell) \in \H(L)$ with all of its entries less than $b^\delta$ in absolute value. In view of Remark~\ref{rem:mmm}, we have an identification
\begin {eqnarray*}
 P(\tilde \Lambda, \Lambda,  \eps) &\xrightarrow{\phantom{blai} \cong \phantom{bla}}& P(\tilde \Lambda|_{L-M}) \times  \TR^\eps \\ \s = (s_1, \dots, s_\ell) &\xrightarrow{\phantom{blablabla}}& \Bigl( \psi^{M}(\s), \prod_{\{i\mid \eps_i = 1\}} T_i^{s_i - s^0_i} \Bigr).
\end {eqnarray*}

Therefore, we can write
\begin {equation}
\label {eq:cdfolded}
 \C^\delta\twist = \bigoplus_{M \subseteq L} \bigoplus_{\s \in  P(\tilde \Lambda|_{L-M})}   \Chain^{-, \delta}(\Hyper^{L-M}, \s)  \otimes \TR^\eps(M),
 \end {equation}
with the differential made of terms
\begin {equation}
\label {eq:Dfolded}
 \D^{\eps, \eps', \delta}_{\eps^0, \s}\twist =  T^{  (\Lambda|_N)_{\orL, \orN}} \cdot  \Phi_{\psi^M(\s)}^{\orN, \delta},
 \end {equation}
where $\eps^0=\eps(M), \eps = \eps(N)$, and the orientation $\orN$ on $N$ is given by $i \in I_-(\orL, \orN) \iff \eps'_i = 1$. By $(\Lambda|_N)_{\orL, \orN}$ we mean the image of $\Lambda_{\orL, \orN} = \sum_{ i \in I_-(\orL, \orN)} \Lambda_i \in H_1(Y-L)$ under the natural map $ H_1(Y-L) \to H_1(Y-N)$.

Our main goal is to show:

\begin {proposition}
\label {prop:foldtruncate}
Fix $\delta > 0$. If we pick the values $m_i$ sufficiently large, then the complexes $\C^\delta$ and $\C^\delta\twist$ are quasi-isomorphic.
\end {proposition}

The proof of Proposition~\ref{prop:foldtruncate} will proceed along the following lines. Recall from Proposition~\ref{prop:deltatruncate} that $\C^\delta$ is quasi-isomorphic to the truncated complex $\C^\delta\bd$. Thus, it suffices to show that $\C^\delta\bd$ and $\C^\delta\twist$ are quasi-isomorphic. The basic idea is that all nontrivial crossover maps go into an acyclic subcomplex of $\C^\delta\twist$; therefore, deleting them will not not change the quasi-isomorphism type of the complex.

In order to make this argument precise, we first introduce the  auxiliary complex we need, and show that it is acyclic.

Note that many of the crossover maps in $\C^\delta\twist$ are automatically zero, by an application of Lemma~\ref{lemma:bdelta}. In particular, when $L$ is a knot, all crossover maps vanish---compare Example~\ref{ex:knott} and Figure~\ref{fig:crknot}. On the other hand, in Figure~\ref{fig:crossover}, some crossover maps have to vanish and others do not: for example, near the top right corner of Figure~\ref{fig:crossover}, in principle we should have crossover maps of the form $\Phi^{-L_1}$ and $\Phi^{-L_2}$; but these are zero, so we did not show them in the figure.

Roughly,  the auxiliary complex $\C^\delta_{\In}\twist$ (to be defined below) will be the target of all crossover maps that are not being automatically set to zero by applying Lemma~\ref{lemma:bdelta}.
In Figure~\ref{fig:crossover}, this complex is indicated by the green regions. (When $L$ is a knot, the auxiliary complex is zero.)

Precisely, we define an {\em incoming crossover point} to be a pair $(\eps^0, \s)$ with $\eps^0 \in \E_\ell$ and $\s=(s_1, \dots, s_\ell) \in P(\tilde \Lambda, \Lambda, \eps^0)$ satisfying the following property: There exist $\eps, \eps' \in \E_\ell$ such that
\begin {itemize}
 \item $\eps^0 - \eps \in \E_\ell$ and $ \eps' \leq \eps;$
\item $\s' = (s'_1, \dots, s_\ell') := \s - \eps' \cdot \Lambda \not \in P(\tilde \Lambda, \Lambda, \eps^0 - \eps)$; 
\item there is no $i=1, \dots, \ell,$ with $\eps_i =1, \eps_i'=0$ and $s'_i < - b^{\delta}$;
\item there is no $i=1, \dots, \ell,$ with $\eps_i =1, \eps_i'=1$ and $s'_i + \sum_{j \neq i} \lk(L_i, L_j) (1-\eps_j) > b^{\delta}$.
\end {itemize}

In other words, we ask for $(\eps^0, \s)$ to specify the target of a (not a priori trivial) crossover map
$$  \D^{\eps, \eps'}_{\eps^0-\eps', \s} : \C^{\eps^0-\eps}_{\s'} \to \C^{\eps^0}_{\s}.$$

We denote by $\In(\tilde \Lambda, \Lambda) $ the set of all incoming crossover points. (In Figure~\ref{fig:crossover}, the set $ \In(\tilde \Lambda, \Lambda) $ is exactly the union of the green regions.) Define the {\em incoming crossover subcomplex} of $\C^\delta\twist$ to be
$$ \C^\delta_{\In}\twist = \bigoplus_{(\eps, \s) \in  \In(\tilde \Lambda, \Lambda)} \C^{\eps, \delta}_\s.$$

\begin {lemma}
\label {lemma:AcyclicSubcx}
The subcomplex $ \C^\delta_{\In}\twist$ is acyclic. 
\end {lemma}

\begin {proof}
Recall that the crossover maps must end near the faces of the hyper-parallelepipeds under consideration. Near each face Lemma~\ref{lemma:bdelta} indicates the presence of a particular kind of quasi-isomorphism among the edge maps $\Phi^{\pm L_i}_\s$. We claim that these quasi-isomorphisms are part of the differential for $ \C^\delta_{\In}\twist$, and cause its terms to cancel out in pairs. (In Figure~\ref{fig:crossover}, the quasi-isomorphisms that produce the cancellation are indicated by the green arrows.)

To make this precise, let  $(\eps^0, \s)$ be an incoming crossover point. The condition $\s - \eps' \cdot \Lambda \not \in P(\tilde \Lambda, \Lambda, \eps^0 - \eps)$ means that $\s - \eps' \cdot \Lambda$ fails to satisfy at least one of the $2\ell$ inequalities describing the hyper-parallelepiped $P(\tilde \Lambda, \Lambda, \eps^0 - \eps)$. This inequality comes from one of the $2\ell$ codimension one faces of the hyper-parallelepiped, and we say that $(\eps^0, \s)$ is an incoming crossover point {\em blocked} by that face. One way to distinguish one of the $2\ell$ codimension one faces is to say which of the $2\ell$ coordinate half-axes of $\rr^\ell$ it intersects. In general, if $P$ is a hyper-parallelepiped like the ones in our construction, the face of $P$ that intersects the positive  $s_i$ axis will be called the $[i, +]$ face, and the face  that intersects the negative  $s_i$ axis will be called $[i, -]$ face. 

Suppose an inequality that fails for $\s - \eps' \cdot \Lambda$ is the one defining the $[1, +]$ face of $P(\tilde \Lambda, \Lambda, \eps^0 - \eps)$. Since $\s$ does satisfy the inequality defining the $[1, +]$ face of $P(\tilde \Lambda, \Lambda, \eps^0)$, it must be that $\s$ lies close to those faces---where by ``close'' we mean a distance comparable to the entries of $\Lambda$ and much smaller than $b^\delta$ or the $m_i$'s. (An example is when $\s$ is in the rightmost vertical green triangle in Figure~\ref{fig:crossover}, and $\eps^0 = (0,1)$ or $(1,1)$; there, $\eps = (0,1)$ and $\eps' = (0,0)$, so $\Phi^{+L_2}_{[\s']}$ is the crossover map.) Suppose further that $\eps^0_1=0$, and let $\hat{\eps}^0 \in \E_\ell$ be obtained from $\eps^0$ by changing the entry $\eps^0_1$ from $0$ to $1$. Changing $\eps^0$ into $\hat{\eps}^0$ does not change the $[1, +]$ face of $P(\tilde \Lambda, \Lambda, \eps^0)$ at all. The same is true for the $[1, +]$ face of $P(\tilde \Lambda, \Lambda, \eps^0 - \eps)$. Therefore, $(\s, \eps^0)$ is an incoming crossover point blocked by the $[1, +]$ face if and only if $(\s, \hat{\eps}^0)$ is. Moreover, since we are near the $[1, +]$ face, the value $s_1$ is bigger than $b^\delta$, so Lemma~\ref{lemma:bdelta} says that  
$$\Phi^{+L_1}_\s: \C^{\eps^0, \delta}_{\s} \longrightarrow \C^{\hat{\eps}^0, \delta}_{\s }$$  is a quasi-isomorphism. (This is shown as one of the green arrows in Figure~\ref{fig:crossover}.) Consider the subcomplex of $ \C^\delta_{\In}\twist$ consisting of $\C^{\eps^0, \delta}_\s$ for all $(\eps^0, \s) \in  \In(\tilde \Lambda, \Lambda)$ blocked by the $[1, +]$ face. This subcomplex is acyclic, because it consists of two hypercubes of chain complexes (corresponding to $\eps^0_1 =0$ and $\eps^0_1=1$, respectively), connected by a chain map which has quasi-isomorphisms $\Phi^{+L_1}_\s$ along the edges. We can kill off this acyclic subcomplex and obtain a new complex quasi-isomorphic to $ \C^\delta_{\In}\twist$, which does not involve any $(\s, \eps^0)$ blocked by the $[1, +]$ face. Similarly, from this new complex we can kill off the terms blocked by the $[2, +]$ face, then the terms blocked by the $[3, +]$ face, and so on, all without changing the quasi-isomorphism type.

An analogous elimination procedure works for the terms blocked by faces of the type $[i, -]$. Indeed, suppose an inequality that fails for $\s - \eps' \cdot \Lambda$ is the one defining the $[1, -]$ face of $P(\tilde \Lambda, \Lambda, \eps^0 - \eps)$. Since $\s$ does satisfy the corresponding inequality for $P(\tilde \Lambda, \Lambda, \eps^0)$, it must be that $\s$ lies close to those faces. (An example of this situation is when $\s$ is in one of the two leftmost vertical green triangles in Figure~\ref{fig:crossover}, and $\eps^0 = (0,1)$ or $(1,1)$; there, $\eps = \eps' = (0,1)$, so $\Phi^{-L_2}_{[\s']}$ is the crossover map.) Suppose further that $\eps^0_1=0$, and let $\hat{\eps}^0 \in \E_\ell$ be obtained from $\eps^0$ by changing the entry $\eps^0_1$ from $0$ to $1$. This time, changing $\eps^0$ into $\hat{\eps}^0$  results in translating the $[1, -]$ face of $P(\tilde \Lambda, \Lambda, \eps^0)$ by the vector $\Lambda_1$. The same happens to the $[1, -]$ face of $P(\tilde \Lambda, \Lambda, \eps^0 - \eps)$. Therefore, $(\hat{\eps}^0, \s + \Lambda_1)$ is also an incoming crossover point blocked by the $[1, -]$ face.  Further, since $s_1 \ll -b^\delta$, Lemma~\ref{lemma:bdelta} says that 
$$\Phi^{-L_1}_\s: \C^{\eps^0, \delta}_{\s} \longrightarrow \C^{\hat{\eps}^0, \delta}_{\s + \Lambda_1}$$ 
is a quasi-isomorphism. Consequently, the terms in $ \C^\delta_{\In}\twist$ that are blocked by the $[1, -]$ face cancel each other out in pairs (using $\Phi^{-L_1}_\s$). By a similar argument, we can kill off the remaining terms blocked by the $[2, -]$ face, and so forth. We conclude that the whole of $ \C^\delta_{\In}\twist$ is acyclic.
\end {proof}

It is worthwhile to note that we can similarly define an {\em outgoing crossover point} to be a pair $(\eps^0, \s)$ with $\eps^0 \in \E_\ell$ and $\s=(s_1, \dots, s_\ell) \in P(\tilde \Lambda, \Lambda, \eps^0)$ satisfying the following property: There exist $\eps, \eps' \in \E_\ell$ such that
\begin {itemize}
 \item $\eps^0 + \eps \in \E_\ell$ and $ \eps' \leq \eps;$
\item $\s + \eps' \cdot \Lambda \not \in P(\tilde \Lambda, \Lambda, \eps^0 + \eps)$; 
\item there is no $i=1, \dots, \ell,$ with $\eps_i =1, \eps_i'=1$ and $s_i  > b^{\delta}$;
\item there is no $i=1, \dots, \ell,$ with $\eps_i =1, \eps_i'=0$ and $s_i + \sum_{j \neq i} \lk(L_i, L_j) (1-\eps_j^0) < - b^{\delta}$.
\end {itemize}

We denote by $ \Out(\tilde \Lambda, \Lambda) $ the set of all outgoing crossover points. (In Figure~\ref{fig:crossover}, the set $ \Out(\tilde \Lambda, \Lambda) $ is the union of the red regions.) Define the {\em outgoing crossover quotient complex} of $\C^\delta\twist$ to be
$$ \C^\delta_{\Out}\twist = \bigoplus_{(\eps, \s) \in  \Out(\tilde \Lambda, \Lambda)} \C^{\eps, \delta}_\s,$$
with the differential induced from $ \C^\delta_{\Out}\twist$. Roughly,  $\C^\delta_{\Out}\twist$ is the domain of all crossover maps that are not being automatically set to zero by applying Lemma~\ref{lemma:bdelta}. One can show that $\C^\delta_{\Out}\twist$ is acyclic by an argument similar to the one used for $ \C^\delta_{\In}\twist$. However, we will not need to use the complex $\C^\delta_{\Out}\twist$ in the proof below.

\begin {proof}[Proof of Proposition~\ref{prop:foldtruncate}]
Let
$$ \C^\delta\dtwist = \bigoplus_{\substack{
\eps \in \E_\ell, \ \s \in P(\tilde \Lambda, \Lambda, \eps) \\
(\eps, \s) \not\in \In(\tilde \Lambda, \Lambda)}} \C^{\eps, \delta}_\s$$
be the quotient complex of $\C^\delta\twist$ obtained from by quotienting out the subcomplex $\C^\delta_{\In}\twist$. Since $\C^\delta_{\In}\twist$ are acyclic, it follows that $\C^\delta\twist$ and $\C^\delta\dtwist$ are quasi-isomorphic.

Equally well, we can think of $ \C^\delta_{\In}\twist$ as a subcomplex of the truncation $\C^\delta\bd$ (the one without crossover maps). Thus, $\C^\delta\dtwist$ is a quotient of $\C^\delta\bd$, and quasi-isomorphic to it. We deduce that $\C^\delta\bd$ and $\C^\delta\twist$ are quasi-isomorphic.  The conclusion then follows from Proposition~\ref{prop:deltatruncate}.
\end {proof}

\section {A general surgery exact sequence}
\label {sec:exact}

Theorem~\ref{thm:surgery} is a generalization of  the corresponding result for knots, Theorem 1.1 in \cite{IntSurg}. The key component of the proof of that result is the existence of a long exact sequence:
\begin {equation}
\label {eq:extri}
 \dots \lra \HFm(Y_n(K)) \lra \HFm(Y_{m+n}(K)) \lra \oplus^m \HFm(Y) \lra \dots, \end {equation}
see \cite[Theorem 3.1]{IntSurg}, stated there with $\iHF^+$ instead of $\HFm$. 
(see also  \cite[Theorem 6.2]{RatSurg}, for 
a generalization to rationally null-homologous knots inside three-manifolds).
Here $n, m \in \zz$ are surgery coefficients of the knot $K \subset Y$ with $m > 0$, and $Y$ an integral homology sphere. The fact that we work with the completed version $\HFm$, see Section~\ref{sec:conventions}, allows the triangle to exist for $\HFm$ just as for $\iHF^+$. 

As a first step towards the proof of Theorem~\ref{thm:surgery}, in this section we will discuss a broader generalization of \eqref{eq:extri}, in which $K$ is an arbitrary knot in a closed, oriented three-manifold $Y$. In Sections~\ref{sec:4c} and ~\ref{sec:grtwisted} we establish a few useful facts about cobordism maps in Heegaard Floer homology. We then proceed to state the general exact sequence. Its proof in Section~\ref{sec:LES} is a straightforward adaptation of Theorem 3.1 of~\cite{IntSurg}, as long as one does not keep track of gradings and decompositions into $\spc$ structures. We deal with these two last issues at length in Sections~\ref{sec:case1}, \ref{sec:case2} and \ref{sec:x}, for the exact sequence with vertically truncated complexes.

\subsection {Four-dimensional cobordisms}
\label {sec:4c}

We spell out  here a way of describing homology classes, cohomology classes, and $\spc$ structures on four-dimensional cobordisms (given by two-handle additions) in terms of surgery. This will be useful to us when discussing the surgery long exact sequence later in this section, as well as when discussing surgery maps in Section~\ref{sec:beyond}. 

Let $(\orL, \Lambda)$ be a framed $\ell$-component link inside an
integral homology three-sphere $Y$. Recall that the 
 space of $\spc$
structures on $Y_{\Lambda}(L)$ is identified with the quotient $\H(L)/H(L,
\Lambda)$  (see Section~\ref{subsec:surgery}).
Also,
$H_1(Y_\Lambda(L)) \cong H^2(Y_{\Lambda}(L))$ is identified with $ \zz^\ell/ H(L,
\Lambda)$ and $H_2(Y_\Lambda(L)) \cong H^1(Y_\Lambda(L))$ with $H(L,
\Lambda)^\perp$
(see Section~\ref{sec:gradings}).

Let 
$L'   \subseteq L$ be an $\ell'$-component sublink, with the orientation induced from $\orL$. Denote by
$$H(L, \Lambda|_{L'}) \subset H_1(Y-L) \cong \zz^\ell$$ 
the sublattice generated by the framings $\Lambda_i$, for $L_i \subseteq L'$.

Let $W_{\Lambda}(L', L)$ be the  cobordism from $Y_{\Lambda|_{L'}}(L')$ to $Y_{\Lambda}(L)$ given by surgery on $L - L'$ (framed with the restriction of $\Lambda$). Further, let $W_{\Lambda}(L) = W_{\Lambda}(\emptyset, L)$, so that   
$$W_{\Lambda}(L)=W_{\Lambda|_{L'}}(L') \cup W_{\Lambda}(L', L).$$

\begin {lemma}
\label {lem:HH}

(a) There is a natural identification
\begin {equation}
\label {eq:homology}
H_2(W_{\Lambda}(L', L)) \cong H(L, \Lambda|_{L'})^\perp = \{ \vs \in \zz^\ell | \vs \cdot \Lambda_i = 0, \forall i \text{ with } L_i \subseteq L' \},
\end {equation}
under which the intersection product 
$$  H_2(W_{\Lambda}(L', L)) \otimes H_2(W_{\Lambda}(L', L)) \to \zz$$
is given by $ \vs \otimes \vs' \to \vs^t  \Lambda  \vs'$.

(b) There is a natural identification
\begin {equation}
\label {eq:cohomology}
H^2(W_{\Lambda}(L', L)) \cong \zz^\ell/H(L, \Lambda|_{L'}),
\end {equation}
under which the natural projection 
$$\pi^{L, L'} :  \zz^\ell/H(L, \Lambda|_{L'}) \longrightarrow \zz^\ell/H(L, \Lambda)$$
corresponds to restriction to $H^2(Y_{\Lambda}(L))$, while restriction to the coordinates corresponding to $L'$
$$  \zz^\ell/H(L, \Lambda|_{L'}) \to \zz^{\ell'}/H(L', \Lambda|_{L'}) $$
corresponds to restriction to $H^2(Y_{\Lambda|_{L'}}(L'))$.

(c) Under the identifications \eqref{eq:homology} and \eqref{eq:cohomology}, the evaluation map
$$ H^2(W_{\Lambda}(L', L)) \otimes H_2(W_{\Lambda}(L', L)) \to \zz$$
corresponds to the usual scalar multiplication of vectors in $\zz^\ell$. Also, the composition of Poincar\'e duality with the natural map in the long exact sequence of a pair:
$$ H_2(W_{\Lambda}(L', L)) \xrightarrow{\cong} H^2(W_{\Lambda}(L', L), \del W_{\Lambda}(L', L)) \longrightarrow H^2(W_{\Lambda}(L', L)) 
$$
corresponds to the multiplication $\vs \to \Lambda  \vs$. \end {lemma}

\begin {proof}
We start by proving the claims in the case when $L' = \emptyset$, i.e. $W_{\Lambda}(L', L) = W_{\Lambda}(L)$. Choose Seifert surfaces $F_i \subset Y$ for each link component $L_i$. (Of course, these surfaces may intersect each other.) Let $\hat F_i$ be the surface obtained by capping off $F_i$ in $W_{\Lambda}(L)$, using the core of the respective two-handle. Note that the homology class $[\hat F_i]$ is independent of our choice of $F_i$. Further, since $Y$ is an integral homology sphere, the classes $[\hat F_i], i=1, \dots, \ell$, form a basis of $H_2( W_{\tilde \Lambda}(L))$, so we have the desired identification
\begin {equation}
\label {eq:H_2}
H_2( W_{\Lambda}(L)) \cong \zz^\ell.
\end {equation}
Note that the intersection form on $H_2(W_{\Lambda}(L))$ in this basis is the framing matrix $\Lambda$.

We also have another identification:
\begin {equation}
\label {eq:H^2}
H^2( W_{ \Lambda}(L)) \cong \zz^\ell,
\end {equation}
obtained by sending a cohomology class $c$ to $(\langle c, [\hat F_1]\rangle, \dots, \langle c, [\hat F_\ell]\rangle)$. The claims in (c) are then easy to check for $W_{\Lambda}(L)$. 

Let us now consider the general case of $L' \subseteq L$. For part (a), observe that
$$ H_3(W_{\Lambda}(L), W_{\Lambda}(L', L)) \cong H_3(W_{\Lambda|_{L'}}(L'), Y_{\Lambda|_{L'}}(L')) = 0,$$
because $W_{\Lambda|_{L'}}(L')$ consists of two-handle additions only. Hence, the long exact sequence in homology for the pair $(W_{\Lambda}(L), W_{\Lambda}(L', L))$ reads
$$ 0 \to H_2(W_{\Lambda}(L', L)) \to  H_2(W_{\Lambda}(L)) \xrightarrow{f} H_2(W_{\Lambda}(L), W_{\Lambda}(L', L)) \cong H^2(W_{\Lambda|_{L'}}(L')).$$
Using \eqref{eq:H_2} and \eqref{eq:H^2} we can view $f$ as a map from $\zz^\ell$ to $\zz^{\ell'}$. From part (c) for $W_{\Lambda}(L)$ we get that $f$ is given in matrix form by the restriction of $\Lambda$ to the components of $L'$. The identification \eqref{eq:homology} follows.

For part (b), use the identifications \eqref{eq:H^2} for $W_{\Lambda}(L)$ and $W_{\Lambda|_{L'}}(L')$, as well as the commutative diagram
$$ \begin {CD}
0 @>>> H^2(W_{\Lambda}(L), W_{\Lambda}(L', L)) @>>> H^2(W_{\Lambda}(L)) @>>> H^2(W_{\Lambda}(L', L)) @>>> 0 \\
@. @VV{\cong}V @VVV  @VVV \\
0 @>>> H^2(W_{\Lambda|_{L'}}(L') , Y_{\Lambda|_{L'}}(L')) @>>> H^2(W_{\Lambda|_{L'}}(L')) @>>> H^2( Y_{\Lambda|_{L'}}(L')) @>>> 0
\end {CD}$$

Part (c) follows from the respective statements for $W_{\Lambda}(L)$.
\end {proof}

We also have a description of the space of $\spc$ structures on the cobordism  $W_{\Lambda}(L', L)$:

\begin {lemma}
\label {lem:spc4}
There is a natural identification:
\begin {equation}
\label {eq:idspc}
\spc ( W_{\Lambda}(L', L)) \cong \H(L)/H(L, \Lambda|_{L'})
\end {equation}
under which the natural projection 
$$\pi^{L, L'} :  \H(L)/H(L, \Lambda|_{L'}) \longrightarrow \H(L)/H(L, \Lambda)$$
corresponds to restricting the $\spc$ structures to $Y_{\Lambda}(L)$, while the map
$$ \psi^{L-L'}: \H(L)/H(L, \Lambda|_{L'}) \to \H(L')/H(L', \Lambda|_{L'})$$
corresponds to restricting them to $Y_{\Lambda|_{L'}}(L')$. Further, the first Chern class map
$$ c_1 : \spc ( W_{\Lambda}(L', L)) \to H^2( W_{\Lambda}(L', L)) \cong \zz^\ell / H(L, \Lambda|_{L'})$$
is given by
\begin {equation}
\label {eq:c1s}
c_1([\s]) = [2\s - (\Lambda_1 + \dots + \Lambda_\ell)].
\end {equation}
\end {lemma}

\begin {proof}
The space of $\spc$ structures on $W_{\Lambda}(L)$ is identified 
with the space of relative $\spc$ structures on $(Y, L)$, and hence with $\H(L)$ via the formula \eqref{eq:c1s} (see Remark~\ref{rem:h1}). Similarly, we have an identification $\spc(W_{\Lambda|_{L'}}(L')) \cong \H(L')$. Moreover, there is a commutative diagram
$$ \begin {CD}
0 @>>> \spc(W_{\Lambda}(L), W_{\Lambda}(L', L)) @>>> \spc(W_{\Lambda}(L)) @>>> \spc(W_{\Lambda}(L', L)) @>>> 0 \\
@. @VV{\cong}V @VVV  @VVV \\
0 @>>> \spc(W_{\Lambda|_{L'}}(L') , Y_{\Lambda|_{L'}}(L')) @>>> \spc(W_{\Lambda|_{L'}}(L')) @>>> \spc( Y_{\Lambda|_{L'}}(L')) @>>> 0
\end {CD}$$
where the two horizontal rows are short exact sequences. The conclusion easily follows. \end {proof}

\subsection{Gradings, cobordism maps, and twisted coefficients}
\label {sec:grtwisted}
We now discuss the general grading properties for Heegaard Floer complexes with twisted coefficients and cobordism maps between them. Let $Y$ be a closed, oriented 3-manifold with a $\spc$ structure $\ux$ and a $\Field[H^1(Y; \zz)]$-module $M$. A pointed, admissible Heegaard diagram for $Y$ (with a complete set of paths for the generators, as in \cite[Definition 3.12]{HolDisk}) gives rise  to a twisted Heegaard Floer complex $ {\CFm}(Y, \ux; M)$, see \cite[Section 8.1]{HolDiskTwo}. 
(Here we use the completed version, and delete the usual underline from notation for simplicity.) 
This complex admits a relative $\zz/\delt(\ux, M)\zz$-grading, where
\begin {equation}
\label {eq:twistedgrading}
 \delt(\ux, M) = \gcd_{ \{\xi \in H_2(Y; \zz)| PD(\xi) \cdot m = m, \forall m \in M \}} \langle c_1(\ux), \xi\rangle.
 \end {equation}

This is true because the ambiguity in the grading difference between two generators is the Maslov index of periodic domains for which the corresponding classes in $H^1(Y; \zz)$ act trivially on $M$; compare \cite[Section 8.1]{HolDiskTwo}.

Next, we set up cobordism maps with twisted coefficients, following \cite[Section 8.2]{HolDiskTwo}, but in slightly more generality. 

Consider a pointed, admissible triple Heegaard diagram $(\Sigma, \alphas, \betas, \gammas, w)$ which represents a cobordism $X= X_{\alpha, \beta, \gamma}$ with boundaries $Y_{\alpha, \beta}, Y_{\beta, \gamma}$ and $Y_{\alpha, \gamma}$. Suppose we are given an $\ff[H^1(Y_{\alpha, \beta}; \zz)]$-module $M_{\alpha, \beta}$ and an $\ff[H^1(Y_{\beta, \gamma}; \zz)]$-module $M_{\beta, \gamma}$.  Given a $\spc$ structure $\t$ on $X$, we denote by $\underline{\spc}(X; \t)$ the space of relative $\spc$ structures on $X$ representing $\t$. Note that $\underline{\spc}(X; \t)$ has a natural action of $H^1(Y_{\alpha, \beta}; \zz)  \times H^1(Y_{\beta, \gamma}; \zz)  \times H^1(Y_{\alpha, \gamma}; \zz)$. The $\spc$ structure $\t$ induces a $\ff[H^1(Y_{\alpha, \gamma}; \zz)]$-module 
$$\{M_{\alpha, \beta} \otimes M_{\beta, \gamma}\}^\t = \frac{(m_{\alpha, \beta}, m_{\beta, \gamma}, \underline \t) \in M_{\alpha, \beta} \times M_{\beta, \gamma} \times \underline{\spc}(X; \t)}{(m_{\alpha, \beta}, m_{\beta, \gamma}, \underline \t) \sim (h_{\alpha, \beta} \cdot m_{\alpha, \beta}, h_{\beta, \gamma} \cdot m_{\beta, \gamma}, (h_{\alpha, \beta} \times h_{\beta, \gamma} \times 0) \cdot \underline \t) },$$ 
where $h_{\alpha, \beta}$ and $h_{\beta, \gamma}$ are arbitrary elements of $H^1(Y_{\alpha, \beta}; \zz)$ and $H^1(Y_{\beta, \gamma}; \zz)$, respectively.

Suppose we are also given a $\ff[H^1(Y_{\alpha, \gamma}; \zz)]$-module $M_{\alpha, \gamma}$ and a module homomorphism
$$ \zeta: \{M_{\alpha, \beta} \otimes M_{\beta, \gamma}\}^\t \longrightarrow M_{\alpha, \gamma}.$$

For simplicity, let us write $M$ for the triple $(M_{\alpha, \beta}, M_{\beta, \gamma}, M_{\alpha, \gamma})$. We then have a cobordism map with twisted coefficients:
$$ {f^{-}_{\alpha, \beta, \gamma; \t, M, \zeta}}: {\CFm}(Y_{\alpha, \beta}, \t|_{Y_{\alpha, \beta}}; M_{\alpha, \beta}) \otimes  {\CFm}(Y_{\beta, \gamma}, \t|_{Y_{\beta, \gamma}}; M_{\beta, \gamma}) \to {\CFm}(Y_{\alpha, \gamma}, \t|_{Y_{\alpha, \gamma}}; M_{\alpha, \gamma}), $$
given by
$$ f^{-}_{\alpha, \beta, \gamma; \t, M, \zeta}(m_{\alpha, \beta}\x \otimes m_{\beta \gamma}\y) = \sum_{\z \in \Ta \cap \Tg} \sum_{\bigl \{\phi \in \pi_2(\x, \y, \z) \big \vert \substack{\mu(\phi)=0 \\ \t_w(\phi) = \t}\bigr \}} (\# \M(\phi)) \cdot U^{n_w(\phi)} \zeta(m_{\alpha, \beta}\otimes m_{\beta \gamma} \otimes \bar\t_w(\phi)) \z,$$
with the usual notations in Heegaard Foer theory; compare \cite[Equation (10)]{HolDiskTwo}.

Set
\begin {equation}
\label {eq:deltam}
 \delt(\t, M) = \gcd \bigl(\delt(\t|_{Y_{\alpha, \beta}}, M_{\alpha, \beta}),  \delt( \t|_{Y_{\beta, \gamma}}, M_{\beta, \gamma}), \delt(\t|_{Y_{\alpha, \gamma}}, M_{\alpha, \gamma})  \bigr). 
 \end {equation}

Note that the Floer complexes 
$${\CFm}(Y_{\alpha, \beta}, \t|_{Y_{\alpha, \beta}}; M_{\alpha, \beta}),  {\CFm}(Y_{\beta, \gamma}, \t|_{Y_{\beta, \gamma}}; M_{\beta, \gamma}),  {\CFm}(Y_{\alpha, \gamma}, \t|_{Y_{\alpha, \gamma}}; M_{\alpha, \gamma})$$ all admit relative $\zz/\delt(\t, M)\zz$-gradings.

\begin {lemma}
\label{lem:twistedmaps}
(a) The map $ f^{-}_{\alpha, \beta, \gamma; \t, M, \zeta}$ preserves the relative $\zz/\delt(\t, M)\zz$-gradings.

(b) Let $\t'$ be another $\spc$ structure on $X$ with the same restrictions to $\del X$ as $\t$, and let $\zeta': \{M_{\alpha, \beta} \otimes M_{\beta, \gamma}\}^{\t'} \longrightarrow M_{\alpha, \gamma}$ be a module homomorphism. Set $\t' - \t = u \in H^2(X, \del X; \zz)$. Then, for every pair of homogeneous elements $\x \in {\CFm}(Y_{\alpha, \beta}, \t|_{Y_{\alpha, \beta}}; M_{\alpha, \beta})$ and $\y \in {\CFm}(Y_{\beta, \gamma}, \t|_{Y_{\beta, \gamma}}; M_{\beta, \gamma})$, we have
\begin {equation}
\label {eq:grshift}
 \gr {f^{-}_{\alpha, \beta, \gamma; \t', M, \zeta'}}(\x \otimes \y)  - \gr {f^{-}_{\alpha, \beta, \gamma; \t, M, \zeta}}(\x \otimes \y)  \equiv \bigl( c_1(\t) \smallsmile  u + u \smallsmile u \bigr) [X] \  \  (\mod \delt(\t, M)). 
\end {equation}
\end {lemma}

\begin {proof}
(a) The ambiguity in the grading shift comes from doubly periodic domains with trivial module actions. The contribution of these domains to the Maslov index is given by their pairing with $c_1(\t)$.

(b) This follows from the formula for the Maslov index of a triply periodic domain in \cite[Section 5]{SarkarIndex}.  
\end {proof}

Of course, Equation~\eqref{eq:twistedgrading} and Lemma~\ref{lem:twistedmaps} apply equally way to truncated Floer complexes $\CFmd$ instead of $\CFm$, and to triple Heegaard diagrams with several basepoints.

\begin {remark} A particular example of cobordism map is the untwisted one corresponding to a two-handle addition, as in \cite{HolDiskFour}. In this case $Y_{\beta, \gamma}$ is a connected sum of $S^1 \times S^2$'s, and we consider only triangles with one vertex at a representative for the top-degree homology generator $\y = \Theta^\can_{\beta, \gamma}$. Supposing further that $c_1(\t)$ has torsion restrictions to $Y_{\alpha, \beta}$ and $Y_{\alpha, \gamma}$, the respective Floer homology groups have absolute $\qq$-gradings compatible with their relative $\zz$-gradings. Let $W$ be the cobordism from $Y_{\alpha, \beta}$ and $Y_{\alpha, \gamma}$ obtained from $X$ by filling in the other boundary component with three-handles. It is shown in \cite{HolDiskFour} that the cobordism map shifts absolute grading by
\begin {equation}
\label {eq:shift}
\frac{c_1(\t)^2 - 2\chi(W) - 3\sigma(W)}{4},
\end {equation}
where $\chi$ and $\sigma$ denote Euler characteristic and signature, respectively. If we further suppose that $c_1(\t')$ has torsion restrictions to $Y_{\alpha, \beta}$ and $Y_{\alpha, \gamma}$, we can then view Equation \eqref{eq:grshift} as a simple consequence of the formula \eqref{eq:shift}, with $\delt(\t, M) = 0$.
\end {remark}

\subsection{A long exact sequence}
\label {sec:LES}
In this section we sketch the construction of the surgery long exact
sequence from \cite[Section 3]{IntSurg} and \cite[Section
6.1]{RatSurg}. It is stated there for rationally null-homologous
knots, but the construction can be generalized to works for arbitrary knots inside
three-manifolds. In the more general setting, one of  the Floer
complexes may appear with genuinely twisted coefficients; see
Proposition~\ref{prop:leseq} below, as well as its graded refinements: Propositions~\ref{prop:gr+}, ~\ref{prop:gr-}, \ref{prop:gr0}, \ref{prop:gr00}, \ref{prop:gr002}, and \ref{prop:gr00twisted}. (A slightly
different generalized surgery sequence with twisted coefficients was
proved by Fink \cite{Fink}.)

We work in the setting of Section~\ref{sec:4c}, with $(\orL , \Lambda)$ being a framed $\ell$-component link inside an integral homology three-sphere $Y$. Let $\bar \Lambda$ be the framing on $L$ obtained from $\Lambda$ by adding $m_1 > 0$ to the surgery coefficient $\lambda_1$ of the first component $L_1$; that is, the corresponding framing vectors are $\bar \Lambda_1 = \Lambda_1 + m_1 \tau_1$,
where here  $\tau_1$ is a meridian for $L_1$,  and $\bar \Lambda_i = \Lambda_i$ for $i \neq 1$. Also, we let $L' = L - L_1$ and denote by $\Lambda'$ the restriction of the framing $\Lambda$ to $L'$, with framing vectors $\Lambda'_2, \dots, \Lambda'_\ell$. We can view $K=L_1$ as a knot inside the three-manifold $Y_{\Lambda'}(L')$, and $Y_{\Lambda}(L)$ as the result of surgery along that knot. (Note that, given any knot $K$ inside a three-manifold $M$, we can find $L = L_1 \cup L' \subset Y = S^3$ such that $M=S^3_{\Lambda'}(L')$ and $K$ corresponds to $L_1$.)

Suppose we have a multi-pointed Heegaard diagram $(\Sigma, \alphas,
\betas, \ws, z_1)$ for $K \subset Y_{\Lambda'}(L')$, with possibly
several free basepoints, but a single basepoint pair $(w_1, z_1)$ on
$K$. Let $g$ be the genus of $\Sigma$ and $k$ the total number of $w$
basepoints, as in Section~\ref{sec:hed}. Moreover, we assume that
$w_1$ and $z_1$ can be connected by a path which crosses $\beta_{g+k-1}$ exactly
once, and which is disjoint from all the other $\alpha$ and $\beta$ curves.

As in \cite[proof of Theorem 3.1]{IntSurg}, we let $\gamma_{g+k-1}$ be a simple, closed curve in $\Sigma$ disjoint from $\beta_1, \dots, \beta_{g+k-2}$ which specifies the $\lambda_1$-framing of $L_1 \subset Y$. We complete this to a $(g+k-1)$-tuple $\gammas$ of attaching curves on $\Sigma$ by taking curves $\gamma_1, \dots, \gamma_{g+k-2}$ which approximate $\beta_1, \dots, \beta_{g+k-2}$ in the sense of Definition~\ref{def:approx}. We define another collection $\deltas$ similarly, only now $\delta_{g+k-1}$ specifies the framing $\lambda_1 + m_1$ on $L_1$. Thus, $(\Sigma, \alphas, \gammas, \ws)$ and $(\Sigma, \alphas, \deltas, \ws)$ are Heegaard diagrams for $Y_{\Lambda}(L)$ and $Y_{\bar \Lambda}(L)$, respectively.

The Heegaard triple $(\Sigma, \alphas, \gammas, \deltas, \ws)$ represents a four-manifold $X_1$ with three boundary components,
$$ Y_{\Lambda}(L), \ L(m_1, 1) \# (\#^{g+k-2} (S^1 \times S^2)), \text{ and } Y_{\bar \Lambda}(L).$$

There is a canonical torsion $\spc$ structure on the manifold $L(m_1, 1) \# (\#^{g+k-2} (S^1 \times S^2))$, see \cite[Definition 3.2]{IntSurg}. We arrange that the Floer homology $\HFm(\T_{\gamma}, \T_{\delta}, \ws)$ in that $\spc$ structure, in the maximal degree with nonzero homology, is represented by a unique intersection point, which we call {\em canonical}. We define a map 
$$f_1^-: \CFm(\Ta, \Tg, \ws) \longrightarrow \CFm(\Ta, \T_{\delta}, \ws)$$ 
by counting holomorphic triangles with one vertex in the canonical intersection point. 

Let us add $g+k-2$ three-handles to $X_1$ to kill off the $S^1 \times S^2$ summands in the middle boundary component, then remove a neighborhood of a path between the first two boundary components. We thus obtain a cobordism $W_1$ from $Y_{\Lambda}(L) \# L(m_1, 1)$ to $Y_{\bar \Lambda}(L)$. It is easy to see that
$$ f_1^{-}(\x) = f_{W_1}^- (\x \otimes \Theta^{\can}_{\gamma, \delta}),$$
where $f_{W_1}^-$ is the map on Floer complexes induced by the cobordism $W_1$ (as in \cite{HolDiskFour}), and $\Theta^\can_{\gamma,\delta}$ is the top degree generator for the Floer homology of $L(m_1, 1)$ (in its canonical $\spc$ structure $\ux^\can$).

Next, we look at the Heegaard triple $(\Sigma, \alphas, \deltas, \betas, \ws)$. This represents a cobordism $X_2$ with three boundary components,
$$ Y_{\bar \Lambda}(L), \#^{g+k-2} (S^1 \times S^2), \text{ and } Y_{ \Lambda'}(L'). $$

By filling in the middle component with three-handles, we obtain a cobordism $W_2$ from $ Y_{\bar \Lambda}(L)$ to $ Y_{ \Lambda'}(L')$. This is simply the reverse of the cobordism $-W_{\bar \Lambda}(L', L)$ in the notation of Section~\ref{sec:4c}. Here, by reverse of a cobordism we mean a reversal in direction, i.e., turning the cobordism around so that we view it as a cobordism from the final to the initial manifold. The minus sign denotes the additional reversal of orientation.

To $X_2$ we associate a cobordism map with twisted coefficients as in Section~\ref{sec:grtwisted}. (Compared to the notation there, we now have $X_2$ instead of $X$, $\deltas$ instead of $\betas$, and $\betas$ instead of $\gammas$.) We let $M_{\alpha, \delta}$ and $M_{\delta, \beta}$ be the modules $\ff$, with trivial action by the respective cohomology groups. Consider the ring 
$$\TR_1 = \ff[\zz/m_1\zz] = \ff[T_1]/(T_1^{m_1} - 1).$$  
We make $\TR_1$ into a $\ff[H^1(Y_{ \Lambda'}(L'))]$-module $M_{\alpha, \beta}$ by letting $h \in H^1(Y_{ \Lambda'}(L'))$ 
act by multiplication by $T_1^{\langle h, [L_1] \rangle}$. In the corresponding Floer complex
$ \CFm(Y_{\Lambda'}(L'); \TR_1)$
each isolated holomorphic strip gets counted in the differential with a coefficient 
$T_1^{n_{w_1}(\phi)-n_{z_1}(\phi)}$, 
where $\phi$ is the homotopy class of the strip.

Given $\t \in \spc(W_2)$, we  abuse notation slightly, letting $\t$ also denote the restriction 
of $\t$ from $W_2$ to $X_2$. We consider the module homomorphism
$$ \zeta: \{M_{\alpha, \delta} \otimes M_{\delta, \beta}\}^\t \longrightarrow M_{\alpha, \beta}$$
given by composing the projection
$$   \{M_{\alpha, \delta} \otimes M_{\delta, \beta}\}^\t \longrightarrow \frac{\underline{\spc}(X_2; \t)}{H^1( Y_{\bar \Lambda}(L)) \times H^1(\#^{g+k-2}(S^1 \times S^2))} $$ 
 with the map
 \begin {equation}
 \label {eq:hor}
  \frac{\underline{\spc}(X_2; \t)}{H^1( Y_{\bar \Lambda}(L)) \times H^1(\#^{g+k-2}(S^1 \times S^2))}  \longrightarrow \TR_1, \ \ [\psi] \to T_1^{n_{w_1}(\psi) - n_{z_1}(\psi)},
  \end {equation}
where $\psi \in {\underline{\spc}(X_2; \t)}$ is viewed as a homotopy class of triangles, compare \cite[Section 8.2.1]{HolDiskTwo}. Note that the map \eqref{eq:hor} is well-defined, because every doubly-periodic domain $\phi$ for $Y_{\bar \Lambda}(L)$ or $\#^{g+k-2}(S^1 \times S^2)$ has $n_{w_1}(\phi) = n_{z_1}(\phi)$.

We define
$$ f_2^- : \CFm(\Ta, \Td, \ws) \longrightarrow \CFm(\Ta, \Tb, \ws; \TR_1)$$
by
$$ f_2^- (\x) = \sum_{\t \in \spc(W_2)} f^-_{\alpha, \delta, \beta; \t, M, \zeta} (\x \otimes \Theta^\can_{\delta, \beta}),$$
where $\Theta^\can_{\delta, \beta}$ is the respective canonical generator. In other words, the map $f_2^-$ counts holomorphic triangles in a class $\psi$ with a coefficient $T_1^{n_{w_1}(\psi) - n_{z_1}(\psi)}$. For future reference, since we will be interested in the grading properties of $f_2^-$, we note that
\begin {equation}
\label {eq:dtm}
\delt(\t, \TR_1) = \gcd\bigl( \delt(\t|_{Y_{\bar \Lambda}(L)}), \delt(\t|_{Y_{ \Lambda'}(L')}, \TR_1)  \bigr) , 
\end {equation}
where we denoted $\delt(\t, \TR_1) = \delt(\t, M)$ for simplicity. Equation~\eqref{eq:dtm} follows from \eqref{eq:deltam} together with the triviality of the modules $M_{\alpha, \delta}$ and $M_{\delta, \beta}$.

Finally, we consider the cobordism $X_3$ corresponding to the Heegaard triple $(\Sigma, \alphas, \betas, \gammas)$. The associated filled-in cobordism $W_3$ from $Y_{ \Lambda'}(L')$ to $Y_{\Lambda}(L)$ is simply $W_{\Lambda}(L', L), $ i.e., surgery on the framed knot $(L_1, \Lambda_1)$.
We set up a map
$$ f_3^- : \CFm(\Ta, \Tb, \ws; \TR_1) \longrightarrow \CFm(\Ta, \Tg, \ws)$$
as follows. (Compare the definition of the map $f_3^+$ in \cite[Section 3]{IntSurg}.) Fix any triangle class $\psi \in \pi_2(\Theta^{\can}_{\gamma, \beta}, \Theta^{\can}_{\beta, \delta}, \Theta^{\can}_{\gamma, \delta})$, where $\Theta^{\can}_{i,j}$ are canonical generators as before. Set 
$$c=n_{w_1}(\psi) - n_{z_1}(\psi).$$
By analyzing triply periodic domains, one can check that the residue of $c$ modulo $m_1$ is independent of the choice of $\psi$. Set
\begin{equation}
\label{eq:f3}
 f_3^-(T_1^{s_1} \cdot \x) = \sum_{\y \in \Ta \cap \Tg} \sum_{\bigl\{\phi \in \pi_2(\x, \Theta^{\can}_{\beta, \gamma}, \y) \big \vert \substack{\mu(\phi)=0 \\ s_1 + n_{w_1}(\phi)-n_{z_1}(\phi) \equiv c (\text{mod  }m_1)} \bigr \}}(\# \M(\phi)) \cdot \U^{n_w(\phi)} \y,
 \end{equation}
where 
$$ \U^{n_{w}(\phi)} := U_1^{n_{w_1}(\phi)} \cdot \dots \cdot U_k^{n_{w_k}(\phi)}.$$

Observe that $f_3^-$ is not a usual cobordism map with twisted coefficients. Rather, its value depends on the exponent $s_1$ in the input $T_1^{s_1}\cdot \x$. When applied to a fixed $T_1^{s_1}\cdot \x$, the map $f_3^-$ acts as a sum of several cobordism maps with twisted coefficients; these are associated to the cobordism $W_3=W_{\Lambda}(L', L)$ and to the $\spc$ structures $\s \in \spc(W_3) \cong \H(L)/H(L, \Lambda|_{L'})$ with some fixed value of $\s$ modulo $(m_1, 0, \dots, 0)$.

We denote by $F_1^-, F_2^-, F_3^-$ the maps induced by $f_1^-, f_2^-, f_3^-$ on homology.

\begin {proposition}
\label {prop:leseq}
For any framed link $(\orL, \Lambda)$ inside an integral homology sphere, there is a long exact sequence
$$ \dots \longrightarrow \HFm(Y_{ \Lambda}(L)) \xrightarrow{F_1^-} \HFm(Y_{\bar \Lambda}(L) ) \xrightarrow{F_2^-} \HFm(Y_{\Lambda'}(L'); \TR_1 ) \xrightarrow{F_3^-} \cdots$$
In fact, the complex $ \CFm(Y_{ \Lambda}(L)) \cong \CFm(\Ta, \Tg, \ws)$ is quasi-isomorphic to the mapping cone of $f_2^-$.
\end {proposition}

The proof of Proposition~\ref{prop:leseq} is given in \cite[Section 6.1]{RatSurg} for the case when $L_1$ is rationally null-homologous inside $Y_{\Lambda'}(L')$ (and for $\iHF^+$ instead of $\HFm$), but it applies equally well to our situation, so we omit it. 

We would like to have a refined statement of Proposition~\ref{prop:leseq}, in which we keep track of the decomposition of $ \CFm(Y_{ \Lambda}(L)) $ into $\spc$ structures, as well as the respective relative gradings on it. However, keeping track of
gradings is possible only if we work with vertically truncated complexes $\CFmd$ (as in Section~\ref{sec:algebra}) instead of $\CFm$. We
write $f_i^\delta$ (resp. $H_i^{\delta}$) for the vertical truncation of $f_i^-$ 
(resp. $H_i^{\delta}$).
The proof of Proposition~\ref{prop:leseq} in fact gives the following more precise version:

\begin {proposition}
\label {prop:leseqd}
Fix $\delta \geq  0$.  For any framed link $(\orL, \Lambda)$ inside an integral homology sphere, the complex $\CFmd(\Ta, \Tg, \ws) = \CFmd(Y_{\Lambda}(L))$ is quasi-isomorphic to the mapping cone
\begin {equation}
\label {eq:mapcone}
\CFmd(Y_{\bar \Lambda}(L)) \xrightarrow{f_2^\delta} \CFmd(Y_{\Lambda'}(L'); \TR_1).
\end {equation}
\end {proposition}

In fact, following \cite{IntSurg} and \cite{RatSurg}, there are two natural quasi-isomorphisms that can be used to prove Proposition~\ref{prop:leseq}. First, we have the quasi-isomorphism
\begin {equation}
\label {eq:qi1} 
\CFmd(Y_{\Lambda}(L)) \xrightarrow{\sim}  Cone \bigl( \CFmd(Y_{\bar \Lambda}(L)) \xrightarrow{f_2^\delta} \CFmd(Y_{\Lambda'}(L'); \TR_1) \bigr) 
\end {equation}
given by a triangle-counting map $f_1^\delta$ to the first term in the mapping cone and a quadrilateral-counting map $H_1^\delta$ (a null-homotopy of $f_2^\delta \circ f_1^\delta$) to the second factor. Second, we have a quasi-isomorphism in the opposite direction
\begin {equation}
\label {eq:qi2}
Cone \bigl( \CFmd(Y_{\bar \Lambda}(L)) \xrightarrow{f_2^\delta} \CFmd(Y_{\Lambda'}(L'); \TR_1) \bigr) \xrightarrow{\sim}  \CFmd(Y_{\Lambda}(L)),
\end {equation}
given by a triangle-counting map $f_3^\delta$ from the second term of the mapping cone, and a quadrilateral-counting map $H_2^\delta$ (a null-homotopy of $f_3^\delta \circ f_2^\delta$)  from the first term.

\subsection{Refinements of Proposition~\ref{prop:leseqd}}

We will discuss below several refinements of
Proposition~\ref{prop:leseqd}, on a case-by-case basis, depending on
$\Lambda$. By choosing $m_1$ judiciously (in particular, sufficiently
large compared to $\delta$), we describe cases where the
quasi-isomorphism~\eqref{eq:qi1} or \eqref{eq:qi2} has good $\spc$
structure decompositions and good grading-preserving properties.

More precisely, we consider the following two cases: 

\begin {itemize}
\item {\bf Case I}:
$\Lambda$ is non-degenerate;

\item {\bf Case II}: $\Lambda$ is degenerate,
and furthermore $\Lambda_1\in Span(\Lambda_2,\dots,\Lambda_\ell)$.
\end {itemize}

(We do not consider the case where $\Lambda$ is degenerate, but
$\Lambda_1\not\in Span(\Lambda_2,\dots,\Lambda_\ell)$. That case
turns out not to be needed for our present applications.)

\medskip

{\bf Case I} is further subdivided as follows.
 Let $h > 0$ be the smallest integer such that the vector
$$  \a = (a_1, \dots, a_\ell) = h \Lambda^{-1} (1, 0, \dots, 0) $$
has all integer coordinates. We choose $m_1 \gg 0$ such that the vector $m_1 \tau_1= (m_1, 0, \dots, 0)$ is in $H(L, \Lambda) = Span(\Lambda_1, \dots, \Lambda_\ell)$.\footnote{By ``span'' we will mean the span of integral vectors over $\zz$, unless we explicitly refer to the ``$\qq$-span.'' } Hence, the value $d = m_1/h$ is an integer, too. Choosing $m_1$ sufficiently large is the same as choosing $d$ sufficiently large. In the case when $a_1 \neq 0$, we impose an additional constraint on our choices of sufficiently large $m_1$: namely, we ask for $m_1$ to be a multiple of $a_1h$, i.e., for $d$ to be a multiple of $a_1$.

It is easy to check that $a_1=0$ if and only if $\Lambda'$ is degenerate, and when $\Lambda'$ is nondegenerate, the sign of $a_1$ is the same as the sign of 
the restriction of $\Lambda$ (viewed as a symmetric bilinear form) to the one-dimensional space $Span_{\qq}(\Lambda_2, \dots, \Lambda_\ell)^\perp \subset \qq^\ell$. (Here, the orthogonal complement is taken with respect to the standard inner product.) 

We thus distinguish three subcases, according to the sign of $a_1$:

\smallskip
\noindent {\bf Subcase I (a):} $a_1 > 0$. 
In this subcase,    the graded refinement of Proposition~\ref{prop:leseqd}
is given in Proposition~\ref{prop:gr+} below.
The model example to keep in mind is that of positive surgery on a knot $K \subset Y$, e.g. 
$Y=S^3$, $n$ is a positive integer, $\Lambda=(n)$, ${\bar \Lambda}=(n+m_1)$, $\Lambda'=\emptyset$,
$a_1=1$, $h=n$.

\smallskip
\noindent{\bf Subcase I (b):} 
$a_1 < 0$. 
In this subcase,    the refinement
is given in Proposition~\ref{prop:gr-} below.
The model example to keep in mind is that of negative surgery on a knot $K \subset Y$, e.g. 
$Y=S^3$, $n$ is a positive integer, $\Lambda=(-n)$, ${\bar \Lambda}=(-n+m_1)$, $\Lambda'=\emptyset$,
$a_1=-1$, $h=n$.

\smallskip
\noindent {\bf Subcase I (c):} $a_1=0$.
In this subcase, the refinement is given in Proposition~\ref{prop:gr0} below.
The model example here is fractional surgery on a knot in $S^3$, as follows.
Let $K\subset Y=S^3$ be a knot, which we can promote to a two-component link $L=L_1 \cup K$, where here
$L_1$ is a meridian for $K$.
In this case, write
$$\Lambda=\left(\begin{array}{ll}
0 & 1 \\
1 & 0 
\end{array}\right), \
{\bar\Lambda}=\left(\begin{array}{ll}
m_1 & 1 \\
1 & 0 
\end{array}\right), \
\Lambda'=(0), \ a_1=0, \ h=1.$$
In this model example, $Y_\Lambda(L)=S^3$, $Y_{\bar \Lambda}(L)=S^3_{-1/m_1}(K)$, and
$Y_{\Lambda'}(L')=S^3_0(K)$. Note now that $Y_{\Lambda'}(L')$ appears with twisted coefficients.
Thus, Proposition~\ref{prop:gr0} can be viewed as a generalized, negative-surgery version of the fractional surgery theorem~\cite[Theorem 9.14]{HolDiskTwo}.

\medskip
In {\bf Case II}, there are two refinements, given in
Propositions~\ref{prop:gr00} and~\ref{prop:gr002}, corresponding to the
two possible quasi-isomorphisms from Equation~\eqref{eq:qi1}
and~\eqref{eq:qi2} respectively. Another case, where the Floer complex
for $Y_\Lambda(L)$ has twisted coefficients, is given in
Proposition~\ref{prop:gr00twisted}.  The model example for Case II 
is that of zero-surgery on a knot in an integral homology
three-sphere: $\Lambda=(0)$, ${\bar \Lambda}=(m_1)$,
$\Lambda'=\emptyset$. In this example, Proposition~\ref{prop:gr00} should be compared
with the integer surgeries exact sequence~\cite[Theorem~9.19]{HolDiskTwo}. 

\medskip

In all these cases, to understand the $\spc$ decompositions, we find it useful to study further the cobordism $W_1$ from $Y_\Lambda(L)$ to $Y_{\bar\Lambda}(L)$, in the manner of Section~\ref{sec:4c}.

\begin{figure}
\begin{center}
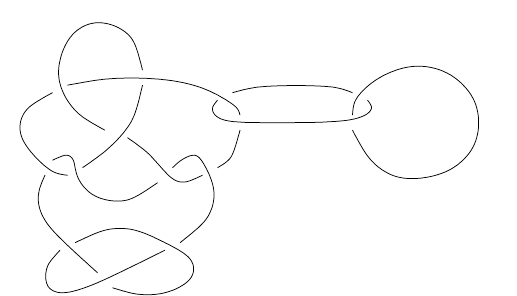
\end{center}
\caption {{\bf The cobordism $W_1$.} By adding an unlinked, unknotted component $L_{\ell+1}$ with framing $m$ to the original link $L$ we obtain a surgery presentation for $Y_{\Lambda}(L) \# L(m_1,1)$. The cobordism $W_1$ is then given adding a two-handle along the unknot component $L_{l+2}$ with framing $0$, which has linking number one with both $L_1$ and $L_{\ell+1}$. 
}
\label{fig:kirby}
\end{figure}

The cobordism $W_1$ consists of a single two-handle addition. A Kirby calculus picture for it is shown in Figure~\ref{fig:kirby}. If we denote by $L^+$ the $(\ell+2)$-component link $L \cup L_{\ell+1} \cup L_{\ell+2}$ from the figure, and by $\Lambda^+$ its given framing (also shown in the figure), we can express our cobordism as
$$ W_1 = W_{\Lambda^+}(L^+ - L_{\ell+2}, L^+),$$
in the notation used in Section~\ref{sec:4c}. In matrix form, the framing $\Lambda^+$ for $L^+$ is
\begin {equation}
\label {eq:lambda+}
 \left (
\begin{array}{cccc|cc} 
 & & & & 0 &  1\\
 & \Lambda& & & 0 & 0\\ 
 & & & & \vdots & \vdots\\ 
 & & & & 0 & 0\\
 \hline 
 0&0 & \ldots& 0 & m_1& 1\\ 
 1& 0& \ldots & 0 & 1& 0
\end {array}
\right )
\end {equation}

We denote the rows of $\Lambda^+$ by $\Lambda_i^+, \ i=1, \dots, \ell+2$.

\begin {definition}
\label {def:linked}
Two $\spc$ structures $\ux$ on $Y_{\Lambda}(L)$ and $\bar \ux$ on $Y_{\bar \Lambda}(L)$ are called {\em linked} if there exists a $\spc$ structure on $W_1$ which restricts to $\ux \# \ux^\can$  on $Y_{\Lambda}(L) \# L(m_1, 1)$, and to $\bar \ux$ on $Y_{\bar \Lambda}(L)$.
\end {definition}

Note that we can use Lemma~\ref{lem:spc4} to describe $\spc$ structures and restriction maps concretely. We have:

\begin {lemma}
\label {lem:linked}
Two $\spc$ structures $\ux \in \spc(Y_{\Lambda}(L)) \cong \H(L)/H(L, \Lambda)$ and $\bar \ux \in \spc(Y_{\bar \Lambda}(L)) \cong \H(L)/H(L, \bar \Lambda)$ are linked if and only if there exists $\s \in \H(L)$ and $j \in \zz$ such that $\ux=\s + H(L,\Lambda)$ and $\bar\ux=\s+jm_1\tau_1 + H(L,\bar\Lambda)$. 
\end {lemma} 

\begin {proof}
Suppose we have a $\spc$ structure 
$$\t \in \spc(W_1) \cong \H(L^+)/Span(\Lambda^+_1, \dots, \Lambda^+_\ell, \Lambda^+_{\ell+1}),$$
and let $\s^+ = (s_1^+, \dots, s_{\ell}^+, s_{\ell+1}^+, s_{\ell+2}^+) \in \H(L^+)$ be a representative of $\t$. Using Lemma~\ref{lem:spc4}, the condition that the restriction of $\t$ to $Y_{\Lambda}(L) \# L(m_1, 1)$ is $\ux \# \ux^\can$ translates into 
$$ \ux = \s + H(L, \Lambda),  \ \text{where} \ \s=(s_1^+ - 1/2, s_2^+, \dots, s_{\ell}^+) \in \H(L)$$
and
$$ \psi^{L-L_{\ell+1}}(\s^+) = 0 \in \spc(L(m_1,1)) \cong \zz/m_1$$
i.e.
$$ s^+_{\ell+1} = qm_1 + 1/2 $$
 for some $ q \in \zz$.

On the other hand, the restriction of $\t$ to $Y_{\bar \Lambda}(L)$ is
$$  \s^+ + H(L^+, \Lambda^+) \in \bigl(\H(L^+)/ H(L^+, \Lambda^+)\bigr).$$

Note that we have two different expressions for $\spc(Y_{\bar \Lambda}(L))$ as a quotient, one coming from the cobordism $W_1$ and one from $W_{\bar \Lambda}(L)$. They are related by the isomorphism
\begin{eqnarray*} 
 \H(L^+)/H(L^+, \Lambda^+) &\longrightarrow&  \H(L)/H(L, \bar \Lambda)  \\
\s^+ + H(L^+, \Lambda^+) &\longrightarrow& \bar \s + H(L, \bar \Lambda),
\end {eqnarray*}
where
$$
\bar  \s = (\bar s_1, \dots, \bar s_\ell) = (s_1^+-s^+_{\ell+1} + m_1 s^+_{\ell+2}, s^+_2, \dots, s^+_\ell).
$$
The conclusion follows by taking $j = s^+_{\ell+2} - q$. 
\end {proof}

\subsection {Refinements of Proposition~\ref{prop:leseqd}, Case I: $\Lambda$ is non-degenerate}
 \label {sec:case1}

\medskip{\bf Subcase I (a):} $a_1 > 0$. (The model example to keep in mind is that of positive surgery on a knot $K \subset Y$, for which the grading issues have been studied in \cite[Sections 4.5-4.6]{IntSurg}.) 

In this situation $Y_{\Lambda}(L), Y_{\bar \Lambda}(L)$, and $Y_{\Lambda'}(L')$ are all rational homology spheres, and therefore the corresponding Floer chain complexes all admit absolute $\qq$-gradings (and relative $\zz$-gradings) in each $\spc$ structure, cf. \cite{HolDiskFour}. 

In fact, we can rephrase this subcase as being about positive surgery on the knot $L_1$ inside the rational homology sphere $Y_{\Lambda'}(L')$. As such, it was fully treated in \cite[Section 6.3]{RatSurg} (see also~\cite[Sections 4.5-4.6]{IntSurg}).  However, we present a slightly different proof here, which will be easier to adapt to the other cases below, and which gives us the opportunity to establish some notation. 

We aim to understand the quasi-isomorphism $(f_1^\delta, H_1^\delta)$ from \eqref{eq:qi1}. We start by studying the map $f_1^\delta$, which corresponds to the cobordism $W_1$ from Figure~\ref{fig:kirby}.

\begin {lemma}
\label {lem:w1-}
If $a_1 > 0$, then the cobordism $W_1$ is negative definite. In fact, a generator $\Sigma_1$ of $H_2(W_1; \zz) \cong \zz$ satisfies $\Sigma_1 \cdot \Sigma_1 = -m_1(a_1d + 1)$. 
\end {lemma}

\begin {proof}
Using the matrix form \eqref{eq:lambda+} for $\Lambda^+$ and the identification \eqref{eq:homology}, we get that the generator of $H_2(W_1; \zz)$ is given in coordinates by
\begin {equation}
\label {eq:sigma1}
\Sigma_1 =(a_1d, \dots, a_\ell d, 1, -m_1).
\end {equation}
The conclusion then follows from Lemma~\ref{lem:HH} (a). 
\end {proof}

Using Lemma~\ref{lem:linked}, we can specify the relation between the $\spc$ structures on $W_1$ and its two boundary components. Precisely, since $m_1 \tau_1= (m_1, 0, \dots, 0) \in H(L, \Lambda)$, we deduce that $\spc$ structures $\ux \in \spc(Y_{\Lambda}(L)) \cong \H(L)/H(L, \Lambda)$ and $\bar \ux \in \spc(Y_{\bar \Lambda}(L)) \cong \H(L)/H(L, \bar \Lambda) $ are linked if and only if
$$ \pi(\bar \ux) =  \ux,$$
where $$\pi: \Bigl( \H(L)/H(L, \bar \Lambda) \Bigr) \to \Bigl( \H(L)/H(L,  \Lambda) \Bigr)$$
is the natural projection. Note that the projection makes sense because $H(L , \bar \Lambda) = Span (\Lambda_1 + m_1  \tau_1, \Lambda_2, \dots, \Lambda_\ell)$ is a subset of $H(L, \Lambda)$. 

This observation implies that the map $f_1^\delta$ is a direct sum of the maps
$$ f_{1, \ux}^\delta : \CFmd(Y_{\Lambda}(L), \ux) \to \bigoplus_{\{\bar \ux \in \spc(Y_{\bar \Lambda}(L)) | \pi(\bar \ux) = \ux\} }\CFmd(Y_{\bar \Lambda}(L), \bar \ux), $$
for $\ux \in \spc(Y_{\Lambda}(L))$. Further, the map $f_{1, \ux}^\delta$ is itself a sum of maps $f_{1, \ux; \t}^{\delta}$, one for each $\spc$ structure $\t \in \spc (W_1)$ with $\psi^{L_{\ell+2}}(\t) = \ux \# \ux^\can$.

Let us now turn our attention to the map 
$$f_2^\delta: \CFmd(Y_{\bar \Lambda}(L)) \to \CFmd(Y_{\Lambda'}(L'); \TR_1) \cong \bigoplus_{i=0}^{m_1-1} T_1^i \cdot \CFmd(Y_{\Lambda'}(L')),$$ whose mapping cone appears in \eqref{eq:qi1}. This is simply a twisted coefficient map associated to the cobordism $W_2$, the reverse of $-W_{\bar \Lambda}(L', L)$ in the notation of Section~\ref{sec:4c}.

\begin {lemma}
\label {lem:w2-}
The cobordism $W_2$ is negative definite. Its second homology is generated by a class $\Sigma_2$ with $\Sigma_2 \cdot \Sigma_2 = -  a_1h (a_1d + 1)$.
\end {lemma}

\begin {proof}
It is easier to think about the cobordism $W_{\bar \Lambda}(L', L)$, for which we can apply Lemma~\ref{lem:HH}. Indeed, by part (a) of that Lemma, the generator of the second homology is the vector $\a$, with $\a^t \bar \Lambda \a = a_1h(a_1d + 1)$. The change in sign in the final answer is due to the fact that in $W_2$, the orientation is reversed.
\end {proof}

Using Lemma~\ref{lem:spc4} we can relate $\spc$ structures on the two boundaries of $W_2$. By also keeping track of the powers of $T_1$, we obtain that $f_2^{\delta}$ is a  sum of maps
\begin {equation}
\label {eq:f2ut}
 f^{\delta}_{2, \bar \ux; \t} : \CFmd(Y_{\bar \Lambda}(L), \bar \ux) \to T_1^{\t} \cdot \CFmd(Y_{\Lambda'}(L'), \psi^{L_1}(\t)),
 \end {equation}
for $\t \in \spc(W_2) \cong \H(L)/H(L, \Lambda|_{L'})$ restricting to $\bar \ux \in \spc(Y_{\bar \Lambda}(L))$. Here,  $T_1^\t$ denotes $T_1^{\a \cdot (\t - \t_0)}$, where $\t_0$ is a fixed $\spc$ structure on $W_2$. Note that the expression $\a \cdot (\t - \t_0)$ is well-defined because $\a \cdot \vs = 0$ for all $\vs \in H(L, \Lambda|_{L'})$. (In fact, as can be seen from the proof of Lemma~\ref{lem:w2-}, multiplication with $\a$ represents evaluation on the homology generator.)

Taking the sum of all maps $f^{\delta}_{2, \bar \ux; \t}$ over all $\t$ (with $\bar \ux$ fixed), we obtain a map
$$  f^{\delta}_{2, \bar \ux} : \CFmd(Y_{\bar \Lambda}(L), \bar \ux) \to \bigoplus_{\{\t \in  \spc(W_2) | \pi^{L, L'}(\t) = \bar \ux)} T_1^{\t} \cdot \CFmd(Y_{\Lambda'}(L'), \psi^{L_1}(\t)).$$ 
 
 By a slight abuse of notation, for $\ux \in \spc(Y_{\Lambda}(L))$, we set 
 $$f^{\delta}_{2, \ux} = \sum_{\{\bar \ux \in \spc(Y_{\bar \Lambda}(L)) | \pi(\bar \ux) = \ux\} } f^{\delta}_{2, \bar \ux}.$$
  
Note that each term $T_1^i \cdot  \CFmd(Y_{\Lambda'}(L'), \ux')$, for $i=0, \dots, m_1-1$ and $\ux' \in \spc (Y_{\Lambda'}(L'))$, appears in the target of exactly one of the maps $f^{\delta}_{2, \ux}$. Indeed, two $\spc$ structures $\t, \t'$ on $W_2$ have the same term in the target if and only if $\t - \t' \in \zz^\ell /H(L, \Lambda|_{L'})$ has a representative of the form $(dh/a_1, 0, \dots, 0)$. Since we have chosen $d$ to be a multiple of $a_1$, and $(h, 0, \dots, 0)$ is in the span of $\Lambda_1, \dots, \Lambda_\ell$, we get that $\t$ and $\t'$ have the same reduction modulo $H(L, \Lambda)$, i.e. they correspond to the same $\ux$. We conclude that the map $f_2^{\delta}$ is the direct sum of all $f^{\delta}_{2, \ux}$, for $ \ux \in  \spc(Y_{ \Lambda}(L))$. 
    
The last map that appears in \eqref{eq:qi1} is the null-homotopy $H_1^{\delta}$. This also splits as a direct sum of maps $H_{1, \ux}^{\delta}$ over $\ux \in \spc(Y_{\Lambda}(L))$, where $H_{1, \ux}^{\delta}$ denotes the corresponding null-homotopy of $f^{\delta}_{2, \ux} \circ f^{\delta}_{1, \ux}$. Moreover, each $H_{1, \ux}^{\delta}$ is a sum of maps $H_{1, \ux; \t}^{\delta}$, over $\spc$ structures on $W_1 \cup W_2$ restricting to $\ux \# \ux^\can$ on the boundary component $Y_{\Lambda}(L) \# L(m_1, 1)$.
 
The quasi-isomorphism \eqref{eq:qi1} can then be viewed as a direct sum itself. Precisely, the summand 
 $\CFmd(Y_{\Lambda}(L), \ux)$ is quasi-isomorphic to the mapping cone of
 \begin {equation}
 \label {eq:qi1u}
 \bigoplus_{\{\bar \ux \in \spc(Y_{\bar \Lambda}(L)) | \pi(\bar \ux) = \ux\} }\CFmd(Y_{\bar \Lambda}(L), \bar \ux) \xrightarrow{f_{2, \ux}^\delta} \bigoplus_{\{\t \in  \spc(W_2) | \pi(\pi^{L, L'}(\t)) = \ux\}} T_1^\t \cdot \CFmd(Y_{\Lambda'}(L'), \psi^{L_1}(\t)) 
  \end {equation}
via the map $ (f_{1, \ux}^\delta, H_{1, \ux}^{\delta})$.

Note that the direct summands appearing in \eqref{eq:qi1u} are
absolutely $\qq$-graded and relatively $\zz$-graded (in a compatible
way). 

We seek to prove:
\begin {proposition}
\label {prop:gr+}
Fix $\delta > 0$. 
Suppose that $\Lambda$ is nondegenerate, $a_1 > 0$, and $m_1$ (a multiple of $a_1h$) is sufficiently large.
Then for every $\ux \in \spc(Y_{\Lambda}(L))$,  there is a relative $\Z$-grading on $Cone (f_{2, \ux}^\delta)$
such that the quasi-isomorphism
$$  (f_{1, \ux}^\delta, H_{1, \ux}^{\delta}): \CFmd(Y_{\Lambda}(L), \ux) \xrightarrow{\sim}  Cone (f_{2, \ux}^\delta) $$
respects the relative $\zz$-gradings on the two sides.
\end {proposition}

For this, we need a few lemmas, for which we make use of the absolute $\qq$-gradings on our complexes:

\begin {lemma}
\label {lem:24}
Let $\orL \subset Y$ and $\Lambda$ be as above, and fix $\delta > 0$. Then, there are constants $C_1$ and $C_2$ (depending only on $\orL, \Lambda$ and $\delta$) such that for all sufficiently large $m_1$ (chosen as specified above, i.e. a multiple of $a_1h$) and for all $\bar \ux \in \spc(Y_{\bar \Lambda}(L)) \cong \H(L)/H(L, \bar \Lambda)$, we have
\begin {equation}
\label {eq:101}
 \max \gr  \CFmd(Y_{\bar \Lambda}(L), \bar \ux) -  \min \gr  \CFmd(Y_{\bar \Lambda}(L), \bar \ux) \leq C_1,
 \end {equation}
\begin {equation}
\label {eq:102}
\left | \max \gr  \CFmd(Y_{\bar \Lambda}(L), \bar \ux) - \min_{\{\s \in \H(L)| [\s]=\bar \ux\}}  \frac{|(2\s-\bar \Lambda_1) \cdot \a|^2}{4a_1 h(a_1d+1)}   \right |   \leq C_2.
 \end {equation}
\end {lemma}

\begin {proof}
This is an analogue of Corollary 2.4 from \cite{IntSurg}, and has a similar proof, so we only sketch the argument. (See also \cite[Lemma 4.6 and the proof of Proposition 4.2]{RatSurg}.) The first inequality follows from 
the large surgeries theorem (\cite[Theorem 4.1]{RatSurg}), which gives a relatively graded identification of $  \CFmd(Y_{\bar \Lambda}(L), \bar \ux)$ with a generalized Floer complex of $L_1$ inside $Y_{\Lambda'}(L')$. This latter complex is independent of $m_1$. 

For the second inequality,
note that if $\t \in \spc(W_2) \cong \H(L)/H(L, \Lambda|_{L'})$, then, using Lemma~\ref{lem:w2-} and the formula $\Sigma_2 = \a$ for the homology generator in standard coordinates, we get 
$$ \langle c_1(\t), [\Sigma_2] \rangle = (2\s - \bar \Lambda_1) \cdot \a,$$
and
\begin {equation}
\label {eq:c1t}
 c_1(\t)^2 = - \frac{|(2\s-\bar \Lambda_1) \cdot \a|^2}{a_1 h(a_1d+1)} ,
 \end {equation}
where $\s \in \H(L)$ is any representative of $\t$. Consider the map induced on (truncated) Floer homologies by the cobordism $W_2$ in the $\spc$ structure $\t$ for which $c_1(\t)^2$ is maximized along all $\t$ that restrict to $\bar \ux$.
The large surgeries theorem~\cite[Theorem 4.1]{RatSurg} also
identifies this map with a standard inclusion map between generalized Floer complexes. Since the target graded group $  \CFmd(Y_{ \Lambda'}(L'), \psi^{L_1}(\t))$ is independent of $m_1$, inequality \eqref{eq:102} follows from the formula \eqref{eq:shift} for absolute grading shifts.
\end {proof}

\begin {lemma}
\label {lem:w22}
Fix $\orL, \Lambda, \delta$ as above, and a constant $C_3 \in \rr$. Then, there is a constant $b$ with the following property. For all sufficiently large $m_1$ (divisible by $a_1h$), for any fixed $\bar \ux \in \spc(Y_{\bar \Lambda}(L)) \cong \H(L)/H(L, \bar \Lambda)$, there are at most two $\spc$ structures $\t$ on $W_2$ whose restriction to $ Y_{\bar \Lambda}(L)$ is $\bar \ux$, and with the property that
\begin {equation}
\label {eq:103}
\max \gr  \ \CFmd(Y_{\bar \Lambda}(L), \bar \ux) \geq C_3 - \frac{c_1(\t)^2}{4};
\end {equation}
these are the $\spc$ structures $\t= \t^{\pm}_{\bar \ux}$ with
$$  \langle c_1(\t^{\pm}_{\bar \ux}), [\Sigma_2] \rangle = (2\s \pm \bar \Lambda_1) \cdot \a,$$
where $\s \in \H(L)$ is any representative of $\bar \ux$ which satisfies
\begin {equation}
\label {eq:104}
 - h(a_1d+1)/2 \leq \s \cdot \a < h(a_1d + 1)/2.
 \end {equation}
All other $\spc $ structures $\t$ restricting to $\bar \ux$ satisfy the inequality
$$ c_1(\t)^2 \leq -4m_1.$$
Moreover, if there is no representative $\s$ of $\bar \ux$ satisfying $|\s \cdot \a| < b$, then there is a unique $\spc$ structure $\t$ (restricting to $\bar \ux$) that satisfies \eqref{eq:103}, namely the one for which $|\langle c_1(\t), [\Sigma_2] \rangle |$ is minimal.
\end {lemma}

\begin {proof}
Since $\bar \Lambda_1 \cdot \a = h(a_1d+1)$ and $\Lambda_i \cdot \a = 0$ for $i > 1$, each $\bar \ux \in \H(L)/H(L, \bar \Lambda)$ has indeed 
a representative $\s$ satisfying \eqref{eq:104}, and the value of $\s\cdot \a$ is independent of that representative.

We use formula~\eqref{eq:c1t} for $c_1(\t)^2$, as well as the inequalities \eqref{eq:101} and \eqref{eq:102} to verify the statements of the lemma. For more details in a special case, see \cite[Lemma 4.4]{IntSurg}.
\end {proof}

\begin {lemma}
\label {lem:w11}
Fix a constant $C_0$. For all sufficiently large $m_1$ (divisible by $a_1h$),  the following statement holds. Each $\spc$ structure $\bar \ux$ over $Y_{\bar \Lambda}(L)$ has at most one extension $\t$ to $W_1$ whose restriction to $Y_{\Lambda}(L) \# L(m_1,1)$ is $\pi(\bar \ux) \# \ux^\can$ and for which
\begin {equation}
\label{eq:c00}
 C_0 \leq c_1(\t)^2 + m_1.
 \end {equation}
Further, if such a $\spc$ structure $\t = \t^\circ_{\bar \ux}$ exists, then
\begin {equation}
\label{eq:c01}
 c_1(\t^\circ_{\bar \ux})^2 = - \frac{4d (\s \cdot \a)^2}{h(a_1d +1)},
 \end {equation}
where $\s \in \H(L)$ is the representative of $[\s] = \bar \ux \in \H(L)/H(L, \bar \Lambda)$ for which the absolute value $|\s \cdot \a|$ is minimal. 
\end {lemma}

\begin {proof} (Compare \cite[Lemma 4.7]{IntSurg} and \cite[Lemma 6.7]{RatSurg}.) Since both boundaries of $W_1$ are rational homology spheres, for any $\t \in \spc(W_1)$ we can write $c_1(\t) = \alpha \cdot \text{PD}(\Sigma_1)$, for some $\alpha \in \qq$. Here, $\Sigma_1$ is the homology generator from Lemma~\ref{lem:w1-}. Using the computation  $\Sigma_1 \cdot \Sigma_1=-m_1(a_1d+1)$ from that lemma,  inequality~\eqref{eq:c00} becomes
\begin {equation}
\label {eq:alfa}
 |\alpha| \leq \sqrt{\frac{m_1 - C_0}{m_1(a_1d+1)}}.
\end {equation}
(Recall that $d=m_1/h$.)
Restriction to the boundaries determines $\t$ up to addition of PD$(\Sigma_1)$, i.e. $\bar \ux$ determines $\alpha$ up to the addition of an even integer. Since the right
hand side of \eqref{eq:alfa} becomes very small when $d$ gets large, in each equivalence class mod $2\zz$ there is at most one $\alpha$ satisfying \eqref{eq:alfa}. 
Inequality~\eqref{eq:c00} follows.

To establish \eqref{eq:c01}, we investigate in more detail the possible values of $ \langle c_1(\t), [\Sigma_1] \rangle$, over all $\t$ which have fixed restrictions $\bar \ux$ and $\pi(\bar \ux) \# \ux^\can$ to the boundaries. Given such a $\spc$ structure 
$$\t \in \spc(W_1) \cong \H(L^+)/Span(\Lambda^+_1, \dots, \Lambda^+_\ell, \Lambda^+_{\ell+1}),$$
let $\s^+ = (s_1^+, \dots, s_{\ell}^+, s_{\ell+1}^+, s_{\ell+2}^+) \in \H(L^+)$ be a representative of $\t$. Recall from the proof of Lemma~\ref{lem:linked}, that 
\begin {equation}
\label {eq:106}
 s^+_{\ell+1} = m_1 + \frac{1}{2},
 \end {equation}
 and the restriction of $\t$ to $Y_{\bar \Lambda}(L)$ is $\bar \ux = [\s] \in \H(L)/ H(L, \bar \Lambda)$, 
where
\begin {equation}
\label {eq:105}
 \s = (s_1, \dots, s_\ell) = (s_1^+-s^+_{\ell+1} + m_1 s^+_{\ell+2}, s^+_2, \dots, s^+_\ell).
 \end {equation}

Using \eqref{eq:106}, \eqref{eq:105}, and the formula \eqref{eq:sigma1} for the homology generator $\Sigma_1$, we get
\begin {eqnarray*}
 \langle c_1(\t), \Sigma_1 \rangle &=& (2\s^+ - \Lambda_{\ell+2}^+) \cdot (a_1d, \dots, a_\ell d, 1, -m_1)\\
&=& 2d\s \cdot \a + 2m_1(a_1d+1)(1-s^+_{\ell+2})\\
&=& 2d(\s + \bar \Lambda_1(1-s^+_{\ell+2})) \cdot \a.  
\end {eqnarray*}

Since the equivalence class $\bar \ux = [\s] \in \H(L)/H(L, \bar \Lambda)$ is unchanged by the addition of a multiple of $\bar \Lambda_1$, by re-labelling $\t$ we conclude that the possible values of $\langle c_1(\t), \Sigma_1 \rangle$ (when $\bar \ux$ is fixed) are exactly given by
$$ \langle c_1(\t), \Sigma_1 \rangle= 2d \s \cdot \a,$$
where $\s \in \H(L)$ is a representative of $\t$. Hence,
$$\alpha = \frac{2\s \cdot \a }{h(a_1d+1)} \in \qq,$$
and
$$ c_1(\t)^2 = - \frac{4d (\s \cdot \a)^2}{h(a_1d +1)}.$$
Thus, $|\alpha|$ is small if and only if $|\s \cdot \a|$ is small, and \eqref{eq:c01} follows. 
\end {proof}

\begin {proof}[Proof of Proposition~\ref{prop:gr+}] 
We start by equipping the domain of the 
map $f_{2,\ux}^\delta$ with a relative $\Z$-grading such that $f_{1,\ux}^\delta$ respects the gradings.
We do this as follows.
Lemma~\ref{lem:24} supplies a constant $C_0$ with the property that for all $\t \in\spc(W_1)$ with
 $\psi^{L_{\ell+2}}(\t)=\ux\#\ux^\can$, the map $f_{1,\ux;\t}^\delta$  is zero unless Inequality~\eqref{eq:c00}
is satisfied. Now, Lemma~\ref{lem:w11} shows that for each
$\bar\ux\in \spc(Y_{\bar\Lambda}(L))$ such that $\pi(\bar\ux)=\ux$,
there is at most one extension $\t=\t_{\bar\ux}^\circ$ of $\bar\ux$ to $W_1$ 
satisfying ~\eqref{eq:c00}. Now we choose a relative grading on the domain of $f_{2,\ux}^\delta$ so that
each map $f_{1,\ux;\t_{\bar\ux}^\circ}$ is grading preserving.

We grade the range of $f_{2,\ux}^\delta$ so that this map (thought of as differential of the mapping cone) drops grading by one.
To do this, we need to check that, for any given term
$$ T_1^i \cdot \CFmd (Y_{\Lambda'}(L'), \ux') $$ in the target of
$f^{\delta}_{2, \ux}$, the compositions $f^{\delta}_{2, \ux; \t_2}
\circ f^{\delta}_{1, \ux; \t_1}$ which hit that summand induce the
same shift in grading. Lemmas~\ref{lem:w22} and \ref{lem:w11} imply
that, if we fix the target, there are at most two pairs $(\t_1, \t_2)$
for which the corresponding compositions are nonzero. In fact, in many
cases there is at most one such pair. If two pairs with nonzero
compositions (with the same domain and target) exist, they are of the
form
$$ (\t_1, \t_2) = (\t_{\bar \ux}^\circ, \t_{\bar \ux}^+) \ \text{ and }  \ (\t_1', \t_2')=(\t_{\bar \ux+ \Lambda_1}^\circ, \t_{\bar \ux+  \Lambda_1}^-),$$
where $\bar \ux \in \spc (Y_{\bar \Lambda}(L))$ admits a representative $\s$ with $|\s \cdot \a| < b$. 

 Let $\s$ be the representative of $\bar \ux$ with $|\s \cdot \a|$ minimal. Using the formula \eqref{eq:shift}, checking that $f^{\delta}_{2, \ux; \t_2} \circ f^{\delta}_{1, \ux; \t_1}$ and $f^{\delta}_{2, \ux; \t_2'} \circ f^{\delta}_{1, \ux; \t_1'}$ shift degree by the same amount is equivalent to showing that
\begin {equation}
\label {eq:cccc}
 c_1(\t_1)^2 + c_1(\t_1')^2 = c_1(\t_2)^2 + c_1(\t_2)^2.
\end{equation}

Indeed, using Lemma~\ref{lem:w11} we get
\begin {eqnarray*}
c_1(\t_1)^2 - c_1(\t_1')^2 &=& \frac{4d| (\s+ \Lambda_1) \cdot \a|^2 }{h(a_1d+1)} - \frac{4d|\s \cdot \a|^2 }{h(a_1d+1)} \\
&=& \frac{4d(\Lambda_1 \cdot \a)((2\s + \Lambda_1) \cdot \a)}{h(a_1d+1)} \\
&=& \frac{4d(2\s + \Lambda_1) \cdot \a}{a_1d+1}.
\end {eqnarray*}

Using Lemma~\ref{lem:w22}, we get
\begin {eqnarray*}
c_1(\t_2')^2 - c_1(\t_2)^2 &=& \frac{(2\s + \bar \Lambda_1) \cdot \a}{a_1h(a_1d+1)} - \frac{(2\s + 2\Lambda_1 - \bar \Lambda_1) \cdot \a}{a_1h(a_1d+1)} \\
&=& \frac{4\bigl( (\bar \Lambda_1 - \Lambda_1) \cdot \a\bigr ) \cdot \bigl( (2\s + \Lambda_1) \cdot \a \bigr)}{a_1h(a_1d+1)} \\
&=& \frac{4d(2\s + \Lambda_1) \cdot \a}{a_1d+1}.
\end {eqnarray*}

Equation~\eqref{eq:cccc} is therefore satisfied, and we conclude that  $f^\delta_{1, \ux}$ preserves the relative grading. The other component of the quasi-isomorphism under consideration is the null-homotopy $H^\delta_{1, \ux}$. To check that it is grading-preserving, note that it is a sum of terms $H^{\delta}_{1, \ux; \t}$ over certain $\spc$ structures on $W_1 \cup W_2$. The inequalities \eqref{eq:101}, \eqref{eq:102}, \eqref{eq:103} and \eqref{eq:c00} imply that the only nonzero terms correspond to $\spc$ structures of the form $(\t_1, \t_2) = (\t^\circ_{\bar \ux}, \t^\pm_{\bar \ux})$; compare \cite[Lemma 4.8]{IntSurg}. The grading shifts are one less than to the corresponding shifts of the compositions  $f^{\delta}_{2, \ux; \t_2} \circ f^{\delta}_{1, \ux; \t_1}$. Hence, since $f^\delta_{1, \ux}$ preserves the relative grading, so does $H^\delta_{1, \ux}$. 
\end {proof}

\medskip{\bf Subcase I (b):} $a_1 < 0$.
The model example is that of negative surgery on a knot $K \subset Y$, which was discussed in \cite[Section 4.7]{IntSurg}. 

In this situation again $Y_{\Lambda}(L), Y_{\bar \Lambda}(L)$, and $Y_{\Lambda'}(L')$ are rational homology spheres, so we can use absolute gradings and computations of grading shifts. The computations in Lemmas~\ref{lem:w1-} and \ref{lem:w2-} still hold, with the important difference that now $\Sigma_1 \cdot \Sigma_1 = -m_1(a_1d+1)$ is positive, hence $W_1$ is positive definite.
Since we need to work with maps induced by negative definite cobordisms, we consider instead the quasi-isomorphism
in the other direction (i.e. Equation~\eqref{eq:qi2} rather than~\eqref{eq:qi1}), and show that
it preserves relative grading. To this end, instead of $W_1$ we use the negative definite cobordism $W_3 = W_{\Lambda}(L', L)$ from $Y_{\Lambda'}(L')$ to $Y_{\Lambda}(L)$, which gives rise to the map $f_3^{\delta}$. 

The quasi-isomorphism \eqref{eq:qi2} is a direct sum of quasi-isomorphisms $(H^{\delta}_{2, \ux}, f^\delta_{3, \ux})$, one for each $\ux \in \spc (Y_{\Lambda}(L))$. 

\begin {proposition}
\label {prop:gr-}
Fix $\delta > 0$. If $\Lambda$ is nondegenerate, $a_1 < 0$, and $m_1$ (a multiple of $a_1h$) is sufficiently large, then for every $\ux \in \spc(Y_{\Lambda}(L))$, the quasi-isomorphism
$$  (H_{2, \ux}^\delta, f_{3, \ux}^{\delta}): Cone (f_{2, \ux}^\delta) \xrightarrow{\sim} \CFmd(Y_{\Lambda}(L), \ux) $$
respects the relative $\zz$-gradings on the two sides.
\end {proposition}

\begin {proof}[Sketch of the proof]
The quasi-isomorphism under consideration is a sum of several maps, corresponding to $\spc$ structures on $W_2 \cup W_3$. Analogues of Lemmas~\ref{lem:w22} and \ref{lem:w11} show that, if we fix the domain and the target of the maps, there at most two $\spc$ structures on $W_2 \cup W_3$ which give nonzero maps, and a computation similar to that for \eqref{eq:cccc} shows that those two maps shift the absolute gradings in the same way. Compare \cite[Section 4.7]{IntSurg} and \cite[Section 6.4]{RatSurg}.
\end {proof}

\medskip{\bf Subcase I (c): $a_1 = 0$, i.e. $\Lambda'$ is degenerate.} (A model example to keep in mind is surgery on a link of two components $L= L_1 \cup L_2 \subset Y$, where the linking number of $L_1$ and $L_2$ is nonzero, and the framing coefficient of $L_2$ is zero.)

In this situation Lemma~\ref{lem:w1-} still holds, we have $\Sigma_1 \cdot \Sigma_1 = - m_1$, so the cobordism $W_1$ is negative definite. Consequently, we choose to look at the quasi-isomorphism \eqref{eq:qi1}, just like in Subcase I (a). The main difference from that subcase is that now $Y_{\Lambda'}(L')$ has $b_1 =1$, so the respective Floer complex (with twisted coefficients) does not have a relative $\zz$-grading. 

Let us study the relationship between $\spc$ structures on the boundaries of the cobordisms $W_1$ and $W_2$, just as we did in Case I (a). For $W_1$, it is still true that $\ux  \in \spc(Y_{\Lambda}(L) )$ and $\bar \ux \in \spc(Y_{\bar \Lambda}(L))$ are linked if and only if $\pi(\bar \ux) =  \ux$, where 
$$\pi: \Bigl( \H(L)/H(L, \bar \Lambda) \Bigr) \to \Bigl( \H(L)/H(L,  \Lambda) \Bigr)$$
is the natural projection. Moreover, because now $(m_1, 0, \dots, 0) = \bar \Lambda_1 - \Lambda_1$ is in the span of $\Lambda_2, \dots, \Lambda_\ell$, we have $H(L, \Lambda) = H(L, \bar \Lambda)$, and the projection $\pi$ is actually a bijection. Thus, the triangle map $f_1^\delta$ decomposes as a direct sum of maps
$$  f_{1, \ux}^\delta : \CFmd(Y_{\Lambda}(L), \ux) \to \CFmd(Y_{\bar \Lambda}(L), \pi^{-1}(\ux) ), $$
for $\ux \in \spc(Y_{\Lambda}(L))$. Just as before, the map $f_{1, \ux}^\delta$ is itself a sum of maps $f_{1, \ux; \t}^{\delta}$, over $\spc$ structures $\t \in \spc (W_1)$ with the given restriction to the boundary.

With regard to the cobordism $W_2$, we have
\begin {lemma}
\label {lem:w2restriction}
If $\Lambda$ is nondegenerate and $a_1 = 0$, then the natural restriction map
$$ \spc(W_2) \to \spc(\del W_2) \cong \spc(Y_{\bar \Lambda}(L)) \times \spc(Y_{\Lambda'}(L'))$$
is injective. 
\end {lemma}

\begin {proof}
It suffices to prove the analogous statement for second cohomology groups, which is easier because then the restriction map is a homomorphism and we can talk about its kernel. It is also easier to think of the cobordism $W_{\bar \Lambda}(L',L)$, which is the reverse of $W_2$ with the opposite orientation, and apply Lemma~\ref{lem:HH} to it. Suppose $\s \in \zz^\ell$ is such that 
$$ [\s] \in H^2(W_{\bar \Lambda}(L',L)) \cong  \zz^\ell/Span(\Lambda_2, \dots, \Lambda_\ell)$$
has trivial projections to $\zz^\ell/ H(L, \bar \Lambda)$ and to $\zz^{\ell-1}/ H(L', \Lambda')$. We need to show that $[\s]=0$.

By assumption, there exist $b_2, \dots, b_\ell$ such that $\s - b_2\Lambda_2 - \dots- b_\ell\Lambda_\ell$ has only the first coordinate nonzero. Hence $[\s]$ has a representative $\s'=(s'_1, 0, \dots, 0)$. Also by hypothesis, $\s'$ must be in $H(L, \bar \Lambda) = H(L, \Lambda) = Span(\Lambda_1, \dots, \Lambda_\ell)$. But since $(h, 0, \dots, 0)$ is in the span of $\Lambda_2, \dots, \Lambda_\ell$, and $\Lambda$ is nondegenerate, we must have that $\s'$ is also in the span of $\Lambda_2, \dots, \Lambda_\ell$. Therefore, $[\s]=[\s']=0$.
\end {proof}

Recall that $\Lambda \cdot (0, a_2, \dots, a_\ell) = (h, 0, \dots, 0)$, so since $\Lambda$ is nondegenerate, the kernel $H(L', \Lambda')^\perp$ of $\Lambda'$ (which is identified with $H^1(Y_{\Lambda'}(L'))$, cf. Section~\ref{sec:gradings}) must be a copy of $\zz$, generated by $\a' = (a_2, \dots, a_\ell)$. If $\Lambda'_1 = \Lambda_1|_{L'} $ is the vector of linking numbers between $L_1$ and the other components, we have $\Lambda'_1 \cdot \a' = h$, i.e., the generator $\a'$ of $H^1(Y_{\Lambda'}(L'))$ evaluates to $h$ on $[L_1]$. We define an $\ff[H^1(Y_{\Lambda'}(L'))]$-module $$\V_1 = \ff[V_1]/(V_1^d - 1),$$ where $c \in H^1(Y_{\Lambda'}(L'))$ acts on $\V_1$ by multiplication by 
$ V_1^{\langle c, [L_1] \rangle /h}$.

If we equate $T_1^h = V_1$, it is easy to check that for every $\ux' \in  \spc(Y_{\Lambda'}(L'))$ we have an isomorphism
$$ \CFmd(Y_{\Lambda'}(L'), \ux'; \TR_1) \cong \bigoplus_{i=0}^{h-1}  T_1^i \cdot \CFmd(Y_{\Lambda'}(L'), \ux'; \V_1).$$
 
Lemma~\ref{lem:w2restriction} implies that $f_2^\delta$ is a direct sum of maps
$$ f^{\delta}_{2, \ux, \ux'} : \CFmd(Y_{\bar \Lambda}(L), \pi^{-1}( \ux)) \to T_1^{\ux} \cdot \CFmd(Y_{\Lambda'}(L'), \ux'; \V_1) $$
each either zero, or corresponding to the unique $\spc$ structure on $W_2$ that extends $\pi^{-1}(\ux) \in \spc(Y_{\bar \Lambda}(L))$ and $\ux' \in \spc(Y_{\Lambda'}(L'))$. Here, $\ux$ denotes a $\spc$ structure on $Y_{\Lambda}(L)$, and $T_1^{\ux}$ is shorthand for $T_1^{(\ux - \ux_0) \cdot \a}$, where $\ux_0$ is a fixed $\spc$ structure on  $Y_{\Lambda}(L)$. (Observe that the dot product with $\a$ is well-defined modulo $m_1$.)  We let $f^{\delta}_{2, \ux}$ be the sum of $ f^{\delta}_{2, \ux, \ux'}$, over all possible $\ux'$ (fixing $\ux$). 

For $\ux' = [(s_2, \dots, s_{\ell})] \in \spc(Y_{\Lambda'}(L')) \cong \H(L')/H(L', \Lambda')$, using \eqref{eq:twistedgrading} we get
$$ \delt(\ux', \TR_1) = \delt(\ux', \V_1) =  \gcd_{\{\vs \in H(L', \Lambda')^\perp| \Lambda_1' \cdot \vs \equiv 0 (\mod m_1) \}} \sum_{i=2}^\ell 2s_iv_i.$$
Note that $ \Lambda_1' \cdot \vs \equiv 0 (\mod m_1)$, or equivalently $( \Lambda_1' / h)\cdot \vs \equiv 0 (\mod d)$, is the condition on $\vs \in H(L', \Lambda)^\perp \cong H^1(Y_{\Lambda'}(L'))$ to act trivially on the module $\TR_1$ (or, equivalently, on $\V_1$).

Recall that $H(L', \Lambda')^\perp$ is a copy of $\zz$, generated by $\a' = (a_2, \dots, a_\ell)$. For $v \in \zz$, we have $\Lambda_1' \cdot (v \a) = vh$, which is divisible by $m_1 =dh$ if and only if $v$ is divisible by $d$. Hence,
$$\delt(\ux', \V_1) = d \cdot \left | \sum_{i=2}^\ell 2s_ia_i \right |.$$

Note that this is always divisible by $2d$, because the expression inside the absolute value is even (see Remark~\ref{rem:even}). Therefore, in particular, the Floer complexes $\CFmd(Y_{\Lambda'}(L'), \ux'; \V_1) $ have relative $\zz/2d\zz$-gradings. 

Further, for any $\t \in \spc(W_2)$ whose restriction to $Y_{\Lambda'}(L')$ is $\ux'$, Equation~\eqref{eq:dtm}  says that $\delt(\t, \V_1) = \delt(\t, \TR_1) = \delt(\ux', \TR_1)$. Thus, according to Lemma~\ref{lem:twistedmaps} (a), every map $f_{2, \ux, \ux'}^\delta$ preserves the relative $\zz/2d\zz$-gradings.

We equip the mapping cone of $f^{\delta}_{2, \ux}$ with the relative $\zz/2d\zz$-grading which gets decreased by one under $f^{\delta}_{2, \ux}$. Note that the target of $f^{\delta}_{2, \ux}$ is a direct sum of terms corresponding to possible $\ux'$, which a priori have unrelated relative $\zz/2d\zz$-gradings. As part of the mapping cone, however, their direct sum becomes relatively $\zz/2d\zz$-graded.

The second ingredient in \eqref{eq:qi1}, the null-homotopy $H_1^\delta$, is then a direct sum of maps $H_{1, \ux}^\delta$, which represent null-homotopies of $f^{\delta}_{2,\ux} \circ f_{1, \ux}^\delta$. 

\begin {proposition}
\label {prop:gr0}
Fix $\delta > 0$. If $\Lambda$ is nondegenerate and $a_1=0$, then for every $m_1 = dh$ sufficiently large, and for every $\ux \in \spc(Y_{\Lambda}(L))$, the quasi-isomorphism
$$  (f_{1, \ux}^\delta, H_{1, \ux}^{\delta}): \CFmd(Y_{\Lambda}(L), \ux) \xrightarrow{\sim}  Cone (f_{2, \ux}^\delta) $$
respects the relative $\zz/2d\zz$-gradings on the two sides.
\end {proposition}

Instead of Lemma~\ref{lem:w11}, we now have the following:

\begin {lemma}
\label {lem:w1deg}
Fix a constant $C_0$. For all sufficiently large $m_1=dh$,  the following statement holds. Each $\spc$ structure $\bar \ux$ over $Y_{\bar \Lambda}(L)$ has at most two extensions $\t$ to $W_1$ whose restrictions to $Y_{\Lambda}(L) \# L(m_1,1)$  are $\pi(\bar \ux) \# \ux^\can$ and for which
\begin {equation}
\label{eq:000}
 C_0 \leq c_1(\t)^2 + m_1.
 \end {equation}
Further, if two such $\spc$ structures $\t, \t'$ exist, then
\begin {equation}
\label{eq:001}
 c_1(\t)^2 = c_1(\t')^2.
 \end {equation}
\end {lemma}

\begin {proof}
Just as in the proof of Lemma~\ref{lem:w11}, we find that if $\t$ satisfies \eqref{eq:000} and $c_1(\t) = \alpha \cdot \text{PD}(\Sigma_1)$, then
\begin {equation}
\label {eq:alfalec}
 |\alpha| \leq \sqrt{\frac{m_1 - C_0}{m_1(a_1d+1)}} = \sqrt{1-\frac{C_0}{m_1}}.
 \end {equation}

Also, $ \alpha = 2 (\s \cdot \a)/h$ for some representative $\s \in \H(L)$ of $\t$, so $\alpha$ must live in the discrete subset $\frac{1}{h} \cdot \zz \subset \qq$. Hence, if we choose $m_1 > -C_0h/(2h+1)$, the only way inequality \eqref{eq:alfalec} is satisfied is if
$$ |\alpha| \leq 1.$$

If we fix the restriction of $\t$ to $\del W_1$, the value of $\alpha$ is determined up to the addition of even integers. Typically there is at most one $\t$ with $|\alpha| \leq 1$. The only time there are two
such $\spc$ structures $\t$ and $\t'$, the respective values of $\alpha$ are $\pm 1$, and we have
$$ c_1(\t)^2 = c_1(\t')^2 = -m_1,$$
as claimed. \end {proof}

\begin {proof}[Proof of Proposition~\ref{prop:gr0}] Using Lemma~\ref{lem:w1deg} and formula~\eqref{eq:shift}, we find that each $f^{\delta}_{1, \ux}$ is the sum of at most two nonzero maps corresponding to $\spc$ structures on $W_1$, and if there are two such maps, they shift absolute grading by the same amount. Thus, $f^{\delta}_{1, \ux}$ preserves the relative $\zz$-gradings on the two sides (and hence their $\zz/2d\zz$ reductions). From our discussion of $W_2$ we also know that each $f^{\delta}_{2, \ux}$ preserves relative grading. The $\spc$ structures that give nontrivial contributions to $H_{1, \ux}^{\delta}$ are subject to similar constraints, and shift grading by one degree less than the respective maps $f^{\delta}_{2,\ux} \circ f_{1, \ux}^\delta$. \end {proof}

\subsection {Refinements of Proposition~\ref{prop:leseqd}, Case II: $\Lambda_1 \in Span_\qq(\Lambda_2, \dots, \Lambda_{\ell})$} 
\label {sec:case2}
When the framing matrix $\Lambda$ is degenerate, we will discuss refinements only in the case when $\Lambda_1$ is in the $\qq$-span of $\Lambda_2, \dots, \Lambda_\ell$. (A model example for this is zero surgery on a knot in the integral homology sphere $Y$.)  Note that in the present case,
$b_1(Y_{\bar\Lambda}(L))=b_1(Y_{\Lambda'}(L)) = b_1(Y_{\Lambda}(L))-1$.

In this situation we are free to choose whether to consider the quasi-isomorphism in \eqref{eq:qi1}, or the
one in \eqref{eq:qi2}. We will focus on the quasi-isomorphism \eqref{eq:qi1}, which involves the cobordisms $W_1$ and $W_2$. We denote by $h$ the smallest positive integer such that
$$ h \Lambda_1 \in Span(\Lambda_2, \dots, \Lambda_\ell).$$

\begin {lemma}
\label {lem:linear}
If $\Lambda_1$ is in the  $\qq$-span of $\Lambda_2, \dots, \Lambda_\ell$, then the vector $\tau_1=(1,0, \dots, 0)$ is not in the $\qq$-span of $\Lambda_1, \dots, \Lambda_\ell$.
\end {lemma}

\begin {proof}
Let us view the framing matrix $\Lambda$ as a self-adjoint linear operator on $\qq^\ell$ (with the standard inner product). Since $\Lambda_1$ is in the $\qq$-span of  $\Lambda_2, \dots, \Lambda_\ell$, there exists a vector $\vs = (v_1, v_2, \dots, v_{\ell})$ in the kernel of $\Lambda$, withe $v_1 \neq 0$. The kernel is orthogonal to the image of $\Lambda$, so since $\vs \cdot \tau_1= v_1 \neq 0$, the image cannot contain $\tau_1$.
\end {proof}

\begin{lemma}
\label{lem:whatever}
The vector $j\tau_1$ is in $Span(\bar\Lambda_1,\Lambda_2,\dots,\Lambda_\ell)$ if and only if $j$ is divisible by $m_1 h$.
\end{lemma}
\begin{proof}
Writing
$$j\tau_1= v_1(\Lambda_1+m_1\tau_1)+v_2\Lambda_2+\dots+v_\ell \Lambda_\ell,$$
we have that $(j-v_1m_1)\tau_1=v_1\Lambda_1+\dots+v_\ell\Lambda_\ell$. 
This can only happen if $j=v_1 m_1$ (by Lemma~\ref{lem:linear}) and $v_1$ is divisible by $h$ (by the definition of $h$).
\end{proof}

\begin {lemma}
\label {lem:l1'}
Suppose $\Lambda_1$ is in the  $\qq$-span of $\Lambda_2, \dots, \Lambda_\ell$, and let $\Lambda_1'= \Lambda_1|_{L'}$ be the vector of linking numbers between $L_1$ and the other components. Then $j\Lambda_1' \in Span(\Lambda_2', \dots, \Lambda_\ell')$ if and only if $j$ is a multiple of $h$.
\end {lemma}

\begin {proof}
If $i$ is a multiple of $h$, we already know that $i \Lambda_1 \in Span(\Lambda_2, \dots, \Lambda_\ell)$. For the converse, suppose $j\Lambda_1' = \sum v_i \Lambda_2'$ for some $v_i \in \zz, i=2, \dots, \ell$. Then $j\Lambda_1 - \sum v_i \Lambda_i$ is a multiple of $\tau_1$. From Lemma~\ref{lem:linear} we see that it must be zero, so $j$ is a multiple of $h$.
\end {proof}

\begin {lemma}
\label {lem:injj}
If $\Lambda_1$ is in the  $\qq$-span of $\Lambda_2, \dots, \Lambda_\ell$, the natural restriction map
$$ \spc(W_1) \to \spc(\del W_1) \cong \spc (Y_{\Lambda}(L) \# L(m_1, 1)) \times \spc(Y_{\bar \Lambda}(L))$$
is injective.
\end {lemma}

\begin {proof}
We prove the similar statement for second cohomology groups, by looking at the kernel of the correspondinghomomorphism. Suppose $\s=(s_1, \dots, s_{\ell+2}) \in \zz^{\ell+2}$ is such that $$[\s] \in H^2(W_1) \cong \zz^{\ell+2} / Span(\Lambda^+_1, \dots, \Lambda^+_{\ell+1})$$ has trivial projection to the boundary. We need to show that $[\s] =0$.

By hypothesis, the vector $(s_1, \dots, s_{\ell})$ is in the image of $\Lambda$, and $s_{\ell+1}$ is a multiple of $m_1$. By adding suitable multiples of $\Lambda^+_1, \dots, \Lambda^+_{\ell+1}$, 
we can assume that $\s=(0, \dots, 0, s_{\ell+2})$ without changing the equivalence class $[\s]$. Also by hypothesis, there exist integers $b_i, i=1, \dots, \ell+2$ such that
$$ \s = b_1  \Lambda_1^+  \dots + b_{\ell+2} \Lambda_{\ell+2}^+.$$

Restricting attention to the first $\ell$ coordinates of $\s$ we get that $(b_{\ell+2}, 0, \dots, 0)$ is in the span of $\Lambda_1, \dots, \Lambda_{\ell}$. From Lemma~\ref{lem:linear} we see that $b_{\ell+2}=0$. 
Hence $\s$ is in the span on $\Lambda^+_1, \dots, \Lambda^+_{\ell+1}$, as desired.
\end {proof}

Let us define an equivalence relation on $\spc$ structures on $Y_{\Lambda}(L)$ as follows. For $\ux_1, \ux_2 \in \spc(Y_{\Lambda}(L)) \cong \H(L)/H(L, \Lambda)$, we say $\ux_1$ is equivalent to $\ux_2$ if and only if there exist $\s_1, \s_2 \in \H(L)$ with
$$ [\s_1] = \ux_1,  \ [\s_2]=\ux_2, \ \text{and } \s_1 - \s_2 = jm_1\tau_1,$$
for some $j \in \zz$.

We call a $\spc$ structure $\ux$ on $Y_{\Lambda}(L)$ {\em small} if our complex $\CFmd(Y_{\Lambda}(L)) = \CFmd(\Ta, \Tb, \ws)$ has at least one generator in that $\spc$ structure.  Clearly there 
are only finitely many small $\spc$ structures. Further, for $m_1 \gg 0$, using Lemma~\ref{lem:linear} we see that every equivalence class of $\spc$ structures on $Y_{\Lambda}(L)$ contains at most one small structure. Let us pick one representative (in $\spc(Y_{\Lambda}(L))$) from each equivalence class, in such a way that all small structures are picked. We call the chosen $\spc$ structures {\em special}. By construction, for every special  $\ux$ we have
\begin {equation}
\label {eq:unique}
\CFmd(Y_{\Lambda}(L), \ux) = \bigoplus_{\ux_1 \sim \ux} \CFmd(Y_{\Lambda}(L), \ux).
\end {equation}
 
Building up on Definition~\ref{def:linked}, we write $\ux \sim \bar \ux$ if  $\ux \in  \spc(Y_{\Lambda}(L))$
and $\bar \ux \in  \spc(Y_{\bar \Lambda}(L))$ are linked. Further, we say that two $\spc$ structures $\bar \ux$ on $Y_{\bar \Lambda}(L)$ and $\ux'$ on $Y_{\Lambda'}(L')$ are {\em linked} (and write $\bar \ux \sim \ux'$) if there exists a $\spc$ structure on $W_2$ interpolating between the two. We also say that $\ux \in  \spc(Y_{\Lambda}(L))$ is {\em linked} to $\ux' \in \spc(Y_{\Lambda'}(L'))$ (and write $\ux \sim \ux'$) if there exists $\bar \ux \in  \spc(Y_{\bar \Lambda}(L))$ such that $\ux \sim \bar \ux$ and $\bar \ux \sim \ux'$.
  
The following lemma describes how the $\spc$ structures on the three manifolds $Y_{\Lambda}(L), Y_{\bar \Lambda}$ and $Y_{\Lambda'}$ are linked to each other. 

\begin {lemma}
\label {lem:lk}
Suppose that $\Lambda_1$ is in the  $\qq$-span of $\Lambda_2, \dots, \Lambda_\ell$, and $m_1 > 0$ is sufficiently large. Then: 

(a) Every $\ux \in  \spc(Y_{\Lambda}(L))$ is linked to exactly $h$ $\spc$ structures on $Y_{\bar \Lambda}(L)$, and every $\bar \ux \in  \spc(Y_{\bar \Lambda}(L))$ is linked to exactly one special $\spc$ structure on $Y_{ \Lambda}(L)$.

(b) Every $\bar \ux \in  \spc(Y_{\bar \Lambda}(L))$ is linked to exactly  $h$ $\spc$ structures on $Y_{ \Lambda'}(L')$, and every $\ux' \in \spc(Y_{\Lambda'}(L'))$ is linked to exactly $m_1h$ $\spc$ structures on $Y_{\bar \Lambda}(L)$.

(c) Every $\ux \in  \spc(Y_{\Lambda}(L))$ is linked to exactly  $h$ $\spc$ structures on $Y_{ \Lambda'}(L')$, and every $\ux' \in \spc(Y_{\Lambda'}(L'))$ is linked to exactly $m_1$ special $\spc$ structures on $Y_{\Lambda}(L)$.
\end {lemma}

\begin{figure}
\begin{center}
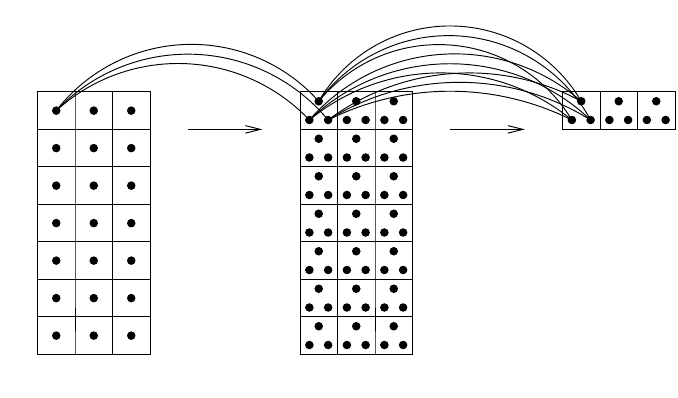
\end{center}
\caption {{\bf Linking of $\spc$ structures.} Each dot in the first rectangle corresponds to a special $\spc$ structure on $Y_{\Lambda}(L)$, and is linked to the three dots (representing $\spc$ structures on $Y_{\bar \Lambda}(L)$) in the corresponding square of the second rectangle. Furthermore, all the dots in the same column of the second rectangle are linked to the three dots (representing $\spc$ structures on $Y_{\Lambda'}(L')$) in the corresponding column of the third rectangle.
}
\label{fig:spinc}
\end{figure}

Before proving the lemma, it is useful to illustrate its content in a particular example, shown graphically in Figure~\ref{fig:spinc}. We consider a framed link $L = L_1 \cup L_2$ with framing matrix 
$$\Lambda = \begin{pmatrix} 1 & 3 \\ 3 & 9 \end{pmatrix},$$ so that $h=3$ and $\H(L) = (\zz + \frac{1}{2})^2$. The set of $\spc$ structures on $Y_{\Lambda}(L)$ is identified with $\H(L)/\langle (1,3) \rangle$. The ones that differ by $(m_1, 0)$ are called equivalent. Thus, there are only $3m_1$ special $\spc$ structures (one from each equivalence class). The set of special $\spc$ structures on $Y_{\Lambda}(L)$ is the left rectangle in the picture (where $m_1 = 7$), with each structure being represented by a black dot. Going down one square in the rectangle (in a cyclical fashion) corresponds to adding the vector $(1,0)$. Going left  corresponds to adding the vector $(0,1)$, and is also done in a cyclical fashion, except when we go from the third to the first column we also move one step up, according to the relation $(0,3) = (-1,0)$. 

The second rectangle in Figure~\ref{fig:spinc} represents the space of all $9m_1$ $\spc$ structures on $Y_{\bar \Lambda}(L)$, where 
$$\bar \Lambda = \begin{pmatrix} m_1 + 1 & 3 \\ 3 & 9 \end{pmatrix}.$$ 
The dots in the same square differ from each other by multiples of $(m_1, 0)$. Going down one square in the rectangle  still means adding $(1,0)$, and going left means adding $(0,1)$. When we go horizontally in a cycle we move up one square as well, according to the relation $(0,3) = (-m_1-1,0)$. 

Finally, the rectangle on the right of Figure~\ref{fig:spinc} represents the set of $\spc$ structures on $Y_{\Lambda'}(L')$, which is simply identified with $\zz/9$. The dots in the same square differ from each other by multiples of $3$. 

The linking of $\spc$ structures between the three rectangles is as shown in Figure~\ref{fig:spinc}, and corresponds to the description in Lemma~\ref{lem:lk}.

\begin {proof}[Proof of Lemma~\ref{lem:lk}]
(a) Consider two $\spc$ structures $\bar \ux_1, \bar \ux_2$ on $Y_{\bar \Lambda}(L)$. By Lemma~\ref{lem:linked}, $\bar \ux_1 $ and $\bar \ux_2$ are linked to a single $\ux \in \spc(Y_{\Lambda}(L))$ if and only if the difference 
$$\bar \ux_1 - \bar \ux_2 \in H^2(Y_{\bar \Lambda}(L)) \cong \zz^\ell/H(L, \bar \Lambda)$$  can be represented (in $\zz^\ell$) by a multiple of $m_1 \tau_1 = \bar \Lambda_1 - \Lambda_1$. The smallest such multiple that lies in the span of $\bar \Lambda_1, \dots, \Lambda_\ell$ (i.e. represents a trivial cohomology class) is $h(\bar \Lambda_1 - \Lambda_1)$. 

In the other direction, consider two $\spc$ structures $\ux_1, \ux_2$ on $Y_{\Lambda}(L)$.  Applying Lemma~\ref{lem:linked} again, we see that $\ux_1$ and $\ux_2$ can be linked to a single $\bar \ux \in \spc(Y_{\bar \Lambda}(L))$ if and only if the difference $$\ux_1 - \ux_2 \in H^2(Y_{\Lambda}(L)) \cong \zz^\ell/H(L, \Lambda)$$ can be represented (in $\zz^\ell$) by a multiple of $m_1 \tau_1$, i.e. $\ux_1$ and $\ux_2$ are equivalent. There is a unique special structure in each equivalence class, by construction.

(b) Apply Lemma~\ref{lem:HH} (b) to $W_{\bar \Lambda}(L', L)$, which is $-W_2$ turned upside down. Two structures $\ux_1', \ux_2' \in \spc(Y_{\Lambda'}(L')) \cong \H(L')/H(L', \Lambda')$ are linked to a single one on $Y_{\bar \Lambda}(L)$
if and only if 
$$ \ux_1' - \ux_2' \in H^2(Y_{\Lambda'}(L')) \cong \zz^{\ell-1}/H(L', \Lambda')$$
can be represented (in $\zz^{\ell-1}$) by a multiple of $\bar \Lambda_1|_{L'} = \Lambda'_1$. From Lemma~\ref{lem:l1'} we see 
that $j \Lambda'_1$ is zero in cohomology if and only if the respective factor $j$ is divisible by $h$. 

In the other direction, two $\spc$ structures $\bar \ux_1, \bar \ux_2$ on $Y_{\bar \Lambda}(L)$ are linked to a single one on $Y_{\Lambda'}(L')$ if and only if $\bar \ux_1 - \bar \ux_2$ can be represented (in $\zz^\ell$) by a multiple of $\tau_1$. By Lemma~\ref{lem:whatever}, the first such multiple that lies in $H(L, \bar \Lambda)$ is $m_1h\tau_1$.

(c) Use the descriptions of linking in (a) and (b). We get that two $\spc$ structures $\ux'_1, \ux_2'$ on $ Y_{\Lambda'}(L')$ are linked to a single one on $Y_{ \Lambda}(L)$  if and only if $\ux'_1 - \ux'_2$ can be represented by a multiple of $\Lambda'_1$, and
two $\spc$ structures $ \ux_1, \ux_2$ on $Y_{ \Lambda}(L)$ are linked to a single one on $Y_{\Lambda'}(L')$ if and only if $ \ux_1 -  \ux_2$ can be represented by a multiple of $\tau_1$. 
\end {proof}

From \eqref{eq:unique} and Lemma~\ref{lem:lk} (a) we see that the cobordism map $f_1^\delta$ splits into the direct sum of the maps
$$ f_{1, \ux}^{\delta}: \CFmd(Y_{\Lambda}(L), \ux) \longrightarrow \bigoplus_{\bar \ux \sim \ux} \CFmd(Y_{\bar \Lambda}(L), \bar \ux) ,$$
over special $\spc$ structures $\ux$. 

Turning our attention to the map $f_2^\delta$, note that the action of $H^1(Y_{\Lambda'}(L'))\cong H(L', \Lambda')^\perp$ on the module $\TR_1 = \ff[T_1]/(T_1^{m_1} -1 )$ is trivial. Indeed, if $\vs' \in \zz^{\ell-1}$ satisfies $\vs' \cdot \Lambda'_i = 0$ for all $i=2, \dots, \ell$, we must also have $\vs' \cdot \Lambda_1' = 0$ (because $\Lambda_1$ is in the $\qq$-span of $\Lambda_2, \dots, \Lambda_\ell$).

Hence, for every $\ux' \in \spc(Y_{\Lambda'}(L'))$, we have a decomposition
$$ \CFmd(Y_{ \Lambda'}(L'),  \ux'; \TR_1) \cong \bigoplus_{i=0}^{m_1 - 1} T_1^i \cdot  \CFmd(Y_{ \Lambda'}(L'),  \ux').$$

\begin{figure}
\begin{center}
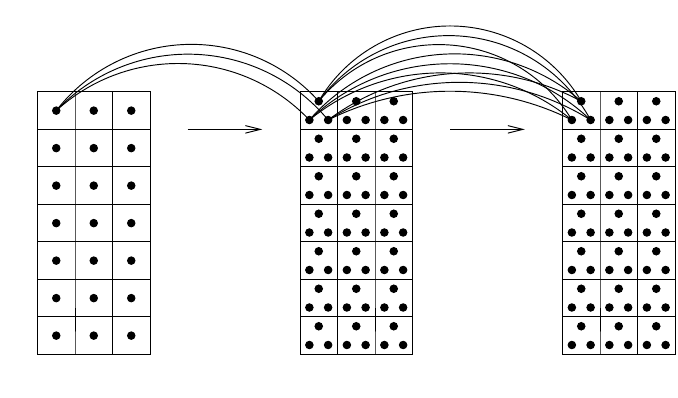
\end{center}
\caption {{\bf Linking of $\spc$ structures, keeping track of the powers of $T_1$.} The map $f_2^\delta$ decomposes as a sum according to powers of $T_1$. Each dot in the rectangle on the right represents a pair $(\ux', i)$, where $\ux' \in \spc(Y_{\Lambda'}(L'))$ and $i \in \zz/_{m_1}$ is an exponent of $T_1$. Going down one square in the rectangle corresponds to multiplication by $T_1$. This way, a dot in the second rectangle is linked (that is, gives a contribution to $f_2^\delta$) only with the three dots in the corresponding square of the third rectangle. Therefore, the pair of maps $(f_1^\delta, f_2^\delta)$ decomposes as a direct sum according to dots in the first rectangle (i.e., special $\spc$ structures $\ux$ on $Y_{\Lambda}(L)$). 
}
\label{fig:spinctwisted}
\end{figure}

Moreover, using the description of linking in Lemma~\ref{lem:lk} (b) (see also Figure~\ref{fig:spinctwisted}), we get that $f_2^\delta$ splits as a direct sum of maps
$$ f_{2, \ux}^{\delta}:  \bigoplus_{\bar \ux \sim \ux} \CFmd(Y_{\bar \Lambda}(L), \bar \ux) \longrightarrow  \bigoplus_{ \ux' \sim \ux} T_1^{\ux}\cdot \CFmd(Y_{ \Lambda'}(L'),  \ux'),$$
over special $\spc$ structures $\ux$ on $Y_{\Lambda}(L)$. Here, $T_1^{\ux}$ denotes $T_1^{(\ux - \ux_0) \cdot \vs}$, where $\ux_0$ is a fixed $\spc$ structure on $Y_{\Lambda}(L)$, and $\vs$ is an arbitrary vector in $H_2(W_2) \cong H(L, \bar \Lambda|_{L'})^\perp$.

Using Lemma~\ref{lem:lk} again, we also obtain a similar decomposition of $H_2^\delta$ into maps $H_{2, \ux}^\delta$, one for each special $\spc$ structure $\ux$ on $Y_{\Lambda}(L)$. Each $\CFmd(Y_{\Lambda}(L), \ux)$ is quasi-isomorphic to the mapping cone of the respective $f_{2, \ux}^\delta$, via the pair $  (f_{1, \ux}^\delta, H_{1, \ux}^{\delta})$.

\begin {proposition}
\label {prop:gr00}
Fix $\delta > 0$ and suppose that $\Lambda_1$ is in the $\qq$-span of $\Lambda_2, \dots, \Lambda_\ell$. Then, for every $m_1 \gg 0$, and for every special $\spc$ structure $\ux$ on $\spc(Y_{\Lambda}(L))$ such that $m_1$ is divisible by $\delt(\ux)$, the mapping cone $Cone (f_{2, \ux}^\delta) $ admits a relative $\zz/\delt(\ux)\zz$-grading, and the quasi-isomorphism
$$  (f_{1, \ux}^\delta, H_{1, \ux}^{\delta}): \CFmd(Y_{\Lambda}(L), \ux) \xrightarrow{\sim}  Cone (f_{2, \ux}^\delta) $$
respects the relative $\zz/\delt(\ux)\zz$-gradings on the two sides.
\end {proposition}

Before proceeding with the proof, we establish a few lemmas. First, note that for $\ux \in \spc(Y_{\Lambda}(L)) \cong \H(L)/H(L, \Lambda)$, Equation~\eqref{eq:dgr} gives
\begin {equation}
\label {eq:dum}
 \delt(\ux) = \gcd_{\{\vs \in \zz^\ell| \Lambda \vs = 0\}} 2\vs \cdot \s,
 \end {equation}
for $[\s]=\ux$. 

\begin{lemma}
\label {lem:baru}
Suppose $\ux \in \spc(Y_{\Lambda}(L))$ and $\bar \ux  \in \spc(Y_{\bar \Lambda}(L))$ are linked. Then, $\delt(\ux)$ divides $\delt(\bar \ux)$.
\end {lemma}

\begin{proof}
If $[\s] = \ux$, then according to the proof of Lemma~\ref{lem:lk} (a) we have $[\s + jm_1 \tau_1] = \bar \ux$ for some $j \in \zz$. We have
$$ \delt(\bar \ux) = \gcd_{\{\vs=(v_1, \dots, v_\ell) \in \zz^\ell| \bar \Lambda \vs = 0\}} 2(\vs \cdot \s + jv_1m_1).$$
From Lemma~\ref{lem:linear} we see that $\bar \Lambda \vs = m_1v_1 \tau_1 + \Lambda \vs = 0$ implies $v_1 = 0$ and $\Lambda \vs = 0$. Hence,  $2(\vs \cdot \s + jv_1m_1) = 2\vs \cdot \s$ is one of the  elements for which we take the greatest common divisor in \eqref{eq:dum}. Hence, $\delt(\ux)$ divides $\delt(\bar \ux)$.
\end {proof}

\begin{lemma}
\label {lem:uprime}
Suppose $\ux \in \spc(Y_{\Lambda}(L))$ and $\ux'  \in \spc(Y_{ \Lambda'}(L'))$ are linked. Then, $\delt(\ux)$ divides $\delt(\ux')$.
\end {lemma}

\begin {proof}
If $\ux = [\s] $ with $\s = (s_1, \dots, s_\ell)$ and we let $\s' = (s_2, \dots, s_\ell)$, then according to the proof of Lemma~\ref{lem:lk} (c), we have $\ux' = [\s' + j \Lambda_1']$ for some $j \in \zz$. We have
$$ \delt(\ux') = \gcd_{\{ \vs' \in \zz^{\ell-1}| \Lambda' \vs' = 0\}} 2 \vs' \cdot (\s' + j \Lambda_1').$$

Note that $\Lambda' \vs' = 0$ implies $\Lambda_1' \cdot \vs' = 0$, so the expression in the $\gcd$ is $2\vs' \cdot \s'$. Moreover, if we let $\vs \in \zz^\ell$ be the vector with first coordinate zero and the others given by $\vs'$, then $\Lambda \vs = (u, 0, \dots, 0)$ for some $u \in \zz$.   Lemma~\ref{lem:linear} implies that $u = 0$, so $\vs$ gets counted in \eqref{eq:dum}, and the conclusion follows.
\end {proof}

\begin{proof}[Proof of Proposition~\ref{prop:gr00}] Let $\ux$ be a special $\spc$ structure on $Y_{\Lambda}(L)$ with $m_1$ divisible by $\delt(\ux)$. Lemma~\ref{lem:twistedmaps} (a), together with Lemmas~\ref{lem:injj} and \ref{lem:baru} imply that we can equip 
$$\bigoplus_{\bar \ux \sim \ux} \CFmd(Y_{\bar \Lambda}(L), \bar \ux) $$
with a relative $\zz/\delt(\ux)\zz$-grading such that $f_{1, \ux}^\delta$ is grading-preserving. Furthermore, Lemma~\ref{lem:uprime} shows that each direct summand $T_1^{\ux}\cdot \CFmd(Y_{ \Lambda'}(L'),  \ux')$ in the target of $f_{2, \ux}^\delta$ is relatively $\zz/\delt(\ux)\zz$-graded. 

We would like to give the mapping cone of $f_{2, \ux}^\delta$ a relative $\zz/\delt(\ux)\zz$-grading such that $f_{1, \ux}^\delta$ is grading-preserving. To be able to do this, we need to check that for any two $\spc$ structures $\t, \t + u$ on $W_1 \cup W_2$ which restrict to $\ux$ on $Y_{\Lambda}(L)$, and to the same $\spc$ structure on $Y_{\Lambda'}(L')$, the contributions to the composition $f_{2, \ux}^\delta \circ 
f_{1, \ux}^\delta$ coming from $\t$ and $\t+u$ shift the relative $\zz/\delt(\ux)\zz$-gradings by the same amount. Using Lemma~\ref{lem:twistedmaps} (b), this is equivalent to showing that
\begin {equation}
\label {eq:cu}
\langle  c_1(\t) \smallsmile u + u \smallsmile u , [W_1 \cup W_2] \rangle \equiv 0 (\mod \delt(\ux)).
\end {equation}
Here $u \in H^2(W_1 \cup W_2, \del(W_1 \cup W_2); \zz)$ is Poincar\'e dual to some class in $H_2(W_1 \cup W_2; \zz)$. 

The cobordism $W_1 \cup W_2$ admits the following description in terms of surgery. Consider the link $L^+$ from Figure~\ref{fig:kirby}, and add an zero-framed, unknotted component $L_{\ell+3}$, which forms a Hopf link (in particular, has linking number one) with $L_{\ell+1}$, and is unlinked with the other components of $L^+$. Call the resulting framed link $(L^{++}, \Lambda^{++})$. Then, in the notation of Section~\ref{sec:4c}, we have
$$ W_1 \cup W_2 = W_{\Lambda^{++}}(L_1 \cup \dots \cup L_{\ell+1}, L^{++}).$$
 
In matrix form, we have
$$ \Lambda^{++} =   \left (
\begin{array}{cccc|ccc} 
 & & & & 0 &  1 & 0\\
 & \Lambda& & & 0 & 0 & 0\\ 
 & & & & \vdots & \vdots &\\ 
 & & & & 0 & 0 & 0\\
 \hline 
 0&0 & \ldots& 0 & m_1& 1 & 1\\ 
 1& 0& \ldots & 0 & 1& 0 &0 \\
 0 & 0 & \ldots & 0 & 1 & 0 & 0
\end {array}
\right )
$$

Using Lemma~\ref{lem:HH} (a), we see that 
$$ \text{PD}(u) = [(v_1, \dots, v_\ell, j, 0, -m_1j)] \in H_2(W_1 \cup W_2)$$
for some $j \in \zz$ and $\vs = (v_1, \dots, v_\ell) \in H(L, \Lambda)^\perp$. We get
$$ \langle u \smallsmile u , [W_1 \cup W_2] \rangle = -m_1 j^2,$$
which is divisible by $\delt(\ux)$ (because we assumed $m_1$ was so).
Also, if $\ux = [\s]$ with $\s= (s_1, \dots, s_\ell) \in \H(L)$, then we must have
$$ \t = [(s_1 + 1/2, s_2, \dots, s_\ell, m_1 + 1, s_{\ell+2}, s_{\ell+3})] \in \H(L^{++})/H(L^{++}, \Lambda^{++}),$$
which gives
$$ \langle  c_1(\t) \smallsmile u , [W_1 \cup W_2] \rangle = 2\s \cdot \vs + 2m_1j(1-s_{\ell+3}).$$
This is divisible by $\delt(\ux)$ because of \eqref{eq:dum} and by our choice of $m_1$. 

Therefore, \eqref{eq:cu} holds, and $f_{1, \ux}^{\delta}$ can be made grading-preserving. The homotopy $H_{1, \ux}^\delt$ is also automatically grading-preserving.
\end {proof}
 
 \begin {corollary}
 \label {cor:small}
Suppose  $\Lambda_1$ is in the $\qq$-span of $\Lambda_2, \dots, \Lambda_\ell$, and fix $\delta > 0$. Then for every $m_1 \gg 0$ suitably chosen, and for every small $\spc$ structure $\ux$ on $\spc(Y_{\Lambda}(L))$, the mapping cone $Cone (f_{2, \ux}^\delta) $ admits a relative $\zz/\delt(\ux)\zz$-grading, and the quasi-isomorphism $ (f_{1, \ux}^\delta, H_{1, \ux}^{\delta})$ is grading-preserving.
\end {corollary}
\begin {proof}
Since there are only finitely many small $\spc$ structures $\ux$, we can choose $m_1$ to be a multiple of $\delt(\ux)$ for all $\ux$ small. We then apply Proposition~\ref{prop:gr00}.
\end {proof}

Proposition~\ref{prop:gr00} and Corollary~\ref{cor:small} give a (partial) grading-preserving decomposition of the quasi-isomorphism $(f_1^\delta, H_1^\delta)$ from \eqref{eq:qi1}. Similar arguments can be applied to the quasi-isomorphism $(H_2^\delta, f_3^\delta)$ from \eqref{eq:qi2}. Indeed, one can check that an analogue of Lemma~\ref{lem:injj} holds for the cobordism $W_3$; i.e., the restriction map
$$ \spc(W_3) \to \spc(\del W_3) \cong \spc(Y_{\Lambda'}(L')) \times \spc (Y_{\Lambda}(L) \# L(m_1, 1)) $$
is injective. Further, the maps $f_3^\delta$ and $H_2^\delta$ split into direct sums of maps $f_{3, \ux}^\delta$ and $H_{2, \ux}^\delta$, according to special $\spc$ structures $\ux$ on $Y_{\Lambda}(L)$. A study of grading differences similar to that in the proof of Proposition~\ref{prop:gr00} yields the following:

\begin {proposition}
\label {prop:gr002}
Fix $\delta > 0$ and suppose that $\Lambda_1$ is in the $\qq$-span of $\Lambda_2, \dots, \Lambda_\ell$. Then, for every $m_1 \gg 0$, and for every special $\spc$ structure $\ux$ on $\spc(Y_{\Lambda}(L))$ such that $m_1$ is divisible by $\delt(\ux)$, the quasi-isomorphism
$$  (H_{2, \ux}^\delta, f_{3, \ux}^{\delta}):Cone (f_{2, \ux}^\delta) \xrightarrow{\sim}  \CFmd(Y_{\Lambda}(L), \ux)$$
respects the relative $\zz/\delt(\ux)\zz$-gradings on the two sides.
\end {proposition}

To conclude this subsection we state a twisted coefficients generalization of Propositions~\ref{prop:gr00} and \ref{prop:gr002}. We work in the following setting. Recall that $(\Sigma, \alphas, \betas, \ws), (\Sigma, \alphas, \gammas, \ws)$ and $(\Sigma, \alphas, \deltas, \ws)$ are the Heegaard diagrams 
for $Y_{\Lambda}(L), Y_{\bar \Lambda}(L)$, and $Y_{\Lambda'}(L')$, respectively. Also, for $i < g+k-1$, the curves $\beta_i, \gamma_i$ and $\delta_i$ approximate one another. For $i=1, \dots, g+k-1$, let us place points $p_i, q_i$ on each side of the curve $\beta_i$, such that they can be joined by an arc that intersects $\beta_i, \gamma_i$ and $\delta_i$ once each, and does not intersect any of the alpha curves. (These are the analogues of $w_1$ and $z_1$.) Let $n_1, \dots, n_{g+k-2}$ be nonnegative integers. Consider the ring
\begin {equation}
\label {eq:sr}
 \SR = \ff[S_1, \dots, S_{g+k-2}]/(S_1^{n_1} - 1, \dots, S_{g+k-1}^{n_{g+k-2}} - 1).
 \end {equation}

We can then construct Floer complexes with twisted coefficients
$$ \CFmd(Y_{\Lambda}(L); \SR), \ \CFmd(Y_{\bar \Lambda}(L); \SR),\  \CFmd(Y_{\Lambda'}(L'); \SR \otimes \TR_1)$$ 
and cobordism maps $f_{1; \SR}^\delta, f_{2; \SR}^\delta, f_{3; \SR}^\delta$ relating them, by counting all pseudo-holomorphic strips and triangles using the coefficients
$$ S_1^{n_{p_1}(\phi) - n_{q_1}(\phi)} \cdot \dots \cdot S_{g+k-2}^{n_{p_{g+k-2}}(\phi) - n_{q_{g+k-2}}(\phi)},$$
where $\phi$ denotes the respective relative homology class. (This is all in addition to the power of $T_1$ coming from the twisting by $\TR_1$.) 

A straightforward generalization of Proposition~\ref{prop:leseqd} shows that $\CFmd(Y_{\Lambda}(L); \SR)$ is quasi-isomorphic to the mapping cone of $f_{2; \SR}^\delta$.  Note that twisting by $\SR$ can have a non-trivial effect on the Floer homology of $Y_{\Lambda}(L)$ only when $b_1(Y_{\Lambda}(L)) > 0$. Supposing further that $\Lambda_1$ is in the $\qq$-span of $\Lambda_2, \dots, \Lambda_\ell$ (as we did in this section), we have twisted coefficients analogues $f_{i, \ux; \SR}^\delta, H_{i, \ux; \SR}^\delta$ of the maps $f_{i, \ux}^\delta$ and the homotopies $H_{i, \ux}^\delta$, respectively.

\begin {proposition}
\label {prop:gr00twisted}
Fix $\delta > 0$ and suppose that $\Lambda_1$ is in the $\qq$-span of $\Lambda_2, \dots, \Lambda_\ell$. Then, for every $m_1 \gg 0$, and for every special $\spc$ structure $\ux$ on $\spc(Y_{\Lambda}(L))$ such that $m_1$ is divisible by $\delt(\ux, \SR)$, the mapping cone $Cone (f_{2, \ux; \SR}^\delta) $ admits a relative $\zz/\delt(\ux, \SR)\zz$-grading, and the quasi-isomorphisms
$$  (f_{1, \ux; \SR}^\delta, H_{1, \ux; \SR}^{\delta}): \CFmd(Y_{\Lambda}(L), \ux; \SR) \xrightarrow{\sim}  Cone (f_{2, \ux; \SR}^\delta) $$
and
$$  (H_{2, \ux; \SR}^\delta, f_{3, \ux; \SR}^{\delta}):Cone (f_{2, \ux; \SR}^\delta) \xrightarrow{\sim}  \CFmd(Y_{\Lambda}(L), \ux; \SR)$$
respect the relative $\zz/\delt(\ux, \SR)\zz$-gradings on the two sides.
\end {proposition}

The proof of Proposition~\ref{prop:gr00twisted} is similar to those of Propositions ~\ref{prop:gr00} and \ref{prop:gr002}.

\subsection {Cobordism maps} \label {sec:x}

Propositions~\ref{prop:gr+}, \ref{prop:gr-}, \ref{prop:gr0} and \ref{prop:gr00} all describe quasi-isomorphisms between chain complexes of the form $\CFmd(Y_{\Lambda}(L), \ux)$ and mapping cones $Cone(f_{2, \ux}^\delta)$. Our goal in this section is to describe commutative diagrams which relate the inclusion of the target of $f_{2, \ux}$ into $Cone(f_{2, \ux}^\delta)$, with a cobordism map from the Floer complex (possibly with twisted coefficients) of $Y_{\Lambda'}(L')$ to the Floer complex of $Y_{\Lambda}(L)$. These diagrams will prove useful in our description of cobordism maps in terms of surgery, in Section~\ref{sec:surgerymaps}.

We state several results, on a case-by-case basis, corresponding to the case analysis 
described above. We will always denote by $f_{3; \t}^\delta$ the map (with twisted coefficients) induced by the cobordism $W_3$ with a $\spc$ structure $\t$. We start with Case I. In Subcase I (a), we have:

\begin {proposition}
\label {prop:cob+}
Fix $\delta > 0$ and $i \in \zz$. Suppose $\Lambda$ is nondegenerate, $a_1 > 0$, and $m_1$ is a sufficiently large multiple of $a_1h$. Suppose $\s_0 \in \H(L)$ is such that $[\s_0]=\t_0 \in \spc(W_2)$ is the base $\spc$ structure used to define the expression $T_1^\t$ in \eqref{eq:f2ut}. If $\s \in \H(L)$ is such that $\a \cdot (\s - \s_0) = i$, let us denote 
$\ux=[\s] \in \spc(Y_{\Lambda}(L)) \cong \H(L)/H(L, \Lambda), \  \ux' = [\psi^{L_1}(\s)] \in \spc(Y_{\Lambda'}(L')) \cong \H(L')/H(L', \Lambda')$, and 
$$\t_i = [\s] \in \spc(W_3) \cong \H(L)/Span(\Lambda_2, \dots, \Lambda_\ell).$$ Then, there is a diagram
$$
\begin{CD}
\CFmd(Y_{\Lambda'}(L'), \ux') @>{f_{3; \t_i}^\delta}>> \CFmd(Y_{\Lambda}(L), \ux)\\
@V{\cong}VV @VV{(f_{1,\ux}^\delta, H_{1, \ux}^\delta)}V \\
T_1^i \cdot \CFmd(Y_{\Lambda'}(L'), \ux') @>>> Cone(f_{2, \ux}^\delta),
\end {CD}
$$
commuting up to chain homotopy. Here, the bottom horizontal arrow is the inclusion into the mapping cone.
\end {proposition}

\begin {proof}
The existence of such a diagram is a consequence of the proof of the fact that $(f_{1,\ux}^\delta, H_{1, \ux}^\delta)$ is a quasi-isomorphism, see the homological algebra Lemma 4.4 in \cite{BrDCov}, compare also \cite[Proof of Theorem 4.2]{IntSurg}. A priori, the top horizontal arrow is a sum of all maps $f_{3; \t}^\delta$, over all $\t \in \spc(W_3)$ such that $\t$ and $\t_i$ have the same restrictions to $\del W_3$, and $\a \cdot (\t - [\s_0]) \equiv i (\mod m_1)$. However, any such $\t$ differs from $\t_i$ by a multiple of $d \cdot \text{PD}(\Sigma_3)$, where $\Sigma_3 = \s$ is the generator of $H_2(W_3)$. Since $i$ and $\delta$ are fixed, if $m_1$ (and hence $d=m_1/h$) is chosen sufficiently large, all $\t_i + dj\cdot \text{PD}(\Sigma_3)$ produce trivial cobordism maps for $j \in \zz, j\neq 0$.
\end {proof}

We have a similar result in Subcase I (b):
\begin {proposition}
\label {prop:cob-}
Fix $\delta > 0$ and $i \in \zz$. Suppose $\Lambda$ is nondegenerate, $a_1 < 0$, and $m_1$ is a sufficiently large multiple of $a_1h$. Suppose $\s_0, \s \in \H(L)$ (with $\a \cdot (\s - \s_0) = i$), $\ux, \ux'$ and $\t_i$ be as in Proposition~\ref{prop:cob+}.  Then, there is a commutative diagram
$$
\begin{CD}
T_1^i \cdot \CFmd(Y_{\Lambda'}(L'), \ux') @>>> Cone(f_{2, \ux}^\delta),\\
@V{\cong}VV @VV{(H_{2, \ux}^\delta, f_{3,\ux}^\delta)}V \\
\CFmd(Y_{\Lambda'}(L'), \ux') @>{f_{3; \t_i}^\delta}>> \CFmd(Y_{\Lambda}(L), \ux),
\end {CD}
$$
where the top horizontal arrow is the inclusion into the mapping cone.
\end {proposition}

\begin {proof}
This is simpler than Proposition~\ref{prop:gr+} (and the diagram commutes on the nose, rather than only up to chain homotopy), because the relevant quasi-isomorphism \eqref{eq:qi2} already involves $f_{3, \ux}$. The fact that the only contribution to $f_{3, \ux}$ comes from $\t_i$ (for $\delta, i$ fixed and $m_1$ large) is a consequence of the proof of Proposition~\ref{prop:gr-}.
\end {proof}

In Subcase I (c), we have
\begin {proposition}
\label {prop:cob0}
Fix $\delta > 0$. Suppose $\Lambda$ is nondegenerate, $a_1 = 0$, and $m_1$ is a sufficiently large multiple of $h$. Suppose 
$$\t \in \spc(W_3) \cong \H(L)/ Span(\Lambda_2, \dots, \Lambda_\ell)$$ 
has restrictions $\ux'$ to $Y_{\Lambda'}(L')$ and $ \ux $ to $Y_{\Lambda}(L)$.  Then, there is a diagram
$$
\begin{CD}
\CFmd(Y_{\Lambda'}(L'), \ux'; \V_1) @>{f_{3; \t}^\delta}>> \CFmd(Y_{\Lambda}(L), \ux)\\
@V{\cong}VV @VV{(f_{1,\ux}^\delta, H_{1, \ux}^\delta)}V \\
T_1^\ux \cdot \CFmd(Y_{\Lambda'}(L'), \ux'; \V_1) @>>> Cone(f_{2, \ux}^\delta),
\end {CD}
$$
commuting up to chain homotopy, with the bottom horizontal arrow being inclusion into the mapping cone.
\end {proposition}

\begin {proof}
The existence of the commutative diagram follows from the same reasoning as in Proposition~\ref{prop:cob+}, with the top arrow being the sum of cobordism maps over all possible $\spc$ structures on $W_3$ with the given restrictions to the boundary. However, we claim that $\t$ is the unique such $\spc$ structure. For this, it suffices to prove that the restriction
\begin {equation}
\label {eq:w3res}
H^2(W_3) \longrightarrow H^2(\del W_3) \cong H^2(Y_{\Lambda}(L)) \oplus H^2(Y_{\Lambda'}(L'))
\end {equation}
is injective. Suppose $[\vs] \in H^2(W_3) = \zz^\ell/Span(\Lambda_2, \dots, \Lambda_\ell)$ is in the kernel. Because $[\vs]$ becomes trivial when projected to its last $\ell-1$ components, without loss of generality we can assume $\vs = (j, 0, \dots, 0)$ for some $j \in \zz$. We also know that $\vs \in H(L, \Lambda)$, so $j$ must be a multiple of $h$. Since $(h, 0, \dots, 0)$ is in the span of $\Lambda_2, \dots, \Lambda_\ell$ (because $a_1 = 0$), we must have $[\vs] = 0$. This completes the proof.
\end {proof}

Proposition~\ref{prop:cob0} involves the map $f_{3; \t}^{\delta}$, which uses twisted coefficients. From here one can also similarly identify the usual (untwisted) map corresponding to the cobordism $W_3 = W_{\Lambda}(L', L)$.  Indeed, the untwisted complex $\CFmd(Y_{\Lambda'}(L'), \ux')$ can be viewed as a subcomplex of $\CFmd(Y_{\Lambda'}(L'), \ux'; \V_1)$ via the inclusion
\begin {eqnarray}
\label {eq:jota}
\iota: \CFmd(Y_{\Lambda'}(L'), \ux') & \hookrightarrow& \CFmd(Y_{\Lambda'}(L'), \ux'; \V_1), \\
 \x &\mapsto& \x(1 + V_1 + \dots + V_1^{d-1}). \notag
\end {eqnarray}

The composition $f_{3; \t}^{\delta} \circ \iota$ is the untwisted chain map representing the cobordism $W_3$. By pre-composing with $\iota$ in the diagram from Proposition~\ref{prop:cob0}, we get the following corollary:

\begin {proposition}
\label {prop:cob0cor}
Under the hypotheses of Proposition~\ref{prop:cob0}, we have a diagram
$$
\begin{CD}
\CFmd(Y_{\Lambda'}(L'), \ux') @>{f_{3; \t}^\delta \circ \iota}>> \CFmd(Y_{\Lambda}(L), \ux)\\
@V{\cong}VV @VV{(f_{1,\ux}^\delta, H_{1, \ux}^\delta)}V \\
T_1^\ux \cdot \CFmd(Y_{\Lambda'}(L'), \ux') @>>> Cone(f_{2, \ux}^\delta),
\end {CD}
$$
commuting up to chain homotopy, with the bottom horizontal map being the composition of $\iota$ and the inclusion into the mapping cone.
\end {proposition}

Finally, we have the following result in Case II:
\begin {proposition}
\label {prop:cob00}
Fix $\delta > 0$. Suppose $\Lambda_1$ is in the $\qq$-span of $\Lambda_2, \dots, \Lambda_\ell$, and $m_1$ is suitably chosen sufficiently large. Suppose 
$$\t \in \spc(W_3) \cong \H(L)/ Span(\Lambda_2, \dots, \Lambda_\ell)$$ 
has restrictions $\ux'$ to $Y_{\Lambda'}(L')$ and $ \ux $ to $Y_{\Lambda}(L)$, such that $\ux$ is special.  Then:

$(a)$  There is a diagram
$$
\begin{CD}
\CFmd(Y_{\Lambda'}(L'), \ux') @>{f_{3; \t}^\delta}>> \CFmd(Y_{\Lambda}(L), \ux)\\
@V{\cong}VV @VV{(f_{1,\ux}^\delta, H_{1, \ux}^\delta)}V \\
T_1^\ux \cdot \CFmd(Y_{\Lambda'}(L'), \ux') @>>> Cone(f_{2, \ux}^\delta),
\end {CD}
$$
commuting up to chain homotopy, with the bottom horizontal arrow being inclusion into the mapping cone.

$(b)$  There is a diagram
$$
\begin{CD}
T_1^\ux \cdot \CFmd(Y_{\Lambda'}(L'), \ux') @>>> Cone(f_{2, \ux}^\delta),\\
@V{\cong}VV @VV{(H_{2, \ux}^\delta, f_{3,\ux}^\delta)}V \\
\CFmd(Y_{\Lambda'}(L'), \ux') @>{f_{3; \t}^\delta}>> \CFmd(Y_{\Lambda}(L), \ux),
\end {CD}
$$
commuting up to chain homotopy, with the top horizontal arrow being inclusion into the mapping cone.
\end {proposition}

\begin {proof}
For both (a) and (b), just as in the proof of Proposition~\ref{prop:cob0}, it suffices to prove that the restriction map \eqref{eq:w3res} is injective. As before, if $[\vs]$ is in the kernel, we can assume that $\vs = (j, 0, \dots, 0)$ for some $j \in \zz$. Lemma~\ref{lem:linear} then shows that $\vs \in H(L, \Lambda)$ only if $j=0$.
\end {proof}

Let us also state the twisted coefficients generalization of Proposition~\ref{prop:cob00}, which has a similar proof:
\begin {proposition}
\label {prop:cob00twisted}
Let $\delta, \Lambda, \t, \ux, \ux'$ be as in Proposition~\ref{prop:cob00}, and $\SR$ a twisted coefficients ring as in \eqref{eq:sr}. Then:

$(a)$  There is a diagram
$$
\begin{CD}
\CFmd(Y_{\Lambda'}(L'), \ux'; \SR) @>{f_{3; \t, \SR}^\delta}>> \CFmd(Y_{\Lambda}(L), \ux; \SR)\\
@V{\cong}VV @VV{(f_{1,\ux; \SR}^\delta, H_{1, \ux; \SR}^\delta)}V \\
T_1^\ux \cdot \CFmd(Y_{\Lambda'}(L'), \ux'; \SR) @>>> Cone(f_{2, \ux; \SR}^\delta),
\end {CD}
$$
commuting up to chain homotopy, with the bottom horizontal arrow being inclusion into the mapping cone.

$(b)$  There is a diagram
$$
\begin{CD}
T_1^\ux \cdot \CFmd(Y_{\Lambda'}(L'), \ux'; \SR) @>>> Cone(f_{2, \ux; \SR}^\delta),\\
@V{\cong}VV @VV{(H_{2, \ux; \SR}^\delta, f_{3,\ux; \SR}^\delta)}V \\
\CFmd(Y_{\Lambda'}(L'), \ux'; \SR) @>{f_{3; \t; \SR}^\delta}>> \CFmd(Y_{\Lambda}(L), \ux; \SR),
\end {CD}
$$
commuting up to chain homotopy, with the top horizontal arrow being inclusion into the mapping cone.

\end {proposition}

\section {Proof of the surgery theorem for link-minimal complete systems}
\label {sec:proof}

The goal of this section is to prove Theorem~\ref{thm:surgery}. We first do so for the case of the basic systems from Definition~\ref{def:basic}. The proof in that case will be modeled on the proof of the formula for the Heegaard Floer homology of integral surgeries on knots, see \cite{IntSurg}. We will need to combine the arguments in \cite{IntSurg} with the homological algebra from \cite{BrDCov}. At the end we will explain how the statement of Theorem~\ref{thm:surgery} for basic systems implies the statement for all link-minimal complete systems of hyperboxes. 

\subsection {Large surgeries on links}
\label {sec:large}
Let $\orL \subset Y$ be a link in an integral homology three-sphere as in Section~\ref{sec:statement}. 
We let $\tilde \Lambda \gg 0$ be a sufficiently large framing on $L$, meaning that the framing coefficients $\tilde \lambda_i$ on each component are sufficiently large, as for the framing denoted $\tilde \Lambda$ in Sections~\ref{sec:another}. We let $\tilde \Lambda_i \in H_1(Y - L)$ be the induced framings on each component $L_i$, as usual.  Recall that $H(L, \tilde \Lambda)$ denotes the lattice in $H_1(Y-L) \cong \zz^\ell$ generated by all $\tilde \Lambda_i$. 

We use the notation from Section~\ref{sec:4c}. In particular, surgery on the framed link $L$ produces a cobordism $W_{\tilde \Lambda}(L)$ between $Y$ and the surgered manifold $Y_{\tilde \Lambda}(L)$. Since $\tilde \Lambda$ is chosen sufficiently large, the manifold $Y_{\tilde \Lambda}(L)$ is a rational homology three-sphere. Let $W'_{\tilde \Lambda}(L)$ be the cobordism between $Y_{\tilde \Lambda}(L)$ and $Y$, obtained by turning around the cobordism $-W_{\tilde \Lambda}(L)$.

As in the proof of Lemma~\ref{lem:HH}, choose a Seifert surface $F_i \subset Y$ for each link component $L_i$, and let $\hat F_i$ be the surface obtained by capping off $F_i$ in $W_{\tilde \Lambda}(L)$. By a slight abuse of notation, we also denote by $\hat F_i$ the corresponding surface in $W'_{\tilde \Lambda}(L)$. The homology classes $[\hat F_i], i=1, \dots, \ell$, form a basis of $H_2( W'_{\tilde \Lambda}(L))$. As in \eqref{eq:H^2}, we identify $H^2(W'_{\tilde \Lambda}(L))$ with $\zz^\ell$ by sending a cohomology class $c$ to $(\langle c, [\hat F_1]\rangle, \dots, \langle c, [\hat F_\ell]\rangle)$.

Given a $\spc$ structure $\ux$ over $Y_{\tilde \Lambda}(L)$, we can extend it to a $\spc$ structure $\tt$ over $W'_{\tilde \Lambda}(L)$. We can then find $\s \in \H(L)$ such that
$$ c_1(\tt) \equiv 2\s - (\tilde \Lambda_1 + \dots + \tilde \Lambda_\ell)  \mod 2H(L, \tilde \Lambda) . $$ 

The correspondence $\ux \mapsto \s$ determined by the above formula induces an isomorphism
$$ \spc(Y_{\tilde \Lambda}(L)) \longrightarrow \H(L)/H(L, {\tilde \Lambda}).$$

Let $P(\tilde \Lambda)$ be the intersection of the lattice $\H(L)$ with the hyper-parallelepiped with vertices
$$ \zeta + \frac{1}{2}(\pm \tilde \Lambda_1 \pm \tilde \Lambda_2 \pm \dots \pm \tilde \Lambda_\ell),$$
as in Section~\ref{sec:another}. This is a fundamental domain for $\H(L)/H(L, \tilde \Lambda)$, see Equation~\eqref{eq:plambda}. Hence, there is a bijection
$$ \spc(Y_{\tilde \Lambda}(L)) \cong P(\tilde \Lambda),$$
see Section~\ref{subsec:surgery}. From now on we will denote a $\spc$ structure on $Y_{\tilde \Lambda}(L)$ by the corresponding value $\s \in P(\tilde \Lambda) \subset \H(L)$. For $\s \in P(\tilde \Lambda)$, we denote by $\nx_\s$ the $\spc$ structure over $W'_{\tilde \Lambda}(L)$ satisfying
\begin {equation}
\label {eq:nx}
 c_1(\nx_\s) = 2\s - (\tilde \Lambda_1 + \dots + \tilde \Lambda_\ell),
 \end {equation}
 compare Equation~\eqref{eq:c1s}.

Let us choose a basic system $\Hyper$ for $\orL \subset Y$, as in Section~\ref{sec:basic}. Recall that the initial Heegaard diagram $\Hyper^L=(\Sigma, \alphas, \betas, \ws, \zs)$ in the system contains $\ell$ beta curves $\beta_1, \dots, \beta_\ell$ such that the basepoints $w_i, z_i$ lie one on each side of $\beta_i$.  Let $\Ring =\ff[[U_1, \dots, U_\ell]]$.

\begin {theorem}
\label {thm:groups}
For $\tilde \Lambda$ sufficiently large, there exist quasi-isomorphisms of relatively $\zz$-graded complexes of $\Ring$-modules
$$\Psi_{\tilde \Lambda, \s}^- :  \CFm(Y_{\tilde \Lambda}(L), \s) \longrightarrow  \Chain^-(\Hyper^L, \s),$$
for all $\s \in P(\tilde \Lambda)$.
\end {theorem}

\begin {proof}
For each $i = 1, \dots, \ell$, we construct a curve $\delta_i$ by twisting the longitude of $L_i$ $\tilde \lambda_i$ times along $\beta_i$,  in a symmetric way as in \cite[Figure 5]{Knots} (but with the twisting done in the opposite direction). Thus, $\delta_i$ specifies the framing $\tilde \Lambda_i$ of the component $L_i, i=1, \dots, \ell$. We complete this to a full set of attaching circles $\delta$ by taking curves $\delta_{\ell+1}, \dots, \delta_{\ell + g-1}$ that approximate (i.e. are small Hamiltonian translates of) $\beta_{\ell +1}, \dots, \beta_{\ell + g-1}$. The result is a triple Heegaard diagram $(\Sigma, \alphas, \deltas, \betas, \ws)$ for the cobordism $W'_{\tilde \Lambda}(L)$,  such that we have the three-manifolds $Y_{\alpha, \beta} \cong Y, \ Y_{\alpha, \delta} \cong Y_{\tilde \Lambda}(L), Y_{\delta, \beta} \cong \#^{g-1}(S^1 \times S^2)$. 

For $\s \in \H(L)$, we can define a chain map $\Psi_{\tilde \Lambda, \s}^- :  \CFm(Y_{\tilde \Lambda}(L), \s) \longrightarrow  \Chain^-(\Hyper^L, \s)$ by the formula
\begin {equation}
\label {eq:quasiiso}
  \Psi_{\tilde \Lambda, \s}^- (\x) = \sum_{\y \in \Ta \cap \Tb} \sum_{\substack{\phi \in \pi_2(\x, \Theta, \y), \ \mu(\phi) =0 \\ n_{w_i}(\phi) - n_{z_i}(\phi) = A_i(\y) - s_i, \forall i }} \# \M(\phi) \cdot \prod_{i=1}^\ell U_i^{\min(n_{z_i}(\phi), n_{w_i}(\phi))} \cdot  \y, 
  \end {equation}
where $\Theta \in \CFm(Y_{\delta, \beta})$ is the top degree generator in homology.

The proof that $\Psi_{\tilde \Lambda, \s}$ is a quasi-isomorphism for $\tilde \Lambda \gg 0$ and $\s \in P(\tilde \Lambda)$ then proceeds along the same lines as \cite[proof of Theorem 4.4]{Knots}; see also \cite[Theorem 2.3]{IntSurg}, \cite[Theorem 4.1]{RatSurg}. Roughly, the argument is as follows. There are $\ell$ winding regions on the surface $\Sigma$, that is, neighborhoods of the curves $\beta_i \ (i=1, \dots, \ell)$ in which the twisting of the corresponding curves $\delta_i$ takes place. A generator $\x \in \Ta \cap \T_{\delta}$ is said to be {\em supported in the winding regions} if it contains points in all the $\ell$ winding regions. If this is the case, the $i\th$ component of the $\spc$ structure $\s \in P(\tilde \Lambda)$ of $\x$ equals the depth of the respective point of $\x$ inside the $i\th$ winding region, up to the addition of a constant; see \cite[Equation (14)]{Knots}.  A $\spc$ structure $\s$ is said to be supported in the winding regions if all the generators $\x \in \s$ are supported in the winding regions. For such $\s$, the map $\Psi_{\tilde \Lambda, \s}$ is actually an isomorphism of chain complexes, because it is approximated (with respect to area filtrations) by a ``nearest point'' map which is a bijection. 

Note that we have some freedom in choosing the winding region. Indeed, by replacing each $\delta_i \ (i=1, \dots, \ell)$ with an isotopic curve $\delta_i'$ so that the number of twists to the left of $\beta_i$ is changed (i.e.  translating the twists to the right or left of the curve $\beta_i$), we obtain another strongly equivalent triple Heegaard diagram $(\Sigma, \alphas, \betas, \deltas', \ws, \zs)$. If $\delta'_i$ differs from $\delta_i$ by $k_i$ twists, the set of $\spc$ structures supported in the new winding regions is a translate of the old set by $(k_1, \dots, k_\ell)$. We have 
$$(\tilde \lambda_1 - C_1)\cdot (\tilde \lambda_2 - C_2) \cdots (\tilde \lambda_\ell - C_\ell)$$
possibilities for the position of the attaching set $\deltas$, where $C_1, \dots, C_\ell$ are constants (independent of the framing coefficients). For each of these $\deltas$, the number of $\spc$ structures not supported in the respective winding regions is of the order of
$$ \tilde \lambda_1 \cdots  \tilde \lambda_\ell \cdot \sum_{i=1}^\ell \frac{C_i'}{\tilde \lambda_i},$$ 
for some constants $C_i'$. 

Let us choose $2^\ell$ different sets of curves $\deltas^\eps, \ \eps \in \{0,1\}^\ell$, such that $\deltas^\eps_i$ and $\deltas^{\eps'}_i$ differ by $\tilde \lambda_i/2$ twists whenever $\eps_i \neq \eps'_i$ (and are the same curve if $\eps_i = \eps_i'$).  If $\tilde \lambda_1, \dots, \tilde \lambda_i $ are sufficiently large, we see that each $\spc$ structure $\s$ on $Y_{\tilde \Lambda}(L)$ is supported in one of the winding regions for some $\deltas^\eps$. The fact that the new map $\Psi^{\eps, -}_{\tilde \Lambda, \s}$ (using $\deltas^\eps$ instead of $\deltas$) is an isomorphism implies that the original map $\Psi^-_{\tilde \Lambda, \s}$ (which differs from $\Psi^{\eps, -}_{\tilde \Lambda, \s}$ by composition with chain homotopy equivalences) is a quasi-isomorphism. 
\end {proof}

Let $\orL^\ori$ be $L$ with some arbitrary orientation $\ori$, and let $\nx_\s^{\ori}$ be the  $\spc$ structure on $W'_{\tilde \Lambda}(L)$  satisfying
\begin {equation}
\label {eq:cori}
 c_1(\nx^{\ori}_\s) = c_1(\nx_\s) + 2\tilde \Lambda_{\orL, \orL^\ori},
 \end {equation}
where $\tilde \Lambda_{\orL, \orL^\ori}=\sum_{i \in I_-(\orL, \orL^\ori)} \tilde \Lambda_i$ and $I_-(\orL, \orL^\ori)$ is the set of indices describing components of $\orL$ oriented differently in $\orL^\ori$, see Section~\ref{sec:inclusions}.

We  denote by $F^-_{W, \nx}$ the map on Heegaard Floer complexes induced by a particular cobordism $W$ and $\spc$ structure $\nx$, see \cite{HolDiskFour}.

\begin {theorem}
\label {thm:maps}
Fix $\orL^\ori \subset Y$ as above. For any $\tilde \Lambda \gg 0$ and $\s \in P(\tilde \Lambda)$, there is a commutative diagram:
\begin {equation}
\label {eq:cd}
\begin {CD}
\CFm(Y_{\tilde \Lambda}(L), \s) @>{F^-_{W'_{\tilde \Lambda}(L), \nx^{\ori}_\s}}>> \CFm(Y) \\
@V{\Psi^-_{\tilde \Lambda, \s}}VV @VV{\cong}V  \\
 \Chain^-(\Hyper^L, \s) @>{\Pr^{\orL^\ori}_\s}>> \Chain^-(\Hyper^L, p^{\orL^\ori}(\s)),
\end {CD} 
\end {equation}
for all $\s \in P(\tilde \Lambda)$.
\end {theorem}

The proof of Theorem~\ref{thm:maps} is similar to that of \cite[Theorem 2.3]{IntSurg}. In the diagram~\eqref{eq:cd}, we implicitly identified $\Chain^-(\Hyper^L, p^{\orL^\ori}(\s))$ with $\Chain^-(r_{\orL^\ori}(\Hyper^L), \psi^{\orL^\ori}(\s))$,  see Equation~\eqref{eq:red_lmn}. Also, it is worth mentioning that in the proof of Theorem~\ref{thm:maps}, the map $F^-_{W'_{\tilde \Lambda}}(L)$ is defined using the triple Heegaard diagram $(\Sigma, \alphas, \deltas, \betas, \ws^\ori)$, where 
$$ \ws^\ori =  \{ w_i  |  i\in I_+(\orL, \orL^\ori) \} \cup \{ z_i  |  i\in I_-(\orL, \orL^\ori) \}.$$

(This is the set of basepoints that would be denoted $\ws^{\orL, \orL^\ori}$ in  the notation of Section~\ref{sec:basic}.)

There is a more refined version of Theorem~\ref{thm:maps}, as follows. Note that the cobordism $W'_{\tilde \Lambda}(L)$ consists of $\ell$ two-handle additions, which can be composed in any order. Different ways of composing are related by chain homotopies, forming a hypercube of chain complexes. In the refined version that we state below (Theorem~\ref{thm:hypermaps}) we replace the top arrow in \eqref{eq:cd} with this hypercube; this is isomorphic to another hypercube, replacing the bottom arrow in \eqref{eq:cd}. 

Let $\deltas = (\delta_1, \dots, \delta_{g+\ell -1})$ be a set of attaching curves as in the proof of Theorem \ref{thm:groups}. Given $\eps = (\eps_1, \dots, \eps_\ell) \in \E_\ell = \{0,1\}^\ell$, we define a new $(g+\ell-1)$-tuple of attaching circles $\etas^\eps$ by
$$ \eta^\eps_i \approx \begin {cases} 
\delta_i & \text{if } \eps_i = 0,\\
\beta_i & \text{if } \eps_i = 1.
\end {cases} $$

For each $\eps \in \E_\ell$, we denote by $\orL^{\ori,\eps} \subseteq \orL^\ori$ the oriented sublink consisting of those components $L_i$ such that $\eps_i = 1$, all taken with the orientation induced from $\ori$. We let $L^\eps$ be the underlying sublink. Note that the Heegaard diagram $(\Sigma, \alphas, \etas^\eps, \ws^\ori)$ represents the three-manifold $Y_{\tilde \Lambda |_{L - L^\eps}}(L-L^\eps)$. 

Further, for every $\eps < \eps'$, the Heegaard diagram $(\Sigma, \etas^\eps, \etas^{\eps'}, \ws^\ori)$ represents a connected sum of some copies of $S^1 \times S^2$. We can arrange so that the Floer homology $\iHF(\T_{\eta^\eps}, \T_{\eta^{\eps'}})$, in the maximal degree with nontrivial homology (and in the torsion $\spc$ structure), is represented by a unique intersection point, which we denote by $\Theta^\can_{\eps, \eps'}$. Set:
\begin {equation}
\label {eq:teps}
 \Theta_{\eps, \eps'} = \begin {cases} 
\Theta^\can_{\eps, \eps'} & \text{if } \|\eps' - \eps\| =1\\
0 & \text{otherwise.}
\end {cases} 
\end {equation}

For $\eps \leq \eps'$, let 
$$ W'_{\tilde \Lambda}(L-L^\eps, L-L^{\eps'}) \subseteq W'_{\tilde \Lambda}(L)$$
be the cobordism from $Y_{\tilde \Lambda |_{L - L^\eps}}(L-L^\eps)$ to $Y_{\tilde \Lambda |_{L-L^{ \eps'}}}(L-L^{\eps'})$  obtained by reversing the surgery on $L^{\eps' - \eps}$. When $\nx$ is a $\spc$ structure on $ W'_{\tilde \Lambda}(L)$, we keep the same notation $\nx$ for its restriction to 
$ W'_{\tilde \Lambda}(L-L^\eps, L-L^{\eps'})$.

Consider the polygon map
$$ F(\eps, \eps', \ws^\ori): \CFm(\Ta, \T_{\eta^\eps}, \ws^\ori) \to \CFm(\Ta, \T_{\eta^{\eps'}}, \ws^\ori),$$
$$  F(\eps, \eps', \ws^\ori)(\x) = \sum_{\eps = \eps^0 < \dots < \eps^p = \eps'} f(\x \otimes \Theta_{\eps^0, \eps^1} \otimes \dots \otimes \Theta_{\eps^{p-1}, \eps^p}),$$
in the notation of Section~\ref{sec:polygon}, used here for polygon maps between ordinary Floer chain complexes, as in \cite[Section 4.2]{BrDCov}. When $\eps = \eps'$, this is simply the differential $\del$. When $\|\eps' - \eps\|=1$, the map $ F(\eps, \eps', \ws^\ori)$ is a triangle map representing the cobordism $ W'_{\tilde \Lambda}(L-L^\eps, L-L^{\eps'})$ and, as such, it decomposes as a sum of maps according to the $\spc$ structures on that cobordism. For general $\eps < \eps'$, the map  $F(\eps, \eps', \ws^\ori)$ is a higher order chain homotopy relating the different ways of splitting  $ W'_{\tilde \Lambda}(L-L^\eps, L-L^{\eps'})$ into two-handle additions. It still decomposes as a sum of maps 
$$F(\eps, \eps', \ws^\ori, \nx)$$
according to the $\spc$ structures $\nx$ on the cobordism $ W'_{\tilde \Lambda}(L-L^\eps, L-L^{\eps'})$. 

\begin {theorem}
\label {thm:hypermaps}
Fix $ \tilde\Lambda \gg 0$, $\s \in P(\tilde \Lambda)$ and an orientation $\ori$ on $L$ as above. Then, the hypercube with chain groups $$C^\eps = \CFm(\Ta, \T_{\eta^\eps}, \ws^\ori, \psi^{\orL^{\ori,\eps}}(\s))=\CFm(Y_{\tilde \Lambda |_{L - L^\eps}}(L-L^\eps), \psi^{\orL^{\ori,\eps}}(\s))$$ and maps
$$ D^{\eps' -\eps}_\eps = F(\eps, \eps', \ws^\ori, \nx_\s^\ori)$$
is quasi-isomorphic to the hypercube with chain groups 
$$C^\eps = \Chain^-(\Hyper^L, p^{\orL^{\ori,\eps}}(\s)) = \Chain^-(r_{\orL^{\ori,\eps}}(\Hyper^L), \psi^{\orL^{\ori, \eps}}(\s))$$
and maps
$$ D^{\eps' -\eps}_\eps = \begin {cases}
\del & \text{if } \eps = \eps', \\
\Pr^{\orL^\ori_i}_{p^{\orL^{\ori,\eps}}(\s)} & \text{if } \|\eps' - \eps\|=1, \ \orL^{\ori,\eps'} = \orL^{\ori, \eps} \amalg \orL^\ori_i, \\
0 & \text{otherwise.}
\end{cases} $$
\end {theorem}

\begin{proof}
  The maps
  $$\Psi^-_{\tilde \Lambda|_{L-L^{ \eps}}, \psi^{\orL^{\ori,\eps}}(\s)}:
  \CFm(Y_{\tilde \Lambda |_{L - L^\eps}}(L-L^\eps),
  \psi^{\orL^{\ori,\eps}}(\s)) \to
  \Chain^-(r_{\orL^{\ori,\eps}}(\Hyper^L), \psi^{\orL^{\ori,
      \eps}}(\s))$$ given by Equation~\eqref{eq:quasiiso}
  can be generalized to give maps
  (increasing $\eps$) which count higher polygons, and
  where the $U$ powers are counted just as in \eqref{eq:quasiiso}. 
  These form a chain map between the two  hypercubes.
  
  By
  definition, a quasi-isomorphism of hypercubes means that the
  corresponding $\eps$-preserving maps are quasi-isomorphisms for all
  $\eps$. 
  Indeed, the maps $\Psi^-_{\tilde \Lambda|_{L-L^{ \eps}},
    \psi^{\orL^{\ori,\eps}}(\s)}$
  are quasi-isomorphisms by a simple extension of 
  Theorem~\ref{thm:groups}.
\end{proof}

We can change the two hypercubes in Theorem~\ref{thm:hypermaps} by chain homotopy equivalences, and arrive at the following:

\begin {proposition}
\label {prop:hypermaps2}
Fix $ \tilde\Lambda \gg 0$ and $\s \in P(\tilde \Lambda)$ and an orientation $\ori$ as above. Then, the hypercube $\hyp^\ori_{\s}$ with chain groups $$C^\eps = \CFm(\Ta, \T_{\eta^\eps}, \ws, \psi^{\orL^{\ori,\eps}}(\s))=\CFm(Y_{\tilde \Lambda |_{L - L^\eps}}(L-L^\eps), \psi^{\orL^{\ori,\eps}}(\s))$$ and maps
$$ D^{\eps' -\eps}_\eps = F(\eps, \eps', \ws, \nx_\s^\ori)$$
is quasi-isomorphic to the hypercube $\Hyper^\ori_{\s}$ having chain groups 
$$C^\eps = \Chain^-(\Hyper^{L- L^\eps}, \psi^{\orL^{\ori, \eps}}(\s))$$
and maps
$$ D^{\eps' -\eps}_\eps = \Phi^{\orL^{\ori,\eps'-\eps}}_{ \psi^{\orL^{\ori, \eps}}(\s)}= \De^{\orL^{\ori,\eps'-\eps}}_{(p^{\orL^{\ori,\eps'-\eps}} \circ \psi^{\orL^{\ori, \eps}})(\s)} \circ \Pr^{\orL^{\ori,\eps'-\eps}}_{\psi^{\orL^{\ori, \eps}}(\s)},$$
in the notation of Section~\ref{subsec:desublink}.
\end {proposition}

\begin {proof}
Let us compare the first hypercube in Theorem~\ref{thm:hypermaps} with the first hypercube in Proposition~\ref{prop:hypermaps2}. Note that the Heegaard diagrams $(\Sigma, \alphas, \etas^\eps, \ws^\ori)$ and $(\Sigma, \alphas, \etas^\eps, \ws)$ both represent the same three-manifold $Y_{\tilde \Lambda |_{L - L^\eps}}(L-L^\eps)$. Hence, the respective chain complexes are chain homotopy equivalent. In fact, we can describe the chain homotopy equivalence along the lines of Section~\ref{sec:basic}. First, note that $(\Sigma, \alphas, \etas^\eps, \ws^\ori)$ is isotopic to $(\Sigma, \alphas, \etas^{''\eps}, \ws)$, where the collection $\etas^{''\eps}$ differs from $\etas^{\eps}$ by replacing every curve that approximates $\beta_i$ with one that approximates $\beta''_i$. Here $\beta''_i$ is as in Section~\ref{sec:basic}, and similarly we recall that we also have an intermediate curve $\beta_i'$. One can relate $\CFm(\Ta, \T_{\eta^{''\eps}}, \ws)$ to $\CFm(\Ta, \T_{\eta^{'\eps}}, \ws)$ and then to $\CFm(\Ta, \T_{\eta^{\eps}}, \ws)$ via chain homotopy equivalences given by triangle maps with one vertex in the respective canonical generator. It is straightforward to lift these to chain homotopy equivalences between the respective hypercubes. Moreover, we can restrict everything to a $\spc$ structure $\psi^{\orL^{\ori,\eps}}(\s)$, which is left unchanged throughout.

Now let us compare the second hypercube in Theorem~\ref{thm:hypermaps} with the second hypercube in Proposition~\ref{prop:hypermaps2}. The former is reminiscent of a canonical hypercube as in Definition~\ref{def:canhyper}, but it has inclusions rather than identity maps along its edges. Nevertheless, we can construct a chain map very similar to the canonical inclusion from Section~\ref{sec:caninc}, as follows.

For $\eps=(\eps_1, \dots, \eps_\ell) \in \E_\ell$ and $i \in \{0,1, \dots, \ell\}$, we let $\eps[\leq i]$ resp. $\eps[>i]$ be the multi-indices obtained from $\eps$ by changing all entries indexed by $j > i$ (resp. $j \leq i$) into zeros; compare Section~\ref{sec:caninc}. We define an intermediate hypercube $H[i]$ to have chain groups
$$ C[i]^\eps = \Chain^-(\Hyper^{L - L^{\eps[\leq i]}}, (\psi^{\orL^{\ori, \eps[\leq i]}} \circ p^{\orL^{\ori, \eps}})(\s) ) $$
and maps
$$ D[i]_{\eps}^{\eps' - \eps} : \Chain^-(\Hyper^{L - L^{\eps[\leq i]}}, (\psi^{\orL^{\ori, \eps[\leq i]}} \circ p^{\orL^{\ori, \eps}})(\s) ) \to \Chain^-(\Hyper^{L - L^{\eps'[\leq i]}}, (\psi^{\orL^{\ori, \eps'[\leq i]}} \circ p^{\orL^{\ori, \eps'}})(\s) )$$
given by
$$ D[i]_{\eps}^{\eps'-\eps} = \begin {cases}
D^{\orL^{\ori,(\eps'-\eps)}}  \circ \Pr^{\orL^{\ori,(\eps'-\eps)}} &\text{if } \eps[>i] = \eps'[>i],\\
\Pr^{\orL^{\ori,(\eps'-\eps)}}  &\text{if } \eps[\leq i] = \eps'[\leq i] \text{ and } \|\eps'[>i] - \eps[>i] \| =1,\\
0 & \text{otherwise.}
\end {cases} $$
We omitted here the subscripts in the maps $D$ and $\Pr$, as they are uniquely determined by the domains of those maps.

Note that $H[0]$ is the  second hypercube in Theorem~\ref{thm:hypermaps}, while $H[\ell]$ is  the second hypercube in Proposition~\ref{prop:hypermaps2}.

For $i=1, \dots, \ell$, we define chain maps
$$ F[i]: H[i-1] \to H[i]$$
to consist of 
$$ F[i]_{\eps}^{\eps' - \eps} : \Chain^-(\Hyper^{L - L^{\eps[\leq (i-1)]}}, (\psi^{\orL^{\ori, \eps[\leq (i-1)]}} \circ p^{\orL^{\ori, \eps}})(\s) ) \to \Chain^-(\Hyper^{L - L^{\eps'[\leq i]}}, (\psi^{\orL^{\ori, \eps'[\leq i]}} \circ p^{\orL^{\ori, \eps'}})(\s) ),$$
$$ F[i]_{\eps}^{\eps' - \eps} =  \begin {cases}
D^{\orL^{\ori,\eps'[\leq i]-\eps[\leq (i-1)]}}  \circ \Pr^{\orL^{\ori,(\eps'-\eps)[\leq (i-1)]}}  &\text{if } \eps_i =1,\  \eps[>i] = \eps'[>i],\\
\Id &\text{if } \eps = \eps' \text{ and } \eps_i =0,\\
0 & \text{otherwise.}
\end {cases} $$

Note that when $\eps = \eps'$ the map $ F[i]_{\eps}^{\eps' - \eps} $ is either the identity (when $\eps_i = 0$) or an edge map of the form $ D^{\orL_i^\ori}$ (when $\eps_i =1$); in either case, it is a chain homotopy equivalence. One can lift the respective chain homotopies to the level of the hypercubes. This shows that each $F[i]$ is a chain homotopy equivalence of hypercubes. The composition $F[\ell] \circ \dots \circ F[1]$ then represents a chain homotopy equivalence between $H[0]$ and $H [\ell]$.

The claim now follows from Theorem~\ref{thm:hypermaps}.
\end {proof}

\subsection {Iterating the exact triangle} \label{sec:iterate}

Our goal in this section is to present a generalization of the exact triangle \eqref{eq:extri} in the form of a description of $\HFm(Y_{\Lambda}(L))$ for arbitrary surgery on a link $L$ in an integral homology sphere $Y$. This will be based on iterating the more general exact triangle from Proposition~\ref{prop:leseq}.
 
Let $\tilde \Lambda \gg 0$ be a new framing for $L$, as in Section~\ref{sec:large}. We denote by $\tilde \lambda_i$ and $\lambda_i$ the framing coefficients on the component $L_i$, coming from $\tilde \Lambda$ resp. $\Lambda$. For each $i=1, \dots, \ell$, set
$$ m_i = \tilde \lambda_i - \lambda_i \gg 0.$$

We keep all the notation from Section~\ref{sec:large}. In particular, we have chosen a basic system $\Hyper$ for $\orL \subset Y$, and we have collections of curves $\etas^\eps$ for each $\eps \in \E_\ell = \{0,1\}^\ell$.

We now define collections $\etas^\eps$ for all $\eps \in \{0,1,\infty\}^\ell$, such that when $\eps \in \{0,1\}^\ell$ the respective collection coincides with the one already defined. Let $\gamma_i$ be a simple closed curve in $\Sigma$ disjoint from the basepoints and the beta curves, which specifies the framing $\Lambda_i$ of the component $L_i, i=1, \dots, \ell$. (In other words, this is the analogue of $\delta_i$ when we use $\Lambda$ instead of $\tilde \Lambda$.) We complete this to a full set of attaching circles $\gammas$ by taking curves $\gamma_{\ell+1}, \dots, \gamma_{\ell + g-1}$ that approximate (i.e. are small Hamiltonian translates of) $\beta_{\ell +1}, \dots, \beta_{\ell _ g-1}$.  The Heegaard diagram $(\Sigma, \alphas, \gammas, \ws)$ then represents the three-manifold $Y_\Lambda(L)$. 

Given $\eps = (\eps_1, \dots, \eps_\ell) \in \{0,1,\infty\}^\ell$, we define the $(g+\ell-1)$-tuple of attaching circles $\etas^\eps$ by
$$ \eta^\eps_i \approx \begin {cases} 
\delta_i & \text{if } \eps_i = 0,\\
\beta_i & \text{if } \eps_i = 1,\\
\gamma_i & \text{if } \eps_i = \infty.
\end {cases} $$

For every $\eps < \eps'$, the Heegaard diagram $(\Sigma, \etas^\eps, \etas^{\eps'}, \ws)$ represents a connected sum of some copies of $S^1 \times S^2$ and lens spaces. As such, there is a canonical torsion $\spc$ structure on this manifold, see \cite[Definition 3.2]{IntSurg} and Section \ref{sec:LES}. We arrange so that the Floer homology $\HFm(\T_{\eta^\eps}, \T_{\eta^{\eps'}}, \ws)$ in that $\spc$ structure, in the maximal degree with nonzero homology, is represented by a unique intersection point. We denote that point by $\Theta^\can_{\eps, \eps'}$. We then define $\Theta_{\eps, \eps'}$ just as in \eqref{eq:teps}.

Consider the ring 
$$ \TR = \ff [T_1, \dots, T_\ell]/(T_1^{m_1} -1, \dots, T_\ell^{m_\ell} -1).$$

Construct the chain complex with twisted coefficients $\CFm(\Ta, \Tb, \ws; \TR)$, which as a module is $\CFm(\Ta, \Tb, \ws) \otimes_\ff \TR$, and comes equipped with the differential
\begin {equation}
\label {eq:ti}
 \del^- \x = \sum_{\y \in \Ta \cap \Tb} \sum_{\{\phi \in \pi_2(\x, \y) | \mu(\phi)=1\}} \# \M(\phi) \cdot \Bigl( \prod_{i=1}^\ell T_i^{n_{w_i}(\phi)-n_{z_i}(\phi)} U_i^{n_{w_i}(\phi)}\Bigr) \cdot \y. 
\end {equation}
 Since $Y$ is an integral homology sphere, all the periodic domains on the diagram $(\Sigma, \alphas, \betas, \ws)$ are multiples of $\Sigma$. As a consequence, there exists an isomorphism of chain complexes
$$ \CFm(\Ta, \Tb, \ws; \TR) \cong \CFm(\Ta, \Tb, \ws) \otimes_\ff \TR = \oplus^{m_1m_2 \dots m_\ell} \CFm(\Ta, \Tb, \ws),$$
compare \cite[Equation (7)]{IntSurg}. 

More generally, for $\eps \in \{0,1,\infty\}^\ell$, we let $\TR^\eps$ be the polynomial ring in variables $T_i$ for those $i$ such that $\eps_i = 1$, and with relations $T_i^{m_i} = 1$. We consider the chain complex with twisted coefficients 
$$ \CC^\eps = \CFm(\Ta, \T_{\eta^\eps}, \ws; \TR^\eps),$$ constructed as above, but in which we only keep track of the multiplicities $n_{w_i}-n_{z_i}$ (using $T_i$) for those $i$ with $\eps_i = 1$. Let 
$$m^\eps = \prod_{\{i | \eps_i = 1\}} m_i.$$

We then have an identification of chain complexes
$$ \CC^\eps = \CFm(\Ta, \T_{\eta^\eps}, \ws; \TR^\eps) \cong \CFm(\Ta, \T_{\eta^\eps}, \ws)\otimes_\ff  \TR^\eps = \oplus^{m^\eps} \CFm(\Ta, \T_{\eta^\eps}, \ws).$$

For  multi-indices $\eps \leq \eps'$ define linear maps
$$ \DD^{\eps' - \eps}_\eps: \CFm(\Ta, \T_{\eta^\eps}, \ws; \TR^\eps) \to \CFm(\Ta, \T_{\eta^{\eps'}}, \ws; \TR^{\eps'}),$$
$$ \DD^{\eps' - \eps}_\eps (\x) = \sum_p \sum_{\{ \eps = \eps^0 < \dots < \eps^p = \eps'\}} f(\x \otimes \Theta_{\eps^0, \eps^1} \otimes \dots \otimes \Theta_{\eps^{p-1}, \eps^p} ),$$
where $f$ is the polygon map as in Section~\ref{sec:polygon}, which keeps track of the difference in multiplicities at $w_i$ and $z_i$ according to $T_i$, just as in \eqref{eq:ti} above. (Compare Section~\ref{sec:hyperfloer}, \cite[Section 3]{IntSurg} and \cite[Equation (9)]{BrDCov}.) 

The direct sum 
$ \bigoplus_{\eps \in \{0,1,\infty\}^\ell} \CC^\eps $
forms a chain complex with differential $\DD = \sum \DD^{\eps' - \eps}_\eps$. The proof that $\DD^2=0$ is similar to that of \cite[Proposition 4.4]{BrDCov}, with the difference that here we use twisted coefficients as in the proof that $f_2^+ \circ f_1^+ \simeq 0$ in \cite[Section 3]{IntSurg}. Note that $\DD$ is composed of ordinary twisted coefficients maps (and higher homotopies) between Floer complexes. In Section~\ref{sec:LES} (and, in particular, in Proposition~\ref{prop:leseq}) we needed a more complicated map denoted $f_3^-$. In our current context, $f_3^-$ would go from $\CC^\eps$ to $\CC^{\eps'}$ where in $\eps'$ we changed an $\infty$ coordinate of $\eps$ to $0$. Maps of this type do not appear in $\DD$.

Let us consider the quotient complex corresponding to $\eps \in \E_\ell = \{0,1\}^\ell \subset \{0, 1, \infty\}^\ell$. The projection of the differential (which we still denote by $\DD$) turns this quotient complex into a hypercube of chain complexes, which we denote by $\hyp$. 

\begin {proposition}
\label {prop:hyp1}
The complex $\CFm(Y_{\Lambda}(L)) = \CFm(\Ta, \Tb, \ws) = \CC^{(\infty, \dots, \infty)}$ is quasi-isomorphic to the total complex of the hypercube $ \hyp = (\CC^\eps, \DD^\eps )_{\eps \in \E_\ell}$.
\end {proposition}

\begin {proof}
Iterate the quasi-isomorphism from Proposition~\ref{prop:leseq} along the same lines as in the proof of \cite[Theorem 4.1]{BrDCov}. Observe that in order to prove exactness in Proposition~\ref{prop:leseq} we had to use the more complicated map $f_3^-$. However, having established in Proposition~\ref{prop:leseq} that there are quasi-isomorphisms along the edges of our hypercube, we automatically obtain quasi-isomorphisms under iteration. Thus, there is no need to define higher homotopy versions of the map $f_3^-$.
\end {proof}

\begin{remark}
For now, we just regard Proposition~\ref{prop:hyp1} as a statement about ungraded complexes. The same goes for all the quasi-isomorphisms discussed below, until we deal to the question of relative gradings in Section~\ref{sec:rg}.
\end {remark}

We can give an alternate description of the hypercube $\hyp$ as follows. For $\eps \in \E_\ell$, note that 
$$\CC^\eps \cong \CFm(\Ta, \T_{\eta^\eps}, \ws)\otimes_\ff  \TR^\eps = \CFm(Y_{\tilde \Lambda|_{L-L^\eps}}(L-L^\eps)) \otimes_\ff \TR^\eps,$$  
where $L^\eps \subset L$ is the sublink consisting of those components $L_i$ such that $\eps_i = 1$.  Since $\tilde \Lambda \gg 0, $ the manifold $Y_{\tilde \Lambda|_{L-L^\eps}}(L-L^\eps)$ is a rational homology three-sphere.

Consider now $\eps, \eps' \in \E_\ell$ with $\eps \leq \eps'$. Suppose $L^{\eps'} = L^{\eps} \cup L_{i_1} \cup \dots \cup L_{i_p}$ for $p=\|\eps' - \eps\|$. If we set all the variables $T_i$ equal to $1$, the map $\DD_\eps^{\eps' - \eps}$ would simply be (several copies of) the polygon map $F(\eps, \eps', \ws)$ from Section~\ref{sec:large}.  When we keep the variables $T_i$ as they are, they keep track of the $\spc$ structures on the cobordism $W'_{\tilde \Lambda}(L-L^\eps, L-L^{\eps'})$, which is a two-handle attachment relating the manifolds $Y_{\tilde \Lambda|_{L-L^\eps}}(L-L^\eps)$ and $Y_{\tilde \Lambda|_{L-L^{\eps'}}}(L-L^{\eps'})$. Indeed, by the argument used for \cite[Equation (8)]{IntSurg}, there is an identification 
\begin {equation}
\label {eq:kks}
\spc(W'_{\tilde \Lambda}(L-L^\eps, L-L^{\eps'})) \xrightarrow{\cong} \zz^p, \ \ \ \kk \to (k_1, \dots, k_p)
\end {equation}
such that
 \begin {equation}
 \label {eq:lee}
  \DD_\eps^{\eps' - \eps} = \sum_{\kk \in \zz^p} T^\kk \cdot  F(\eps, \eps', \ws,  \kk), 
  \end {equation}
where we formally wrote
 \begin {equation}
  T^\kk = \prod_{j=1}^p T_{i_j}^{k_j}.
  \end {equation} 

\subsection {Modifying the hypercube}

Fix an integer $\delta > 0$. As in \cite{IntSurg}, it is helpful to replace the Floer complexes $\CFm$ by the corresponding complexes $\CFmd$, see Section~\ref{sec:algebra}. Proposition~\ref{prop:hyp1} has the following immediate consequence (which could also be obtained by iterating Proposition~\ref{prop:leseqd}):

\begin {proposition}
\label {prop:hyp1d}
The complex $\CFmd(Y_\Lambda(L)) = \CFmd(\Ta, \T_{\gamma}, \ws)$ is quasi-isomorphic to the total complex of the hypercube
$$ \hyp^\delta = \bigl(\CC^{\eps, \delta} = \CFmd(\Ta, \T_{\eta^\eps}, \ws; \TR^\eps), \DD^{\eps, \delta} \bigr)_{\eps \in \E_\ell}.$$
\end {proposition}

Here is the vertical truncation of Equation~\eqref{eq:lee}:

\begin {equation}
 \label {eq:leed}
  \DD_\eps^{\eps' - \eps, \delta} = \sum_{\kk \in \zz^p}  T^\kk \cdot  F(\eps, \eps', \ws,  \kk)^\delta. 
  \end {equation}

Let $L^{\eps'} = L^{\eps} \amalg M$, with $M= L^{\eps' - \eps} =L_{i_1} \cup \dots \cup L_{i_p}$. Note that among the $\spc $ structures on the cobordism $W'_{\tilde \Lambda}(L-L^\eps, L-L^{\eps'})$ there are some special ones, namely those of the form $$\nx^{\ori}_\s= \nx^{\ori}_\s|_{W'_{\tilde \Lambda}(L-L^\eps, L-L^{\eps'})}$$ (cf. Section~\ref{sec:large}),  for all possible orientations $\ori$ of $L$ and $\spc$ structures $\s$ on $Y_{\tilde \Lambda}(L)$. Note that $\nx^{\ori}_\s$, viewed as a $\spc$ structure on $W'_{\tilde \Lambda}(L-L^\eps, L-L^{\eps'})$, only depends on the restriction of $\ori$ to $M$. Indeed, from Equation~ \eqref{eq:cori} we deduce
$$ c_1(\nx^{\ori}_\s|_{W'_{\tilde \Lambda}(L-L^\eps, L-L^{\eps'})}) = c_1(\nx_\s|_{W'_{\tilde \Lambda}(L-L^\eps, L-L^{\eps'})}) + 2(\tilde\Lambda|_M)_{\orL, \orM^\ori}.$$

To simplify notation, for $\orM \in \Omega(M)$ we denote by $\nx^{\orM}_\s$ the restriction to $ W'_{\tilde \Lambda}(L-L^\eps, L-L^{\eps'})$ of any $\spc$ structure $\nx^{\ori}_\s$ such that the restriction $\orM^\ori$ of $\ori$ to $M$ is $\orM$. Then, with respect to the identification $W'_{\tilde \Lambda}(L-L^\eps, L-L^{\eps'}) \cong \zz^p$ from \eqref{eq:kks}, we have
\begin {equation}
\label {eq:norms}
\nx^{\orM}_\s = \nx_\s + (\tilde \Lambda|_M)_{\orL, \orM}.
\end {equation}

The advantage of using truncated maps comes from the following:
\begin {lemma}
\label {lemma:grshift}
Fix $\delta > 0$. Then, for sufficiently large $\tilde \Lambda \gg 0$ (compared to $\delta$), and for any $\eps, \eps' \in \E_\ell$ with $ \eps < \eps'$, we have 
$$ F(\eps, \eps', \ws,  \kk)^\delta= 0$$ whenever the $\spc$ structure $\kk$ on  $W'_{\tilde \Lambda}(L-L^\eps, L-L^{\eps'})$ is not of the form $\nx^{\orM}_\s$ for any  $\s \in \spc(Y_{\tilde \Lambda}(L))=P(\tilde \Lambda)$ and $\orM \in \Omega(L^{\eps' - \eps})$. Moreover, there is a constant $b^\delta > 0$ such that 
$$F(\eps, \eps', \ws,  \nx^\ori_\s)^\delta= 0$$
whenever $\s = (s_1, \dots, s_\ell) \in P(\tilde \Lambda) \subset \rr^\ell$ admits some $i \in \{1, \dots, \ell\}$ with the
 property that either:
\begin {itemize}
\item $s_i > b^\delta$ and $i \in I_-(\orL, \orM)$, or 
\item $s_i < -b^\delta$ and  $i \in I_+(\orL, \orM)$.
\end {itemize}
\end {lemma}

\begin {proof}
Since both $Y_{\tilde \Lambda|_{L-L^\eps}}(L-L^\eps)$ and $Y_{\tilde \Lambda|_{L-L^{\eps'}}}(L-L^{\eps'})$ are rational homology three-spheres, the respective Floer complexes $\CC^{ \eps, \delta}$ and $\CC^{\eps', \delta}$ admit absolute $\qq$-gradings, see \cite{HolDiskFour}. Because these complexes are vertically truncated, their absolute gradings lie in a finite range. Each map $F(\eps, \eps', \ws,  \kk)^\delta$ shifts the grading by a definite amount. The vanishing of the claimed maps happens because the respective shifts take the range for the initial complex to outside the grading range for the final complex. We leave the verification of the details to the interested reader. In the case $p=1$, this was done in \cite[Lemma 4.4 and Section 4.3]{IntSurg}.
\end {proof}

Observe that the restriction of the $\spc$ structure $\nx^{\orM}_\s$ to the cobordism  $W'_{\tilde \Lambda}(L-L^\eps, L-L^{\eps'})$ only depends on the value $\psi^{L^\eps}(\s) \in P(\tilde \Lambda|_{L-L^\eps})$. (As we recall from Remark~\ref{rem:mmm}, the hyper-parallelepipeds $P(\tilde \Lambda)$ behave well with respect to restriction to sublinks.) Thus, we can write $\nx^{\orM}_{\bar \s}$ for $\nx^{\orM}_\s$ whenever $\bar \s = \psi^{L^\eps}(\s) \in P(\tilde \Lambda|_{L-L^\eps})$.

By Equation~\eqref{eq:leed} and Lemma~\ref{lemma:grshift}, for $\tilde \Lambda \gg 0$ we have
$$ \DD_\eps^{\eps' - \eps, \delta} = \sum_{\orM \in \Omega(L^{\eps' - \eps})} \sum_{\s \in P(\tilde \Lambda|_{L-L^\eps})} T^{ \nx^{\orM}_\s} \cdot  F(\eps, \eps', \ws,  \nx^{\orM}_\s)^\delta. $$

By adjusting the identifications \eqref{eq:kks} if necessary, and taking Equation~\eqref{eq:norms} into account, we can simply write:
$$ \DD_\eps^{\eps' - \eps, \delta} = \sum_{\orM \in \Omega(L^{\eps'-\eps})} \sum_{\s \in P(\tilde \Lambda|_{L-L^\eps})} T^{ (\tilde \Lambda|_{L^{\eps' - \eps}})_{\orL, \orM}} \cdot  F(\eps, \eps', \ws,  \nx^{\orM}_\s)^\delta.$$

We have $T^{\tilde \Lambda_i - \Lambda_i} = 1$ for all $i$, by the definition of the $T_i$ variables. Hence,
$$   T^{ (\tilde \Lambda|_M)_{\orL, \orM}} =  T^{ (\Lambda|_M)_{\orL, \orM}}.$$

To summarize, the total complex $\CC^\delta$ of the hypercube $\hyp^\delta$ is
\begin {equation}
\label {eq:C1}
 \CC^\delta = \bigoplus_{\eps \in \E_\ell} \bigoplus_{\s \in P(\tilde \Lambda|_{L^\eps})} \CFmd(\Ta, \T_{\eta^\eps}, \ws, \s) \otimes \TR^\eps
 \end {equation}
with the differential on each summand being 
\begin {equation}
\label {eq:D1}
\DD^\delta = \sum_{\eps' \geq \eps} \sum_{\orM \in \Omega(L^{\eps' - \eps})} T^{ ( \Lambda|_{L^{\eps' - \eps}})_{\orL, \orM}} \cdot  F(\eps, \eps', \ws,  \nx^{\orM}_\s)^\delta.
 \end {equation}

On the other hand, in Section~\ref{sec:twt} we constructed a folded truncated complex 
$\C^{-, \delta}(\Hyper, \Lambda)\twist$ quasi-isomorphic to $\C^{-, \delta}(\Hyper, \Lambda)$. By rephrasing its description from \eqref{eq:cdfolded}, \eqref{eq:Dfolded}, we have
\begin {equation}
\label {eq:C2}
\C^{-, \delta}(\Hyper, \Lambda)\twist = \bigoplus_{\eps \in \E_\ell} \bigoplus_{\s \in P(\tilde \Lambda|_{L-L^\eps})} \Chain^{-, \delta} (\Hyper^{L-L^\eps}, \s) \otimes \TR^\eps
 \end {equation}
with the differential
\begin {equation}
\label {eq:D2}
 \D^{-, \delta} = \sum_{\eps' \geq \eps} \sum_{\orM \in \Omega(L^{\eps' - \eps})} T^{ ( \Lambda|_{L^{\eps' - \eps}})_{\orL, \orM}} \cdot  \Phi_\s^{L^{\eps' - \eps}, \delta}.
  \end {equation}

Looking at \eqref{eq:D1}, we observe that the maps  $F(\eps, \eps', \ws,  \nx^{\orM}_\s)^\delta$ also appear in the vertical truncations (by $\delta$) of the hypercubes $\hyp^\ori_{\s}$ considered in the statement of Proposition~\ref{prop:hypermaps2}. 
Similarly, looking at \eqref{eq:D2}, we observe that the maps $ \Phi_\s^{L^{\eps' - \eps}, \delta}$  appear  in the vertical truncation of the other hypercube $\Hyper^\ori$ considered in the statement of Proposition~\ref{prop:hypermaps2}. We seek to apply the result of Proposition~\ref{prop:hypermaps2} to obtain a quasi-isomorphism between the hypercubes $\hyp^\delta$ and $\Hyper^\delta$. Before doing so, however, we need a basic result from homological algebra:

\begin{lemma}
\label {lemma:mcones}
Let $A, B, A', B'$ be hypercubes of chain complexes (of the same dimension $d$, with differentials denoted by $\del$), and 
$$F_1, F_2: A \to B, \  \ \ F_1', F_2': A' \to B'$$  be chain maps.  Suppose that, for $i=1,2$, the corresponding $(d+1)$-dimensional hypercubes for the maps $F_i$ and $F_i'$ (that is, the mapping cones $A \xrightarrow{F_i} B$ and $A' \xrightarrow{F_i} B'$) can be related by a quasi-isomorphism consisting of maps $(\Phi_i, \Psi_i, h_i)$ as in the diagram
$$
\xymatrix{
A\ar[d]_{F_i} \ar[r]^{\Phi_i} \ar[dr]^{h_i} & A' \ar[d]^{F_i'} \\
B\ar[r]_{\Psi_i}  & B'
 }
$$
Assume that $\Phi_1$ and $\Psi_1$ are quasi-isomorphisms. Further, suppose that the maps $\Phi_1$ and $\Phi_2$ are chain homotopic, and so are $\Psi_1$ and $\Psi_2$. Then the mapping cones $A \xrightarrow{F_1 + F_2} B$ and $A' \xrightarrow{F_1' + F_2'} B'$, viewed as $(d+1)$-dimensional hypercubes, are quasi-isomorphic as well.
\end {lemma}

\begin {proof}
By hypothesis, we have
$$ F'_1 \Phi_1 - \Psi_1F_1 = \del h_1 + h_1 \del, \ \  F'_2 \Phi_2  - \Psi_2F_2 = \del h_2 + h_2 \del,$$
and there are homotopies $\phi: A \to A', \psi: B \to B'$ such that
$$ \Phi_1 - \Phi_2 = \del \phi + \phi \del, \ \  \Psi_1 - \Psi_2 = \del\psi + \psi \del.$$
It follows that
$$
 (F'_1 + F_2')\Phi_1 - \Psi_1(F_1 + F_2) = \del g + g \del,$$
 where
 $$ g = h_1 + h_2 + F_2'\phi + \psi F_2.$$

Thus, $(\Phi_1, \Psi_1, g)$ form a chain map between the mapping cones $A \xrightarrow{F_1 + F_2} B$ and $A' \xrightarrow{F_1' + F_2'} B'$. This map is a quasi-isomorphism because $\Phi_1$ and $\Psi_1$ are so, compare Definition~\ref{def:quasiiso}. \end {proof}

\begin {proposition}
\label {prop:weaksurgery}
Given a basic system $\Hyper$ for a link $\orL$ in an integral homology sphere $Y$, there is an isomorphism:
$$ H_*(\C^{-, \delta}(\Hyper, \Lambda), \D^{-, \delta}) \cong \HFmd_*(Y_\Lambda(L)).$$ 
\end {proposition}

\begin {proof}
Let us denote by $\hyp^{\ori, \delta}_{\s}, \Hyper^{\ori, \delta}_{\s}$ the vertical truncations by $\delta$ of the hypercubes $\hyp^\ori_{\s}, \Hyper^\ori_{\s}$ considered in Proposition~\ref{prop:hypermaps2}. The result of that proposition implies (after vertical truncation) that, for any orientation $\ori$, the hypercubes $\hyp^{\ori, \delta}_{\s}$ and $\Hyper^{\ori, \delta}_{\s}$ are related by a quasi-isomorphism. Further, if we change the orientation $\ori$, it is easy to see that the respective quasi-isomorphisms are chain homotopic.

The hypercube  $\Hyper^{\delta}= (\C^{-, \delta}(\Hyper, \Lambda) \twist, \D^{-, \delta})$ described in \eqref{eq:C2}, \eqref{eq:D2}, is basically obtained from the hypercubes $\prod_{\s} \Hyper^{\ori, \delta}_\s$ for all possible choices of orientations $\ori \in \Omega(L)$, by gluing those hypercubes along their common parts (corresponding to sublinks). More precisely, this gluing process is an iteration of the one that made an appearance in Lemma~\ref{lemma:mcones}: given two mapping cone hypercubes $A \xrightarrow{F_1} B$ and $A  \xrightarrow{F_2} B$, we replace them by $A \xrightarrow{F_1 + F_2} B$. (See Figure~\ref{fig:glued} for an illustration of the gluing procedure in the case $\ell =2$.) Observe also that the hypercube  $\hyp^\delta$ described in \eqref{eq:C1}, \eqref{eq:D1} is obtained by from the hypercubes $\prod_{\s} \hyp^{\ori, \delta}_\s$ by gluing them using the same process.

\begin{figure}
\begin{center}
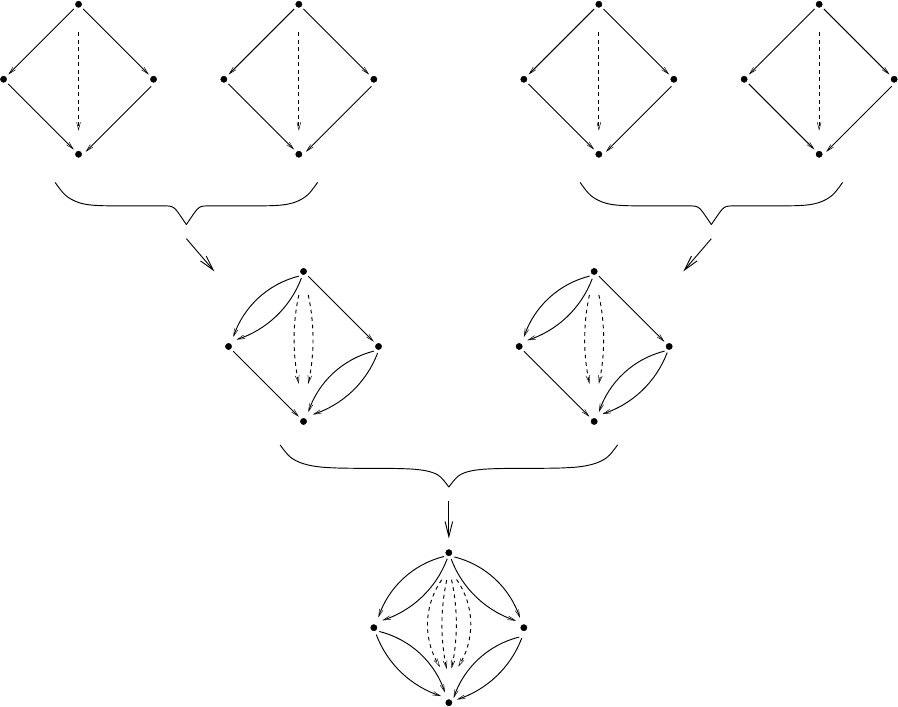
\end{center}
\caption {{\bf Gluing hypercubes in the proof of Proposition~\ref{prop:weaksurgery}.}
Given a two-component link $\orL = L_1 \cup L_2$, each of the four squares in the top row represent one of the hypercubes $\prod_{\s} \Hyper^\ori_{\s}$ (with the factors $\Hyper^\ori_{\s}$ as defined in Proposition~\ref{prop:hypermaps2}) for the four possible orientations $\ori$ of $L$. Combining these hypercubes (by adding up the respective differentials, as shown) yields the hypercube $\Hyper$ at the bottom, which is the one appearing in the statement of the Surgery Theorem \ref{thm:surgery}. The dashed lines represent chain homotopies such as $\Phi^{L_1 \cup L_2}, \Phi^{(-L_1) \cup L_2}$, etc. 
}
\label{fig:glued}
\end{figure}

Applying Lemma~\ref{lemma:mcones} repeatedly,  we can construct a quasi-isomorphism between $(\CC^\delta, \DD^\delta)$ and $(\C^{-, \delta}(\Hyper, \Lambda)\twist, \D^{-, \delta})$. The former complex is quasi-isomorphic to $\CFmd(Y_\Lambda(L))$ by Proposition~\ref{prop:hyp1d}, and the latter to $\C^{-, \delta}(\Hyper, \Lambda)$ by Proposition~\ref{prop:foldtruncate}.
\end {proof}

\subsection {$\spc$ structures} 
\label {sec:spc}
Recall from Section~\ref{subsec:surgery} that the complexes $\C^-(\Hyper, \Lambda)$ and $\CFm(Y_\Lambda(L))$ both break into direct sums of complexes $\C^-(\Hyper, \Lambda, \ux)$ and $\CFm(Y_\Lambda(L), \ux)$ according to the set $\spc(Y_\Lambda(L)) \cong \H(L)/H(L, \Lambda)$. Of course, the same is true for their vertical truncations by $\delta$. We would like to prove that the isomorphism in Proposition~\ref{prop:weaksurgery} preserves these decompositions.

As a warm-up exercise toward this goal, let us investigate to what extent we can make the total complex $\CC^\delta$ of  the hypercube $\hyp^\delta$ break into a direct sum according to $\spc$ structures $\ux \in \H(L)/H(L, \Lambda)$. 
 
We claim that this splitting can be realized when the lattice $H(L, \Lambda)$ is nondegenerate, i.e. the vectors $\Lambda_i$ are linearly independent over $\qq$ or, equivalently, $Y_\Lambda(L)$ is a rational homology three-sphere. Indeed, from \eqref{eq:lee} we see that the complex $\CC^\delta$ breaks into a direct sum according to equivalence classes of $\spc$ structures on $Y_{\tilde \Lambda}(L)$, where two structures in $\spc( Y_{\tilde \Lambda}(L)) \cong \H(L)/H(L, \tilde \Lambda)$ are equivalent if they differ by an element in the $\zz$-span of the vectors $m_i\tau_i = \tilde \Lambda_i - \Lambda_i, i=1, \dots, \ell$. In other words, the decomposition is according to $\H(L)/H(L,\Lambda,  \tilde \Lambda)$, where $H(L, \Lambda, \tilde \Lambda) \subseteq \zz^\ell$ is the lattice generated by all $\Lambda_i$ and $\tilde \Lambda_i$. Since $\Lambda$ is nondegenerate, we can arrange so that $H(L,\Lambda,  \tilde \Lambda) = H(L,\Lambda)$ by simply choosing $m_i \gg 0$ with $m_i \tau_i \in H(L, \Lambda)$, for all $i$. With this choice of $\tilde \Lambda$, our claim about the splitting of $\CC^\delta$ holds true.

When $H(L, \Lambda)$ is degenerate, we need to refine this approach, and settle for a splitting of a complex quasi-isomorphic to $\CC^\delta$, rather than one of $\CC^\delta$ itself. Indeed, since the vectors $\Lambda_i$ do not span $\qq^\ell$, we cannot always find $m_i \in \zz$ such that $m_i \tau_i \in H(L, \Lambda)$. Nevertheless, let us first choose some arbitrary $m_i' \gg 0$, such that the resulting framing $\tilde \Lambda'$ satisfies the conditions in Lemma~\ref{lemma:grshift}. Next, suppose the vectors $\Lambda_i$ span a subspace $\V \subset \qq^\ell$ of dimension $\ell-p$, with $p>0$. Choose $p$ coordinate vectors $\tau_i$ such that their span is complementary to $\V$. Without loss of generality, assume those coordinate vectors are $\tau_i, i=1, \dots, p$. Choose $m_i \gg m_i'$ arbitrarily for $i=1, \dots, p$, and let $H(L, \Lambda)^+ \subset \zz^\ell$ be the lattice spanned  by $\Lambda_i$'s together with $m_1\tau_1, \dots, m_{p}\tau_{p}$. Now choose $m_i \gg m_i'$ for $i=p+1, \dots, \ell$ such that $m_i \tau_i \in H(L, \Lambda)^+$. The result is a framing $\tilde \Lambda \gg \tilde \Lambda'$.

Lemma~\ref{lemma:grshift} says that the complex $\CC^\delta$ (constructed with respect to the framing $\tilde \Lambda)$ is similar in structure to the complex $\C^\delta=\C^{-, \delta}(\Hyper, \Lambda)$, compare Lemma~\ref{lemma:bdelta}. Hence, we can apply the combined truncation procedure from Section~\ref{sec:twt} to $\CC^\delta$, using the smaller framing $\tilde \Lambda'$. The result is a complex $\CC^\delta \langle \langle  \tilde \Lambda' \rangle \rangle$ quasi-isomorphic to $\CC^\delta$. 

The advantage of the complex $\CC^\delta \langle \langle  \tilde \Lambda' \rangle \rangle$ is that it splits as a direct sum according to $\H(L)/H(L, \Lambda)$. Indeed, we know it splits according to $\H(L)/H(L, \Lambda, \tilde \Lambda)$ just like $\CC^\delta$. Moreover, it is supported roughly on the hyper-parallelepiped $P(\tilde \Lambda')$, which is very small compared to $P(\tilde \Lambda)$. The key observation is that, with our choice of the values $m_i$, if two elements in $P(\tilde \Lambda')$ differ by an element in $H(L, \Lambda, \tilde \Lambda)$, they differ by an element in $H(L, \Lambda)$. This is true because  $H(L, \Lambda, \tilde \Lambda)$ is contained in the union of parallel subspaces
\begin {equation}
\label {eq:union}
 \bigcup_{t_1, \dots, t_p \in \zz} (\V + t_1m_1\tau_1 + \dots + t_pm_p\tau_p).
 \end {equation}

Set $\Delta(P(\tilde \Lambda')) = \{\s-\s'|\s,\s' \in P(\tilde \Lambda')\}$. If $\Delta(P(\tilde \Lambda'))$ is sufficiently small compared to the values $m_i$, the only one of the parallel subspaces in the union \eqref{eq:union} that intersects it nontrivially is $\V$ itself. This implies our claim about the decomposition of $\CC^\delta \langle \langle  \tilde \Lambda' \rangle \rangle$ according to $\H(L)/H(L, \Lambda)$.

Observe that similar remarks are applicable to the folded truncated complex $\C^{-, \delta}(\Hyper, \Lambda)\twist$. When $H(L, \Lambda)$ is nondegenerate, it splits according to $\spc$ structures on $Y_{\tilde \Lambda}(L)$, provided the values $m_i$ are chosen suitably. When $H(L, \Lambda)$ is degenerate, we can replace it by a quasi-isomorphic complex 
$$\C^{-, \delta}(\Hyper, \Lambda)\twist \langle \langle  \tilde \Lambda' \rangle \rangle = \C^{-, \delta}(\Hyper, \Lambda) \langle \langle  \tilde \Lambda' \rangle \rangle,$$
which again has the desired splitting.

\subsection {Relative gradings}
\label {sec:rg}

As explained in Sections~\ref{subsec:surgery} and \ref{sec:gradings}, for every $\ux \in \spc(Y_{\Lambda}(L))$, both complexes
$$ \C^{-, \delta}(\Hyper, \Lambda, \ux) \text{ and } \CFmd(Y_\Lambda(L), \ux)$$
admit relative $\zz/\delt(\ux)\zz$ gradings, where $\delt(\ux) \in \zz$ depends on $\ux$. Our goal is to prove the following strengthening of Proposition~\ref{prop:weaksurgery}:

\begin {proposition}
\label {prop:relgr}
Let $\Hyper$ be a basic complete system for an oriented link $\orL$ inside an integral homology sphere $Y$. We then have $\zz/\delt(\ux)\zz$-grading preserving isomorphisms 
\begin {equation}
\label {eq:star'}
H_*(\C^{-, \delta}(\Hyper, \Lambda, \ux), \D^{-, \delta}) \cong \HFmd_*(Y_\Lambda(L), \ux),
\end {equation}
and
\begin {equation}
\label {eq:star2'}
H_*^{\delta \from \delta'}(\C^{-}(\Hyper, \Lambda, \ux)) \cong \bHF^{-, \delta \from \delta'}_*(Y_\Lambda(L), \ux),
\end {equation}
for all $\delta' \geq \delta \geq 0$ and $\ux \in \spc( Y_{\tilde \Lambda}(L))$.
\end {proposition}

\begin {proof}
Recall that the quasi-isomorphism in Proposition~\ref{prop:hyp1d} can be obtained by iterating Proposition~\ref{prop:leseqd}. Indeed, let $\bar \Lambda$ be the framing on $L$ obtained from $\Lambda$ by adding $m_1$ to the coefficient of the first component, as in Section~\ref{sec:LES}. Also, we let $L' = L - L_1$ and denote by $\Lambda'$ the restriction of the framing $\Lambda$ to $L'$,  Consider the ring $\TR_1 = \ff[T_1]/(T_1^{m_1} - 1)$. The iteration process in the proof of Proposition~\ref{prop:hyp1d} starts by applying Proposition~\ref{prop:leseqd} to get that 
$\CFmd(Y_{\Lambda}(L))= \C^{(\infty, \infty, \dots, \infty), \delta}$ (in the notation of Section~\ref{sec:iterate}) is quasi-isomorphic to the mapping cone complex
\begin {equation}
\label {eq:mapcone2}
\CFmd(Y_{\bar \Lambda}(L)) \xrightarrow{f_2^\delta} \CFmd(Y_{\Lambda'}(L'); \TR_1),
\end {equation}
where the left hand side is the Floer complex $\C^{(0, \infty, \dots, \infty), \delta}$, the right hand side is the Floer complex $\C^{(1, \infty, \dots, \infty), \delta}$ with twisted coefficients, and the map $f_2^\delta$ is the triangle-counting map $\D^{(1, 0, \dots, 0), \delta}_{(0, \infty, \dots, \infty)}$ from \eqref{eq:leed}.

The next step in the iteration process will be to show that each of the two sides in \eqref{eq:mapcone2} is itself quasi-isomorphic to a mapping cone (for Floer complexes corresponding to multi-indices in which another one of the $\infty$ components is replaced by $0$ and $1$). These quasi-isomorphisms extend to give a quasi-isomorphism between $\CFm(Y_{\Lambda}(L))$ and a two-dimensional hypercube of complexes. We continue this until we get a quasi-isomorphism between $\CFm(Y_{\Lambda}(L))$ and the total complex of the hypercube $\hyp^\delta$.

Note that we have some freedom in this iteration: we could change the ordering of the components and start with $L_2$ instead of $L_1$, for example. We will choose the ordering as follows. For every $i = 1, \dots, \ell$, denote by $\Lambda^{(i)}$ the restriction of $\Lambda$ to 
$$L_i \cup L_{i+1} \cup \dots \cup L_\ell,$$ 
and by $\Lambda^{(i)}_i, \dots, \Lambda^{(i)}_\ell$ the respective framing vectors. We require that, for every $i=1, \dots, \ell$, either $\Lambda^{(i)}$ is nondegenerate, or else $\Lambda^{(i)}_i$ is in the $\qq$-span of $\Lambda^{(i)}_{i+1}, \dots, \Lambda^{(i)}_\ell$. Note that this can easily be arranged, by choosing which component we call $L_1$ first, then which component we call $L_2$, and so on.

Another degree of freedom in the iteration comes from the direction of the quasi-isomorphisms. In order to apply the results of Section~\ref{sec:exact}, at each step we have to choose a quasi-isomorphism as in \eqref{eq:qi1} or as in \eqref{eq:qi2}, depending on the framing $\Lambda$. At the first step, when we relate $\CFmd(Y_{\Lambda}(L))$ to $Cone(f_2^\delta)$, we choose \eqref{eq:qi1} if $\Lambda$ is nondegenerate and the restriction of $\Lambda$ to $Span(\Lambda_2, \dots, \Lambda_\ell)^{\perp}$ is zero or positive definite. (This corresponds to Cases I (a) and I (c) discussed in Section~\ref{sec:case1}.) We choose \eqref{eq:qi1} if $\Lambda$ is nondegenerate and the restriction of $\Lambda$ to $Span(\Lambda_2, \dots, \Lambda_\ell)^{\perp}$ is negative definite. When $\Lambda$ is degenerate (so, by our choice of ordering, $\Lambda_1$ is in the span of the other framing vectors), we are free to choose either \eqref{eq:qi1} or \eqref{eq:qi2}, compare Case II in Section~\ref{sec:case2}.

At the second step, we need to combine a quasi-isomorphism relating $\CFmd(Y_{\bar \Lambda}(L))$ to a mapping cone, and one relating $\CFmd(Y_{\Lambda'}(L'); \TR_1)$ to a mapping cone. Of course, these quasi-isomorphisms should go in the same direction. We choose the direction according to the same recipe as at the first step. Precisely, if $\bar \Lambda$ is nondegenerate, we choose the direction based on the sign of the restriction of $\bar \Lambda$ to $Span(\bar \Lambda_1, \Lambda_3, \dots, \Lambda_\ell)^{\perp}$. In particular, if $\Lambda$ was nondegenerate to start with, then $\bar \Lambda$ is automatically nondegenerate and, for $m_1$ sufficiently large, the relevant sign is negative if and only if the sign of the restriction of $\bar \Lambda$ to $Span(\tau_1, \Lambda_3, \dots, \Lambda_\ell)^{\perp}$ is negative; or, equivalently, if the sign of the restriction of $\Lambda'$ to $Span(\Lambda_3, \dots, \Lambda_\ell)^{\perp}$ is negative. If $\bar \Lambda$ is degenerate (for $m_1 \gg 0$), then $\Lambda'$ is degenerate also, and we choose the direction of the quasi-isomorphisms arbitrarily.
 
We continue to choose directions this way at the following steps. At step $i$, we look at the framing matrix $$\bar \Lambda^{i} = (\Lambda_1 + m_1 \tau_1, \dots, \Lambda_{i-1} + m_{i-1} \tau_{i-1}, \Lambda_i, \dots, \Lambda_\ell),$$ for $m_1, \dots, m_{i-1} \gg 0$.  If it is degenerate, we choose the direction arbitrarily. If it is nondegenerate, we choose it according to the sign of the restriction of $\bar \Lambda^{i}$ to the subspace
$$ Span(\Lambda_1 + m_1 \tau_1, \dots, \Lambda_{i-1} + m_{i-1} \tau_{i-1}, \Lambda_{i+1}, \dots, \Lambda_\ell)^\perp.$$
 
With these choices of ordering and quasi-isomorphism directions, at every step in the iteration process we can apply one of the refinements of Proposition~\ref{prop:leseqd} discussed in Sections~\ref{sec:case1}-\ref{sec:case2}: namely, Proposition~\ref{prop:gr+}, \ref{prop:gr-}, \ref{prop:gr0}, \ref{prop:gr00} or \ref{prop:gr002}. We first choose $m_1 \gg 0$ in such a way that the respective proposition applies, then we choose $m_2 \gg 0$, and so on. At the second step and later we may need to use the twisted coefficients variant, Proposition~\ref{prop:gr00twisted}. We claim that in the end we obtain the desired isomorphisms \eqref{eq:star'}, one for each $\ux \in \spc(Y_{\Lambda}(L))$.

Indeed, at least in the case when $\Lambda$ and all $\Lambda^{(i)}$'s are nondegenerate, all the complexes appearing in the proof of Proposition~\ref{prop:weaksurgery} decompose according to $\spc(Y_\Lambda(L)) \cong \H(L)/H(L, \Lambda)$. It is straightforward to check that the decompositions correspond to each other under the respective quasi-isomorphisms, and these quasi-isomorphisms preserve the relative $\zz$-gradings. Note that when iterating the exact sequences which give quasi-isomorphisms between mapping cones, at later steps in addition to the maps involved in Propositions~\ref{prop:gr+} and \ref{prop:gr-} we also have certain higher homotopies. However, these decompose into $\spc$ structures and preserve the relative $\zz$-gradings by the same arguments as those used in the discussion of the homotopy $H^\delta_{1, \ux}$ in Proposition~\ref{prop:gr+}, for example.

When $\Lambda=\Lambda^{(1)}$ or one of the other $\Lambda^{(i)}$'s is degenerate, there are two additional complications. First, in order to get a good decomposition into $\spc$ structures we have to replace the complex $\CC^\delta$ by a horizontally truncated, quasi-isomorphic one $\CC^\delta \langle \langle  \tilde \Lambda' \rangle \rangle$, as discussed in Section~\ref{sec:spc}. (When applying Propositions~\ref{prop:gr00} or \ref{prop:gr002}, this corresponds to focusing on a subset of all $\ux \in \spc(Y_{\Lambda}(L))$, and choosing the respective value $m_1$ so that $\delt(\ux)$ divides $m_1$ for all such $\ux$, compare Corollary~\ref{cor:small}. The unused $\spc$ structures $\ux$ give rise to trivial complexes, so we can ignore them.) We then get a decomposition of $\CC^\delta \langle \langle  \tilde \Lambda' \rangle \rangle$ according to $\spc$ structures $\ux$. We have similar decompositions of $\C^\delta  \langle \langle  \tilde \Lambda' \rangle \rangle$, as well as of all the other complexes in the proof of Proposition~\ref{prop:weaksurgery}, provided we truncate them with respect to $\tilde \Lambda'$. The respective truncations are quasi-isomorphic to the original complexes. Putting everything together, we obtain the desired isomorphisms \eqref{eq:star'}, for any $\ux \in \spc( Y_{\tilde \Lambda}(L))$. 

The second complication has to do with the grading-preserving properties of the isomorphisms \eqref{eq:star'}. If $\delt(\ux) = 0$, then when we apply Proposition~\ref{prop:gr0} at a step in the iteration process, the respective quasi-isomorphism preserves only the relative $\zz/2d\zz$-reduction of the relative $\zz$-grading. (Indeed, its target is only $\zz/2d\zz$-graded.) Thus, the resulting isomorphism \eqref{eq:star'} only preserves this $\zz/2d\zz$-grading. However, we can get such an isomorphism for any $d$ in a sequence $\{d_n\}$ with $d_n \to \infty$. Both sides of \eqref{eq:star'} are finite dimensional, relatively $\zz$-graded vector spaces, so if they are related by a $\zz/2d_n\zz$-grading preserving isomorphism for all $d_n$, they must in fact be isomorphic as relatively $\zz$-graded vector spaces. This completes the proof of the claim about the existence of a grading-preserving isomorphism \eqref{eq:star'}.

Finally, as $\delta$ varies, the isomorphisms \eqref{eq:star'} commute with the natural maps between the respective truncations, and we get the isomorphisms \eqref{eq:star2'}. 
\end {proof}

\begin{proof}[Proof of Theorem~\ref{thm:surgery} for basic systems]   Apply  \eqref{eq:star'} and Lemma~\ref{lemma:tors} (b) for $\ux$ torsion, and  \eqref{eq:star2'} and  Lemma~\ref{lemma:nontors} (b) for $\ux$ non-torsion. \end {proof}

\subsection {Link-minimal complete systems}
\label {sec:invariance}

Now that we have established the truth of Theorem~\ref{thm:surgery} for basic 
systems, we are left to do so for all link-minimal complete systems. 

\begin{proof}[Proof of Theorem~\ref{thm:surgery}]   Let $\Hyper$ be a link-minimal complete system of hyperboxes for a link $L$, and let $\Hyper_b$ be a basic system for $L$. According to Proposition~\ref{prop:moves} (b), the system $\Hyper$ can be obtained from $\Hyper_b$ by a sequence of system moved that do not include any index zero/three link stabilizations; i.e., by a sequence of 3-manifold isotopies, index one/two stabilizations and destabilizations, free index zero/three stabilizations, global shifts, elementary enlargements and contractions, and trajectory $\alpha$-slides.

According to Proposition~\ref{prop:moves2}, except for the trajectory $\alpha$-slides, all the other moves induce quasi-isomorphisms (in fact, chain homotopy equivalences) between the compressions of the respective hyperboxes of generalized Floer complexes $\Chain^-(\Hyper^{L', M}, \s)$ and $\Chain^-(\Hyper_b^{L', M}, \s)$. These compressions are the building blocks of the surgery complexes  $\C^-(\Hyper, \Lambda, \ux)$ and $\C^-(\Hyper_b, \Lambda, \ux)$. The quasi-isomorphisms commute with the 
restriction maps relating these building blocks, so by putting them together, we obtain a quasi-isomorphism between the two surgery complexes. 

The effect of a trajectory $\alpha$-slide on the surgery complex of a complete system was studied by Zemke in \cite[Corollary 13.5]{Zemke2}. In general, changing a trajectory $c_w$ by concatenating it with a curve $\gamma$ changes the surgery complex by a formula involving the action of the homology class $[\gamma] \in H_1(\Sigma; \Z)$. Examining the formula for the action of $\gamma$ in \cite[Section 13.3]{Zemke2}, we see that this action vanishes identically on the span of the alpha curves. In our case, $\gamma$ is a concatenation of the form $\tau * \alpha_j * \tau^{-1}$, where $\tau$ and $\alpha_j$ are as in Figure~\ref{fig:aslide} (a). In particular, $\gamma$ is homologous to the curve $\alpha_j$ which bounds a disk in the alpha handlebody, so the action of $\gamma$ on the surgery complex is trivial. Applying  Zemke's result, it follows that the two surgery complexes are quasi-isomorphic.

Since Theorem~\ref{thm:surgery} holds for the basic system $\Hyper_b$, we deduce that it must also hold for $\Hyper$.
\end {proof}

\section {The general link surgery formula}
\label{sec:general}

This section contains the precise statement, as well as the proof, of the link surgery formula (Theorem~\ref{thm:FirstSurgery}). Up to now we have established the formula in the case of link-minimal complete systems. The case of general complete systems requires several new ingredients: explicit resolutions $\Chain^-(\Hyper, \s)$ for the generalized Floer complexes $\Am(\Hyper, \s)$ (compare Section~\ref{sec:alternative}), as well as transition maps that allow one to define maps $\Pr^{\orM}_{\s}$ similar to the ones in Section~\ref{sec:inclusions}; cf. Remark~\ref{rem:notLM}. Moreover, the proof of the general surgery formula is based on an analysis of the effect of an index zero/three link stabilization on a complete system, and uses the material from Section~\ref{sec:03}.

The new constructions are rather involved so, to help the reader, we build up the general statement by first considering some simpler cases. In Section~\ref{sec:knot4} we describe the link surgery formula for a knot represented by a diagram with two $w$ and two $z$ basepoints. In Section~\ref{sec:toylink} we study the simplest version of an index zero/three link stabilization, relating a knot diagram with one basepoint of each type to a diagram with two basepoints of each type. In Section~\ref{sec:thetas} we explain how to define polygon maps between the resolutions $\Chain^-(\Hyper, \s)$, again in the simplest non-trivial example, that of a knot with two $w$ and two $z$ basepoints. In Section~\ref{sec:knots2p} we  describe the link surgery formula in more generality, still for a knot, but now with an arbitrary number of basepoints. In Section~\ref{sec:linktwo} we do yet another warm-up case, that of a link with two components, and with two basepoints of each type per component. Finally, in Section~\ref{sec:gens} we give the statement of the link surgery formula in full generality, and in Section~\ref{sec:proofgen} we prove it.

\subsection{A knot with four basepoints}
\label{sec:knot4}
We consider the case of a knot $K \subset Y$, represented by a Heegaard diagram $\Hyper^K$ with four linked basepoints: $w_1$, $z_1$, $w_2$ and $z_2$. We suppose $\Hyper^K$ is part of a complete system $\Hyper$, that also contains a diagram $\Hyper^{\emptyset}$ for $Y$, and two hyperboxes $\Hyper^{K, K}$ and $\Hyper^{K, -K}$; cf. Example~\ref{ex:completeknot}.

In Section~\ref{sec:chains} we introduced the link Floer complex $\CFLm(\Hyper^K)$, which in our case is freely generated over $\Ring=\Field[[U_1, U_2]]$ by intersection points $\x \in \Ta \cap \Tb$. Let us consider the infinity version of this complex, obtained by inverting the $U$ variables:
$$ \CFLI(\Hyper^K) = (U_1, U_2)^{-1}\CFLm(\Hyper^K).$$ 

As a direct product, the complex $\CFLI(\Hyper^K)$ has generators (over $\Field$) $U_1^{n_1} U_2^{n_2} \x$, with $n_i \in \Z$ and $\x \in \Ta \cap \Tb.$ The differential takes into account intersection with $w$'s, so we write $\CFLI(\Hyper^K)= C\{w\}$. On the other hand, $\CFLm(\Hyper^K)$ only has the generators that satisfy $n_1, n_2 \geq 0$, so we write
$$ \CFLm(\Hyper^K) = C\{ w; n_1 , n_2 \geq 0 \}.$$
For each $s \in \Z$, we also have a subcomplex
$$ \Am(\Hyper^K, s) = C\{ w; n_1 , n_2 \geq 0, A(\x) - n_1 - n_2 \leq s\}.$$

For the surgery formula, we have an inclusion $\Pr^K_{s}$ of $\Am(\Hyper^K, s)$ into $\CFLm(\Hyper^K)= C\{ w; n_1 , n_2 \geq 0\}$; see Section~\ref{sec:inclusions}. We can then use the Heegaard moves specified by the hyperbox $\Hyper^{K, K}$ to relate $\CFLm(\Hyper^K)$ to $\CFm(\Hyper^{\emptyset})$. The composition could play the role of the map $\Phi^{K}_s$ (just as in the case of link-minimal systems).

However, we also need a map $\Phi^{-K}_s$, which would involve a chain map $\Pr^{-K}_s$ from $\Am(\Hyper^K, s)$ to $ C\{ z; n_1 , n_2 \geq 0\}$, a complex with the differential counting $z$'s. It does not seem possible to construct a suitable such map. To deal with this problem, we will replace $\Am(\Hyper^K, s)$ with an explicit resolution, $\Chain^-(\Hyper^K, s)$.

In the case of a knot with two basepoints, the map $\Pr^{-K}_s$ involves the inclusion of the complex $\Am(\Hyper^K, s) = C\{w; n\geq 0, A(\x) - n \leq s\}$ into $C\{w; A(\x) - n \leq s\}$. The latter complex is then identified with $C\{z; n \geq 0\}$ by using the isomorphism $\x \to U^{s-A(\x)}\x$.

Similarly, for four basepoints, we have an inclusion of $\Am(\Hyper^K, s)$ into the intermediate complex 
$$\Cint = C\{ w;  n_2 \geq 0, A(\x) - n_1 - n_2 \leq s\}.$$

We will start by constructing a resolution of $\Cint$, denoted $\Ccint$, and then identify a subcomplex of $\Ccint$ as the resolution $\Chain^-(\Hyper^K, s)$ of $\Am(\Hyper^K,s)$. We will also construct a chain map (called a transition map) from $\Ccint$ to $ C\{ z; n_1 , n_2 \geq 0\}$.

Consider the following elements of $\Cint$:
$$ \x_0 := U_1^{A(\x) - s} \x,  \ \x_1:=U_1^{A(\x)-s-1}U_2 \x, \ \x_2:=U_1^{A(\x)-s-2}U_2^2 \x, \ \dots $$
over all possible $\x \in \Ta \cap \Tb.$ Every element of $\Cint$ can be written as an infinite sum
\begin{equation}
\label{eq:isum}
 \sum_{\x \in \Ta \cap \Tb} \sum_{i=1}^{\infty} U_1^{a_i} U_2^{b_i} \x_{n_i} =  \sum_{\x \in \Ta \cap \Tb} \sum_{i=1}^{\infty} U_1^{A(\x) - s - n_i + a_i} U_2^{b_i + n_i} \x,
 \end{equation}
where $(a_i, b_i, n_i)$ is a sequence (depending on $\x$) consisting of distinct triples of nonnegative integers. Furthermore, in order to avoid infinitely many negative powers of $U_1$ in the sum, we need to impose the condition
\begin{equation}
\label{eq:bdd}
 \exists K >0 \text{ such that } n_i - a_i < K, \text{ for all } i.
 \end{equation}
Observe also, that, if \eqref{eq:bdd} is satisfied, then \eqref{eq:isum} is a well-defined sum. Indeed, because $a_i, n_i, b_i \geq 0$, we cannot have that the pair $(A(\x) - s - n_i + a_i, b_i + n_i)$ is some fixed $(a, b)$ for infinitely many values of $i$.

\begin{remark}
Informally, we can think of the elements $\x_n$ as generators for $\Cint$. If we had worked with polynomials in $U_1, U_2$ instead of power series, then every element of the resulting $\Cint$ would have been a finite linear combination of the $\x_n$. With power series, however, we need to allow infinite sums satisfying \eqref{eq:bdd}. 
\end{remark}

The elements $\x_n$ are not free, as they satisfy relations of the form
$$ U_2 \x_n = U_1 \x_{n+1}.$$

To get a resolution of $\Cint$, we turn $\x_n$ into linearly independent elements, by introducing a new variable $V_2$, so that $\x_n$ corresponds to $V_2^n \x$. We also need to make sure $U_2 V_2^n \x - U_1 V_2^{n+1} \x=  (U_2 + U_1V_2) V_2^n \x$ is in the image of the differential, so we introduce a further variable $Y_2$ so that
$$ Y_2^2=0, \ \ \del Y_2 = U_2 + U_1 V_2.$$
 
Furthermore,  in order to  determine the differential of some $\x \in \Ccint$, note that the differential of the corresponding element $U_1^{A(\x)-s} \x \in \Cint$ consists of terms of the form
$$ U_1^{A(\x)-s} \cdot U_1^{n_{w_1}(\phi)}U_2^{n_{w_2}(\phi)} \y = U_1^{n_{z_1}(\phi) + n_{z_2}(\phi)} \Bigl( U_1^{A(\y)-s-n_{w_2}(\phi)}U_2^{n_{w_2}(\phi)} \y \Bigr).$$
The last expression in parentheses corresponds to $V_2^{{n_{w_2}(\phi)}} \y$. Thus, $\del \x$ in $\Ccint$ should consist of terms of the form
$$ U_1^{n_{z_1}(\phi) + n_{z_2}(\phi)} V_2^{n_{w_2}(\phi)}\y.$$

Here is the precise definition of $\Ccint$. We first consider the subring
$$ \Ringbig =\bigl \{ \sum_{i=1}^{\infty} U_1^{a_i} U_2^{b_i} V_2^{n_i} \mid n_i - a_i \text{ is bounded above} \bigr\} \subset \Field[[U_1, U_2, V_2]].$$
This is a module over $\Field[[U_1, U_2]]$, and contains the polynomial module $\Field[[U_1, U_2]][V_2].$ We will write
$$ \Ringbig =  \Field[[U_1, U_2]][V_2]'$$
where the prime indicates that we allow infinite sums satisfying the boundedness condition \eqref{eq:bdd}.

We define $\Ccint^0$ to be the free complex over $\Ringbig$ with generators $\x \in \Ta \cap \Tb$, and differential
\begin{equation}
\label{eq:delint}
\del \x = \sum_{\y \in \Ta \cap \Tb}  \sum_{\substack{\phi \in \pi_2(\x, \y) \\ \mu(\phi)=1} }  \# (\M(\phi)/\R) \cdot U_1^{n_{z_1}(\phi) + n_{z_2}(\phi)} V_2^{{n_{w_2}(\phi)}} \y.
\end{equation}

\begin{lemma}
The map $\del$ on $\Ccint^0$ satisfies $\del^2=0$.
\end{lemma}

\begin{proof}
The usual arguments in Floer theory apply, with $\del^2$ counting the ends of a one-dimensional moduli space of pseudo-holomorphic disks. One caveat is the presence of disk bubbles. According to the analysis in \cite[Theorem 5.5]{Links}, there are four such bubbles, one going over $w_1$ and $z_1$, one over $z_1$ and $w_2$, one over $w_2$ and $z_2$, and one over $z_2$ and $w_1$. Their contribution to $\del^2 \x$ is
$$ U_1 \x + U_1 V_2 \x + U_1 V_2 \x + U_1  \x =0.$$
This shows that $\del^2 \x=0$. 
\end{proof}

We now let $\Ccint$ be the mapping cone
\begin{equation}
\label{eq:hcone}
Y_2 \cdot \Ccint^0 \xrightarrow{U_2 + U_1V_2} \Ccint^0.
 \end{equation}

Alternatively, we can simply define $\Ccint$ starting from the free dg (differential graded) module over the dg algebra 
$$ \Algbig=\Field[[U_1, U_2]][V_2, Y_2]'/(Y_2^2=0, \del Y_2 = U_2 + U_1 V_2),$$
with generators $\x \in \Ta \cap \Tb$, by introducing the differential on generators as in \eqref{eq:delint}. The homological grading on $\Algbig$ is given by $\gr(U_i)=-2, \gr(V_2)=0, \gr(Y_2)=-1$. The prime again indicates the condition \eqref{eq:bdd} in infinite sums. 

However, most of the time we will view $\Ccint$ as a complex over the ring $\Ring=\Field[[U_1, U_2]]$, rather than over $\Ringbig$ or $\Algbig$.

Define a chain map $P: \Ccint \to \Cint$ by
\begin{equation}
\label{eq:Fmap}
P(V_2^n \x)= U_1^{A(\x) - s -n} U_2^n  \cdot \x, \ \ \ P(V_2^n Y_2 \x)=0,
\end{equation}
and extending it linearly (including to the allowed infinite sums) over $\Ring=\Field[[U_1, U_2]]$.

\begin{lemma}
\label{lem:Fchain}
$P$ is a chain map.
\end{lemma}

\begin{proof}
In the expression $(\del P + P \del)(V_2^n \x)$, for each holomorphic disk in a class $\phi \in \pi_2(\x, \y) $, we obtain a term involving $\y$ with coefficient 
$$  U_1^{A(\x) - s-n}U_2^n \cdot U_1^{n_{w_1}(\phi)}  U_2^{n_{w_2}(\phi)} + U_1^{n_{z_1}(\phi) + n_{z_2}(\phi)}  \cdot U_1^{A(\y) - s - n- n_{w_2}(\phi)} U_2^{n_{w_2(\phi)}+n}.$$
Since $A(\x) - A(\y) = n_{z_1}(\phi) + n_{z_2}(\phi) - n_{w_1}(\phi) - n_{w_2}(\phi)$, these terms cancel.

We also have
\begin{align*}
(\del P + P \del)(V_2^n Y_2 \x) &= P(\del ( V_2^n Y_2\x)) \\ 
&= P(U_2 V_2^n \x + U_1 V_2^{n+1} \x + V_2^n Y_2 (\del \x)) \\
&= U_2 \cdot U_1^{A(\x) - s -n} U_2^{n}  \cdot \x + U_1 \cdot U_1^{A(\x) - s -n-1} U_2^{n+1}  \cdot \x + 0\\
&=0.
\end{align*}
This shows that $P$ is a chain map. 
\end{proof}

\begin{lemma}
\label{lem:CintResolution}
The map $P: \Ccint \to \Cint$ exhibits $\Ccint$ as a resolution of $\Cint$  over $\Ring=\Field[[U_1, U_2]]$; that is, $P$ is a surjective quasi-isomorphism.
\end{lemma}

\begin{proof}
Clearly $P$ is surjective; we need to show that it is also a quasi-isomorphism. 
Let $G$ be the kernel of $P$, so that we have a short exact sequence:
$$ 0 \longrightarrow G \longrightarrow \Ccint \longrightarrow \Cint \longrightarrow 0.$$
By considering the long exact sequence in homology, we see that it suffices to prove that $G$ is acyclic. 

Note that $G$ is freely generated over $\Ringbig = \Field[[U_1, U_2]][V_2]'$ by the elements of the form $(U_2  + U_1 V_2)\x$ and $Y_2\x$. The map $H: G \to G$ given by
$$ H((U_2  + U_1 V_2)\x) = Y_2\x, \ \ \ H(Y_2\x)=0$$
is a chain homotopy between the identity and zero. Hence $G$ is acyclic.
\end{proof}

Next, we define a filtration $\F$ on $\Cint$ by the negative of the exponent of $U_1$, i.e., by setting $$\F(\x) = 0, \ \F(U_1)=-1, \ \F(U_2)=0.$$ Note that $\Am(\Hyper^K, s)$ is the subcomplex of $\Cint$ given by $\F \leq 0$. 

We also define a filtration $\FF$ on $\Ccint$ by setting
\begin{equation}
\label{eq:FF}
 \FF(\x) = -A(\x), \ \FF(U_1)=-1, \ \FF(U_2)=0, \ \FF(V_2)=1, \ \FF(Y_2)=0.
 \end{equation}
Using \eqref{eq:delint} and \eqref{eq:hcone}, it is easy to check that the differential of $\Ccint$ respects the filtration. 

One can also verify that the map $P$ from \eqref{eq:Fmap} is filtered of degree $-s$ with respect to $\FF$ and $\F$. Indeed, we have
$$ \FF(V_2^n \x) = n-A(\x)= \F(U_1^{A(\x)-s-n} U_2^n \x)-s.$$

We define the resolution $\Chain^-(\Hyper^K, s)$ of $\Am(\Hyper^K,s)$ to be the subcomplex of $\Ccint$ corresponding to filtration level $\FF \leq s$. Concretely, $\Chain^-(\Hyper^K, s)$ consists of  infinite sums (with coefficients in $\Ring$, and a suitable boundedness condition) of the elements
\begin{equation}
\label{eq:txna}
 \tx_n^a := U_1^{\max(s+n-A(\x), 0)} V_2^n Y_2^a \x, \ \x \in \Ta \cap \Tb, \ n\in \Z_{\geq 0}, \ a \in \{0,1\}.
 \end{equation}

Let
$$\tilde{P}: \Chain^-(\Hyper^K, s) \to \Am(\Hyper^K,s)$$
be the restriction of $P$ from \eqref{eq:Fmap}. In terms of the elements $\tx_n^0$ and $\tx_n^1$, we have
\begin{equation}
\label{eq:tildeF}
 \tilde{P}(\tx_n^0) = U_1^{\max(A(\x) - s -n, 0)} U_2^n  \cdot \x, \ \ \ \tilde{P}(\tx_n^1) =0.
 \end{equation}

\begin{lemma}
\label{lem:AResolution}
The map $\tilde{P}: \Chain^-(\Hyper^K, s) \to \Am(\Hyper^K, s)$ exhibits $\Chain^-(\Hyper^K, s)$ as a resolution of $\Am(\Hyper^K, s)$. 
\end{lemma}

\begin{proof}
This is a filtered version of Lemma~\ref{lem:CintResolution}. Note that the kernel  $\tilde{G}=\ker(\tilde{P})$ is the subgroup of $G =\ker(P)$ in filtration level $\FF \leq s$. Moreover, the null-homotopy $H: G \to G$ from the proof of Lemma~\ref{lem:CintResolution} is a filtered map; hence, by restricting it to $\tilde{G}$, we get that $\tilde{G}$ is acyclic.
\end{proof}

We now construct a {\em transition map}
$$ \Trans : \Ccint \to C\{z; n_1, n_2 \geq 0\}$$
by setting
\begin{equation}
\label{eq:Trans1}
 \Trans (V_2^n Y_2 \x) = \x
 \end{equation}
and
\begin{equation}
\label{eq:Trans0}
 \Trans(V_2^n \x) =  \sum_{\y \in \Ta \cap \Tb}  \sum_{\substack{\phi \in \pi_2(\x, \y) \\ \mu(\phi)=1} }  \# (\M(\phi)/\R) \cdot U_1^{n_{z_1}(\phi)} \frac{U_1^{n_{z_2}(\phi)} - U_2^{n_{z_2}(\phi)}}{U_1 - U_2} \y.
 \end{equation}
This definition is inspired from \cite{SarkarMoving}. Of course, by $ \frac{U_1^{n_{z_2}(\phi)} - U_2^{n_{z_2}(\phi)}}{U_1 - U_2}$ we actually mean the polynomial
$$ U_1^{n_{z_2}(\phi) - 1}+ U_1^{n_{z_2}(\phi)-2} U_2 + \dots + U_1 U_2^{n_{z_2}(\phi)-2} + U_2^{n_{z_2}(\phi) - 1}.$$
We extend $\Xi$ to be a module map over $\Ring = \ff[[U_1, U_2]]$.

\begin{lemma}
\label{lem:TransChain0}
$\Trans$ is a chain map.
\end{lemma}

\begin{proof}
Note that $\Trans$ is a module map over $\Ringbig=\Field[[U_1, U_2]][V_2]'$, where $V_2$ acts on the target by the identity. Thus, it suffices to check that $\Trans \del + \del \Trans=0$ on generators of the form $\x$ and $Y_2\x$.

Let us first check that
$$\Trans \del (\x) + \del \Trans (\x)=0.$$
 In the expression $\Trans \del (\x)$, a pair of holomorphic disks in classes $\phi \in \pi_2(\x, \a), \psi \in \pi_2(\a, \y)$ contributes
$$ U_1^{n_{z_1}(\psi)} \frac{U_1^{n_{z_2}(\psi)} - U_2^{n_{z_2}(\psi)}} {U_1 - U_2} \cdot U_1^{n_{z_1}(\phi) + n_{z_2}(\phi)}\y.$$

In $\del \Trans (\x)$, the same pair contributes
$$ U_1^{n_{z_1}(\psi)} U_2^{n_{z_2}(\psi)} \cdot U_1^{n_{z_1}(\phi)} \frac{U_1^{n_{z_2}(\phi)} - U_2^{n_{z_2}(\phi)} }{U_1 - U_2} \y.$$

Adding these together we get a total contribution of
\begin{equation}
\label{eq:totch}
U_1^{n_{z_1}(\phi * \psi)} \frac{U_1^{n_{z_2}(\phi * \psi)} - U_2^{n_{z_2}(\phi * \psi)}}{U_1 - U_2} \y.\end{equation}

This only depends on the combined domain $\phi * \psi$. The usual proof that $\del^2=0$ applies here, and shows that adding up all these contributions gives $0$. Note that when $\x = \y$, there are four periodic domains to take care of: one contains $w_1$ and $z_1$, one $w_1$ and $z_2$, one $w_2$ and $z_1$, and the last $w_2$ and $z_2$. Their contributions add up to zero, because the $w$'s do not play a role in the expression \eqref{eq:totch}.

Now let us check that
$$\Trans \del (Y_2 \x) + \del \Trans (Y_2 \x)=0.$$
We have
\begin{align*}
\Trans \del (Y_2 \x) &= \Trans \Bigl(  U_2 \x + U_1V_2 \x + \sum_{\y \in \Ta \cap \Tb}  \sum_{\substack{\phi \in \pi_2(\x, \y) \\ \mu(\phi)=1} }  \# (\M(\phi)/\R) \cdot U_1^{n_{z_1}(\phi) + n_{z_2}(\phi)} V_2^{n_{w_2}(\phi)} Y_2 \y \Bigr)\\
&=  \sum_{\y \in \Ta \cap \Tb}  \sum_{\substack{\phi \in \pi_2(\x, \y) \\ \mu(\phi)=1} }  \# (\M(\phi)/\R) \Bigl ( (U_2+U_1) \cdot U_1^{n_{z_1}(\phi)} \frac{U_1^{n_{z_2}(\phi)} - U_2^{n_{z_2}(\phi)}}{U_1 - U_2} +  U_1^{n_{z_1}(\phi) + n_{z_2}(\phi)} \Bigr ) \y \\
&=  \sum_{\y \in \Ta \cap \Tb}  \sum_{\substack{\phi \in \pi_2(\x, \y) \\ \mu(\phi)=1} }  \# (\M(\phi)/\R) U_1^{n_{z_1}(\phi)} U_2^{n_{z_2}(\phi)} \y.
\end{align*}

On the other hand, we also have
$$
\del \Trans (Y_2 \x) =\del \x = \sum_{\y \in \Ta \cap \Tb}  \sum_{\substack{\phi \in \pi_2(\x, \y) \\ \mu(\phi)=1} }  \# (\M(\phi)/\R) U_1^{n_{z_1}(\phi)} U_2^{n_{z_2}(\phi)} \y.$$
This concludes the proof.
\end{proof}

We now have all the ingredients to write down the surgery formula for $K$ (with four basepoints). This takes the form of a mapping cone, where the initial complex is the direct product of $\Chain^-(\Hyper^K, s)$ over $s \in \Z$, and the final complex is the direct product of copies of some $\CFm(\Hyper^{\emptyset}) = \CFm(S^3)$.

The map $\Phi^K_s$ is defined to be the composition
\begin{equation}
\label{eq:PhiK}
 \Chain^-(\Hyper^K, s) \xrightarrow{\tilde P} \Am(\Hyper^K, s) \hookrightarrow C\{w; n_1, n_2 \geq 0\} \xrightarrow{\sim} \CFm(\Hyper^{\emptyset}),
 \end{equation}
with the last map being an equivalence obtained from the Heegaard moves in $\Hyper^{K, K}$. In fact, we denote the composition of the first two maps in \eqref{eq:PhiK} by
$$ \I^{K}_s:  \Chain^-(\Hyper^K, s) \to C\{w; n_1, n_2 \geq 0\} = \CFLm(\Hyper^K) = \CFm(p^{K}(\Hyper^K)).$$
This plays the role of the inclusion map from Section~\ref{sec:inclusions}. Since in our case it also includes the projection $\tilde P$, we call $\I^K_s$ a {\em projection-inclusion map}.

The final equivalence in \eqref{eq:PhiK}, induced by Heegaard moves, is a descent map
$$ D^K : \CFm(p^{K}(\Hyper^K)) \to \CFm(\Hyper^{\emptyset}),$$
similar to the ones considered in Section~\ref{subsec:desublink}. Thus, we can write 
$$\Phi^K_s = D^K \circ \I^K_s,$$
just as before.

On the other hand, the map $\Phi^{-K}_s$ is defined to be the composition
\begin{equation}
\label{eq:Phi-K}
 \Chain^-(\Hyper^K, s) \hookrightarrow \Ccint \xrightarrow{\Trans} C\{z; n_1, n_2 \geq 0\} \xrightarrow{\sim} \CFm(\Hyper^{\emptyset}),
 \end{equation}
with the last map again obtained from Heegaard moves, this time the ones from the hyperbox $\Hyper^{K, -K}$. In this context, we let $\I^{-K}_s$ be the first inclusion above (into $\Ccint$), and we define the descent map $D^{-K}$ to be the composition of $\Trans$ and the equivalences induced by Heegaard moves. We can still write $\Phi^{-K}_s = D^{-K} \circ \I^{-K}_s$.

Just as in the link-minimal case, by summing up the $\Phi^K_s$ and $\Phi^{-K}_s$ maps (with suitable shifts in $s$), we obtain the surgery complexes $\C^-(\Hyper, n)$ for $n \in \Z$.

\subsection{Invariance under an index zero/three link stabilization}
\label{sec:toylink}
As a warm-up to the general case of invariance, in this section we will study the effect of an index zero/three link stabilization on a minimally-pointed complete system of hyperboxes for a knot $K \subset Y$. 

Let $\bar \Hyper$ denote such a complete system. Recall from Example~\ref{ex:completeknot} that $\bar \Hyper$ consists of:
\begin{itemize}
\item a Heegaard diagram $\bar \Hyper^K$ for $K$ with two basepoints $w_1$ and $z_1$,
\item a Heegaard diagram $\bar \Hyper^{\emptyset}$ for $Y$ with the single basepoint $w_1$, 
\item a hyperbox (a sequence of Heegaard moves with $\Theta$ elements) $\bar \Hyper^{K, K}$ from $p^{K}(\bar \Hyper^K)$ to $\bar \Hyper^{\emptyset}$, 
\item another hyperbox, $\bar \Hyper^{K, -K}$, relating $p^{-K}(\bar \Hyper^K)$ to $\bar \Hyper^{\emptyset}$. 
\end{itemize}
From $\bar \Hyper$ we can construct a surgery complex $\C^-(\bar \Hyper, n)$ for each $n \in \Z$. This is a free complex over $\Field[[U_1]]$.

\begin{figure}
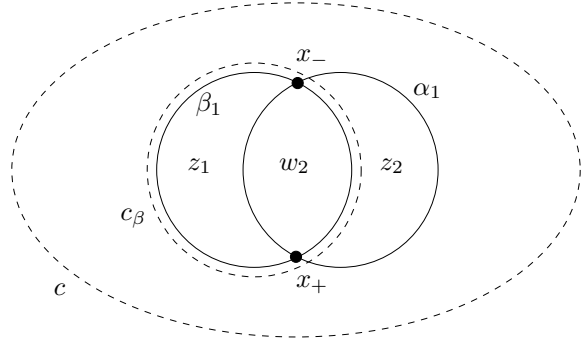

\begin{center}

\end{center}
\caption {An index zero/three link stabilization.}
\label{fig:stabilize}
\end{figure}

We now do a stabilization as in Figure~\ref{fig:stabilize}, by introducing two new curves $\alpha_1$ and $\beta_1$, and two new basepoints $w_2$ and $z_2$ in the neighborhood of the existing $z_1$. We let $\Hyper$ be the complete system obtained from $\bar \Hyper$ by this stabilization; cf. Section~\ref{sec:moves}.

Thus, in the complete system $\Hyper$ we have:
\begin{itemize}
\item  the diagram $\Hyper^K$ obtained from $\bar \Hyper^K$ by the link stabilization,
\item the diagram $\Hyper^{\emptyset}$ obtained from $\bar \Hyper^{\emptyset}$ by a free index zero/three stabilization introducing $\alpha_1$, $\beta_1$ and $w_2$,
\item a hyperbox $\Hyper^{K, K}$ from $p^{K}(\Hyper^K)$ to $\Hyper^{\emptyset}$ obtained from $\bar \Hyper^{K, K}$ by doing free index zero/three stabilizations at each step,
\item a hyperbox $\Hyper^{K, -K}$ from $p^{-K}(\Hyper^K)$ to $\Hyper^{\emptyset}$ obtained $\bar \Hyper^{K, -K}$ by doing free index zero/three stabilizations at each step, and then the move shown in Figure~\ref{fig:StabTriple0}.
\end{itemize}
From $\Hyper$ we can form a surgery complex $\C^-(\Hyper, n)$ as explained in Section~\ref{sec:knot4}. This is a free complex over $\Field[[U_1, U_2]]$.

In Section~\ref{sec:hd03} we described the behavior of holomorphic disks under an index zero/three link stabilization; cf. Proposition~\ref{prop:ZemkeDisks}. For our calculations, we will choose almost complex structures as in that proposition, that is, with sufficiently stretched necks along the curves $c$ and $c_{\beta}$ from Figure~\ref{fig:StabTriple0}.

\begin{proposition}
For every $n \in \Z$, the surgery complexes $\C^-(\bar \Hyper, n)$ and $\C^-(\Hyper, n)$ are quasi-isomorphic over $\Field[[U_1]]$.
\end{proposition}

\begin{proof}
Let $\bar C\{w; n \geq 0\}=\CFm(p^{K}(\bar \Hyper^K))$ and $\bar C\{z; n \geq 0 \}=\CFm(p^{-K}(\bar \Hyper^K))$ be the complexes obtained from the diagram $\bar \Hyper^K$ by using the $w_1$ resp. the $z_1$ basepoint. Then, as in the proof of Theorem~\ref{thm:LinkInvariance} or in the discussion from Section~\ref{sec:hd03}, we can identify $C\{w; n_1, n_2 \geq 0\}=\CFm(p^{K}(\Hyper^K))$ with the mapping cone
\begin{equation}
\label{eq:icecream}
  \bar C\{w; n \geq 0\}_-[[U_2]] \xrightarrow{U_1 - U_2} \bar C\{w; n \geq 0\}_+[[U_2]],
  \end{equation}
where the subscripts $-$, $+$ indicate the point in $\alpha_1 \cap \beta_1$ that is part of the respective generators ($x_-$ or $x_+$). Hence, by Corollary~\ref{cor:rhonormal}, we have a chain homotopy equivalence (and hence a quasi-isomorphism)
$$ \rho: C\{w; n_1, n_2 \geq 0\} \to \bar C\{w; n \geq 0\},$$
acting on the domain of \eqref{eq:icecream} by zero, and on the target by $U_1^{n_1} U_2^{n_2} \x \to U_1^{n_1 + n_2} \x$. 

By a slight abuse of notation, we will denote by $\rho$ all quasi-isomorphisms of a similar form. For example, there are such equivalences $\rho$ relating the complexes $\CFm$ for all the diagrams in the hyperbox $\Hyper^{K, K}$ to the corresponding ones in $\bar \Hyper^{K, K}$. Furthermore, according to Proposition~\ref{prop:StabPolygonH}, these equivalences commute with the descent maps induced by Heegaard moves. Therefore, we have a commutative diagram
\begin{equation}
\label{eq:commuteplus0}
\begin{CD} 
C\{w; n_1, n_2 \geq 0\}=\CFm(p^{K}(\Hyper^K)) @>{D^K}>> \CFm(\Hyper^{\emptyset}) \\
@V{\rho}VV @VV{\rho}V  \\
\bar C\{w; n \geq 0\}=\CFm(p^{K}(\bar \Hyper^K)) @>{D^K}>> \CFm(\bar \Hyper^{\emptyset}).
\end{CD}
\end{equation}

Theorem~\ref{thm:LinkInvariance}(b) gives yet another equivalence
$$  \Am(\Hyper^K, s) \xrightarrow{\sim} \Am(\bar \Hyper^K, s).$$
We pre-compose this with the map $\tilde P$ from \eqref{eq:tildeF}, which is a quasi-isomorphism  according to Lemma~\ref{lem:AResolution}. The result is a quasi-isomorphism
\begin{equation}
\label{eq:rhoA}
 \rho: \Chain^-(\Hyper^K, s) \to  \Am(\bar \Hyper^K, s)= \Chain^-(\bar \Hyper^K, s).
 \end{equation}
By construction, we get a commutative diagram
$$\begin{CD}
 \Chain^-(\Hyper^K, s) @>{\I^K_s}>> \CFm(p^{K}(\Hyper^K)) \\
@V{\rho}VV @VV{\rho}V  \\
 \Chain^-(\bar \Hyper^K, s) @>{\I^K_s}>> \CFm(p^{K}(\bar \Hyper^K)).
\end{CD}
$$
Combining this with \eqref{eq:commuteplus0}, we obtain a commutative diagram
\begin{equation}
\label{eq:commuteplus}
\begin{CD} 
 \Chain^-(\Hyper^K, s) @>{\Phi^K}>> \CFm(\Hyper^{\emptyset}) \\
@V{\rho}VV @VV{\rho}V  \\
  \Chain^-(\bar \Hyper^K, s) @>{\Phi^K}>> \CFm(\bar \Hyper^{\emptyset}).
\end{CD}
\end{equation}

We aim to construct a similar diagram involving $\Phi^{-K}$ instead of $\Phi^K$:
\begin{equation}
\label{eq:commuteminus}
\xymatrix{
 \Chain^-(\Hyper^K, s) \ar[r]^{\Phi^{-K}} \ar[d]_{\rho} \ar[dr] & \CFm(\Hyper^{\emptyset})  \ar[d]^{\rho}  \\
  \Chain^-(\bar \Hyper^K, s) \ar[r]^{\Phi^{-K}} & \CFm(\bar \Hyper^{\emptyset}).
}
\end{equation}
This time, the maps $\rho$ will not commute with $\Phi^{-K}$ on the nose, but only up to chain homotopy. Still, by introducing the chain homotopy as a diagonal map in the diagram above, we can form a chain map from the cone of the first row to the cone of the second row. This map is a quasi-isomorphism because the vertical maps are.

Recall that  the $\Phi^{K}$ and $\Phi^{-K}$ maps are the building blocks in the surgery complex. Thus, once we have \eqref{eq:commuteplus} and \eqref{eq:commuteminus}, we can form a quasi-isomorphism between the surgery complexes $\C^-(\Hyper, n)$ and $\C^-(\bar \Hyper, n)$. 

We are left to construct \eqref{eq:commuteminus}. By Proposition~\ref{prop:ZemkeDisks}, we can identify $C\{z; n_1, n_2 \geq 0\}$ with
\begin{equation}
\label{eq:zcone}
 \bar C\{z; n \geq 0 \}^{U_1 \to U_2}_+[[U_1]] \xrightarrow{U_1 - U_2} \bar C\{z; n \geq 0 \}^{U_1 \to U_2}_-[[U_1]],
 \end{equation}
and hence, by Corollary~\ref{cor:rhovariant}, we have an equivalence
$$ \rho: C\{z; n_1, n_2 \geq 0\} \to \bar C\{z; n \geq 0\}.$$

We can do the same with the complexes appearing in the hyperboxes $\Hyper^{K, -K}$ versus the ones in $\bar \Hyper^{K, -K}$. We obtain a diagram that commutes up to some diagonal chain homotopy:
\begin{equation}
\label{eq:commuteminus0}
\xymatrix{
C\{z; n_1, n_2 \geq 0\}=\CFm(p^{-K}(\Hyper^K)) \ar[r]^-{\sim} \ar[d]_{\rho} \ar[dr] & \CFm(\Hyper^{\emptyset}) \ar[d]^{\rho} \\
\bar C\{z; n \geq 0\}=\CFm(p^{-K}(\bar \Hyper^K)) \ar[r]^-{\sim} & \CFm(\bar \Hyper^{\emptyset}),
}
\end{equation}
where the horizontal maps are equivalences induced by the Heegaard moves in the respective hyperboxes. The fact that these maps commute with the projections $\rho$, up to chain homotopy, follows from Propositions~\ref{prop:PolyDestab03b} and \ref{prop:PolyDestabMove}. (The latter proposition applies for the final map in the stabilized hyperbox, the one shown in Figures~\ref{fig:StabTriple0} and \ref{fig:StabTriple}.) 

We claim that there is also a diagram
\begin{equation}
\label{eq:commuteminus1}
\xymatrix{
 \Chain^-(\Hyper^K, s) \ar[r]^-{\I_s^{-K}} \ar[d]_{\rho} &  \Ccint \ar[r]^-{\Xi} \ar[d]^{\rho} \ar[dr]^{Z}   & C\{z; n_1, n_2 \geq 0\}  \ar[d]^{\rho}\\
  \Chain^-(\bar \Hyper^K, s) \ar[r]^-{\I_s^{-K}} & C\{z; n \geq 0\} \ar[r]^{\id} & \bar C\{z; n \geq 0\}.
}
\end{equation}
where the left square in the diagram commutes, and the right square commutes up to the chain homotopy defined by $Z$, that is:
\begin{equation}
\label{eq:checkZ}
 \rho+ \rho \circ \Xi = \del \circ Z  + Z \circ  \del.
\end{equation}
In the diagram \eqref{eq:commuteminus1}, the two maps that we still have to define are the middle vertical map 
\begin{equation}
\label{eq:midrho}
\rho: \Ccint \to C\{z; n \geq 0\}, 
\end{equation}
and the chain homotopy $Z$.

We first note that, in view of the formula \eqref{eq:delint} and the description of the holomorphic disks in Proposition~\ref{prop:ZemkeDisks}, we can identify $\Ccint^0$ with the cone
\begin{equation}
\label{eq:Ccint0id}
  \bar C\{z; n \geq 0\}^{U_1 \to U_2}_-[[U_1]][V_2]' \xrightarrow{V_2 - 1} \bar C\{z; n \geq 0\}^{U_1 \to U_2}_+[[U_1]][V_2]'.
  \end{equation}
 The same argument as in Lemma~\ref{lem:ap1} shows that $ \Ccint^0$ is chain homotopy equivalent to $C\{z; n \geq 0\}^{U_1 \to U_2}[[U_1]]$ over $\Field[[U_1, U_2]]$, via the projection that takes 
 \begin{equation}
 \label{eq:somep}
 \x \times x_+ \to \x, \ \ \x \times x_- \to 0, \ \ V_2 \to 1,
 \end{equation}
 for $\x \in \Ta \cap \Tb$, where $\Ta$ and $\Tb$ are the totally real tori coming from the diagram $\bar \Hyper^{K}$.
 
Furthermore, $\Ccint$ is obtained from $\Ccint^0$ by taking the cone of $U_2 + U_1 V_2$. Consequently, using \eqref{eq:somep}, we get that $\Ccint$ is chain homotopy equivalent to the cone of $U_2 + U_1$ on $C\{z; n \geq 0\}^{U_1 \to U_2}[[U_1]]$. In turn, this is chain homotopy equivalent to $C\{z; n \geq 0\}$ over $\Field[[U_1]]$, via the usual map taking the domain of the cone to $0$, and projecting the target with $U_2 \mapsto U_1$; cf. Lemma~\ref{lem:ap2}.

The combined homotopy equivalence from $\Ccint$ to $C\{z; n \geq 0\}$ is the desired map \eqref{eq:midrho}. We can write it as follows. The generators of $\Ccint$ over $\Field[[U_1]]$ are 
$$U_2^{m} V_2^n Y_2^a (\x \times r), \ \ m, n \geq 0, \ a \in \{0,1\}, \ \x \in \Ta \cap \Tb,\ r \in \{x_+,x_-\}.$$
 Then, we have
\begin{equation}
\label{eq:rhomidformula}
 \rho(U_2^{m} V_2^n Y_2^a (\x \times r)) = \begin{cases}
U_1^m \x & \text{if } a=0, r=x_+, \\
0 & \text{otherwise.}
\end{cases}.
\end{equation}

As for the homotopy $Z$ from \eqref{eq:commuteminus1}, we let it be the $\Field[[U_1]]$-module map defined on generators by
\begin{equation}
\label{eq:Zorro}
 Z(U_2^m V_2^n Y_2^a (\x \times r)) = \begin{cases}
mU_1^{m-1}\x & \text{if } a=0, r=x_-,\\
0 & \text{otherwise.}\end{cases}
\end{equation}

As an aside, observe that the Alexander grading of $\x \times r$ is 
\begin{equation}
\label{eq:Alex}
A(\x \times r) = \begin{cases}
A(\x) & \text{if  } r=x_+,\\
A(\x)-1& \text{if } r=x_-.\end{cases}
\end{equation} 

Let us now check that the left square in \eqref{eq:commuteminus1} commutes. In preparation for that, note that the generators of $\Chain^-(\Hyper^K, s)$ over $\Field[[U_1]]$ can be written
$$ U_2^m(\widetilde{\x \times r})^a_n, \  \x \in \Ta \cap \Tb, \ r \in \{x_+,x_-\}, \ n \in \Z_{\geq 0}, \ a \in \{0,1\}.$$
The left vertical map $\rho: \Chain^-(\Hyper^K, s) \to \Chain^-(\bar \Hyper^K, s)$ in \eqref{eq:commuteminus1} is constructed in \eqref{eq:rhoA}, by composing $\tilde P$ from \eqref{eq:tildeF} with a projection taking $U_2^m (\x \times r)$ to $U_1^m \x$ if $r=x_+$ and to zero if $r=x_-$. Therefore, we have
$$ \rho(U_2^m(\widetilde{\x \times r})^a_n)= \begin{cases}
U_1^{\max(A(\x)-s-n, 0)+m+n} \x &\text{if } a=0, r=x_+\\
0 & \text{otherwise.} 
\end{cases}$$

The first map in the bottom row in \eqref{eq:commuteminus1} is the inclusion $\I^{-K}_s$ defined in \eqref{eq:proj}, by the formula $\x \mapsto U_1^{\max(s-A(\x), 0)}\x$. However, that formula was based on the alternative description of $\Chain^-(\bar \Hyper^K)$ from Section~\ref{sec:alternative}, whereas here we use its description as a subcomplex of $\bar C\{w; n \geq 0\}$. The translation between the two is given by the inverse of Equation~\eqref{eq:mapsto}, namely $\x \mapsto U_1^{-\max(A(\x)-s, 0)}\x$. Composing these two formulas we obtain:
$$ \I^{-K}_s(\x) = U_1^{s-A(\x)}\x.$$

Therefore, we have
\begin{equation}
\label{eq:zabrze}
(\I_s^{-K}\circ \rho)(U_2^m(\widetilde{\x \times r})^a_n) = \begin{cases}
U_1^{\max(s-A(\x)+n, 0)+m} \x &\text{if } a=0,r=x_+,\\
0 & \text{otherwise.} 
\end{cases}
\end{equation}

With regard to the first map in the top row of \eqref{eq:commuteminus1}, by \eqref{eq:txna}, we have
\begin{equation}
\label{eq:Isk}
\I_s^{-K}(U_2^m(\widetilde{\x \times r})^a_n)= U_1^{\max(s+n-A(\x \times r), 0)} U_2^m V_2^n Y_2^a (\x \times r).
\end{equation}

Using \eqref{eq:rhomidformula} and \eqref{eq:Alex}, we get that $\I_s^{-K}\circ \rho$ is given by the same formula as \eqref{eq:zabrze}, so the left square in \eqref{eq:commuteminus1} commutes. In fact, the middle vertical map $\rho: \Ccint  \to C\{z; n \geq 0\}$, given by\eqref{eq:rhomidformula}, is filtered (with respect to the filtration $\FF$ from \eqref{eq:FF} on the domain and the Alexander filtration on the target). Its restriction to filtration levels $\leq s$ is exactly the left vertical $\rho$ map in \eqref{eq:commuteminus1}.

Next, we check that the right square in \eqref{eq:commuteminus1} commutes up to the homotopy $Z$, i.e., that Equation~\eqref{eq:checkZ} is satisfied.

To compute the top horizontal map $\Trans$ in this square, we use the formulas \eqref{eq:Trans1} and \eqref{eq:Trans0}, as well as the description of holomorphic disks in $\Hyper^K$ from Proposition~\ref{prop:ZemkeDisks}. With our choice of almost complex structure, the holomorphic disks in $\Hyper^K$ in classes $\phi$ with $n_{z_2}(\phi) \neq 0$ are:
\begin{itemize}
\item the disks from $\x \times x_+$ to $\x \times x_-$, which contain $z_2$ once and no other basepoints, and 
\item disks from $\x \times x_{\pm}$ to $\y \times x_{\pm}$, which corresponds to disks in the destabilized diagram from $\x$ to $\y$, and go over $z_2$ as many times as they went over $z_1$ in the destabilized diagram.
\end{itemize}
It is helpful to introduce the following notation. If the differential on $\bar C\{z; n \geq 0\}$ is given by
\begin{equation}
\label{eq:delxz}
 \del \x = \sum_{\y \in \Ta \cap \Tb} c(\x, \y) U_1^{n(\x, \y)} \y,
 \end{equation}
we set
$$ \del^{\Xi}\x =  \sum_{\y \in \Ta \cap \Tb} c(\x, \y) \frac{U_1^{n(\x, \y)}- U_2^{n(\x, \y)}}{U_1 - U_2} \y,$$
and
$$ \del'\x=  \sum_{\y \in \Ta \cap \Tb} c(\x, \y) n(\x, \y) U_1^{n(\x, \y)-1} \y.$$
Then, we have
$$ \Trans(U_2^m V_2^n Y_2^a (\x \times r)) = \begin{cases}
U_2^m(\x \times r)&\text{if } a=1,\\
U_2^m(\x \times x_- + (\del^{\Xi}\x)\times x_+) &\text{if } a=0, r=x_+,\\
(\del^{\Xi}\x) \times x_-  & \text{if }  a=0, r=x_-. 
\end{cases}$$

Moreover, in view of the identification \eqref{eq:zcone}, the rightmost vertical map in \eqref{eq:commuteminus1} is given by
$$\rho(U_2^m (\x \times r)) = \begin{cases}
U_1^m \x &\text{if } r=x_-,\\
0 &\text{if } r=x_+.
\end{cases}.$$

Therefore, we have
\begin{equation}
\label{eq:chorzow}
(\rho \circ \Xi )(U_2^mV_2^n Y_2^a (\x \times r)) = \begin{cases}
U_1^{m} \x &\text{if } a=1,r=x_- \text{ or } a=0,r=x_+\\
U_1^m (\del'\x) & \text{if } a=0,r=x_-,\\
0 & \text{if } a =1, r=x_+.
\end{cases}
\end{equation}

Combining this with \eqref{eq:rhomidformula}, we get a formula for the left hand side of \eqref{eq:checkZ}:
\begin{equation}
\label{eq:katowice}
U_2^mV_2^n Y_2^a (\x \times r) \mapsto \begin{cases}
U_1^{m} \x &\text{if } a=1,r=x_-,\\
U_1^m (\del'\x) & \text{if } a=0,r=x_-,\\
0 & \text{if } r=x_+. 
\end{cases}
\end{equation}

Let us compare this with the right hand side of \eqref{eq:checkZ}, which is the commutator $\del \circ Z + Z \circ \del$. Recall that $Z$ is given by Equation \eqref{eq:Zorro}. Note that, since $\Ccint$ is the mapping cone of $U_2 + U_1V_2$ on $\Ccint^0$, the differential $\del$ on $\Ccint$ consists of two parts: one coming from $\Ccint^0$, and one taking $Y_2$ to $U_2 + U_1V_2$. Further, in view of the identification \eqref{eq:Ccint0id}, the first part consist of a term coming from the differential on $\bar C \{z; n \geq 0\}^{U_1 \to U_2}$, and a term taking $\x \times x_-$ to $(V_2-1)\x \times x_+$. With this in mind, we can compute the effect of $\del \circ Z + Z \circ \del$ on a generator $U_2^mV_2^n Y_2^a (\x \times r)$. 

When $r=x_+$, note that $\del$ cannot take a generator containing $x_+$ into one containing $x_-$. Since $Z$ kills off all generators of the form $\y \times x_+$, we have
$$(\del \circ Z + Z \circ \del)(U_2^mV_2^n Y_2^a (\x \times x_+))=0,$$
in agreement with \eqref{eq:katowice}.

When $r=x_-$ and $a=0$, if we write $\del \x \in \bar C\{z; n \geq 0\}$ as in \eqref{eq:delxz}, then by using \eqref{eq:Ccint0id}, we obtain
$$\del(U_2^mV_2^n (\x \times x_-)) = (V_2-1)U_2^mV_2^n (\x \times x_+) + \sum_{\y \in \Ta \cap \Tb} c(\x, \y) U_2^{m+n(\x, \y)} V_2^n (\y \times x_-)$$
and hence
$$ (Z \circ \del)(U_2^mV_2^n (\x \times x_-)) = \sum_{\y \in \Ta \cap \Tb} c(\x, \y)(m+n(\x, \y)) U_1^{m+n(\x, \y)-1} \y.$$
On the other hand, we have
$$ (\del \circ Z)(U_2^mV_2^n (\x \times x_-)) = \sum_{\y \in \Ta \cap \Tb} c(\x, \y)m U_1^{m+n(\x, \y)-1} \y.$$
From here we get 
$$(\del \circ Z + Z \circ \del)(U_2^mV_2^n (\x \times x_+))=U_1^m (\del'\x),$$
as we expected from \eqref{eq:katowice}.

Finally, when $r=x_-$ and $a=1$, we have
\begin{multline*}
\del(U_2^mV_2^nY_2 (\x \times x_-)) = (U_2 + U_1V_2)U_2^mV_2^n(\x \times x_-) + (V_2-1)U_2^mV_2^nY_2 (\x \times x_+) \\
+ \sum_{\y \in \Ta \cap \Tb} c(\x, \y) U_2^{m+n(\x, \y)} V_2^n Y_2 (\y \times x_-),
\end{multline*}
so
$$
  (Z \circ \del)(U_2^mV_2^nY_2 (\x \times x_-)) = ((m+1)+m)U_1^m\x +  0 + 0 = U_1^m \x,$$
 since we work modulo two. We also have
$$ (\del \circ Z)(U_2^mV_2^nY_2 (\x \times x_-))=0.$$
Therefore,
$$(\del \circ Z + Z \circ \del)(U_2^mV_2^n Y_2 (\x \times x_+))= U_1^m \x,$$
again in agreement with \eqref{eq:katowice}.

This proves \eqref{eq:checkZ}, i.e., that the right square in  \eqref{eq:commuteminus1} commutes up to the given chain homotopy. 

Putting the diagrams \eqref{eq:commuteminus0} and \eqref{eq:commuteminus1} side by side, we obtain \eqref{eq:commuteminus}, and the proof is complete.
\end{proof}

\subsection{Maps between resolutions}
\label{sec:thetas}
Recall that one essential ingredient in the link surgery formula are the hyperboxes of strongly equivalent Heegaard diagrams. These are defined in Section~\ref{sec:hse}, for arbitrary (not necessarily link-minimal) diagrams. The construction uses $\Theta$ chain elements in generalized Floer complexes of the form $\Chain^-(\Tb, \Tbp, \zero) \cong \Am(\Tb, \Tbp, \zero)$ or $\Chain^-(\Ta, \Tap, \zero) \cong \Am(\Ta, \Tap, \zero)$, associated to link-minimal diagrams of an unlink in a connected sum of $(S^1 \times S^2)$'s. Furthermore, in Section~\ref{sec:hyperfloer}, we showed that a hyperbox of strongly equivalent Heegaard diagrams produces one of generalized Floer complexes of the form $\Am(\Ta, \Tb, \s)$; this used the polygon maps between $\Am$ complexes constructed in Section~\ref{sec:polygon}. 

For the link surgery formula with diagrams that are not link-minimal, we will also need to see how the same $\Theta$ chain elements produce hyperboxes made of the {\em resolutions} $\Chain^-(\Ta, \Tb, \s)$ of generalized Floer complexes. To do so, it suffices to construct polygon maps between these resolutions. We will explain this here in the simplest non-trivial case, for two strongly equivalent diagrams representing a knot $K \subset Y$, with two $w$ and two $z$ basepoints:
$$ \he = (\Sigma, \alphas, \betas, w_1, w_2, z_1, z_2), \ \ \ \he'=(\Sigma,  \alphas, \betas', w_1, w_2, z_1, z_2).$$
We assume that $\betas$ and $\betas'$ are strongly equivalent collections of curves, so together they form a (link-minimal, admissible) diagram representing a two-component unlink in a connected sum of $(S^1 \times S^2)$'s. Pick a cycle
$$ \Theta \in \Am(\Tb, \Tbp, \zero).$$
We can write $\Theta$ as a linear combination of intersection points,
$$ \Theta = \sum_{\a \in \Tb \cap \Tbp} c(\a) U_1^{u_1(\a)} U_2^{u_2(\a)} \a,$$
for $c(\a) \in \{0,1\}$ and $u_1(\a), u_2(\a) \geq 0$.
 
There are two Alexander filtrations on $\Tb \cap \Tbp$, call them $A_1$ and $A_2$, and we are interested in their sum, $A=A_1 + A_2$. The fact that $\Theta$ is in bi-filtration level $\leq \zero$ means that
$$A_1(\a) - u_1(\a) \leq 0,  A_2(\a) - u_2(\a) \leq 0\ \text{for all $\a$ with} \ c(\a)=1.$$
and hence
\begin{equation}
\label{eq:Aau}
 A(\a) - u_1(\a) - u_2(\a) \leq 0 \ \text{for all $\a$ with} \ c(\a)=1.
 \end{equation}

The $\Theta$ element induces a triangle map (cf. Section~\ref{sec:polygon}):
$$
 \Am(\Ta, \Tb, s) \to \Am(\Ta, \Tbp, s),$$

$$ \x \mapsto \sum_{\y \in \Ta \cap \Tbp} \sum_{\a \in \Tb \cap \Tbp} \sum_{{\substack{\phi \in \pi_2(\x, \a, \y) \\ \mu(\phi)=0}} } \#\M(\phi) \cdot c(\a) U_1^{n_{w_1}(\phi) + u_1(\a)} U_2^{n_{w_2}(\phi) + u_2(\a)} \y.
$$

Consider the resolutions $\Chain^-(\Ta, \Tb, s)$ of $\Am(\Ta, \Tb, s)$ and  $\Chain^-(\Ta, \Tbp, s)$ of  $\Am(\Ta, \Tbp, s)$, constructed as in Section~\ref{sec:knot4}. We claim that $\Theta$ also induces a chain map 
$$ \Chain^-(\Ta, \Tb, s) \to \Chain^-(\Ta, \Tbp, s).$$
 
To construct it, let us first define an analogue of the complex $\Ccint^0$ from Section~\ref{sec:knot4} for the $\beta$-$\beta'$ diagram. This complex, denoted $\Ccint^0(\Tb, \Tbp)$, is free over $\Ringbig = \Field[[U_1, U_2]][V_2]'$ with generators $\a \in \Tb \cap \Tbp$, and differential similar to \eqref{eq:delint}, counting holomorphic disks with coefficients
$$U_1^{n_{z_1}(\phi) + n_{z_2}(\phi)} V_2^{{n_{w_2}(\phi)}}.$$
 
Recall that in \eqref{eq:Fmap} we defined a map $P: \Ccint \to \Cint$, whose restriction to $\Ccint^0$  takes
$$ \x \mapsto U_1^{A(\x)-s} \x, \ \ V_2 \to U_1^{-1}U_2.$$

The idea is to ``pull back'' the element $\Theta$ under a map of this form (with $s=0$), that is, to consider
$$ \tTheta =  \sum_{\a \in \Tb \cap \Tbp} c(\a) U_1^{u_1(\a)+u_2(\a)-A(\a)} V_2^{u_2(\a)} \a \in \Ccint^0(\Tb, \Tbp).$$
Note that the condition \eqref{eq:Aau} implies that $U_1$ has nonnegative exponents. Further, it is easy to check that since $\Theta$ is a cycle, so is $\tTheta$.
 
There is a filtration $\FF$ on $ \Ccint^0(\Tb, \Tbp)$ similar to the one from Section~\ref{sec:knot4}, i.e.,
$$ \FF(\x) = -A(\x), \ \FF(U_1)=-1, \ \FF(U_2)=0, \ \FF(V_2)=1.$$
Note that all the terms in the expression for $\tTheta$ are in filtration degree $\FF \leq 0$.

We have a triangle map
$$ f: \Ccint(\Ta, \Tb) \otimes \Ccint^0(\Tb, \Tbp) \to \Ccint(\Ta, \Tbp)$$
given by
\begin{equation}
\label{eq:tTheta1}
f(\x \otimes \a) =  \sum_{\y \in \Ta \cap \Tbp}  \sum_{{\substack{\phi \in \pi_2(\x, \a, \y) \\ \mu(\phi)=0}} } \#\M(\phi) \cdot  U_1^{n_{z_1}(\phi) + n_{z_2}(\phi)} V_2^{n_{w_2}(\phi)} \y,
 \end{equation}
and equivariant with respect to the actions of $U_1, U_2, V_2$ and $Y_2$ (with $Y_2$ acting by $1$ on $\Ccint^0(\Tb, \Tbp)$). This map is filtered with respect to the filtrations $\FF$.

By plugging in the $\tTheta$ element and restricting to $\FF \leq 0$, we obtain the desired map $ \Chain^-(\Ta, \Tb, s) \to \Chain^-(\Ta, \Tbp, s)$, given by the formula
\begin{equation}
\label{eq:tTheta2}
 \x \to   \sum_{\a \in \Tb \cap \Tbp} c(\a) U_1^{u_1(\a)+u_2(\a)-A(\a)} V_2^{u_2(\a)} f(\x \otimes \a).
 \end{equation}

Higher polygon maps between the resolutions are constructed in a similar manner.

\subsection{A knot with many basepoints}
\label{sec:knots2p}
We now consider a knot $K$ represented by a Heegaard diagram $\Hyper^K$ with $2p$ basepoints, denoted (in order) $w_1, z_1, w_2, z_2, \dots, w_p, z_p$.

Just as for a knot with four basepoints we had to construct a resolution of $\Cint = C\{ w;  n_2 \geq 0, A(\x) - n_1 - n_2 \leq s\}$, we now start by constructing a resolution $\Ccint$ of the intermediate complex
$$ \Cint = C\{ w;  n_2, \dots, n_p \geq 0, A(\x) - n_1 - \ldots - n_p \leq s\},$$
with elements of the form
\begin{equation}
\label{eq:genCint}
\x_{\ns} := U_1^{A(\x) - s - |\ns|} U_2^{n_2} \dots U_p^{n_p} \x, \ \ \x\in \Ta \cap \T, \ \ns=(n_2, \dots n_p)  \in (\Z_{\geq 0})^{p-1},
\end{equation}
where we wrote $|\ns| = n_2 + \dots + n_p$. Every element of $\Cint$ can be written as an infinite sum
$$\sum_{\x \in \Ta \cap \Tb} \sum_{i=1}^{\infty} U_1^{a_i} U_2^{b^2_i} \dots U_n^{b^n_i} \x_{\ns_i},$$
for some sequence of distinct tuples $(a_i, b^2_i, \dots, b^n_i, \ns_i = (n_{i,2}, \dots, n_{i,p}))_{i \geq 1}$  satisfying a boundedness condition similar to \eqref{eq:bdd}, namely
 $$ \exists K > 0 \text{ such that } |\ns_i | - a_i < K, \text{ for all } i.$$

In the construction of $\Ccint$, we turn $\x_{\ns}$ into linearly independent elements over $\Field[[U_1, \dots, U_p]]$, by introducing new variables $V_2, \dots, V_p$, so that $\x_{\ns}$ corresponds to $V_2^{n_2}\dots V_p^{n_p} \x$. More precisely, we proceed as follows. Given a monomial $\ms$ in several variables, we will write $e_{\ms}(X)$ for the exponent of a variable $X$ in $\ms$. We define the subring $\Ringbig \subset  \Field[[U_1, U_2, \dots, U_p, V_2, \dots, V_p]]$ to consist of infinite sums of monomials $\ms$ such that 
$$ e_{\ms}(V_2) + \dots + e_{\ms}(V_p) - e_{\ms}(U_1) $$
is bounded above (over all $\ms$ in a sum). Note that $\Ringbig$ is a module over $\Field[[U_1, U_2]]$, and contains the ring $\Field[[U_1, U_2, \dots, U_p]][V_2, \dots, V_p]$. We will write
$$ \Ringbig=\Field[[U_1, U_2, \dots, U_p]][V_2, \dots, V_p]'.$$

Let $\Ccint^0$ be the complex freely generated over $\Ringbig$
by the intersection points $\x \in \Ta \cap \Tb$, with differential 
\begin{equation}
\label{eq:delknotp}
  \del \x = \sum_{\y \in \Ta \cap \Tb}  \sum_{\substack{\phi \in \pi_2(\x, \y) \\ \mu(\phi)=1} }  \# (\M(\phi)/\R) \cdot U_1^{n_{z_1}(\phi) + \ldots + n_{z_p}(\phi)} V_2^{n_{w_2}(\phi)} \dots  V_p^{n_{w_p}(\phi)} \y.\end{equation}
The fact that it is a complex ($\del^2=0$) follows from \cite[Theorem 5.5]{Links}: There are $2p$ boundary degenerations to account for, each containing a pair of basepoints, as follows:
$$ (w_1, z_1), (z_1, w_2), (w_2, z_2), \dots, (w_p, z_p), (z_p, w_1).$$
Since we keep track of all $z$ basepoints through the same variable $U_1$, we see that the contributions from boundary degenerations cancel in pairs.

Next, we construct $\Ccint$ from $\Ccint^0$ by taking iterated mapping cones for maps of the form $U_i + U_1 V_i$, that is, we introduce new variables $Y_2, \dots, Y_p$ satisfying 
$$Y_i^2=0, \ \ \ \del Y_i = U_i + U_1 V_i.$$
We then take the tensor product
$$\Ccint = \Ccint^0 \otimes_{\Ringbig} \bigotimes_{i=2}^p \Ybig_i,$$
where $\Ybig_i$ is the free $\Ringbig$-complex with two generators $1$ and $Y_i$, and with differential $ \del Y_i = U_i + U_1 V_i$.

In other words, $\Ccint$ is the dg module over the dga
$$ \Ringbig^Y = \Field[[U_1, \dots, U_p]][V_2, \dots, V_p, Y_2, \dots, Y_p]'/(Y_i^2=0, \del Y_i = U_i + U_1 V_i)$$
generated by $\x \in \Ta \cap \Tb$, with differential \eqref{eq:delknotp}.

We define a projection map $P: \Ccint \to \Cint$ by
\begin{equation}
\label{eq:Fxepsns}
P(V_2^{n_2}\dots V_p^{n_p} Y_2^{a_2} \dots Y_p^{a_p}\x) = \begin{cases}
 U_1^{A(\x) - s - |\ns|} U_2^{n_2} \dots U_p^{n_p} \x  & \text{for } a_2=\dots =a_p = 0, \\
 0 & \text{otherwise,}
\end{cases}
\end{equation}
and extending it linearly (to the allowed infinite sums) over $\Ring = \Field[[U_1, U_2, \dots, U_p]].$

\begin{lemma}
\label{lem:FCh}
The map $P$ is a chain map, and exhibits $\Ccint$ as a resolution of $\Cint$.
\end{lemma}

\begin{proof}
The fact that $P$ commutes with the differential follows as in the proof of Lemma~\ref{lem:Fchain}. 

To show that $P$ is a quasi-isomorphism, observe that $P$ is surjective, and denote its kernel by $G$. We seek to prove that $G$ is acyclic.

Note that $G$ is generated over $\Ringbig$ by the elements that contain non-trivial powers of some $Y_i$'s, as well as by the sums $(U_i + U_1 V_i)\x$, corresponding to relations between the generators $\x_\ns$ of $\Cint$; cf. \eqref{eq:genCint}. 

Consider an area filtration $\A$ on the generators of $G$, as follows. Let $\Pi$ be the space of periodic domains $\phi$ on the Heegaard diagram $\Hyper^K$ that satisfy 
$$ n_{z_1}(\phi) + \ldots + n_{z_p}(\phi)=0, \ n_{w_2}(\phi)= \dots = n_{w_p}(\phi)=0.$$

Equip the Heegaard surface with an area form such that all domains in $\Pi$ have zero area. Pick a base generator $\x_b \in \Ta \cap \Tb$, and set $\A(\x_b)=0$. Then, for any other generator of the form $U_1^{m} V_2^{n_2} \dots V_p^{n_p}\x$, pick a class $\phi \in \pi_2(\x, \x_b)$ with 
$$n_{z_1}(\phi) + \ldots + n_{z_p}(\phi)=m, \ n_{w_2}(\phi)=n_2,  \dots, n_{w_p}(\phi)=n_p$$
and let $\A(U_1^{m} V_2^{n_2} \dots V_p^{n_p}\x)$ be the area of $\phi$. Finally, let $\A$ be unchanged when we multiply by variables $Y_i$ or $U_i$ for $i \geq 2$. 

In view of Equation~\eqref{eq:delknotp}, we see that $\A$ defines a filtration on the complex $\Ccint$, and hence on its subcomplex $G$. The associated graded of $\A$ on $\Ccint$  is a direct sum, over $\x \in \Ta \cap \Tb$, of Koszul complexes for the regular sequence
$$ (U_2+U_1V_2,  U_3 + U_1 V_3, \dots, U_p + U_1 V_p)$$
over the ring $\Ringbig$. The homology of such a Koszul complex is concentrated in degree zero, with respect to the grading given by the sum of the exponents of the $Y_i$ variables. It follows that the same is true for the associated graded for $\A$ on $G$. Moreover, in degree zero, all the generators $(U_i + U_1 V_i)\x$ of $G$ are in the image of the differential $\del$ on the associated graded (because $\del Y_i = U_i + U_1 V_i$). We conclude that the associated graded of $G$ is acyclic, and from here that $G$ is acyclic.
\end{proof}

Next, as in the case of four basepoints, we define a filtration $\F$ on $\Cint$ by the negative of the exponent of $U_1$. Then $\Am(\Hyper^K, s)$ is the subcomplex of $\Cint$ given by $\F \leq 0$. 

We also define a filtration $\FF$ on $\Ccint$ by setting
$$ \FF(\x) = -A(\x), \ \FF(U_1)=-1, \ \FF(U_i)=0, \ \FF(V_i)=1, \ \FF(Y_i)=0,\ i \geq 2.$$
One can check that the map $P$ from \eqref{eq:Fxepsns} is filtered of degree $-s$ with respect to $\FF$ and $\F$. 

We define a resolution $\Chain^-(\Hyper^K, s)$ of $\Am(\Hyper^K,s)$ to be the subcomplex of $\Ccint$ corresponding to $\FF \leq s$. Thus, $\Chain^-(\Hyper^K, s)$ consists of infinite sums (with $\Ring$ coefficients, and a suitable boundedness condition) of the elements
$$ \tx_\ns^\as := U_1^{\max(s+|\ns|-A(\x), 0)} V_2^{n_2} \dots V_p^{n_p} Y_2^{a_2} \dots Y_p^{a_p} x,$$ 
where we denote $\ns = (n_2, \dots, n_p)$, $\as = (a_2, \dots, a_p)$. 

We let 
$$\tilde{P}: \Chain^-(\Hyper^K, s) \to \Am(\Hyper^K,s)$$
to be the restriction of $P$ from \eqref{eq:Fxepsns}. In terms of the generators $\tx_{\ns}^{\as}$, we have
$$ \tilde{P}(\tx_{\ns}^\as) = \begin{cases}
  U_1^{\max(A(\x) - s - |\ns|, 0)} U_2^{n_2} \dots U_p^{n_p} \x & \text{for } a_2 =\dots=a_p=0, \\
 0 & \text{otherwise.}
\end{cases}
$$
 
By the same arguments as in Lemma~\ref{lem:FCh}, we can show that $\tilde{P}$
exhibits $\Chain^-(\Hyper^K, s)$ as a resolution of $\Am(\Hyper^K, s)$. 

As in Section~\ref{sec:knot4}, we define the map $\Phi^K_s$ in the surgery formula as the composition
$$ \Chain^-(\Hyper^K, s) \xrightarrow{\tilde P} \Am(\Hyper^K, s) \hookrightarrow C\{w; n_1, \dots, n_p \geq 0\} \xrightarrow{\sim} \CFm(\Hyper^{\emptyset}),$$

Recall also that in Section~\ref{sec:knot4} we defined the map $\Phi^{-K}_s$ as the composition
$$ \Chain^-(\Hyper^K, s) \hookrightarrow \Ccint \xrightarrow{\Trans} C\{z; n_1, n_2 \geq 0\} \xrightarrow{\sim} \CFm(\Hyper^{\emptyset}),$$
with the last map again obtained from Heegaard moves. 

In the general case, we have a similar sequence, but instead of the transition map $\Trans$ we have a composition of $p-1$ maps
\begin{equation}
\label{eq:transitions}
 \Ccint = \Cc_0 \xrightarrow{\Trans_0} \Cc_1  \xrightarrow{\Trans_1} \dots  \xrightarrow{\Trans_{p-2}} \Cc_{p-1} = C\{z; n_1, \dots, n_p \geq 0\}.
 \end{equation}

Here, for $0 \leq j \leq p-1$, we consider the dga
$$ \Ringbig^Y_j = \Field[[U_1, \dots, U_p]][V_{j+2}, \dots, V_p, Y_{j+2}, \dots, Y_p]'/(Y_i^2=0, \del Y_i = U_i + U_1 V_i),$$
where the prime refers to allowing infinite sums of monomials $\ms$ in the given variables, provided that
$$ e_{\ms}(V_{j+2}) + \dots + e_{\ms}(V_p) - e_{\ms}(U_1)$$
is bounded above. The complex $\Cc_j$ is defined as the dg module over 
$ \Ringbig^Y_j$ generated by $\x \in \Ta \cap \Tb$, with differential 
$$ \del \x = \sum_{\y \in \Ta \cap \Tb}  \sum_{\substack{\phi \in \pi_2(\x, \y) \\ \mu(\phi)=1} }  \# (\M(\phi)/\R) \cdot U_1^{n_{z_1}(\phi)} \dots U_{j}^{n_{z_{j}}(\phi)}  U_{j+1}^{n_{z_{j+1}}(\phi) + \dots + n_{z_p}(\phi)} V_{j+2}^{ n_{w_{j+2}}(\phi)} \dots V_p^{n_{w_p}(\phi)}.
$$

The transition maps $\Trans_j: \Cc_j \to \Cc_{j+1}$ are constructed as follows. Recall that for $p=2$ we had Equations~\eqref{eq:Trans1} and \eqref{eq:Trans0}. In general, we set
\begin{equation}
\label{eq:Transjdiff}
 \Trans_j(V_{j+2}^n Y_{j+2}\x) =  \x.
 \end{equation}
and
\begin{multline}
\label{eq:Transjsame}
 \Trans_j(V_{j+2}^n \x) =  \sum_{\y \in \Ta \cap \Tb}  \sum_{\substack{\phi \in \pi_2(\x, \y) \\ \mu(\phi)=1} }  \# (\M(\phi)/\R) \cdot U_1^{n_{z_1}(\phi)} \dots U_{j+1}^{n_{z_{j+1}}(\phi)} \frac{U_{j+1}^{n_{z_{j+2}}(\phi)+ \dots + n_{z_p}(\phi)} - U_{j+2}^{n_{z_{j+2}}(\phi)+ \dots + n_{z_p}(\phi)}}{U_{j+1} - U_{j+2}} \\ \cdot 
 V_{j+3}^{ n_{w_{j+3}}(\phi)} \dots V_p^{n_{w_p}(\phi)} \y,
\end{multline}
then extend these maps equivariantly with respect to the action of the variables $U_1, \dots, U_p$,$V_{j+3}, \dots, V_p$, $Y_{j+3}, \dots, Y_p$.

\begin{lemma}
\label{lem:TransChain1}
$\Trans_j$ is a chain map.
\end{lemma}

\begin{proof}
This is entirely similar to the proof of Lemma~\ref{lem:TransChain0}, with the variables $U_{j+1}$ and $U_{j+2}$ playing the roles of $U_1$ and $U_2$. Furthermore, all holomorphic disks get counted with an additional factor of 
$$U_1^{n_{z_1}(\phi)} \dots U_{j}^{n_{z_{j}}(\phi)} V_{j+3}^{ n_{w_{j+3}}(\phi)} \dots V_p^{n_{w_p}(\phi)}.$$
\end{proof}

\subsection{A link with two components and eight basepoints}
\label{sec:linktwo}
As the next warm-up for the general case of setting up the hyperboxes, let us consider a link $L = L_1 \cup L_2$ with two components. Suppose this is represented by a complete system of hyperboxes $\Hyper$, where the initial diagram $\Hyper^L$ has four basepoints on each component:
$$ w_{1,1}, z_{1,1}, w_{1,2}, z_{1,2} \ \text{on} \ L_1,$$
and
$$ w_{2,1}, z_{2,1}, w_{2,2}, z_{2,2} \ \text{on} \ L_2.$$

The system $\Hyper$ consists of several diagrams and hyperboxes, as listed in Example~\ref{ex:hyperboxes2}. We seek to construct resolutions $\Chain^-(\Hyper^M, \s)$ of $\Am (\Hyper^M, \s)$, for each sublink $M \subseteq L$, and also maps 
$$\Phi_{\s}^{\orM}:  \Chain^-(\Hyper^{L'}, \s) \to \Chain^-(\Hyper^{L' - M}, \psi^{\orM}(\s)) \text{ for } M \subseteq L' \subseteq L,$$
associated to the hyperboxes $\Hyper^{\orL', \orM}$. 

The complexes $\Chain^-(\Hyper^{L_1}, \s)$, $\Chain^-(\Hyper^{\emptyset})=\CFm(\Hyper^{\emptyset})$ and the chain map associated to $\Hyper^{\orL_1, \pm \orL_1}$ are constructed just as in Section~\ref{sec:knot4}, with $L_1$ playing the role of $K$, except that when counting holomorphic disks and triangles, we also keep track of the basepoints on $L_2$, with coefficients
$$ U_{2,1}^{n_{w_{2,1}}(\phi)} U_{2,2}^{n_{w_{2,2}}(\phi)}.$$
Similar remarks apply to $\Chain^-(\Hyper^{L_2}, \s)$ and the map associated to $\Hyper^{\orL_2, \pm \orL_2}$.

Let us construct $\Chain^-(\Hyper^L, \s)$. By analogy with what we did in Section~\ref{sec:knot4}, we first define $\Ringbig = \Field[[U_{1,1}, U_{1,2}, U_{2,1}, U_{2,2}]][V_{1,2}, V_{2,2}]'$ to consist of power series made of monomials $\ms$ such that
$$ e_{\ms}(V_{1,2}) - e_{\ms}(U_{1,1}) \ \text{ and } \ e_{\ms}(V_{2,2}) - e_{\ms}(U_{2,1})$$
are bounded above. We let $\Ccint^0(\Hyper^L)$ to be freely generated by $\x \in \Ta \cap \Tb$ over $\Ringbig$, with differential counting holomorphic disks with coefficients
$$ U_{1,1}^{n_{z_{1,1}}(\phi) + n_{z_{1,2}}(\phi)} V_{1,2}^{n_{w_{1,2}}(\phi)} U_{2,1}^{n_{z_{2,1}}(\phi) + n_{z_{2,2}}(\phi)} V_{2,2}^{n_{w_{2,2}}(\phi)}.$$

The complex $\Ccint(\Hyper^L)$ is obtained from $\Ccint^0(\Hyper^L)$ by adjoining variables $Y_{1,2}, Y_{2,2}$ with relations
\begin{equation}
\label{eq:Y12}
 Y_{1,2}^2=0, \ \del Y_{1,2} = U_{1,2} + U_{1,1} V_{1,2},
 \end{equation}
 \begin{equation}
\label{eq:Y22}
 Y_{2,2}^2=0, \ \del Y_{2,2} = U_{2,2} + U_{2,1} V_{2,2}.
 \end{equation}

There are now two filtrations, $\FF_1$ and $\FF_2$, on $\Ccint(\Hyper^L)$, given by
$$ \FF_1(\x) = -A_1(\x), \ \FF_1(U_{1,1})=-1, \ \FF(U_{1,2})=0, \ \FF(V_{1,2})=1, \ \FF(Y_{1,2})=0$$
$$ \FF_1(U_{2,1})= \FF(U_{2,2})= \FF(V_{2,2})= \FF(Y_{2,2})=0$$
and similarly for $\FF_2$. For $\s = (s_1, s_2)$, we let $\Chain^-(\Hyper^L, \s)$ be the subcomplex of $\Ccint(\Hyper^L)$ in bi-filtration degrees $\F_ 1\leq s_1, \F_2 \leq s_2.$

Next, consider the one-dimensional hyperboxes $\Hyper^{L, \pm L_1}$. The corresponding maps
$$ \Phi^{\pm L_1}_{\s} : \Chain^-(\Hyper^L, \s) \to \Chain^-(\Hyper^{L_2}, \psi^{\pm L_1}(\s))$$
are constructed as in Section~\ref{sec:knot4}, except that at each step, we keep track of the basepoints on $L_2$ through coefficients
$$U_{2,1}^{n_{z_{2,1}}(\phi) + n_{z_{2,2}}(\phi)} V_{2,2}^{n_{w_{2,2}}(\phi)}$$
and we also keep the variable $Y_{2,2}$ with relations \eqref{eq:Y22}. Note that, as part of the construction, the $\Theta$ elements that are part of the data in $\Hyper^{L, \pm L_1}$ are used to define an equivalence from $\Chain^-(r_{\pm L_1}(\Hyper^L), \psi^{\pm L_1}(\s))$ to $\Chain^-(\Hyper^{L_2}, \psi^{\pm L_1}(\s))$. This equivalence is the composition of several triangle  maps, induced by the $\Theta$ elements as explained in Section~\ref{sec:thetas}.

The maps associated to $\Hyper^{L, \pm L_2}$ are defined in the same way.

The last ingredient in the surgery formula are chain homotopies associated to the four two-dimensional hyperboxes in $\Hyper$. We leave the constructions corresponding to $\Hyper^{L, L_1 \cup L_2}$, $\Hyper^{L, L_1 \cup -L_2}$ and $\Hyper^{L, (-L_1) \cup L_2}$ as an exercise for the reader, and focus on the most complicated case, that of the hyperbox $\Hyper^{L, -L_1 \cup -L_2}$. (These constructions are all particular cases of the general definition, which will be given in Section~\ref{sec:gens}.)

We seek to construct the map $\Phi^{-L}_{\s}$, where $\s=(s_1, s_2)$. This will be a composition 
$$ \Phi^{-L}_{\s} = D^{-L}_{p^{-L}(\s)} \circ \I^{-L}_{\s}.$$
Here, $\I^{-L}_{\s}$ is just the inclusion of $\Chain^-(\Hyper^L, \s)$ into $\Ccint(\Hyper^L)$. The descent map $D^{-L}_{\s}$ will be obtained by compressing a certain hyperbox, which we proceed to describe. Note that, in this case, $p^{-L}(\s) = (-\infty, -\infty)$ does not depend on $\s$, so we can simply write $D^{-L}$ or $D^{-L}_{\s}$.

As part of the data we have the two-dimensional hyperbox $\Hyper^{L, -L}$ of strongly equivalent Heegaard diagrams (representing $\emptyset \subset S^3$). This hyperbox starts at $r_{-L}(\Hyper^L)$, the diagram obtained from $\Hyper^L$ by deleting the $w$ basepoints, and ends at $\Hyper^{\emptyset}$. Suppose that this hyperbox (rectangle) is of size $\dd = (d_1, d_2)$, and let us denote by $\Hyper^{L, -L}_{\eps}$ the diagram at position $\eps = (\eps_1, \eps_2)$, with $\eps_i \in \{0, \dots, d_i\}$.

We also have the one-dimensional hyperbox $\Hyper^{L, -L_1}$ of size $d_1$, consisting of strongly equivalent Heegaard diagrams that represent $L_2$. This hyperbox starts at $r_{-L_1}(\Hyper^L)$, the diagram obtained from $\Hyper^L$ by deleting the $w$ basepoints on $L_1$ only. Further, by the compatibility relation between hyperboxes (cf. Definition~\ref{def:comp1}), when we apply $r_{-L_2}$ to the hyperbox $\Hyper^{L, -L_1}$, we should obtain the first row in the two-dimensional hyperbox $\Hyper^{L, -L}$. Thus, if we denote by $\Hyper^{L, -L_1}_{\eps_1}$ the diagram at position $\eps_1$, we have $r_{-L_2}(\Hyper^{L, -L_1}_{\eps_1}) = \Hyper^{L, -L}_{\eps_1, 0}$.

There is also a similar hyperbox $\Hyper^{L, -L_2}$, of size $d_2$, in the other direction. (Compare Example~\ref{ex:hyperboxes2} and Figure~\ref{fig:completeL}.) Its reduction with respect to $L_1$ is the first column in $\Hyper^{L, -L}$.

\begin{figure}
\begin{center}
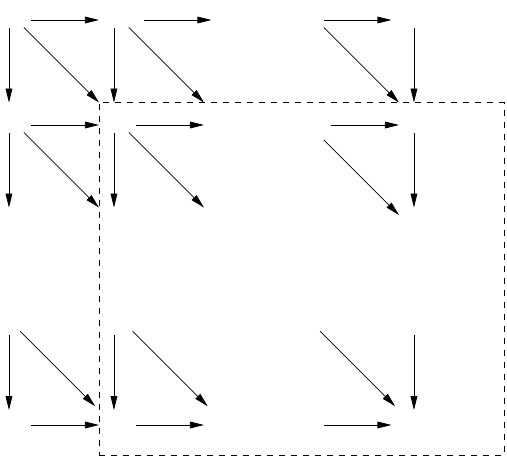
\end{center}
\caption {{\bf A hyperbox of chain complexes of size $(d_1 + 1, d_2 + 1)$.} This is the hyperbox $\Cc(\Hyper^{L,-L})$. The initial vertex is $\Ccint(\Hyper^L)$, and the rectangle bounded by dashed lines is the hyperbox obtained by applying $\CFm$ to the diagrams in $\Hyper^{L, -L}$.}
\label{fig:Enlarged}
\end{figure}

Note that from $\Hyper^{L, -L}$ we could get a hyperbox of chain complexes of size $(d_1, d_2)$, by applying $\CFm$ to each diagram $\Hyper^{L, -L}_\eps$, as in Section~\ref{sec:hyperfloer}. However, we want to include transition maps into this hyperbox, so we add an additional row and column in the beginning. The resulting hyperbox, denoted $\Cc(\Hyper^{L, -L})$ and pictured in Figure~\ref{fig:Enlarged}, is of size $(d_1 + 1, d_2 + 1)$. At its vertices we have complexes
$$ C^{\eps_1, \eps_2}= \begin{cases}
\Ccint(\Hyper^L) & \text{if } \eps_1=\eps_2=0,\\
\Ccint(\Hyper^{L, - L_1}_{\eps_1-1}) & \text{if } \eps_1 \neq 0, \eps_2=0,\\
\Ccint(\Hyper^{L, -L_2}_{\eps_2-1}) & \text{if } \eps_1 = 0, \eps_2 \neq 0,\\
\CFm(\Hyper^{L, -L}_{\eps_1-1, \eps_2-1}) & \text{if } \eps_1,\eps_2 \neq 0.
\end{cases}$$

We need to specify the chain maps along the horizontal and vertical segments in the hyperbox, and the chain homotopies along diagonals. In the bottom-right rectangle of size $(d_1, d_2)$, they are the usual maps in $\CFm(\Hyper^{L, -L})$, given by counting holomorphic triangles or quadrilaterals, with some vertices fixed at $\Theta$ elements that are part of the data in $\Hyper^{L, -L}$. Note that these holomorphic triangles and quadrilaterals, in classes $\phi$, are counted with coefficients
$$U_{1,1}^{n_{z_{1,1}}(\phi)} U_{1,2}^{n_{z_{1,2}}(\phi)} U_{2,1}^{n_{z_{2,1}}(\phi)} U_{2,2}^{n_{z_{2,2}}(\phi)}.$$

We are left to specify the chain maps and chain homotopies in the first row and column of the hyperbox $\Cc(\Hyper^{L, -L})$. 

\begin{figure}
\begin{center}
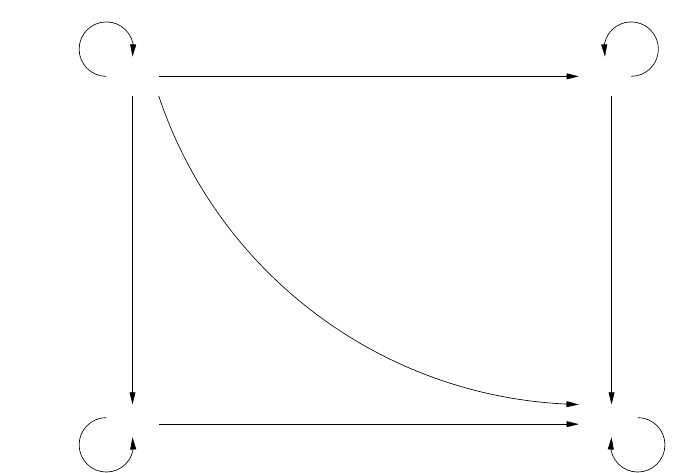
\end{center}
\caption {{\bf The first square in the hyperbox $\Cc(\Hyper^{L, -L})$.} Each arrow is a map counting holomorphic disks. We indicated the coefficients that keep track of the basepoints $z_{1,2}$ and $z_{2,2}$.}
\label{fig:Link4}
\end{figure}

Let us start by looking at the top-left square, pictured in Figure~\ref{fig:Link4}. Observe that, at its vertices, we have chain complexes associated to the same diagram $\Hyper^L$, in which we keep track of the basepoints in various ways. With regard to the basepoints $z_{1,1}$ and $z_{2,1}$, at all four vertices, the coefficients in the differential pick up factors of
$$ U_{1,1}^{n_{z_{1,1}}(\phi)} U_{2,1}^{n_{z_{2,1}}(\phi)}.$$
On the other hand, with regard to $z_{1,2}$ and $z_{2,2}$, these are counted differently, depending on the vertex. If we set
$$ m_1 := n_{z_{1,2}}(\phi), \ \ m_2 := n_{z_{2,2}}(\phi),$$
then $m_1$ appears as an exponent of either $U_{1,1}$ or $U_{1,2}$, and $m_2$ appears as an exponent of either $U_{2,1}$ or $U_{2,2}$ as indicated in Figure~\ref{fig:Link4}. We also have the $w$ basepoints, which appear as exponents of $V$ variables. In $C^{0,0}=\Ccint(\Hyper^L)$ we have
$$ V_{1,2}^{n_{w_{1,2}}(\phi)} V_{2,2}^{n_{w_{2,2}}(\phi)},$$
then in $C^{1,0}$ we set $V_{1,2}=1$ (and keep $V_{2,2}$ as before). In $C^{0,1}$ we set $V_{2,2}=1$, and in $C^{1,1}$ both $V$ variables are set to $1$. A similar pattern holds for the $Y$ variables: in $C^{0,0}$ we have both $Y_{1,2}$ and $Y_{2,2}$, with relations \eqref{eq:Y12} and \eqref{eq:Y22}; in $C^{1,0}$ we drop $Y_{1,2}$, in $C^{0,1}$ we drop $Y_{2,2}$, and in $C^{1,1}$ we drop both.

With regard to the chain maps along the edges of the square in Figure~\ref{fig:Link4}, these are transition maps similar to the map $\Trans$ from Section~\ref{sec:knot4}. For example, the top horizontal map
$$ D^{1,0}_{0,0}: C^{0,0} \to C^{1,0}$$
is defined as follows. This map is equivariant with respect to the action of the $U$ and $V$ variables, as well to the action of $Y_{2,2}$ (but not $Y_{1,2}$, which is a variable that does not appear in the target complex $C^{1,0}$). The map  $D^{1,0}_{0,0}$ takes
\[
Y_{1,2} \x \mapsto \x, \]
\[
\x \mapsto \sum_{\y \in \Ta \cap \Tb}  \sum_{\substack{\phi \in \pi_2(\x, \y) \\ \mu(\phi)=1} }  \# (\M(\phi)/\R) \cdot U_{1,1}^{n_{z_{1,1}}(\phi)} \frac{U_{1,1}^{n_{z_{1,2}}(\phi)} - U_{1,2}^{n_{z_{1,2}}(\phi)}}{U_{1,1} - U_{1,2}}  U_{2,1}^{n_{z_{2,1}}(\phi) + n_{z_{2,2}}(\phi)} \y.
\]
The bottom horizontal map, $D^{1,0}_{0,1}$, is defined in the same way, with the only difference that when counting holomorphic disks, the quantity $n_{z_{2,2}}(\phi)$ appears as an exponent of $U_{2,2}$ instead of $U_{2,1}$. The two vertical chain maps have similar definitions, with the roles of $L_1$ and $L_2$ (and the corresponding variables) switched.

Finally, the chain homotopy along the diagonal of the square,
$$ D^{1,1}_{0,0} : C^{0,0} \to C^{1,1}$$
is taken to be equivariant with respect to the $U$ and $V$ variables, to take the $Y_{1,2}$ and $Y_{2,2}$ variables to zero, and to satisfy
$$ \x \mapsto  \sum_{\y \in \Ta \cap \Tb}  \sum_{\substack{\phi \in \pi_2(\x, \y) \\ \mu(\phi)=1} }  \# (\M(\phi)/\R) \cdot U_{1,1}^{n_{z_{1,1}}(\phi)} \frac{U_{1,1}^{n_{z_{1,2}}(\phi)} - U_{1,2}^{n_{z_{1,2}}(\phi)}}{U_{1,1} - U_{1,2}}  U_{2,1}^{n_{z_{2,1}}(\phi)}  \frac{U_{2,1}^{n_{z_{2,1}}(\phi)} - U_{2,2}^{n_{z_{2,2}}(\phi)}}{U_{2,1} - U_{2,2}}  \y.$$

\begin{lemma}
\label{lem:D2first}
With these definitions, the total differential $D$ on the first square in the hyperbox $\Cc(\Hyper^{L, -L})$ satisfies $D^2=0$.
\end{lemma}

\begin{proof}
The one interesting case to check is when we apply $D^2$ to a generator $\x \in C^{0,0}$ (without any $Y$ variables). Then, $\y \in \Ta \cap \Tb$ (viewed as an element in $C^{1,1}$) appears in $D^2(\x)$ with a sum of coefficients that come from pairs of holomorphic disks, in classes $\phi \in \pi_2(\x, \a), \psi \in \pi_2(\a, \y)$, over all $\a \in \Ta \cap \Tb$. For each pair $(\phi, \psi)$, there are four such contributions, from each of the four terms in
\begin{equation}
\label{eq:Dsquare}
 D^2=D^{0,1}_{1,0} \circ D^{1,0}_{0,0} + D^{1,0}_{0,1} \circ D^{0,1}_{0,0} + D^{1,1}_{0,0} \circ D^{0,0}_{0,0} + D^{0,0}_{1,1} \circ D^{1,1}_{0,0}.
 \end{equation}
In other words, we go from the initial to the final vertex in Figure~\ref{fig:Link4} by either following the edges, or along the diagonal combined with a differential (self-map) on either the initial or the final vertex.

We claim that these terms cancel in pairs. To see this, we need to track the powers of the $U$ variables in each of them. In all the terms we have factors of the form
$$ U_{1,1}^{n_{z_{1,1}}(\phi * \psi)} U_{2,1}^{n_{z_{2,1}}(\phi * \psi)}.$$
We will ignore these, and focus on the appearances of $n_{z_{1,2}}$ and $n_{z_{2,2}}$. These are indicated in Figure~\ref{fig:Link4}, for each map. Since we now have two different classes $\phi$ and $\psi$, let us denote
$$m_1 := n_{z_{1,2}}(\phi), \ \ m_2 := n_{z_{2,2}}(\phi),$$
$$m'_1 := n_{z_{1,2}}(\psi), \ \ m'_2 := n_{z_{2,2}}(\psi).$$
For each pair of disks in classes $(\phi, \psi)$, we get the following contributions to the coefficient of $\y$ in $\del^2 \x$, corresponding to the four terms in \eqref{eq:Dsquare}:
\begin{align*}
& \frac{U_{1,1}^{m_1} - U_{1,2}^{m_1}}{U_{1,1} - U_{1,2}} U_{2,1}^{m_2}  U_{1,2}^{m'_1} \frac{U_{2,1}^{m'_2}- U_{2,2}^{m'_2}}{U_{2,1}-U_{2,2}} + U_{1,1}^{m_1}\frac{U_{2,1}^{m_2}- U_{2,2}^{m_2}}{U_{2,1}-U_{2,2}}  \frac{U_{1,1}^{m'_1}-U_{1,2}^{m'_1}}{U_{1,1}-U_{1,2}} U_{2,2}^{m'_2}    \\
+ \  & U_{1,1}^{m_1} U_{2,1}^{m_2} \frac{U_{1,1}^{m'_1}-U_{1,2}^{m'_1}}{U_{1,1}-U_{1,2}} \frac{U_{2,1}^{m'_2}- U_{2,2}^{m'_2}}{U_{2,1}-U_{2,2}} +  \frac{U_{1,1}^{m_1}-U_{1,2}^{m_1}}{U_{1,1}-U_{1,2}} \frac{U_{2,1}^{m_2}- U_{2,2}^{m_2}}{U_{2,1}-U_{2,2}} U_{1,2}^{m'_1} U_{2,2}^{m'_2}.
\end{align*}
We obtain a fraction whose denominator is $(U_{1,1}-U_{1,2})(U_{2,1}-U_{2,2})$, and whose numerator is a sum of $16$ terms. Twelve of those terms cancel one another in pairs, and the remaining four terms can be combined into the product
$$ (U_{1,1}^{m_1+m_1'} - U_{1,2}^{m_1+m_1'})(U_{2,1}^{m_2+m'_2}- U_{2,2}^{m_2+m'_2}).$$
Observe that $m_1 + m_1'= n_{z_{1,2}}(\phi * \psi)$ and $m_2+m_2' = n_{z_{2,2}}(\phi * \psi)$. Thus, we are left with a count of pairs of holomorphic disks, where the exponents depend only on the combined class $\phi * \psi$. By the usual arguments in Floer theory, these pairs of holomorphic disks form the boundary of a one-dimensional moduli space of holomorphic disks, in the class $\phi * \psi$. Since the coefficients only depend on $\phi * \psi$, it follows that the relevant count of pairs of disks is zero. The contributions from the boundary degenerations cancel out in pairs, because the $w$'s do not appear in our expressions. We conclude that $D^2=0$, as desired.
\end{proof}

This completes the description of the top-left square in the hyperbox $\Cc(\Hyper^{L, -L})$. We now proceed to describe a typical square in the first row of that hyperbox. (The squares in the first column will be described similarly.)

\begin{figure}
\begin{center}
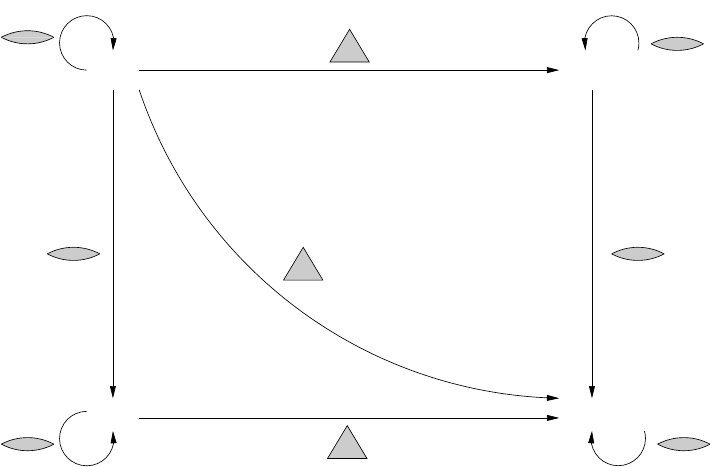
\end{center}
\caption {{\bf A further square in the first row of the hyperbox $\Cc(\Hyper^{L, -L})$.} Each arrow is a map counting holomorphic disks or triangles, as indicated by the nearby picture. We also wrote the coefficients that involve the quantities $m=n_{z_{2,2}}(\phi)$ and $u_{\a}=u_{2,2}(\a)-A_{2,2}(\a)$.}
\label{fig:Link4e}
\end{figure}

Consider the square in the first row that starts at $C^{\eps_1,0}$. This square is pictured in Figure~\ref{fig:Link4e}. It consists of data coming from two, strongly equivalent Heegaard diagrams for $L_2$:  namely, $\Hyper^{L, -L_1}_{\eps_1}$ and $\Hyper^{L, -L_1}_{\eps_1+1}$. Without loss of generality, let us suppose these differ by a strong equivalence among the beta curves; we denote the sets of curves on $\Hyper^{L, -L_1}_{\eps_1}$ by $(\alphas, \betas)$, and those on $\Hyper^{L, -L_1}_{\eps_1+1}$ by $(\alphas, \betas')$. Moreover, note that both of these diagrams have the same six basepoints: the free basepoints $z_{1,1}$ and $z_{1,2}$, and the four basepoints on $L_2$.

With this notation, the two complexes on the top edge of Figure~\ref{fig:Link4e} are
$$ C^{\eps_1, 0} = \Ccint(\Ta, \Tb), \ \ C^{\eps_1+1, 0} = \Ccint(\Ta, \Tbp).$$
The respective differentials count holomorphic disks with coefficients
\begin{equation}
\label{eq:topedge}
U_{1,1}^{n_{z_{1,1}}(\phi)} U_{1,2}^{n_{z_{1,2}}(\phi)} U_{2,1}^{n_{z_{2,1}}(\phi) + n_{z_{2,2}}(\phi) }  V_{2,2}^{n_{w_{2,2}}(\phi)},
\end{equation} 
and also have terms from tensoring with a mapping cone complex $\Ybig_{2,2}$ (with a variable $Y_{2,2}$). 

By contrast, the complexes on the bottom edge of Figure~\ref{fig:Link4e} are
$$ C^{\eps_1, 1} = \CFm(\Ta, \Tb), \ \ C^{\eps_1+1, 1} = \CFm(\Ta, \Tbp).$$
They don't involve the basepoints $w_{1,2}$ and $w_{2,2}$, and have no variables $V_{2,2}$ or $Y_{2,2}$. Their differentials count holomorphic disks with coefficients
\begin{equation}
\label{eq:bottomedge}
U_{1,1}^{n_{z_{1,1}}(\phi)} U_{1,2}^{n_{z_{1,2}}(\phi)} U_{2,1}^{n_{z_{2,1}}(\phi)} U_{2,2}^{n_{z_{2,2}}(\phi)}.
\end{equation}

The vertical maps in Figure~\ref{fig:Link4e} are transition maps defined by the following formulas:
\begin{equation}
\label{eq:verticalY}
Y_{2,2} \x \mapsto \x, 
\end{equation}
\begin{equation}
\label{eq:verticalnoY}
\x \mapsto \sum_{\y \in \Ta \cap \Tb}  \sum_{\substack{\phi \in \pi_2(\x, \y) \\ \mu(\phi)=1} }  \# (\M(\phi)/\R) \cdot U_{1,1}^{n_{z_{1,1}}(\phi)} U_{1,2}^{n_{z_{1,2}}(\phi)} U_{2,1}^{n_{z_{2,1}}(\phi)} \frac{U_{2,1}^{n_{z_{2,2}}(\phi)} - U_{2,2}^{n_{z_{2,2}}(\phi)}}{U_{2,1} - U_{2,2}} \y.
\end{equation}

The horizontal maps in Figure~\ref{fig:Link4e} are triangle maps induced by the strong equivalence from $(\alphas, \betas)$ to $(\alphas, \betas')$. They are specified by the choice of a cycle 
$$ \Theta \in \Am(\Tb, \Tbp, \zero)$$
that represents the maximal degree element in the homology of $\Am(\Tb, \Tbp, \zero)$. Note that the diagram formed by the $\beta$ and $\beta'$ curves represents an unlink in a connected sum of $(S^1 \times S^2)$'s, and this diagram is link minimal. The diagram contains the basepoints $z_{1,1}$ and $z_{1,2}$ from $L_1$, and all the four basepoints coming from $L_2$. The unlink consists of two components, one specified by $w_{2,1}$ and $z_{2,1}$, the other by $w_{2,2}$ and $z_{2,2}$. As such, the intersection points $\a \in \Tb \cap \Tbp$ admit two Alexander gradings $A_{2,1}$ and $A_{2,2}$, whose sum is denoted $A_2$.

We can write $\Theta$ as a linear combination of intersection points:
$$ \Theta = \sum_{\a \in \Tb \cap \Tbp} c(\a) \prod_{i_1, i_2 \in \{1,2\}} U_{i_1,i_2}^{u_{i_1,i_2}(\a)}  \a,$$
for $c(\a) \in \{0,1\}$ and some nonnegative integers $u_{i_1, i_2}(\a)$.

By analogy with the construction of the triangle maps in Section~\ref{sec:thetas}, cf. Equations~\eqref{eq:tTheta1} and \eqref{eq:tTheta2}, we define the map along the top edge of Figure~\ref{fig:Link4e} by the formula
\begin{multline}
 \x \mapsto \sum_{\y \in \Ta \cap \Tbp} \sum_{\a \in \Tb \cap \Tbp} \sum_{\substack{\phi \in \pi_2(\x, \a, \y) \\ \mu(\phi)=0} }  \# \M(\phi) \cdot c(\a) U_{1,1}^{n_{z_{1,1}}(\phi)+u_{1,1}(\a)} U_{1,2}^{n_{z_{1,2}}(\phi)+u_{1,2}(\a)}\\
 \cdot U_{2,1}^{n_{z_{2,1}}(\phi) + n_{z_{2,2}}(\phi) + u_{2,1}(\a) + u_{2,2}(\a)-A_2(\a)
  } 
 V_{2,2}^{n_{w_{2,2}}(\phi)+u_{2,2}(\a)} \y.
\end{multline}

On the other hand, for the map along the bottom edge, we use the image of $\Theta$ under the inclusion into the complex that uses the $z$ basepoints only; cf. Remark~\ref{rem:switch}. Using \eqref{eq:proj3}, this element can be written
$$  \sum_{\a \in \Tb \cap \Tbp} c(\a) \prod_{i \in \{1,2\}} U_{1,i}^{u_{1,i}(\a)} \prod_{i \in \{1,2\}} U_{2,i}^{u_{2,i}(\a)-A_{2,i}(\a)}  \a$$
By analogy with \eqref{eq:bottomedge}, we let
\begin{multline}
 \x \mapsto \sum_{\y \in \Ta \cap \Tbp} \sum_{\a \in \Tb \cap \Tbp} \sum_{\substack{\phi \in \pi_2(\x, \a, \y) \\ \mu(\phi)=0} }  \# \M(\phi) \cdot c(\a) U_{1,1}^{n_{z_{1,1}}(\phi)+u_{1,1}(\a)} U_{1,2}^{n_{z_{1,2}}(\phi)+u_{1,2}(\a)} \\
 \cdot U_{2,1}^{n_{z_{2,1}}(\phi) + u_{2,1}(\a)-A_{2,1}(\a)} U_{2,2}^{n_{z_{2,2}}(\phi)  + u_{2,2}(\a) -A_{2,2}(\a)}   \y.
 \end{multline}

Finally, for the chain homotopy along the diagonal of Figure~\ref{fig:Link4e}, by analogy with \eqref{eq:verticalY} and \eqref{eq:verticalnoY}, we let
\begin{equation}
\label{eq:diagonalY}
Y_{2,2} \x \mapsto 0, 
\end{equation}
\begin{multline}
\label{eq:diagonalnoY}
\x \mapsto \sum_{\y \in \Ta \cap \Tb}  \sum_{\substack{\phi \in \pi_2(\x, \y) \\ \mu(\phi)=1} }  \# (\M(\phi)/\R) \cdot c(\a) U_{1,1}^{n_{z_{1,1}}(\phi)} U_{1,2}^{n_{z_{1,2}}(\phi)} U_{2,1}^{n_{z_{2,1}}(\phi)+  u_{2,1}(\a)-A_{2,1}(\a)} \\
\cdot \frac{U_{2,1}^{n_{z_{2,2}}(\phi)+u_{2,2}(\a) -A_{2,2}(\a)} - U_{2,2}^{n_{z_{2,2}}(\phi)+u_{2,2}(\a) -A_{2,2}(\a)}}{U_{2,1} - U_{2,2}} \y.
\end{multline}

\begin{lemma}
With these definitions, the total differential $D$ on the square pictured in Figure~\ref{fig:Link4e} satisfies $D^2=0$.
\end{lemma}

\begin{proof}
Just as in Lemma~\ref{lem:D2first}, the interesting case to check is when we apply $D^2$ to a generator $\x$ in the initial complex $C^{\eps_1, 0}$ (without the $Y_{2,2}$ variable). The result is a 
sum of counts of pairs consisting of a holomorphic bigon and a holomorphic triangle. There are two cases, according to whether the triangle or the bigon comes first in the composition.

Let us consider the case when we first have a triangle (in a class $\phi$), and then a bigon (in a class $\psi$). Such contributions come from going along the top edge, and then the right edge in Figure~\ref{fig:Link4e}, and also from going along the diagonal and then using the differential on the final complex $C^{\eps_1+1, 1}$. In the coefficients of some $\y$ in $D^2\x$, we always obtain factors of the form
$$ c(\a) U_{1,1}^{n_{z_{1,1}}(\sigma)+u_{1,1}(\a)} U_{1,2}^{n_{z_{1,2}}(\sigma)+u_{1,2}(\a)} U_{2,1}^{n_{z_{2,1}}(\sigma) + u_{2,1}(\a)-A_{2,1}(\a)},$$
where $\sigma$ is the combined class $\phi * \psi$. These factors get multiplied by factors having to do with the basepoint $z_{2,2}$, as shown in Figure~\ref{fig:Link4e}. If we denote $m=n_{z_{2,2}}(\phi)$, $m'=n_{z_{2,2}}(\psi)$, and $u_{\a}=u_{2,2}(\a)-A_{2,2}(\a)$, the sum of the latter factors is
$$U_{2,1}^{m+u_{\a}} \frac{U_{2,1}^{m'} - U_{2,2}^{m'}}{U_{2,1}-U_{2,2}} + \frac{U_{2,1}^{m+u_{\a}} - U_{2,2}^{m+u_{\a}}}{U_{2,1}-U_{2,2}} U_{2,2}^{m'}= \frac{U_{2,1}^{m+m'+u_{\a}} - U_{2,2}^{m+m'+u_{\a}}}{U_{2,1}-U_{2,2}}.$$ 
Note that $m+m'=n_{z_{2,2}}(\phi * \psi)$, so the result only depends on the combined class $\sigma= \phi * \psi$. 

The case when we first have a bigon and then a triangle is similar, and we get the same coefficients in front of $\y$, with dependence only on $\sigma$.

Thus, it is natural to consider the one-dimensional moduli space of triangles with vertices at $\x, \a$ and $\y$, in the class $\sigma$, and weigh it by the same factor
\begin{multline}
c(\a) U_{1,1}^{n_{z_{1,1}}(\sigma)+u_{1,1}(\a)} U_{1,2}^{n_{z_{1,2}}(\sigma)+u_{1,2}(\a)} U_{2,1}^{n_{z_{2,1}}(\sigma) + u_{2,1}(\a)-A_{2,1}(\a)} \\
\cdot \frac{U_{2,1}^{n_{z_{2,2}}(\sigma)+u_{2,2}(\a)-A_{2,2}(\a)} - U_{2,2}^{n_{z_{2,2}}(\sigma)+u_{2,2}(\a)-A_{2,2}(\a)}}{U_{2,1}-U_{2,2}}.
\end{multline}

The boundary of this moduli space represents pairs of bigons and triangles, including the two cases described above. The other possibility involves a triangle with vertices at $\x, \b$ and $\y$ (for some $\b \in \Tb \cap \Tbp$), combined with a bigon from $\a$ to $\b$. However, the fact that $\Theta$ is a cycle in $\Am(\Tb, \Tbp, \zero)$ implies that those contributions sum up to zero. It follows that the coefficient of any $\y$ in $D^2\x$ vanishes.
\end{proof}

This completes the description of the hyperbox $\Cc(\Hyper^{L, -L})$. By compressing this hyperbox, we obtain the descent map
$$ D^{-L}: \Ccint(\Hyper^L) \to \CF(\Hyper^{\emptyset}).$$
After pre-composing this with the inclusion $\I^{-L}_{\s}: \Chain^-(\Hyper^L, \s) \to \Ccint(\Hyper^L)$, we arrive at the map
$$ \Phi_{\s}^{-L}:  \Chain^-(\Hyper^L, \s) \to  \CF(\Hyper^{\emptyset}),$$
which can be used as a chain homotopy in the surgery complex.

\subsection{The general statement}
\label{sec:gens}
By building up on the particular cases presented in Sections~\ref{sec:knot4}, \ref{sec:knots2p} and \ref{sec:linktwo}, we are now ready to give the general statement of the link surgery formula, Theorem~\ref{thm:FirstSurgery}. In other words, we plan to define the complex $\C^-(\Hyper, \Lambda)$, where $\Hyper$ is an arbitrary complete system of hyperboxes for the link $\orL \subset Y$, and $\Lambda$ is a framing. Let $L_1, \dots, L_{\ell}$ be the components of $\orL$. 

The surgery complex $\C^-(\Hyper, \Lambda)$ has roughly the same form as in the link-minimal case presented in Section~\ref{subsec:surgery}:
 $$\C^-(\Hyper, \Lambda) = \bigoplus_{M \subseteq L} \prod_{\s \in \H(L)}  \Chain^-(\Hyper^{L - M}, \psi^{M}(\s) ),$$
 with differential
  \begin {eqnarray*} 
\D^-(\s, \x) &=& \sum_{N \subseteq L - M} \sum_{\orN \in \Omega(N)} (\s + \Lambda_{\orL, \orN}, \Phi^{\orN}_{\psi^{M}(\s)}(\x)) \\
&\in&  \bigoplus_{N \subseteq L - M} \bigoplus_{\orN \in \Omega(N)}  \Chain^-(\Hyper^{L-M-N}, \psi^{M \cup \orN} (\s))   \subseteq \C^-(\Hyper, \Lambda).
\end {eqnarray*}

Our job is to define the complexes $\Chain^-(\Hyper^{L - M}, \psi^{M}(\s) )$ and the maps $\Phi^{\orN}_{\psi^{M}(\s)}$ in the general setting. It suffices to describe $\Chain^-(\Hyper^L, \s)$ and $\Phi^{\orM}_{\s}$; then, the same definitions will apply to the complexes for the sublinks of $L$, and to the maps having these as domains.

To define $\Chain^-(\Hyper^L, \s)$, we denote the basepoints on the component $L_i$, in order, by
$$ w_{i,1}, z_{i,1}, w_{i,2}, z_{i,2}, \dots, w_{i, p_i}, z_{i, p_i}.$$
We also allow some free basepoints, denoted
$$ w_{\circ, 1}, w_{\circ, 2}, \dots, w_{\circ, p}.$$

Consider first the dga
$$ \Ringbig^{Y} = \Field[[(U_{i,j})_{\substack{1\leq i \leq \ell, \\ 1 \leq j \leq p_i}} (U_{\circ,j})_{1\leq j \leq p}]][(V_{i,j}, Y_{i,j})_{\substack{1\leq i \leq \ell, \\ 2 \leq j \leq p_i}}]'/(Y_{i,j}^2=0, \del Y_{i,j} = U_{i,j} + U_{i,1}V_{i,j}),$$
where the prime refers to allowing infinite sums of monomials $\ms$ such that, for each $i=1, \dots, \ell$, we have that
$$e_{\ms}(V_{i,2}) + \dots + e_{\ms}(V_{i, p_i}) - e_{\ms}(U_{i, 1})$$
is bounded above. 

We define the complex $\Ccint(\Hyper^{L})$ as the dg module over $ \Ringbig^{Y}$ 
generated by $\x \in \Ta \cap \Tb$, with differential 
\begin{equation}
\label{eq:Ccintgen}
 \del \x = \sum_{\y \in \Ta \cap \Tb}  \sum_{\substack{\phi \in \pi_2(\x, \y) \\ \mu(\phi)=1} }  \# (\M(\phi)/\R) \cdot \prod_{i=1}^\ell  U_{i,1}^{n_{z_{i,1}}(\phi) + \ldots + n_{z_{i,p_i}}(\phi)} V_{i,2}^{n_{w_{i,2}}(\phi)} \dots  V_{i,p_i}^{n_{w_{i,p_i}}(\phi)} \cdot \prod_{j=1}^p U_{\circ,j}^{n_{w_{\circ, j}}(\phi)} \y.
 \end{equation}

This complex admits filtrations $\FF_i$, for $1 \leq i \leq \ell$, such that
$$ \FF_i(\x) = -A_i(\x), \ \FF_i(U_{i,1})=-1, \ \FF(V_{i,j})=1, \ j \geq 2,$$
and $\FF_i$ takes all the other $U$, $V$ and $Y$ variables to zero (i.e, the action of those variables respects $\FF_i$).

For $\s=(s_1, \dots, s_\ell) \in \H(L)$, we define $\Chain^-(\Hyper^L, \s)$ to be the filtered part of $\Ccint(\Hyper^{L})$ given by
$$ \FF_i \leq s_i, \ i=1, \dots, \ell.$$

Let $\Ring$ be the power series ring freely generated over $\Field$ by the $U$ variables (without $V$'s or $Y$'s). The complex $\Chain^-(\Hyper^L, \s)$ admits an $\Ring$-linear projection $\tilde{P}$ to $\Am(\Hyper^L, \s)$, given by 
$$ \tilde{P}\Bigl(\prod_{i,j} V_{i,j}^{n_{i,j}}  Y_{i,j}^{a_{i,j}} \x \Bigr) = 
\begin{cases}
\prod_{i=1}^{\ell} U_{i,1}^{A_i(\x) - s_i - \sum_{j=2}^{p_i} n_{i,j}} U_{i,2}^{n_{i,2}} \dots U_{i,p_i}^{n_{i,p_i}} \x  & \text{if all } a_{i,j} = 0, \\
 0 & \text{otherwise.}\end{cases}$$
The proofs of Lemmas~\ref{lem:AResolution} and \ref{lem:FCh} extend to our setting, showing that $\tilde{P}$ exhibits $\Chain^-(\Hyper^L, \s)$ as a resolution of $\Am(\Hyper^L, \s)$, over the ring $\Ring$.

Next, let $M \subseteq L$ a sublink, with an orientation $\orM$. As usual, $I_+(\orL, \orM)$ (resp. $I_-(\orL, \orM)$) denotes the set of indices $i$ such that the component $L_i$ is in $M$ and its orientation induced from $\orL$ is the same as (resp. opposite to) the one induced from $\orM$. Let
$$ I(\orL, \orM) = I_+(\orL, \orM) \cup I_-(\orL, \orM).$$
Also, as in Section~\ref{sec:complete}, we set
$$ M_{\pm} = \bigcup_{i \in I_{\pm}(\orL, \orM)} L_i,$$
so that $M=M_+ \amalg M_-$.

We seek to define maps
$$ \Phi^{\orM}_{\s}: \Chain^-(\Hyper^L, \s) \to \Chain^-(\Hyper^{L-M}, \psi^{\orM}(\s)).$$

Similarly to what we did in Section~\ref{sec:statement}, these maps will be constructed as compositions 
$$ \Phi^{\orM}_{\s}= D^{\orM}_{p^{\orM}(\s)} \circ \I^{\orM}_{\s},$$
where $\I^{\orM}_{\s}$ is a projection-inclusion and $D^{\orM}_{p^{\orM}(\s)}$ is a descent map; compare \eqref{eq:PhiK} and \eqref{eq:Phi-K}.

Let us first define the map 
$$\I^{\orM}_{\s}:   \Chain^-(\Hyper^L, \s) \to  \Chain^-(\Hyper^L,\orM, p^{\orM}(\s)).$$

Here, the target $\Chain^-(\Hyper^L,\orM, p^{\orM}(\s))$ is a complex that looks like a resolution $\Chain^-$ with respect to $L-M$, like $\CFm=C\{w\}$ with respect to $M_+$, and like $\Ccint$ with respect to $M_-$. Precisely, we first consider $\Ccint(r_{M_+}(\Hyper^L))$, a complex constructed similarly to $\Ccint(\Hyper^L)$, except that for $i \in I_+(\orL, \orM)$ we do not introduce variables $V_{i,j}$ and $Y_{i,j}$; instead, for such $i$, we use  $n_{w_{i,j}}(\phi)$ as exponents for $U_{i,j}$ in the formula for the differential:
\begin{multline}
 \del \x = \sum_{\y \in \Ta \cap \Tb}  \sum_{\substack{\phi \in \pi_2(\x, \y) \\ \mu(\phi)=1} }  \# (\M(\phi)/\R) \cdot \prod_{i \not \in I_+(\orL, \orM)}  U_{i,1}^{n_{z_{i,1}}(\phi) + \ldots + n_{z_{i,p_i}}(\phi)} V_{i,2}^{n_{w_{i,2}}(\phi)} \dots  V_{i,p_i}^{n_{w_{i,p_i}}(\phi)} \\
\cdot    \prod_{i \in I_+(\orL, \orM)} U_{i,1}^{n_{w_{i,1}}(\phi)} \cdots U_{i,p_i}^{n_{w_{i,p_i}}(\phi)} 
  \prod_{j=1}^p U_{\circ,j}^{n_{w_{\circ, j}}(\phi)} \y.
 \end{multline}
The complex $\Ccint(r_{M_+}(\Hyper^L))$ has filtrations $\FF_i$ for $i \not \in I_+(\orL, \orM)$, where we use the convention that the Alexander gradings $A_i$ of $\x \in \Ta \cap \Tb$ are the ones from the diagram $\Hyper^L$ (rather than those from $r_{M_+}(\Hyper^L)$, where $\s$ becomes $\psi^{M_+}(\s)$).

We let $\Chain^-(\Hyper^L,\orM, p^{\orM}(\s))$ be the subcomplex of $\Ccint(r_{M_+}(\Hyper^L))$ in filtration degrees
$$ \FF_i \leq s_i \text{ for } L_i \not\subseteq M.$$

The map $\I^{\orM}_{\s}$ is the composition of a projection similar to $\tilde P$, but taken only with respect to the indices $i \in I_+(\orL, \orM)$, and the natural inclusion into $\Chain^-(\Hyper^L,\orM, p^{\orM}(\s))$. Specifically, we set
\begin{multline}
\label{eq:Ims}
\I^{\orM}_{\s}\Bigl(\prod_{i\in I_+(\orL, \orM)} \prod_j V_{i,j}^{n_{i,j}}  Y_{i,j}^{a_{i,j}} \x \Bigr) = \\
\begin{cases}
\prod_{i \in I_+(\orL, \orM)} U_{i,1}^{A_i(\x) - s_i - \sum_{j=2}^{p_i} n_{i,j}} U_{i,2}^{n_{i,2}} \dots U_{i,p_i}^{n_{i,p_i}} \x  & \text{if all } a_{i,j} = 0 \text{ for } i \in I_+(\orL, \orM), \\
 0 & \text{otherwise.}\end{cases}
 \end{multline}
We then extend $\I^{\orM}_{\s}$ to be equivariant with respect to the action of the variables $V_{i,j}$ and $Y_{i,j}$ for $i \not \in  I_+(\orL, \orM)$, as well as to that of all the $U$ variables.

We are left to describe the descent map 
\begin{equation}
\label{eq:descentgeneral}
D^{\orM}_{p^{\orM}(\s)}: \Chain^-(\Hyper^L,\orM, p^{\orM}(\s))\to  \Chain^-(\Hyper^{L-M}, \psi^{\orM}(\s)).
\end{equation}
This is associated to the hyperbox $\Hyper^{\orL, \orM}$ in the system $\Hyper$. Recall that the hyperbox $\Hyper^{\orL, \orM}$ consists of diagrams for $L-M$, such that the initial diagram is $r_{\orM}(\Hyper^L)$ and the final one is $\Hyper^{\orL, \orM}$, which is surface isotopic to $\Hyper^{L-M}$. The dimension of $\Hyper^{\orL, \orM}$ is the number of components of $M$, denoted $m$; we label its coordinate directions by the indices $i \in \{1, \dots, \ell\}$ such that $L_i \subseteq M$, and order them accordingly.

Let $d_i$ be the side length of the hyperbox $\Hyper^{\orL, \orM}$ in the direction $i$. To construct the descent map $D^{\orM}_{p^{\orM}(\s)}$, we will compress a hyperbox of chain complexes associated to $\Hyper^{\orL, \orM}$, of the same dimension $m$ but possibly of larger size. The models are:
\begin{itemize}
\item
the pre-composition with the $p-1$  transition maps $\Trans_j$ in Section~\ref{sec:knots2p}, to construct the descent map for a knot with $2p$ basepoints of each type;
\item the construction of the hyperbox $\Cc(\Hyper^{L, -L})$ in Section~\ref{sec:linktwo}, where we added one extra row and column at the beginning; cf. Figure~\ref{fig:Enlarged}.
\end{itemize}
In the general setting, we construct an $m$-dimensional hyperbox of chain complexes
$$C=\Cc(\Hyper^{\orL, \orM}, \psi^{\orM}(\s))$$ whose sides are labeled by $i\in I(\orL, \orM)$, and such that the side length in direction $i$ is
\begin{equation}
\label{eq:di}
 d'_i = \begin{cases}
d_i &\text{if } i\in I_+(\orL, \orM),\\
d_i + p_i - 1 &\text{if } i\in I_-(\orL, \orM)
\end{cases}.
\end{equation}
(Recall that $p_i$ denotes the number of basepoints of each type on the component $L_i$.) We let $\dd'$ be the vector with $d'_i$ as entries.

At a vertex $\eps=(\eps_i)_{i\in I(\orL, \orM)}$ in the hyperbox, we place a complex $C^{\eps}$ defined as follows. Let
$$ I(\eps) = \{ i \in I_-(\orL, \orM) \mid \eps_i < p_i - 1 \}$$
and
$$ M_{\eps}= \bigcup_{i \in I(\eps)} L_i \subseteq M_-.$$
(For example, in Section~\ref{sec:linktwo}, the part of the hyperbox $\Cc(\Hyper^{L, -L})$ corresponding to $\eps$ such that $I(\eps)=\emptyset$ is the dashed rectangle in Figure~\ref{fig:Enlarged}.)

We split the information in $\eps$ into two parts, $\eps^<$ and $\eps^>$, corresponding to $i \in I(\eps)$ and $i \not \in I(\eps)$, respectively. Specifically, $\eps^<$ is a vector whose entries are indexed by $i \in I(\eps)$, such that
$$ \eps^<_i = \eps_i.$$

The second vector, $\eps^>$, has entries indexed by $i \in I(\orL, \orM) - I(\eps)$, such that
$$ \eps^>_i = \begin{cases}\eps_i & \text{if } i \in I_+(\orL, \orM), \\
\eps_i - (p_i -1) & \text{if } i \in I_-(\orL, \orM) - I(\eps).
\end{cases}$$
In the second case, we subtract $p_i-1$ so that we are in agreement with the indexing of the entries in the hyperbox $\Hyper^{\orL, \orM}$. Indeed, note that the values of $\eps^>_i$ vary between $0$ and $d_i$.

Consider the hyperbox $\Hyper^{\orL, \orM-M_{\eps}}$, whose reduction at $M_{\eps}$ is identified with a sub-hyperbox of $\Hyper^{\orL, \orM}$; cf. Definitions~\ref{def:comp1} and \ref{def:precomplete}. The complex $C^{\eps}$ will be associated to the diagram at position $\eps^>$ in the hyperbox $ \Hyper^{\orL, \orM-M_{\eps}}.$
We write
\begin{equation}
\label{eq:CCmix}
 C^{\eps} = \Cc_{\eps^<}(\Hyper^{\orL, \orM-M_{\eps}}_{\eps^>}, M_{\eps}, \psi^{\orM}(\s)).
 \end{equation}
This is a Floer complex whose construction is a mix of those in Sections~\ref{sec:knots2p} and \ref{sec:linktwo}. Roughly speaking, $C^{\eps}$ looks like
\begin{itemize}
\item the resolution $\Chain^-(\cdot, \psi^{\orM}(\s))$ in the directions $i$ such that $L_i \subseteq L-M$, i.e., $i \not \in I(\orL, \orM)$,
\item the complexes $\Cc_j$ from \eqref{eq:transitions}, with $j=\eps_i$, in the directions $i \in I(\eps)$, 
\item the complex $\CFm=C\{z\}$ in the directions $i \in I_-(\orL, \orM) - I(\eps)$,
\item the complex $\CFm=C\{w\}$ in the directions $i \in I_+(\orL, \orM)$. 
\end{itemize}
We will construct $C^{\eps}$ as a subcomplex of a complex denoted
$$ \Cc_{\eps^<}(\Hyper^{\orL, \orM-M_{\eps}}_{\eps^>}, M_{\eps}),$$
which is similar, except it looks like $\Ccint$ in the directions $i\not \in I(\orL, \orM)$.

Let us make these definitions fully rigorous. First, recall that we have a base ring
$$ \Ring = \Field[[(U_{i,j})_{\substack{1\leq i \leq \ell, \\ 1 \leq j \leq p_i}} (U_{\circ,j})_{1\leq j \leq p}]].$$
We introduce the dga
$$ \Ringbig^Y_{\eps^<}(\orL, \orM, M_{\eps})$$
generated over $\Ring$ (in the sense of infinite sums with boundedness conditions) by new variables $V_{i,j}, Y_{i,j}$, for indices $i, j$ such that
\begin{equation}
\label{eq:condo}
\bigl(i \not\in I(\orL, \orM), \ 2 \leq j \leq p_i \bigr ) \ \text{ or } \ \bigl(i \in I(\eps), \eps_i +2 \leq j \leq p_i \bigr)
\end{equation}
with relations
$$ Y_{i,j}^2=0, \ \ \del Y_{i,j} = U_{i,j} + U_{i,1}V_{i,j}.$$
We let $ \Cc_{\eps^<}(\Hyper^{\orL, \orM-M_{\eps}}_{\eps^>}, M_{\eps})$ be the complex freely generated over $\Ringbig^Y_{\eps^<}(\orL, \orM, M_{\eps})$ by $\x \in \Ta \cap \Tb$, with differential
\begin{multline}
\label{eq:multidel}
 \del \x = \sum_{\y \in \Ta \cap \Tb}  \sum_{\substack{\phi \in \pi_2(\x, \y) \\ \mu(\phi)=1} }  \# (\M(\phi)/\R) \cdot \prod_{i \not \in I(\orL, \orM)}  U_{i,1}^{n_{z_{i,1}}(\phi) + \ldots + n_{z_{i,p_i}}(\phi)} V_{i,2}^{n_{w_{i,2}}(\phi)} \dots  V_{i,p_i}^{n_{w_{i,p_i}}(\phi)} \\
\cdot \prod_{i \in I(\eps)}  U_{i,1}^{n_{z_{i,1}}(\phi)} \cdots U_{i,\eps_i}^{n_{z_{i,\eps_i}}(\phi)}
U_{i,\eps_i+1}^{n_{z_{i,\eps_i+1}}(\phi) + \ldots + n_{z_{i,p_i}}(\phi)} V_{i,\eps_i+2}^{n_{w_{i,\eps_i+2}}(\phi)} \dots  V_{i,p_i}^{n_{w_{i,p_i}}(\phi)}\\
\cdot  \prod_{i \in I_-(\orL, \orM) - I(\eps)} U_{i,1}^{n_{z_{i,1}}(\phi)} \cdots U_{i,p_i}^{n_{z_{i,p_i}}(\phi)} 
   \prod_{i \in I_+(\orL, \orM)} U_{i,1}^{n_{w_{i,1}}(\phi)} \cdots U_{i,p_i}^{n_{w_{i,p_i}}(\phi)} 
  \prod_{j=1}^p U_{\circ,j}^{n_{w_{\circ, j}}(\phi)} \y.
 \end{multline}

This complex admits filtrations $\FF_i$, for $i \not\in I(\orL, \orM)$, such that
$$ \FF_i(\x) = -A_i(\x), \ \FF_i(U_{i,1})=-1, \ \FF_i(V_{i,j})=1, \ j \geq 2,$$
and $\FF_i$ takes all the other $U$, $V$ and $Y$ variables to zero. (Once again, the Alexander gradings $A_i$ are taken to be the ones induced from the diagram $\Hyper^L$.)

We then define $C^{\eps} = \Cc_{\eps^<}(\Hyper^{\orL, \orM-M_{\eps}}_{\eps^>}, M_{\eps}, \psi^{\orM}(\s))$ to be the subcomplex of $\Cc_{\eps^<}(\Hyper^{\orL, \orM-M_{\eps}}_{\eps^>}, M_{\eps})$ given by
$$ \FF_i \leq s_i, \ i \not\in I(\orL, \orM).$$

In particular, when $\eps = \zero$, observe that $M_{\zero} = M_-$, $I(\zero)=I_-(\orL, \orM)$ and (by the compatibility conditions between hyperboxes) $\Hyper^{\orL, \orM-M_{\zero}}_{\zero}=r_{M+}(\Hyper^L).$ Therefore, we have identifications
$$ \Cc_{\zero}(\Hyper^{\orL, \orM-M_{\zero}}_{\zero}, M_{\zero}) \cong \Ccint(r_{M_+}(\Hyper^L))$$
and
\begin{equation}
\label{eq:Cinitial}
C^{\zero}=\Cc_{\zero}(\Hyper^{\orL, \orM-M_{\zero}}_{\zero}, M_{\zero}, \psi^{\orM}(\s)) \cong \Chain^-(\Hyper^L,\orM, p^{\orM}(\s)).
\end{equation}

At the other extreme, for $\eps$ such that $I(\eps)=\emptyset$, we have $M_{\eps} = \emptyset$ and
$$
 C^\eps \cong \Chain^-(\Hyper^{\orL, \orM}_{\eps^>}, \psi^{\orM}(\s)).
 $$
Thus, the final complex in the hyperbox $\Cc(\Hyper^{\orL, \orM}, \psi^{\orM}(\s))$, with $\eps= \dd'$, is
\begin{equation}
\label{eq:Cfinal}
C^{\dd'} \cong\Chain^-(\Hyper^{\orL, \orM}(M), \psi^{\orM}(\s)).
 \end{equation}

For general $\eps, \eps'$ such that $\eps \leq \eps'$ and $\eps, \eps'$ are neighbors, we now proceed to define the maps
$$ D^{\eps' - \eps}_{\eps}: C^{\eps} \to C^{\eps'}$$
that are part of the hyperbox $\Cc(\Hyper^{\orL, \orM}, \psi^{\orM}(\s))$. 

Since $\eps$ and $\eps'$ are neighbors, for any $i$ we must have $\eps'_i = \eps_i$ or $\eps'_i=\eps_i+1$. Observe that $I(\eps') \subseteq I(\eps)$. When an index $i$ belongs to $I(\eps) - I(\eps')$, it means that $\eps_i = p_i-2$ and $\eps'_i = p_i -1$, so $(\eps')^>$ has $0$ in position $i$. Thus, while $(\eps')^>$ may have more entries that $\eps^>$, all these additional entries are zero; let us denote by $(\eps^>)'$ the vector obtained from $(\eps')^>$ by deleting these entries. In the same manner, $(\eps')^<$ may have fewer entries than $\eps^<$, in which case the additional entries of $\eps^<$ are $p_i-2$; we let $(\eps^<)'$ be the vector obtained from $(\eps^<)'$ by introducing new entries of $p_i-1$ in those positions. Thus, $(\eps^>)'$ has the same length as $\eps^>$, and $(\eps^<)'$ has the same length as $\eps^<$.

Note that (by the compatibility between hyperboxes) $\Hyper^{\orL, \orM-M_{\eps'}}_{(\eps')^>}$, the diagram used to define $C^{\eps'}$, is a reduction of a diagram in the hyperbox $\Hyper^{\orL, \orM-M_{\eps}}$:
$$ \Hyper^{\orL, \orM-M_{\eps'}}_{(\eps')^>} = r_{M_{\eps} - M_{\eps'}} (\Hyper^{\orL, \orM-M_{\eps}}_{(\eps^>)'}).$$

Thus, the complex $C^{\eps'}$ can be interpreted as being associated to the diagram $\Hyper^{\orL, \orM-M_{\eps}}_{(\eps^>)'}.$ Indeed, from the definitions we see that
\begin{align*}
 C^{\eps'} &=  \Cc_{(\eps')^<}(\Hyper^{\orL, \orM-M_{\eps'}}_{(\eps')^>}, M_{\eps'}, \psi^{\orM}(\s)) \\
 &= \Cc_{(\eps')^<}(r_{M_{\eps} - M_{\eps'}} (\Hyper^{\orL, \orM-M_{\eps}}_{(\eps^>)'}), M_{\eps'}, \psi^{\orM}(\s)) \\
 &=  \Cc_{(\eps^<)'}(\Hyper^{\orL, \orM-M_{\eps}}_{(\eps^>)'}, M_{\eps}, \psi^{\orM}(\s)).
 \end{align*}

Similar remarks apply to the intermediate indices between $\eps$ and $\eps'$. Specifically, we can write the corresponding complexes as
$$ \Cc_{\nu}(\Hyper^{\orL, \orM-M_{\eps}}_{\eta}, M_{\eps}, \psi^{\orM}(\s))$$
for all $\nu, \eta$ with
$$ \eps^< \leq \nu \leq (\eps^<)', \ \ \ \eps^> \leq \eta \leq (\eps^>)'.$$
These complexes have almost the same definition as the one for index $\eps$, except that changing $\eps^>$ into $\eta$ results in a change in the underlying Heegaard diagram (a different vertex in the hyperbox $\Hyper^{\orL, \orM-M_{\eps}}$), whereas changing $\eps^<$ into $\nu$ results in a change in the exponents of $U_{i,j}$ for some $i\in I(\eps)$ in the formula \eqref{eq:multidel}; namely, $\eps_i$ gets replaced by $\nu_i$, which may equal $\eps_i + 1$.

In light of this discussion, to define $ D^{\eps' - \eps}_{\eps}$ we will focus on the hyperbox
$$\Hyper^{\orL, \orM-M_{\eps}},$$
and specifically on the sub-hypercube in that hyperbox corresponding to indices $\eta$ such that
$$ \eps^> \leq \eta \leq (\eps^>)'.$$
We denote this hypercube by
$$ \Hyper^{\orL, \orM-M_{\eps}}(\eps^>, (\eps^>)').$$
This is a $v$-dimensional hypercube of strongly equivalent Heegaard diagrams, with 
$$v=\| (\eps^>)' - \eps^{>} \|.$$ 

Recall that in the link-minimal case, the maps $ D^{\eps' - \eps}_{\eps}$ were defined in Section~\ref{sec:hyperfloer} as counts of holomorphic $k$-gons, where $k \leq (\| \eps' - \eps \|+2)$. In our setting, they will be given by counts of holomorphic $k$-gons for $k \leq v+2$, coming from the above hypercube.

The formula in the link-minimal case was:
$$ \De^{\eps'-\eps}_{\eps}(\x) =\sum_{l, q} \sum_{\{\eps^{'\alpha} = \gamma^0 > \dots > \gamma^l = \eps^\alpha\}}  \sum_{\{\eps^\beta = \zeta^0 < \dots < \zeta^q = \eps^{'\beta}\}}  f(\Theta^\alpha_{\gamma^0, \gamma^1} \otimes \dots \otimes \Theta^\alpha_{\gamma^{l-1}, \gamma^l} \otimes \x \otimes \Theta^\beta_{\zeta^0, \zeta^1} \otimes \dots \otimes \Theta^\beta_{\zeta^{q-1}, \zeta^q} ).
$$

In general, our hypercube still has bipartition maps (cf. Section~\ref{sec:hse}), so we can talk about indices
$$ (\eps^>)^{\alpha} \leq (\eps^>)^{'\alpha}, \ \ \  (\eps^>)^{\beta} \leq (\eps^>)^{'\beta}.$$
For any intermediate $\gamma$ between (and including) $(\eps^>)^{ \alpha}$ and $(\eps^>)^{'\alpha}$,  as part of the hyperbox we have a collection of curves $\alphas^\gamma$, such that all these collections are strongly equivalent. In particular, we write $\alphas =\alphas^{(\eps^>)^{\alpha}}$ and $\alphas' =\alphas^{(\eps^>)^{'\alpha}}.$  Furthermore, for any two intermediate indices $\gamma, \gamma'$ with
$$  (\eps^{>})^{\alpha} \leq \gamma \leq \gamma' \leq (\eps^>)^{'\alpha}$$
we have a chain element
$$ \Theta^{\alpha}_{\gamma', \gamma} \in \Am(\T_{\alpha^{\gamma'}}, \T_{\alpha^\gamma}, \zero).$$
We can write this element explicitly in terms of intersection points:
$$ \Theta^{\alpha}_{\gamma', \gamma} = \sum_{\a^{\alpha}_{\gamma' \gamma} \in  \T_{\alpha^{\gamma'}} \cap \T_{\alpha^\gamma}} c(\a^{\alpha}_{\gamma' \gamma}) \prod_{i,j} U_{i,j}^{u_{i,j}(\a^{\alpha}_{\gamma' \gamma})} \prod_j U_{\circ, j}^{u_{\circ,j}(\a^{\alpha}_{\gamma' \gamma})} \a^{\alpha}_{\gamma' \gamma},$$
with coefficients $c(\a^{\alpha}_{\gamma' \gamma}) \in \{0,1\}$ and exponents $u_{i,j}(\a^{\alpha}_{\gamma' \gamma}), u_{\circ,j}(\a^{\alpha}_{\gamma' \gamma}) \geq 0.$ 

Note that $ \Theta^{\alpha}_{\gamma', \gamma}$ are part of the data in the hyperbox  $\Hyper^{\orL, \orM - M_{\eps}}$, where the $U_{i,j}$ basepoints keep track of the basepoints $z_{i,j}$ for $i \in I_-(\orL, \orM) - I(\eps)$, and of the basepoints $w_{i,j}$ for other $i$. 

Further, the $\alpha^{\gamma'}$-$\alpha^{\gamma}$ diagram is link-minimal, and therefore the intersection points $\a^{\alpha}_{\gamma' \gamma}$ admit Alexander gradings $A_{i,j}$ for $L_i \not \subseteq M-M_{\eps}$, that is, for $i \not \in I(\orL, \orM)$ as well as for $i \in I(\eps)$, and for all corresponding $j$. We will need to pull back the $\Theta$ elements as in Section~\ref{sec:thetas}, so that the corresponding exponents are in terms of the $z$ basepoints rather than the $w$ basepoints. It helps to introduce the notation
$$ \bu_{i,j}(\a^{\alpha}_{\gamma' \gamma}) := u_{i,j}(\a^{\alpha}_{\gamma' \gamma}) - A_{i,j}(\a^{\alpha}_{\gamma' \gamma}).$$

The hypercube $ \Hyper^{\orL, \orM-M_{\eps}}(\eps^>, (\eps^>)')$ also contains collections of beta curves and corresponding $\Theta$ elements. For $\zeta$ between $ (\eps^>)^{\beta}$ and $(\eps^>)^{'\beta}$ we have a collection of curves $\betas^{\zeta}$, and for 
$$  (\eps^{>})^{\beta} \leq \zeta \leq \zeta' \leq (\eps^>)^{'\beta}$$
we have a chain element
$$ \Theta^{\beta}_{\zeta, \zeta'} \in \Am(\T_{\beta^{\zeta}}, \T_{\beta^{\zeta'}}, \zero)$$
of the form
$$ \Theta^{\beta}_{\zeta, \zeta'} = \sum_{\a^{\beta}_{\zeta \zeta'} \in  \T_{\beta^{\zeta}} \cap \T_{\beta^{\zeta'}}} c(\a^{\beta}_{\zeta \zeta'}) \prod_{i,j} U_{i,j}^{u_{i,j}(\a^{\beta}_{\zeta \zeta'})} \prod_j U_{\circ, j}^{u_{\circ,j}(\a^{\beta}_{\zeta \zeta'})} \a^{\beta}_{\zeta \zeta'}.$$
Similar remarks to those for the alpha curves apply here. In particular, we set
$$ \bu_{i,j}(\a^{\beta}_{\gamma' \gamma}) := u_{i,j}(\a^{\beta}_{\zeta \zeta'}) - A_{i,j}(\a^{\beta}_{\zeta \zeta'}).$$

The map $D_{\eps}^{\eps'-\eps}: C^{\eps} \to C^{\eps'}$ are defined as follows. For an intersection point $ \x \in \Ta \cap  \Tb$, we set
\begin{equation}
\label{eq:bigformula}
D_{\eps}^{\eps'-\eps}(\x) = \sum_{E} \#\M(\phi) \cdot \c(E) \U^{u(E)} \y,
\end{equation}
where the summation is over data $E$ consisting of:
\begin{itemize}
\item$l, q \geq 0, $\\
\item $ (\eps^>)^{'\alpha} = \gamma^0 > \dots > \gamma^l = (\eps^>)^{\alpha}$,\\
\item $ (\eps^>)^{\beta} = \zeta^0< \dots < \zeta^q = (\eps^>)^{'\beta},$\\
\item $ \a^{\alpha}_{\gamma^k \gamma^{k+1}} \in \T_{\alpha_{\gamma^k}} \cap \T_{\alpha_{\gamma^{k+1}}} \text{ for } k=0, \dots, l-1, $\\
\item $ \a^{\beta}_{\zeta^k \zeta^{k+1}} \in \T_{\beta_{\zeta^k}} \cap \T_{\beta_{\zeta^{k+1}}} \text{ for } k=0, \dots, q-1, $\\
\item $ \y \in  \Tap \cap  \Tbp$,\\
\item $ \phi \in \pi_2(\a^{\alpha}_{\gamma^0\gamma^1}, \dots, \a^{\alpha}_{\gamma^{l-1}\gamma^l} , \x, \a^{\beta}_{\zeta^0\zeta^1}, \dots, \a^{\beta}_{\zeta^{q-1}\zeta^q}, \y) \text{ with } \mu(\phi) = 1-l-q.$
\end{itemize}
In the formula \eqref{eq:bigformula}, by $\# \M(\phi)$ we mean the count of holomorphic polygons in the class $\phi$ (where we divide by the action of $\rr$ if we count bigons). Moreover, the coefficient $\c(E) \in \{0,1\}$ is 
$$ \c(E) := \prod_{i,j} c(\a^{\alpha}_{\gamma^0 \gamma^1}) \cdots c(\a^{\alpha}_{\gamma^{l-1}\gamma^l}) \cdot c(\a^{\beta}_{\zeta^0\zeta^1}) \cdots c(\a^{\beta}_{\zeta^{q-1}\zeta^q})$$
whereas the factor $\U^{u(E)}$ is a product of powers of $U_{i,j}$ that will be described soon. Let us first introduce the shorthand notation
$$
 n^{\a}_{w_{i,j}}(\phi) := n_{w_{i,j}}(\phi) + \sum_{k=0}^{l-1} u_{i,j}(\a^{\alpha}_{\gamma^k \gamma^{k+1}}) + \sum_{k=0}^{q-1} u_{i,j} (\a^{\beta}_{\zeta^k \zeta^{k+1}}),$$
$$ n^{\a}_{z_{i,j}}(\phi) := n_{z_{i,j}}(\phi) + \sum_{k=0}^{l-1} \bu_{i,j}(\a^{\alpha}_{\gamma^k \gamma^{k+1}}) + \sum_{k=0}^{q-1} \bu_{i,j} (\a^{\beta}_{\zeta^k \zeta^{k+1}}),$$
$$  n^{\a}_{w_{\circ,j}}(\phi) := n_{w_{\circ,j}}(\phi) + \sum_{k=0}^{l-1} u_{\circ,j}(\a^{\alpha}_{\gamma^k \gamma^{k+1}}) + \sum_{k=0}^{q-1} u_{\circ,j} (\a^{\beta}_{\zeta^k \zeta^{k+1}}).$$
Then, by analogy with the $U$-powers in the differentials \eqref{eq:Transjsame}, \eqref{eq:diagonalnoY} and \eqref{eq:multidel}, we set
\begin{multline}
\label{eq:multiU}
\U^{u(E)} :=  \prod_{i \not \in I(\orL, \orM)}  U_{i,1}^{n^{\a}_{z_{i,1}}(\phi) + \ldots + n^{\a}_{z_{i,p_i}}(\phi) } V_{i,2}^{n^{\a}_{w_{i,2}}(\phi)} \dots  V_{i,p_i}^{n_{w_{i,p_i}}(\phi)} \\
\cdot \prod_{\{i \in I(\eps)\mid \nu_i = \eps_i \}}  U_{i,1}^{n^{\a}_{z_{i,1}}(\phi)} \cdots U_{i,\eps_i}^{n^{\a}_{z_{i,\eps_i}}(\phi)}
U_{i,\eps_i+1}^{n^{\a}_{z_{i,\eps_i+1}}(\phi) + \ldots + n^{\a}_{z_{i,p_i}}(\phi)} V_{i,\eps_i+2}^{n^{\a}_{w_{i,\eps_i+2}}(\phi)} \dots  V_{i,p_i}^{n^{\a}_{w_{i,p_i}}(\phi)}\\
\cdot \prod_{\{i \in I(\eps)\mid \nu_i = \eps_i+1 \}}  U_{i,1}^{n^{\a}_{z_{i,1}}(\phi)} \cdots U_{i,\eps_i+1}^{n^{\a}_{z_{i,\eps_i+1}}(\phi)}
\frac{U_{i,\eps_i+1}^{n^{\a}_{z_{i,\eps_i+2}}(\phi) + \ldots + n^{\a}_{z_{i,p_i}}(\phi)} - U_{i,\eps_i+2}^{n^{\a}_{z_{i,\eps_i+2}}(\phi) + \ldots + n^{\a}_{z_{i,p_i}}(\phi)} }{U_{i, \eps_i + 1} - U_{i, \eps_i + 2}} V_{i,\eps_i+3}^{n^{\a}_{w_{i,\eps_i+3}}(\phi)} \dots  V_{i,p_i}^{n^{\a}_{w_{i,p_i}}(\phi)}
\\
\cdot  \prod_{i \in I_-(\orL, \orM) - I(\eps)} U_{i,1}^{n^{\a}_{z_{i,1}}(\phi)} \cdots U_{i,p_i}^{n^{\a}_{z_{i,p_i}}(\phi)} 
\cdot   \prod_{i \in I_+(\orL, \orM)} U_{i,1}^{n^{\a}_{w_{i,1}}(\phi)} \cdots U_{i,p_i}^{n^{\a}_{w_{i,p_i}}(\phi)} 
 \cdot  \prod_{j=1}^p U_{\circ,j}^{n^{\a}_{w_{\circ, j}}(\phi)}.
\end{multline}
We have described the values of $D^{\eps'-\eps}_{\eps}$ at intersection points $\x$. To define the map in full, we extend it to be equivariant with respect to the action of all the $U_{i,j}, V_{i,j}, Y_{i,j}$ and $U_{\circ, j}$ variables that appear in the construction of the target $C^{\eps'}$. In view of the conditions \eqref{eq:condo}, note that the domain $C^{\eps}$ may have a few more variables, namely 
$$V_{i,\eps_i+2}, Y_{i, \eps_i + 2} \text{ for } i \in I(\eps) \text{ such that } \eps'_i = \eps_i + 1.$$

We let $D^{\eps'-\eps}_{\eps}$ be equivariant with respect to $V_{i, \eps_i+2}$, where we let this variable act on the target by $1$. With regard to $Y_{i, \eps_i + 2}$, this behaves somewhat like in an identity hyperbox, with identity maps on the edges and zeros on the diagonals. Precisely, we set 
$$D^{\eps'-\eps}_{\eps}(Y_{i, \eps_i + 2} \x) = 0\  \text{ when } \|\eps' - \eps\| > 1,$$
whereas when $\|\eps' - \eps\| =1$ there is a unique $i$ with $\eps'_i = \eps_i + 1$, and (for the extra variable $Y$ to exist) we must have $i \in I(\eps)$; that is,
$$ (\eps^<)'_i = \eps^<_i + 1, \ \  (\eps^>)'_i = \eps^>_i.$$

In this particular case and for that particular $i$, by analogy with \eqref{eq:Transjdiff} and \eqref{eq:Transjsame}, the formula for $D^{\eps'-\eps}_{\eps}(Y_{i, \eps_i + 2} \x)$ differs from the one for $D^{\eps'-\eps}_{\eps}( \x)$ by replacing 
$$U_{i,1}^{n^{\a}_{z_{i,1}}(\phi)} \cdots U_{i,\eps_i+1}^{n^{\a}_{z_{i,\eps_i+1}}(\phi)}
\frac{U_{i,\eps_i+1}^{n^{\a}_{z_{i,\eps_i+2}}(\phi) + \ldots + n^{\a}_{z_{i,p_i}}(\phi)} - U_{i,\eps_i+2}^{n^{\a}_{z_{i,\eps_i+2}}(\phi) + \ldots + n^{\a}_{z_{i,p_i}}(\phi)} }{U_{i, \eps_i + 1} - U_{i, \eps_i + 2}} V_{i,\eps_i+3}^{n^{\a}_{w_{i,\eps_i+3}}(\phi)} \dots  V_{i,p_i}^{n^{\a}_{w_{i,p_i}}(\phi)}$$
on the third line of the formula \eqref{eq:multiU}, with $1$.

We have now completed the description of the hyperbox $\Cc(\Hyper^{\orL, \orM}, \psi^{\orM}(\s))$. This plays the role of the hyperbox $\Chain^-(\Hyper^{\orL, \orM}, \psi^{\orM}(\s))$ from the link-minimal case presented in Section~\ref{subsec:desublink}. Just as there, we let
$$  \hat\De^{\orM}_{\s} : \Chain^-(\he^L, \orM, \s) \to \Chain^-(\he^{\orL, \orM}(M), \psi^{\orM}(\s)),$$
be the longest diagonal map in the hypercube obtained from $\Cc(\Hyper^{\orL, \orM}, \psi^{\orM}(\s))$ by compression; compare \eqref{eq:Cinitial} and \eqref{eq:Cfinal}. After identifying $\he^{\orL, \orM}(M)$ with $\he^{L-M}$, we get the desired descent map $D^{\orM}_{p^{\orM}(\s)}$ as in \eqref{eq:descentgeneral}. By pre-composing it with the corresponding projection-inclusion map, we get 
$$ \Phi^{\orM}_{\s}= D^{\orM}_{p^{\orM}(\s)} \circ \I^{\orM}_{\s}.$$
These are the maps used to define the surgery complex $\C^-(\Hyper, \Lambda)$. The fact that the differential on this complex squares to zero is an extension of Proposition~\ref{prop:phiphi}, with a similar proof.

\subsection{Proof of the theorem}
\label{sec:proofgen}
In this section we prove Theorem~\ref{thm:FirstSurgery} for arbitrary complete systems $\Hyper$. In Section~\ref{sec:proof} we have already given the proof in the case where $\Hyper$ is link-minimal. Recall that any complete system $\Hyper$ can be obtained from a basic one by a sequence of the system moves listed in Section~\ref{sec:moves}; cf. Proposition~\ref{prop:moves}. Thus, to prove Theorem~\ref{thm:FirstSurgery} it suffices to check that the quasi-isomorphism type of the surgery complex $\C^-(\Hyper, \Lambda)$ is unchanged by these moves. Invariance under most of these moves (3-manifold isotopies, index one/two stabilizations and destabilizations, free index zero/three stabilizations, global shifts, and elementary enlargements and contractions) follows by the same arguments as in the link-minimal case; see Section~\ref{sec:invariance}.

It remains to establish the following. 

\begin{proposition}
Let $\bHyper$ be a complete system of hyperboxes for a link $\orL \subset Y$, and let $\Hyper$ be the system obtained from $\bHyper$ by an index zero/three link stabilization. Then, for any framing $\Lambda$, the surgery complexes $\C^-(\bHyper, \Lambda)$ and $\C^-(\Hyper, \Lambda)$ are quasi-isomorphic over $\ff[[U]].$
\end{proposition}

\begin{proof}
We keep the notation from Section~\ref{sec:gens} for everything in the construction of $\C^-(\Hyper, \Lambda)$. Furthermore, we assume that the stabilization from $\bHyper$ to $\Hyper$ is done on the link component $L_1$, and introduces the points labeled $w_{1,2}$ and $z_{1,2}$ in the neighborhood of an existing $z_{1,1}$; cf. Figure~\ref{fig:stabilize2}. \footnote{Strictly speaking, we lose some generality by choosing the stabilization on $L_1$, because the compression procedure from Section~\ref{sec:hyperco} involves an ordering of the components. However, this issue will not be relevant here, and the arguments below apply equally well to all other components, with only cosmetic modifications.} Thus, the basepoints on $L_1$ in $\bHyper$ are $w_{1,1}, z_{1,1}, w_{1,3}, z_{1,3}, \dots, w_{1,p_1}, z_{1, p_1}.$

\begin{figure}
\begin{center}
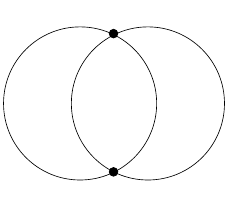
\end{center}
\caption {An index zero/three link stabilization.}
\label{fig:stabilize2}
\end{figure}

Note that $\C^-(\Hyper, \Lambda)$ is a complex over the ring
$$ \Ring = \ff[[(U_{i,j})_{1\leq i \leq \ell, \  1 \leq j \leq p_i}, (U_{\circ,j})_{1\leq j \leq p}]]$$
and $\C^-(\bHyper, \Lambda)$ is a complex over the smaller ring $\bRing$, in which we do not have the variable $U_{1, 2}$. We will in fact prove that the complexes $\C^-(\bHyper, \Lambda)$ and $\C^-(\Hyper, \Lambda)$ are quasi-isomorphic over $\bRing$. This will imply that they are quasi-isomorphic over $\ff[[U]]$, where $U$ is any of the variables in $\bRing$. 

The proof will be very similar to the one in Section~\ref{sec:toylink}. 

Let us start by comparing the complexes $\Chain^-(\Hyper^L, \s)$ and $\Chain^-(\bHyper^L, \s)$. 
These are subcomplexes of $\Ccint(\Hyper^L)$ and $\Ccint(\bHyper^L)$, respectively. 
Recall that $\Ccint(\Hyper^L)$ was described in Section~\ref{sec:gens} as a dg module over the dga 
$$ \Ringbig^{Y} = \Field[[(U_{i,j})_{\substack{1\leq i \leq \ell, \\ 1 \leq j \leq p_i}} (U_{\circ,j})_{1\leq j \leq p}]][(V_{i,j}, Y_{i,j})_{\substack{1\leq i \leq \ell, \\ 2 \leq j \leq p_i}}]'/(Y_{i,j}^2=0, \del Y_{i,j} = U_{i,j} + U_{i,1}V_{i,j}).$$
Similarly, $\Ccint(\bHyper^L)$ is a dg module over a dga $\bRingbig^Y$; compared to $\Ringbig^Y$, in $\bRingbig^Y$ we do not have the variables $U_{1,2}, V_{1,2}$ and $Y_{1,2}$. Further, by looking at the formula \eqref{eq:Ccintgen} for the differential on $\Ccint(\Hyper^L)$, and its analogue for $\Ccint(\bHyper^L)$, we deduce that $\Ccint(\Hyper^L)$ is obtained from $\Ccint(\bHyper^L)$ by first considering the cone
\begin{equation}
\label{eq:coneCint}
\Ccint(\bHyper^L)^{U_{1,1} \to U_{1,2}}_-[[U_{1,1}]][V_{1,2}]' \xrightarrow{V_{1,2} -1}  \Ccint(\bHyper^L)^{U_{1,1} \to U_{1,2}}_+[[U_{1,11}]][V_{1,2}]'
\end{equation}
and then taking the cone of $U_{1,2} + U_{1,1} V_{1,2}$, i.e., introducing the variable $Y_{1,2}$ with $\del Y_{1,2}=U_{1,2} + U_{1,1} V_{1,2}$. We are using here the description of the holomorphic disks in the stabilized diagram; cf. Proposition~\ref{prop:ZemkeDisks}, and the subscripts $-$ and $+$ refer to containing the points $x_-$ resp. $x_+$ from Figure~\ref{fig:stabilize2}.  The situation is entirely similar to the one in Section~\ref{sec:toylink}, where $\Ccint$ was obtained from $\bar C\{z ; n \geq 0\}$ in this way. As explained there, this implies that the complexes are homotopy equivalent over $\bRing$. In our case, the equivalence is given by the projection
$$ \rho: \Ccint(\Hyper^L) \to \Ccint(\bHyper^L)$$
 \begin{equation}
\label{eq:rhomidformula2}
 \rho(U_{1,2}^{m} V_{1,2}^n Y_{1,2}^a (\x \times r)) = \begin{cases}
U_{1,1}^m \x & \text{if } a=0, r=x_+, \\
0 & \text{otherwise.}
\end{cases}
\end{equation}
Compare Equation~\eqref{eq:rhomidformula}. The map $\rho$ is filtered with respect to the filtrations $\FF_i$, and as such it induces an equivalence between the subcomplexes in filtration levels $\FF_i \leq s_i$: 
$$\rho: \Chain^-(\Hyper^L, \s) \to \Chain^-(\bHyper^L, \s).$$

Similar remarks apply when we compare the complexes $\Chain^-(\Hyper^{L-M}, \s)$ and $\Chain^-(\bHyper^{L-M}, \s)$ for sublinks $M$ such that $L_1 \not \subseteq M.$ We get that these complexes are equivalent via a projection $\rho$, defined by the same formula \eqref{eq:rhomidformula2}.

When $L_1 \subseteq M$, the situation is simpler, as the basepoints $z_{1,1}$ and $z_{1,2}$ disappear. The diagram $\Hyper^{L-M}$ is obtained from $\bHyper^{L-M}$ by a free index zero/three stabilization. Hence, we can identify $\Chain^-(\Hyper^{L-M}, \s)$ with the cone
$$\Chain^-(\bHyper^{L-M}, \s)_+[[U_{1,2}]] \xrightarrow{U_{1,2} -U_{1,1}}  \Chain^-(\bHyper^{L-M}, \s)_-[[U_{1,2}]].$$
Thus, we can still construct an equivalence
$$ \rho: \Chain^-(\Hyper^{L-M}, \s) \to \Chain^-(\bHyper^{L-M}, \s),$$
this time given by
\begin{equation}
\label{eq:rhomidformula3}
 \rho (U_{1,2}^{m}  (\x \times r)) = \begin{cases}
U_{1,1}^m \x & \text{if }  r=x_+, \\
0 & \text{otherwise.}
\end{cases}
\end{equation}

We can view the  surgery complexes $\C^-(\Hyper, \Lambda)$ and $\C^-(\bHyper, \Lambda)$ as $\ell$-dimensional hypercubes (as in Section~\ref{sec:prelim}). We aim to construct a chain map between these hypercubes, in the sense of Definition~\ref{def:chmap}, by starting from the equivalences
$$\rho: \Chain^-(\Hyper^{L-M}, \s) \to \Chain^-(\bHyper^{L-M}, \s),$$
and then adding diagonal maps as in the diagram
\begin{equation}
\label{eq:Phirho}
\xymatrix{
\Chain^-(\Hyper^L, \s) \ar[r]^-{\Phi^{\orM}_{\s}} \ar[d]_{\rho} \ar[dr] & \Chain^-(\Hyper^{L-M}, \psi^{\orM}(\s)) \ar[d]^{\rho} \\
\Chain^-(\bHyper^L, \s) \ar[r]^-{\Phi^{\orM}_{\s}} &\Chain^-(\bHyper^{L-M}, \psi^{\orM}(\s)).
}
\end{equation}

For these to form a chain map between hypercubes, we need for example that when $M$ is a knot (so that $\Phi^{\orM}_{\s}$ are chain maps), the projections $\rho$ commute with $\Phi^{\orM}_{\s}$ up to the chain homotopy given by the diagonal. In general, we need the diagonal maps to be higher chain homotopies, i.e., to satisfy the relations \eqref{eq:DF}.

Observe that, if we can construct such a chain map, then it would be filtered with respect to the filtration given by the value of $-\|\eps\|$ in the hypercube. The maps on the associated graded would be just the projections $\rho$, which are quasi-isomorphisms over $\bRing$. This would imply that the whole map from $\C^-(\Hyper, \Lambda)$ and $\C^-(\bHyper, \Lambda)$ is a quasi-isomorphism, and complete the proof of the proposition.

We are left to find diagonal maps as in \eqref{eq:Phirho}. Recall that each $\Phi^{\orM}_{\s}$ is the composition of a projection-inclusion $\I^{\orM}_{\s}$ and a descent map $D^{\orM}_{p^{\orM}(\s)}$. We will first show that the projections $\rho$ commute with the projection-inclusions (on the nose), and then that they commute with the descent maps up to chain homotopies that fit into suitable hypercubes.

Thus, with regard to the projection-inclusions, we seek commutative diagrams of the form
\begin{equation}
\label{eq:commuterhoint}
\xymatrix{
\Chain^-(\Hyper^L, \s) \ar[r]^-{\I^{\orM}_{\s}} \ar[d]_{\rho} & \Chain^-(\Hyper^L, \orM, p^{\orM}(\s) )\ar[d]^{\rho} \\
\Chain^-(\bHyper^L, \s) \ar[r]^-{\I^{\orM}_{\s}} & \Chain^-(\bHyper^L, \orM, p^{\orM}(\s)).
}
\end{equation}
Here, the map $\rho$ on the right hand side is defined just as before, through formulas of the type \eqref{eq:rhomidformula3} if $L_1 \subseteq M_+$, or \eqref{eq:rhomidformula2} otherwise. 

To check that \eqref{eq:commuterhoint} commutes, we use the formula \eqref{eq:Ims} for the projection-inclusion maps. If $L_1 \not \subseteq M_+$, the verification is straightforward, since both vertical maps are given by the same formula \eqref{eq:rhomidformula2}. If $L_1 \subseteq M_+$, i.e., $1 \in I_+(\orL, \orM)$, then we use \eqref{eq:rhomidformula2} for the left vertical map and \eqref{eq:rhomidformula3} for the right vertical map. We calculate
\begin{multline*}
\I^{\orM}_{\s}\Bigl(\rho\bigl(U_{1,2}^{m} \cdot \prod_{i\in I_+(\orL, \orM)} \prod_j V_{i,j}^{n_{i,j}}  Y_{i,j}^{a_{i,j}}  (\x \times r)\bigr)\Bigr) \\
= \begin{cases}
 \I^{\orM}_{\s}(U_{1,1}^{m}\cdot \prod_{j=3}^{p_1}V_{i,j}^{n_{i,j}} Y_{i,j}^{a_{i,j}} \cdot \prod_{i \in I_+(\orL, \orM), i\neq 1}  \prod_jV_{i,j}^{n_{i,j}} Y_{i,j}^{a_{i,j}} \x) & \text{if } a_{1,2}=0, r=x_+, \\
0 & \text{otherwise,}
\end{cases}
\end{multline*}
which gives
$$ U_{1,1}^{m+A_1(\x) - s_1 - \sum_{j=3}^{p_1} n_{1,j}} U_{1,3}^{n_{1,3}} \dots U_{1,p_1}^{n_{1,p_1}} \cdot \prod_{i \in I_+(\orL, \orM), i \neq 1} U_{i,1}^{A_i(\x) - s_i - \sum_{j=2}^{p_i} n_{i,j}} U_{i,2}^{n_{i,2}} \dots U_{i,p_i}^{n_{i,p_i}} \x $$
if $r=x_+ \text{ and all } a_{i,j} = 0 \text{ for } i \in I_+(\orL, \orM),$ and gives zero otherwise.
 This is the same answer as
\begin{multline*}
\rho\Bigr (\I^{\orM}_{\s} \bigl(U_{1,2}^{m} \cdot \prod_{i\in I_+(\orL, \orM)} \prod_j V_{i,j}^{n_{i,j}}  Y_{i,j}^{a_{i,j}}  (\x \times r)\bigr)\Bigr)  \\
=\begin{cases}
\rho\Bigl(U_{1,2} ^{m} \cdot \prod_{i \in I_+(\orL, \orM)} U_{i,1}^{A_i(\x) - s_i - \sum_{j=2}^{p_i} n_{i,j}} U_{i,2}^{n_{i,2}} \dots U_{i,p_i}^{n_{i,p_i}} (\x \times r)\Bigr)  & \text{if all } a_{i,j} = 0 \text{ for } i \in I_+(\orL, \orM), \\
 0 & \text{otherwise.}\end{cases}
\end{multline*}
Hence, the diagram \eqref{eq:commuterhoint} commutes.

Next, we look at descent maps of the form \eqref{eq:descentgeneral}. When $M$ is a knot, we seek a diagram
\begin{equation}
\label{eq:commutedescent}
\xymatrix{
\Chain^-(\Hyper^L, \orM, p^{\orM}(\s)) \ar[r]^-{D^{\orM}_{p^{\orM}(\s)}} \ar[d]_{\rho} \ar[dr] & \Chain^-(\Hyper^{L-M}, \psi^{\orM}(\s)) \ar[d]^{\rho} \\
\Chain^-(\bHyper^L, \orM, p^{\orM}(\s)) \ar[r]^-{D^{\orM}_{p^{\orM}(\s)}} & \Chain^-(\bHyper^{L-M}, \orM, \psi^{\orM}(\s)).
}
\end{equation}
that commutes up to the homotopy given by the diagonal map. For general $M$, we need higher chain homotopies that fit into a chain map between hypercubes.

Recall from Section~\ref{sec:gens} that the map 
$$D^{\orM}_{p^{\orM}(\s)}: \Chain^-(\Hyper^L, \orM, p^{\orM}(\s)) \to\Chain^-(\Hyper^{L-M}, \psi^{\orM}(\s)) $$ 
is the largest diagonal in the compression of a hyperbox 
$$ C = \Cc(\Hyper^{\orL, \orM}, \psi^{\orM}(\s)).$$
Similarly, the corresponding map $D^{\orM}_{p^{\orM}(\s)}$ coming from $\bHyper$ is obtained by compressing a hyperbox
$$ \bC = \Cc(\bHyper^{\orL, \orM}, \psi^{\orM}(\s)).$$

The relation between the hyperboxes of chain complexes $\Hyper^{\orL, \orM}$ and $\bHyper^{\orL, \orM}$ was described in Section~\ref{sec:moves}. When $L_1 \not \subseteq M_-,$ these two hyperboxes have the same size, and the diagrams in  $\Hyper^{\orL, \orM}$ are obtained from the corresponding ones in $\bHyper^{\orL, \orM}$ by an index zero/three link stabilization. As such, there are equivalences $\rho$ that relate the corresponding complexes in the hyperboxes $\Hyper^{\orL, \orM}$ and $\bHyper^{\orL, \orM}$. These equivalences commute (up to chain homotopy) with the polygon maps in the hyperboxes, by Propositions~\ref{prop:StabPolygon} and \ref{prop:PolyDestab03a}. Thus, the maps $\rho$ commute with $D^{\orM}_{p^{\orM}(\s)}$ up to chain homotopy, and we obtain commuting diagrams of the form \eqref{eq:commutedescent}.

The more interesting case is when $L_1 \subseteq M_-$. Then, as explained in Section~\ref{sec:moves}, the hyperbox $\Hyper^{\orL, \orM}$ is obtained from $\bHyper^{\orL, \orM}$ by stabilizing all diagrams as in Figure~\ref{fig:variantfree}, and then also increasing the length of the hyperbox by one in the direction $i=1$, by adding the move shown in Figure~\ref{fig:StabTriple0}. Furthermore, recall from Equation~\eqref{eq:di} that the length of the hyperbox $C$ in the direction $i=1$ is $d_1 + p_1-1$, where $d_1$ is the length of $\Hyper^{\orL, \orM}$ in that direction. A similar formula holds for $\bC$ and $\bHyper^{\orL, \orM}$. Since there is one extra pair of basepoints in $\Hyper^{\orL, \orM}$ than in $\bHyper^{\orL, \orM}$, we deduce that the first side length of $C$ is two more than that of $\bC$. 

To be able to compare the hyperboxes, we replace $\bC$ by a new hyperbox $\tC$, which is obtained  from $\bC$ by two elementary enlargements, as in Section~\ref{sec:ele} both in the direction $i=1$, with one done at the beginning and one at the end. The beginning one has to do with the additional basepoint (and thus additional transition maps) present in $C$, and the ending one with the move from Figure~\ref{fig:StabTriple0}.

Now, the hyperboxes $C$ and $\tC$ have the same size. We claim that we can construct a chain map $$Z: C \to \tC$$ between these hyperboxes, in the sense of Definition~\ref{def:chmap}, such that the maps $Z^\zero_\eps=\rho: C^{\eps} \to \tC^{\eps}$ are chain homotopy equivalences over $\bar \Ring$. (Recall that a chain map between hyperboxes consists of maps $Z^\zero_\eps$ that preserve the index $\eps$, as well as higher diagonal maps $Z^{\eps'-\eps}_{\eps}: C^{\eps} \to \tC^{\eps'}$ when $\eps' > \eps$ are neighbors.) After compression, our map $Z$ would produce a chain map between the resulting hypercubes. Furthermore, by Lemma~\ref{lemma:ci}, elementary enlargements leave the compressed hypercubes unchanged. We would thus get the desired chain map between the compressions of $C$ and $\bC$, providing diagonal maps in  \eqref{eq:commutedescent}.

We start by constructing the quasi-isomorphisms $\rho: C^{\eps} \to \tC^{\eps}$. Recall from Section~\ref{sec:gens} that $C^{\eps}$ looks like
\begin{itemize}
\item the resolution $\Chain^-(\cdot, \psi^{\orM}(\s))$ in the directions $i \not \in I(\orL, \orM)$,
\item the complexes $\Cc_j$ from \eqref{eq:transitions}, with $j=\eps_i$, in the directions $i \in I(\eps)$, 
\item the complex $\CFm=C\{z\}$ in the directions $i \in I_-(\orL, \orM) - I(\eps)$,
\item the complex $\CFm=C\{w\}$ in the directions $i \in I_+(\orL, \orM)$. 
\end{itemize}
A similar description applies to $\tC^{\eps}$. In our case, since $L_1 \subseteq M_-$, we have either $1 \in I(\eps)$ or $1 \in  I_-(\orL, \orM) - I(\eps)$, according to whether $\eps_1 < p_1-1$ or not. 

Observe that $\eps_1$ can take values from $0$ to $d_1'=d_1+p_1-1.$ 
If $0 < \eps_1 < d_1'$, then (regardless of whether $\eps_1$ is less or greater than $p_1-1$), the
basepoint $w_{2,1}$ does not play a role in the complex $C^{\eps}$, and we are in the situation of Figure~\ref{fig:variantfree}. Thus, as in \eqref{eq:secondcone}, we can identify $C^{\eps}$ with the cone
$$ 
 (\tC^{\eps})^{U_{1,1} \to U_{1,2}}_+[[U_{1,1}]] \xrightarrow{U_{1,1} - U_{1,2}}  (\tC^{\eps})^{U_{1,1} \to U_{1,2}}_-[[U_{1,1}]].
 $$
We get an equivalence $\rho: C^{\eps} \to \tC^{\eps}$ given by
$$ \rho (U_{1,2}^{m}  (\x \times r)) = \begin{cases}
U_{1,1}^m \x & \text{if }  r=x_-, \\
0 &  \text{if }  r=x_+.
\end{cases}
$$

On the other hand, when $\eps_1=0$, then the complex $C^{\eps}$ is obtained from $\tC^{\eps}$ by taking cones on $V_{1,2}-1$ and $U_{1,2} + U_{1,1}V_{1,2}$, as in \eqref{eq:coneCint}, and therefore we have an equivalence $\rho$ given by the formula~\eqref{eq:rhomidformula2}:
$$  \rho(U_{1,2}^{m} V_{1,2}^n Y_{1,2}^a (\x \times r)) = \begin{cases}
U_{1,1}^m \x & \text{if } a=0, r=x_+, \\
0 & \text{otherwise.}
\end{cases}
$$

Lastly, when $\eps_1=d_1'$, then $C^{\eps}$ is constructed from the Heegaard diagram at the end of the move in Figure~\ref{fig:StabTriple0}, i.e., using the $\alpha_1$ and $\gamma_1$ curves. This is  just a free index zero/three stabilization of the diagram that gives $\tC^{\eps}$, and therefore we can identify $C^{\eps}$ with the cone
$$  \tC^{\eps}_+[[U_{1,2}]] \xrightarrow{U_{1,1} - U_{1,2}}  \tC^{\eps}_-[[U_{1,2}]].$$
We obtain an equivalence $ \rho: C^{\eps} \to \tC^{\eps}$ given by
$$ \rho (U_{1,2}^{m}  (\x \times r)) = \begin{cases}
U_{1,1}^m \x & \text{if }  r=y_-, \\
0 & \text{if }  r=y_+.
\end{cases}
$$

We now construct the other components of the chain map $Z$, i.e., $Z^{\eps'-\eps}_{\eps}$  for $\eps < \eps'$.  

When $\eps_1 > 0$, we simply set $Z^{\eps'-\eps}_{\eps} = 0$. If we denote by $D^{\eps'-\eps}_{\eps}$ and $\tD^{\eps'-\eps}_{\eps}$ the maps in the hyperboxes $C$ and $\tC$, we claim that the diagrams
$$\xymatrix{
C^{\eps} \ar[r]^-{D^{\eps'-\eps}_{\eps}} \ar[d]_{\rho} & C^{\eps'} \ar[d]^{\rho} \\
\tC^{\eps} \ar[r]^-{\tD^{\eps'-\eps}_{\eps}} & \tC^{\eps'} }$$
commute on the nose. Indeed, when $1 < \eps_1 < d_1'-1$, this  follows from the study of holomorphic polygons in the stabilized diagram; cf. Proposition~\ref{prop:PolyDestab03b}. For $\eps_1=d_1'-2$ and $\eps'_1= d_1'-1$, we can use the description of holomorphic polygons under the move in Figure~\ref{fig:StabTriple0}, that is, Proposition~\ref{prop:PolyDestabMove}. When $\eps_1=\eps'_1=d_1'-1$, we use the similar description of a free index zero/three stabilization; see Proposition~\ref{prop:StabPolygonH}.
 
When $\eps_1=0$, we again set  $Z^{\eps'-\eps}_{\eps} = 0$ (for $\eps < \eps'$), with one exception: if $\eps'$ differs from $\eps$ only in position $i=1$, where $\eps'_1=1$ (in other words, if we have $\eps' =\eps+\tau_1$, in the notation of Section~\ref{sec:hyperv}.) When $\eps_1=0$ and $\eps' =\eps+\tau_1$, we set
$$
 Z^{\eps'-\eps}_{\eps}(U_{1,2}^m V_{1,2}^n Y_{1,2}^a (\x \times r)) = \begin{cases}
mU_{1,1}^{m-1}\x & \text{if } a=0, r=x_-,\\
0 & \text{otherwise.}\end{cases}
$$
This is exactly analogous to the situation described in Section~\ref{sec:toylink}, where we had the formula \eqref{eq:Zorro} for the diagonal map $Z$. Checking that $Z^{\eps'-\eps}_{\eps}$ fits into a chain map between the hypercubes in $C$ and $\tC$ is entirely similar to the proof of Equation~ \eqref{eq:checkZ} in Section~\ref{sec:toylink}.
 
 This completes the construction of the chain map between $C$ and $\tC$. From here we obtain the desired quasi-isomorphism between the surgery complexes $\C^-(\Hyper, \Lambda)$ and $\C^-(\bHyper, \Lambda)$.
\end{proof}

\section {Beyond the surgery theorem}
\label {sec:beyond}

We discuss here several extensions of Theorem~\ref{thm:FirstSurgery}. 

\subsection {Maps induced by surgery} \label{sec:surgerymaps} We work in the setting of Section~\ref{sec:proofgen}, with $\orL \subset Y$ being a link in an integral homology three-sphere, and $\Hyper$ a complete system of hyperboxes for $\orL$. Let $L' \subseteq L$ be a sublink, with the orientation induced from $\orL$. The hyperboxes $\Hyper^{\orL'', \orM}$ with $M \subseteq L'' \subseteq L'$ form a complete system of hyperboxes for $\orL'$, which we denote by $\Hyper|_{L'}$. 

Following the notation from Section~\ref{sec:4c}, we let $W_{\Lambda}(L', L)$ be the  cobordism from $Y_{\Lambda|_{L'}}(L')$ to $Y_{\Lambda}(L)$ given by surgery on $L - L'$ (framed with the restriction of $\Lambda$). Recall that in Lemma~\ref{lem:spc4} we established an identification:
$$\spc ( W_{\Lambda}(L', L)) \cong \H(L)/H(L, \Lambda|_{L'})$$
such that the natural projection 
$$\pi^{L, L'} : \bigl( \H(L)/H(L, \Lambda|_{L'})\bigr) \longrightarrow \bigl(\H(L)/H(L, \Lambda)\bigr)$$
corresponds to restricting the $\spc$ structures to $Y_{\Lambda}(L)$, and the map
$$ \psi^{L-L'}: \bigl(\H(L)/H(L, \Lambda|_{L'})\bigr) \to \bigl(\H(L')/H(L', \Lambda|_{L'})\bigr) $$
corresponds to restricting them to $Y_{\Lambda|_{L'}}(L')$.

\begin {lemma}
\label {lemma:bij}
Fix an equivalence class $\t \in \H(L)/H(L, \Lambda|_{L'})$. Then the map
\begin {equation}
\label {eq:psimap}
\psi^{L-L'} : \{\s \in \H(L) \mid [\s] = \t \} \longrightarrow \{\s' \in \H(L') \mid  [\s'] = \psi^{L-L'}(\t)\} 
 \end {equation}
is always surjective, and it is injective if and only if every component $L_j \subseteq L-L'$ is rationally null-homologous inside the surgered manifold $Y_{\Lambda|_{L'}} (L')$.
\end {lemma}

\begin {proof}
Using affine identifications between spaces of $\spc$ structures and second cohomology, surjectivity and injectivity of \eqref{eq:psimap} are equivalent to the same conditions for the restriction map: 
\begin {equation}
\label {eq:kermap}
\ker\bigl (H^2(W_{\Lambda}(L)) \to H^2(W_{\Lambda}(L', L)) \bigr) \longrightarrow \ker\bigl (H^2(W_{\Lambda|_{L'}}(L')) \to H^2(Y_{\Lambda|_{L'}}(L')\bigr).
\end {equation}

Using long exact sequences for pairs, we can rewrite \eqref{eq:kermap} as
$$
\im\bigl (H^2(W_{\Lambda}(L), W_{\Lambda}(L', L)) \to H^2(W_{\Lambda}(L))  \bigr) \longrightarrow \im\bigl (H^2(W_{\Lambda|_{L'}}(L') ,Y_{\Lambda|_{L'}}(L') ) \to H^2(W_{\Lambda|_{L'}}(L')) \bigr).
$$

Surjectivity of the last map follows directly from the existence of a commutative diagram:
$$\begin {CD}
H^2(W_{\Lambda}(L), W_{\Lambda}(L', L)) @>>> H^2(W_{\Lambda}(L))  \\
@V{\cong}VV @VVV \\
H^2(W_{\Lambda|_{L'}}(L') ,Y_{\Lambda|_{L'}}(L')) @>>> H^2(W_{\Lambda|_{L'}}(L')).
\end {CD}$$

Injectivity of \eqref{eq:kermap} is equivalent to injectivity of the combined restriction map 
\begin {equation}
\label {eq:kermap3}
 H^2(W_{\Lambda}(L)) \longrightarrow H^2(W_{\Lambda}(L', L)) \oplus H^2(W_{\Lambda|_{L'}}(L')).
 \end {equation}

In turn, using the Mayer-Vietoris sequence in cohomology, injectivity of \eqref{eq:kermap3} is equivalent to surjectivity of the map
\begin {equation}
\label {eq:kermap4}
H^1(W_{\Lambda}(L', L)) \oplus H^1(W_{\Lambda|_{L'}}(L'))\longrightarrow H^1(Y_{\Lambda|_{L'}}(L')).
\end {equation}

Since $H^1(W_{\Lambda|_{L'}}(L')) = 0$ and $H^1$ groups have no torsion, we can rephrase surjectivity of \eqref{eq:kermap4} as injectivity of the dual map with rational coefficients:
$$ H_1(Y_{\Lambda|_{L'}}(L'); \qq) \longrightarrow H_1(W_{\Lambda}(L', L); \qq).$$

In turn, this is equivalent to the vanishing of the boundary map
\begin {equation}
\label {eq:kermap5}
H_2( W_{\Lambda}(L', L), Y_{\Lambda|_{L'}}(L'); \qq)  \longrightarrow H_1(Y_{\Lambda|_{L'}}(L'); \qq).
\end {equation}

Since the domain of \eqref{eq:kermap5} is generated by the cores of the $2$-handles attached along link components $L_j \subseteq L-L'$, the conclusion follows.
\end {proof}

Observe that, for every equivalence class $\t \in \H(L)/H(L, \Lambda|_{L'})$, 
$$ \widetilde{\C}^-(\Hyper, \Lambda)^{L', \t} := \bigoplus_{\{M \mid L - L' \subseteq M \subseteq L\} } \prod_{\{\s \in \H(L)| [\s] = \t\}}  \Chain^-(\Hyper^{L - M}, \psi^{M}(\s)),$$
is a subcomplex of $\C^-(\Hyper, \Lambda, \pi^{L, L'}(\t)) \subseteq \C^-(\Hyper, \Lambda)$. 

If the injectivity condition in Lemma~\ref{lemma:bij} holds, so that the map \eqref{eq:psimap} is a bijection, then the complex $ \widetilde{\C}^-(\Hyper, \Lambda)^{L', \t} $ is isomorphic to 
$$ \C^-(\Hyper|_{L'}, \Lambda|_{L'}, \psi^{L-L'}(\t))  = \bigoplus_{M' \subseteq L'} \prod_{\{\s' \in \H(L')| [\s'] = \psi^{L-L'}(\t)\}} \Chain^-(\Hyper^{L' - M'}, \psi^{M'}(\s')). $$
Indeed, the isomorphism is induced by taking $M$ to $M' = M - (L-L')$ and $\s$ to $\s'=\psi^{L - L'}(\s)$. 

In the general case, even if the map \eqref{eq:psimap} is not bijective, we still have an injective chain map
$$ \iota: \C^-(\Hyper|_{L'}, \Lambda|_{L'}, \psi^{L-L'}(\t)) \to \widetilde{\C}^-(\Hyper, \Lambda)^{L', \t}$$
given by summing up all the natural identifications 
$$\Chain^-(\Hyper^{L' - M'}, \psi^{M'}(\s')) \xrightarrow{\phantom{bla} \cong \phantom{bla} } \Chain^-(\Hyper^{L - M}, \psi^{M}(\s))$$ such that $M' = M - (L-L')$, $[\s] = \t$ and $\s'=\psi^{L - L'}(\s)$. 

We denote by $\C^-(\Hyper, \Lambda)^{L', \t}$ the image of the inclusion $\iota$. This is a subcomplex of $\C^-(\Hyper, \Lambda, \pi^{L, L'}(\t))$ isomorphic to $\C^-(\Hyper|_{L'}, \Lambda|_{L'}, \psi^{L-L'}(\t))$. When the injectivity condition in Lemma~\ref{lemma:bij} holds, $\C^-(\Hyper, \Lambda)^{L', \t}$ coincides with the previously defined $\widetilde{\C}^-(\Hyper, \Lambda)^{L', \t}$.

\begin {example}
Suppose $L = L_1 \cup L_2$ is a two-component link, with $\lk(L_1, L_2) = 1$. Consider $(0,0)$-surgery on $L$, so that $\Lambda = \begin{pmatrix}
0 & 1\\
1 & 0
\end {pmatrix}$. Let $L' = L_1$. Pick $s' \in \zz$ and set $$\t= [(s' + \tfrac{1}{2}, \tfrac{1}{2})] \in \H(L)/H(L, \Lambda|_{L'}) \cong (\zz + \tfrac{1}{2})^2/(0,1).$$ Then $\C^-(\Hyper|_{L'}, \Lambda|_{L'}, \psi^{L-L'}(\t))$ is the complex associated to $0$-surgery on $L_1$, and consists of two generalized Floer complexes:
$$ \Chain := \Chain^-(\Hyper^{L_1}, s'), \ \ \ \Bain := \Chain^-(\Hyper^{\emptyset}),$$
related by the sum of the maps $\Phi^{+L_1}_{s'}$ and $\Phi^{-L_1}_{s'}$. On the other hand, $\widetilde{\C}^-(\Hyper, \Lambda)^{L', \t}$ is the direct product of infinitely many copies of $\Chain$ and $\Bain$, related by maps $\Phi^{+L_1}_{s'}$ and $\Phi^{-L_1}_{s'}$ arranged in a zigzag, as in Figure~\ref{fig:zigzag}. The inclusion $\iota$ is similar to the one defined in $\eqref{eq:jota}$ from Section~\ref{sec:x}, except now the target has infinitely many factors.

\begin{figure}
\begin{center}
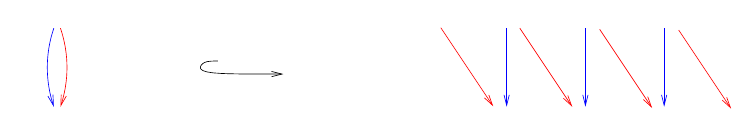
\end{center}
\caption {{\bf The inclusion of $\C^-(\Hyper|_{L'}, \Lambda|_{L'}, \psi^{L-L'}(\t))$ into $ \widetilde{\C}^-(\Hyper, \Lambda)^{L', \t}$.} The blue and red  arrows are copies of the maps $\Phi^{+L_1}_{s'}$ and $\Phi^{-L_1}_{s'}$, respectively.
The inclusion $\iota$ is given by sending each generator $\x \in \Chain$ to $(\cdots, \x, \x, \x, \cdots)$ in the infinite direct product on the right hand side; and similarly for generators $\y \in \Bain$. The image of $\iota$ is denoted $\C^-(\Hyper, \Lambda)^{L', \t}$.} 
\label{fig:zigzag}
\end{figure}

\end {example}

Theorem~\ref{thm:FirstSurgery} implies that the homology of $ \C^-(\Hyper|_{L'}, \Lambda|_{L'}, \psi^{L-L'}(\t))$, and hence also the homology of $\C^-(\Hyper, \Lambda)^{L', \t}$, are isomorphic to 
$$ \HFm_{*}(Y_{\Lambda|_{L'}}(L'),  \t|_{Y_{\Lambda|_{L'}}(L')} ).$$

In \cite{HolDiskFour}, the authors associated a map $F^-_{W, \t}$ to any cobordism $W$ between connected three-manifolds, and $\spc$ structure $\t$ on that cobordism. In the case when the cobordism $W$ consists only of two-handles (i.e. is given by integral surgery on a link), the following theorem gives a way of looking at the map $F^-_{W, \t}$ in terms of complete systems of hyperboxes:

\begin {theorem}
\label {thm:Cobordisms}
Let $\orL \subset Y$ be a link in an integral homology three-sphere, $L' \subseteq L$, a sublink,  $\Hyper$ a complete system of hypercubes for $\orL$, and $\Lambda$ a framing of $L$. Then, for any $\t \in \spc ( W_{\Lambda}(L', L)) \cong \H(L)/H(L, \Lambda|_{L'})$, the following diagram commutes:
$$\begin {CD}
H_*(\C^-(\Hyper, \Lambda)^{L', \t}) @>{\phantom{F^-_{W_{\Lambda}(L', L), \t} \otimes \Id  }}>> H_*(\C^-(\Hyper, \Lambda, \pi^{L, L'}(\t))) \\
 @V{\cong}VV @VV{\cong}V \\
\HFm_{*}(Y_{\Lambda|_{L'}}(L'), \t|_{Y_{\Lambda|_{L'}}(L')})  @>{F^-_{W_{\Lambda}(L', L), \t}}>> \HFm_{*}(Y_{\Lambda}(L), \t|_{Y_{\Lambda|_{L}}(L)}).
\end {CD}$$
Here, the top horizontal map is induced from the inclusion of chain complexes, while the two vertical isomorphisms are the ones from Theorem~\ref{thm:FirstSurgery}. 
\end {theorem}

\begin {proof}
We first discuss the proof in the case when $\Hyper$ is a basic system. The argument is similar to the one in \cite[Theorem 4.2]{IntSurg}; basically, one has to keep track of the surgery maps all throughout the arguments in Section~\ref{sec:proof}. The key point is to find commutative diagrams that relate the inclusion maps of sub-hypercubes of $\hyp^\delta$ to cobordism maps between the respective Floer complexes. This is done by applying the results of Section~\ref{sec:x} (precisely, Propositions~\ref{prop:cob+}, \ref{prop:cob-}, \ref{prop:cob0cor} and \ref{prop:cob00twisted}) repeatedly, as we follow the iteration process in the proof of Proposition~\ref{prop:relgr}. 

We emphasize that, even though in the proof of Proposition~\ref{prop:relgr} we have chosen a particular ordering of the components of $L$ (such that if $\Lambda$ is degenerate, $\Lambda_1$ is in the span of the other framing vectors, etc.), regardless of the ordering we can apply the results of Section~\ref{sec:x} to obtain an identification of the cobordism maps associated to $W_{\Lambda}(L',L)$ for {\em any} sublink $L' \subset L$. For the sake of concreteness, we explain how this works in the case of a link of two components $L = L_1 \cup L_2$.

Suppose that the framing matrix $\Lambda$ is degenerate, so we had to choose the ordering of the components such that $\Lambda_1$ is a multiple of $\Lambda_2$. (This is a constraint when $\Lambda_1$ is the zero vector.) To pick a particular situation, let us assume that $\Lambda_1$ is indeed zero (so, in particular, $L_1$ and $L_2$ have zero linking number), whereas the surgery coefficient of $L_2$ is positive. In the iteration process that leads to Proposition~\ref{prop:hyp1d} we have a diagram of maps:

\begin {equation}
\label {eq:great}
\begin {CD}
\CFmd(Y_{\Lambda_1}(L_1); \TR_2) @>{g_1^\delta}>> \CFmd(Y_{\Lambda_1+ m_1\tau_1}(L_1); \TR_2) @>{g_2^\delta}>>  \CFmd(Y; \TR) \\
@VV{j^\delta_3}V @VV{k_3^\delta}V @VV{l^\delta_3}V \\
\CFmd(Y_{\Lambda}(L)) @>{f^\delta_1}>> \CFmd(Y_{\Lambda+ m_1\tau_1}(L)) @>{f^\delta_2}>>  \CFmd(Y_{\Lambda_2}(L_2); \TR_1) \\
@. @VV{k_1^{\delta}}V @VV{l_1^{\delta}}V \\
@. \CFmd(Y_{\tilde \Lambda}(L)) @>{h_2^\delta}>> \CFmd(Y_{ \Lambda_2 + m_2 \tau_2}(L_2); \TR_1) \\
@. @VV{k_2^\delta}V @VV{l_2^\delta}V \\
@. \CFmd(Y_{ \Lambda_1 + m_1\tau_1}(L_1); \TR_2) @>{g_2^\delta}>> \CFmd(Y; \TR). 
\end {CD}
\end {equation}
We have denoted $\TR_1= \ff[T_1]/(T_1^{m_1} -1), \TR_2 = \ff[T_2]/(T_2^{m_2} - 1)$, and we have dropped the chain homotopies (corresponding to diagonals in the squares above) from notation for simplicity. Note that the bottom two rows in \eqref{eq:great} form the complex $\CC^\delta$. 

In the proof of Proposition~\ref{prop:relgr}, we first look at a quasi-isomorphism $(f_1^\delta, H_1^\delta)$ from $\CFmd(Y_{\Lambda}(L))$ to the mapping cone $Cone(f_2^{\delta})$, and then at a quasi-isomorphism from the latter to the mapping square $\CC^\delta$. Let $\ux$ be a $\spc$ structure on $Y_{\Lambda}(L)$. An application of Proposition~\ref{prop:gr00} shows that the first quasi-isomorphism above decomposes into a direct sum of several quasi-isomorphisms; one of them relates $\CFmd(Y_{\Lambda}, \ux)$ to a mapping cone $Cone(f_{2, \ux}^\delta)$, in a grading-preserving fashion. Then,  
a double application of Proposition~\ref{prop:gr+} gives a grading-preserving quasi-isomorphism from 
$Cone(f_{2, \ux}^\delta)$ to a direct summand $\CC^\delta_{\ux} \subset \CC^\delta$. 

Our new claim is that under the composition of these quasi-isomorphisms, the inclusion of the subcomplex $Cone(g_{2, \ux}^\delta)$ (resp. $Cone(l_{2, \ux}^\delta)$) into $\CC^\delta_{\ux}$ corresponds to a cobordism map (coming from a unique $\spc$ structure) from $Y_{\Lambda_1}(L_1)$ to $Y_{\Lambda}(L)$ (resp. from $Y_{\Lambda_2}(L_2)$ to  $Y_{\Lambda}(L)$).

Indeed, in the case of $Cone(g_{2, \ux}^\delta)$, a double application of Proposition~\ref{prop:cob+} gives a commutative diagram between its inclusion into $\CC^\delta_{\ux}$ and a map (consisting of $k_{3,\ux}^\delta, l_{3, \ux}^\delta$ and a diagonal chain homotopy) from $Cone(g_{2, \ux}^\delta)$ to $Cone(f_{2, \ux}^\delta)$. Then, another commutative diagram relates the latter map to a cobordism map from $Y_{\Lambda_1}(L_1)$ to $Y_{\Lambda}(L)$, using a double application of Proposition~\ref{prop:gr00twisted}.

In the case of $Cone(l_{2, \ux}^\delta)$, a double application of Proposition~\ref{prop:gr+} gives a commutative diagram between its inclusion into $\CC^\delta_{\ux}$ and the inclusion of a summand of $\CFmd(Y_{\Lambda_2}(L_2); \TR_1)$ into $Cone(f_{2, \ux}^\delta)$. Then, another commutative diagram relates the latter inclusion to a cobordism map from $Y_{\Lambda_1}(L_2)$ to $Y_{\Lambda}(L)$, by applying Proposition~\ref{prop:cob00}.

It is straightforward to extend this argument to links of several components (and arbitrary sublinks). This leads to a proof of Theorem~\ref{thm:Cobordisms} in the case when the complete system $\Hyper$ is basic.

For general complete systems, note that the quasi-isomorphisms used in the proof of Theorem~\ref{thm:FirstSurgery} in Sections~\ref{sec:invariance} and \ref{sec:proofgen} respect the inclusion maps. We obtain the desired commutative diagram, except that a priori, the bottom row is a more general cobordism map than the one considered in \cite{HolDiskFour}. More precisely, it counts holomorphic triangles between multi-pointed Heegaard diagrams (for the respective three-manifolds) that may have more than one basepoint; the original cobordism maps $F^-_{W_{\Lambda}(L', L), \t}$ as defined in \cite{HolDiskFour}, were going between singly pointed diagrams. Nevertheless, after some handleslides, isotopies and index one/two stabilizations and destabilizations, we can arrange so that the multi-pointed Heegard triple diagrams involved are all obtained from basic ones by a sequence of index zero/three (free and link) stabilizations. The fact that the bottom row can be identified with $F^-_{W_{\Lambda}(L', L), \t}$ then follows from Propositions~\ref{prop:StabPolygonH} and \ref{prop:PolyDestab03a}.
\end {proof}

\subsection {Other versions} The chain complex $\C^-(\Hyper, \Lambda, \ux)$ from Section~\ref{subsec:surgery} was constructed so that the version of Heegaard Floer homology appearing in Theorem~\ref{thm:Surgery} is $\HFm$. We now explain how one can construct similar chain complexes $\hat \C(\Hyper, \Lambda, \ux), \C^+(\Hyper, \Lambda, \ux)$ and $\C^{\infty}(\Hyper, \Lambda, \ux)$, corresponding to the theories $\widehat{\iHF}$,  $\iHF^+$ and $\HFinf$. 

The chain complex $\hat \C(\Hyper, \Lambda, \ux)$ is simply obtained from $\C^-(\Hyper, \Lambda, \ux)$ by setting one of the variables $U_i$ equal to zero. Its homology computes $\widehat{\iHF}(Y_{\Lambda}(L), \ux)$.

The chain complex  $\C^\infty(\Hyper, \Lambda, \ux)$ is obtained from  $\C^-(\Hyper, \Lambda, \ux)$ by inverting all the $U_i$ variables. It is a module over the ring of Laurent semi-infinite polynomials 
$$\Ring^{\infty} = \Field[[U_1, \dots, U_k; U_1^{-1},\dots, U_k^{-1}] = (U_1, \dots, U_k)^{-1} \Ring. $$
In other words, $\Ring^{\infty}$ consists of those power series in $U_i$'s that are sums of monomials with degrees bounded from below.

Note that $C^-(\Hyper, \Lambda, \ux)$ is a subcomplex of $\C^\infty(\Hyper, \Lambda, \ux)$. We denote the respective quotient complex by $\C^+(\Hyper, \Lambda, \ux)$.  Theorems~\ref{thm:FirstSurgery} and \ref{thm:Cobordisms} admit the following extension:

\begin {theorem}
\label {thm:AllVersions}
Fix a complete system of hyperboxes $\Hyper$ for an oriented, $\ell$-component link $\orL$
in an integral homology three-sphere $Y$, and fix a framing $\Lambda$ of $L$. Pick $\ux \in \spc(Y_{\Lambda}(L)) \cong \H(L)/H(L, \Lambda)$. Then, there are vertical isomorphisms and horizontal long exact sequences making the following diagram commute:
$$\begin {CD}
\cdots \to  H_*(\C^-(\Hyper, \Lambda, \ux)) @>>> H_*(\C^{\infty}(\Hyper, \Lambda, \ux))  @>>> H_*(\C^+(\Hyper, \Lambda, \ux)) \to \cdots\\
@VV{\cong}V  @VV{\cong}V @VV{\cong}V \\
\cdots \to \HFm_{*}(Y_\Lambda(L), \ux) @>>>  \HFinf_{*}(Y_\Lambda(L), \ux)  @>>>  \iHF^+_{*}(Y_\Lambda(L), \ux)  \to \cdots
 \end {CD}$$

 Furthermore, the maps in these diagrams behave naturally with respect to cobordisms, in the sense that there are commutative diagrams analogous to those in Theorem~\ref{thm:Cobordisms}, involving the cobordism maps $F^-_{W_{\Lambda(L', L)}, \tt}, F^{\infty}_{W_{\Lambda(L', L)}, \tt}, F^+_{W_{\Lambda(L', L)}, \tt}$.
\end {theorem}

\begin {proof}
Inverting the $U_i$ variables is an exact operation on modules, see for example \cite[Proposition 3.3]{AtiyahM}. Hence the quasi-isomorphisms relating $\C^-(\Hyper, \Lambda, \ux)$ and  $\CFm(Y_\Lambda(L), \ux) $ induce similar ones between the respective infinity versions. The five lemma then implies that the resulting maps between the plus versions are quasi-isomorphisms as well. Naturality with respect to the cobordism maps is clear from the construction.
\end {proof}

\begin{remark}
For the plus version, at least when $c_1(\ux)$ is torsion, we can replace direct products with direct sums in the link surgery formula, and thus obtain a formula closer in spirit to \cite{IntSurg}. Indeed, when $c_1(\ux)$ is torsion, we have a relative $\Z$-grading on the surgery complex $\C^+(\Hyper, \Lambda, \ux))$, and only finitely many terms in the direct product are nonzero in each grading level. Therefore, if we are interested in a fixed grading level, direct sums and direct products give the same answer. This implies that the plus surgery complexes constructed with direct sums are the same as those with direct products, overall. (The $U$ actions also coincide, because they coincide in each grading level.)

By contrast, recall that, for $\CFm$, we do need to use direct products, as explained in Section~\ref{sec:unknot}.
\end{remark}

\subsection {Mixed invariants of closed four-manifolds}

Let us recall the definition of the closed four-manifold invariant from \cite{HolDiskFour}. Let $X$ be a closed, oriented four-manifold with $b_2^+(X) \geq 2$. By deleting two four-balls from $X$ we obtain a cobordism $W$ from $S^3$ to $S^3$. We can cut $W$ along a three-manifold $N$ so as to obtain two cobordisms $W_1, W_2$ with $b_2^+(W_i) > 0$; further, the manifold $N$  can be chosen such that $\delta H^1(N; \Z) \subset H^2(W; \Z)$ is trivial. (If this is the case, $N$ is called an {\em admissible cut}.) Let $\tt$ be a $\spc$ structure on $X$ and $\tt_1, \tt_2$ its restrictions to $W_1, W_2$. In this situation, the cobordism maps
$$ F^-_{W_1, \tt_1} : \HFm(S^3) \to \HFm(N, \tt|_N)$$
and
$$ F^+_{W_2, \tt_2}: \iHF^+(N, \tt|_N) \to \iHF^+(S^3)$$
factor through $\iHF_{\operatorname{red}}(N, \tt|_N)$, where
$$ \iHF_{\operatorname{red}} = \coker(\HFinf \to \iHF^+) \cong \ker (\HFm \to \HFinf).$$ By composing them we obtain the mixed map
$$ F^{\operatorname{mix}}_{W, \tt}: \HFm(S^3) \to \iHF^+(S^3),$$
which changes degree by the quantity
$$ d(\tt) = \frac{c_1(\tt)^2 - 2\chi(X) - 3\sigma(X)}{4}.$$

Let $\Theta_-$ be the maximal degree generator in $\HFm(S^3)$.  Clearly the map $ F^{\operatorname{mix}}_{W, \tt}$ can be nonzero only when $d(\tt)$ is even and nonnegative. If this is the case, the value 
\begin {equation}
\label {eq:mixedOS}
\Phi_{X, \tt} = U^{d(\tt)/2} \cdot F^{\operatorname{mix}}_{W, \tt}(\Theta_-) \in \iHF^+_{0}(S^3) \cong \ff 
\end {equation}
is an invariant of the four-manifold $X$ and the $\spc$ structure $\tt$. It is conjecturally the same as the Seiberg-Witten invariant.

\begin {remark}
In \cite[Section 9]{HolDiskFour}, the mixed invariant was defined as a map
$$ \ff[U] \otimes \Lambda^*\bigl (H_1(X)/\operatorname{Tors} \bigr) \to \ff.$$
We only discuss here the value of this map at $1$, which is exactly $\Phi_{X, \tt}$ as defined in \eqref{eq:mixedOS}. 
\end {remark}

The following definition was sketched in the Introduction:

\begin {definition}
Let  $X$ be a closed, oriented four-manifold with $b_2^+(X) \geq 2$. A {\em cut link presentation  for $X$} consists of a link $L \subset S^3$, a decomposition of $L$ as a disjoint union
$$ L = L_1 \amalg L_2 \amalg L_3,$$
and a framing $\Lambda$ for $L$ (with restrictions $\Lambda_i$ to $L_i, i=1, \dots, 3$)  
with the following properties:
\begin {itemize}
\item $S^3_{\Lambda_1}(L_1)$ is a connected sum of $m$ copies of $S^1 \times S^2$, for some $m \geq 0$. We denote by $W_1$ the cobordism from $S^3$ to $\#^m (S^1 \times S^2)$ given by $m$ one-handle attachments;
\item $S^3_{\Lambda_1 \cup \Lambda_2 \cup \Lambda_3} (L_1 \cup L_2 \cup L_3)$ is a connected sum of $m'$ copies of $S^1 \times S^2$, for some $m' \geq 0$. We denote by $W_4$ the cobordism from $\#^{m'} (S^1 \times S^2)$ to $S^3$ given by $m'$ three-handle attachments;
\item If we denote by $W_2$ resp. $W_3$ the cobordisms from $S^3_{\Lambda_1}(L_1)$ to $S^3_{\Lambda_1\cup \Lambda_2}(L_1 \cup L_2)$, resp. from $S^3_{\Lambda_1\cup \Lambda_2}(L_1 \cup L_2)$ to $S^3_{\Lambda_1 \cup \Lambda_2 \cup \Lambda_3} (L_1 \cup L_2 \cup L_3)$, given by surgery on $L_2$ resp. $L_3$ (i.e. consisting of two-handle additions), then 
$$ W = W_1 \cup W_2 \cup W_3 \cup W_4$$
is the cobordism from $S^3$ to $S^3$ obtained from $X$ by deleting two copies of $B^4$;
\item The manifold $N=S^3_{\Lambda_1\cup \Lambda_2}(L_1 \cup L_2)$ is an admissible cut for $W$, i.e. $b_2^+(W_1 \cup W_2) > 0, b_2^+(W_3 \cup W_4) > 0$, and $\delta H^1(N) =0$ in $H^2(W)$.
\end {itemize}
\end {definition}

\begin {lemma}
Any closed, oriented four-manifold $X$ with $b_2^+(X) \geq 2$ admits a cut link presentation.
\end {lemma}
\begin {proof}
Start with a decomposition $W = W' \cup_N W''$ along an admissible cut. Split $W'$ into three cobordisms 
$$ W'= W_1' \cup W_2' \cup W_3'$$ 
such that $W_i'$ consists of $i$-handle additions only. It is easy to check that the decomposition
$$ W = (W_1' \cup W_2') \cup (W_3' \cup W'')$$
is still along an admissible cut. Next, split the cobordism $W_3' \cup W''$ into
$$ W''_1 \cup W_2'' \cup W_3'' ,$$
such that  $W_i'$ consists of $i$-handle additions only. Finally, adjoin the one-handles from $W_1''$ to $W_1' \cup W_2'$ and rearrange the handles to obtain a decomposition
$$ W_1' \cup W_2' \cup W_1'' = W_1 \cup W_2,$$
where $W_i, i=1,2$ consists of $i$-handle additions only. If we set $W_3 = W_2''$ and $W_4 = W_3''$, we obtain a decomposition along admissible cut of the form:
$$ W = (W_1 \cup W_2) \cup (W_3 \cup W_4),$$

We can then find a framed link $L = L_1 \cup L_2 \cup L_3$ such that surgery on $L_1$ produces the same $3$-manifold as at the end of the cobordism $W_1$ (made of one-handles), whereas surgery on $L_2$ and $L_3$ is represented by the cobordisms $W_2$ and $W_3$ (made of two-handles), respectively.
\end {proof}

\begin {definition}
Let  $X$ be a closed, oriented four-manifold with $b_2^+(X) \geq 2$. A {\em hyperbox presentation} $\Gamma$ for $X$ consists of a cut link presentation $(L = L_1 \cup L_2 \cup L_3, \Lambda)$ for $X$, together with a complete system of hyperboxes for $L$.
\end {definition}

The four-manifold invariant $\Phi_{X, \tt}$ can be expressed in terms of a hyperbox presentation $\Gamma$ for $X$ as follows. Using Theorem~\ref{thm:AllVersions}, we can express the maps $F^-_{W_2, \tt|_{W_2}}$ and $F^+_{W_3, \tt|_{W_3}}$ in terms of counts of holomorphic polygons on a symmetric product of the surface. We can combine these maps using their factorization through $\iHF_{\operatorname{red}}$, and obtain a mixed map
$$ F^{\mix}_{W_2 \cup W_3, \tt|_{W_2 \cup W_3}} : \HFm(\#^m (S^1 \times S^2)) \to \iHF^+(\#^{m'} (S^1 \times S^2)).$$

On the other hand, by composing the maps induced on homology by natural inclusions of chain complexes (of the kind used in Theorem~\ref{thm:Cobordisms}), via factoring through a reduced group
$$ \ker \bigl( H_*(\C^-(\Hyper, \Lambda)^{L_1 \cup L_2, \t|_{W_3}}) \to H_*(\C^\infty(\Hyper, \Lambda)^{L_1 \cup L_2, \t|_{W_3}}) \bigr)$$
we can construct a map
$$ F^{\mix}_{\Gamma, \tt}: H_*(\C^-(\Hyper, \Lambda)^{L_1, \tt|_{W_2 \cup W_3}}) \to H_*(\C^+(\Hyper, \Lambda)^{L_1 \cup L_2 \cup L_3, \tt|_{\#^{m'} (S^1 \times S^2)}}).$$

Theorem~\ref{thm:AllVersions} implies that $F^{\mix}_{\Gamma, \tt}$ is the same as 
$ F^{\mix}_{W_2 \cup W_3, \tt|_{W_2 \cup W_3}}$,
up to compositions with isomorphisms on both the domain and the target. Note, however, that at this point we do not know how to identify elements in the domains (or targets) of the two maps in a canonical way. For example, we know that there is an isomorphism
\begin {equation}
\label {eq:isoV}
H_*(\C^-(\Hyper, \Lambda)^{L_1, \tt|_{W_2 \cup W_3}}) \cong \HFm(\#^m (S^1 \times S^2)) ,
\end {equation}
but it may be difficult to pinpoint what the isomorphism is, in terms of $\Hyper$. Nevertheless, the good news is that $ \HFm(\#^m (S^1 \times S^2))$ has a unique maximal degree element $\Theta_{\max}^m$ . We can identify what $\Theta_{\max}^m$ corresponds to on the left hand side of ~\eqref{eq:isoV} by simply computing degrees. Let us denote the respective element by 
$$\Theta_{\max}^\Gamma \in H_*(\C^-(\Hyper, \Lambda)^{L_1, \tt|_{W_2 \cup W_3 }}). $$

The following proposition says that one can decide whether $\Phi_{X, \tt} \in \ff$ is zero or one from information in the hyperbox presentation $\Gamma:$

\begin {theorem}
\label {thm:Mixed}
Let $X$ be a closed, oriented four-manifold $X$ with $b_2^+(X) \geq 2$, with a $\spc$ structure $\tt$ with $d(\tt) \geq 0$ even. Let $\Gamma$ be a hyperbox presentation for $X$. Then $\Phi_{X, \tt} = 1$ if and only if  $U^{d(\tt)/2} \cdot F^{\mix}_{\Gamma, \tt} (\Theta_{\max}^\Gamma)$ is nonzero.
\end {theorem}
 
\begin {proof}
We have
$$ \Phi_{X, \tt} = U^{d(\tt)/2} \cdot F^{\operatorname{mix}}_{W, \tt}(\Theta_-) = F^+_{W_4, \tt|_{W_4}} (U^{d(\tt)/2} \cdot F^{\mix}_{W_2 \cup W_3, \tt|_{W_2 \cup W_3}} (F^-_{W_1, \tt|_{W_1}}(\Theta_-))).$$
By the definition of the one-handle addition maps from \cite[Section 4.3]{HolDiskFour},
$$F^-_{W_1, \tt|_{W_1}}(\Theta_-) = \Theta_{\max}^m.$$ 

Note that $  U^{d(\tt)/2}\cdot F^{\mix}_{W_2 \cup W_3, \tt|_{W_2 \cup W_3}} (\Theta_{\max}^m)$ lies in the minimal degree $k$ for which 
$$\iHF^+_k(S^3_{\Lambda_1 \cup \Lambda_2 \cup \Lambda_3} (L_1 \cup L_2 \cup L_3)) = \iHF^+_k(\#^{m'}(S^1 \times S^2))$$ is nonzero, namely $k=-m'/2$. There is a unique nonzero element in the Floer homology in that degree, which is taken to $1$ by the three-handle addition map $F^+_{W_4, \tt|_{W_4}}$, see  \cite[Section 4.3]{HolDiskFour}.

We deduce from here that $\Phi_{X, \tt} = 1$ if and only if $ U^{d(\tt)/2}\cdot F^{\mix}_{W_2 \cup W_3, \tt|_{W_2 \cup W_3}} (\Theta_{\max}^m)$ is nonzero.  The claim then follows from the fact that the maps $F^{\mix}_{\Gamma, \tt}$ and   $ F^{\mix}_{W_2 \cup W_3, \tt|_{W_2 \cup W_3}}$ are the same up to pre- and post-composition with isomorphisms.
\end {proof}

\subsection {The link surgeries spectral sequence} Our goal here will be to explain how the link surgeries spectral sequence from \cite[Section 4]{BrDCov} can be understood in terms of complete systems of hyperboxes for links in $S^3$.

We recall the main result from \cite[Section 4]{BrDCov}. Let $M = M_1 \cup \dots \cup  M_\ell$ be a framed $\ell$-component link in a 3-manifold $Y$. For each $\eps = (\eps_1, \dots, \eps_\ell) \in \E_\ell = \{0,1\}^\ell$, we let $Y(\eps)$ be the $3$-manifold obtained from $Y$ by doing $\eps_i$-framed surgery on $M_i$ for $i=1, \dots, \ell$. 

When $\eps'$ is an immediate successor to $\eps$ (that is, when $\eps < \eps'$ and $\|\eps' - \eps\| = 1$), the two-handle addition from $Y(\eps)$ to $Y(\eps')$ induces a map on Heegaard Floer homology
$$ F^-_{\eps < \eps'} : \HFm(Y(\eps)) \longrightarrow \HFm (Y(\eps')). $$

The following is the link surgery spectral sequence (Theorem 4.1 in \cite{BrDCov}, but phrased here in terms of $\HFm$ rather than $\widehat{\iHF}$ or $\iHF^+$):

\begin {theorem}[Ozsv\'ath-Szab\'o]
\label {thm:OSspectral}
There is a spectral sequence whose $E^1$ term is $\bigoplus_{\eps \in \E_\ell} \HFm(Y(\eps))$, whose $d_1$ differential is obtained by adding the maps $F^-_{\eps < \eps'}$ (for $\eps'$ an immediate successor to $\eps$), and which converges to $E^{\infty} \cong \HFm(Y)$.
\end {theorem}

To relate this to the constructions in this paper, we represent $Y(0,\dots, 0)$ itself as surgery on a framed link $(L', \Lambda')$ inside $S^3$. Let $L'_1, \dots, L'_{\ell'}$ be the components of $L'$. There is another framed link $(L=L_1 \cup \dots \cup L_\ell, \Lambda)$ in $S^3$, disjoint from $L'$, such that surgery on each component $L_i$ (with the given framing) corresponds exactly to the 2-handle addition from $Y(0, \dots, 0)$ to $Y(0, \dots, 0, 1, 0, \dots, 0)$, where the $1$ is in position $i$. For $\eps \in \E_\ell$, we denote by $L^\eps$ the sublink of $L$ consisting of those components $L_i$ such that $\eps_i = 1$.

Let $\Hyper$ be a complete system of hyperboxes for the link $L' \cup L \subset S^3$. As mentioned in Section~\ref{sec:surgerymaps}, inside the surgery complex $\C^-(\Hyper, \Lambda' \cup \Lambda)$ (which is an $(\ell'+\ell)$-dimensional hypercube of chain complexes) we have various subcomplexes corresponding to surgery on the sublinks on $L' \cup L$. We will restrict our attention to those sublinks that contain $L'$, and use the respective subcomplexes to construct a new, $\ell$-dimensional hypercube of chain complexes $\C^-(\Hyper, \Lambda' \cup \Lambda \hey L)$ as follows.

At a vertex $\eps \in \E_\ell$ we put the complex 
$$ \C^-(\Hyper, \Lambda' \cup \Lambda \hey L)^\eps = \C^-(\Hyper|_{L' \cup L^\eps}, \Lambda' \cup \Lambda|_{L^\eps}).$$

Consider now an edge from $\eps$ to $\eps' = \eps +\tau_i$ in the hypercube $\E_\ell$. The corresponding complex $\C^-(\Hyper|_{L' \cup L^\eps}, \Lambda' \cup \Lambda|_{L^\eps})$ decomposes as a direct product over all $\spc$ structures $\s$ on $Y(\eps) = S^3(L' \cup L^\eps, \Lambda' \cup \Lambda|_{L^\eps})$. As explained in Section~\ref{sec:surgerymaps}, each factor $\C^-(\Hyper|_{L' \cup L^\eps}, \Lambda' \cup \Lambda|_{L^\eps}, \s)$ admits an inclusion into $\C^-(\Hyper|_{L' \cup L^{\eps'}}, \Lambda' \cup \Lambda|_{L^{\eps'}})$ as a subcomplex. In fact, there are several such inclusion maps, one for each $\spc$ structure $\t$ on the 2-handle cobordism from $Y(\eps)$ to $Y(\eps')$ such that $\t$ restricts to $\s$ on $Y(\eps)$. Adding up all the inclusion maps on each factor, one obtains a combined map
$$G^-_{\eps < \eps'} : \C^-(\Hyper|_{L' \cup L^\eps}, \Lambda' \cup \Lambda|_{L^\eps}) \longrightarrow \C^-(\Hyper|_{L' \cup L^{\eps'}}, \Lambda' \cup \Lambda|_{L^{\eps'}}).$$

We take $G^-_{\eps < \eps'}$ to be the edge map in the hypercube of chain complexes $\C^-(\Hyper, \Lambda' \cup \Lambda \hey L)$. Since the edge maps are just sums of inclusions of subcomplexes, they commute on the nose along each face of the hypercube. Therefore, in the hypercube $\C^-(\Hyper, \Lambda' \cup \Lambda \hey L)$ we can take the diagonal maps to be zero, along all faces of dimension at least two. 

This completes the construction of $\C^-(\Hyper, \Lambda' \cup \Lambda \hey L)$. As an $\ell$-dimensional hypercube of chain complexes, its total complex admits a filtration by $-\|\eps\|$, which induces a spectral sequence; we will refer to the filtration by $-\|\eps\|$ as the {\em depth filtration} on 
$\C^-(\Hyper, \Lambda' \cup \Lambda \hey L)$.

\begin {theorem}
\label{thm:SpectralSequence}
Fix a complete system of hyperboxes $\Hyper$ for an oriented link $\orL' \cup \orL$
in $S^3$, and fix framings $\Lambda$ for $L$ and $\Lambda'$ for $L'$.  Suppose $\Hyper$ has $k$ basepoints of type $w$ and $p$ colors, and that $L$ has $\ell$ components $L_1, \dots, L_\ell$. Let $Y(0,\dots,0) = S^3_{\Lambda'}(L')$, and let $Y(\eps)$ be obtained from $Y(0,\dots,0)$ by surgery on the components $L_i \subseteq L$ with $\eps_i = 1$ (for any $\eps \in \E_\ell$). Then, there is an isomorphism between the link surgeries spectral sequence from Theorem~\ref{thm:OSspectral}, and the spectral sequence associated to the depth filtration on $\C^-(\Hyper, \Lambda' \cup \Lambda \hey L)$.   
\end{theorem}

\begin {proof}
Theorem~\ref{thm:FirstSurgery} gives a quasi-isomorphism between $\CFm (Y(1, \dots,1))  = \CFm(S^3_{\Lambda' \cup \Lambda}(L' \cup L)) $ and the surgery complex $\C^-(\Hyper, \Lambda' \cup \Lambda)$.  Let us summarize the main steps in the construction of this quasi-isomorphism. It suffices to construct the quasi-isomorphism at the level of the vertical truncations $\C^{-, \delta \from \delta'}(\Hyper, \Lambda' \cup \Lambda)$ and  $\CFmdd (Y(1, \dots,1)) $. (In the proof, we used vertical truncations by $\delta$ for torsion $\spc$ structures and by $\delta \from \delta'$ for non-torsion $\spc$ structures. However, it is clear that always truncating by $\delta \from \delta'$ works, too.)
We consider a basic system $\Hyper_b$ for $L' \cup L$. By relating $\Hyper$ to $\Hyper_b$ via some moves on complete systems, we find a quasi-isomorphism between $\Cdd(\Hyper, \Lambda' \cup \Lambda)$ and $\Cdd(\Hyper_b, \Lambda' \cup \Lambda) $. By a version of the Large Surgeries Theorem, the hypercube $\Cdd(\Hyper_b, \Lambda' \cup \Lambda)$ is shown to be quasi-isomorphic to a hypercube $\hyp^{\delta \from \delta'}$, in which at each vertex we have the truncated 
 Floer complex of some surgery on a sublink of $L' \cup L$, possibly with twisted coefficients. Finally, iterating a variant of the surgery exact triangle, we obtain a quasi-isomorphism between $\CFmdd(Y(1,\dots,1))$  and the total complex of $\hyp^{\delta \from \delta'}$. 
    
 Recall that from the $(\ell + \ell')$-dimensional hypercube  $\C^-(\Hyper, \Lambda' \cup \Lambda)$ we construct an $\ell$-dimensional hypercube $\C^-(\Hyper, \Lambda' \cup \Lambda \hey L)$. In the new hypercube, at each vertex $\eps \in \E_\ell$ we have an $\|\eps\|$-dimensional sub-hypercube of the original $\C^-(\Hyper, \Lambda' \cup \Lambda)$; along the edges we have corresponding inclusion maps, and along higher-dimensional faces the diagonal maps are trivial. We can apply an analogous procedure to the truncated hypercubes $\Cdd(\Hyper, \Lambda' \cup \Lambda), \Cdd(\Hyper_b, \Lambda' \cup \Lambda)$ and $\hyp^{\delta \from \delta'}$, and obtain $\ell$-dimensional hypercubes from them; we denote these by $\Cdd(\Hyper, \Lambda' \cup \Lambda \hey L), \Cdd(\Hyper_b, \Lambda' \cup \Lambda \hey L)$ and $\hyp^{\delta \from \delta'} \hey L$. The quasi-isomorphisms between $(\ell + \ell')$-dimensional hypercubes constructed in the proof of Theorem~\ref{thm:FirstSurgery} all preserve the corresponding depth filtrations, and thus induce quasi-isomorphisms between respective sub-hypercubes. As a consequence, we can construct filtered quasi-isomorphisms (with respect to the depth filtration) between $\Cdd(\Hyper, \Lambda' \cup \Lambda \hey L), \Cdd(\Hyper_b, \Lambda' \cup \Lambda \hey L) $ and $(\hyp^{\delta \from \delta'} \hey L) $. This implies that the corresponding spectral sequences (induced by the depth filtrations on $\ell$-dimensional hypercubes) 
 are isomorphic.
 
 It remains to find an isomorphism between the spectral sequence induced by the depth filtration on  
 the hypercube $\hyp^{\delta \from \delta'} \hey L$, and the $\delta \from \delta'$ truncation of the link surgeries spectral sequence from Theorem~\ref{thm:OSspectral}. Let us first explain how this is done in the simplest case, when $\ell =1$ so that the link $L$ has a single component $L_1 = K$. We further assume that $K$ has linking number zero with each component of $L'$. We will drop the truncation symbol $\delta \from \delta'$ from notation for simplicity (in fact, in this step of the argument the quasi-isomorphism exists also at the level of untruncated complexes). The spectral sequence from Theorem~\ref{thm:OSspectral} is simply associated to the depth filtration on a one-dimensional hypercube, which
 is the mapping cone
 \begin {equation}
 \label {eq:cone1}
  \CFm(S^3_{\Lambda'}(L')) \longrightarrow \CFm(S^3_{\Lambda' \cup \Lambda} (L' \cup K)).
 \end {equation} 
   
On the other hand, the hypercube $\hyp \hey L$ is a mapping cone
\begin {equation}
\label {eq:cone2}
 (\hyp \hey L)^0 \longrightarrow  (\hyp \hey L)^1,
 \end {equation}
where $(\hyp \hey L)^0 $ is $\CFm(S^3_{\Lambda'}(L'))$ and $(\hyp \hey L)^1$ is itself a mapping cone
\begin {equation}
\label {eq:cone3}
 \CFm(S^3_{\Lambda' \cup (\Lambda+m)} (L' \cup K)) \longrightarrow \oplus^m \CFm(S^3_{\Lambda'}(L')),
 \end {equation}
 for some $m \gg 0$. 
 
 The surgery exact triangle (Proposition~\ref{prop:leseq}) says that the second term in \eqref{eq:cone1} is quasi-isomorphic to the second term in \eqref{eq:cone2}. In fact, we can take a quasi-isomorphism  given by a triangle-counting map from $\CFm(S^3_{\Lambda' \cup \Lambda} (L' \cup K))$ to the first term in \eqref{eq:cone3}, and a quadrilateral-counting map to the second term in \eqref{eq:cone1}. We can extend this quasi-isomorphism to one between the mapping cone \eqref{eq:cone1} and the mapping cone \eqref{eq:cone2}, by taking the identity between their first terms $\CFm(S^3_{\Lambda'}(L'))$, and also adding a diagonal map from $\CFm(S^3_{\Lambda'}(L'))$ to the mapping cone \eqref{eq:cone2}. This diagonal map consists of a quadrilateral-counting map from $\CFm(S^3_{\Lambda'}(L'))$ to the first term in \eqref{eq:cone3}, and a pentagon-counting map from  $\CFm(S^3_{\Lambda'}(L'))$ to the second term in \eqref{eq:cone3}. This produces a chain map between \eqref{eq:cone1} and \eqref{eq:cone2}, which is a quasi-isomorphism because it is so on the level of the associated graded of the depth filtrations.
 
Let us now discuss how this construction can be generalized to the case when the link $L$ has an arbitrary number $\ell$ of connected components. The link surgeries spectral sequence from Theorem~\ref{thm:OSspectral} is associated to  the depth filtration on an $\ell$-dimensional hypercube, where at each vertex we have a Heegaard Floer complex $\CFm(Y(\eps))$, along the edges we have triangle-counting maps (producing the cobordism maps $F^-_{\eps < \eps'}$ on the $E^1$ page), and along the higher-dimensional faces we have higher polygon-counting maps. We can construct a filtered quasi-isomorphism between this hypercube and $\hyp \hey L$ as follows: between corresponding vertices we use the quasi-isomorphisms given iterating the surgery exact triangle (Proposition~\ref{prop:leseq}), and then we complete this to a chain map between the hypercubes, by adding further polygon-counting maps along the diagonals. A filtered quasi-isomorphism between $\ell$-dimensional hypercubes produces an isomorphism between the respective spectral sequences, and this completes the proof.
\end {proof}

\clearpage \section{The surgery theorem applied to grid diagrams with free markings}
\label {sec:grid}

In this section we state a variant of the surgery theorem, in terms of counts of polygons on grid diagrams with free markings. The proof involves applying Theorem~\ref{thm:FirstSurgery} to a special kind of complete system, associated to the grid.

\subsection {Grid diagrams with free markings}
Toroidal grid diagrams are a particular kind of Heegaard diagrams for a link in $S^3$. In \cite{MOS}, \cite{MOST}, they have been used to give combinatorial descriptions to link Floer complexes.

\begin{figure}
\begin{center}
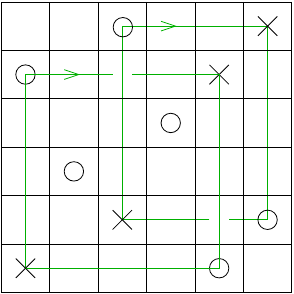
\end{center}
\caption {{\bf A grid diagram for the Hopf link, with two free markings.} The link is drawn in green.}
\label{fig:gridfree}
\end{figure}

We introduce here a slightly more general concept, that of a {\em toroidal grid diagram with free markings}. (An example is shown in Figure~\ref{fig:gridfree}.) Such a diagram $G$ consists of a torus $\TT$, viewed as a square in the plane with the opposite sides identified, and split into $n$ annuli (called rows) by $n$ horizontal circles $\alpha_1, \dots, \alpha_n$, and into $n$ other annuli (called columns) by $n$ vertical circles $\beta_1, \dots, \beta_n$. Further, we are given several markings on the torus, of two types: $X$ and $O$, such that:
\begin {itemize}
\item each row and each column contains exactly one $O$ marking;
\item each row and each column contains at most one $X$ marking;
 \item if the row of an $O$ marking contains no $X$ markings, then the column of that $O$ marking  contains no $X$ markings either. An $O$ marking of this kind is called a {\em free marking}. 
  \end {itemize}
 It follows that $G$ contains exactly $n$ $O$ markings and $n-q$ $X$ markings, where $q$ is the number of free markings. A marking that is not free is called {\em linked}. The number $n$ is called  the {\em grid number} of $G$. 
  
Given $G$ as above, we draw horizontal arcs between the (linked) markings in the same row, and vertical arcs between the markings in the same column. Letting the vertical arcs be overpasses whenever they intesect the horizontal arcs, we then obtain a planar diagram for a link $\orL \subset S^3$, which we orient so that all horizontal arcs go from an $O$ to an $X$. We denote by $\ell$ the number of components of $L$. 
 
A toroidal grid is a type of multi-pointed Heegaard diagram for $\orL$, with the horizontal circles being the alpha curves $\alpha_1, \dots, \alpha_n$, and the vertical curves being the beta curves $\beta_1, \dots, \beta_n$. We let $\S = \S(G) = \Ta \cap \Tb$, where $\Ta = \alpha_1 \times \dots \times \alpha_n, \Tb = \beta_1 \times \dots \times \beta_n$ are tori in the symmetric product $\Sym^n(\TT)$.

Following the notation of \cite{MOST}, we use $X$ and $O$ markings to play the role of the basepoints $z_i$ and $w_i$, respectively. In particular, we change the notation from Section~\ref{sec:chains} and write $X_i, O_i$ for $n_{z_i}, n_{w_i}$, respectively. Further, we let $\Xs_i$ and $\Os_i$ be the sets of $X$ and $O$ markings on a link component $L_i$. We write $\Xs = \cup_i \Xs_i$ for the set of all $X$ markings, and $\Os$ for the set of all $O$ markings (including the linked ones). We also fix an ordering of the basepoints in each set $\Xs_i, \Os_i$. Note that the free markings are exactly the free basepoints; we arrange so that they are denoted $O_{n-q+1}, \dots, O_n$.

The arguments in \cite{MOS} show that isolated holomorphic disks in $\Sym^n(\TT)$ with boundaries on $\Ta, \Tb$ are in a natural one-to-one correspondence to empty rectangles on the grid $G$. For $\x, \y \in \S$, we let $\EmptyRect(\x, \y)$ be the set of empty rectangles between $\x$ and $\y$; compare \cite{MOST}. Then, the Floer chain complex $\CFLm(G)=\CFLm(\Ta, \Tb)$ is the free module over $\Ring = \Field[[U_1, \dots, U_n]]$ generated by $\S$, and endowed with the differential:
\begin{equation}
\label {eq:du}
\partial \x = \sum_{\y \in \S} \sum_{r\in \EmptyRect(\x, \y)} U_1^{O_1(r)} \cdots  U_n^{O_n(r)} \y.
\end {equation}
We have Alexander gradings $A_i$ on $\S$ for $i=1, \dots, \ell$, which produce filtrations on $\CFLm(G)$. For  $\s = (s_1, \dots, s_\ell) \in \bH(L)$, we denote by $\Am(G, \s) \subset \CFLm(G)$ the subcomplex given by $A_i \leq s_i$, $i=1, \dots, \ell$.

We have a resolution $\Chain^-(G, \s)$ of $\Am(G, \s)$ constructed as in Section~\ref{sec:gens}. Specifically, let us re-label the basepoints on each component $L_i \subseteq L$ as
$$ O_{i,1}, X_{i,1}, O_{i,2}, X_{i,2}, \dots, O_{i, p_i}, X_{i, p_i}.$$
We also give the alternate name $U_{i,j}$ to the $U$ variable associated to $O_{i,j}$.

Consider the complex $\Ccint(G)$, which is the dg module over the dga
$$ \Ringbig^{Y} = \Field[[(U_{i,j})_{\substack{1\leq i \leq \ell, \\ 1 \leq j \leq p_i}} (U_{j})_{n-q+1\leq j \leq p}]][(V_{i,j}, Y_{i,j})_{\substack{1\leq i \leq \ell, \\ 2 \leq j \leq p_i}}]'/(Y_{i,j}^2=0, \del Y_{i,j} = U_{i,j} + U_{i,1}V_{i,j}),$$
generated by $\x \in \S$, with differential 
$$
 \del \x = \sum_{\y \in \S}  \sum_{r\in \EmptyRect(\x, \y)}  \prod_{i=1}^\ell  U_{i,1}^{X_{i,1}(r) + \ldots + X_{i,p_i}(r)} V_{i,2}^{O_{i,2}(r)} \dots  V_{i,p}^{O_{i,p}(r)} \cdot \prod_{j=n-q+1}^n U_{j}^{O_{ j}(r)} \y.
$$

This complex admits filtrations $\FF_i$, for $1 \leq i \leq \ell$, as in Section~\ref{sec:gens}. For $\s=(s_1, \dots, s_\ell) \in \H(L)$, we define $\Chain^-(\Hyper^L, \s)$ to be the filtered part of $\Ccint(\Hyper^{L})$ given by $ \FF_i \leq s_i$, $i=1, \dots, \ell.$

If $M \subseteq L$ is a sublink with an orientation $\orM$, we have a complex $\Chain^-(G, \orM, p^{\orM}(\s))$ and a projection-inclusion map
\begin {equation}
\label {eq:Proj}
 \Pr^{\orM}_\s : \Chain^-(G, \s) \to \Chain^- (G, \orM, p^{\orM}(\s)),
\end {equation}
defined as in \eqref{eq:Ims}.

We also introduce the following notation: $\Xs^{\orM} \subseteq \Xs$ will be the subset consisting of the $X$ basepoints on $L - M$, and 
$$ \Os^{\orM} := \Bigl( \Os - \bigcup_{i \in I_-(\orL, \orM)} \Os_i \Bigr) \cup  \bigcup_{i \in I_-(\orL, \orM)} \Xs_i.$$ 

Thus, $(\TT, \alphas, \betas, \Os^{\orM}, \Xs^{\orM})$ is the reduction $r_{\orM}(G)$, in the sense of Definition~\ref{def:reduce}.

\subsection {Handleslides over a set of markings}
\label {sec:DestabSeveral}
From now on we will specialize to the situation in which $G$ is a toroidal grid diagram {\em with at least one free marking.}

Consider a subset $\zed= \{Z_1, \dots, Z_k \} \subseteq \Xs \cup \Os$ consisting only of linked markings. We say that $\zed$ is {\em consistent} if, for any $i$, at most one of the sets $\zed \cap \Os_i$ and $\zed \cap \Xs_i$ is nonempty. From now on we shall assume that $\zed$ is consistent.

We let $L(\zed) \subseteq L$ be the sublink consisting of those components $L_i$ such that at least one of the markings on $L_i$ is in $\zed$. We orient $L(\zed)$ as $\orL(\zed)$, such that a component $L_i$ is given the orientation coming from $\orL$ when $\zed \cap \Os_i \neq \emptyset$, and is given the opposite orientation when  $\zed \cap \Xs_i \neq \emptyset$.

Let us define a new set of curves $\betas^\zed = \{\beta_j^{\zed}|  j=1, \dots, n\}$ on
the torus $\TT$. Let $j_i$ be the index corresponding to the vertical
circle $\beta_{j_i}$ just to the left of a marking $Z_i \in \zed$. We let
$\beta_{j_i}^{\zed}$ be a circle encircling $Z_i$ and intersecting
$\beta_{j_i}$, as well as the $\alpha$ curve just below $Z_i$, in two
points each; in other words, $\beta_{j_i}^{\zed}$ is obtained from
$\beta_j$ by handlesliding it over the vertical curve just to the
right of $Z_i$. For those $j$ that are not $j_i$ for any $Z_i \in \zed$, we 
let $\beta_j^{\zed}$ be a curve isotopic to $\beta_j$ and intersecting it in two points.

\begin {remark}
Our assumption on the existence of a free marking is crucial here, because it ensures that $\betas^\zed$ is a good set of attaching curves. Indeed, since $\zed$ has strictly fewer than $n$ markings, at least one beta curve in $\betas^\zed$ remains homologically nontrivial in $H_1(\TT)$.  
\end {remark}

We denote 
$$\mathbb{T}_{\beta}^{\zed} = \beta_1^\zed \times \dots \times \beta_n^\zed \subset \Sym^n(\TT).$$

Observe that
$$ (\TT, \alphas, \betas^\zed, \Os^{\orL(\zed)}, \Xs^{\orL(\zed)})$$
is a multi-pointed Heegaard diagram representing the link $\orL - L(\zed)$. For each $\s \in \bH(L)$, we can associate to this diagram a resolved Floer complex 
$$\Chain^-(\Ta, \Tb^\zed, \psi^{\orL(\zed)}(\s)),$$ as in Section~\ref{sec:gens}. For simplicity, we omit the basepoints from the notation for this chain complex. The map $ \psi^{\orL(\zed)} : \bH(L) \longrightarrow \bH(L-L(\zed))$ is as in Section~\ref{sec:reduction}.

More generally, let $\orM$ be any sublink of $L$ containing $\orL(\zed)$, such that the restriction to $L(\zed)$ of the orientation on $\orM$ coincides with $\orL(\zed)$. We can then consider a diagram
$$ (\TT, \alphas, \betas^\zed, \Os^{\orM}, \Xs^{\orM})$$
and resolved Floer complexes
$$\Chain^-(\Ta, \Tb^\zed, \psi^{\orM}(\s)).$$
By writing $\psi^{\orM}(\s)$ we implicitly mean that the basepoints are taken to be those in $\Os^{\orM}$ and $\Xs^{\orM}$.

\begin{figure}
\begin{center}
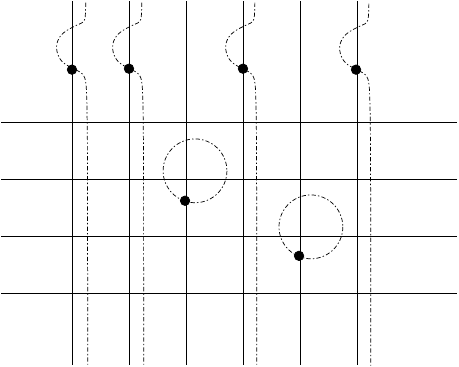
\end{center}
\caption {{\bf A new collection of curves.} We show here a part of a grid diagram, with the horizontal segments lying on  curves in $\alphas$ and the straight vertical segments lying on curves in $\betas$. The interrupted curves (including the two circles) represent curves in $\betas^\zed$, where $\zed$ consists of the two markings $Z_1$ and $Z_2$. The maximal degree intersection point $\Theta^\can_{\emptyset, \zed}$ is represented by the black dots.}
\label{fig:betaZ}
\end{figure}

When we have two collections of markings $\zed, \zed'$ such that $\zed \cup \zed'$ is consistent, we will require that $\beta_i^\zed$ and $\beta_i^{\zed'}$ intersect in exactly two points. Hence, there is always a unique maximal degree intersection point $\Theta^\can_{\zed, \zed'} \in \mathbb{T}_{\beta}^{\zed} \cap \mathbb{T}_{\beta}^{\zed'}$. See Figure~\ref{fig:betaZ}.

Of course, for the above Floer complexes to be well-defined, we need to make sure that the underlying Heegaard diagrams are admissible. This follows from the stronger admissibility result in Lemma~\ref{lem:admissiblegrid} below.

\subsection {Handleslides for a sublink}
\label {sec:desublink2}

Let $M \subseteq \orL$ be a sublink, endowed with an arbitrary orientation $\orM$. We seek to define a descent map for the sublink $\orM$, along the lines of Section~\ref{subsec:desublink}.

Set
$$\zed(\orM) = \bigcup_{i \in I_+(\orL, \orM)} \Os_i \cup \bigcup_{i \in I_-(\orL, \orM)} \Xs_i. $$ 

We proceed to construct a hyperbox $\Hyper^{\orL, \orM}_G$ for the pair $(\orL, \orM)$ (in the sense of Definition \ref{def:pair}), as follows. Order the components of $M$ according to their ordering as components of $L$:
$$ M = L_{i_1} \cup \dots \cup L_{i_{m}}, \ \ i_1 < \dots < i_{m}.$$
For $j=1, \dots, m$, let us denote $M_j = L_{i_j}$ for simplicity, and equip $M_j$ with the orientation $\orM_j$ induced from $\orM$. Then $\zed(\orM_j)$ is either $\Os_{i_j}$ or $\Xs_{i_j}$. In either case, we have an ordering of its elements, so we can write 
$$ \zed(\orM_j) = \{Z^{\orM_j}_1, \dots, Z^{\orM_j}_{d_{i_j}} \},$$
where $d_{i_j}$ is the cardinality of $\zed(\orM_j)$.

The hyperbox $\Hyper^{\orL, \orM}_G$ will be $m$-dimensional, of size 
$$\dd^{M} = (d_{i_1}, \dots, d_{i_m}).$$ It will be a pure $\beta$-hyperbox, i.e. the $\alpha$  curves remain fixed: they are the horizontal circles in the grid diagram $G$. For each multi-index $\eps = (\eps_1, \dots, \eps_m) \in \E(\dd^M)$, we let $\zed(\orM)^\eps \subseteq \zed(\orM)$ be the collection of markings 
$$ \zed(\orM)^{\eps} = \bigcup_{j=1}^m \{Z^{\orM_j}_1, \dots, Z^{\orM_j}_{\eps_j} \}.  $$

We then set 
$$\betas^{\eps} = \betas^{\zed(\orM)^\eps}$$
be the collection of beta curves handleslid at the points of $\zed(\orM)^\eps$. For each $\eps$, consider the Heegaard diagram $$\Hyper^{\orL, \orM}_{G, \eps} = (\TT, \alphas, \betas^\eps, \Os^{\orM}, \Xs^{\orM}).$$ This diagram represents the link $\orL - M$. 

According to Definitions~\ref{def:alphabeta} and \ref{def:pair}, to have a well-defined hyperbox for a pair, we need admissibility for the Heegaard multi-diagrams corresponding to each hypercube in the hyperbox. In fact, it is convenient to consider the bigger multi-diagram that includes the curve collections $\betas^\zed$ for all subsets $\zed \subseteq \zed(\orM)$. Once we prove that that diagram is admissible, it will follow that all its sub-diagrams (such as the ones relevant for our hypercubes) are admissible as well.

\begin{lemma}
\label{lem:admissiblegrid}
Consider the Heegaard multi-diagram consisting of the curve collections $\alphas$ and $\betas^{\zed}$ on $\TT$, for all possible $\zed \subseteq \zed(\orM)$, and with basepoints given by $\Os^{\orM}, \Xs^{\orM}$. Then, this multi-diagram is admissible, in the sense of Definition~\ref{def:mpd}. 
\end{lemma}

\begin{proof}
The markings in $\Xs^{\orM}$ are irrelevant for admissibility, so we can ignore them. There is one marking in $\Os^{\orM}$ in each row and in each column of the grid $G$. Without loss of generality, let us assume that $\orM$ is oriented as in $\orL$, so that $\Os^{\orM} = \Os=\{O_1, \dots, O_n\}$. 

Recall from Definition~\ref{def:mpd} that a multi-periodic domain is a $\zz$-linear combinations of regions, such that its boundary is a $\zz$-linear combination of the curves in the multi-diagram, and such that its multiplicity at any $O_i$ is zero. Admissibility means there are no non-trivial multi-periodic domains that have only nonnegative local multiplicities. By way of contradiction, suppose we had such a domain $D$.

Let us mark an asterisk in each square of the original grid diagram $G$. When we construct the new beta curves $\beta^{\zed}_j$, we make sure that each beta curve encircling some $O_i$ does not go around an asterisk, and also that whenever we isotope a beta curve to obtain a new beta curve intersecting it at two points, these isotopies do not cross the asterisks. See Figure~\ref{fig:asterisks}.

\begin{figure}
\begin{center}
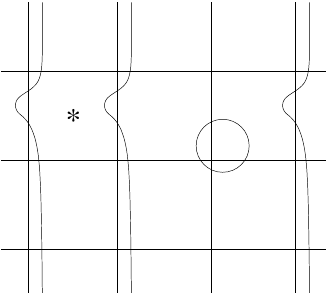
\end{center}
\caption {{\bf Asterisks in the Heegaard multi-diagram.}}
\label{fig:asterisks}
\end{figure}

Let $a_j \in \zz$ be the coefficient with which the curve $\alpha_j$ appears in $\del D$. This means that, whenever we have two regions in the multi-diagram separated by a segment in $\alpha_j$, the difference between the multiplicities of $D$ on these two regions equals $a_j$. Observe that there are two such regions, such that one on top of $\alpha_j$ contains an $O$ marking. The multiplicity of $D$ is supposed to be zero at that marking and nonnegative at the region below. It follows that $a_j \leq 0$, and therefore the multiplicity increases (or stays constant) whenever we cross an alpha curve downward. In particular, if we look at the multiplicities of $D$ at the regions containing asterisks in a given column, they should be nondecreasing in cyclic order, which means that they are constant.

We can do a similar argument with the beta curves, by letting $b_j$ be the sum of the coefficients in $\del D$ of $\beta_j$ and all of its Hamiltonian translates. In this case $b_j$ is the difference in multiplicities as we cross all the translates of $\beta_j$. By looking at the region containing the $O$ marking in each column, we find that $b_j$ are all of the same sign. This shows that the multiplicities of $D$ at the regions containing asterisks in a given row are constant. 

Combining these two results, we deduce that $D$ has the same multiplicity at all the regions containing asterisks. However, we know there is at least one free $O$ marking in $G$, so that no handleslides are done over that marking. Therefore, the region containing that marking has an asterisk, and $D$ has multiplicity zero there. It follows that $D$ has multiplicity zero at all asterisks.

This means that $D$ is a sum of disjoint periodic domains, each supported in the area bounded by the Hamiltonian translates of the same curve: either a vertical beta curve, or a curve encircling an $O$ marking; cf. Figure~\ref{fig:shaded}. We can study these domains separately.

Let us consider the neighborhood of a vertical beta curve, as on the left hand side of Figure~\ref{fig:shaded}. (The other case is similar.) Let $c_1, \dots, c_p$ be the coefficients in $\del D$ of each of the beta curves in the picture. As we cross a beta curve (say, in the directions shown by arrows in Figure~\ref{fig:shaded}), the multiplicity of $D$ should changed by the corresponding $c_i$. Note that the outside region at the top of the picture has an asterisk and hence multiplicity zero. Therefore, we have $c_i \geq 0$ for all $i$. Let $m_0 \leq m_1 \leq \dots \leq m_{p}$ be the multiplicities at the regions shown in the figure, as we go from left to right. We must have $m_0 = m_p=0$, because those regions contain asterisks. Also, $m_i = m_{i-1} + c_i$, which implies that 
$$ 0 =m_0 \leq m_1 \leq \dots \leq m_{p-1} \leq m_p=0.$$

\begin{figure}
\begin{center}
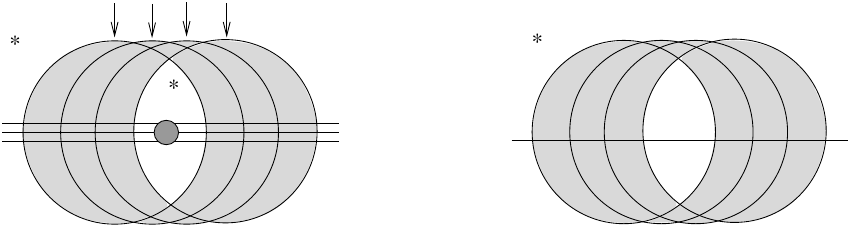
\end{center}
\caption {{\bf The Hamiltonian translates of a beta curve.} On the left we have the neigborhood of a vertical beta curve on the grid, shown here as a circle. The other circles are Hamiltonian translates of the same curve. The horizontal lines are alpha curves, and the disk in the middle is a handle that is connected to the exterior of the picture to produce the torus $\TT$. The shaded area is the potential support of a periodic domain, in the proof of Lemma~\ref{lem:admissiblegrid}. On the right we show a similar neighborhood of a beta curve encircling an $O$ marking.}
\label{fig:shaded}
\end{figure}

This shows that $m_i=0$ and $c_i=0$ for all $i$, which implies that the domain is trivial. We arrived at a contradiction, and therefore we have admissibility. 
\end{proof}

To finish the construction of the hyperbox $\Hyper^{\orL, \orM}_G$, whenever $\eps, \eps' \in \E(\dd^{M})$ with $\eps < \eps'$ are neighbors, we need to pick chain elements $ \Theta_{\eps, \eps'} \in \Am(\T_{\beta^\eps}, \T_{\beta^{\eps'}}, \zero)$. We choose
$$ 
 \Theta_{\eps, \eps'} = \begin {cases} 
\Theta^\can_{\zed(\orM)^{\eps}, \zed(\orM)^{\eps'}} & \text{if } \|\eps' - \eps\| =1\\
0 & \text{otherwise.}
\end{cases}
$$
The proof of Equation \eqref{eq:compatibility} for these chains is similar to those of \cite[Lemma 9.7]{HolDisk} and \cite[Lemma 4.3]{BrDCov}. 

As in Section~\ref{sec:gens}, for each $\s \in \bH(L)$, we have an associated hyperbox of generalized Floer complexes 
$$ \Cc(\Hyper^{\orL, \orM}_G, \psi^{\orM}(\s)),$$
of size $\dd^{'M}:=(d'_{i_1}, \dots, d'_{i_m})$, where
$$ d'_{i_j} = \begin{cases}
d_{i_j} &\text{if } i_j\in I_+(\orL, \orM),\\
2d_{i_j} - 1 &\text{if } i_j\in I_-(\orL, \orM)
\end{cases}.$$
(Note that in our setting, the side length $d_{i_j}$ of the hyperbox $\Hyper_G^{\orL, \orM}$ is the same as the number $p_{i_j}$ of $w$ basepoints on the respective link component.)

At each vertex $\eps \in \E(\dd^{'M})$ we have the Floer complex 
$$C^\eps = \Cc_{\eps^<}(\Hyper^{\orL, \orM-M_{\eps}}_{G, \eps^>}, M_{\eps}, \psi^{\orM}(\s)),$$
and along the faces we have linear maps $\De^{\eps' -\eps}_{\eps}: C^{\eps} \to C^{\eps'}$. 

We compress the hyperbox $\Cc(\Hyper^{\orL, \orM}_{G}, \psi^{\orM}(\s))$ as in  Section~\ref{sec:compression2}, and define 
\begin {equation}
\label {eq:dezed}
\hat \De^{\orM}_{p^{\orM}(\s)} : \Chain^-(G, \orM, p^{\orM}(\s))  \cong \Chain^-(\Ta, \Tb, \orM, \psi^{\orM}(\s)) \to \Chain^-(\Ta, \Tb^{\zed(\orM)}, \psi^{\orM}(\s))
\end {equation}
to be the longest diagonal map in the resulting hypercube.

For example, when $M = L_i$ is a single component with orientation $\orM = \orL$, the map $\hat \De^{\orM}_{p^{\orM}(\s)}$ is a composition of the triangle maps corresponding to handleslides over the basepoints in $\zed(\orM)$, in the given order. When $M$ has several components, it is a sum of compositions of more complicated polygon maps, corresponding to chain homotopies (of higher order) between compositions of the handleslide maps.

\begin{remark}
The formulas used to define $\De^{\eps' -\eps}_{\eps}$ in Section~\ref{sec:gens} simplify somewhat in the setting of grids. Indeed, all the hypercubes of strongly Heegaard diagrams are of pure $\beta$ type, meaning that the curve collection $\alphas$ stays constant; that is, the bipartition maps take only the value $\beta$. Further, the chain elements $\Theta^{\beta}_{\zeta, \zeta'}$ consist of a single intersection point (in Alexander gradings zero) when $\|\zeta' -\zeta \|=1$, and are zero otherwise. Thus, we have
$$ c(\a^{\beta}_{\gamma' \gamma})=\begin{cases} 
1 &\text{if } \a = \Theta^{\beta}_{\zeta, \zeta'},\\
0 &\text{otherwise}.
\end{cases}$$
and
$$  \bu_{i,j}(\a^{\beta}_{\gamma' \gamma}) := u_{i,j}(\a^{\beta}_{\zeta \zeta'}) = 0.$$
From here, we find that, for $\a = \Theta^{\beta}_{\zeta, \zeta'}$ (which is the only intersection point with a nonzero contribution), the quantities $n^{\a}$ that appear in \eqref{eq:multiU} are the same as the corresponding basepoint multiplicities $n$ (with the same subscripts).
\end{remark}

\subsection {The handleslid complex}
\label {sec:handleslid}

For $O_j \in \Os_i$, we let $U_{j'}$ be the variable corresponding to
the row exactly under the row through $O_j$. Given a sublink 
$M \subseteq L$, we define a chain complex 
\begin {equation}
\K (M) := \bigotimes_{M_i \subseteq M} \bigotimes_{j \in \Os_i} \Bigl( \Ring \xrightarrow{U_j - U_{j'} } \Ring \Bigr).
\end {equation}

Given $\s \in \bH(L)$ and an orientation $\orM$ of $M$, we seek to describe the complex $\Chain^-(\Ta, \Tb^{\zed(\orM)}, \psi^{\orM}(\s))$ explicitly.  Before stating the result, let us introduce some terminology. We call the Heegaard diagram $$\Hyper^{\orL, \orM}_{G, \dd^{M}} = (\TT, \alphas, \betas^{\zed(\orM)}, \Os^{\orM}, \Xs^{\orM})$$ the {\em handleslid diagram}. It contains two kinds of beta curves: those that approximate curves in the original collection $\betas$, and those that encircle one of the points in $\zed(\orM)$. We call the former curves {\em old}, and the latter {\em recent}. Each recent curve bounds a disk in $\TT$, which we call a {\em marked disk}. The curves in $\alphas$ are also of two types: those that intersect a recent $\beta^{\zed(\orM)}$ curve, and those that do not. We call the first alpha curves {\em special}, and the latter {\em normal}. The connected components of $\TT \setminus \alphas$ are called rows, and a row is called {\em special} if at least one of its boundaries is  special. If both of its alpha boundaries are special, the row is called {\em very special}.

Observe that each recent beta curve intersects a unique (special) alpha curve, and  it does so in two points. Hence, each $\x \in   \Ta \cap \Tb^{\zed(\orM)}$ must contain exactly one of these two points.

Note that there is a {\em quasi-destabilized grid diagram} $G^{L - M}$ obtained from $G$ by
eliminating all rows and columns on which $M$ is supported; compare Section~\ref{sec:triangles}. The diagram $G^{L-M}$ represents the link $L - M$, with the orientation induced from $\orL$.

\begin {proposition}
\label {lemma:pear}
For any $\s \in \bH(L)$, for a suitable choice of almost complex structure on $\Sym^n(\TT)$, there is an isomorphism of chain complexes of $\Ring$-modules:
\begin {equation}
\label {eq:psim}
\Psi^{\orM}_{\psi^{\orM}(\s)} : \Chain^-(\Ta, \Tb^{\zed(\orM)}, \psi^{\orM}(\s)) \to \Chain^-(G^{L-M}, \psi^{\orM}(\s))[[\{U_i \}_{L_i \subseteq M }]] \otimes_{\Ring} \K(M).
\end {equation}
Here, the square brackets mean adjoining the variables $U_i$, in order to make the respective group into an $\Ring$-module. 
\end {proposition}

Observe that on the right hand side of \eqref{eq:psim} we have the complex associated to the quasi-destabilized grid diagram $G^{L-M}$, for which holomorphic disks in the symmetric product exactly correspond to empty rectangles on the grid, see \cite{MOS}. In proving Proposition~\ref{lemma:pear}, the challenge will be to show that the complex on the left hand side of \eqref{eq:psim} has the same description. The proof will occupy the rest of this subsection.

The curves $\alphas$ and $\betas^{\zed(\orM)}$ split the torus $\TT$ into a number of regions. These can be rectangles, bigons, hexagons, and octagons. The bigons and hexagons appear near the markings in $\zed( \orM)$, while the octagons only appear when two markings in $\zed(M)$ lie in adjacent rows and adjacent columns of $G$. The hexagons and octagons are always supported in special rows.

Let $\phi \in \pi_2(\x, \y)$ be a homology class of Whitney disks in $\Sym^n(\TT)$ with boundaries on $\Ta$ and $\Tb^{\zed(\orM)}$, for $\x , \y \in  \Ta \cap \Tb^{\zed(\orM)}$. The class $\phi$ has an associated domain, see Section~\ref{sec:domains}. Recall that the Maslov index $\mu(\phi)$ can be calculated in terms of the domain $\D=\D(\phi)$ using Lipshitz's formula \eqref{eq:lipshitz1}:
\begin {equation}
\label {eq:lipshitz}
 \mu(\phi) = \sum_{x \in \x} n_x(\D) + \sum_{y \in \y} n_y(\D) + e(\D).
\end {equation} 

Suppose $\D$ is written as a linear combination of regions
 $$\D = \sum a_i R_i, \ \ a_i \in \zz.$$

For any region $R$, we set $n(R, \x, \y)= n_\x(R) + n_\y(R)$ and $\mu(R, \x, \y) = n(R, \x, \y) + e(R)$. Since the Euler measure and vertex multiplicities are additive, Equation~\eqref{eq:lipshitz} gives
\begin {equation}
\label {eq:L2}
 \mu(\phi) = \sum a_i \mu(R_i, \x, \y).
 \end {equation}

The differential on $\Chain^-(\Ta, \Tb^{\zed(\orM)}, \psi^{\orM}(\s))$ involves counts of holomorphic disks in homology classes $\phi$ with $\mu(\phi) = 1$. The following are two necessary conditions for a domain $\D=\D(\phi)$ of index one to admit a holomorphic representative:
\begin {enumerate}
\item All local multiplicities $a_i$ of $\D$ must be nonnegative;
\item The support of the domain $\D$ (that is, the closure of the union of the regions $R_i$ such that $a_i \neq 0$) must be connected. 
\end {enumerate}

The first condition is a consequence of the principle of positivity of intersection for holomorphic objects, see \cite[Lemma 3.2]{HolDisk}. The second condition needs to be satisfied because if the support were disconnected, the class $\phi$ would be the sum of two domains $\phi_1$ and $\phi_2$, with each $\phi_i$ required to admit pseudo-holomorphic representatives generically. Hence $\mu(\phi_1), \mu(\phi_2) \geq 1$, which would contradict the fact  that $ \mu(\phi_1) + \mu(\phi_2) = \mu(\phi) = 1$.
 
A domain $\D$ that satisfies the two conditions above is called {\em positive} and {\em connected}.  
  
For our handleslid Heegaard diagram, we can completely characterize the index one domains 
that are positive and connected. Indeed, given a class $\phi \in \pi_2(\x, \y)$ that has such a domain, let us understand the quantities $\mu(R_i, \x, \y)$ that appear with nonzero multiplicity $a_i > 0$ in Equation~\eqref{eq:L2}:

\begin {itemize}
\item  If $R_i$ is a rectangle, then $e(R_i) = 0$ and $\mu(R_i, \x, \y) \in [0,1]$, with the exact value depending on the number of vertices in $\x, \y$ among the corners of $R_i$. 

\item If $R_i$ is a bigon, then $e(R_i) = 1/2$ and $\mu(R_i, \x, \y) = 1$. Indeed, note that $\x$ contains exactly one of the two points of intersection between  the respective new beta curve and a special alpha curve. Hence $\x$ contributes $1/4$ to the quantity $n(R_i, \x, \y)$. The same goes for $\y$, while the Euler measure of the bigon is $1/2$. 

\item If $R_i$ is a hexagon, then $e(R_i) = -1/2$ and $\mu(R_i, \x, \y) \in [0, 1/2]$.  Indeed, there is a contribution of $1/2$ to the vertex multiplicity coming from points on the recent beta curve, and there may be an additional contribution from intersections between old beta curves and normal alpha curves. We denote this additional contribution by $n'(R_i, \x, \y)$. It is at most $1/2$.

\item  If $R_i$ is an octagon, then $e(R_i) = -1$ and $\mu(R_i, \x, \y) = 0$, because the vertex multiplicity contribution is always one, coming from points on the two recent beta curves on the boundary of the octagon.
\end {itemize} 
 
Hence:
\begin {equation}
\label {eq:L3}
 \mu(\phi) = \sum_{ \text{bigons} R_i} a_i + \sum_{ \text{rectangles} R_i} a_i n(R_i, \x, \y)  + \sum_{ \text{hexagons} R_i}  a_i n'(R_i, \x, \y)  \geq 0.
 \end {equation} 
   
If $\mu(\phi) = 1$, it follows that we can have at most one bigon in the support of $\D$, and that the bigon (if it exists) must appear with multiplicity one. Further, no points of $\x \cap \y$ can be contained in the interior of $D$. 

\begin{figure}
\begin{center}
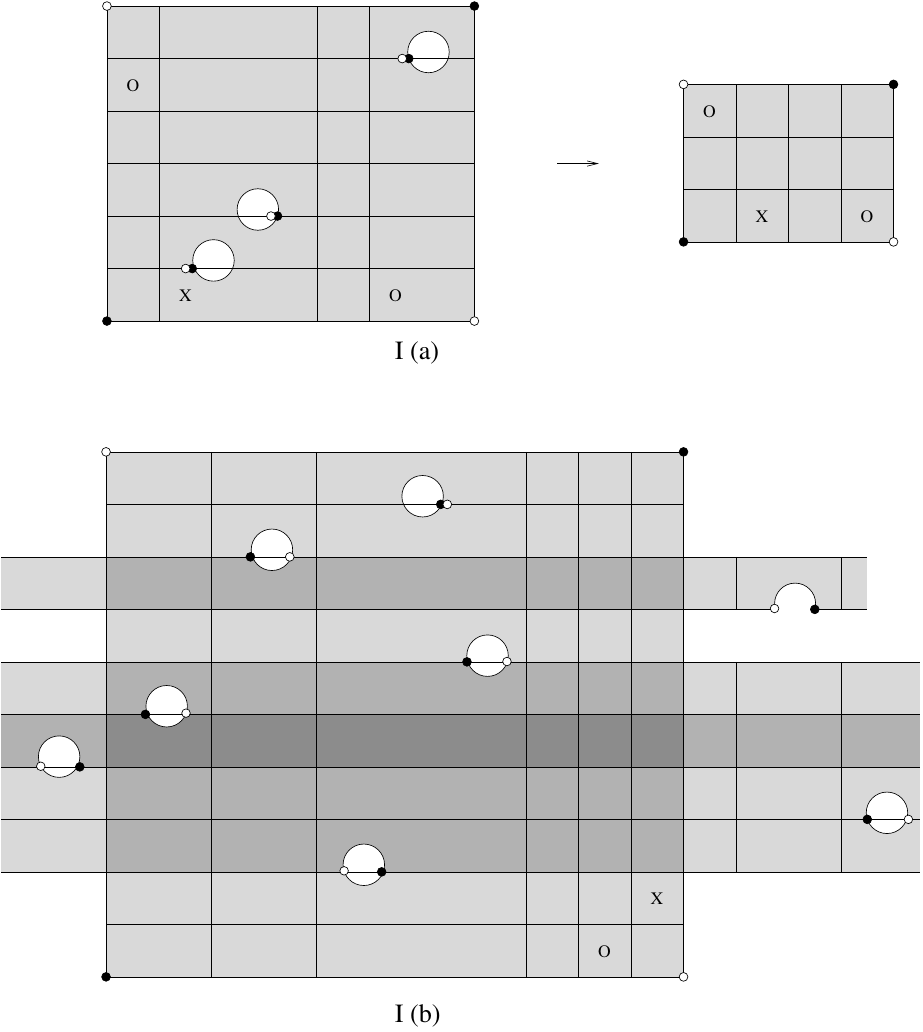
\end{center}
\caption {{\bf Positive connected domains of index one: Case I.} In each subcase, the domain is shaded, with darker shading corresponding to higher local multiplicity. In subcase (a) we only have multiplicities zero and one, while in subcase (b) we can arbitrarily large multiplicities; in the figure, they go up to multiplicity three. The black dots represent components of the initial generator $\x$, and the white dots components of $\y$. There can be various markings ($X$ and $O$) in the domains, though not in the rows containing markings in $\zed(\orM)$. In subcase (a), apart from the domain in the handleslid diagram $(\TT, \alphas, \betas^{\zed(M)})$ on the left hand side, we also show the corresponding domain in the quasi-destabilized  grid diagram $G^{L-M}$, on the right hand side. }
\label{fig:CaseI}
\end{figure}

\begin{figure}
\begin{center}
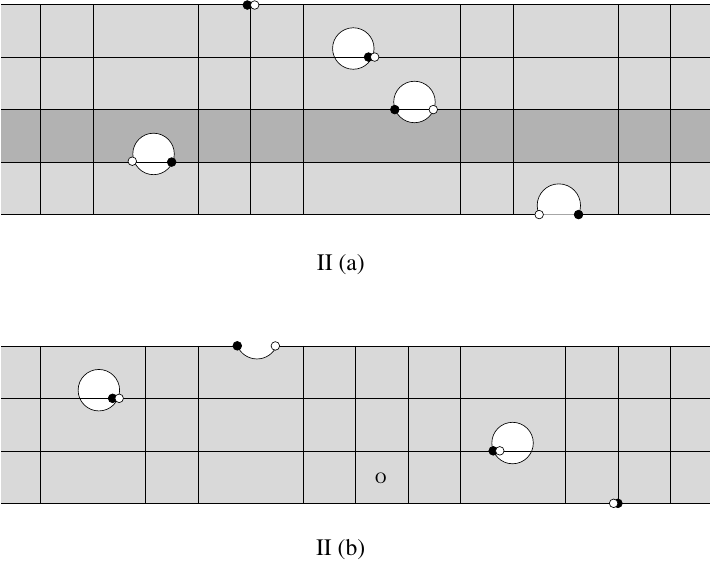
\end{center}
\caption {{\bf Positive connected domains of index one: Case II.} We use the same conventions as in Case I. In subcase (a) we show a domain of width four and having local multiplicity two in the darkly shaded regions. In (b) we show a domain of width three with only zero and one local multiplicities. Of course, in each subcase the widths and multiplicities can be arbitrarily large, just as we saw in Case I (b).}
\label{fig:CaseII}
\end{figure}

\begin{figure}
\begin{center}
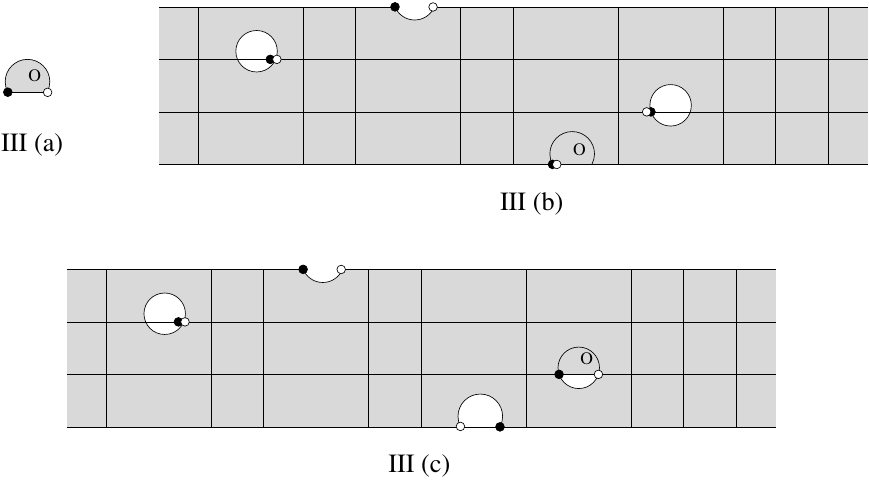
\end{center}
\caption {{\bf Positive connected domains of index one: Case III.} Same conventions as before.  In subcases (b) and (c) the widths and multiplicities can be arbitrarily large, just as in Case II.}
\label{fig:CaseIII}
\end{figure}

\begin{figure}
\begin{center}
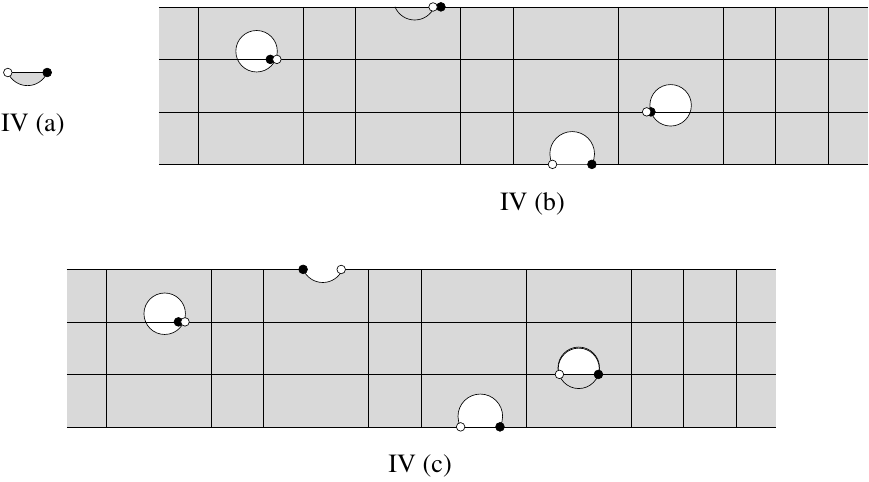
\end{center}
\caption {{\bf Positive connected domains of index one: Case IV.} Same conventions as before.  In subcases (b) and (c) the widths and multiplicities can be arbitrarily large, just as in Cases II and III.}
\label{fig:CaseIV}
\end{figure}

We distinguish several cases for $\phi \in \pi_2(\x, \y)$ positive, connected and of index one. These cases are shown in Figures \ref{fig:CaseI}-\ref{fig:CaseIV}.
\medskip

I.  There exists a normal alpha curve $\alpha_j$ such that $\alpha_j \cap \x \neq \alpha_j \cap \y$. Then the two points $\alpha_j \cap \x$ and  $\alpha_j \cap \y$ are on different (old) beta curves, so there must exist another normal alpha curve $\alpha_k \neq \alpha_j$ with   $\alpha_k \cap \x \neq \alpha_k \cap \y$.   We get a contribution of at least $1/4$ to $n(\D)$ from each of the four intersection points on $\alpha_j$ and $\alpha_k$. Since $\mu(\phi)=1$, there can be no other contribution. It follows that the support of $\D$ has no bigons.  One possibility is that $\D$ is a big rectangle with vertices the four intersection points $\alpha_j \cap \x, \alpha_j \cap \y, \alpha_k \cap \x, \alpha_k \cap \y$, with all the marked disks in the interior of the rectangle removed. This is Case I (a)  shown in Figure~\ref{fig:CaseI}. Note that there is a one-to-one correspondence between such rectangles and empty rectangles in the quasi-destabilized  grid diagram $G^{L-M}$.
  
 The other possibility is that in addition to a big rectangle, $\D$ also contains some additional rectangle, hexagon, or octagon regions. These have to be supported in some very special rows. Further, all non-bigon regions in those rows are contained in the domain $\D$, because the multiplicity of regions in a very special row can only change when we pass the vertical curves through  $\alpha_j \cap \x$ and  $\alpha_j \cap \y$. Thus, the domain consists of a big rectangle as in Case I (a), plus one or more very special rows intersecting the big rectangle, minus the bigons in those rows. The multiplicities in the special rows cannot change by more than one as we pass a special alpha curve, because the support of $\D$ contains no bigons. A domain of this type is said to be in Case I (b). An example is shown in Figure~\ref{fig:CaseI}. 
 
 \medskip
II.  We have $\alpha_j \cap \x = \alpha_j \cap \y$ for all normal alpha curves $\alpha_j$, and the support of $\D$ contains no bigons. Then $\D$ is  supported in some adjacent special rows. Further, there are no old beta curves such that when we pass them the multiplicity of $\D$ gets changed. Hence, if the support of $\D$ contains a non-bigon region in a special row, it must contain all the non-bigon regions in that row, with the same multiplicity.  Further, one of the boundaries of the support must be a normal alpha curve, which produces all nonzero contributions to $\mu(\phi)$ in Equation \eqref{eq:L3}. There are two subcases, II (a) and II (b), according to whether the normal alpha curve is the top or the bottom boundary of the support. In both cases, the multiplicity in the row bounded by the normal alpha curve must be one. The multiplicities in the other rows can be higher, but they are constant in each row and cannot change by more than one as we pass a special alpha curve. See Figure~\ref{fig:CaseII}.

\medskip
III. We have $\alpha_j \cap \x = \alpha_j \cap \y$ for all normal alpha curves $\alpha_j$, and the support of $\D$ includes a bigon which contains a marking in $\zed(\orM)$. Then the domain $\D$ is  supported in some adjacent special rows, and the bigon gives the only nonzero contributions to $\mu(\phi)$ in Equation \eqref{eq:L3}. Further, if $\D$ contains a rectangle, hexagon or octagon in a special row, it must contains all the rectangles, hexagons, and octagons in that row, with the same multiplicity. There are three subcases III (a) - (c), according to whether: (a) $\D$ is the bigon itself; (b) $\D$ contains some special rows, but only above the bigon; (c) $\D$ contains special rows both below and above the bigon. See Figure~\ref{fig:CaseIII}.
   
  \medskip 
IV.  We have $\alpha_j \cap \x = \alpha_j \cap \y$ for all normal alpha curves $\alpha_j$, and the support of $\D$ includes a bigon which does not contain a marking in $\zed(\orM)$.  This is similar to Case III, in that we have the same observations about the support consisting of some adjacent special rows, and three subcases: (a) $\D$ is the bigon itself; (b) $\D$ contains some special rows, but only  below the bigon; (c) $\D$ contains special rows both below and above the bigon. See Figure~\ref{fig:CaseIV}.
   
 \medskip 
 
When a domain $\D$ is in one of the cases above, we say that it is of the corresponding {\em type}. For example, a domain in Case II (b) is called of type II (b). 
 
 \begin {definition}
If $\D$ is a domain on the handleslid grid, the number of rows that have nontrivial intersection with its support is  called the {\em width} of $\D$.    
\end {definition}   
   
We would like to count the number of holomorphic representatives for each class $\phi$ with a positive, connected domain $\D$ of index one. One should note that many of the domains in Cases I-IV above are decomposable, meaning that we can write them as a sum $\D_1 + \D_0$, with $\D_i$ of index $i \in \{0,1\}$ being positive and connected. This means that the number of holomorphic representatives for $\D$ may depend on the choice of almost complex structure on the symmetric product. We will explain how to make a suitable choice soon.

When counting holomorphic representatives, it is convenient to use Lipshitz's cylindrical formulation of Heegaard Floer homology \cite{LipshitzCyl}, as in Section~\ref{sec:lipshitz}. In his setting, we need to choose an almost complex structure on $W = \TT \times [0,1] \times \rr$, rather than on the symmetric product.

Suppose 
$$\zed(\orM) = \{Z_1, \dots, Z_k\}.$$

The labeling of the markings is as follows. Choose a normal alpha curve $\alpha_1$ on the handleslid diagram. Go vertically down around the torus, starting at $\alpha_1$ and ending at $\alpha_1$, and number the markings in $\zed(\orM)$ as $Z_1, \dots, Z_k$, in the order in which they are encountered.

\begin{figure}
\begin{center}
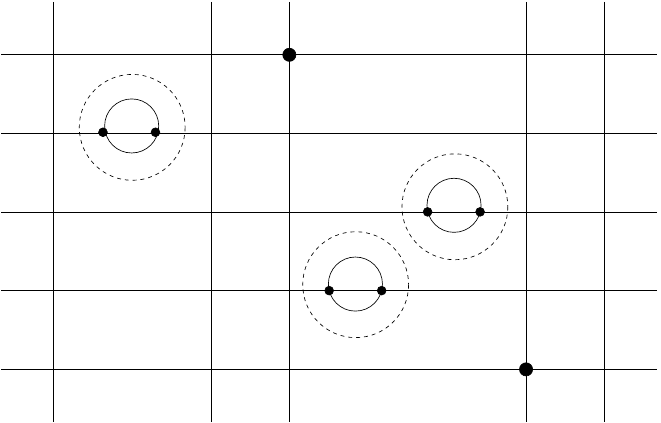
\end{center}
\caption {{\bf Degenerations.}  We choose the complex structure by stretching the neck along the dashed curves. The bigger black dots are components of the generator $\bar \x$ in Lemma~\ref{lemma:barx}.}
\label{fig:3rows}
\end{figure}

For each $Z_i \in \zed(M)$, there is a recent beta curve $\beta_{j_i}^{\zed(\orM)}$, which is the boundary of the marked disk containing  $Z_i$. Choose a slightly bigger curve around this beta curve, and introduce a long cylindrical neck of length $T_i$ there, for $T_i \gg 0$. We can then choose an almost complex structure $J(T_1, \dots, T_k)$ on $W$ which only depends on $T_i$ on the neck around $Z_i$, and is split on these necks. This is the analogue of the ``stretching the neck'' process from Section~\ref{sec:gluing}, only now we consider flow lines rather than triangles. In fact, we can view the handleslid diagram as a special connected sum of the quasi-destabilized grid diagram $G^{L-M}$ and $k$ genus zero diagrams $\Sphere_1, \dots, \Sphere_k$. Each $\Sphere_i$ is obtained from the diagram in Figure~\ref{fig:ess} by deleting the gamma curve. See Figure~\ref{fig:3rows}. 

Given a marking $Z_i$, we denote by $M_i$ the marked disk containing it, and by $x_i$ and $y_i$ the two points of intersection of $\beta_{j_i}^{\zed(\orM)}$ with the corresponding alpha curve. The point $x_i$ is on the left of $y_i$ inside the marked disk.

\begin {lemma}
\label {lemma:degz}
 Let $\phi \in \pi_2(\x, \y)$ be a homology class of Whitney disks in the handleslid diagram. Suppose $\mu(\phi)=1$ and the domain $\D = \D(\phi)$ is positive and connected. Suppose $\alpha_j$ is a special alpha curve that intersects a recent beta curve $\beta_j$. Let $M_i$ be the marked disk whose boundary is $\beta_j$, and let $Z_i$ be the marking inside. Suppose $\alpha_j \cap \x = x_i, \ \alpha_j \cap \y = y_i$, and that $\D$ has multiplicity zero everywhere inside the marked disk $M_i$.  Denote by $r$  the row whose upper boundary is $\alpha_j$, and denote by $b$ the bigon 
 $r \cap M_i$.  If the row $r$ contains a marking $Z_{i+1} \in \zed(M)$, fix the corresponding neck-length $T_{i+1}$. If the class $\phi$ admits holomorphic representatives for a sequence of almost complex structures $J(T^{(n)}_1, \dots, T_k^{(n)})$ such that  $T_i^{(n)} \to \infty $ (where $T_i^{(n)}$ corresponds to $Z_i$), and $T_{i+1}^{(n)} = T_{i+1}$ is fixed, then the support of the domain $\D$ is exactly the union of all the regions in the row $r$, except the bigon $b$.
 \end {lemma}
   
\begin {proof}   
For simplicity, let us first assume that all neck-lengths $T_s^{(n)} = T_s$ are fixed for $s \neq i$, and $T_i^{(n)} \to \infty$.

Let $m_1$ and $m_2$ the multiplicities of $\D$ on each side of $\alpha_j$, near the connected sum curve. By an analogue of Proposition~\ref{prop:converges}, in the limit $T_i \to \infty$ the presumed holomorphic representatives of $\phi$ degenerate into a broken flow line $v$ on $\Sphere_i$, and the union of a broken flow line $v'$ on the other side (i.e. on the diagram $G'$ obtained by deleting $\beta_j$) with an annoying $\alpha$-degeneration with domain $\Per$. Let $\psi, \psi'$ be the homology classes of $v$ and $v'$, respectively. Lemmas~\ref{lemma:mumu} and \ref{lemma:tris} can be applied to homology classes of flow lines just as well as to triangles. Therefore,
$$ \mu(\phi) = \mu(\psi) + \mu(\psi') + \mu(\Per) - m_1 - m_2$$   
and
$$ \mu(\psi) = m_1 + m_2.$$

In principle, the latter equality should also involve the multiplicities of $\D$ inside the two bigons that form $M_i$, but those multiplicities are zero by hypothesis. We get:
$$ 1 = \mu(\phi) = \mu(\psi') + \mu(\Per).$$

Further, we must have $m_2 = m_1 + 1$ because of the acceptability condition \eqref{eq:corners}. This means that $\del \Per$ contains $\alpha_j$ with multiplicity $-1$. In particular, $\Per \neq 0$, so $\mu(\Per) \geq 1$, see Lemma~\ref{lemma:peru} for $\mu(\Per)$.  We deduce that $\mu(\psi') = 0$. Since $\psi'$ admits holomorphic representatives (for an almost complex structure which can be chosen to be generic away from the neck), it must be the case that $\psi'$ is trivial. As for the remaining domain $\Per$, it is positive, of index one, and satisfies $m_2 = m_1 + 1$. The only possibility is that it is the row $r$, viewed in the quasi-destabilized diagram $G'$, compare Equation~\eqref{eq:muannoy}. The claim follows, provided that the neck-lengths $T_i$ stay fixed for $s \neq i$.

Now allow the other neck-lengths to vary as well, except $T_{i+1}$. In the limit some of them go to infinity, so the holomorphic representatives of $\phi$ degenerate into broken flow lines on a number of surfaces $\Sphere_s$, plus the union of a broken flow line and $\alpha$-boundary degenerations on a diagram $G'$ obtained by deleting several recent  beta curves. Nevertheless, most of the arguments above still apply, with some notational modifications. The broken flow line on $G'$ must be trivial, so all we are left with is some $\alpha$-boundary degenerations on $G'$. The main difference is that now, instead of Equation~\eqref{eq:muannoy}, we have:
$$ \mu(\Per) = \sum_\alpha m(\Per, \alpha),$$ 
where the sum is over all curves $\alpha$ that did not intersect any of the deleted beta curves, and the quantity $m(\Per, \alpha)$ is the sum of the multiplicities of $\Per$ in the two rows that have $\alpha$ as part of their boundary.

It follows that the inequality $\mu(\Per) \geq 0$ in Lemma~\ref{lemma:peru} still holds true, but equality can happen for some nonzero domains on $G'$, namely those composed of very special rows such that the normal beta curves intersecting their boundaries have been deleted. As such, a priori there are more possibilities for the positive domain $\Per$ of index one that represents the $\alpha$-boundary degenerations. However, the fact that the recent beta curve just below $\beta_j$ (if it exists) has not been deleted, together with the fact that $\Per$ is connected and satisfies $m_2 = m_1 +1$, suffice to determine $\Per$ uniquely as being the row below $\alpha_j$.
\end {proof}   

Lemma~\ref{lemma:degz} suggests the following way of choosing the neck-lengths $T_1, \dots, T_k$ in the construction of the almost complex structure $J(T_1, \dots, T_k)$. We first choose $T_k \gg 0$ such that the conclusion of the Lemma holds true for all possible classes $\phi$ satisfying the hypothesis, with $i=k$. This is possible because there are a finite number of such classes, and the row $r$ below $Z_k$ does not contain any marking in $\zed(\orM)$. Next, choose $T_{k-1} \gg T_k$ such that the conclusion of the Lemma holds for all classes with $i=k-1$. Iterate this procedure until we get to $T_1$, so that   
$$T_1 \gg T_2 \gg \dots \gg T_k \gg 0.$$   
   
Lemma~\ref{lemma:degz} implies that, for this choice of almost complex structure, the domains of type I (b), III (c) and IV (c) have holomorphic disk counts equal to zero (mod $2$), so they do not contribute to the differential. The same goes for domains of type II (b) or III (b) that have width larger than one.

On the other hand, the holomorphic disk counts are one (mod $2$) for: 
\begin {itemize}
\item simple bigons, that is, domains of types III (a) and  IV (a);
\item domains of type I (a), see \cite[Lemma 3.11]{CubeResolutions};
\item annular domains of type  II (b) and width one, see \cite[proof of Lemma 3.4]{HolDiskTwo};
\item annular domains of type  III (b) and width one, see \cite[Lemma 9.4]{HolDisk}. Indeed, for these the conformal angle of the beta part of the upper boundary is very small, because of our choice of neck-lengths. This means that there exists a choice of cut length on the lower boundary that makes the annulus holomorphic. 
 \end {itemize}  
 
This leaves us with the domains of types II (a) and IV (b). It is more difficult to establish the holomorphic disk counts for these. Instead, we will settle for computing their net effect on the differential $\del$ on the handleslid complex.
   
More precisely, let us split the markings in $\zed(\orM)$ into equivalence classes, with the equivalence relation being generated by the relations $Z \sim Z'$ if $Z$ and $Z'$ lie in adjacent rows. Then, an equivalence class consists of all the markings in $\zed(\orM)$ that lie between two normal alpha curves. 

Let us consider such an equivalence class. Without loss of generality, say it is composed of the markings $Z_1=O_1, \dots, Z_p=O_p$, for some $p \leq k$. (See Figure~\ref{fig:3rows} for the case $p=3$.) Let $\bar G$ be the grid diagram obtained by deleting the recent beta curves $\beta_1, \dots, \beta_p$ that encircle $O_1, \dots, O_p$, as well as the special alpha curves $ \alpha_1, \dots, \alpha_p$ that intersect $\beta_1, \dots, \beta_p$. Given a generator $\x \in \Ta \cap \Tb^{\zed(\orM)}$, let $\bar \x$ be the corresponding generator in $\bar G$ obtained from $\x$ by deleting its components on $\alpha_i$, for $i=1, \dots, p$. Define $I(\x)$ to be the Maslov index of $\bar \x$ in $\bar G$. Then $I$ defines a filtration on the handleslid complex. Indeed, if $\y$ contributes a nonzero term to the differential $\del \x$, we must have $I(\x) \geq I(\y)$, with equality if and only if $\bar \x = \bar \y$.

Thus, the associated graded of the handleslid complex $\Chain^-(\Ta, \Tb^{\zed(\orM)}, \psi^{\orM}(\s))$ splits as a direct sum of
terms $C(\bar \x)$, where $C(\bar \x)$ is generated by all possible $\x$ with the fixed reduction $\bar \x$ on $\bar G$. Pick such a direct summand $C(\bar \x)$. 

\begin{lemma}
\label {lemma:barx}
For the almost complex structure $J(T_1, \dots, T_k)$ chosen as before (with $T_1 \gg T_2 \gg \dots \gg T_k \gg 0$), the complex $C(\bar \x)$ is isomorphic to the tensor product of the complexes
$$ \Ring \xrightarrow{U_i - U_{i+1}} \Ring,$$
for $i=1, \dots, p$.
\end {lemma}

\begin {proof}
Recall that the intersection points between $\alpha_i$ and $\beta_i$ are denoted $x_i$ and $y_i$, compare Figure~\ref{fig:3rows}. The generators of $C(\bar \x)$ can be written as $\x^A$, where $A \subseteq \{ 1, \dots, p\}$ is the set of indices $i$ such that $x_i \in \x^A$. 

The contributions to the differential on $C(\bar \x)$ come from domains completely supported in the rows containing $O_1, \dots, O_{p+1}$. By our previous observations, if their contribution is nonzero, the domains must be of types II (b) of width one, III (b) of width one, III (a), IV (a), II (a), or IV (b). We already know that the first four cases in this series produce nonzero contributions mod $2$. More precisely, the first three cases give a term $(U_i - U_{i+1}) \x^{A \setminus \{i\}}$ in $\del \x^A, $ whenever $i \in A$. 

The challenge is to figure out the contributions from domains of type II (a) and IV (b). Lemma~\ref{lemma:degz}
 implies that nonzero contributions can come only from domains that have no local multiplicities bigger than one. Hence, the differential on the complex $C(\bar \x)$ takes the form
 $$ \del \x^A = \sum_{i \in A} (U_i - U_{i+1}) \cdot \x^{A \setminus \{i\}} + \sum_{i \not 
 \in A} n^A_i \cdot \x^{A \cup \{i\} },$$
 for some values $n^A_i \in \{0, 1\}$. These values are constrained by the requirement that $\del^2 = 0$. In fact, it suffices to look at the coefficient of $\x^A$ in $\del^2 (\x^A):$
$$\sum_{i  \in A} n^{A \setminus \{i\}}_i \cdot (U_i - U_{i+1}) + \sum_{i \not \in A} n^A_i  \cdot (U_i - U_{i+1}) = 0.$$

Since the values $U_i - U_{i+1}$ for $i=1, \dots, p$, are linearly independent, we deduce that  $n^A_i = 0$ for all possible $A$ and $i \not \in A$. This completely determines the complex $C(\bar \x)$. 
\end {proof}
   
\begin{figure}
\begin{center}
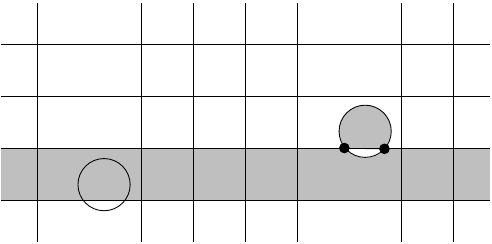
\end{center}
\caption {{\bf The origin of the complexes $\K(M)$.}  The figure shows
  part of a grid diagram with some arcs on the $\alpha$ and
  $\beta^{\{Z_1,Z_2,\dots, Z_k\}}$ curves drawn. There are two intersection
  points (marked as bullets) between the alpha curve below the marking
  $Z_1$, and the corresponding beta curve. There are two differentials
  going from the left to the right generator: a bigon containing $Z_1$
  and an annulus containing $Z_2$, both drawn shaded in the
  diagram. This produces a factor of the form $\Ring
  \xrightarrow{U_{j} - U_{j'} } \Ring$ in the definition of the
  complexes $\K(M)$.  }
\label{fig:annuli}
\end{figure}

\begin {proof}[Proof of Proposition~\ref{lemma:pear}]   
For our choice of almost complex structure, all the domains that contribute to the differential on the handleslid complex are either of type I (a), or appear in the differential on a summand $C(\bar \x)$ of an associated graded. 

The holomorphic count for all domains of type I (a) is always $\pm 1$. As previously noted, such domains are in one-to-one correspondence with empty rectangles on the quasi-destabilized grid $G^{L-M}$. The result now follows by combining this observation with Lemma~\ref{lemma:barx}. See Figure~\ref{fig:annuli}.
  \end {proof}

\subsection {The surgery theorem}
\label {sec:surgery}

By composing the maps
\eqref{eq:Proj}, \eqref{eq:dezed} and \eqref{eq:psim}, we construct a
map
\begin {align}
\Phi_\s^{\orM} &: \Chain^- (G, \s) \longrightarrow \Chain^-(G^{L-M}, \psi^{\orM}(\s))[[\{U_{i,j} \}_{L_i \subseteq M }]] \otimes_{\Ring} \K(M), 
\notag \\
\Phi_\s^{\orM} &= \Psi_{\psi^{\orM}(\s)}^{\orM} \circ \hat\De_{p^{\orM}(\s)}^{\orM} \circ \Pr^{\orM}_\s ,
\label{eq:phis2}\end {align}
defined for any $\s \in  \H(L)$.

This is the analogue of the map $\Phi_\s^{\orM}$ from \eqref{eq:Phi}. There is also an analogue of Proposition~\ref{prop:phiphi}: 

\begin {proposition}
For any $\orM$ and $\s \in \bH(L)$, we have 
\begin {equation}
\label {eq:PhiPhi}
 \sum_{\orM_1 \amalg \orM_2 = \orM} \Phi^{\orM_2}_{\psi^{\orM_1} (\s)} \circ \Phi^{\orM_1}_{\s} = 0,
 \end {equation} 
where $\orM_1$ and $\orM_2$ are only considered with the orientations induced from $\orM$. 
\end {proposition}

\begin{proof} This is similar to the proof of Proposition~\ref{prop:phiphi}, except we have the additional factors $\Psi_{\psi^{\orM}(\s)}^{\orM}$ in \eqref{eq:phis2}. These factors commute with the descent maps $\hat\De$ corresponding to disjoint sublinks. Indeed, if $\orM_1 \cap \orM_2 = \emptyset$, then Proposition~\ref{prop:polydegen} implies that 
$$ \Psi^{\orM_1}_{\s'} \circ \hat \De^{\orM_2}_{\s'}  = \hat \De^{\orM_2}_{\s'} \circ  \Psi^{\orM_1}_{\s'},$$
where $\s' = p^{\orM_2}(\psi^{\orM_1}(\s)).$

The maps $\Psi^{\orM_1}$ also commute with the inclusion maps $\Pr^{\orM_2}$ corresponding to disjoint sublinks. With these observations in mind, the proof of Proposition~\ref{prop:phiphi} can be easily adapted to the present context. 
\end {proof}

Using the new maps $\Phi^{\orM}_\s$, we can exactly mimic Section~\ref{subsec:surgery}, and state a version of the surgery theorem for grid diagrams instead of complete systems of hyperboxes.

We transfer most of the notation from Section~\ref{subsec:surgery}. In particular, we have a framing $\Lambda$ for the link $\orL$. For a sublink $N \subseteq L$, recall that $\Omega(N)$ is the set of all
possible orientations on $N$. For $\orN \in \Omega(N)$, the quantity $ \Lambda_{\orL, \orN} \in H_1(S^3 - L)$ is the sum of the framings of those components of $N$ oriented differently in $\orN$ as in $\orL$.

We consider the $\Ring$-module
$$ \C^-(G, \Lambda) = \bigoplus_{M \subseteq L} \prod_{\s \in \H(L)} \Bigl( \Chain^-(G^{L - M}, \psi^{M}(\s))[[\{ U_{i,j}\}_{ L_i \subseteq M}]] \Bigr)  \otimes_\Ring \K(M),$$
where $\psi^{M}$ simply means $\psi^{\orM}$ with $\orM$ being the orientation induced from the one on $\orL$.

We equip $\C^-(G, \Lambda)$ with a boundary operator $\D^-$ as follows. 

For $\s \in \H(L)$ and $\x \in \bigl( \Chain^-(G^{L - M}, \psi^{M}(\s))[[\{ U_{i,j}\}_{ L_i
  \subseteq M}]] \bigr) \otimes_\Ring \K(M)$, we set
 \begin {align*} 
\D^-(\s, \x) &= \sum_{N \subseteq L - M} \sum_{\orN \in \Omega(N)} (\s + \Lambda_{\orL, \orN}, \Phi^{\orN}_\s(\x)) \\
&\in  \bigoplus_{N \subseteq L - M} \bigoplus_{\orN \in \Omega(N)} \Bigl( \Chain^-(G^{L-M-N}, \psi^{M \cup \orN} (\s)) [[\{ U_{i,j}\}_{ L_i \subseteq M \cup N}]]\Bigr)  \otimes_\Ring \K(M \cup N) \subset \C^-(G, \Lambda).
\end {align*}

According to \eqref{eq:PhiPhi}, $\C^-(G, \Lambda)$ is a
chain complex. Just like the complex $\C^-(\Hyper, \Lambda)$ from Section~\ref{subsec:surgery}, it splits into a direct product of complexes $\C^-(G, \Lambda, \ux)$, according to $\ux \in \spc(S^3_\Lambda(L)) \cong \H(L)/H(L, \Lambda)$.
The complexes $\C^-(G, \Lambda, \ux)$ admit relative $\zz/\delt(\ux)\zz$-gradings, constructed just as for $\C^-(\Hyper, \Lambda, \ux)$.

\begin {theorem}
\label {thm:Surgery}
Fix a grid diagram $G$ (with at least one free marking) for an oriented link $\orL$
in $S^3$, and fix a framing $\Lambda$ of $L$.  Then, for every $\ux \in \spc(S^3_\Lambda(L))$, we have an isomorphism of relatively graded $\Field[[U]]$-modules: 
\begin {equation}
\label {eq:surgery}
 H_*(\C^-(G, \Lambda, \ux), \D^-) \cong  \HFm_{*}(S^3_\Lambda(L), \ux).
  \end {equation}
\end {theorem}

\subsection {A complete system associated to the grid}
\label {sec:completegrid}

Theorem~\ref{thm:Surgery} is a consequence of Theorem~\ref{thm:FirstSurgery}, albeit not an immediate one. A first guess would be to construct a complete system of hyperboxes $\Hyper_G$ out of a grid diagram $G$ (with at least one free marking) as follows: as noted in Section~\ref{sec:desublink2}, the grid $G$ gives rise to hyperboxes $\Hyper^{\orL, \orM}_G$ for the pairs $(\orL, \orM)$. One is tempted to try to include these into a complete system. Unfortunately, this is not possible: in a complete system, the diagrams $\Hyper^{\orL, \orM}(M)$ and $\Hyper^{\orL, - \orM}(M)$ are required to be surface isotopic, 
while for a grid $G$, typically the diagrams $\Hyper^{\orL, \orM}_{G, \dd^{M}}$ and $\Hyper^{\orL, -\orM}_{G, \dd^{M}}$ are not surface isotopic. For example, if $\orM$ has the orientation induced from $\orL$, then in the diagram $\Hyper^{\orL, \orM}_G$ the new beta curves are obtained by handlesliding over the $O$ markings on $M$, whereas in  $\Hyper^{\orL, -\orM}_G$ they are obtained by handlesliding over the $X$ markings on $M$. 

Nevertheless, we can still associate to a grid $G$ a complete system of hyperboxes $\BHyper_G$, such that the (compressed) maps on Floer complexes coming from $\BHyper_G^{\orL, \orM}$ are essentially the same as those coming from $\Hyper^{\orL, \orM}_G$. 

In fact, if $\orM$ is the orientation induced from $\orL$, we define $\BHyper_G^{\orL, \orM}$ to be exactly the same as $\Hyper^{\orL, \orM}_G$. On the other hand, if $\orM$ has some components with a different orientation, we need to change $\Hyper^{\orL, \orM}_G$ using some more beta handleslides.

\begin{figure}
\begin{center}
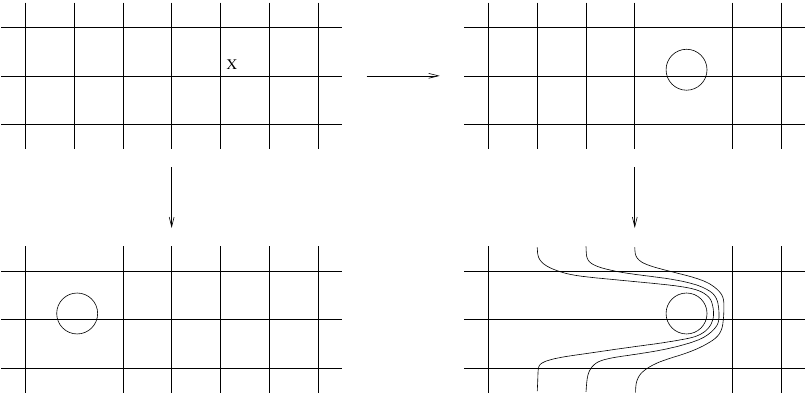
\end{center}
\caption {{\bf Relating the handleslides over $O$ and over $X$.}
The top left picture shows a part of a grid diagram, with an $O$ and $X$ markings in the same row. The vertical arrow on the left describes a handleslide over the $O$. The top arrow describes a handleslide over the $X$, as well as relabeling the $X$ as an $O$. The vertical arrow on the right shows a sequence of handleslides, the end result being a diagram surface isotopic with the one on the bottom left. 
}
\label{fig:ox}
\end{figure}

Before constructing $\BHyper^{\orL, \orM}_G$ in general, let us make an observation. Let $O_j, X_j$ be two markings on the same row in the grid  $G$. Then the two collections of beta curves $\betas^{O_j}$ and $\betas^{X_j}$ are strongly equivalent, being obtained from each other through handleslides of the beta curves between $O_j$ and $X_j$, as in Figure~\ref{fig:ox}. 

\begin{figure}
\begin{center}
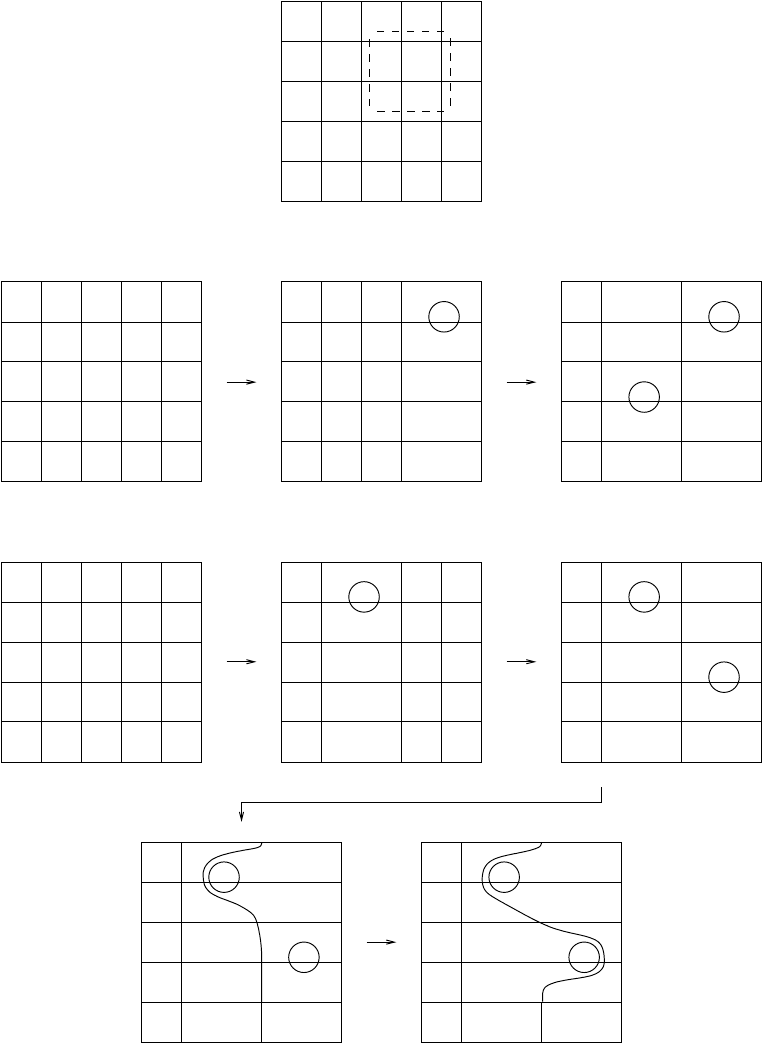
\end{center}
\caption {{\bf Some hyperboxes in a complete system for the Hopf link.}
We show a grid diagram $G$ for the Hopf link $\orL$ (with one free marking) and two one-dimensional hyperboxes in the complete system $\BHyper_G$ associated to the grid. At the top we have the grid $G$, with one component $L_1$ shown by the dashed lines. The second row is the hyperbox $\BHyper_G^{\orL, \orL_1}$ of size $(2)$, where $\orL_1$ has  the orientation induced from $\orL$. The two bottom rows show the hyperbox  $\BHyper_G^{\orL, -\orL_1}$, of size $(4)$.  We describe each hyperbox as a sequence of Heegaard diagrams, from the initial to the final vertex. The two initial diagrams are obtained from the grid $G$ by reduction at $\orL_1$, resp. $-\orL_1$. The two final diagrams are surface isotopic.
}
\label{fig:hopfgrid}
\end{figure}

For simplicity, let us first consider the case when $\orM = -\orL_i$ is a single component of $L$, oriented oppositely from its orientation in $\orL$. Then $\Hyper^{\orL, -L_i}_G$ is a one-dimensional hyperbox, consisting of a sequence of moves from the reduction $r_{- L_i}(G)$ to the diagram $\Hyper^{\orL, -\orL_i}_{G, \dd^{L_i}}$ obtained from it by handleslides over all $X_j \in \Xs_i$. In the new one-dimensional hyperbox $\BHyper^{\orL, -\orL_i}_G$, we add at the end of this sequence a series of handleslides as in Figure~\ref{fig:ox}, so that we end at the diagram $\Hyper^{\orL, \orL_i}_{G, \dd^{L_i}}$, which is $r_{\orL_i}(G)$ handleslid over all $O_{i,j} \in \Os_i$. See Figure~\ref{fig:hopfgrid} for an example. There are unique choices for the intersection points that play the role of $\Theta$-chain elements. 

More generally, for $L_i \subset L'$, we define the one-dimensional hyperboxes $\BHyper^{\orL', \orL_i}_G$ and $\BHyper^{\orL', \orL_i}_G$ to consist of similar moves as in the case $L' = L$, except starting at the diagram handleslid over all the $O$ markings that are on the components of $L-L'$.
 
Next, to define the higher-dimensional hyperboxes $\bar \Hyper^{\orL', \orM}_G$ for all $L'$ and $\orM$, note that we already know the sequence of moves on their edges. The moves in one direction of the hyperbox are independent of those in a different direction, so it is straightforward to fill in the hyperbox by combining these moves. The $\Theta$-chain elements are taken to be zero on higher dimensional faces, just as we did for the hyperboxes $\Hyper^{\orL, \orM}_G$ from Section~\ref{sec:desublink2}.  The result is a complete system of hyperboxes $\BHyper_G$ for $\orL$. Indeed, for the associated good set of trajectories, we can take horizontal paths on the grid joining each $O$ marking to the $X$ marking in the same row.

\begin {proof}[Proof of Theorem~\ref{thm:Surgery}]
In light of Theorem~\ref{thm:FirstSurgery}, it suffices to show that the chain complexes $\C^-(G, \Lambda, \ux)$ and $\C^-(\BHyper_G, \Lambda, \ux)$ are chain homotopy equivalent. Indeed, for $\s \in \bH(L)$, let us look at the hyperbox of Floer complexes $\Cc(\BHyper_{G}^{\orL, \orM}, \psi^{\orM}(\s))$ associated to a hyperbox $\BHyper_G^{\orL, \orM}$ in $\BHyper_G$, as in Section~\ref{sec:gens}. This contains as a sub-hyperbox $\Cc(\Hyper_{G}^{\orL, \orM}, \psi^{\orM}(\s))$. 

In fact, $\Cc(\BHyper_{G}^{\orL, \orM}, \psi^{\orM}(\s))$ is obtained from $\Cc(\Hyper_{G}^{\orL, \orM}, \psi^{\orM}(\s))$ by adding polygon maps corresponding to handleslides over marked disks as on the right of Figure~\ref{fig:ox}. According to Proposition~\ref{prop:polydegen}, if we choose suitable almost complex structures, these extra higher polygon maps are the same as the corresponding polygon maps in the quasi-destabilized diagrams. In the quasi-destabilized diagrams, the handleslides from 
the right hand side of Figure~\ref{fig:ox} are nothing more than curve isotopies. On the level of hyperboxes of Floer complexes, isotopies have (up to chain homotopy) the effect of identity shifts---compare Lemma~\ref{lemma:id} and the discussion of elementary enlargements in Section~\ref{sec:movesh}. It follows that $\Cc(\BHyper_{G}^{\orL, \orM}, \psi^{\orM}(\s))$ is obtained from $\Cc(\Hyper_{G}^{\orL, \orM}, \psi^{\orM}(\s))$ by elementary enlargements, in the sense of Section~\ref{sec:ele}. Elementary enlargments leave unchanged the corresponding compressed hypercubes, see Lemma~\ref{lemma:ci}. Hence, the maps involved in the complexes $\C^-(G, \Lambda, \ux)$ and $\C^-(\BHyper_G, \Lambda, \ux)$ are the same, up to compatible chain homotopy equivalences.
\end {proof}

\begin {remark}
Theorems~\ref{thm:Cobordisms}, \ref{thm:AllVersions} and \ref{thm:Mixed} also admit straightforward adaptations to the case of grids, with the respective complexes all being expressed in terms of hyperboxes $\Hyper^{\orL, \orM}_G$. The exact statements are given in \cite{MOT}.
\end {remark}

\newpage
\bibliographystyle{custom}
\bibliography{biblio}

\end{document}